\documentclass{pspum-l}
\usepackage[T1]{fontenc}
\usepackage[utf8]{inputenc}

\usepackage{bm}

 \usepackage[style=ext-alphabetic,doi=false,isbn=false,maxalphanames=99,maxnames=99,articlein=false,backend=biber]{biblatex}
\DeclareFieldFormat{language}{}
\renewbibmacro*{language}{}

\AtEveryBibitem{%
  \clearfield{language}%
  \ifboolexpr{
    test {\ifentrytype{article}} 
      or test {\ifentrytype{book}}
            or test {\ifentrytype{incollection}}
  }{%
    \clearfield{url}%
  }{}%
  \iffieldundef{eprint}
    {}
    {\clearfield{url}}
} 
  \DefineBibliographyStrings{english}{page = {}, pages = {}}
\DefineBibliographyStrings{french}{page = {}, pages = {}}
\DefineBibliographyStrings{german}{page = {}, pages = {}}

\DeclareMathOperator{\CoCh}{CoCh}
\newcommand{\unif}{\varpi}
\newcommand{\alg}{\mathrm{alg}}

\newcommand{\smG}{\mathrm{sm}\,G}
\newcommand{\smH}{\mathrm{sm}\,H}

\newcommand{\smcg}{\mathrm{c.g. sm}}
\newcommand{\Ghat}{\widehat{G}}
\newcommand{\unip}{\operatorname{unip}}
\newcommand{\Sel}{\operatorname{Sel}}

\usepackage{etex}
\usepackage{fixltx2e}
\usepackage{tikz}
\usepackage{tikz-cd}
\usepackage{adjustbox}
\usepackage{enumitem}
\usepackage{amsmath,calligra,mathrsfs}

\newcommand{\sRHom}{R\mathcal{H}om}
\newcommand{\sHom}{\mathcal{H}om}

\newcommand{\lambdau}{\underline{\lambda}}
\newcommand{\crys}{\mathrm{crys}}

\newcommand{\semis}{{\operatorname{ss}}}

\newcommand{\chibar}{\overline{\chi}}

\newcommand{\thetabar}{\overline{\theta}}

\newcommand{\sigmacomp}{\sigma^{\operatorname{co}}}
\newcommand{\sigmasigmacomp}{\sigma|\sigmacomp}
\newcommand{\Ad}{\mathrm{Ad}}
\newcommand{\ad}{\mathrm{ad}}
\newcommand{\T}{\mathbf{T}}
\newcommand{\tame}{\mathrm{tame}}
\newcommand{\Pro}{\operatorname{Pro}}
\newcommand{\IndCoh}{\operatorname{IndCoh}}

\newcommand{\cg}{\operatorname{c.g.}}
\newcommand{\Dfp}{\mathop{D}_{\mathrm{f.p.}}\nolimits}
\newcommand{\DfpI}{\mathop{D}_{\mathrm{f.p.},I}\nolimits}
\newcommand{\qc}{\mathrm{qc}}
\newcommand{\coh}{\mathrm{coh}}
\newcommand{\Pj}{\mathbf{P}}

\newcommand{\fA}{\mathfrak{A}}
\newcommand{\PS}{\operatorname{PS}}
\newcommand{\tP}{\widetilde{P}}
\newcommand{\tH}{\widetilde{{H}}}
\newcommand{\cont}{\mathrm{cont}}

\newcommand{\gA}{{\mathfrak{A}}}
\newcommand{\odd}{\mathrm{odd}}

\DeclareMathOperator{\Bun}{Bun}
\newcommand{\ps}{\operatorname{ps}}
\newcommand{\Pone}{\mathbf{P}^1}

\newcommand{\Zptimes}{\Zp^\times}

\newcommand{\Qptimes}{\Qp^\times}

\newcommand{\Iw}{I}

\usepackage{amssymb,amsmath,amsfonts,amsthm,epsfig,amscd,amsxtra,latexsym,euscript}
\usepackage{stmaryrd}
\usepackage[all,cmtip,poly]{xy}
\usepackage{color}
\usepackage{mathtools}

\newcommand{\sm}{{\mathrm{sm}}}
\newcommand{\triv}{\mathrm{triv}}

\newcommand{\pro}{\operatorname{Pro}}

\newcommand{\red}{\operatorname{red}}
\newcommand{\quoteslim}[1]{\text{``}\varprojlim_{#1}\text{''}}
\newcommand{\colim}{\operatorname{colim}}

\newcommand{\fB}{\mathfrak{B}}

\newcommand{\set}[1]{\left\{#1\right\}}

\newcommand{\univ}{{\operatorname{univ}}}
\newcommand{\To}{\longrightarrow}
\newcommand{\isoto}{\stackrel{\sim}{\To}}
\newcommand{\cotimes}{\, \widehat{\otimes}}

\newcommand{\pr}{\operatorname{pr}}
\newcommand{\fg}{\operatorname{fg}}
\newcommand{\fgadm}{\operatorname{fg.adm}}

\newcommand{\Rep}{\operatorname{Rep}}

\newcommand{\JH}{\operatorname{JH}}
\newcommand{\sigmabar   }{\overline{\sigma}}
\newcommand{\Vbar   }{\overline{V}}   
\newcommand{\iso}{\buildrel \sim \over \longrightarrow}

\newcommand{\id}{\operatorname{id}}

\newtheorem{theorem}[subsubsection]{Theorem}
\newtheorem{thm}[subsubsection]{Theorem}
\newtheorem{expectedthm}[subsubsection]{Expected Theorem}
\newtheorem{lemma}[subsubsection]{Lemma}
\newtheorem{lem}[subsubsection]{Lemma}

\newtheorem{cor}[subsubsection]{Corollary}
\newtheorem{conj}[subsubsection]{Conjecture}

\newtheorem{expectation}[subsubsection]{Expectation}
\newtheorem{conc}[subsubsection]{Conclusion}

\newtheorem{prop}[subsubsection]{Proposition}

\theoremstyle{definition}
\newtheorem{df}[subsubsection]{Definition}
\newtheorem{defn}[subsubsection]{Definition}

\theoremstyle{remark}
\newtheorem{remark}[subsubsection]{Remark}
\newtheorem{rem}[subsubsection]{Remark}

\newtheorem{example}[subsubsection]{Example}

\newcommand{\numequation}{\addtocounter{subsubsection}{1}\begin{equation}}
  \newcommand{\nummultline}{\addtocounter{subsubsection}{1}\begin{multline}}
    \newcommand{\numalign}{\addtocounter{subsubsection}{1}\begin{align}}
\newcommand{\anumequation}{\addtocounter{subsection}{1}\begin{equation}}
\newcommand{\anummultline}{\addtocounter{subsection}{1}\begin{multline}}
  \let\originaltheequation\theequation
\renewcommand{\theequation}{\arabic{section}.\arabic{subsection}.\arabic{subsubsection}}

\newcommand{\A}{\mathbf{A}}
\newcommand{\C}{\CC}
\newcommand{\F}{\FF}

\newcommand{\Q}{\QQ}

\newcommand{\Z}{\ZZ}

\newcommand{\D}{\cD}

\newcommand{\m}{\frakm}
\newcommand{\p}{\frakp}

\newcommand{\X}{\frakX}

\newcommand{\CC}{{\mathbf C}}

\newcommand{\FF}{{\mathbf F}}
\newcommand{\GG}{{\mathbf G}}

\newcommand{\NN}{{\mathbf N}}

\newcommand{\QQ}{{\mathbf Q}}

\newcommand{\ZZ}{{\mathbf Z}}

\newcommand{\bA}{\ensuremath{\mathbf{A}}}

\newcommand{\bQ}{\ensuremath{\mathbf{Q}}}

\newcommand{\cA}{{\mathcal A}}
\newcommand{\cB}{{\mathcal B}}
\newcommand{\cC}{{\mathcal C}}
\newcommand{\cD}{{\mathcal D}}
\newcommand{\cE}{{\mathcal E}}
\newcommand{\cF}{{\mathcal F}}
\newcommand{\cG}{{\mathcal G}}
\newcommand{\cH}{{\mathcal H}}
\newcommand{\cI}{{\mathcal I}}

\newcommand{\cK}{{\mathcal K}}
\newcommand{\cL}{{\mathcal L}}
\newcommand{\CL}{{\mathcal{L}}}
\newcommand{\cM}{{\mathcal M}}

\newcommand{\cO}{{\mathcal O}}

\newcommand{\cR}{{\mathcal R}}

\newcommand{\cU}{{\mathcal U}}
\newcommand{\cV}{{\mathcal V}}
\newcommand{\cW}{{\mathcal W}}
\newcommand{\cX}{{\mathcal X}}
\newcommand{\cY}{{\mathcal Y}}
\newcommand{\cZ}{{\mathcal Z}}
\newcommand{\tF}{{\widetilde F}}

\newcommand{\frakm}{\mathfrak{m}}

\newcommand{\frakp}{\mathfrak{p}}

\newcommand{\frakX}{\mathfrak{X}}

\newcommand{\Fbar}{\overline{\F}}
\newcommand{\Qbar}{\overline{\Q}}

\newcommand{\Fp}{\F_p}
\newcommand{\Fq}{\F_q}

\newcommand{\Fpbar}{\Fbar_p}

\newcommand{\Fpbartimes}{\Fpbar^{\times}}

\newcommand{\Zp}{\Z_p}

\newcommand{\Qp}{\Q_p}

\newcommand{\Qpbar}{\Qbar_p}

\DeclareMathOperator{\Aut}{Aut}
\DeclareMathOperator{\Coh}{Coh}
\DeclareMathOperator{\QCoh}{QCoh}
\DeclareMathOperator{\LocSys}{LocSys}

\DeclareMathOperator{\End}{End}
\DeclareMathOperator{\REnd}{REnd}
\DeclareMathOperator{\Ext}{Ext}

\DeclareMathOperator{\Fil}{Fil}
\DeclareMathOperator{\Fun}{Fun}
\DeclareMathOperator{\Gal}{Gal}
\DeclareMathOperator{\GL}{GL}

\DeclareMathOperator{\Gr}{Gr}

\DeclareMathOperator{\Hom}{Hom}
\DeclareMathOperator{\RHom}{RHom}

\DeclareMathOperator{\Maps}{Maps}

\DeclareMathOperator{\im}{im}
\DeclareMathOperator{\Ind}{Ind}

\DeclareMathOperator{\Mod}{Mod}
\DeclareMathOperator{\LMod}{LMod}
\DeclareMathOperator{\RMod}{RMod}
\DeclareMathOperator{\Mor}{Mor}

\DeclareMathOperator{\PGL}{PGL}
\DeclareMathOperator{\Pic}{Pic}

\DeclareMathOperator{\SL}{SL}

\DeclareMathOperator{\Spec}{Spec}

\DeclareMathOperator{\Spf}{Spf}

\DeclareMathOperator{\Sym}{Sym}

\DeclareMathOperator{\WD}{WD}

\newcommand\cInd{c\text{-}\!\Ind}
\newcommand{\ab}{\mathrm{ab}}

\newcommand{\Fr}{\mathrm{Fr}}
\newcommand{\Frob}{\mathrm{Frob}}

\newcommand{\HT}{\mathrm{HT}}
\newcommand{\Id}{\mathrm{Id}}
\newcommand{\ladm}{\mathrm{l.adm}}
\newcommand{\nr}{\mathrm{nr}}

\newcommand{\St}{\mathrm{St}}

\newcommand{\ur}{\mathrm{ur}}

\newcommand{\rhobar}{\overline{\rho}}

\newcommand{\et}{\text{\'{e}t}}

\newcommand{\op}{\mathrm{op}}

\newcommand{\into}{\hookrightarrow}

\newcommand{\Ga}{\GG_a}
\newcommand{\Gm}{\GG_m}
\newcommand{\Gmhat}{\widehat{\GG}_m}

\newcommand{\etabar}{\overline\eta}

\newcommand{\varepsilonbar}{\overline{\varepsilon}}
\newcommand{\rbar}{\overline{r}}
\newcommand{\Res}{\operatorname{Res}}

\usepackage[bookmarksopen,bookmarksdepth=2]{hyperref}
\usepackage{xurl} %

\begin{document}

\title[Categorical $p$-adic Langlands] 
{%
An introduction to the categorical $p$-adic
  Langlands program}

\author[M. Emerton]{Matthew Emerton}\email{emerton@uchicago.edu}
\address{Department of Mathematics, University of Chicago,
  5734 S.\ University Ave., Chicago, IL 60637, USA}
\author[T. Gee]{Toby Gee} \email{toby.gee@imperial.ac.uk} \address{Department of
  Mathematics, Imperial College London,
  London SW7 2AZ,~UK}
\author[E. Hellmann]{Eugen Hellmann}\email{e.hellmann@uni-muenster.de}
\address{Mathematisches Institut, Westf\"alische Wilhelms Universit\"at M\"unster, Einsteinstra\ss e 62, 48149 M\"unster,~Germany}
\thanks{M.E.\ was supported in part by the
  NSF grants 
  DMS-1902307, DMS-1952705, and DMS-2201242. T.G.\ was 
  supported in part by an ERC Advanced grant. This project has received funding from the European Research Council (ERC) under the European Union's Horizon 2020 research and innovation programme
  (grant agreement No. 884596). M.E.\ and T.G.\ were
  both supported in part by the Simons Collaboration on Perfection in Algebra, Geometry and Topology. E.H. was supported by Germany's Excellence Strategy EXC 2044-390685587 ``Mathematics M\"unster:
Dynamics--Geometry--Structure" and by the CRC 1442 Geometry: Deformations and Rigidity of the DFG}
\subjclass[2020]{Primary 11R39}

\begin{abstract}We give an introduction to the ``categorical'' approach to the
  $p$-adic Langlands program, in both the ``Banach'' and ``analytic'' settings. 
\end{abstract}

\maketitle

\setcounter{tocdepth}{2}
\tableofcontents

\section{Introduction}
\label{sec: introduction}
The aim of these notes is to discuss some of the $p$-adic aspects of the 
Langlands program, and especially the emerging  ``categorical''
perspective. We begin with an overview, first of some history and then of our goals;
we refer the reader to Section~\ref{sec:notation} for any unfamiliar notation.

\subsection{A rapid overview of the Langlands program}
\label{sec: a rapid overview of the Langlands program}
The Langlands program began with Langlands's celebrated letter to
Weil (reproduced e.g.\ in~\cite[\S 6]{MR3309997}),
the contents of which were elaborated on in several subsequent writings of
Langlands, including his article ``Problems in the theory
of automorphic forms''~\cite{MR0302614} and
his Yale lectures on Euler products~\cite{MR0419366}.
The article focuses on the {\em functoriality} conjecture for automorphic
representations and its consequences, 
while the Yale lectures explain his construction of automorphic
$L$-functions, and discuss the problem of proving their expected
analytic properties (analytic continuation and functional equation).\footnote{Since
the standard $L$-functions associated to automorphic representations
on $\GL_d$ coincide with the $L$-functions constructed by Godement--Jacquet \cite{MR342495} and 
Tamagawa \cite{MR144928}, which are known to admit analytic continuations and functional equations,
the expected properties would follow from functoriality.   Langlands's Yale lectures
discuss another approach, via constant terms of Eisenstein series.}

In his Yale lectures Langlands already touches, if only tangentially,  on the idea
that the automorphic $L$-functions may include all motivic $L$-functions.  
Subsequently, this idea developed into the {\em reciprocity} conjecture.
(See~\cite{MR546619} for one articulation of this conjecture.)
Roughly speaking, the reciprocity conjecture articulates a correspondence between
certain automorphic representations (those which are {\em algebraic})
of $G(\A_F)$ (for a connected reductive algebraic group
$G$ over a number field $F$),
and motives (with coefficients in $\Qbar$, say) over $F$ whose motivic
Galois group is closely related to to the $C$-group\footnote{The $C$-group,
introduced in~\cite{MR3444225}, 
is a refinement  of the $L$-group introduced by Langlands, which is
better adapted to the problem of relating  automorphic forms  and Galois
representations;  we %
recall the  definition below.}
 ${}^cG$ of $G$.  %
Conjectures on Galois representations (especially the Fontaine--Mazur
conjecture~\cite{MR1363495}) suggest that such motives in turn may
be identified with  compatible systems %
of $\ell$-adic %
Galois representations $\Gal_F \to {}^cG(\Qbar_{\ell})$. \footnote{As far as we know, such an expectation can only be made precise for~$G=\GL_n$, as it is unclear what the precise definition of a
  compatible system of Galois representations should be for a general~$G$; see e.g.\ \cite[\S 6]{MR4018263}.} %

There are many subtleties involved in trying to formulate a precise reciprocity
conjecture.  For example,
in the case that $G$ is not some  $\GL_d$, one has to
worry about $L$-packets; and one should restrict to automorphic representations
which are not {\em anomalous} in the sense of~\cite{MR546619} 
(or else replace $\Gal_F$  by  $\Gal_F \times \SL_2$ and work
in the framework of Arthur's conjectures regarding non-tempered
endoscopy~\cite{MR1021499}).
We refer to \cite{MR3444225} for a more technical discussion
of the conjecture, and to \cite{MR4274535} for a more thorough historical overview.

We note that there is a relationship between reciprocity and functoriality:
namely, since functoriality is tautologically true for (compatible
systems of) Galois representations, cases of reciprocity can be used
to deduce cases of functoriality.  Also, since $L$-groups and $C$-groups involve
Galois groups in their definition, representations of Galois groups
into algebraic groups can sometimes be related to $L$-homomorphisms
of $L$-groups  or  $C$-homomorphisms of $C$-groups, and thus some
cases of reciprocity can be subsumed into functoriality.  (This is
Langlands's original perspective on the Artin conjecture, as explained
in~\cite{MR0302614}.)   If one replaces the Galois group $\Gal_F$  by the hypothetical
Langlands group~$L_F$, and compatible systems of $\Gal_F$-representations
by representations of~$L_F$, then one can also extend the reciprocity conjecture
to non-algebraic automorphic representations, or (using the ``$L_F$-form'' of
the $L$-group or $C$-group) entirely subsume reciprocity into functoriality.
(A version of this last-mentioned  perspective was already adopted at times by
Langlands, by using the ``Weil  group form'' of the $L$-group.)

\subsection{The \texorpdfstring{$p$}{p}-adic perspective}
\label{sec: the p-adic perspective}
Our intention in these notes is not so much to focus on reciprocity in the manner
described above (the relationship between automorphic representations
and motives, or compatible systems of Galois representations),
but rather to fix a prime $p$ and consider the relationship
between automorphic representations and $p$-adic Galois
representations (for this fixed choice of~$p$).   At first,
this doesn't much change the problem, since (at least for
representations valued in $\GL_d$) a compatible system
of semisimple Galois representations is determined by any one of
its members, %
and the Fontaine--Mazur conjecture gives a purely Galois-theoretic
condition for a $p$-adic Galois representation to be motivic. 

But focusing on a particular  prime $p$ brings to the fore certain
aspects of the theory of automorphic forms and Galois representations
which are absent in the more motivic perspective on reciprocity.
For example, since at least the work of Ramanujan, it has been
known that automorphic forms can satisfy interesting congruences
modulo powers of~$p$.  Since the work of Swinnerton--Dyer~\cite{MR0406931}
and Serre~\cite{MR0466020}, %
it has been understood that these congruences are related to (or are
manifestations of, if one prefers) analogous
congruences between $p$-adic Galois representations.
Extending the notion of congruences between automorphic forms,
one is led to take various $p$-adic completions of spaces
of automorphic forms to obtain the notion of $p$-adic automorphic forms,
with associated Galois representations that need not be motivic.
Related to this, one has notions of mod $p$ automorphic forms, to which one
might associate mod $p$ Galois representations; and reciprocity
conjectures have been formulated in this context, famously
by Serre~ \cite{SerreDuke} in the context of classical modular forms,
and more recently in greater generality by others
(e.g.~\cite{BuzzardDiamondJarvis}). %

We do not intend at all to survey these developments; rather, our aim here 
is simply to indicate that the $p$-adic perspective emerged naturally over
a long period of time, and has led to a natural collection
of problems and concerns:  e.g.\ how  to arrive at {\em some} notion
of automorphic form which admits $p$-adic integral or mod $p$ coefficients,
and which allows for genuinely $p$-adic objects?
And how to phrase reciprocity in a manner which allows for $p$-adic integral
or mod $p$ Galois representations, which will be associated to the $p$-adic integral
or  mod $p$ objects  that  one introduces on the automorphic side?

We make one last remark on the $p$-adic theory for now:  once one allows
oneself to $p$-adically interpolate automorphic forms, one sees
that automorphic forms naturally lie in families (e.g.\ Hida families,
Coleman families, \ldots; see \cite{MR2906349} for a survey).   And
Galois representations also lie in families (e.g.\ via Mazur's deformation
theory~\cite{MR1012172}, which he was at least partly motivated to develop in
response to Hida's theory of ordinary families of $p$-adic modular forms).
These phenomena of continuous families of objects
--- on both the automorphic and Galois side ---
are not a feature of reciprocity in its more motivic formulation.
(Cuspidal automorphic representations are rigid objects, and so  are motives,
if one has fixed a particular number field as the field of definition.)
Before passing to our next topic, we note that families of Galois
representations {\em are} a feature of the {\em geometric} Langlands program.
We will return to this point below.

\subsection{Proofs of reciprocity}%
\label{sec: proofs of reciprocity}
Some of the very first results on reciprocity were proved in the case
of representations with solvable image, by combining class field theory
(which implies reciprocity in the abelian context) 
with techniques more classically automorphic in nature (see in
particular the work of Langlands and Tunnell, \cite{LanglandsGL2, Tunnell}).   %
In~1995, 
a breakthrough in our understanding of reciprocity was achieved with
the proof by Wiles~\cite{WilesFLT} and Taylor--Wiles~\cite{MR1333036} of the modularity of semistable
elliptic curves over~$\Q$, a result which was soon improved to handle
all elliptic curves over~$\Q$~\cite{BCDT}.
Building on these methods, reciprocity has since been proved in many
further interesting contexts; see Calegari's recent
survey~\cite{calegari2021reciprocity} for some of the highlights of
the last 25 years.

The crux of the method of~\cite{WilesFLT, MR1333036} and the many
subsequent results that build on their ideas is to prove the
modularity --- or, more generally, automorphy --- of a $p$-adic Galois
representation. More precisely, the Taylor--Wiles method gives a way
to deduce from the modularity of a single $p$-adic Galois
representation~$\rho_1$ the modularity of another $p$-adic Galois
representation~$\rho_2$ which is congruent to~$\rho_1$ modulo~$p$;
this is typically referred to as {\em modularity lifting}.  (We ``lift''
modularity from the mod $p$ reduction of~$\rho_1$ --- which, since
it coincides with the mod $p$ reduction of~$\rho_2$, is inherited from
that of~$\rho_2$ --- to the the modularity of the $p$-adic 
representation $\rho_1$ itself.)
It turns out that there are many such congruences between Galois
representations, and modularity can be propagated from a single
representation to many others in this way (either using a fixed
prime~$p$, or using several different primes~$p$). Of course, it is
necessary to have an initial supply of representations whose
modularity is already known, such as CM forms, or 2-dimensional Artin
representations with solvable image (the latter being used in Wiles'
proof of Fermat's Last Theorem~\cite{WilesFLT}).

Modularity lifting theorems are also known as ``$R=\T$'' theorems,
because they are proved by identifying a $\Zp$-algebra $R$ (a ``Galois
deformation ring'') parameterizing Galois representations congruent to
a fixed representation $\rho_1$ with another $\Zp$-algebra $\T$, a
``Hecke algebra'', the endomorphism algebra of an appropriate space of ($p$-adic) modular
forms. Over the last 25 years it has become apparent that such
theorems should hold in great generality, although finding appropriate
definitions of $R$ and $\T$ so that we literally have $R=\T$ is still
something of an art. In particular, it has become clear (following in
particular Skinner--Wiles \cite{MR1793414},
  Wake--Wang-Erickson \cite{MR3738092}, and Newton--Thorne
  \cite{https://doi.org/10.48550/arxiv.1912.11265}) that in general
  $R$ should be taken to be a so-called pseudodeformation ring, i.e.\
  a deformation ring for pseudorepresentations, rather than a
  deformation ring for literal Galois representations.

  Roughly speaking, pseudorepresentations capture the information
  given by characteristic polynomials of representations; for this
  reason they are also known as pseudocharacters. They were originally
  introduced for~$\GL_2$ by Wiles~\cite{MR969243}, and were
  considerably developed for~$\GL_d$ by many authors (we highlight in
  particular Chenevier's notion of a determinant~\cite{MR3444227},
  which showed how to define them for arbitrary primes~$p$). A theory
  valid for general reductive groups was introduced by Vincent
  Lafforgue~\cite{MR3787407}.

  From Lafforgue's point of view, as
  recently made precise (in the setting of local Galois representations) by Fargues--Scholze~\cite[\S
  VIII]{fargues--scholze}, a pseudodeformation ring~$R$ is the ring of
  global functions on a moduli stack of $L$-parameters, %
  so that the
  moduli space of pseudorepresentations is a coarse moduli space\footnote{Here and below,
we use 
the expression ``coarse moduli space'' in its usual informal manner; since stacks of $L$-parameters
are not Deligne--Mumford stacks, they do not admit coarse moduli spaces
in the technical sense, but rather {\em adequate moduli spaces} in the
sense of~\cite{MR3272912}.}
 for
  the moduli stack of $L$-parameters; and for general reductive
  groups, the Hecke algebra~$\T$ is replaced with a larger algebra,
  the algebra of ``excursion operators'' (which agrees with the usual
  Hecke algebra for~$\GL_d$).
  
From this perspective, it is natural to wonder whether it is possible
to upgrade the putative equality $R=\T$ to a statement on the actual
moduli stack of $L$-parameters, rather than on its coarse moduli
space; for example, one could ask whether spaces of $p$-adic modular
forms admit an interpretation in terms of sheaves on the stack of
$L$-parameters. The possibility of such statements, and their local
analogues, is the main subject of these notes; but before discussing this possibility
further, we turn to the geometric Langlands program, which was
also a substantial motivation for Vincent Lafforgue's work.

\subsection{Unramified geometric Langlands and the Fargues--Scholze
  conjecture}\label{subsec: FS}Very roughly, the unramified geometric Langlands
correspondence for a connected reductive group~$G$ and a curve~$X$
takes the form of an equivalence of categories \numequation\label{eqn:
  rough form of geometric Langlands}D(\textrm{const.}
\Bun_G(X))\cong \QCoh(\LocSys_{\Ghat}(X)),\end{equation}where the
left hand side is a derived category of ``constructible sheaves'' on
the stack $\Bun_G$ of $G$-bundles on~$X$, and the right hand side is a
derived category of ``quasicoherent sheaves'' on a stack of
``$\Ghat$-local systems'' on~$X$; see for
example~\cite{MR3666025}.

The recent work of Fargues--Scholze \cite{fargues--scholze} geometrizes
the classical local Langlands correspondence by ``taking $X$ to be the
Fargues--Fontaine curve''. (One motivation for doing this is that the
fundamental group of the Fargues--Fontaine curve is~$\Gal_F$.)
Slightly more precisely (but still only approximately), if $F/\Qp$ is
a finite extension, and $X$ is the Fargues--Fontaine curve for~$F$,
then for each prime $\ell\ne p$, one expects an equivalence of the
approximate form~\eqref{eqn: rough form of geometric Langlands}, where
the left hand side is replaced by an appropriate derived category of
(solid) $\ell$-adic sheaves on $\Bun_G$, and the right hand side is a
derived category of (Ind-)coherent sheaves on a moduli stack of
$\Ghat$-valued $\ell$-adic representations of~$\Gal_F$.

The stack~$\Bun_G$ of $G$-bundles on the Fargues--Fontaine curve
admits a stratification into locally closed substacks, with one of the
open strata being a quotient stack $[\cdot/G(F)]$. The $\ell$-adic sheaves
on this quotient stack correspond to smooth $\ell$-adic
representations of $G(F)$, so by considering an appropriate
pushforward of sheaves from this quotient stack to~$\Bun_G$, an
equivalence of the form ~\eqref{eqn: rough form
  of geometric Langlands} implies in particular that there is a fully
faithful functor \numequation\label{eqn: conjectural fully faithful
  functor}D(\mathrm{sm.}\ G(F))\to
\QCoh(\LocSys_{\Ghat})\end{equation}where the left hand side is a
derived category of $\ell$-adic representations of $G(F)$, and the
right hand side is a derived category of quasicoherent sheaves on the
stack of $\Ghat$-valued $\ell$-adic representations of~$\Gal_F$.

Specialising to the case $G=\GL_d$, we have the following rough
conjecture, where we now allow the case $\ell=p$. %
\begin{conj}
\label{conj:main intro}
If $F$ is a finite extension  of~$\Q_{p}$,
if $\cO$ is the ring of integers in a finite extension of~$\Q_{\ell}$
{\em (}$p$ and $\ell$ each  denoting some fixed prime{\em )},
and if $d$ is a  positive  integer, then the category of smooth
$\GL_d(F)$-representations on torsion $\cO$-modules
admits a fully faithful embedding into the category of 
quasicoherent sheaves on~$\cX$, an appropriate moduli stack parameterizing
$d$-dimensional $\ell$-adic representations of the absolute Galois group~$\Gal_F$.
\end{conj}

We have described this  conjecture as ``rough'' for several reasons:
one should be precise about the stack $\cX$ appearing in the statement of
the conjecture --- when  $\ell \neq p$, it should be understood to be
the stack of Weil--Deligne representations introduced in
\cite{zhu2020coherent} (see also~\cite{dat2020moduli}
and~\cite{fargues--scholze} for alternative constructions of the
underlying classical stacks),
while in the case $\ell = p$ it should be one of the moduli stacks
of \'etale $(\varphi,\Gamma)$-modules
constructed in~\cite{emertongeepicture} (and one should {\em a priori}
consider these stacks as derived stacks --- although in the $\ell \neq p$
context, it is shown in
\cite{zhu2020coherent} that the resulting stacks are
in fact classical, and the same has been shown by Min~\cite{min2024classicalityderivedemertongeestack} in the
$\ell = p$ context, using the results of B\"ockle--Iyengar--Pa\v{s}k\={u}nas~\cite{BIP});
one should also be
precise about what sort of categories are being considered --- the envisaged
embedding will not be compatible  with  the natural abelian category
structure on source and target, so one should work with appropriate derived
(or, probably better, stable infinity) categories (where ``appropriate''
is also doing some work --- e.g.\ one should work with 
the triangulated/stable $\infty$-category
of Ind-coherent complexes on $\cX$, and then
make an analogous modification to
the triangulated category of representations as well,
as in \cite[\S 4.1]{zhu2020coherent}).

We also aim for a version of the conjecture where smooth $p$-power torsion representations of ${\rm GL}_d(F)$ are replaced by locally analytic representations (i.e.~the kind of $p$-adic representations
relevant to %
the theory of overconvergent $p$-adic automorphic forms, Coleman families, \dots), and the stack $\cX$ is replaced by a stack of equivariant vector bundles on the Fargues--Fontaine curve (which can also be interpreted as a stack of $(\varphi,\Gamma)$-modules, now with coefficients in 
a Robba ring; see Theorem~\ref{EHequivVBphiGamma} 
 below). As there is not even an informal way to incorporate this case into the above rough conjecture (that simultaneously treats the cases $\ell\neq p$ and $\ell=p$) we leave the formulation of the conjecture in the $p$-adic locally analytic
case to the body of the notes. We only remark that the conjecture in this case can be seen as an overconvergent version of the $p$-adic limit
of Conjecture \ref{conj:main intro}
(given by passing from torsion $\cO$-modules to complete torsion-free $\cO$-modules)
in the $\ell=p$ case.
  
\begin{rem}
  \label{rem: who first thought of this}As far as we know, the possibility of a conjecture along the lines of
Conjecture~\ref{conj:main intro} with $\ell=p$ was first raised by
Michael Harris  (see~\cite[Question~4.7]{MR3530163}). We first learned that a $p$-adic
version of the conjecture of Fargues--Scholze in the form~\eqref{eqn:
  rough form of geometric Langlands} should hold  from
Laurent Fargues in the summer of 2016.  Peter Scholze explained some
related conjectures to T.G.\ in 2018, and a rough formulation of a 
version of Conjecture~\ref{conj: Banach functor} was discussed by M.E.,
T.G., E.H.\ and Scholze at the Hausdorff school in Bonn in 2019. 

\end{rem}

\begin{rem}
  \label{rem: history in the case l=p}We motivated
  Conjecture~\ref{conj:main intro} by the work of Fargues--Scholze. In
  fact, in the case $\ell\ne p$ various
   conjectures along
  the lines of Conjecture~\ref{conj:main intro} were independently
  proposed by Ben-Zvi--Chen--Helm--Nadler, Fargues--Scholze, E.H., and Zhu.
  The history of such
  conjectures is discussed in
  the introductions to the papers \cite{fargues--scholze,
    zhu2020coherent, hellmann2020derived, benzvi2020coherent}. It is
  natural to guess that a similar conjecture should hold in the case
  $\ell=p$, but it is less clear what the precise formulation should
  be. We propose two precise conjectures in Section~\ref{sec:
  categorical conjectures}, one in the ``Banach'' case (focusing on representations on torsion $\cO$-modules and in the limit on lattices in Banach space representations), and one in the
``analytic'' setting (focusing on locally analytic representations), and explain some relationships between 
them. However, it does not seem that either conjecture implies the
other, and we do not know whether they admit a common refinement.

The two cases correspond to two classical notions of $p$-adic
automorphic forms, or $p$-adic modular forms. Serre defined $p$-adic modular forms as limits of $p$-adic modular forms with torsion coefficients (or more precisely: as the $p$-adic limits of the reductions modulo powers of $p$ of the $q$-expansions of modular forms). These spaces of $p$-adic modular forms certainly belong to the ``Banach case", whereas the subspace of overconvergent modular forms, as introduced by Katz, belong to the ``analytic case".
In the discussion in these notes the space of $p$-adic modular forms is replaced by the completed cohomology of modular curves (or more general towers of Shimura varieties), as in Sections~\ref{sec:patching local-global-compatibility} and \ref{subsubsec: cohomology of Shimura}, and the space of overconvergent modular forms is then replaced by the space of locally analytic vectors underlying the completed cohomology, which is used to construct eigenvarieties, see Section~\ref{sec:eigenvar}.

\end{rem}%

\begin{rem}
  \label{rem:why only fully faithful functors}The reader might
  ask why we are only looking for fully faithful functors,
  rather than equivalences of categories as in the work of
  Fargues--Scholze described above. It does seem reasonable to hope
  that there are such equivalences of categories (even in the global
  setting, as speculated in \cite[\S 7]{fargues--scholze-IHES}).
  One difficulty is to even define  categories of  $p$-adic
  sheaves on $\Bun_G$. These difficulties are addressed in
  the recent work of Lucas
  Mann~\cite{https://doi.org/10.48550/arxiv.2206.02022} and in joint work of Lucas Mann with Johannes Anschütz and Arthur-César le Bras \cite{anschütz20246functorformalismsolidquasicoherent} that will be extended in forthcoming work of Johannes Anschütz, Arthur-César le Bras, Juan Esteban Rodriguez Camargo and Peter Scholze on analytic syntomification.  This work
  should then play an important role 
  in formulating a precise conjecture;
 but 
  nevertheless, at the time of writing there is no such precise
  conjecture in the literature (let alone a theorem --- and we note
  that as far as we are aware, it is not expected that the
  Fargues--Scholze construction of the spectral action will go over to
  the $p$-adic setting in any simple manner). In contrast, we can make
  precise conjectures about the fully faithful embeddings of
  categories of representations of $\GL_n(F)$, and even prove theorems
  in the case of~$\GL_2(\Qp)$. Furthermore, at least for~$\GL_2(\Qp)$
  and~$\GL_2(\Q_{p^2})$ (and their inner forms), we produce
  semiorthogonal decompositions of the categories of (Ind-)coherent
  sheaves, which we expect to correspond to a semiorthogonal
  decomposition of the appropriate category of sheaves on~$\Bun_G$
  induced by the closure relations on~$\Bun_G$ itself. (See
  Section~\ref{subsec:semi-decomp} for some results in this
  direction.)
\end{rem}%

\begin{remark}
\label{rem:why only GL}
There is no particular reason to restrict to the group
$\GL_d$ in the statement  of Conjecture~\ref{conj:main intro},
rather than considering dual groups or $L$-groups  or $C$-groups
of more general reductive groups over~$F$ (as in
\cite{benzvi2020coherent, hellmann2020derived, zhu2020coherent}). 
However, since we are interested in the $\ell = p$ case,
we will for the most part focus on the case of~$\GL_d$, since the relevant
``stacks of parameters'' have been constructed in this case.
(Since the literature on Taylor--Wiles patching, which sits in
the background as an important source of intuition 
and motivation, also focuses on this case,
this provides another reason for us to focus on it  in this paper as
well.)%
\end{remark}

\subsection{Some differences between the \texorpdfstring{$\ell$}{l}-adic and \texorpdfstring{$p$}{p}-adic
  cases}%
\label{sec:some-diff-betw}While it is possible to formulate
$\ell$-adic and $p$-adic versions of Conjecture~\ref{conj:main intro}
in a somewhat uniform fashion, there are significant differences
between the two cases. %

That there should be differences is not surprising from the point of view of
representation theory. Indeed, already for $\GL_2(F)$, the
classification of irreducible smooth mod~$\ell$ representations is
well-understood for $\ell\ne p$, and can be formulated uniformly for
any~$F$ (these results are due to Vign\'{e}ras,
see~\cite{MR1026328}). In contrast, for $\ell=p$, there is only a
classification when $F=\Qp$ (due to Barthel--Livn\'e and Breuil
\cite{BarthelLivneDuke, BreuilGL2I}), while %
for $F$ a non-trivial unramified extension of~$\Qp$ %
many complications arise; for example, there are infinite families of admissible irreducible
representations not lying in the principal series (constructed by Breuil--Pa\v{s}k\={u}nas
\cite{BreuilPaskunas}), and (absolutely) irreducible representations
which are not admissible (constructed by Le~\cite{MR4078693}).
Furthermore, supersingular irreducible
admissible representations are never finitely presented (see Schraen's paper
\cite{MR3365778} and Wu's~\cite{MR4280498}; here we use
``finitely presented'' in the sense of Definition~\ref{def:finitely presented}).

These differences and difficulties are mirrored on the Galois side of
the correspondence. For $\ell\ne p$ the moduli stacks of Galois
representations are $0$-dimensional ind-algebraic stacks of a rather
mild type: the ind-structure arises from there being infinitely many
connected components (controlled by the conductors of the Galois
representations), and each connected component is a quotient of an
affine scheme by a reductive group. Accordingly these stacks admit
many global functions, and have coarse moduli spaces which can be
described in terms of pseudorepresentations (at least away from a
small number of bad primes; this has been proved by Fargues--Scholze~\cite[Thm.\ VIII.3.6]{fargues--scholze}).

In the Banach
case for $\ell=p$, the moduli stacks we consider are those defined
in~\cite{emertongeepicture}. These are also ind-algebraic stacks, but
their dimension grows with~$[F:\Qp]$, and the ind-structure is much
more complicated: indeed, they are naturally formal algebraic stacks
with only finitely many connected components (arising just %
from a condition on central characters). Furthermore, we do not expect them
to admit any interesting global functions, and beyond the case
of~$\GL_2(\Qp)$ (and the case of tori) we do not expect them to have
non-trivial associated moduli spaces.

In the analytic case the moduli stacks are defined in these notes (a more detailed discussion is going to appear in the forthcoming paper \cite{HellmannHernandezSchraen}) as the moduli stacks of equivariant vector bundles on the Fargues--Fontaine curve (or equivalently as moduli stacks of $(\varphi,\Gamma)$-modules over the Robba ring). 
As far as the dimension is concerned the analytic case behaves similar to the Banach case: the (expected) dimension grows with $[F:\Q_p]$ and should coincide with the dimension of the generic fiber of the formal stacks in the Banach case (and in fact conjecturally this generic fiber should be closely related to an open substack of the analytic moduli stack). However, these stacks now should not have an ind-structure but rather should be rigid analytic Artin stacks. Though there are a few more interesting global functions on these stacks we expect that they all arise from the theory of Hodge--Tate weights and hence the algebra of global functions still loses much of the %
information contained in the moduli stack.

\begin{rem}
  \label{rem: already a mess for GL2 beyond Qp}Another difference
  between the cases $\ell=p$ and $\ell\ne p$ is the dependence on the
  field~$F$. In the case $\ell\ne p$, the (in general conjectural)
  functors admit a qualitatively uniform description over all~$F$, and
  it does not seem to be any easier to establish instances of the
  conjecture for $F=\Qp$ than for a general~$F$. In contrast, if
  $\ell=p$ then we only know how to construct a functor for
  $\GL_2(\Qp)$ in the Banach case, in which case it shares some qualitative similarities
  with the case $\ell\ne p$ (for example, it is ``not very derived'',
  with the only derived phenomena relating to finite-dimensional
  representations; and supersingular irreducible admissible 
  representations roughly correspond to skyscraper sheaves supported
  on irreducible Galois representations).

  However already for~$\GL_2(\Q_{p^2})$ we expect it is essential to
  consider the functor on the derived level, and supersingular irreducible
  admissible representations no longer correspond to
  skyscraper sheaves, but instead should correspond to complexes with
  non-coherent~$H^{-1}$. (See Section~\ref{subsec: Banach Qp2 stuff}.)
\end{rem}

\begin{rem}
  \label{rem: analogy between Betti de Rham and ell p Langlands}The
  differences between the $\ell$-adic and $p$-adic settings are
  somewhat reminiscent of the differences between the Betti and de Rham
  versions of geometric Langlands (see for example~\cite{MR3821166}
  for the former and~\cite{MR3364744} for the latter), although for
  reasons of ignorance we are
  unsure of how seriously one should take this analogy.
\end{rem}

\subsection{Motivation from Taylor--Wiles patching}
\label{subsec:patched-modules}%
Fargues has emphasised that the conjectures described in
Section~\ref{subsec: FS} did not arise by attempting to ``copy and
paste'' from the geometric Langlands program to the setting of the
Fargues--Fontaine curve, but rather from a study of $p$-adic period
morphisms, the $p$-adic geometry of Shimura varieties and
Rapoport--Zink spaces, and the classical local Langlands
correspondence (in particular the phenomena of $L$-packets and
endoscopy). Similarly, while it may be possible to arrive at the
conjectures presented here by an appropriate procedure of ``setting
$\ell=p$ in the Fargues--Scholze conjecture'', we instead came to them via %
a roundabout process over
many years, largely motivated by considerations coming from the theory of $p$-adic
automorphic forms and the Taylor--Wiles patching method. %

As well as its importance in shaping our ideas, the Taylor--Wiles
patching method has suggested several important properties of the
conjectural functors as in Conjecture~\ref{conj:main intro}, even in
the case $\ell\ne p$. For example, the patching method suggests that the
``integral kernel'' of such a functor should be a genuine sheaf (rather than a
complex thereof), a property that was not at all apparent in the first
versions of the constructions of Ben-Zvi--Chen--Helm--Nadler and
Zhu. (There is also a strong connection between the Taylor--Wiles
method and the classical $p$-adic local Langlands correspondence
for~$\GL_2(\Qp)$, for which see~\cite{MR3732208}.) %

In brief, the connection between Conjecture~\ref{conj:main intro} and
the Taylor--Wiles method is as follows. A typical example of a smooth representation of $\GL_d(F)$ is 
a compact induction $\cInd_{\GL_d(\cO_F)}^{\GL_d(F)} V,$ where $V$
is a continuous (i.e.\ smooth) representation of $\GL_d(\cO_F)$ on
a literally finite $\cO$-module~$V$.   One envisages that the embedding
of Conjecture~\ref{conj:main intro}
will map such representations to genuine coherent sheaves (as opposed
to complexes of such), say $\mathfrak A(V)$,
and that these coherent sheaves will be ``globalizations''
over~$\cX$ of the patched modules $M_{\infty}(V)$ that arise in the
Taylor--Wiles--Kisin patching method,
in the sense that if $R_{\infty}$
is one of the rings usually denoted this way in the theory of patching,
and if $\Spf R_{\infty} \to \cX$  is the natural versal morphism,
then the pullback of $\mathfrak A(V)$ along this 
map coincides with~$M_{\infty}(V)$.  %
In particular, the conjecture entails
that patched modules should be ``purely local'' (where now ``local'' is 
understood to pertain to the local Galois-theoretic situation,
in opposition to the global automorphic/Galois-theoretic choice that was made  
in order to perform patching).
Ascertaining the truth or otherwise of this ``pure locality'' 
in general is well-known to be a major
open problem in the arithmetic theory of automorphic forms.

We note that, building upon \cite{Gpatch}, the patching strategy has been extended to the ``analytic case'' in \cite{MR3623233}. Using roughly the same formulation as above one can
associate patched modules %
to locally analytic representations, %
and in particular to locally analytic principal series representations, i.e.\  representations
of the form ${\rm Ind}_B^G(\delta)^{\rm an}$. The corresponding coherent sheaf should then be the pullback of the evaluation of the conjectural ``analytic'' functor on such a representation %
(and is related to the fiber over $\delta$ of the sheaf of overconvergent $p$-adic automorphic forms of finite slope on eigenvarieties in a similar way to
that in which the usual patched modules are related to modules of classical automorphic forms). 
The question of ``pure locality'' of patched modules in this analytic context has again 
arisen as a major open problem in the theory.

We can also
proceed in the opposite direction.  Namely, in light of the fact 
that the patched modules~$M_{\infty}(V)$ are constructed from the
cohomology of Shimura varieties (with coefficients in a local system
corresponding to~$V$), it seems natural to ask whether the cohomology of
Shimura varieties can in some way be computed in terms of the
(conjectural) functors of Conjecture~\ref{conj:main intro}. This
should indeed be the case, as we explain in
Section~\ref{sec:localtoglobal}  (and as will be more thoroughly explained in the paper in
preparation of M.E.\ and Zhu~\cite{emertonzhu}). In parallel to this, Section \ref{sec:eigenvar} gives a precise computation of the sheaves of overconvergent automorphic forms on eigenvarieties (that can be defined in terms of completed cohomology of Shimura varieties) in terms of the conjectural functor in the analytic case.

\subsection{The ``Galois to automorphic'' direction}\label{subsection: connection to
  Breuil program}A variant of the expectation that patched modules are
``purely local'' is the hope that it should be possible to canonically
associate an admissible Banach representation of $\GL_n(F)$ to each
$n$-dimensional $p$-adic representation of $\Gal_F$; and that this
association should be realised in the cohomology of Shimura
varieties. This expectation is borne out by a good deal of evidence in
the case of Shimura curves, the most recent example being the spectacular results on
GK-dimensions due to Breuil, Herzig, Hu, Morra, and Schraen
\cite{https://doi.org/10.48550/arxiv.2009.03127,
  https://doi.org/10.48550/arxiv.2102.06188}.

While our focus is mostly on functors from smooth representations to
sheaves on stacks of Galois representations, our conjectures also entail the
existence of adjoint functors from sheaves to representations of
$p$-adic groups. In
particular, applying these to ``skyscraper'' sheaves (that is, to
sheaves of finite length, supported at particular Galois
representations), we are able to show that our conjectures are
consistent with the above-mentioned
hope. (See e.g.\ Remark~\ref{rem: we expect a unique supersingular pi with skyscraper H0}.)

\subsection{Topics we omit}For reasons of time and ignorance, we
do not discuss many recent (and even not so recent) developments in the $p$-adic
Langlands program, which we nonetheless expect have important
connections to the ideas discussed here. (In revising these notes in early 2025, we have resisted the temptation to try to significantly update their contents, and the following list reflects a list
of topics omitted at the time of the original lectures; it could be significantly lengthened now!) In particular we say
nothing about the cohomology of the Drinfeld upper half plane and
related $p$-adic locally symmetric spaces, as studied by Colmez,
Dospinescu, Le Bras, Niziol and Pan \cite{MR4073863, https://doi.org/10.48550/arxiv.2204.11214, MR4255044, MR3702670,
  MR3641878} (the interested reader could see
\url{https://mathoverflow.net/questions/432932/how-does-the-cohomology-of-the-lubin-tate-drinfeld-tower-fit-into-categorical-p}
for an explanation of the expected connections); about Scholze's $p$-adic
Jacquet--Langlands functor \cite{scholze};  or about the coherent cohomology of Shimura
varieties, and in particular about any relationship to Boxer and
Pilloni's higher Hida and Coleman theories~\cite{pilloniHidacomplexes,
  https://doi.org/10.48550/arxiv.2002.06845,
  https://doi.org/10.48550/arxiv.2110.10251} or to Diamond and
Sasaki's geometric Serre weight conjectures~\cite{Diamond_2022}.

  \subsection{A brief guide to the notes}
\label{sec:guide-notes}
These notes cover many topics, and we hope that it is possible to read
many of the sections more or less independently from each other,
referring back for definitions where necessary. After setting up some notation in Section~\ref{sec:notation}, in Section~\ref{sec: TW
  patching} we give an idiosyncratic introduction to the Taylor--Wiles
method, guided by its connections to the classical $p$-adic local
Langlands correspondence for $\GL_2(\Qp)$, and its role in motivating
our conjectures. In
Section~\ref{sec:moduli-stacks-phi-gamma} we very briefly recall the
definitions and basic properties of the stacks of
$(\varphi,\Gamma)$-modules from~\cite{emertongeepicture}. The rather
longer Section~\ref{sec:analytic-case} explores analogues of these
constructions in the ``analytic'' context of
$(\varphi,\Gamma)$-modules over Robba rings. 

Section~\ref{sec:
  categorical conjectures} states our main (local) $p$-adic Langlands
conjectures, in both the Banach and analytic settings, and explains
their relationships to other topics in the literature, e.g.\ the
Breuil--M\'ezard conjecture. We then explain what is known about these
conjectures in Section~\ref{sec: known cases}; as well as the
relatively straightforward 
case of~$\GL_1$, we explain forthcoming work of Andrea Dotto with
M.E.\ and T.G.\ which proves the Banach conjecture for $\GL_2(\Qp)$,
and we explain some evidence for the Banach conjecture for $\GL_2(F)$
for $F/\Qp$ unramified. %
Section~\ref{sec: categoric LL for l not p} briefly recalls some of
the results and expectations in the case $\ell\ne p$. This material is
used in Section~\ref{sec: global stacks and cohomology of Shimura
varieties}, which explains conjectures on the (completed) cohomology
of Shimura varieties and eigenvarieties.

The  appendices recall various results on (infinity-) category theory
and the representation theory of $p$-adic analytic groups, and establish
others for which we do not know of a reference. %

\begin{rem}
  \label{rem: technicalities and condensed mathematics}The categorical
  $p$-adic Langlands correspondence is still in its infancy, and (as
  these notes make plain) it is unclear to us what the ultimate form
  of many of the conjectures we make should be. There is also a
  considerable amount of technical machinery needed to get off the
  ground. While we have endeavoured to make precise and correct
  statements, our aim has been where possible to explain and
  illustrate phenomena rather than give a detailed development of all
  of this machinery. In particular, we for the most part treat
  infinity categories and derived algebraic geometry as a black box.

  It seems extremely likely that the future development of the theory
  will rely on condensed mathematics, and in particular the theory of
  solid rings (and modules, sheaves\dots). It is also likely that some
  of the definitions and constructions would be streamlined by
  phrasing them in these terms. However, due to our own limitations, %
  we have for the most part not attempted to do this.
\end{rem}

\subsection{Acknowledgements}We would like to thank the organisers
(Pierre-Henri Chaudouard,  Wee Teck Gan, Tasho Kaletha, and Yiannis
Sakellaridis) and participants of the IHES Summer School on the
Langlands Program for the opportunity and encouragement to write these notes, and their patience with our revisions.

We would like to thank the many
colleagues, collaborators and friends with whom we have discussed the
ideas presented in these notes. In particular we would like to thank
Johannes Ansch\"utz, David Ben-Zvi, Roman Bezrukavnikov, George Boxer, Christophe
Breuil, Frank Calegari, Ana Caraiani, Pierre Colmez, Gabriel
Dospinescu, Andrea Dotto, Laurent Fargues, Tony Feng, Dennis
Gaitsgory, David Geraghty, Ian Grojnowski, David Hansen, Michael
Harris, David Helm, Christian Johansson, Arthur-C\'esar Le Bras, Bao Le Hung, Brandon
Levin, Jacob Lurie, Lucas Mann, David Nadler, James Newton, Lue Pan, Vytautas Pa\v{s}k\={u}nas, Vincent Pilloni, Alice
Pozzi, Juan Esteban Rodriguez Camargo, Joaquin Rodrigues Jacinto, David
Savitt, Peter Scholze, Benjamin Schraen, Jack Sempliner, Sug Woo Shin,
James Timmins, Pol van Hoften, Akshay
Venkatesh, Carl Wang-Erickson, and Xinwen Zhu.

We are very grateful to the referees for their many helpful comments and corrections. We would also like to thank John Bergdall, Pierre Colmez, Fred Diamond, Andrea Dotto, Elmar Große-Klönne,
David Hansen, Claudius Heyer, Heejong Lee, Vytautas Pa\v{s}k\={u}nas, and David Savitt for their
comments on earlier versions of these notes. %

\section{Notation}
\label{sec:notation}
\subsection{Fields}
We fix an algebraic closure~$\Qbar$ of~$\Q$, and an algebraic
closure~$\Qpbar$ of~$\Qp$ for each prime~$p$. If~$F/\Qp$ is a finite
extension, write~$\Gal_F$ for the absolute Galois
group~$\Gal(\Qpbar/F)$. Write~$I_F$ for the inertia subgroup
of~$\Gal_F$, $W_F$ for the Weil group, and $k_F$ for the residue field of the ring of
integers~$\cO_F$ of~$F$. We normalise local class field theory so that a uniformizer
corresponds to a geometric Frobenius element.

If $M$ is a number
field, we write $\Gal_M:=\Gal(\Qbar/M)$ for the absolute Galois group
of~$M$. If~$S$ is a set of places of~$M$ then we write  $M(S)/M$ for
the maximal extension unramified outside of~$S$, and write
$\Gal_{M,S}:=\Gal(M(S)/M)$. We fix embeddings $\Qbar\into\Qpbar$ for
each~$p$ and $\Qbar\into\C$, so that if~$v$ is a place of~$M$, there
is a homomorphism $\Gal_{M_v}\to\Gal_M$.

Let $\cO$ denote the ring of integers in a fixed finite extension
$L$ of $\Q_p$, let~$k$ be the residue field of~$\cO$,
and let $\varpi$ denote a uniformizer of~$\cO$. (These will be the coefficients of our various representations.)

\subsection{\texorpdfstring{$p$}{p}-adic Hodge theory}
If $\rho$ is a de Rham representation of $\Gal_F$ on an
$L$-vector space~$W$, then we will write $\WD(\rho)$ for the
corresponding Weil--Deligne representation of $W_F$ (see Section~\ref{EHsub:sub:deRham} for more details of this construction), and if $\sigma:F \into L$ is a continuous
embedding of fields then we will write $\HT_\sigma(\rho)$ for the
multiset of Hodge--Tate weights of $\rho$ with respect to $\sigma$,
which by definition contains $i$ with multiplicity
$\dim_{L} (W \otimes_{\sigma,F} \widehat{\overline{F}}(i))^{\Gal_F} $; %
for example, if~$\varepsilon$ denotes the $p$-adic cyclotomic character,
then $\HT_\sigma(\varepsilon)=\{ -1\}$.

Suppose that $L$ contains the images of all continuous embeddings
$F\into\overline{L}$. By a \emph{$d$-tuple of labeled Hodge--Tate
  weights~$\underline{\lambda}$}, \index{labeled Hodge--Tate weights}
\index{$\underline{\lambda}$} we mean a tuple of integers
$\{\lambda_{\sigma,i}\}_{\sigma:F\into L,1\le i\le d}$ with
$\lambda_{\sigma,i}\ge \lambda_{\sigma,i+1}$ for all~$\sigma$ and all
$1\le i\le d-1$. We will also refer to~$\underline{\lambda}$ as a
\index{Hodge type} \index{inertial type} \emph{Hodge type}. We say
that~$\lambdau$ is \emph{regular} if
$\lambda_{\sigma,i}> \lambda_{\sigma,i+1}$ for all~$\sigma$ and all
$1\le i\le d-1$.
We say that a de Rham representation~$\rho$ has Hodge type~$\underline{\lambda}$ (or
labeled Hodge--Tate weights~$\underline{\lambda}$) if for
each~$\sigma:F\into L$ we
have~$\HT_\sigma(\rho)=\{\lambda_{\sigma,i}\}_{1\le i\le d}$.

By an \emph{inertial type~$\tau$} we mean a
representation $\tau:I_F\to\GL_d(L)$ which extends to a representation
of~$W_F$ with open kernel (so in particular, $\tau$ has finite image).
We say that a de Rham representation~$\rho$ has inertial type~$\tau$ if~$\WD(\rho)|_{I_F}\cong\tau$.

\subsection{Hodge--Tate weights vs. highest weights}
Under the $p$-adic local Langlands correspondence,
Hodge--Tate weights on the Galois-theoretic side of the correspondence
will be related to highest weights on the representation theoretic side,
and it is useful to have notation adapted to either side of the correspondence. As usual, the Hodge--Tate weights differ from the highest weights by a ``$\rho$-shift''.
 
\begin{defn}
  \label{defn: lambda versus xi}If~$\lambdau$ is a regular Hodge type,
  then we write~$\xi_{\sigma,i}=i-1-\lambda_{\sigma,d+1-i}$, so that
$\xi_{\sigma,1}\ge\dots\ge\xi_{\sigma,d}$.  %
 \end{defn}

\begin{defn}
  \label{defn: W lambda}If~$\lambdau$ is a regular Hodge type,
  then we view each
$\xi_\sigma:=(\xi_{\sigma,1},\dots,\xi_{\sigma,d})$ as a dominant
weight of the algebraic group~$\GL_d$ (with respect to the upper
triangular Borel subgroup). We let~$W_{\lambdau}$ be
the corresponding $\cO$-representation of~$\GL_d(\cO_F)$, defined as
follows: for each~$\sigma:F\into L$, we write~$W_{\xi_\sigma}$
for the
algebraic $\cO_F$-representation of~$K$ (more precisely, the dual Weyl
module) of highest
weight~$\xi_\sigma$. Then we define
\[W_{\lambdau}:=\otimes_{\sigma}W_{\xi_\sigma}\otimes_{\cO_F,\sigma}\cO,\]and
write $L(\xi)=W_{\lambdau}[1/p]$, a representation of~$\GL_d(F)$
(or of the algebraic group $\Res_{F/\Q_p} \GL_d$).
Note then that, despite the notation, the highest weight of ~$W_{\lambdau}$ corresponds to the~$\xi_{\sigma}$, rather than the~$\lambda_{\sigma}$.
\end{defn}

\begin{remark}
The preceding notation may not be optimal, but serves our purposes.
The ``$W$'' stands for (dual) {\em Weyl}, and emphasizing the dependence 
on  $\lambdau$ (rather than on~$\xi$) facilitates the comparison with
other considerations of $p$-adic Hodge theory.  
The notation $L(\xi)$ is chosen because it is traditional 
in the theory of category~$\cO$.  Indeed, for any (not necessarily dominant, not necessarily
integral) weight~$\xi$ (of the Lie algebra of $\Res_{F/\Q_p} \GL_d$),
we use $M(\xi)$ to denote the Verma module of highest weight~$\xi$,
and $L(\xi)$ to denote the unique simple quotient of~$M(\xi)$. 
When $\xi$ is dominant integral, one finds that $L(\xi)$ is the usual
highest weight representation~$W_{\lambdau}[1/p]$,  explaining our 
choice of notation in this case.
\end{remark}

\section{Taylor--Wiles patching as a motivation for categorical
  \texorpdfstring{$p$}{p}-adic local Langlands}\label{sec: TW patching}
In this section we review some of the history of the Taylor--Wiles
patching method and of the $p$-adic local
Langlands correspondence for~$\GL_2(\Qp)$, emphasising where possible
the connections between them. We do not pretend to give a thorough
overview, and we view everything through the lens of the
categorical Langlands program.  The reader may wish
to consult Calegari's survey~\cite{calegari2021reciprocity} for a more
mainstream account of developments in the Taylor--Wiles method since
its inception. %
\subsection{A very brief introduction to patching}\label{sec: very
  brief introduction to patching}%
Patching was first introduced by Taylor and Wiles as a technique for proving
modularity lifting theorems~\cite{MR1333036}.  It was then further developed by
Kisin, who showed that it provides a mechanism for relating local
and global aspects of the theory of $p$-adic Galois representations~\cite{MR2600871,MR2392362,KisinFM,MR2827797,geekisin}.

We assume throughout this section that ~$p>2$. We begin by placing ourselves
in the following simple but illustrative context.
Fix a continuous, absolutely irreducible  representation $\rbar:
\Gal_{\Q,\{p,\infty\}} \to \GL_2(k)$ with determinant $\varepsilonbar^{-1}$, and  suppose that $\rhobar :=
\rbar_{| \Gal_{\Q_p}}$ and $\rbar|_{\Gal_{\Q(\zeta_p)}}$ are also
absolutely irreducible. We assume that~$\rbar$ is modular (equivalently, by Serre's conjecture \cite{MR2551763}, we assume that~$\rbar$ is odd).
For any finite set $S$ of primes containing $p$ and $\infty$, we
may consider the global deformation space $\cX_S(\rbar)$
of $\rbar$ over $\cO$, which parameterizes deformations
of $\rbar$ to representations of $\Gal_{\Q,S}$ over complete
Noetherian local
$\cO$-algebras, with fixed determinant~$\varepsilon^{-1}$. %
We may also consider the local deformation
space $\cX(\rhobar)$, which parameterizes deformations
of $\rhobar$ to representations of $\Gal_{\Q_p}$ over complete Noetherian local
$\cO$-algebras, again with fixed determinant~$\varepsilon^{-1}$.
These are each the $\Spf$ of a %
complete Noetherian local $\cO$-algebra,
denoted $R_S(\rbar)$ and $R(\rhobar)$ respectively.
Restriction of Galois representations induces a morphism
\numequation
\label{eqn:deformation restriction}
\cX_S(\rbar) \to \cX(\rhobar),
\end{equation}
which is a finite morphism (see \cite{MR4077579} and~\cite{MR3294389}). %

The reduced tangent space to the closed point of $\cX_S(\rbar)$
(resp.\ $\cX(\rhobar)$)  is equal
to the global Galois cohomology group
$H^1\bigl(G_{\mathbf Q,S}, \Ad^0(\rbar)\bigr)$
(resp.\ to the local Galois cohomology group
$H^1\bigl(G_{\mathbf Q_p},\Ad^0(\rhobar)\bigr)$),
and the restriction map between deformation spaces induces the restriction map
\numequation
\label{eqn:tangent restriction}
H^1\bigl(G_{\mathbf Q,S}, \Ad^0(\rbar)\bigr) \to
H^1\bigl(G_{\mathbf Q_p},\Ad^0(\rhobar)\bigr)
\end{equation}
on reduced tangent spaces, which is given by the usual pullback of
cohomology classes.

In order to simplify the situation,
we will assume that~(\ref{eqn:deformation restriction})
is in fact a closed immersion when $S = \{p,\infty\};$
equivalently, we assume that the   Selmer group %
$\Sel\bigl(G_{\mathbf Q,\{p,\infty\}}, \Ad^0(\rbar)\bigr)$ vanishes (where, by definition,
this Selmer group is the kernel of~\eqref{eqn:tangent restriction}; i.e.\ it is defined
by requiring that the cohomology classes vanish locally at~$p$).
Under this assumption (and the running assumption that $\rbar|_{\Gal_{\Q(\zeta_p)}}$ is absolutely irreducible),
the formula of Greenberg--Wiles (\cite[Thm.\ 2.18]{MR1605752}) shows that %
$$
\dim_k %
H^1\bigl(G_{\mathbf Q, \{p,\infty\}}, \Ad^0(\rbar)(1) \bigr)
=
\dim_k
\Sel\bigl(G_{\mathbf Q,\{p,\infty\}}, \Ad^0(\rbar)\bigr)
+1=1.$$
(The $H^1$ appearing 
on the left-hand side of this equation is, in our context, the ``dual Selmer group''
that appears in the Greenberg--Wiles formula.)
The Taylor--Wiles method then requires us to consider primes $q \equiv 1 \bmod p$
such that the restriction morphism
$
 H^1\bigl(G_{\mathbf Q, \{p,\infty\}}, \Ad^0(\rbar)(1) \bigr)
\to H^1\bigl(\Gal_{\Q_q},\Ad^0(\rbar)(1)\bigr)$
is non-zero (equivalently, injective),
and such that $\rbar(\Frob_q)$ has distinct eigenvalues;
these are the so-called {\em Taylor--Wiles primes}.
Another application of the Greenberg--Wiles formula then shows
that
$$\Sel\bigl(\Gal_{\Q, \{p,q,\infty\}}, \Ad^0(\rbar)\bigr) = 0$$
for such  a Taylor--Wiles prime~$q$ (where again the Selmer condition is defined
by requiring that the cohomology classes vanish locally at~$p$),
and thus that
the restriction  morphism of deformation spaces
\numequation\label{eqn: restriction morphism with q added}\cX_{\{p,q,\infty\}}(\rbar) \to \cX(\rhobar)\end{equation}
is again a closed immersion. 

Write $\rbar_q:=\rbar|_{\Gal_{\Q_q}}$.
Under our
hypotheses that $q\equiv 1\bmod p$ and that $\rbar(\Frob_q)$ has
distinct eigenvalues,
the deformation problem for $\rbar_q$ is representable (even though $\rbar_q$ has 
larger-than-scalar endomorphisms).  Furthermore, the deformation ring $R(\rbar_q)$
takes a
very simple form: namely, $R(\rbar_q)$ is  of the form
 $\cO[[x,y]]/((1+x)^{p^n}-1)$, where $p^n$ is the largest power of~$p$
 dividing $(q-1)$. %
(More precisely, one
checks that all deformations of~$\rbar_q$ are %
a direct sum of characters, each deforming one of unramified characters of which
$\rbar_q$ is the direct sum. Since the determinant of our deformation is fixed,
either one of these characters determines the other, and so, fixing
one of the two summands of $\rbar_q$ once and for all, we see that $R(\rbar_q)$
coincides with the deformation ring of this fixed unramified character. 
The variables $x,y$ then correspond respectively to a generator of the
tame inertia group and to a Frobenius element respectively.)

We now consider the restriction morphism of deformation spaces 
$$\cX_{\{p,q,\infty\}}(\rbar) \to \cX(\rbar_q)$$
(where of course $\cX(\rbar_q) :=\Spf R(\rbar_q)$).
By the discussion of the preceding paragraph, we can rewrite this as a morphism
\numequation
\label{eqn:restriction to q}
\cX_{\{p,q,\infty\}}(\rbar)  \to \Spf \cO[[x]]/((1+x)^{p^n}-1),
\end{equation}
whose fibre over the locus $x =0$ is equal to~$\cX_{\{p,\infty\}}(\rbar)$;
in this way we may regard 
$\cX_{\{p,q,\infty\}}(\rbar)$
as a kind of thickening\footnote{Since if we look modulo the uniformizer of~$\cO$,
then we find that $\Spec k[[x]]/(1+x)^{p^n} -  1)  = \Spec k[x]/x^{p^n}$ is a nilpotent thickening
of~$\Spec k$.}
of~$\cX_{\{p,\infty\}}(\rbar)$.
(This interpretation of the situation would be optimal
if the morphism~\eqref{eqn:restriction to q} were flat.  This is not always the case, but what
comes out of the Taylor--Wiles method is that $\cX_{\{p,q,\infty\}}(\rbar)$ 
supports a faithful module which is flat over~$\cO[[x]]/((1+x)^{p^n} -1),$
namely the completed homology group $\tH_1(q)$ introduced below,
so the situation is always fairly close to the optimal one.)

Now it turns out that $\cX_{\{p,\infty\}}(\rbar)$ has relative dimension $2$
over~$\cO$, while $\cX(\rhobar)$ has relative dimension~$3$ over~$\cO$.
Thus the closed immersion~(\ref{eqn:deformation restriction})
has codimension~$1$, and the closed immersion~\eqref{eqn: restriction morphism with q added}
then realizes $\cX_{\{p,q,\infty\}}(\rbar)$ as a thickening (in the 
``$x$ direction'', if  we refer to~\eqref{eqn:restriction to q})  
of $\cX_{\{p,\infty\}}(\rbar)$
along the transverse direction to $\cX_{\{p,\infty\}}(\rbar)$
in~$\cX(\rhobar)$.

A key point is that for any given power
$p^n$ of $p$, we can find plenty of Taylor--Wiles primes $q$ (ultimately via
a \v{C}ebotarev argument), such that $q \equiv 1 \bmod p^n$.
Thus we can thicken up~$\cX_{\{p,q,\infty\}}(\rbar)$ as  much
as we like inside~$\cX(\rhobar)$, and
the rather surprising idea at the heart of the patching method is 
to try to realize 
$\cX(\rhobar)$ as some kind of limit of all these thickenings~$\cX_{\{p,q,\infty\}}(\rbar)$.

To make this precise, we first 
regard each deformation ring $R(\rbar_q)$ as an
$S_\infty:=\cO[[x]]$-algebra (this gives a uniform meaning to the ``$x$  direction''
referred to above  as $q$ varies), 
and then try to glue together the various formal schemes %
$\cX_{\{p,q,\infty\}}(\rbar)$ 
over $S_{\infty}$ %
as $q$ varies in a manner compatible with their embeddings into~$\cX(\rhobar)$;
equivalently, we try to glue together the various
$\bigl(S_{\infty},R(\rhobar)\bigr)$-bialgebras~$R_{p,q,\infty}(\rbar)$.
As part of the patching process,
we simultaneously glue together 
various (completed) homology
groups.

Recall that we define the completed (first) homology of the modular curve as  %
\numequation\label{eqn:completed homology modular curves}\tH_1 :=  \varprojlim_{r} H_{1}\bigl(Y(p^r), \cO) \cong
\varprojlim_r H^1_{c}\bigl(Y(p^r),\cO\bigr).
\end{equation}
Here $Y(p^r)$ is the modular curve of level $p^r$ (or more precisely a congruence quotient of ${\rm PGL}_2$-symmetric spaces, which is a union of connected components of the usual modular curves for ${\rm GL}_2$). 

This is acted on by a Hecke algebra $\T$ generated by the prime-to-$p$
Hecke operators, as well as by $\PGL_2(\Q_p)$ and $\Gal_{\Q,\{p,\infty\}}$.
All of these actions commute with one another. There is a maximal ideal $\mathfrak m$ of $\T$ corresponding
to the modular Galois representation $\rbar$.   We can localize
 $\tH_1$ %
at $\mathfrak m$ to get a $\T_{\mathfrak m}$-module
$\tH_{1,\mathfrak m}$.
Similarly, for each~$q$ as above we have the completed homology
group
$$\tH_1(q) := \varprojlim_{r} H_{1}\bigl(Y(p^rq), \cO)
\cong \varprojlim_r H^1_c\bigl(Y(p^r q),\cO),$$
which has actions of  $\PGL_2(\Q_p)$, of
$\Gal_{\Q,\{p,q,\infty\}}$, of prime-to-$pq$ Hecke operators, and of
certain Hecke operators at~$q$.  We can and do localize these at~$\m$
(which is now extended to incorporate a Hecke operator at~$q$; the precise 
extension corresponds to the choice that was made above of one of the two summands of $\rbar_q$
when giving the explicit description of $R(\rbar_q)$) 
to get  $S_\infty$- and $R(\rhobar)$-modules $\tH_{1,\mathfrak
  m}(q)$. In fact, since the action of $S_\infty$ on $\tH_{1,\mathfrak
  m}(q)$ is by definition via the composite of our chosen morphism $S_\infty\to
R(\rbar_q)$ and the natural morphism $R(\rbar_q)\to
R_{\{p,q,\infty\}}(\rbar)$, we can use the formal smoothness of
$S_\infty$ and the closed immersion~\eqref{eqn: restriction morphism
  with q added} to produce a morphism $S_\infty\to R(\rhobar)$,
compatible with the actions of both rings on $\tH_{1,\mathfrak
  m}(q)$.

Following \cite[\S9]{scholze}, we make an ultraproduct construction to
produce an $R(\rhobar)$-module $M_\infty$ with an action of
$\PGL_2(\Q_p)$, with the property that %
\[M_\infty\otimes_{R(\rhobar)}R_{\{p,\infty\}}(\rbar)=\tH_{1,\mathfrak
  m}.\]   By
construction (and local-global compatibility at the primes~$q$),
$M_\infty$ is a finite projective $S_\infty[[K]]$-module, where $K=\PGL_2(\Zp)$. %

In general one cannot assume that the Selmer group $\Sel\bigl(G_{\mathbf
  Q,\{p,\infty\}}, \Ad^0(\rbar)\bigr)$ vanishes (equivalently, we
cannot assume that~\eqref{eqn:deformation restriction} is a closed immersion) and the Taylor--Wiles
method goes as follows. Let $g:=\dim_k\Sel\bigl(G_{\mathbf
  Q,\{p,\infty\}}, \Ad^0(\rbar)\bigr)+1$. The 
 \v{C}ebotarev argument alluded to above shows that for each~$n$ there
 is a set $Q_n$ of Taylor--Wiles primes with the properties that
 \numequation\label{eqn: TW primes don't change Selmer group}\dim_k\Sel\bigl(G_{\mathbf
  Q,\{p,\infty\}\cup Q_n},\Ad^0(\rbar)\bigr)=\dim_k\Sel\bigl(G_{\mathbf
  Q,\{p,\infty\}},\Ad^0(\rbar)\bigr),\end{equation} and
 for each $q\in Q_n$, we have $q \equiv 1 \bmod p^n$. We  write
 $S_\infty:=\cO[[x_1,\dots,x_{g}]]$ and
 $R_\infty:=R(\rhobar)[[y_1,\dots,y_{g-1}]]$ for formal variables
 $x_i,y_i$. By~\eqref{eqn: TW primes don't change Selmer group}, we may lift ~\eqref{eqn:deformation
   restriction} and~\eqref{eqn: restriction morphism with q added} to (non-canonical) closed immersions
 $\cX_{\{p,\infty\}}(\rbar)\to\Spf R_\infty$ and
 $\cX_{\{p,\infty\}\cup Q_n}(\rbar)\to\Spf R_\infty$, and straightforward generalizations of the constructions
 explained above in the case $g=1$ yield a finite projective
 $S_\infty[[K]]$-module $M_\infty$ with an action of~$R_\infty$, and
 we can again produce a morphism $S_\infty\to R_\infty$ compatible
 with the actions of each ring on~$M_\infty$. By construction, we
 again have \numequation\label{eqn: going back from Minfty to H}M_\infty\otimes_{R_\infty}R_{\{p,\infty\}}(\rbar)=\tH_{1,\mathfrak
  m}.\end{equation}

The patched module~$M_\infty$ can be thought of as being a way to
extract the local at~$p$ information from completed
(co)homology. Its construction depends on many choices, but the
perspective first taken in~\cite{Gpatch} is that it should be purely
local in an obvious sense, and that the same construction for more
general Shimura varieties should give a candidate for a $p$-adic local
Langlands correspondence beyond the case of~$\GL_2(\Qp)$. This
expectation has yet to be verified, but it is one that has been
important for us in finding the conjectures presented in these notes,
and the modules~$M_\infty$ appear in the statements of our conjectures
(see in particular Remark~\ref{rem: conjectural compatibility with TW
  patching}).

\subsubsection{Patching at finite level}
If $V$ is any finitely generated $\cO$-module with a continuous action
of~$K$, then we set \numequation\label{eqn: formula for
  Minfty V}M_{\infty}(V) :=
M_{\infty}\otimes_{\cO[[K]]}V.\end{equation} Then $V \mapsto M_{\infty}(V)$ is an
exact functor, and~\eqref{eqn: going back from Minfty to H} implies
that \[M_\infty(V)\otimes_{R_\infty}R_{\{p,\infty\}}(\rbar)=H^1_c(V)_{\mathfrak
  m},\]
where $H^1_c(V)_{\mathfrak m}$ denotes the (localized at~$\m$) compactly supported
cohomology group with
coefficients in the local system determined by~$V$.
(Since $\m$ is non-Eisenstein, the natural map  $H^1_c(V)_{\m} \to H^1(V)_\m$ is an isomorphism,
but we have used the $H^1_c$ notation to emphasize that we are working with {\em homology}.)

\begin{rem}
  \label{rem: Kisin and EGS perspective on patching and link to Linfty and A}
  As we indicate below, the definition just given of $M_{\infty}(V)$ in
  terms of~$M_\infty$ is historically backwards. Indeed the original
  Taylor--Wiles method fixed a~$V$ and constructed an associated patched ring --- 
denoted~$R_{\infty}$
  in the original literature, but which we might denote $R_{\infty}(V)$ ---
  which in retrospect can be interpreted as that quotient of the ring $R_{\infty}$ introduced
  above which acts faithfully on~$M_\infty(V)$.  In the contexts originally considered
  by Taylor and Wiles, the patched module $M_{\infty}(V)$ is cyclic over~$R_{\infty}$,
  and hence free over~$R_{\infty}(V)$, and so only appears implicitly (in the use
  of mod $p$ multiplicity one results that show that certain Hecke modules arising 
  from the cohomology of modular curves are free over the Hecke algebras that act on them). 
  The patched modules $M_{\infty}(V)$ themselves were introduced (independently)
  by Diamond~\cite{MR1440309} and Fujiwara~\cite{Fujiwaradefo}, and
  it was the work of Kisin \cite{KisinFM} (as interpreted
  by~\cite{MR3323575}) that then recast patching in terms of the exact
  functor $V \mapsto M_{\infty}(V)$. That this functor is 
  determined by a single ``large'' patched module~$M_{\infty}$ was observed in~\cite{Gpatch}.
  (Note though that the formula~\eqref{eqn: formula for Minfty V}
   does not literally appear in~\cite{Gpatch},
   where instead  one  finds essentially equivalent (but visually more involved)
   formulas  involving $\Hom$s and duals.)

  The reader may wish to compare the formula~\eqref{eqn: formula for Minfty V} 
  for the patching functor to Conjecture~\ref{conj: Banach functor}~\eqref{item: Linfty is a sheaf},
  which suggests the mnemonic
  ``$M_{\infty} = M_{\infty}(\cO[[K]])$''. The jump from considering
  the functor $V \mapsto M_{\infty}(V)$ to the conjectures of these
  notes is roughly to move from considering representations of
  $\PGL_2(\Z_p)$ to representations of~$\PGL_2(\Qp)$, by replacing~$V$
  by the compact induction $\cInd_{\PGL_2(\Zp)}^{\PGL_2(\Qp)}V$, and
  then to move from modules over Galois deformation rings to sheaves
  on stacks of $(\varphi,\Gamma)$-modules.

  Actually, in the context of patching, the passage from $\PGL_2(\Z_p)$ to  $\PGL_2(\Q_p)$
  was one of the main contributions of~\cite{Gpatch}, and so there is a natural progression
  of ideas from~\cite{MR3323575}
  (patching as a functor on $\PGL_2(\Z_p)$-representations)
  to~\cite{Gpatch}
  (patching as a functor on $\PGL_2(\Q_p)$-representations)
  to Expected Theorem~\ref{expectedthm: DEG results} below
  (patching is promoted to a functor from $\PGL_2(\Q_p)$-representations to
  (complexes of) sheaves on a stack of $(\varphi,\Gamma)$-modules).
\end{rem}

As already indicated above, the patching method can be extended from the case of modular
curves to more general Shimura varieties and related congruence quotients.  However,
the case of~$\PGL_2(\Qp)$ considered here is particularly nice, because a $p$-adic local 
Langlands correspondence has been constructed in this case (in fact
for $\GL_2(\Q_p)$).
The idea of a $p$-adic local Langlands correspondence was motivated by
Taylor--Wiles patching, but the correspondence itself was constructed by quite different
techniques. We review some of the history in the rest of this section.

\subsection{\texorpdfstring{$p$}{p}-adic Langlands: an overview and history}
Following preliminary and for the most part unpublished investigations by many researchers (including in
approximate historical order Fontaine--Langlands, Serre, Mazur,
Harris, Vign\'eras, and Schneider--Teitelbaum), the investigation of a
possible $p$-adic Langlands correspondence, relating $p$-adic
representations of the group $\GL_2(\Q_p)$ to $2$-dimensional $p$-adic
representations of the local Galois group $\Gal_{\Q_p}$, was instigated
by Christophe Breuil~\cite{BreuilGL2I, BreuilGL2II}.  He was motivated
to a considerable extent by various considerations that arose out of
the proof of modularity lifting theorems \cite{WilesFLT, MR1333036,
  MR1639612, BCDT}; in particular, he was directly motivated by
considerations arising from the Breuil--M\'ezard conjecture
(see~\cite[\S 1.4]{MR2493214} for an account of these motivations).
We begin our overview by discussing some ideas and methods related
to these theorems, before turning to a discussion of the correspondence
itself, its application to the study of the Breuil--M\'ezard
conjecture, and its relationship to completed (co)homology via
local-global compatibility.

\subsubsection{Taylor--Wiles--Kisin patching}
\label{subsubsec:TWK patching}
Taylor--Wiles patching, as introduced in \cite{MR1333036} and briefly
reviewed above,
is a process that constructs modules over local deformation
rings out of the (co)homology of Shimura varieties (thought of as
Hecke modules). %

The patching method was originally developed to prove modularity
lifting theorems, showing that certain $p$-adic Galois representations
(in the original applications, those associated to elliptic curves)
correspond to modular forms. Originally the patching arguments
required as inputs a number of deep results about congruences between
modular forms, including mod~$p$ multiplicity one theorems, level
raising theorems, cases of the weight part of Serre's conjecture, and
Serre's modularity conjecture itself. Subsequently, the situation has
reversed, and such results can now be deduced from patching
arguments.

While patching arguments are still used to prove that Galois
representations are automorphic, they are now also used to study the
internal properties of Galois representations and of congruences
between automorphic representations, and it is this application which
interests us. We will not attempt to give a full account of these
shifts, and will instead discuss the contemporary
perspective on the patching construction, but in brief we note that
the dependence on  mod~$p$ multiplicity one theorems was removed
independently by Diamond~\cite{MR1440309} and
Fujiwara~\cite{Fujiwaradefo}, level raising was addressed by Taylor's
``Ihara avoidance'' argument~\cite{tay}, the weight part of Serre's
conjecture by T.G.~\cite{geebdj}, and Serre's conjecture
was proved by Khare--Wintenberger~\cite{MR2551763} and
Kisin~\cite{MR2551765}. Many of these developments, and much of what
we explain below, rely crucially on Kisin's improvement of the
Taylor--Wiles method~\cite{MR2600871}.

We maintain the notation introduced above, except
that, to avoid clutter, we write~$R_p$ for~$R(\rhobar)$ from now on; 
and we allow the fixed determinant 
in our global deformation problem
to be an arbitrary odd character, rather than just $\varepsilon^{-1}$,
so that the fixed determinant in the corresponding local
deformation problem is then allowed to be an arbitrary character.
The point of allowing an arbitrary character as the local determinant is that this then
allows us in what follows to consider representations of $\GL_2(\F_p)$, $\GL_2(\Z_p)$, 
or $\GL_2(\Q_p)$ which have a fixed but arbitrary central character --- which is the usual
setting for $p$-adic local Langlands. 
For example,
by carrying out patching on modular curves for $\GL_2$ rather than~$\PGL_2$,
we can obtain a version of the exact functor $V \mapsto M_{\infty}(V)$
described above which maps
the category of continuous $\GL_2(\Z_p)$-representations
on finitely generated $\cO$-modules with appropriately prescribed central character to the category
of finitely generated modules over~$R_\infty$, where~$R_\infty$ is a
power series ring over~$R_p$.
(Here and in the ensuing discussion, we always assume that the fixed determinant
of our deformation problem and the central characters of our $\GL_2$-representations
have been chosen compatibly; since this issue of fixing determinants and central
characters is primarily a technical one, we will try not to belabour it.)

There are certain %
$\GL_2(\Z_p)$-representations on which the evaluation
of $M_{\infty}$ is of particular interest, namely those that
are lattices
in a {\em locally algebraic type}, %
and those that are irreducible representations of $\GL_2(\F_p)$ over $k$.
Evaluating the patching functor $M_{\infty}$ on the
first kind of representation relates
patching to the study of modularity lifting; and
this provided the original motivation for the patching construction.
Evaluating the patching functor on the second kind of representation
establishes a relationship between patching and the weight
	part of Serre's conjecture.

To be more precise, recall that by a locally algebraic type~$\sigma$,
we mean an irreducible $\Qpbar$-representation of~$\GL_2(\Zp)$ obtained by
tensoring an algebraic representation with a finite-dimensional representation corresponding to an
inertial type~$\tau$ via Henniart's inertial local Langlands
correspondence~\cite[App.\ A]{BreuilMezard}; we say that a lift
of~$\rhobar$ is of type~$\sigma$ if it has Hodge--Tate weights
corresponding to the highest weight of the algebraic representation, and inertial
type corresponding to~$\tau$. Suppose that $\sigma$ is a locally algebraic type,
and let $\sigma^{\circ}$ be a $\GL_2(\Z_p)$-invariant lattice in
$\sigma$.
If $R_p(\sigma)$ denotes the quotient of $R_p$
parameterizing lifts of type $\sigma$, then (as a consequence of
local-global compatibility) the $R_p$-action
on $M_{\infty}(\sigma^{\circ})$ factors through $R_p(\sigma)$,
and a modularity lifting theorem for lifts of type $\sigma$ (i.e.\ a theorem
proving that lifts of $\rbar$ of type $\sigma$ are necessarily
modular) can be interpreted in terms of patching as saying
that the action of the corresponding quotient $R_\infty(\sigma)$ of~$R_\infty$ on $M_{\infty}(\sigma)$ should be faithful.

As for the weight part of Serre's conjecture,
this can be rephrased
as the statement that $M_{\infty}(\sigmabar)$ is non-zero
for precisely those Serre weights $\sigmabar$ of $\GL_2(\F_p)$
that lie in a specified set $W(\rhobar)$ (where by definition a Serre
weight is an irreducible $k$-representation of~$\GL_2(\Fp)$). %

\subsubsection{An observation}
\label{subsubsec:observation}
If $V$ is a $\GL_2(\Z_p)$-representation on a finite rank
free $\cO$-module (e.g.\ a lattice in a locally algebraic type),
and $\Vbar := V/\unif V$ is the associated residual representation,
then $M_{\infty}(V)$ is $\cO$-torsion free (by exactness),
and so $M_{\infty}(V) \neq 0$ iff $M_{\infty}(\Vbar) \neq 0$
iff $M_{\infty}(\sigmabar) \neq 0$ for some Jordan--H\"older
	factor $\sigmabar$ of $\Vbar$.
Thus, if the weight part of Serre's conjecture is known,
then the question of whether or not $M_{\infty}(V)$ is non-zero
can be answered via an examination of the set of Jordan--H\"older
factors of $\Vbar$.
This simple observation  motivated much of the
further development of the subject.

\subsubsection{The Breuil--M\'ezard conjecture}
\label{subsubsec:BM}
As already noted, in order to prove a modularity lifting theorem
for Galois representations of locally algebraic type $\sigma$,
one must show that $R_\infty(\sigma)$ 
acts faithfully on $M_{\infty}(\sigma^{\circ})$, for some (equivalently,
any) invariant lattice $\sigma^{\circ}$ contained in $\sigma$.
It follows from the patching construction that $M_{\infty}(\sigma^{\circ})$
is maximal Cohen--Macaulay over $R_\infty(\sigma)$;
so if $R_p(\sigma)$ is a domain, then~$R_\infty(\sigma)$ (which is
then also a domain) necessarily acts faithfully
on $M_{\infty}(\sigma^{\circ})$. 
Perhaps the simplest way to establish that $R_p(\sigma)$
is a domain is to show that it is formally smooth,
and the arguments of \cite{MR1639612} and \cite{BCDT} were devoted to establishing
such a smoothness result (for appropriately chosen~$\sigma)$.
In analyzing the conditions introduced in these papers related to
ascertaining such results, Breuil and M\'ezard were
led to formulate their eponymous conjecture \cite{BreuilMezard}.
This conjecture describes the Hilbert--Samuel multiplicity
of $R_p(\sigma)/\unif$ in terms of the Serre weights of $\rhobar$
that appear as constituents of $\sigma^{\circ}/\unif$.

A key breakthrough in understanding the Breuil--M\'ezard conjecture
was made by Kisin (see e.g.\ \cite{MR2827797}),
who strengthened the observation
of~(\ref{subsubsec:observation}) to note that it follows from
the weight part of Serre's conjecture, the exactness of the
patching functor, and the maximal Cohen--Macaulay nature of
patched modules, that the formula of the
Breuil--M\'ezard conjecture precisely describes the Hilbert--Samuel
multiplicity of $R/\varpi$, where
$R$ denotes that quotient of $R_\infty(\sigma)$
that acts faithfully on $M_{\infty}(\sigma^{\circ})$.
Thus the Breuil--M\'ezard conjecture, when combined
with the weight part of Serre's conjecture,
is seen to be essentially equivalent to the faithfulness
of the action of $R_p(\sigma)$ on $M_{\infty}(\sigma^{\circ})$,
i.e.\ to a modularity lifting theorem.
(See \cite[\S 5]{emertongeerefinedBM} for a more precise formulation
of this equivalence in a more general setting.)

By combining this observation with an argument using the $p$-adic
Langlands correspondence, Kisin was in fact able to prove the
Breuil--M\'ezard conjecture (and so also to deduce corresponding
modularity lifting theorems) \cite{KisinFM}.  We recall more of
the details of this argument below.
For now, we merely remark that,
even without knowing whether or not $R_p(\sigma)$ acts faithfully
on $M_{\infty}(\sigma^{\circ})$, there is an obvious
inequality of Hilbert--Samuel multiplicities
$e_{R_p(\sigma)/\varpi} \geq e_{R/\varpi}$, and so Kisin's observation
always yields a weak form of the Breuil--M\'ezard conjecture,
in which equality is replaced by an inequality in one
direction.

\subsubsection{A locally analytic variant}
Above we discussed an approach to studying Serre weights via patching functors. There is a related story in the analytic setup which is discussed in \cite{MR3319547} and \cite{MR4028517}.
The question of which Serre weights are associated to a residual
Galois representation $\bar\rho$  can be understood as the question of
which irreducible ${\rm GL}_2(\Z_p)$-representations (more generally
$K={\rm GL}_n(\mathcal{O}_F)$-representations) embed into the
(conjectural) mod $p$  representation $\Pi(\bar\rho)$. (Note that here we are
adopting the more traditional perspective on $p$-adic and mod~$p$ local Langlands;
the connection between this perspective and the categorical conjecture
emphasised in these notes is discussed in Remark~\ref{rem: we expect a unique supersingular pi with skyscraper H0}.)
In some analogy with this question, Breuil \cite[Conj.\ 6.1]{MR3319547} has formulated a conjecture about the irreducible constituents of locally analytic principal series representation that embed into the $p$-adic local Langlands correspondence $\Pi(\rho)$ of a $p$-adic Galois representation $\rho$.
While the question about Serre weights is motivated by the global question about the weights of modular forms whose associated Galois representation has fixed reduction modulo $p$, Breuil's socle conjecture is motivated by the global question about  (the weights of) non-classical overconvergent $p$-adic modular forms with prescribed associated Galois representation.

In \cite{MR4028517} many cases of Breuil's socle conjecture are proven using a similar strategy (patching functors and a multiplicity formula in the spirit of the Breuil-M\'ezard conjecture) as discussed above: 
Instead of studying $V\mapsto M_\infty(V)$ as a functor on certain locally algebraic $K$-representations, one studies a related functor on 
certain locally analytic representations (or the BGG category $\mathcal{O}$ and a functor from that category to locally analytic representations). 
In this case the functor produces coherent sheaves on the rigid analytic generic fiber of ${\rm Spf}R_\infty$ that are supported on closed subvarieties defined in terms of trianguline Galois representations (see also Section \ref{subsubsec:loc an padic Langlands} below for a discussion of trianguline representations).
For these trianguline deformation spaces there is a locally analytic analogue of the Breuil-M\'ezard conjecture formulated in \cite[Conj.\ 4.3.4]{MR4028517}.
Roughly this conjecture predicts that the multiplicities of irreducible constituents in locally analytic representations that are parabolically  induced  from a character $\delta$ match the multiplicities of certain cycles in trianguline deformation spaces whose parameters correspond to $\delta$. In fact this conjecture can be proven in many cases \cite[Theorem 4.3.8]{MR4028517} by relating it to multiplicity formulas in geometric representation theory.

\subsubsection{The $p$-adic local Langlands correspondence for
  $\GL_2(\Q_p)$}\label{subsubsec:p-adic LL}
The $p$-adic local Langlands correspondence, associating
admissible continuous unitary $p$-adic Banach space
representations $\Pi(\rho)$
of $\GL_2(\Q_p)$ to two dimensional continuous representations
$\rho: \Gal_{\Q_p} \to \GL_2(L)$, %
was
first proposed by Breuil for those representations $\rho$ that
are potentially semistable
with distinct Hodge--Tate weights \cite{BreuilGL2I, BreuilGL2II, BreuilSemistable}.

More precisely, Breuil proposed in this context that $\Pi(\rho)$ should be
constructed as a certain completion of the locally algebraic representation
$\Pi(\rho)_{\alg}$ obtained by tensoring together the smooth representation
of $\GL_2(\Q_p)$ attached to $\WD(\rho)$ (the Weil--Deligne representation
underlying the Dieudonn\'e module of the potentially semistable
representation $\rho$) with a finite dimensional algebraic representation
of $\GL_2(\Q_p)$ encoding the Hodge--Tate weights of $\rho$. (The
representation~$\Pi(\rho)_\alg$ contains a locally algebraic
type~$\sigma$, and Breuil's proposal was motivated by the
Breuil--M\'ezard conjecture.)
Furthermore,
in the case when $\rho$ is crystabelline or semistabelline
(i.e.\ becomes crystalline,
or semistable, upon restriction to a finite abelian extension of $\Q_p$),
he made his proposal completely explicit, in the sense that he made
a specific proposal as to what $\Pi(\rho)$ should be:
namely, in the irreducible crystabelline case,
he proposed that $\Pi(\rho)$ should be the universal unitary completion
(in the sense of~\cite{MR2181093})
of $\Pi(\rho)_{\alg}$, while in the irreducible genuinely semistabelline  (i.e.\ non-crystabelline)
case, he proposed that $\Pi(\rho)$ should be a certain
completion
of $\Pi(\rho)_{\alg}$ depending on the $\CL$-invariant
of~$\rho$.
(As we elaborate on below, there is still an inexplicit aspect to this proposal,
in so far as it is not obvious how to actually compute 
these universal unitary completions in any concrete fashion.)

The motivation for proposing the universal unitary
completion as the candidate for $\Pi(\rho)$
in the irreducible crystabelline case is that in this
case there is no ``extra'' $p$-adic Hodge theoretic data carried
by $\rho$ besides that already carried by $\Pi(\rho)_{\alg}$
(that is, there is no $\CL$-invariant or the like in this case; more
precisely, there is up to isomorphism a unique weakly admissible
filtration with fixed Hodge--Tate weights on $D_{\rm cris}(\rho)$ such
that the corresponding Galois representation is irreducible). The
representation $\Pi(\rho)_{\rm alg}$ determines (and is determined by)
$D_{\rm cris}(\rho)$ and the Hodge--Tate weights, 
and the universal unitary completion is the only evident choice
of completion available that doesn't depend on any additional
choices. %

The major initial difficulty in investigating
Breuil's conjectured correspondence
was in establishing any non-trivial properties of his proposed completions;
it was not even evident that these completions were non-zero.
One exception to this situation was in the case when $\rho$
was the restriction to $\Gal_{\Q_p}$ of a global two-dimensional representation
$r: \Gal_{\Q} \to \GL_2(L)$ arising from a Hecke eigenform.
In this case the theory of completed cohomology, together with
local-global compatibility for the classical Langlands correspondence,
gives rise to a non-zero unitary completion of $\Pi(\rho)_{\alg}$.
In the semistabelline case, Breuil was furthermore able to show, using
methods arising from the theory of $p$-adic $L$-functions
(and we remind the reader that the notion of $\CL$-invariant had
its origin in that same theory)
that this completion factored through his proposed completion
of $\Pi(\rho)_{\alg}$ \cite{breuil-l-invariant2}.

Inspired by the approach of Kato, Perrin-Riou and
others to the theory of $p$-adic $L$-functions via
Galois cohomology, explicit reciprocity laws, and $p$-adic Hodge theory,
and his own reinterpretation and further developments
of these ideas in terms of $(\varphi,\Gamma)$-modules,
Colmez was then able to make a major breakthrough, and establish,
in the irreducible and genuinely semistabelline case, that Breuil's proposed
completion $\Pi(\rho)$ is non-zero, topologically irreducible, and admissible.
(See~\cite{MR2642407} for the culmination of Colmez's investigations in
this direction.)
Soon thereafter, using similar techniques, Berger and Breuil
were able to prove the analogous result in the crystabelline case~\cite{MR2642406}.  A key consequence of their result is that, for crystabelline $\rho$,
the universal unitary completion $\Pi(\rho)$ of $\Pi(\rho)_{\alg}$
is the unique non-zero unitary completion of $\Pi(\rho)_{\alg}$.

Since the theory of $(\varphi,\Gamma)$-modules applies to arbitrary
representations $\rho$, not just potentially semistable ones, Colmez's results
led him to propose that the $p$-adic local Langlands correspondence
should exist for arbitrary continuous representations
$\rho: \Gal_{\Q_p} \to \GL_2(L)$.  (In the case when $\rho$ is not potentially
semistable with distinct Hodge--Tate weights, one no longer has
a representation $\Pi(\rho)_{\alg}$ available from which to even
conjecturally construct $\Pi(\rho)$ as some kind of completion.)
Reverse-engineering his $(\varphi,\Gamma)$-module approach
to studying Breuil's completions, Colmez discovered the
first of his two functors, the functor $\Pi \mapsto V(\Pi)$
mapping  topologically finite length admissible unitary continuous
Banach $\GL_2(\Q_p)$-representations (satisfying
a finite length condition on their reductions mod $\varpi$,
which was shown by Pa\v{s}k\={u}nas to be superfluous \cite{MR3150248})
to continuous representations of $\Gal_{\Q_p}$ \cite[\S IV]{MR2642409}.

In the light of his investigations of the Breuil--M\'ezard conjecture, Kisin then made the crucial suggestion that Colmez's
proposed
general $p$-adic local Langlands correspondence should be regarded
as taking place over local Galois deformation spaces, and that
deformation-theoretic tools could play a role in its construction.
Following this suggestion, Colmez was able to extend
his $(\varphi,\Gamma)$-module methods to construct
$\Pi(\rho)$ for general $\rho$ \cite{MR2642409}.
The functors are related via a natural isomorphism
$V\bigl(\Pi(\rho)\bigr) \cong \rho.$
Furthermore, he was able to prove that
when $\rho$ is potentially semi-stable with distinct Hodge--Tate weights,
the representation $\Pi(\rho)$ is indeed a completion of $\Pi(\rho)_{\alg}$,
as Breuil had originally proposed.  In the cases when $\rho$ is
crystabelline or semistabelline, this followed from the results
already mentioned due to him and Berger--Breuil (indeed, the
general construction of $\Pi(\rho)$ is motivated by the constructions
in those cases).  However, in general (i.e.\ in cases
that are not crystabelline, or semistabelline, or, more generally,
trianguline), the action of the entirety of $\GL_2(\Qp)$
on $\Pi(\rho)$ is constructed
in a very indirect way (it is only the action of the upper
triangular Borel subgroup of $\GL_2(\Q_p)$ that can be
seen directly in terms of $(\varphi,\Gamma)$-modules) by deformation-theoretic
methods (or, one could say, by rigid analytic interpolation),
and so the nature of the locally algebraic vectors in $\Pi(\rho)$
is not obvious.  Colmez's proof that they contain (and in fact
coincide with) $\Pi(\rho)_{\alg}$ ultimately rests on a comparison
with the global context of completed cohomology, mediated
via the work of one of us (M.E.) on local-global compatibility.

The theory of the $p$-adic local Langlands correspondence for
$\GL_2(\Q_p)$ was completed by Pa\v{s}k\={u}nas,
who showed that the functor $\Pi \mapsto V(\Pi)$ induces
a bijection between isomorphism classes of supersingular
topologically irreducible admissible unitary continuous
Banach space representations of $\GL_2(\Q_p)$, and
irreducible continuous two-dimensional representations of $\Gal_{\Q_p}$
\cite{MR3150248}.
In fact, more was established: essentially, $V$ induces an equivalence
between the category of topologically
finite length admissible
Banach space representations of $\GL_2(\Q_p)$,
and the category
whose objects are deformations of two-dimensional continuous
representations of $\Gal_{\Q_p}$.  Furthermore, this equivalence can
be described Morita-theoretically: there are certain injective objects
in the category of $\GL_2(\Q_p)$-representations, whose endomorphism
rings may be identified with Galois deformation rings, that mediate
the equivalence.  Since these injective objects (or, equivalently,
their topological duals, which are then projective objects in an
appropriate category) play an important role in the work of Andrea
Dotto, M.E.\ and T.G.\ which we explain in
Sections~\ref{sec:banach-case-gl_2qp}--\ref{subsec:examples},
we will discuss them further below.

\subsubsection{The locally analytic $p$-adic Langlands correspondence}\label{subsubsec:loc an padic Langlands}
 As we will discuss further below
  below, the Banach space representations occurring in the $p$-adic
  local Langlands correspondence for~$\GL_2(\Qp)$ are ultimately
  constructed by passage to the limit of constructions modulo~$p^n$
  (and then by inverting~$p$), and in particular they naturally fit
  into the ``Banach'' context in these notes.
  
  On the other hand,
  an important tool for studying Banach space representations such as
  $\Pi(\rho)$ is to study their space of analytic vectors $\Pi(\rho)^{\rm
    an}$ (in the sense of \cite[\S7]{MR1990669}), which are locally
  analytic representations of ${\rm GL}_2(\Q_p)$.
  This goes in the opposite direction to the consideration of mod $p^n$ representations,
  since $\Pi(\rho)^{\rm an}$ forgets about the Banach space (and hence  $p$-adic integral) structure 
  on~$\Pi(\rho)$.
  On the  other hand, Colmez's
  construction shows that the locally analytic representation
  $\Pi(\rho)^{\rm an}$ can be constructed directly from the rigid
  analytic $(\varphi,\Gamma)$-module $D_{\rm rig}(\rho)$ over the Robba ring $\mathcal{R}$ associated to~$\rho$,
  so that it still  makes sense to study the $p$-adic local Langlands
  correspondence in this context.
  Motivated by this observation Colmez \cite{MR3522263} extended
  the $p$-adic local Langlands correspondence for ${\rm GL}_2(\Q_p)$
  to the setup of (possibly non-\'etale) rank two
  $(\varphi,\Gamma)$-modules over the Robba ring and locally analytic
  representations of ${\rm GL}_2(\Q_p)$. This naturally fits into the
  ``analytic'' context in these notes.
  
 We point out that usually $p$-adic Hodge-theoretic properties of the
 representation $\rho$ are not so easy to read off from the associated
 $(\varphi,\Gamma)$-module $D(\rho)$, but it is much easier to read them
 off from the rigid analytic variant $D_{\rm rig}(\rho)$. Hence
 Hodge-theoretic information, e.g.~the Hodge filtration on $D_{\rm
   dR}(\rho)$ (assuming that $\rho$ is de Rham) should be connected to
 properties of the locally analytic representation $\Pi(\rho)^{\rm
   an}$ (and it should not be easy to understand this directly in
 terms of the Banach space representation $\Pi(\rho)$). We now
 elaborate on this point.
 
 One of the original motivations in the $p$-adic Langlands program was to improve the understanding of local-global compatibility in the cohomology of modular curves (or more general Shimura varieties).
We will discuss the $p$-adic local global compatibility involving completed (co-)homology of modular curves in more detail in Section \ref{sec:patching local-global-compatibility} below,
but for now we note that
 to each system of Hecke eigenvalues that occurs in the cohomology
 $H^1(Y(N),\mathbf{Q}_p)$ of the modular curve $Y(N)$ there corresponds a continuous $2$-dimensional representation $\rho$ of ${\rm Gal}_{\Q}$ on a $\Q_p$-vector space.
The restriction $\rho|_{{\rm Gal}_{\Q_\ell}}$ to decomposition groups
at primes $\ell\neq p$ is completely determined (at least up to
Frobenius semisimplicity, which is expected to be automatic) %
by the eigenvalues of the Hecke operators at $\ell$. However the action of the Hecke operators at $p$ only determines the Weil--Deligne representation associated to the de Rham representation  $\rho|_{{\rm Gal}_{\Q_p}}$. The missing piece in determining $\rho|_{{\rm Gal}_{\Q_p}}$ is the Hodge filtration, which is mysteriously determined by the Hecke action away from $p$ but not captured by the Hecke action at $p$.

One of
Breuil's initial ideas, then, was that $p$-adic local global compatibility involving completed cohomology and the $p$-adic local Langlands correspondence should make it possible to read off the Hodge filtration from the ${\rm GL}_2(\Q_p)$-representation on the corresponding
Hecke eigenspace in
completed cohomology $\widetilde H^1(Y(N),\mathbf{Q}_p)$ (as defined in~\eqref{eqn:completed cohom def} below)
without using the action of Hecke operators away from $p$. In fact, the Hodge filtration should be determined by the locally analytic vectors in this eigenspace in $\widetilde H^1(Y(N),\mathbf{Q}_p)$. 
A generalization of this philosophy %
is visible for example in Breuil's socle conjecture (see Section~\ref{sec:companionpoints}), where the relative position of the Hodge filtration with respect to Frobenius stable flags is determined by the socle of a certain locally analytic representation. 

Passing from a (classical) $(\varphi,\Gamma)$-module $D$ to its variant $D_{\rm rig}$ over the Robba ring $\mathcal{R}$ has the advantage that the $(\varphi,\Gamma)$-module can become reducible, even though $D$ was not reducible. For example this happens if $D$ comes from a crystalline or semi-stable (or more generally crystabelline or semistabelline) representation. 
The $(\varphi,\Gamma)$-modules that are completely reducible are called trianguline and were studied extensively by Colmez \cite{MR2493219}. In terms of the $p$-adic local Langlands correspondence a representation $\rho$ is trianguline (i.e.~$D_{\rm rig}(\rho)$ is trianguline) if and only if the locally analytic vectors $\Pi(\rho)^{\rm an}$  contain (irreducible constituents of) locally analytic principal series representations.
In this case $\Pi(\rho)^{\rm an}$ is typically an extension of two parabolically induced representations, and the extension class of these representations is determined by the extension class of $(\varphi,\Gamma)$-modules of rank~$1$ defined by $D_{\rm rig}(D)$ (and vice versa).
This implies in particular that in the case of crystalline or semi-stable representations $\rho$ the Hodge filtration on the Weil--Deligne representation ${\rm WD}(\rho)$ is determined by the locally analytic vectors $\Pi(\rho)^{\rm an}$, as it should be according to the above discussion.

\subsubsection{The mod $p$ local Langlands correspondence
	for $\GL_2(\Q_p)$}
In our discussion so far, we have set the $p$-adic local Langlands
correspondence in the context of {\em $p$-adic} representations
of $\GL_2(\Q_p)$.  However, from the beginning of the theory,
 mod $p$ considerations have played a fundamental role.
In his paper \cite{BreuilGL2I} (and building
on earlier work of Barthel and Livn\'e \cite{BarthelLivneDuke}),
Breuil classified the absolutely
irreducible admissible smooth representations of $\GL_2(\Q_p)$
over a finite field~$k$,
and showed {\em grosso modo} that they are of the form
$\bigl(\cInd_{\GL_2(\Z_p)Z}^{\GL_2(\Q_p)} \sigmabar\bigr)
\otimes_{\cH(\sigmabar),x}k,$
where $Z$ denotes the centre of $\GL_2(\Q_p)$,  $\cH(\sigmabar) =
\End(\bigl(\cInd_{\GL_2(\Z_p)Z}^{\GL_2(\Q_p)} \sigmabar\bigr)$
is the relevant Hecke algebra,
which admits an explicit isomorphism $\cH(\sigmabar) \cong k[S^{\pm 1}, T]$,
and where $x: \cH(\sigmabar) \to k$ is some given character of
$\cH(\sigmabar).$
The key contribution
of \cite{BreuilGL2I} is to show that such representations
are irreducible in the so-called {\em supersingular} case,
i.e.\ the case when $x(T) = 0$ (and it is this statement
that is special to the context of $\GL_2(\Q_p)$).

Building on this classification, Breuil was able to define
a {\em mod $p$ semisimple local Langlands correspondence},
by attaching to each semisimple two-dimensional representation
$\rhobar: \Gal_{\Q_p} \to \GL_2(k)$ a corresponding
semisimple admissible smooth
representation $\kappa(\rhobar)$ of $\GL_2(\Q_p)$ over $k$.
Breuil further conjectured that this correspondence was
compatible with the conjectured $p$-adic local Langlands correspondence,
in the sense that if $\rho: \Gal_{\Q_p} \to \GL_2(L)$ is continuous,
then one should have $\overline{\Pi(\rho)}^{\semis} \cong \kappa(\rhobar^{\semis}).$
(Here, on either side of the conjectured isomorphism,
the notation $\overline{(\text{--})}^{\semis}$ indicates
that we reduce an invariant lattice modulo the uniformizer $\varpi$,
and then semisimplify.
Note also that on either side,
the reductions are only uniquely determined
after semi-simplifying, and so it is natural to pass to semi-simplifications
when making such a comparison between the $p$-adic and mod $p$ settings.)

Subsequently, Berger \cite{MR2642408}
was able to obtain a description of Breuil's
 mod~$p$ correspondence in terms of $(\varphi,\Gamma)$-modules,
and using this, was able to verify Breuil's conjectured compatibility
in the context of the paper \cite{MR2642406}.
The same $(\varphi,\Gamma)$-module
arguments apply to show this compatibility naturally extends to the
general $p$-adic local Langlands correspondence of \cite{MR2642409}.

Colmez was subsequently able to define a fully fledged
{\em mod $p$ local Langlands correspondence}, i.e.\ to define
the correspondence in such a fashion that, in the case
when $\rhobar$ is reducible but indecomposable, it takes into
account the corresponding extension class of characters.
More precisely, for each
representation $\rhobar: \Gal_{\Q_p} \to \GL_2(k)$, one may define
an associated representation
$\kappa(\rhobar)$ of $\GL_2(\Q_p)$, with the property that
$V\bigl(\kappa(\rhobar)\bigr) = \rhobar$.
If $\rhobar$ is {\em not} a twist of an extension
of the mod $p$ cyclotomic character by the trivial character,
then $\kappa(\rhobar)$ can in fact be characterized by this last
property, together with the condition that it admits no
finite-dimensional $\GL_2(\Q_p)$-invariant subobject or
quotient.  (See e.g.\ \cite[Thm.~3.3.2]{emerton2010local}.)

In the exceptional case,
when $\rhobar$ {\em is} a twist of an extension
of the mod $p$ cyclotomic character by the trivial character,
then the ``correct'' definition of $\kappa(\rhobar)$ (by which we mean
the one that is compatible with the correspondence that we explain in Section~\ref{sec:banach-case-gl_2qp}) %
is {\em not} compatible with semi-simplification.
In fact, even in this case, it is possible to alter the definition
of $\kappa(\rhobar)$ so that it is so compatible, and Colmez
did indeed adopt this alternative definition; but this alternative choice
of $\kappa(\rhobar)$ fails to
satisfy %
compatibility with deformations, in the form  discussed in more detail
below. (See~\cite[Rem.\ 6.22]{MR3732208} for a more precise
description of the difference between Colmez's definition and ours.) %

\subsubsection{$p$-adic local Langlands over deformation space
and proofs of the Breuil--M\'ezard conjecture}
\label{subsubsec:padicLL over deformation space}
As already remarked, the construction of the  $p$-adic local
Langlands correspondence for~$\GL_2(\Qp)$ is intimately related to describing the
variation of the correspondence over deformation spaces
of Galois representations, and we now turn to discussing
the various forms this description has taken in the literature.
The first key input to all of these
descriptions is that Colmez's functor $V$ is well-defined on
the mod $p$ and $p$-adic integral level, so that it can be utilized
to study the $p$-adic local Langlands correspondence
in a deformation-theoretic context.

In \cite{emerton2010local}, a deformation-theoretic
description of the correspondence is given for deformations
of a representation $\rhobar: \Gal_{\Q_p} \to \GL_2(k)$ that
is {\em not} a twist of an extension of cyclotomic character by the
trivial character; in what follows we assume
for simplicity that $\rhobar$ also admits only trivial endomorphisms,
so that we can and do take $R_p$ to be the universal deformation
$\cO$-algebra of~ $\rhobar$ with some fixed determinant (rather than
having to work with framed deformations). For simplicity of
exposition, we furthermore revert to our earlier assumption that this fixed determinant
is~$\varepsilon^{-1}$; correspondingly, we will now assume that all of the
representations of~$\GL_2$ that we consider have trivial central
character, and accordingly we will regard them as representations of~$\PGL_2$.

In this case, the arguments of Colmez \cite{MR2642409}
and Kisin~\cite{MR2642410} can be encapsulated as follows: there exists
a representation  $\pi^{\univ}$ %
of $\PGL_2(\Q_p)$ on an
$\mathfrak m_{R_p}$-adically complete and orthonormalizable $R_p$-module,
deforming the smooth representation $\kappa(\rhobar)$,
with the property that $V(\pi^{\univ}) \cong \rho^{\univ}$,
the universal deformation of $\rhobar$.
The $\PGL_2(\Q_p)$-representation attached by the $p$-adic local Langlands
correspondence to any particular deformation of $\rhobar$ is then
obtained by specializing $\pi^{\univ}$ over this deformation.

In \cite{MR3150248}, a deformation-theoretic
description of the correspondence is given that is more flexible
and more conceptual than the formulation just recalled.
In order to make the description as simple as possible,
we assume that $\rhobar$ %
is {\em not} a twist of an extension of the trivial
character by the cyclotomic character.
(Note that this is different to the exceptional case considered
above; the forbidden extension has the characters in the opposite
order.)
We then
let $\tP$ denote a projective envelope of the Pontryagin dual %
      $\kappa(\rhobar)^{\vee}$
in a suitable category of representations of $\PGL_2(\Q_p)$ on compact
$\cO$-modules. 
Let a superscript ${}^d$ denote Schikhof duality on pseudocompact
$\cO$-torsion free $\cO$-modules, i.e.\ 
 $\Hom_\cO^\cont(-,\cO)$. %
One of the main results of \cite{MR3150248}
is that there is a canonical isomorphism  %
\numequation\label{eqn: VP identifying deformation ring with
  endomorphism of projective cover}R_p \cong \End(\tP),\end{equation}%
so that $\tP$ may be naturally regarded as ``living over''
$\Spf R_p$, %
and that the consequent $R_p$-structure on $V(\tP^d)$ realizes
$V(\tP^d)$ as the universal deformation $r^{\univ}$ of $\rhobar$.
If $\rho: \Gal_{\Q_p} \to \GL_2(L)$ admits a lattice deforming $\rhobar$,
corresponding to a morphism $x:R_p \to \cO,$
then the unitary Banach space representation $\Pi(\rho)$
is determined via the isomorphism
$$\Pi(\rho) \cong L\otimes_{\cO} (\tP\otimes_{R_p,x} \cO)^d.$$ 

If $\rhobar$ is not a twist of an extension of the cyclotomic character
by the trivial character (so that we are in neither of the exceptional
cases considered above), then it is furthermore proved in \cite{MR3150248}
that $\tP$ is topologically flat
over $R_p$ (in the sense that the functor $-\widehat{\otimes}_{R_p}\tP$ on
pseudo-compact $R_p$-modules is exact).  In this case, the representation
$\pi^{\univ}$ associated to $\rhobar$ as above %
is related
to $\tP$ via the isomorphism (of $R_p[\PGL_2(\Q_p)]$-modules)
$$ \pi^{\univ} \cong \Hom^{\cont}_{R_p}(\tP,R_p).$$
	
In \cite{KisinFM}, Kisin gives quite a different description
of the $p$-adic local Langlands correspondence, not over
all of $\Spf R_p$, but over the locus $\Spf R_p(\sigma)$
parameterizing lifts of a fixed locally algebraic type $\sigma$.
For deformations $\rho$ of $\rhobar$ that are parameterized by $\Spf R_p(\sigma)$,
Colmez's proof of Breuil's conjecture describing
$\Pi(\rho)$ as a completion of the locally algebraic representation
$\Pi(\rho)_{\alg}$ shows that
the restriction to $\Spf R_p(\sigma)[1/p]$ of the $p$-adic local
Langlands correspondence, which we denote
by $\pi(\sigma)$, %
may be described as a certain
completion of
$\cInd_{\PGL_2(\Z_p)}^{\PGL_2(\Q_p)} \sigma.$
This description of $\pi(\sigma)$
allows Kisin to deduce the Breuil--M\'ezard conjecture
for~$\rhobar$, as we now recall.

A choice of Jordan--H\"older filtration
in the reduction $\sigmabar$ of some chosen $\PGL_2(\Z_p)$-invariant
lattice $\sigma^{\circ}$ in $\sigma$ determines
a corresponding filtration on $\overline{\pi(\sigma)}$
(the reduction mod $\varpi$ of some invariant lattice in $\pi(\sigma)$)
by completions of representations of the form
$\cInd_{\PGL_2(\Z_p)}^{\PGL_2(\Q_p)} \sigmabar_i$,
where $\sigmabar_i$ runs through the Jordan--H\"older factors of
$\sigmabar$. For each~$\sigmabar_i$, we let~$\sigma_i$ be an algebraic
$L$-representation of~$\PGL_2(\Qp)$
with~$\sigma^\circ_i\otimes_{\cO}k\cong\sigmabar_i$, and
set~$R_p(\sigmabar_i)=R_p(\sigma_i)\otimes_{\cO}k$. %
Now the explicit classification of the irreducible representations
of $\PGL_2(\Q_p)$, and the fact that Colmez's functor $V$ is compatible
with reduction, shows that the only $\sigmabar_i$ that can contribute
a non-zero completion are those that appear as Serre weights of $\rhobar$,
and that if a given Serre weight $\sigmabar_i$ contributes a completion,
this completion must be the $\mathfrak m_x$-adic  completion
of $\cInd_{\PGL_2(\Z_p)}^{\PGL_2(\Q_p)} \sigmabar_i,$
where $\mathfrak m_x$ is the kernel of the character $x: \cH(\sigmabar_i) \to k$
for which
$\bigl(\cInd_{\PGL_2(\Z_p)}^{\PGL_2(\Q_p)} \sigmabar_i\bigr) \otimes_{
	\cH(\sigma),x} k$ is a constituent of $\kappa(\rhobar)$.
Furthermore, this completion, when thought of as an $R_p(\sigma)$-module,
in fact is supported on $R_p(\sigmabar_i)$, %
the $\mathfrak m_x$-adic completion of $R_p(\sigmabar_i)$ being
explicitly identified with the $\mathfrak m_x$-adic completion
of $\cH(\sigmabar_i)$. %
	       	In this way, one obtains an inequality
of Breuil--M\'ezard type, but in the opposite direction
to the one given by the argument described in~(\ref{subsubsec:BM}) above.
Combining these two inequalities gives an equality, which
is the Breuil--M\'ezard conjecture.
As noted in~(\ref{subsubsec:BM}), this conjecture
then implies that each patched module $M_{\infty}(\sigma^{\circ})$
is faithful over the ring $R_p(\sigma)$,
giving rise to a modularity lifting theorem for global deformations
of type $\sigma$.  This is the approach to the Fontaine--Mazur
conjecture taken in \cite{KisinFM}.

In the argument that is sketched in the preceding paragraph,
the reason that only an inequality (rather than the precise equality
of the Breuil--M\'ezard conjecture) is obtained is more-or-less
because completions (being a form of projective limit) are not {\em
	a priori} right exact, so that is it not clear {\em a priori}
that all the possible completions
of the representations
$\cInd_{\PGL_2(\Z_p)}^{\PGL_2(\Q_p)} \sigmabar_i$
that might occur in $\overline{\pi(\sigma)}$ actually do so occur.

Issues of exactness of this kind are ultimately of a homological nature,
and so it is perhaps not surprising that having stronger homological
information available allows the argument to be strengthened.
This was achieved by Pa\v{s}k\={u}nas in the paper~\cite{paskunasBM}.
Using the construction of $\tP$ as described above, together
with Colmez's result on locally algebraic vectors, one finds
that $R_p(\sigma)$ is precisely the support of the $R_p$-module
$\Hom_{\PGL_2(\Z_p)}(\sigma^{\circ}, \tP^{\vee})$ (for any choice
of invariant lattice $\sigma^{\circ}$ in $\sigma$).  In fact,
it is more convenient to work with the $R_p$-module %
\numequation\label{eqn: M of sigma}M(\sigma^{\circ}) := \tP\otimes_{\cO[[\PGL_2(\Z_p)]]}\sigma^{\circ},\end{equation}
since this turns out to be a coherent $R_p$-module.
Since $\tP$ is projective as a $\PGL_2(\Q_p)$-representation,
it is in particular projective as a $\PGL_2(\Z_p)$-representation,
and so the functor $M(\text{--})$ implicitly defined by the
preceding formula is {\em exact}.   Arguing with this functor
one can essentially combine the two prongs of Kisin's argument
into a single argument that yields the Breuil--M\'ezard
equality directly.

\subsubsection{Local-global compatibility}\label{sec:patching local-global-compatibility}
In \cite{emerton2010local} a local-global compatibility relating
the $p$-adic local Langlands correspondence to completed cohomology
is proved, giving a $p$-adic analogue of the classical local-global
compatibility \cite{MR0354617, delignetoPS, MR705677}
relating cohomology at infinite level
and the classical local Langlands correspondence.
In order to state the result, we fix 
as above a modular irreducible
residual representation $\rbar:\Gal_{\Q} \to \GL_2(k)$, setting~$\rhobar := \rbar|_{\Gal_{\Qp}}$,
and
work at a tame level $N$ divisible by the prime-to-$p$
Artin conductor of~$\rbar$.

Similar to the completed homology (\ref{eqn:completed homology modular curves}) we let $\tH^1 := \tH^1_{\et}\bigl( Y(N), \cO \bigr)$ denote the $p$-adically
completed \'etale cohomology of the ``modular curves'' (again: more
precisely, congruence quotients of  $\PGL_2$-symmetric spaces, which
are unions of the connected components of the usual modular curves for $\GL_2$) %
of tame level~$N$; 
      by
      definition we have
\numequation
\label{eqn:completed cohom def}
\tH^1:=\varprojlim_s\varinjlim_n H^1_\et\bigl(Y(Np^n),\cO/\varpi^s\bigr),
\end{equation}
and~$\tH^1$ is naturally equipped with commuting actions of the groups $\Gal_{\Q}$ and $\PGL_2(\Q_p)$,
and of a Hecke algebra $\T$ (generated by the usual Hecke
operators~$S_l,T_l$ at primes~$l\nmid Np$).  The representation $\rbar$
gives rise to a maximal ideal $\mathfrak m$ of $\T$,
and we may form the corresponding localisation $\tH^1_{\mathfrak m}$
of $\tH^1$.  There is also a universal pro-modular deformation
$\rho^{\univ}$ of $\rbar$ over $\T_{\mathfrak m},$
and if we set $\rhobar := \rbar_{|\Gal_{\Q_p}}$, then restriction
from $\Gal_{\Q}$ to $\Gal_{\Q_p}$ induces a morphism
$R_p \to \T_{\mathfrak m}$.

Assume now that $\rhobar$ is not a twist of an extension of the cyclotomic
character by the trivial character; the main result of \cite{emerton2010local}
then states that there is an isomorphism (of $\T_{\mathfrak m}[
\Gal_{\Q}\times \PGL_2(\Q_p)]$-modules)
\numequation\label{eqn:local global compatibility H^1}
\tH^1_{\mathfrak m} \cong ( \rho^{\univ} \otimes_{R_p}
\pi^{\univ})\ \cotimes_{R_p} O, %
\end{equation}
where $\pi^{\univ}$ is the orthonormalizable
$\PGL_2(\Q_p)$-representation over $R_p$ discussed in Section~\ref{subsubsec:padicLL over deformation space} above, and where $O$ is a $p$-torsion free cofinitely generated
$\T_{\mathfrak m}$-module\footnote{I.e.\ the $\cO$-dual of a finitely generated
$\T_{\mathfrak m}$-module.} of ``oldforms''
related to the chosen tame level (regarded as an $R_p$-module via the map
$R_p\to \T_{\mathfrak m}$).\footnote{If we work at full level $N$, and then pass to a colimit
over all levels~$N$, the resulting limiting space of oldforms acquires
an action of $\GL_2(\A_f^p)$, and can be described in terms of the
local Langlands correspondence
in families of \cite{MR3250061} and~\cite{MR3867634}.  This is related to the $\ell \neq  p$
case of categorical local Langlands that we discuss briefly in
Section~\ref{sec: categoric LL for l not p},
and which plays a key role in the global conjecture  (Conjecture~\ref{conj:cohomology} below).
For simplicity of exposition,  though, we don't say anything more about it here.}
If $\rhobar$ is furthermore not a twist of an extension of the trivial
character by the cyclotomic character, then passing to continuous $\cO$-duals, %
and taking into account the relationship between $\pi^{\univ}$ and $\tP$
in this case,
we may rephrase~\eqref{eqn:local global compatibility H^1} as an isomorphism %
\numequation\label{eqn: local global compatibility GL2Qp}(\tH_1)_{\mathfrak m} \cong \tP\cotimes_{R_p}
(\rho^{\univ})^{\vee}\otimes_{R_p} O^*, \end{equation} %
where $\tH_1 = \Hom_{\cO}(\tH^1,\cO)$ denotes the completed homology of
the modular curves as in (\ref{eqn:completed homology modular curves}) at tame level $N$,
$(\rho^{\univ})^{\vee}$
denotes the $R_p$-dual of~$\rho^{\univ}$,
and $O^*$ denotes the $\cO$-dual of~$O$, which is now a finitely generated $\T_{\mathfrak m}$-module.
The proof of the isomorphism is a combination of a Morita-theoretic argument,
similar to the more sophisticated such arguments that are developed
in~\cite{MR3150248}, and an interpolation argument, related to
the ``capture'' arguments that appear in \cite{MR3272011}.
The key input is a combination of classical local-global compatibility
and the result of Berger--Breuil (that $\Pi(\rho)$ is the unique
unitary completion of $\Pi(\rho)_{\alg}$ when $\rho$ is crystabelline):
classical local-global compatibility gives rise to copies
of $\Pi(\rho)_{\alg}$ inside $\tH^1$, and then the result of Berger--Breuil
shows that the closures of these in $\tH^1$
are necessarily isomorphic to~$\Pi(\rho)$.

Local-global compatibility gives another approach to the Fontaine--Mazur
conjecture.  Namely, if $x: \T_{\mathfrak m} \to L$
is a homomorphism corresponding to a pro-modular deformation
$r$ of $\rbar$, with kernel $\mathfrak p_x$,
and if $\rho := r_{| \Gal_{\Q_p}}$,
then passing to $\mathfrak p_x$-torsion in the local-global
compatibility isomorphism \eqref{eqn: local global compatibility GL2Qp} yields an isomorphism
\numequation\label{eqn: padic local global compatibility }
\tH^1[\mathfrak p_x]\otimes_{\cO} E \cong r\otimes \Pi(\rho)
\end{equation}
(where we have omitted the contribution from the oldforms; if we choose the
tame level to coincide with the prime-to-$p$ conductor of~$r$, then their 
contribution will indeed be trivial).
If $\rho$ is potentially semi-stable with distinct Hodge--Tate weights,
so that by Colmez's results $\Pi(\rho)_{\alg} \neq 0$,
then one concludes that the system of Hecke eigenvalues
given by $x$ appears in $\tH^1_{\alg}.$
The comparison between locally algebraic vectors in completed cohomology
and the classical cohomology of local systems on modular curves
then shows that $x$, and thus $r$, arises from a classical
modular form. 

The Fontaine--Mazur conjecture for deformations of~$\rbar$ therefore
reduces to the problem of showing that all such deformations are
pro-modular. Under standard hypotheses on~$\rbar$ needed for the
Taylor--Wiles method, this can be accomplished by proving a ``big
$R=\T$'' theorem, identifying~$\T_{\m}$ with the universal
deformation $\cO$-algebra for deformations of~$\rbar$ which are
minimally ramified away from~$p$. Such theorems were first proved
in~\cite{MR1854117}, using the infinite fern of Gouv\^ea and Mazur
together with modularity lifting theorems in small weight (where the
deformation rings~$R_p(\sigma)$ can be described explicitly via
Fontaine--Laffaille theory). Moreover, we note  that (\ref{eqn: padic local global compatibility }) and the construction of the $p$-adic Langlands correspondence imply that $\tH^1[\mathfrak p_x]\otimes_{\cO} L$ is non-zero for every $x: \T_{\mathfrak m} \to L$,
which is an {\em a priori} stronger statement 
than the fact that the associated deformations~$r$ are pro-modular.

The results of~\cite{emerton2010local} therefore give an
alternative proof of the results of~\cite{KisinFM} (with slightly
different hypotheses); while both proofs rely on Taylor--Wiles--Kisin
patching and the $p$-adic Langlands correspondence, they differ in
that the approach of~\cite{KisinFM} is to use the $p$-adic Langlands
correspondence to  prove enough about the ring~$R_p(\sigma)$ to
deduce that
the patched $R_\infty(\sigma)$-module~$M_\infty(\sigma^\circ)$ is
faithful, while the approach of~\cite{emerton2010local} can be viewed
as computing the functor~$M_\infty$ in terms of the $p$-adic Langlands
correspondence. %

The rephrasing~\eqref{eqn: local global compatibility GL2Qp} of the
main result of~\cite{emerton2010local} in terms of completed homology
is highly suggestive of two things. Firstly, from the optic of these
notes, where we replace Galois deformation rings with moduli stacks of
$L$-parameters, it suggests that we can compute the completed
cohomology of modular curves by pulling back some canonical sheaf of
$\PGL_2(\Qp)$-representations from a stack of local parameters to a
stack of global parameters; we pursue this suggestion in Section~\ref{sec: global stacks and cohomology of Shimura
varieties}. Secondly, it suggests the 
possibility of using the Taylor--Wiles patching method to extract~$\tP$ from the completed homology~$\tH_1$;
this possibility was realized in~\cite{MR3732208}. In particular, one
finds that the patching functor $M_\infty$ ``is'' the functor~$M$
implicitly described by~\eqref{eqn: M of sigma}; precisely, it coincides with
the extension of scalars of $M$ along the morphism $R_p \to R_{\infty}$. (Note, though, that
it seems  hopeless to construct the sheaf~$L_\infty$ of
Section~\ref{sec:banach-case-gl_2qp} in this way, or from any other
global construction.)

\subsubsection{Comparison to the proof of classical local
  Langlands}\label{subsubsec: comparison to classical LL}Our construction of the $p$-adic local Langlands correspondence
via global means may be compared with the construction
of the classical local Langlands correspondence, which also
proceeds by realising a local representation
globally~\cite{ht,MR1738446,scholze_local_LL}. However, as discussed
in~\cite[\S 1.3]{Gpatch}, there are some significant differences
between our construction and the classical one. Classically, one first uses the classification
of Weil--Deligne representations, and the parallel classification
of irreducible smooth representations of $\GL_n(\Q_p)$~\cite{BZ77ENS,Zel80}
to reduce to constructing a bijection between irreducible $n$-dimensional Weil--Deligne
representations and cuspidal representations of $\GL_n(\Q_p)$.
A trace formula argument allows one to realize any cuspidal
representation (up to twist)
as the local component of a cuspidal automorphic representation,
and the associated Weil--Deligne representation is then constructed
in the cohomology of a Shimura variety (in particular, this construction
{\em presumes} that local-global compatibility will hold, just as
that of~\cite{MR3732208} does).  In order to verify that the correspondence
so constructed is in fact purely local, one ultimately relies
on the prescription
for the local Langlands correspondence in terms of $L$- and $\varepsilon$-
factors.

In the $p$-adic context, most of these steps aren't available (at
least at present).   For example,
the only known classification of admissible unitary Banach
representations of $\GL_2(\Q_p)$ is that of~\cite{MR3150248}, which is {\em in terms of} the
$p$-adic local Langlands correspondence, rather than being independent of it.
Furthermore,
we can't expect to realize a local $p$-adic Galois representation
globally, even if we require it to be potentially semistable of some
prescribed weights and type: the family of such representations is
one-dimensional (thus uncountable), while there are only countably many modular
representations of $\Gal_{\Q}$.  This is why Taylor--Wiles patching
is so crucial to the approach of~\cite{Gpatch, MR3732208}: it provides a mechanism for relating
the local and the global contexts in the $p$-adic setting, in the
absence of the Bernstein--Zelevinsky classification/trace formula
argument that relates the local and global contexts in the classical
setting.

Perhaps the biggest difference between the $p$-adic and classical contexts
is that, in the $p$-adic context,
we don't have any {\em a priori} prescription for the correspondence
in terms analogous to the classical prescription in terms of $L$- and $\varepsilon$-
factors; the closest analogue is the requirement that there should be a
compatibility between the $p$-adic and classical correspondence mediated
by the passage to locally algebraic vectors (which should parallel the
passage from a potentially semistable Dieudonn\'e module to the underlying
Weil--Deligne representation on the Galois side).  However, it is only
in the principal series case, for the group $\GL_2(\Q_p)$, that this requirement
has the possibility of providing enough information to determine the
$p$-adic correspondence.  This is why the consideration of the (unramified) principal series
case is so crucial in all known proofs of a $p$-adic local Langlands
correspondence for~$\GL_2(\Qp)$. %

\subsubsection{Patching and the analytic setup}%
In a series of papers \cite{MR3623233}, \cite{MR3660309} and
\cite{MR4028517} it was shown that the patching construction can
be used not only to prove results about automorphic forms and Galois
representations (for example, to prove the existence of an automorphic
form whose associated Galois representations has prescribed local properties),
but also to prove results about overconvergent $p$-adic automorphic forms and their associated Galois representations. 

A key observation in the theory of overconvergent $p$-adic automorphic
forms is that the Hecke eigenvalues of overconvergent $p$-adic automorphic eigenforms nicely vary in families and form a so-called eigenvariety $Y$. The $p$-adic automorphic forms (of which the eigenvariety parameterizes the Hecke eigenvalues) themselves form a coherent sheaf $\mathcal{M}$ on this eigenvariety. 
In particular we note that in the theory of eigenvarieties, spaces of (overconvergent finite slope) $p$-adic automorphic forms can be described as the (dual of) global sections of a coherent sheaf on a rigid analytic space, that is closely related to the generic fiber of a deformation space of Galois representations. 

Characterizing the support $Y$ of this coherent sheaf $\mathcal{M}$ purely in terms of Galois representations can be regarded as an \emph{overconvergent version} of the Fontaine--Mazur conjecture: it amounts to giving a necessary and sufficient condition on a Galois representation $\rho$ that ensures that $\rho$ is associated to a  overconvergent $p$-adic automorphic form of finite slope.
Kisin \cite{MR1992017} constructed a rigid analytic subvariety (the \emph{finite slope space} of loc.\ cit.) of (the product of $\mathbf{G}_m$ with) a deformation space of $2$-dimensional ${\rm Gal}_{\Q}$-representations that is defined in purely Galois-theoretic terms (asking that the Galois representations are trianguline\footnote{The name \emph{trianguline} was only invented later by Colmez using a different but equivalent definition as the one in Kisin's paper.} at $p$ and unramified almost everywhere) and that conjecturally equals the eigencurve of Coleman--Mazur. This construction in fact inspired the definition of the trianguline variety in \cite{MR3623233} and the following overconvergent version of the Fontaine--Mazur conjecture (that we state slightly informally, and
which was first proposed by Kisin in the case $d = 2$; see \cite[11.7(2), 11.8]{MR1992017}): an irreducible  $d$-dimensional odd $p$-adic global Galois representation is associated to an overconvergent $p$-adic automorphic form of finite slope on ${\rm GL}_d$ if any only it is unramified almost everywhere and trianguline at places dividing $p$ (recall that $\rho:\Gal_{\Q}\to\GL_d(\Qpbar)$ is odd if and only if the trace of~$\rho(c)$ is~$0$ or~$\pm 1$, where~$c$ denotes a complex conjugation).
(In the case $d = 2$, this conjecture was proved under auxiliary technical hypotheses 
in~\cite[Thm.~1.2.4]{emerton2010local}; a slightly weaker version, under slightly
different hypotheses, is established in~\cite[Cor., p.~642]{KisinFM}.) %

However, these constructions and conjectures only give a
characterization of the support $Y$ of the coherent sheaf
$\mathcal{M}$. We refer to Section~\ref{sec:eigenvar} for a precise
conjecture characterizing the coherent sheaf $\mathcal{M}$ purely in
terms of local Galois representations and a global-to-local map of
stacks of Galois representations. The idea that a characterization
like this should exist originates in the insight that the patching
construction can also be carried out for eigenvarieties.

More precisely, in a similar way as $M_\infty$  over  $R_\infty$ (and $M_\infty(\sigma)$
over $R_\infty(\sigma)$) is related (via pullback along a local-global
map) to a space of automorphic forms over a global Galois deformation
ring, there is (in the setup of a definite unitary group $G$ over a
CM field $F$ satisfying some assumptions necessary for the patching
construction) a \emph{patched} version \cite{MR3623233} of the
eigenvariety $X_p(\bar\rho)$ and a coherent sheaf $\mathcal{M}_\infty$
on $X_p(\bar\rho)$ such that $\mathcal{M}$ and $Y$ may be
reconstructed from $\mathcal{M}_\infty$ and $X_p(\bar\rho)$ via
pullback along a local-global map. The main idea of the construction
is to mimic the construction of eigenvarieties of
M.E.~\cite{MR2207783} (using a locally analytic Jacquet functor), but
using a big patched module $M_\infty$ (respectively its dual
$\Pi_\infty$, which is a big Banach space representation of $G(F_p)$), instead of the space of $p$-adic automorphic forms.

The key idea of \cite{MR3660309} and \cite{MR4028517} is to use Galois-theoretic tools (i.e.~the geometry of the space $\mathfrak{X}_{G,\rm tri}$ from Section~\ref{sec:drinf-style-comp}) to analyze the local geometry of the space $X_p(\bar\rho)$, which in turn gives information about the behavior of the sheaf $\mathcal{M}_\infty$.

\section{Moduli stacks of  \texorpdfstring{$(\varphi,\Gamma)$}{(ϕ,Γ)}-modules: the Banach case}
\label{sec:moduli-stacks-phi-gamma}%
We now review the (sometimes conjectural) constructions of stacks of
$(\varphi,\Gamma)$-modules that we will work with. We work in two
different settings: in this section we cover the ``Banach'' case treated
in~\cite{emertongeepicture}, which (conjecturally) relates to smooth $p$-adic
representations of~$G(\Qp)$. In Section~\ref{sec:analytic-case} we discuss the ``analytic'' case, relating
to locally analytic representations of~$G(\Qp)$. %

In this section we very briefly recall some of
the main results of~\cite{emertongeepicture}. The reader may find it
useful to refer to~\cite{emerton2020moduli} for a more extended
introduction to and overview of~\cite{emertongeepicture}, and to
Sections~\ref{sec:banach-case-gl_2qp} and~\ref{subsec: Banach Qpf stuff} below for a more explicit
description of the stacks in some simple cases.

The stacks that we consider below are formal algebraic stacks in the sense
of~\cite{Emertonformalstacks}; colloquially, they are ``the formal
analogue of Artin stacks'', in the sense that they admit smooth covers
by formal schemes (or more precisely by formal algebraic spaces). We
note that they are in general not $p$-adically formal; equivalently,
their special (i.e.\ mod $p$) fibres are genuinely formal, rather than algebraic. We
refer the reader to~\cite[App.\ A]{emertongeepicture} for a basic
overview of formal algebraic stacks and their properties.

\subsection{An overview of~\texorpdfstring{\cite{emertongeepicture}}{[EG23]}}\label{subsec:
   overview of EG book}
We have our fixed finite extension $F/\Qp$, and we work with
coefficient rings which are $\cO$-algebras, where $\cO=\cO_L$ is
sufficiently large. We let $\cX_{d}$ denote
the moduli stack of projective \'etale $(\varphi,\Gamma)$-modules
over~$F$ of rank~$d$, in the sense of \cite[Defn.\
3.2.1]{emertongeepicture}.

\begin{rem}
  \label{rem:different phi Gamma theories}There are several different
  possible descriptions of the objects
  that are parameterized by~$\cX_d$. For example, one can consider
  them as ~$(\varphi,\Gamma)$-modules with respect to either the full
  cyclotomic extension of~$\Gal_F$, or its~$\Zp$ subextension; but it is
  also sometimes helpful to think of them as \'etale $(\Gal_F,\varphi)$-modules
  over~$W(\C^\flat)$, in the sense of~\cite[Defn.\
  2.7.3]{emertongeepicture}. %
  It can also be described in terms of (prismatic) Laurent $F$-crystals, see~\cite{min2024classicalityderivedemertongeestack}.
\end{rem}

By its very definition (and by Fontaine's theory of
$(\varphi,\Gamma)$-modules \cite{MR1106901}), the groupoid of $\Fpbar$-points of~$\cX_d$,
which coincides with the groupoid of $\Fpbar$-points
of the underlying reduced substack $\cX_{d,\red}$,
is naturally equivalent to the groupoid of continuous representations
$\rhobar:\Gal_F\to\GL_d(\Fpbar)$.  More generally, if $A$ is any finite $\cO$-algebra,
then the groupoid $\cX_d(A)$ is canonically equivalent
to the groupoid of continuous representations $\Gal_F \to \GL_d(A)$. It
follows easily that universal Galois lifting rings are versal rings
to~$\cX_d$ at finite type points (that is to say, at $\Fpbar$-points; note that
these automatically admit representatives over finite extensions of~$\Fp$),
and so we can think of~$\cX_d$ as an
algebraization of Mazur's Galois deformation rings.

\begin{rem}
  \label{rem: comparison to CWE stacks}If $A$ is not a finite
  $\cO$-algebra (e.g.\ if $A=k[X]$), then there is typically no
  Galois representation corresponding to an object
  of~$\cX_d(A)$. There is a maximal substack~$\cX_d^{\Gal}$ of~$\cX_d$ over which the
  universal $(\varphi,\Gamma)$-modules comes from a
  $\Gal_F$-representation; this stack was originally constructed and
  studied by Wang-Erickson \cite{MR3831282}. %

  The inclusion $\cX_d^{\Gal}\subset\cX_d$ induces a bijection on
  $\Fpbar$-points (and is versal at such points), and one can roughly
  imagine that $\cX_d^{\Gal}$ is the union of the formal completions
  of $\cX_d$ at its closed $\Fpbar$-points (which correspond to
  semisimple $\Gal_F$-representations). Thus although the stack
  $\cX_d^{\Gal}$ is the more obvious definition of a moduli stack of
  ``Langlands parameters'', the stack $\cX_d$ has richer geometric
  structure, and although for example the restriction maps from the
  stacks of global Galois representations to the stacks~$\cX_d$ factor
  through the stacks $\cX_d^{\Gal}$, we anticipate that it is the
  stack $\cX_d$ (or rather, its derived categories of quasicoherent sheaves)
  that is the natural setting for the $p$-adic local Langlands
  correspondence. 
\end{rem}

In the $p$-adic Langlands program it is often important to impose
conditions on the $\Gal_F$-representations under consideration, such as
demanding that they be crystalline of fixed Hodge--Tate weights, or
more generally of being potentially crystalline or semistable of a
fixed inertial type. For any\footnote{For the following construction,
we need to assume that $L$ is ``large enough'', e.g.\ so that
$\tau$ is defined over~$L$.  We largely suppress this point in what follows.}
Hodge type $\lambdau$ and inertial
type~$\tau$, there is a corresponding substack
$\cX_d^{\crys,\lambdau,\tau}$ (and similarly $\cX_d^{\semis,\lambdau,\tau}$),
which is characterized as being the unique closed substack of
$(\cX_{d})_{\cO}$ which is flat over $\cO$ and whose groupoid of
$A$-valued points, for any finite flat $\cO$-algebra~$A$, is
equivalent to the groupoid of $A$-valued $\Gal_F$-representations
satisfying the corresponding $p$-adic Hodge theoretic property (of
being potentially crystalline or semistable of Hodge type~$\lambdau$
and inertial type~$\tau$). In the
below if~$\tau$ is trivial then we drop it from the notation.

These substacks are constructed in~\cite{emertongeepicture} using (among
other techniques) Kisin's results on Breuil--Kisin modules and
crystalline representations~\cite{KisinCrys}; they are probably
more naturally defined in terms of (log-) prismatic $F$-crystals,
following~\cite{bhatt2021prismatic}.

The following theorem summarises the main results
of~\cite{emertongeepicture}.
\begin{thm}
  \label{thm: main results of EG moduli}\leavevmode
  \begin{enumerate}
  \item $\cX_d$ is a 
	Noetherian formal algebraic
	stack.
      \item Each substack $\cX_d^{\crys,\lambdau,\tau}$ and
        $\cX_d^{\semis,\lambdau,\tau}$ is a $p$-adic formal algebraic
        stack.
              \item   The underlying reduced substack~$\cX_{d,\red}$ of~$\cX_d$
      	{\em (}which is an algebraic stack{\em )}
	is of finite type
	over $\Fp$, and is equidimensional of dimension
	$[F:\Q_p] d(d-1)/2$. 
      \item The special fibre of each of $\cX_d^{\crys,\lambdau,\tau}$ and
        $\cX_d^{\semis,\lambdau,\tau}$ 	{\em (}which is an algebraic
        stack{\em )} is equidimensional of dimension equal to that of
        the flag variety determined by~$\lambdau$. In particular,
        if~$\lambdau$ is regular, then each of these stacks is
        equidimensional of dimension  	$[F:\Q_p] d(d-1)/2$. 

      \item The irreducible components of 
	$\cX_{d,\red}$
	admit a natural labelling by Serre
        weights.
           
  \end{enumerate}
\end{thm}
\begin{rem}%
  \label{rem: what does labelling of components by Serre weights
    mean?}The notion of a Serre weight is recalled in Section~\ref{subsubsec:
  geometric BM Banach}. The precise meaning of the labelling of the irreducible
  components of $\cX_{d,\red}$ by Serre weights is somewhat involved;
  see \cite[\S 5.5]{emertongeepicture}. Roughly, the description is as
  follows. Each irreducible component of $\cX_{d,\red}$ admits a dense
  open substack, whose $\Fpbar$-points correspond to
  $\Gal_F$-representations which are maximally nonsplit of niveau one;
  that is, they admit a unique filtration by characters. The ordered
  tuple of inertial
  weights of these characters determines the irreducible component,
  and these ordered tuples correspond to Serre weights, because
  Serre weights are by definition (isomorphism classes of) irreducible
  $\Fpbar$-representations of $\GL_d(k)$, and as such are indexed by
  their highest weights. (There is one subtlety in making this
  identification, which is that it is necessary to distinguish Serre
  weights whose highest weights are congruent modulo $p-1$; this
  manifests itself on the Galois side as a distinction between peu and
  tr\`es ramifi\'ee representations.)
\end{rem}

\begin{rem}
  \label{rem: coincidence of dimensions suggests BM}The two
  appearances of $[F:\Qp]d(d-1)/2$ in Theorem~\ref{thm: main results
    of EG moduli} --- as the dimension of the underlying reduced
  substack ~$\cX_{d,\red}$, and of the dimension of the special fibres
  of potentially crystalline/semistable substacks of regular Hodge
  type --- is not a coincidence. As well as being an important
  ingredient in the proof of Theorem~\ref{thm: main results of EG
    moduli}, it underlies the Breuil--M\'ezard conjecture (see Section~\ref{subsubsec:
  geometric BM Banach}). %
\end{rem}

\begin{rem}
  \label{rem: don't make dimension of whole EG stack precise}As far as
  we are aware, there is no general dimension theory for formal
  algebraic stacks in the literature. It is, however, reasonable to
  think that such a theory exists, and that this dimension can be
  computed from the dimensions of the versal rings. This being the
  case, the results of~\cite{BIP} (which in particular compute the
  dimensions of unrestricted local Galois lifting rings) imply that
  ~$\cX_d$ will be equidimensional of dimension~$1+[F:\Qp]d^2$. Since
  its underlying reduced substack $\cX_{d,\red}$ is equidimensional of
  dimension $[F:\Qp]d(d-1)/2$, we see that $\cX_d$ should have
  $1+[F:\Qp]d(d+1)/2$ directions which are purely formal (including
  the $p$-adic direction). In
  particular, it is (provably) never a $p$-adic formal algebraic stack. %
\end{rem}

\begin{rem}
  \label{rem: expect that X is lci}Again, as far as
  we are aware, there is no notion in the literature of what it means
  for a formal
  algebraic stack to be lci. We anticipate that reasonable definitions
  exist, and that $\cX_d$ is lci; again, this should follow from the
  results of~\cite{BIP} (which show that the unrestricted local Galois
  lifting rings are lci).
(A version of this expectation is realized by Min in~\cite{min2024classicalityderivedemertongeestack}, which shows that a natural derived variant of~$\cX_d$ is in fact classical.)
\end{rem}

\subsection{Fixed determinant variant}
It is often convenient to consider a variant of these stacks where we
fix the  determinant. If~$\chi:\Gal_F\to\cO^\times$ is a character, then
we let  $\cX_d^\chi(A)$
be the groupoid of pairs $(D,\theta)$ where~$D$ is a rank~$d$
projective \'etale $(\varphi,\Gamma)$-module with $A$-coefficients,
and $\theta$ is an identification of~$\wedge^dD$ with~$\chi$. These
stacks (and their potentially crystalline and semistable variants) are
studied in an appendix to~\cite{DEGcategoricalLanglands}, where the analogue of
Theorem~\ref{thm: main results of EG moduli} is established (the only
difference being that in the case of ``Steinberg'' Serre weights,
there can be multiple corresponding irreducible components, indexed by
the $d$th roots of unity in~$\Fpbar$). In particular, the underlying
reduced substack of $\cX_d^\chi$ is again equidimensional of dimension $[F:\Qp]d(d-1)/2$.

\section{Moduli stacks of \texorpdfstring{$(\varphi,\Gamma)$}{(ϕ,Γ)}-modules: the analytic case}
\label{sec:analytic-case}

Similarly to Section~\ref{sec:moduli-stacks-phi-gamma}, we want to define and study stacks of $(\varphi,\Gamma)$-modules on the category of rigid analytic spaces. 
It is sometimes more convenient to think of these objects in terms of equivariant vector bundles on the Fargues--Fontaine curve. Most of the results in this section will appear in the forthcoming paper \cite{HellmannHernandezSchraen}.
We start by recalling the setup and the comparison with $(\varphi,\Gamma)$-modules.

\subsection{ \texorpdfstring{$(\varphi,\Gamma)$}{(ϕ,Γ)}-modules over the Robba ring}

\label{sec:varphi-gamma-modules-Robba}
Recall that we have fixed a finite extension $F$ of $\mathbf{Q}_p$ and an algebraic closure $\bar F$ of $F$. We moreover fix a compatible system $\epsilon_{p^n}$ of $p^n$-th roots of unity and write $F_n=F(\epsilon_{p^n})$ and $F_\infty=\bigcup_n F_n\subset \bar F$. We moreover write $\Gamma={\rm Gal}(F_\infty/F)\supset \Gamma_n={\rm Gal}(F_\infty/F_n)$.
We moreover note that the $p$-adic completions of $\bar F$ and of $F_\infty$ are both perfectoid fields. 

Associated with the field $F$ (or rather with the cyclotomic extension
$F_\infty$ and the algebraic closure $\bar F$) we can define two
versions of the Fargues--Fontaine curve. Let $\varpi^\flat\in\mathcal{O}_{\hat F_\infty^\flat}$ denote the choice of a pseudo-uniformizer. We set 
\begin{align*}
Y_{F_\infty}&={\rm Spa}(W(\mathcal{O}_{\hat F_{\infty}^\flat}),W(\mathcal{O}_{\hat F_{\infty}^\flat}))\backslash V(p[\varpi^\flat]),\\
Y_{\bar F}&={\rm Spa}(W(\mathcal{O}_{\hat{ \bar F}^\flat}),W(\mathcal{O}_{\hat {\bar F}^\flat}))\backslash V(p[\varpi^\flat]).
\end{align*}
These adic spaces over $\mathbf{Q}_p$ come equipped with an automorphism $\varphi$ and a continuous action of $\Gamma$, respectively $\Gal_F$, that commutes with~$\varphi$. 
 The Fargues--Fontaine curves $X_{F_\infty}$ and $X_{\bar F}$ are the quotients 
 \begin{align*}
 X_{F_\infty}&=Y_{F_\infty}/\varphi^\mathbf{Z},\\
  X_{\bar F}&=Y_{\bar F}/\varphi^\mathbf{Z},
 \end{align*}
  which inherit a continuous action of $\Gamma$, respectively $\Gal_F$.   
 
 Note that the curves $X_{F_\infty}$ and $X_{\bar F}$ have a distinguished point $\infty$ defined by the kernel of the usual $\theta$ maps of Fontaine, see e.g.~\cite[2.1.2]{MR3917141}:
 \begin{align*}
 \theta_{F_\infty}: W(\mathcal{O}_{\hat F_\infty^\flat})&\longrightarrow \mathcal{O}_{\hat F_\infty},\\
  \theta_{\bar F}: W(\mathcal{O}_{\hat {\bar F}^\flat})&\longrightarrow \mathcal{O}_{\hat {\bar F}}.
  \end{align*}
It is a direct consequence of the construction that the point $\infty$ is a fixed point for the action of $\Gamma$, respectively of $\Gal_F$-action. 

We define the following categories of (semi-)linear algebra objects: 
\begin{enumerate}
\item[-] the category ${\rm VB}_{X_{F_\infty}}^\Gamma$ of $\Gamma$-equivariant vector bundles on $X_{F_\infty}$,
\item[-] the category ${\rm VB}_{X_{\bar F}}^{\Gal_F}$ of $\Gal_F$-equivariant vector bundles on $X_{\bar F}$,
\item[-] the category ${\rm VB}^{\varphi,\Gamma}_{Y_{F_\infty}}$ of $(\varphi,\Gamma)$-equivariant vector bundles on $Y_{F_\infty}$,
\item[-] the category ${\rm VB}^{\varphi,\Gal_F}_{Y_{\bar F}}$ of $(\varphi,\Gal_F)$-equivariant vector bundles on $Y_{\bar F}$.
\end{enumerate}
In fact in all the above categories we ask that the action of $\Gamma$ resp.~$\Gal_F$ on the vector bundles is continuous in the evident sense (that we will not spell out explicitly).
By the work of Berger, Fargues--Fontaine and Kedlaya, all these various categories are known to be canonically equivalent. 
\begin{thm}\label{EHequivVB}
\hfill\\(i) The  projection $X_{\bar F}\rightarrow X_{F_\infty}$ is
equivariant with respect to the canonical projection
$\Gal_F\rightarrow \Gamma$ and induces via pullback an equivalence of categories 
$${\rm VB}_{X_{F_\infty}}^\Gamma\rightarrow {\rm VB}_{X_{\bar F}}^{\Gal_F}.$$
The corresponding statement for $Y_{\bar F}\rightarrow Y_{F_\infty}$ also holds true.\\
\noindent (ii) The projection $Y_{F_\infty}\rightarrow X_{F_\infty}$ is $\Gamma$-equivariant and induces, via pullback, an equivalence of categories
$${\rm VB}_{X_{F_\infty}}^{\Gamma}\rightarrow {\rm VB}_{Y_{F_\infty}}^{\varphi,\Gamma}.$$
The corresponding statement for $Y_{\bar F}\rightarrow X_{\bar F}$ also holds true.
\end{thm} 
The second part of the theorem is rather formal. The first part follows from \cite[Thm. 9.3.1]{MR3917141}
\begin{rem}
 In Scholze's world of diamonds one has $X_{F_\infty}/\Gamma=X_{\bar F}/\Gal_F$, giving a more geometric interpretation of the first equivalence. 
 \end{rem}

Furthermore there are obvious compatibilities between these equivalences. 
There is a further description of these categories in terms of the Beauville--Laszlo gluing
lemma~\cite{MR1320381}. This gives rise to Berger's category of $B$-pairs \cite{MR2377364}.
We set 
\begin{align*}
B_e&=\Gamma(X_{\bar F}\backslash\{\infty\},\mathcal{O}_{X_{\bar F}}),&  B_{{\rm dR}}^+&=\hat{\mathcal{O}}_{X_{\bar F},\infty}, & B_{{\rm dR}}&=B_{\rm dR}^+[1/t].
\end{align*}
Each of these rings is equipped with a continuous $\Gal_F$-action.   Using the Beauville--Laszlo gluing lemma, we can identify ${\rm VB}_{X_{\bar F}}^{\Gal_F}$ with the category of triples $(M,\Lambda, \xi)$ consisting of a finite projective $B_{e}$-module $M$ and a finite projective $B^+_{{\rm dR}}$-module $\Lambda $ equipped with continuous semi-linear $\Gal_F$-actions and a $\Gal_F$-equivariant isomorphism $$\xi:M\otimes_{B_{e}}B_{{\rm dR}}\isoto\Lambda\otimes_{B_{{\rm dR}}^+}B_{{\rm dR}}.$$

There is yet another description of ${\rm VB}_{X_{\bar F}}^{\Gal_F}\cong{\rm VB}_{X_{F_\infty}}^\Gamma$ using $(\varphi,\Gamma)$-modules over imperfect period rings. 
We write $\mathcal{R}_F$ for the Robba ring
$$\mathcal{R}_F=\lim_{\substack{\longrightarrow\\ r\rightarrow 1}}\lim_{\substack{\longleftarrow\\ s\rightarrow 1}}B^{[r,s]}_F,$$
where $B^{[r,s]}_F$ is the ring of rigid analytic functions on the closed annulus $\mathbf{B}^{[r,s]}_{F_{\infty,0}}$ of inner radius $r$ and outer radius $s$ in the closed unit disc over $F_{\infty,0}$, the maximal unramified subextension of $F_\infty$ (which in general can be strictly larger than $F_0$). 
For $s>r\gg 0$ there is a canonical way to define a continuous action
of $\Gamma$ on $B^{[r,s]}_F$ and to define a ring
homomorphism $$\varphi:B^{[r,s]}_F\rightarrow
B_F^{[r^{1/p},s^{1/p}]}$$ commuting with the $\Gamma$-action. These
data induce a ring homomorphism $\varphi:\mathcal{R}_F\rightarrow
\mathcal{R}_F$ and a continuous $\Gamma$-action on $\mathcal{R}_F$
commuting with $\varphi$ (we refer e.g.~to  \cite[2.2]{MR3230818} for
the details of the definition of the $\varphi$- and  $\Gamma$-actions, and for the variants of the ring $\mathcal{R}_F$ with coefficients in an affinoid algebra below).

\begin{defn}
The category ${\rm Mod}^{\varphi,\Gamma}_{\mathcal{R}_F}$ is the category of finite projective $\mathcal{R}_F$-modules together with continuous $\Gamma$-action and an isomorphism $\varphi^\ast D\rightarrow D$ commuting with the $\Gamma$-action.
\end{defn}
\begin{rem}\label{EH rem peropdrings}
We recall the relation of $(\varphi,\Gamma)$-modules over
$\mathcal{R}_F$ with the $(\varphi,\Gamma)$-modules used in Section \ref{sec:moduli-stacks-phi-gamma}. The $(\varphi,\Gamma)$-modules of Section \ref{sec:moduli-stacks-phi-gamma} are defined over a $p$-adically complete ring $\mathbf{A}_F$.
This ring  contains a subring $\mathbf{A}_F^\dagger\subset \mathbf{A}_F$ of \emph{overconvergent elements} that is stable under the action of $\varphi$ and $\Gamma$. Moreover $\mathbf{A}_F^\dagger$ embeds (equivariantly for $\varphi$ and $\Gamma$) into the $\mathbf{Q}_p$-algebra $\mathcal{R}_{F}$.
By a theorem of Cherbonnier--Colmez \cite{MR1645070} every \'etale
$(\varphi,\Gamma)$-module over $\mathbf{B}_F=\mathbf{A}_F[1/p]$ is overconvergent, where by definition the \'etale $(\varphi,\Gamma)$-modules over $\mathbf{B}_F$ are those that come by base change from a $(\varphi,\Gamma)$-module over $\mathbf{A}_F$. More
precisely, the scalar extension from $\mathbf{B}_F^\dagger=\mathbf{A}_F^\dagger[1/p]$ to
$\mathbf{B}_F$ induces an equivalence between the corresponding
categories of \'etale $(\varphi,\Gamma)$-modules.  Moreover, the extension of
scalars from $\mathbf{B}_F^\dagger$ to $\mathcal{R}_F$ induces a fully
faithful embedding of the category of \'etale 
$(\varphi,\Gamma)$-modules over $\mathbf{B}^\dagger_F$ to the category of $(\varphi,\Gamma)$-modules over $\mathcal{R}_{F}$. By definition the
$(\varphi,\Gamma)$-modules over $\mathcal{R}_{F}$ in the essential
image of this functor are called \'etale.

By Kedlaya's slope
filtration theorem \cite{MR2493220} the \'etale $(\varphi,\Gamma)$-modules
are precisely the objects that are semi-stable of slope zero, in the
following sense:
If $D$ is
a $(\varphi,\Gamma)$-module $D$ of rank~$r$,
then by the classification 
of the $(\varphi,\Gamma)$-modules of rank~$1$ (see Section~\ref{sec: GL1} below),
we may write $\bigwedge^r D  = \mathcal{R}_F(\delta)$,
where $\delta$ is a continuous character of $F^\times$.
We then define the {\em slope} of $D$  to be 
$${\rm slope}(D)=\dfrac{{\rm val}_p\bigl(\delta(\varpi_F)\bigr)}{r},$$
where $\varpi_F$ is a choice of a uniformizer in $F$. 
We furthermore say that $D$ is {\em  semi-stable}
if ${\rm slope}(D')\geq {\rm slope}(D)$ for all $(\varphi,\Gamma)$-stable subobjects $D'\subseteq D$. 
\end{rem}

\begin{thm}\label{EHequivVBphiGamma} There is an equivalence of categories
$${\rm Mod}^{\varphi,\Gamma}_{\mathcal{R}_F}\longrightarrow {\rm VB}^{\varphi,\Gamma}_{Y_{F_\infty}}.$$
\end{thm}
\begin{proof}[Proof sketch.]
  The rough idea behind this theorem is the following: any
  $(\varphi,\Gamma)$-module over $\mathcal{R}_F$ admits a model over
  some closed annulus $\mathbf{B}^{[r,s]}_{F_{\infty,0}}$. Using the
  choice of the coordinate function $[\varpi^\flat]$ we can map a
  corresponding closed annulus in $Y_{F_\infty}$ to
  $\mathbf{B}^{[r,s]}_{F_{\infty,0}}$, equivariantly for the action of
  $\Gamma$ (and compatibly with $\varphi$ in an evident way). The
  pullbacks of the $\varphi$-bundles with $\Gamma$-action on the
  various $\mathbf{B}^{[r,s]}_{F_{\infty,0}}$ to the closed annuli in
  $Y_{F_\infty}$ spread out to a $\Gamma$-equivariant
  $\varphi$-bundle on $Y_{F_\infty}$ (as $\varphi$ is an automorphism
  of $Y_{F_\infty}$ that shifts the radii of closed annuli). The hard
  part of the theorem is the essential surjectivity, sometimes referred to as {\em deperfection}. It is
  proved using the additional structure on the cyclotomic tower given
  by Tate's normalized traces. A proof of the theorem can be given by
  combining \cite[10.1.2]{MR3917141} and \cite[Theorem
  2.2.7]{MR2377364}.  We refer to \cite[4.2]{MR2493221} for a variant
  of the descent argument using Tate's normalized traces.
\end{proof}

We will be interested in families of these (equivalent) objects with coefficients in rigid analytic spaces. More precisely we want to define a stack $\mathfrak{X}_d$ whose ${\rm Sp}(\mathbf{\Q}_p)$-valued points are given by the groupoid of rank~$d$ objects in ${\rm VB}_{X_{\bar F}}^{\Gal_F}$ or any of the equivalent categories described above. 
We briefly describe the variants of the above categories over an affinoid rigid space ${\rm Sp}(A)$, which suffices to define our stacks. 

Let $A$ be an affinoid algebra topologically of finite type over $\mathbf{Q}_p$. By a result of Kedlaya \cite{MR3519123} the Fargues--Fontaine curves $X_{F_\infty}, X_{\bar F}$ are strongly Noetherian, and the same applies (locally) to their coverings $Y_{F_\infty}, Y_{\bar F}$.  
Hence the fiber products $$X_{F_\infty,A}:=X_{F_\infty}\times_{{\rm Sp}(\mathbf{Q}_p)}{\rm Sp}(A)$$ etc.\ are well-defined in the category of adic spaces, and come equipped with a canonical continuous action of $\Gamma$, respectively $\Gal_F$.  

Similarly to the absolute case considered above we can define 
\begin{enumerate}
\item[-] the category ${\rm VB}_{X_{F_\infty,A}}^\Gamma$ of $\Gamma$-equivariant vector bundles on $X_{F_\infty,A}$,
\item[-] the category ${\rm VB}_{X_{\bar F},A}^{\Gal_F}$ of $\Gal_F$-equivariant vector bundles on $X_{\bar F,A}$,
\item[-] the category ${\rm VB}^{\varphi,\Gamma}_{Y_{F_\infty},A}$ of $(\varphi,\Gamma)$-equivariant vector bundles on $Y_{F_\infty,A}$,
\item[-] the category ${\rm VB}^{\varphi,\Gal_F}_{Y_{\bar F},A}$ of $(\varphi,\Gal_F)$-equivariant vector bundles on $Y_{\bar F,A}$.
\end{enumerate}
Again we note that the $\Gamma$- respectively $\Gal_F$-action on the vector bundles is supposed to be continuous in the evident way.
Moreover, there is a variant with coefficients of the category of $(\varphi,\Gamma)$-modules over the imperfect Robba ring $\mathcal{R}_{F}$. We define 
$$\mathcal{R}_{F,A}=\lim_{\substack{\longrightarrow\\ r\rightarrow 1}}\lim_{\substack{\longleftarrow\\ s\rightarrow 1}}B_F^{[r,s]}\widehat{\otimes}_{\mathbf{Q}_p}A$$ which is again equipped with a continuous $\Gamma$-action and an endomorphism $\varphi$ commuting with the $\Gamma$-action. Then the category ${\rm Mod}^{\varphi,\Gamma}_{\mathcal{R}_{F,A}}$ is defined to be the category of finite projective $\mathcal{R}_{F,A}$-modules $D$ together with an isomorphism $\varphi_D:\varphi^\ast D\rightarrow D$ and a continuous $\Gamma$-action commuting with $\varphi_D$.

The conclusions of Theorems \ref{EHequivVB} and \ref{EHequivVBphiGamma} still hold true for the categories with coefficients in $A$, i.e.\ there are equivalences of categories
\[
\begin{xy}
\xymatrix{
{\rm VB}_{X_{F_\infty,A}}^\Gamma \ar^-{\sim}[d]\ar^-{\sim}[r]& {\rm VB}_{Y_{F_\infty,A}}^{\varphi,\Gamma}\ar^-{\sim}[d]\ar^-{\sim}[r] &  {\rm Mod}^{\varphi,\Gamma}_{{\mathcal{R}}_{F,A}}\ar[l]\\
{\rm VB}_{X_{\bar F,A}}^{\Gal_F} \ar^-{\sim}[r]& {\rm VB}_{Y_{\bar F,A}}^{\varphi,\Gal_F}
}
\end{xy}
\]
induced by similar functors as in the case without coefficients. Moreover, given a morphism $A\rightarrow B$ of affinoid algebras, there are obvious base change functors
$$(-)\widehat{\otimes}_AB:{\rm VB}_{X_{F_\infty,A}}^\Gamma\longrightarrow {\rm VB}_{X_{F_\infty,B}}^\Gamma$$
etc.\ and the above equivalences of categories are obviously compatible with the base change from $A$ to $B$. 
\begin{rem}
It is also worth noticing that the description of these equivalent categories using the Beauville--Laszlo gluing lemma generalizes to families. This allows us to modify a given family of equivariant vector bundles over $X_{\bar F,A}$ by a given family of invariant lattices over the completion of $X_{\bar F,A}$ along $\{\infty\}\times{\rm Sp}\,A$.
\end{rem}
\bigskip

\subsubsection{Rigid analytic Artin stacks}
\label{subsec:rigid stacks}
We introduce basic notions in the framework of Artin stacks on a category of rigid analytic spaces to set the ground for the definition of the stacks of rigid analytic $(\varphi,\Gamma)$-modules. Denote by ${\rm Rig}_L$ the category of rigid analytic spaces over a fixed base field $L$ that is a finite extension of~$\mathbf{Q}_p$. 
We will equip ${\rm Rig}_L$ with the Tate-fpqc topology \cite[2.1]{MR2524597}. The coverings in this topology are generated by the usual (admissible) Tate coverings and the morphisms ${\rm Sp}(A)\rightarrow {\rm Sp}(B)$ of rigid spaces for faithfully flat maps $B\rightarrow A$ of affinoid algebras. 
With respect to this topology all representable functors are sheaves, and coherent sheaves satisfy descent \cite[Theorem 4.2.8]{MR2266885}. In the following a stack on ${\rm Rig}_L$ will be a category fibered in groupoids over ${\rm Rig}_L$ that satisfies descent for the Tate-fpqc topology. 
As in the case of Artin stacks on schemes, we have to start with the definition of an analogue of algebraic spaces.

\begin{defn}
A quasi-analytic space is a sheaf $\mathcal{F}$ on ${\rm Rig}_L$ such that the diagonal $\mathcal{F}\rightarrow \mathcal{F}\times\mathcal{F}$ is representable and such that there exists an \'etale surjection $U\rightarrow \mathcal{F}$ from a representable sheaf $U$ onto $\mathcal{F}$.
\end{defn}
\begin{rem}
\hfill\\(i) In fact below we do not need this level of generality: by a result of Conrad and Temkin \cite[Theorem 1.2.2]{MR2524597} every separated quasi-analytic space is representable by a rigid analytic space. Below we will only meet separated quasi-analytic spaces. \\
(ii) The PhD thesis of Evan Warner~\cite{MR4257201} also develops a
theory of Artin stacks on adic spaces. The main difference to our set
up is that Warner works with all strongly Noetherian adic spaces (and
uses the \'etale topology), whereas we allow only rigid analytic
spaces (and use the Tate-fpqc topology); moreover Warner uses a
different terminology. We mainly restrict ourselves to rigid analytic
spaces as, at least for the time being, the theory of
$(\varphi,\Gamma)$-modules with coefficients \cite{MR3230818} and in
particular the finiteness results about their cohomology, is limited
to coefficients in affinoid algebras of classical rigid analytic
geometry. %
\end{rem}

\begin{defn}
A rigid analytic Artin stack is a stack $\mathfrak{X}$ on ${\rm Rig}_L$ such that the diagonal $\Delta_\mathfrak{X}:\mathfrak{X}\rightarrow \mathfrak{X}\times_L\mathfrak{X}$ is representable by quasi-analytic spaces and such that there exists a rigid analytic space $U$ and a smooth surjection $U\rightarrow \mathfrak{X}$. 
\end{defn}

Given a rigid analytic Artin stack $\mathfrak{X}$ we can define as usual its category of coherent sheaves as
$${\rm Coh}(\mathfrak{X})=\lim_{\substack {\longleftarrow\\ {\rm Sp}(A)\rightarrow \mathfrak{X}}} {\rm Coh}({\rm Sp}(A)),$$
where the limit is taken over all maps of affinoid rigid spaces ${\rm Sp}(A)$ to $\mathfrak{X}$. As usual, given a smooth surjection $U\rightarrow \mathfrak{X}$ from an (affinoid) rigid analytic space $U$ to $\mathfrak{X}$ this can be computed as the limit of the diagram
$${\rm Coh}(U)\rightrightarrows{\rm Coh}(U\times_{\mathfrak{X}}U)\substack{\rightarrow\\[-1em] \rightarrow \\[-1em] \rightarrow}\dots.$$
Similarly the stable %
$\infty$-category of coherent sheaves on $\mathfrak{X}$ can be defined as the homotopy limit of $\infty$-categories
$$\mathbf{D}_{\rm coh}(\mathfrak{X})=\lim_{\substack {\longleftarrow\\ {\rm Sp}(A)\rightarrow \mathfrak{X}}} \mathbf{D}_{\rm coh}({\rm Sp}(A)),$$
and for a given covering it can be computed as the homotopy limit of the diagram of $\infty$-categories
$$\mathbf{D}_{\rm coh}(U)\rightrightarrows \mathbf{D}_{\rm coh}(U\times_{\mathfrak{X}}U)\substack{\rightarrow\\[-1em] \rightarrow \\[-1em] \rightarrow}\dots.$$

For a morphism of rigid analytic spaces $f:\mathfrak{X}\rightarrow \mathfrak{Y}$, or more generally for a morphism between rigid analytic Artin stacks, the definition implies that we have a canonical derived pullback $Lf^\ast: \mathbf{D}_{\rm coh}(\mathfrak{Y})\rightarrow \mathbf{D}_{\rm coh}(\mathfrak{X})$.
If the map $f$ is proper we also have (derived) pushforward $Rf_\ast: \mathbf{D}_{\rm coh}(\mathfrak{X})\rightarrow \mathbf{D}_{\rm coh}(\mathfrak{Y})$ which is right adjoint to $Lf^\ast$.
More precisely, for general $f$ the pushforward $Rf_\ast$ can be defined on the subcategory of $\mathbf{D}_{\rm coh}(\mathfrak{X})$ consisting of sheaves with proper support over $\mathfrak{Y}$.
For the rigid analytic Artin stacks that we consider there always exists a dualizing complex $\omega_\mathfrak{X}$, and the
internal
duality $$\mathbf{D}_\mathfrak{X}(-)=R\mathcal{H}om_{\cO_\mathfrak{X}}(-,\omega_{\mathfrak{X}})$$
allows us to use 
\[f^!=\mathbf{D}_\mathfrak{X}\circ Lf^\ast\circ\mathbf{D}_\mathfrak{Y}:\mathbf{D}_{\rm coh}(\mathfrak{Y})\rightarrow \mathbf{D}_{\rm coh}(\mathfrak{X})\]
 as an ad hoc definition of the upper shriek pullback $f^!$. 
 
\begin{rem}\label{EHrem:sheavesnotsolid}
The definitions above are a bit restrictive, as we only defined categories of coherent and not of quasicoherent sheaves. This is mainly due to the fact that quasicoherent sheaves on rigid analytic spaces are (for topological reasons) a bit delicate. In the world of condensed rings and modules of Clausen--Scholze it is possible to develop a theory of a (derived) category of quasicoherent sheaves together with a full six-functor formalism that gives a more uniform approach than the above ad hoc definition of $f^!$. 
\end{rem}

We mention two sources of rigid analytic Artin stacks: 
Given an Artin stack $\mathbf{X}$ on the category ${\rm Sch}_L$ of
$L$-schemes such that $\mathbf{X}$ is (locally) of finite type over~
$L$ (i.e.\ $\mathbf{X}$ can be covered by a scheme $\mathbf{U}$
(locally) of finite type over $L$), its analytification can be defined as the sheafification of 
$$\mathbf{X}^{\rm an}({\rm Sp}(A))=\mathbf{X}({\rm Spec}(A)).$$
This is easily  seen to be the correct definition in the case $\mathbf{X}=\mathbf{A}^n$ and hence in the case of affine schemes (of finite type). 
In particular: for a presentation $\mathbf{V}=\mathbf{U}\times_\mathbf{X} \mathbf{U}\rightrightarrows \mathbf{U}\rightarrow \mathbf{X}$ the analytification of $\mathbf{X}$ can be defined as the stack quotient of $\mathbf{V}^{\rm an}\rightrightarrows \mathbf{U}^{\rm an}$.
Given $\mathbf{X}$ we have an obvious analytification functor on the (derived) category of coherent sheaves which commutes with $Lf^\ast, f^!, Rf_\ast$ whenever they are defined.

On the other hand, if $\mathcal{X}$ is a formal algebraic stack over ${\rm Spf}(\mathcal{O}_L)$ in the sense of \cite[5.]{Emertonformalstacks}, then we can define its generic fiber as
\numequation\label{EHdef:genfibofstack}\mathcal{X}_\eta^{\rm
    rig}({\rm Sp}(A))=\lim_{\substack{\longrightarrow\\ Y}}\mathcal{X}(Y),
\end{equation}
where the limit is taken over all formal models (in the sense of Raynaud) $Y$ of ${\rm Sp}(A)$.
\begin{prop}\label{EH prop genfibofformalstacks}
Let $\mathcal{X}$ be a formal algebraic stack topologically of finite type, with diagonal topologically of finite type. Then $\mathcal{X}_\eta^{\rm rig}$ is a rigid analytic Artin stack.
\end{prop}

Here we say that a formal algebraic stack $\mathcal{X}$ is topologically of finite type if it admits a smooth surjection from a formal scheme that is locally of the form ${\rm Spf}\,A$ for a complete $\mathbf{Z}_p$-algebra $A$ topologically of finite type (with a similar definition of the diagonal being topologically of finite type).
The idea to prove Proposition \ref{EH prop genfibofformalstacks} is to reduce to the affine case where $\mathcal{X}$ can be written as the coequalizer of formally smooth maps
$$
{\rm Spf}\, B \rightrightarrows {\rm Spf}\,A
$$
for two complete $\mathbf{Z}_p$-algebras $A,B$ topologically of finite type. The rigid analytic generic fiber of $\mathcal{X}$ can then be identified with the coequalizer of 
$$
({\rm Spf}\, B)_\eta^{\rm rig} \rightrightarrows ({\rm Spf}\,A)_\eta^{\rm rig},
$$
where $(-)_\eta^{\rm rig}$ denotes the rigid analytic generic fiber of a formal scheme in the sense of Berthelot.

\bigskip

\subsubsection{Rigid analytic stacks of $(\varphi,\Gamma)$-modules} We now return to the main objective of this subsection and define stacks of $(\varphi,\Gamma)$-modules.

\begin{defn}
For $d\geq 1$ the category fibered in groupoids $\mathfrak{X}_{{\rm GL}_d}=\mathfrak{X}_d$ over ${\rm Rig}_L$ is the groupoid that assigns to an affinoid space ${\rm Sp}\,A$ the groupoid of $(\varphi,\Gamma)$-modules of rank~$d$ over $\mathcal{R}_{F,A}$.
\end{defn}
Using the fact that a $(\varphi,\Gamma)$-module over $\mathcal{R}_{F,A}$ admits a model over an annulus over ${\rm Sp}\,A$, descent for vector bundles on rigid analytic spaces implies:
\begin{prop}
The category fibered in groupoids $\mathfrak{X}_d$ is a stack.
\end{prop}
\begin{rem}
  \label{rem: comparison to Banach stack}For the expected relationship
  between~ $\mathfrak{X}_d$ and the stack~$\cX_d$ of Section
  ~\ref{sec:moduli-stacks-phi-gamma}, see
  \eqref{EHgenfibEGstacktoanalytic} and the surrounding discussion.
\end{rem}
\begin{remark}\label{EHrem:propertiesphiGammacohom}
By the above discussion the stack $\mathfrak{X}_d$ can be described equivalently as the stack of $\Gal_F$-equivariant vector bundles on the Fargues--Fontaine curve and in the following we will freely choose whichever description fits better with our purposes. Hence there are some obvious similarities with the classical theory of (stacks of) vector bundles on an algebraic curve. We list the most important similarities and differences:
\begin{enumerate}
\item[(i)] There is a cohomology theory $$R\Gamma_{\Gal_F}(X,-):=R\Gamma(\Gal_F,R\Gamma({\rm Sp}\,A\times X_{\bar F},-))$$ for $\Gal_F$-equivariant vector bundles on ${\rm Sp}\,A\times X_{\bar F}$. In terms of $(\varphi,\Gamma)$-modules this cohomology theory coincides with the usual $(\varphi,\Gamma)$ cohomology that is computed by the Herr complex \cite[Def. 2.2.3]{MR3230818}. 
More precisely, given a $(\varphi,\Gamma)$-module $D$ over
$\mathcal{R}_{F,A}$, the corresponding Herr complex is
\[\begin{xy}
\xymatrix{
 R\Gamma_{\varphi,\Gamma}(D)=[D^\Delta \ar[rr]^{(\varphi-1,\gamma-1)} && D^\Delta\oplus D^\Delta \ar[rr]^{(\gamma-1)\oplus(1-\varphi)} && D^\Delta],
}
\end{xy}\]
where $\Delta\subset\Gamma$ is the torsion subgroup and where $\gamma\in\Gamma/\Delta$ is a topological generator.
This complex is (quasi-isomorphic to) a perfect complex of $A$-modules concentrated in degrees $0,1,2$, see \cite[Theorem 4.4.5]{MR3230818}.
\item[(ii)] As a consequence of (i) an analogue of Grothendieck's cohomology and base change theorem for coherent sheaves holds true in the context of equivariant vector bundles on the Fargues--Fontaine curve, respectively in the context of $(\varphi,\Gamma)$-modules.  More precisely, for a given $(\varphi,\Gamma)$-module $D$ over $\mathcal{R}_{F,A}$ and a point $x\in {\rm Sp}\, A$ the base change map
\[{\rm bc}^i_x:H^i_{\varphi,\Gamma}(D)\otimes_Ak(x)\longrightarrow H^i_{\varphi,\Gamma}(D\otimes_{\mathcal{R}_{F,A}}\mathcal{R}_{F,k(x)})\]
is surjective if and only if it is an isomorphism; and if ${\rm bc}_x^i$ is surjective, then $H^i_{\varphi,\Gamma}(D)$ is a vector bundle in a neighborhood of $x$ if and only if ${\rm bc}_x^{i-1}$ is surjective. 
\item[(iii)] There is no analogue of an ample line bundle: for a given equivariant vector bundle $\mathcal{V}$ the Euler characteristic formula (see op.~cit.)
$$\sum_{i\geq 0}(-1)^i \dim H_{\Gal_F}^i(X_{\bar F},\mathcal{V})=-{\rm rk}\,\mathcal{V} [F:\mathbf{Q}_p]$$ implies that there is always a non-vanishing $H^1$ (unless of course $\mathcal{V}=0$). In particular we cannot twist away higher cohomology. (Without taking the $\Gal_F$-equivariance into account, there is an analogue of an ample line bundle on $X_{\bar F}$. This implies that we can at least twist the vector bundle $\mathcal{V}$ such that its $H^2$ vanishes.)
\item[(iv)] Local Tate duality \cite[Theorem 4.4.5 (3)]{MR3230818} asserts that given a $(\varphi,\Gamma)$-module over $\mathcal{R}_{F,A}$ there is a canonical isomorphism
\[R\Gamma_{\varphi,\Gamma}(D)\isoto R{\rm Hom}_{A}(R\Gamma_{\varphi,\Gamma}(D^\vee(\varepsilon)),A)[-2].\]
Here $D^\vee={\rm Hom}_{\mathcal{R}_{F,A}}(D,\mathcal{R}_{F,A})$ and $\varepsilon$ is the cyclotomic character.
\item[(v)] There is a theory of Hodge--Tate--Sen weights, which we explain in Section \ref{sec:HTS-wts} below. It introduces substantial differences in the discussion of $P$-structures, for a parabolic subgroup $P\subseteq {\rm GL}_d$, in Section \ref{sec:drinf-style-comp} below.
\end{enumerate}
\end{remark}

\begin{conj}\label{EHConj:Artinstack}
The stack $\mathfrak{X}_d$ is a rigid analytic Artin stack of dimension $d^2[F:\mathbf{Q}_p]={\rm dim}\ {\rm Res}_{F/\mathbf{Q}_p}{\rm GL}_d$,
which is  furthermore
normal and a local complete intersection. 
\end{conj}

We now describe some partial results in the direction of this conjecture.

\begin{thm}\label{EHthm:chartsforXd}
\hfill\\(i) The diagonal $\Delta_{\mathfrak{X}_d}:\mathfrak{X}_d\rightarrow \mathfrak{X}_d\times\mathfrak{X}_d$ of $\mathfrak{X}_d$ is representable by rigid spaces.\\
(ii) Let $L'$  be a finite extension of $L$ and let $x\in\mathfrak{X}_d(L')$ be an $L'$-valued point. Then there exists an open neighborhood $\mathfrak{U}_x$ of $x$ in $\mathfrak{X}_d$ which is a rigid analytic Artin stack.
Moreover, if $d=2$, then $\mathfrak{U}_x$ has the expected dimension and is normal and a local complete intersection.  
\end{thm}
\begin{proof}[Sketch of Proof.]
(i) This basically follows from the fact that $R\Gamma_{\varphi,\Gamma}(D)$ is universally computed by a perfect complex. More precisely, given two families $D_1,D_2$ of $(\varphi,\Gamma)$-modules over an affinoid rigid analytic space ${\rm Sp}\,A$, we have to show that the functor
$$\underline{{\rm Isom}}(D_1,D_2):{\rm Sp}\,B\longmapsto {\rm Isom}_{\varphi,\Gamma}(D_1\widehat{\otimes}_AB,D_2\widehat{\otimes}_AB)$$
on the category of affinoid spaces over ${\rm Sp}\,A$ is representable by a rigid analytic space. 

By looking at the induced map on the top exterior powers it follows that $\underline{\rm Isom}(D_1,D_2)$ is an open subfunctor of
\begin{align*}
\underline{{\rm Hom}}(D_1,D_2)=\underline{H}^0_{\varphi,\Gamma}(D_1^\vee\otimes D_2):{\rm Sp}\,B\longmapsto &{\rm Hom}_{\varphi,\Gamma}(D_1\widehat{\otimes}_AB,D_2\widehat{\otimes}_AB)\\=&H^0_{\varphi,\Gamma}((D_1^\vee\otimes D_2)\widehat{\otimes}_AB).
\end{align*}
Let us write $D=D_1^\vee\otimes D_2$. Then there is a complex of vector bundles $\mathcal{E}^\bullet$ concentrated in degrees $[0,2]$ together with an isomorphism
$\underline{H}^0_{\varphi,\Gamma}(D)\cong \underline{{\rm ker}}(\mathcal{E}^0\rightarrow \mathcal{E}^1)$, where the functor on the right hand side maps ${\rm Sp}\,B$ to ${\rm ker}(\mathcal{E}^0\otimes_AB\rightarrow \mathcal{E}^1\otimes_AB)$. 
As $\mathcal{E}^0$ and $\mathcal{E}^1$ are vector bundles, this functor is representable by the preimage of the zero section under the induced map on the corresponding geometric vector bundles over ${\rm Sp}\,A$. 

(ii) We sketch the proof in the case $d=2$. Obviously we can reduce to the case $L'=L$. Let $D$ be a $(\varphi,\Gamma)$-module of rank~$2$ over $\mathcal{R}_{F,L}$. In order to construct a chart, i.e.~a smooth surjection from a rigid analytic space to a neighborhood of $D$, we distinguish two possible cases:
\begin{enumerate}
\item[(a)] There exists a character $\delta:F^\times\rightarrow L^\times$ such that $D'=D\otimes \mathcal{R}_{F,L}(\delta)$ is \'etale. 
\item[(b)] No such character exists.
\end{enumerate}

In case (a) we can reduce to the case of Galois representations and use the generic fiber of a (framed) Galois deformation ring in order to construct a chart locally at $D$.
In case (b) we want to reduce to the case of Galois representations as well: 
after possibly replacing $D$ by its twist by a character, we construct an extension
\numequation\label{EHeqn auxiliaryextensionofD}
0\rightarrow D\rightarrow D'\rightarrow \mathcal{R}_{F,L}(\delta)\rightarrow 0
\end{equation}
such that $D'$ is the $(\varphi,\Gamma)$-module attached to a $3$-dimensional Galois representation $\rho'$. 
Then $\rho'$ defines a point in the generic fiber $({\rm Spf} R_{\bar\rho'})^{\rm rig}$ of the (framed) deformation space of some residual $3$-dimensional Galois representation $\bar\rho'$, and (locally around $\rho'$) the space $({\rm Spf} R_{\bar\rho'})^{\rm rig}$ contains a closed subspace $Y$ over which we can deform the extension (\ref{EHeqn auxiliaryextensionofD}). Mapping the universal Galois representation on $Y$ to the subobject in this extension defines the desired morphism to $\mathfrak{X}_2$ which we will prove to be smooth at $\rho'$.
We now explain the argument in cases (a) and (b) in more detail:

Case (a): As $D'$ is \'etale there exists a Galois representation $\rho':\Gal_F\rightarrow {\rm GL}_2(L)$ with $D_{\rm rig}(\rho')=D'$. After conjugating $\rho'$ by some element of ${\rm GL}_2(L)$ if necessary, we may assume that $\rho'$ takes values in ${\rm GL}_2(\mathcal{O})$ (i.e.~we choose a stable lattice in the Galois representation). Let $\bar\rho':\Gal_F\rightarrow {\rm GL}_2(k)$ denote the reduction of $\rho'$ and let $R_{\bar\rho'}$ denote the universal framed deformation ring of $\bar\rho'$. Then $\rho'$ defines a point of the rigid analytic generic fiber $({\rm Spf}\,R_{\bar\rho'})^{\rm rig}$ and we are left to show that the morphism
\begin{align*}
({\rm Spf}\,R_{\bar\rho'})^{\rm rig}\longrightarrow & \mathfrak{X}_2\\
\rho\longmapsto &D_{\rm rig}(\rho)\otimes \mathcal{R}(\delta^{-1})
\end{align*}
that maps $\rho'$ to $D$ is smooth in a neighborhood of $\rho'$.  This follows from the infinitesimal lifting criterion together with the fact that the complete local ring of $({\rm Spf}\,R_{\bar\rho'})^{\rm rig}$ at $\rho'$ pro-represents the universal framed deformation functor of $\rho':\Gal_F\rightarrow{\rm GL}_2(L)$, see \cite[Lemma 2.3.3, Prop. 2.3.5]{MR2600871}.

\noindent Case (b): After again twisting with some rank $1$ object we may assume that all slopes of subobjects of $D$ are non-negative. By Lemma \ref{EHlemexistenceofetaleextension} below there exists a character $\delta:F^\times\rightarrow L^\times$ and a non-split extension 
$$0\longrightarrow D\longrightarrow D'\longrightarrow \mathcal{R}_{F,L}(\delta)\longrightarrow 0$$
such that $D'$ is \'etale and such that ${\rm Ext}^0_{\varphi,\Gamma}(\mathcal{R}_{F,L}(\delta),D)={\rm Ext}^2_{\varphi,\Gamma}(\mathcal{R}_{F,L}(\delta),D)=0$, as well as ${\rm dim}\ {\rm Hom}_{\varphi,\Gamma}(D',\mathcal{R}_{F,L}(\delta))=1$.
 Again we choose an $\mathcal{O}$-lattice $\rho'$ in the Galois representation associated to $D'$ and consider $\rho'$ as a point of $({\rm Spf}\,R_{\bar\rho'})^{\rm rig}$ the rigid analytic generic fiber of the universal framed deformation ring of the reduction $\bar\rho'$ of $\rho'$. 

We construct a locally closed subset $Y\subseteq ({\rm Spf}\,R_{\bar\rho'})^{\rm rig}$ containing $\rho'$ and a smooth morphism $Y\rightarrow \mathfrak{X}_2$ as follows:
Let $\tilde D$ denote the universal $(\varphi,\Gamma)$-module over $({\rm Spf}\,R_{\bar\rho'})^{\rm rig}$ and let  $Y_1$ denote the scheme-theoretic support of $H^2_{\varphi,\Gamma}(\tilde D(\delta^{-1}\varepsilon))$. 
As 
$$H^2_{\varphi,\Gamma}(\tilde D(\delta^{-1}\varepsilon))\otimes k(\rho')=H^2_{\varphi,\Gamma}(D'(\delta^{-1}\varepsilon))=({\rm Hom}_{\varphi,\Gamma}(D',\mathcal{R}_{F,L}(\delta)))^\vee$$
is of dimension $1$, we can find a neighborhood $Y_2\subseteq Y_1$ of $\rho'$ and a lift  $\tilde f\in \Gamma(Y_2,H^2_{\varphi,\Gamma}(\tilde D(\delta^{-1}\omega)))$ of the dual basis to a chosen basis vector $$f^\vee\in{\rm Hom}_{\varphi,\Gamma}(D',\mathcal{R}_{F,L}(\delta)).$$
After replacing $Y_2$ by a smaller open neighborhood of $\rho'$ if necessary, we may assume that $H^2_{\varphi,\Gamma}(\tilde D(\delta^{-1}\varepsilon))|_{Y_2}$ is free of rank~$1$ with basis $\tilde f$. 
Localizing further we find that there is a neighborhood $Y\subseteq Y_2$ of $\rho'$ together with a surjection 
$$\tilde f^\vee:\tilde D|_Y\longrightarrow \mathcal{R}_{F,Y}(\delta)$$
that specializes to the chosen surjection $f^\vee:D'\rightarrow \mathcal{R}_{F,L}(\delta)$ at $\rho'$ and we moreover may assume that at all points $y\in Y$ we have 
\numequation\label{EHtechnicalExtvanishing}
{\rm Ext}^0(\mathcal{R}_{F,k(y)}(\delta),\tilde D\otimes k(y))={\rm Ext}^2(\mathcal{R}_{F,k(y)}(\delta),\tilde D\otimes k(y))=0.
\end{equation}
We claim that the map $Y\rightarrow \mathfrak{X}_2$ defined by $\tilde D|_Y\mapsto {\rm ker}\, \tilde f^\vee$ is smooth. Again it is enough to check this on the complete local ring of a point $y\in Y$. 
The complete local ring at $y$ pro-represents the functor of deformations of $\tilde D_y=\tilde D\otimes k(y)$ together with a deformation of the map $f_y^\vee:\tilde D_y\rightarrow \mathcal{R}_{F,k(y)}(\delta)$. 
By our choices of $Y$ this functor may be identified with the functor of deformations of $D_y={\rm ker}\, f_y^\vee$ together with a deformation of the extension
$$0\longrightarrow D_y\longrightarrow \tilde D_y\longrightarrow \mathcal{R}_{F,L}(\delta)\longrightarrow 0.$$
 By the vanishing of the Ext-groups (\ref{EHtechnicalExtvanishing}) this functor is formally smooth over the deformation functor of $D_y$.

We are left to prove the statement about the dimension and the claims on the local structure. This is Proposition \ref{EHproplocalstructureXd} below.
\end{proof}

\begin{lem}\label{EHlemexistenceofetaleextension}
Let $D$ be a $(\varphi,\Gamma)$-module of rank~$2$ over
$\mathcal{R}_{F,L}$ and assume that all slopes of subobjects of $D$
are non-negative and that $D\not\cong\mathcal{R}_{F,L}(\delta)^{\oplus
2}$ for any character $\delta:F^\times\rightarrow L^\times$.
Let $\eta:F^\times\rightarrow L^\times$ be a continuous character such that 
\begin{enumerate}
\item[-] for all embeddings $\tau:F\hookrightarrow L$ one has $({\rm wt}_\tau(\eta)+\mathbf{Z})\cap {\rm wt}_\tau(D)=\emptyset$, where we refer to Section {\em  \ref{sec:HTS-wts}} below for the definition of the weights. 
\item[-] the slope of $\eta$ is  ${\rm slope}(\eta)=-{\rm slope}(D)$.
\end{enumerate}
Then there exists a non-split extension 
\numequation\label{EHetaleextension}
0\longrightarrow D \longrightarrow D' \longrightarrow \mathcal{R}_{F,L}(\eta)\longrightarrow 0
\end{equation}
such that $D'$ is \'etale and moreover
\begin{equation*}
\begin{aligned}
{\rm Ext}^0_{\varphi,\Gamma}(\mathcal{R}_{F,L}(\eta),D)={\rm Ext}^2_{\varphi,\Gamma}(\mathcal{R}_{F,L}(\eta),D)&=0,\\
\dim_L {\rm Hom}_{\varphi,\Gamma}(D',\mathcal{R}_{F,L}(\eta))&=1.
\end{aligned}
\end{equation*}
\end{lem}
\begin{rem}
Given $D$ as in the lemma, a character $\eta$ satisfying the assumptions in the lemma obviously always exists. 
\end{rem}
\begin{proof}
As  $D(\eta^{-1})$ and $D^\vee(\eta\omega)$ have non-integral Hodge--Tate--Sen weights, the vanishing of the Ext-groups is a direct consequence of the identifications
\begin{align*}
{\rm Ext}^0_{\varphi,\Gamma}(\mathcal{R}_{F,L}(\eta),D)&=H^0_{\varphi,\Gamma}(D(\eta^{-1}))\\
{\rm Ext}^2_{\varphi,\Gamma}(\mathcal{R}_{F,L}(\eta),D)&=H^0_{\varphi,\Gamma}(D^\vee(\eta\omega)^\vee.
\end{align*}

Choose any non-split extension (\ref{EHetaleextension}). The claim on the Hodge--Tate weights also implies the claim
on the dimension of ${\rm
  Hom}_{\varphi,\Gamma}(D',\mathcal{R}_{F,L}(\eta))$, so we only need
to show that~$D'$ can be chosen to be \'etale.  By the condition on the slope of $\eta$ we automatically have ${\rm slope}(D')=0$ and by Kedlaya's slope filtration theorem \cite{MR2493220} we are left to show that for all (saturated) subobjects $D''\subseteq D$ we have ${\rm slope}(D'')\geq 0$. 
By assumption this is true if $D''\subseteq D$. It therefore suffices to
show that we can choose $D'$ such that all saturated
$(\varphi,\Gamma)$-stable subobjects of $D'$ are contained in $D$. Assume to the contrary that that there is some $D''$ such that the map $D''\rightarrow \mathcal{R}_{F,L}(\eta)$ is non-zero. 
We distinguish the cases ${\rm rk}\, D''=1$ and ${\rm rk}\,D''=2$. 

If ${\rm rk}\,D''=1$, then $D''\rightarrow \mathcal{R}_{F,L}(\eta)$ is an isomorphism after inverting $t$, which implies that that the extension (\ref{EHetaleextension}) is split after inverting $t$. 
However, the assumption $({\rm wt}_\tau(\eta)+\mathbf{Z})\cap {\rm wt}_\tau(D)=\emptyset$ implies that the canonical map 
$${\rm Ext}^1_{\varphi,\Gamma}(\mathcal{R}_{F,L},D)\longrightarrow {\rm Ext}^1_{\varphi,\Gamma}(\mathcal{R}_{F,L}[1/t],D[1/t])$$
is an isomorphism, which is a contradiction, as we chose (\ref{EHetaleextension}) to be non-split.

If ${\rm rk}\, D''=2$, then $E=D''\cap D\subseteq D$ is a $(\varphi,\Gamma)$-submodule of rank~$1$. By assumption, $D$ has only finitely many saturated subobjects of rank~$1$. Each such subobject $E$ gives rise to a short exact sequence 
\numequation\label{EHDasanextension}
0\longrightarrow \mathcal{R}_{F,L}(\delta_1)\longrightarrow D\longrightarrow \mathcal{R}_{F,L}(\delta_2)\longrightarrow 0
\end{equation}
and our assumptions on the Hodge--Tate--Sen weights implies the vanishing of certain ${\rm Ext}^0$ and ${\rm Ext}^2$ groups and hence a short exact sequence
\begin{align*}
0\rightarrow {\rm Ext}^1_{\varphi,\Gamma}(\mathcal{R}_{F,L}(\eta),\mathcal{R}_{F,L}(\delta_1))&\rightarrow {\rm Ext}^1_{\varphi,\Gamma}(\mathcal{R}_{F,L}(\eta),D)\rightarrow \\ &\rightarrow {\rm Ext}^1_{\varphi,\Gamma}(\mathcal{R}_{F,L}(\eta),\mathcal{R}_{F,L}(\delta_2))\rightarrow 0.
\end{align*}
Moreover, as in the rank $1$ case, the condition $({\rm wt}_\tau(\eta)+\mathbf{Z})\cap {\rm wt}_\tau(D)=\emptyset$ implies that inverting $t$ induces an isomorphism between this exact sequence and the exact sequence
\begin{align*}
0\rightarrow {\rm Ext}^1_{\varphi,\Gamma}(\mathcal{R}_{F,L}(\eta)[1/t],\mathcal{R}_{F,L}(\delta_1)[1/t])&\rightarrow {\rm Ext}^1_{\varphi,\Gamma}(\mathcal{R}_{F,L}(\eta)[1/t],D[1/t])\rightarrow \\ &\rightarrow {\rm Ext}^1_{\varphi,\Gamma}(\mathcal{R}_{F,L}(\eta)[1/t],\mathcal{R}_{F,L}(\delta_2)[1/t])\rightarrow 0.
\end{align*}
As $D$ has only finitely many subobjects, we can choose the extension $$D'\in{\rm Ext}^1_{\varphi,\Gamma}(\mathcal{R}_{F,L}(\eta),D)$$ such that for all possibilities to write $D$ as an extension (\ref{EHDasanextension}), the extension $D'$ is not in the image of ${\rm Ext}^1_{\varphi,\Gamma}(\mathcal{R}_{F,L}(\eta),\mathcal{R}_{F,L}(\delta_1))$. 

In this case, if $D''\subseteq D'$ is a saturated subobject with $D''\cap D=E=\mathcal{R}_{F,L}(\delta_1)$ and such that $D''\rightarrow \mathcal{R}_{F,L}(\eta)$ is non-zero, then $$\mathcal{R}_{F,L}(\eta)[1/t]\cong (D''/E)[1/t]\subseteq (D'/E)[1/t]$$
implies that the image $(D'/E)[1/t]$ of $D'$ in ${\rm Ext}^1_{\varphi,\Gamma}(\mathcal{R}_{F,L}(\eta)[1/t],\mathcal{R}_{F,L}(\delta_2)[1/t])$ vanishes. This contradicts our choice of $D'$ above. \end{proof}
The existence of local charts implies the existence of versal rings at rigid analytic points of $\mathfrak{X}_d$.
\begin{prop}\label{EHproplocalstructureXd}
Let $d=2$ and $x\in\mathfrak{X}_d$ and let $R$ be a versal ring at
$x$ such that ${\rm Spf}\,R\rightarrow \mathfrak{X}_d$ is formally smooth of relative dimension $m$. Then $R$ is a local complete intersection of dimension
$m+d^2[F:\Q_p]$ and is normal.
\end{prop}
\begin{proof}[Sketch of proof.]
As any two versal rings are smoothly equivalent, we may use any choice of a versal ring to prove the assertion. 
We use the notation from the proof of Theorem \ref{EHthm:chartsforXd}. 
If $x$ is as in case (a) of that proof, then $R$ can be chosen to be
isomorphic to the complete local ring of $({\rm
  Spf}\,R_{\bar\rho'})^{\rm rig}$ at the point $\rho'$. The result
then follows from~\cite{BIP}.

In case~(b), where $x$ is not the twist of an \'etale
$(\varphi,\Gamma)$-module,  standard deformation theory arguments give
us a presentation of a versal ring $R$ such that ${\rm Spf}\,
R\rightarrow \mathfrak{X}_d$ is formally smooth of dimension $m$ and
such that every irreducible component of $R$ has dimension at least
$m+[F:\mathbf{Q}_p]d^2$, and in order to show that~$R$ is lci, we need to show that for all irreducible
components equality holds. By the existence of local charts (i.e.\ Theorem~\ref{EHthm:chartsforXd}~(2)), we may view $R$ as the complete local ring
of a point $y\in Y={\rm Sp}\,A$, where $f:Y={\rm Sp}\, A\rightarrow
\mathfrak{X}_d$ is smooth of relative dimension $m$. %
 After restricting to an open subset of $Y$ we may assume that all irreducible components of $Y$ contain $y$ and we are left to show that all irreducible components of $Y$ have dimension at most $m+[F:\mathbf{Q}_p]d^2$. 
 Let $Z$ be such an irreducible component, then generically on $Z$, i.e.~on some Zariski-open subset $U_Z\subset Z$, the restriction $f_Z:Z\rightarrow \mathfrak{X}_d$ is formally smooth of relative dimension $m$.
 If $U_Z$ contains a point $z$ such that the $(\varphi,\Gamma)$-module $f_Z(z)$ is étale up to twist, then the complete local ring of $U_Z$ at $z$ has dimension $m+[F:\mathbf{Q}_p]d^2$ by the discussion of case (a), and hence $Z$ has dimension $m+[F:\mathbf{Q}_p]d^2$. 
As a $(\varphi,\Gamma)$-module of rank $2$ that is not étale up to twist is necessarily reducible (using Kedlaya's slope filtration theorem \cite{MR2493220} that asserts that being étale up to twist is the same as being semi-stable), this means that we have to rule out the possibility that (generically on an irreducible component $Z$) all points $z$ have the property that the $(\varphi,\Gamma)$-module $f(z)$ is reducible.

As $Z$ is of dimension at least $m+[F:\mathbf{Q}_p]d^2$ and is (generically) smooth of relative dimension $m$ over $\mathfrak{X}_d$, we can rule out this possibility by proving that the dimension of space of reducible $(\varphi,\Gamma)$-modules of rank $2$ is too small. 
We sketch the argument on the level of deformation spaces\footnote{Strictly speaking we have to show something slightly stronger, namely that the subset of the rigid analytic space $Z$ where the $(\varphi,\Gamma)$-modules are pointwise reducible is a countable union of Zariski-closed subspaces of dimension less or equal to $m+[F:\mathbf{Q}_p]\tfrac{d(d+1)}{2}<m+[F:\mathbf{Q}_p]d^2$. Compare Theorem \ref{EHThm propernessofcompactific} and Remark \ref{EH:rem density of XB} below for these statements.}.
Given an extension (\ref{EHDasanextension}) we
consider the groupoid of deformations $\mathfrak{D}$ of this short
exact sequence and denote by $S$ a versal ring to $\mathfrak{D}$ such
that ${\rm Spf}\,S\rightarrow \mathfrak{D}$ is formally smooth of
dimension $m$. Then it is enough to show that all irreducible
components of $S$ are of dimension strictly less than
$m+[F:\mathbf{Q}_p]d^2$. In the case $d=2$ this can easily be checked
by a direct computation (compare also the discussion of the  stack of
$B$-bundles in \ref{sec:drinf-style-comp} below).

For normality we need to check in addition (if $F=\Q_p$) that for a generic choice of a reducible $(\varphi,\Gamma)$-module $$0\rightarrow \mathcal{R}_{F,L}(\delta_1)\rightarrow D\rightarrow \mathcal{R}_{F,L}(\delta_2)\rightarrow 0$$ the stack $\mathfrak{X}_2$ is smooth at $D$. Again this follows from a tangent space computation (resp.~a computation of Ext-groups) together with the fact that generically on the stack of $B$-bundles, the ratio $\delta_1/\delta_2$ is very regular in the sense that 
\[{\rm Ext}^0_{\varphi,\Gamma}(\mathcal{R}_{F,L}(\delta_2),\mathcal{R}_{F,L}(\delta_1))={\rm Ext}^2_{\varphi,\Gamma}(\mathcal{R}_{F,L}(\delta_2),\mathcal{R}_{F,L}(\delta_1))=0.\qedhere\]
\end{proof}
\begin{rem}
Basically the same strategy (using not just the Borel $B$ but all parabolic subgroups) should settle the claim on the dimension and the local structure for general $d$. Only the computation of the dimension of spaces/stacks of extensions becomes more involved and is not finished yet. We note that similar computations of dimensions of Ext-groups also show up in the proof of the main theorem of~\cite{BIP}.
\end{rem}

We elaborate briefly on the connection with the integral theory presented in Section  \ref{subsec:  overview of EG book}. 
Let ${\rm Spf}\,A\rightarrow \mathcal{X}_d$ be an $A$-valued point, for some $\mathbf{Z}_p\langle T_1,\dots, T_m\rangle\twoheadrightarrow A$. This $A$-valued point defines a $(\varphi,\Gamma)$-module $D$ over the $p$-adically complete ring $\mathbf{A}_{F,A}=A\widehat{\otimes}_{\mathbf{Z}_p}\mathbf{A}_F$ (compare Remark \ref{EH rem peropdrings} for the notation) and we would like to be able to (functorially) associate a $(\varphi,\Gamma)$-module $D_{\rm rig}$ over $\mathcal{R}_{F,A[1/p]}$ to $D$. 
As there is no map from $\mathbf{A}_{F,A}$ to $\mathcal{R}_{F,A[1/p]}$ one has to show first that $D$ is overconvergent, i.e.~that it admits a canonical model $D^\dagger$ over the subring $\mathbf{A}^\dagger_{F,A}=A\widehat{\otimes}_{\mathbf{Z}_p}\mathbf{A}^\dagger_F$ (for an appropriately completed tensor product). If $D$ is the $(\varphi,\Gamma)$-module associated to a family of ${\rm Gal}_F$-representations over $A$, then this is proven by Berger--Colmez \cite[Cor. 4.2.7]{MR2493221}. The general case is a theorem of Gal Porat~\cite{galoverconvergent}:
\begin{theorem}\label{EHConj:overconvergence}
Let $A$ be a $p$-adically complete $\mathbf{Z}_p$-algebra topologically of finite type over $\mathbf{Z}_p$. Then every \'etale $(\varphi,\Gamma)$-module over $\mathbf{A}_{F,A}$ is overconvergent.
\end{theorem}
As a consequence of this theorem, we can construct a morphism from the rigid analytic generic fiber of the stack $\mathcal{X}_d$ defined in Section  \ref{subsec:  overview of EG book} to the stack $\mathfrak{X}_d$.
\begin{cor}
There is a morphism 
\numequation\label{EHgenfibEGstacktoanalytic}
\pi_d:\mathcal{X}_{d,\eta}^{\rm rig}\longrightarrow \mathfrak{X}_d
\end{equation}
 given by mapping a $(\varphi,\Gamma)$-module $D\in\mathcal{X}_d({\rm Spf}\,A)$ (for some $p$-adically complete $\mathbf{Z}_p$-algebra $A$ that is topologically of finite type) to $D^\dagger\widehat{\otimes}_{\mathbf{A}^\dagger_{F,A}}\mathcal{R}_{F,A[1/p]}$.
\end{cor}

We expect that the morphism $\pi_d$ is formally \'etale, but not
representable by rigid analytic spaces (nor by quasi-analytic spaces).
The reason is that the morphism is basically given by forgetting the étale lattice in a $(\varphi,\Gamma)$-module over the Robba ring, and hence the fibers of the morphism parametrize the choices of such a lattice.  
In the case of ${\rm GL}_1$ this implies that the fibers look like the stack quotient $\mathbf{G}_m/(\widehat{\mathbf{G}}_m)^{\rm rig}$, where $\widehat{\mathbf{G}}_m$ is the formal multiplicative group over ${\rm Spf}\,\mathbf{Z}_p$, see Section \ref{sec: GL1} for more details. 

\begin{rem}
We point out that the proof of Theorem \ref{EHthm:chartsforXd} only
gives charts locally around rigid analytic points. As a set theoretic
cover of a rigid analytic space by admissible open subsets is not
necessarily an admissible cover, this is not enough to give a full
proof of Conjecture \ref{EHConj:Artinstack}. Extending the argument in
the above proof to all points of the corresponding adic space (not
just those points $x$ such that $L'=k(x)$ is finite over $L$) would
prove Conjecture \ref{EHConj:Artinstack}. The given argument would
directly generalize if one had an analogue of Kedlaya's slope
filtration theorem for $(\varphi,\Gamma)$-modules over
$\mathcal{R}_{F,L'}$, where $L'$ is allowed to be the (completed)
residue field of any point in an adic space of finite type over
$\mathbf{Q}_p$ and if moreover one could prove that the map $\pi_d$ (that exists due to \cite{galoverconvergent}) is formally \'etale. While \'etaleness is not expected to be difficult to prove, we have no idea whether the generalization of Kedlaya's theorem to general coefficient fields is true or not.\end{rem}

\subsubsection{Hodge--Tate--Sen weights and global sections} 
\label{sec:HTS-wts}
The theory of Hodge--Tate--Sen weights for representations of $\Gal_F$ on finite dimensional $\mathbf{Q}_p$-vector spaces generalizes to families of $(\varphi,\Gamma)$-modules, respectively to families of $\Gal_F$-equivariant vector bundles on $X_{\bar F}$. 
Let $A$ denote an affinoid algebra and let $D_A$ be a $(\varphi,\Gamma)$-module $D_A$ over $\mathcal{R}_{F,A}$. We write $\mathcal{E}_A=\mathcal{E}_{F_\infty}(D_A)$ for the corresponding $\Gamma$-equivariant vector bundle on $X_{F_\infty,A}$. 
We then consider the special fiber
$$\widehat{D}_{\rm Sen}(D_A)=\mathcal{E}_A|_{{\rm Sp}\,A\times\{\infty\}}=D_A\widehat{\otimes}_{\mathcal{R}_{F,A}}(\widehat{ F}_\infty\widehat{\otimes}_{\mathbf{Q}_p}A)$$
where the tensor product on the right hand side is induced by a certain (explicit) map $\mathcal{R}_{F}\rightarrow \widehat{ F}_\infty$. The map $\mathcal{R}_F\rightarrow \widehat{F}_\infty$ factors through $F_\infty\subset\widehat{F}_\infty$ and the module $\widehat{D}_{\rm Sen}(D_A)$ has a decompleted version
$$D_{\rm Sen}(D_A)=D_A\otimes_{\mathcal{R}_{F,A}}(A\otimes_{\Q_p}F_\infty)$$
which is a finite projective $A{\otimes}_{\mathbf{Q}_p}F_\infty$-module that carries a continuous action of the $1$-dimensional $p$-adic Lie group $\Gamma$. In fact $D_{\rm Sen}(D_A)$ can even be defined over $A\otimes_{\mathbf{Q}_p}F_n$ for sufficiently large $n$ depending on $D_A$, i.e.~there exists some $n\gg 0$ and a (unique) $A\otimes_{\mathbf{Q}_p}F_n$-module $D_{{\rm Sen},n}(D_A)$ with semi-linear continuous $\Gamma$-action such that 
$$D_{\rm Sen}(D_A)=D_{{\rm Sen},n}(D_A)\otimes_{F_n}F_\infty.$$ 
On the decompleted modules $D_{\rm Sen}(D_A)$ and $D_{{\rm Sen},n}(D_A)$ the $\Gamma$-action is locally analytic, but we remind the reader that it is not smooth, i.e.~it is not a descent datum to an $A\otimes_{\mathbf{Q}_p}F$-module.
 
As the action of $\Gamma_n\subset \Gamma$ on $D_{{\rm Sen},n}(D_A)$ is $A\otimes_{\mathbf{Q}_p}F_n$-linear and locally analytic one can derive this action,  defining a morphism
$${\rm Lie}\,\Gamma={\rm Lie}\,\Gamma_n\longrightarrow {\rm Lie}({\rm GL}(D_{{\rm Sen},n}(D_A)))={\rm End}_{A\otimes_{\mathbf{Q}_p}F_n}(D_{{\rm Sen},n}(D_A)).$$

The cyclotomic character induces a trivialization $\mathbf{Z}_p\cong {\rm Lie}\,\Gamma$ and hence the image of $1\in\mathbf{Z}_p$ defines (after extending scalars from $F_n$ back to $F_\infty$) an $A\otimes_{\mathbf{Q}_p}F_\infty$-linear endomorphism 
$$\Theta:D_{\rm Sen}(D_A)\longrightarrow D_{\rm Sen}(D_A)$$ 
that is independent of the choice of $n$ in the preceding discussion. 
This endomorphism commutes with the $\Gamma$-action, so the coefficients of its characteristic polynomial take values in $F\otimes_{\mathbf{Q}_p}A$. 
\begin{remark}
We caution the reader that, while the characteristic polynomial of $\Theta$ always has coefficients in $A\otimes_{\mathbf{Q}_p}F$, this does not imply that the $\Gamma$-action can be used to (functorially) descend the pair $(D_{\rm Sen}(D_A),\Theta)$ to a finite projective $A\otimes_{\Q_p}F$-module together with an $A\otimes_{\Q_p}F$-linear endomorphism. The reason is that a continuous semi-linear $\Gamma$-action is not a descent datum --- only a smooth semi-linear $\Gamma$-action is. 
Given this statement one might try to twist the $\Gamma$-action by the inverse of 
$$\Gamma\xrightarrow{\epsilon_{\rm cyc}} \mathbf{Z}_p^\times\xrightarrow{\rm log}\mathbf{Z}_p\rightarrow (A\otimes_{\mathbf{Q}_p}F_\infty)^\times,$$
where $\epsilon_{\rm cyc}$ is the cyclotomic character and the last map is given by $1\mapsto {\rm exp}(\Theta)$, in order to obtain a smooth action and  hence a descent datum. However, this only makes sense if ${\rm exp}(\Theta)$ is convergent which is in general not the case. Convergence can always be assured by replacing $\Theta$ with $p^n\Theta$ for $n\gg 0$, but this way we can only twist the restriction of the $\Gamma$-action to $\Gamma_n\subset \Gamma$ and obtain a descent datum from $F_\infty$ to $F_n$. This descent yields exactly the module $D_{{\rm Sen},n}(D_A)$ over $A\otimes_{\mathbf{Q}_p}F_n$ from the discussion above. 
\end{remark}
If $A$ is an $L$-algebra and $|{\rm Hom}(F,L)|=[F:\mathbf{Q}_p]$, then
$F\otimes_{\mathbf{Q}_p} A\cong \prod_{\tau:F\hookrightarrow L}A_\tau$, where
$A_\tau=A$. For each embedding $\tau$ we obtain an endomorphism
$\Theta_\tau$ of the finite projective $A\otimes_{\tau,F}F_\infty$-module $$D_{{\rm
    Sen},\tau}(D_A)=D_{\rm Sen}(D_A)\otimes_{A\otimes_{\Q_p}F}A_\tau$$ with characteristic polynomial $\chi_{\Theta_\tau}$ with coefficients in $A$. The zeros of this polynomial (i.e.~the eigenvalues of $\Theta_\tau$) are called the Hodge--Tate--Sen weights of $D$ labeled by $\tau$. If $A=L$ is a field containing all the zeros of $\chi_{\Theta_\tau}$, then we write ${\rm wt}_\tau(D)\subseteq L$ for the set of all these zeros.

We write $\mathfrak{gl}_n/\!/{\rm GL}_n$ for the GIT quotient of the Lie algebra $\mathfrak{gl}_n$ by the adjoint action of ${\rm GL}_n$.
Mapping a $(\varphi,\Gamma)$-module $D_A$ over $\mathcal{R}_{F,A}$ to the coefficients of the characteristic polynomial $\chi_\Theta$ of $\Theta$ defines a map
\numequation\label{EHeqnweightmap}
\omega_d:\mathfrak{X}_d\longrightarrow \big({\rm Res}_{F/\mathbf{Q}_p} (\mathfrak{gl}_n/\!/{\rm GL}_n)\big)_L={\rm WT}_{d,L}.
\end{equation}
If $|{\rm Hom}(F,L)|=[F:\mathbf{Q}_p]$, then of course ${\rm WT}_{d,L}=\prod_{\tau:F\rightarrow L}{\rm WT}_{\tau,d,L}$ and each ${\rm WT}_{\tau,d,L}=(\mathfrak{gl}_n/\!/{\rm GL}_n)_L$ is an affine space.

We finish this subsection by conjecturing that, up to functions pulled back via the map induced by the determinant ${\rm GL}_d\rightarrow {\rm GL}_1$, all global functions on $\mathfrak{X}_d$ come by pullback along the map $\omega_{d}$. For the precise formulation we need to discuss the determinant map:
taking the top exterior power of an equivariant vector bundle defines a map 
\[{\rm det}:\mathfrak{X}_d\longrightarrow \mathfrak{X}_1.\]
By the discussion in Section~\ref{sec: GL1} below (see in particular (\ref{EHeqn X1analytic})) the stack $\mathfrak{X}_1$ can be written as the stack quotient $\mathcal{T}/\mathbf{G}_m$, where $\mathcal{T}$ is the space of continuous characters of $F^{\times}$ that is equipped with the trivial action of the rigid analytic group $\mathbf{G}_m$. As the $\mathbf{G}_m$-action on $\mathcal{T}$ is trivial $\mathfrak{X}_1=\mathcal{T}\times \ast/\mathbf{G}_m$ admits a canonical projection to $\mathcal{T}$, in particular pullback along this map induces an isomorphism $\Gamma(\mathfrak{X}_1,\mathcal{O}_{\mathfrak{X}_1})\cong \Gamma(\mathcal{T},\mathcal{O}_{\mathcal{T}})$.
By abuse of notation we will often also write 
$${\rm det}:\mathfrak{X}_d\longrightarrow \mathcal{T}$$
for the composition of ${\rm det}$ with that projection. 
Similarly to the above construction, the top exterior power defines a map ${\rm WT}_{d,L}\rightarrow {\rm WT}_{1,L}$ and $\omega_1$ induces a map $\mathcal{T}\rightarrow {\rm WT}_{1,L}$. 
\begin{conj}\label{EHconjglobalsections}
The maps $\omega_d$ and ${\rm det}$ induce an isomorphism 
\[\begin{aligned}
\Gamma(\mathfrak{X}_d,\mathcal{O}_{\mathfrak{X}_d})&\cong \Gamma({\rm WT}_{d,L},\mathcal{O}_{{\rm WT}_{d,L}})\widehat{\otimes}_{\Gamma({\rm WT}_{1,L},\mathcal{O}_{{\rm WT}_{1,L}})} \Gamma(\mathfrak{X}_1,\mathcal{O}_{\mathfrak{X}_1})\\
& \cong \Gamma({\rm WT}_{d,L},\mathcal{O}_{{\rm WT}_{d,L}})\widehat{\otimes}_{\Gamma({\rm WT}_{1,L},\mathcal{O}_{{\rm WT}_{1,L}})} \Gamma(\mathcal{T},\mathcal{O}_{\mathcal{T}}).
\end{aligned}\]
\end{conj}
As the connected components of a rigid analytic space, or a stack on rigid analytic spaces, are in bijection with the primitive idempotents in the global sections of the structure sheaf Conjecture \ref{EHconjglobalsections} gives a description of the connected components of $\mathfrak{X}_d$. Let $$\mathcal{W} =({\rm Spf}\,\mathcal{O}[[\mathcal{O}_F^\times]])^{\rm rig}$$ denote the rigid analytic space of continuous characters of $\mathcal{O}_F^\times$. Then the choice of a uniformizer of $F$ induces an isomorphism $\mathcal{T}\cong \mathcal{W}\times \mathbf{G}_m$, and we obtain the following description of the connected components of $\mathfrak{X}_d$.
\begin{cor}
Assume Conjecture~{\em \ref{EHconjglobalsections}}. Then the connected components of $\mathfrak{X}_d$ are in bijection {\em (}via the determinant{\em )} with the connected components of $\mathfrak{X}_1$ and hence with the connected components of $\mathcal{W}$.
\end{cor}
Note that the connected components of $\mathcal{W}$ are in turn in bijection with the characters $\mu(F)\rightarrow\mathcal{O}^\times$, where $\mu(F)\subset F^\times$ denotes the subgroup of roots of unity.
We point out that that a similar result for the components of universal deformation rings was proven by B\"ockle--Iyengar--Pa\v{s}k\={u}nas~\cite{BIP}. 

\subsection{Stacks of de Rham objects} \label{EHsub:sub:deRham}
The stack of $(\varphi,\Gamma)$-modules $\mathfrak{X}_d$ has closed substacks defined in terms of $p$-adic Hodge theory, similarly to the closed substacks $\cX_d^{{\rm ss},\lambdau,\tau}\subset\cX_d$.
In the case of $(\varphi,\Gamma)$-modules over the Robba ring these closed substacks can be studied in terms of (filtered) Weil--Deligne representations. In the conjectural relation of locally analytic representations with sheaves on the stack $\mathfrak{X}_d$ these closed substacks of de Rham objects (of given weight and inertial type) will play a role in the comparison with the smooth categorical Langlands conjectures.
We will now define these ``de Rham loci".

Let $F'$ be a finite Galois extension of $F$ with Galois group ${\rm Gal}(F'/F)$ and write $F'_0$ respectively $F_0$ for the maximal unramified subextension of $F'$ respectively $F$. Moreover, let $\sigma$ denote the lift of the Frobenius $x\mapsto x^p$  to $F'_0$.
We consider the stack ${\rm Mod}_{d,\varphi,N,F'/F}$ on the category of $L$-schemes that maps an $L$-algebra $A$ to the groupoid of finite projective $A\otimes_{\mathbf{Q}_p} F'_0$-modules $D$ of rank~$d$ together with
\begin{enumerate}
\item[-] an ${\rm id}\otimes\sigma$-linear automorphism $\varphi_D:D\isoto D$.
\item[-] an $A\otimes_{\mathbf{Q}_p}F'_0$-linear endomorphism $N:D\rightarrow D$ satisfying $$N\circ \varphi_D=p\varphi_D\circ \sigma^\ast N.$$
\item[-] an action of ${\rm Gal}(F'/F)$ on $D$, commuting with $\varphi_D$ and $N$, such that $$g((a\otimes f)\cdot d)=(a\otimes g(f))\cdot g(d)$$
for all $a\in A$, $f\in F'_0$, $d\in D$ and $g\in{\rm Gal}(F'/F)$. 
\end{enumerate}

We further recall the stack of $d$-dimensional Weil--Deligne
representations ${\rm WD}_{d,F}$ of $F$ which is an Artin stack on the
category of $L$-schemes, see \cite{dat2020moduli}, \cite[\S 3]{zhu2020coherent}.
Recall that for a given $L$-scheme $S$ the groupoid ${\rm WD}_{d,F}(S)$ parameterizes vector bundles $\mathcal{E}$ on $S$ together with a smooth action $\rho$ of the Weil group $W_F$ and an endomorphism $N:\mathcal{E}\rightarrow \mathcal{E}$ satisfying $Ng=q^{|\!|g|\!|}gN$. Given $F'$ we write ${\rm WD}_{d,F'/F}\subset {\rm WD}_{d,F}$ for the open and closed substack of ${\rm WD}_{d,F}$ consisting of those Weil--Deligne representations $(\mathcal{E},\rho,N)$ such that the restriction of $\rho$ to the inertia subgroup $I_{F'}$ of $F'$ is trivial. 

We note that Fontaine's construction of the Weil--Deligne representation associated to a $(\varphi,N,{\rm Gal}(F'/F))$-module \cite{MR1293977} (see also \cite[Prop. 4.1]{MR2359853}) implies that these two stacks become isomorphic over a large enough field~$L$: 
\begin{lem}\label{EHlemWDtophiNmodules}
Assume that $L$ is large enough such that $[F'_0:\mathbf{Q}_p]=|{\rm Hom}(F'_0,L)|$ and fix an embedding $F'_0\hookrightarrow L$. Then there is an isomorphism of stacks 
\[{\rm WD}_d:{\rm Mod}_{d,\varphi,N,F'/F}\isoto {\rm WD}_{d,F'/F}.\]
\end{lem}
We point out that the isomorphism is non-canonical, i.e.~it depends on the choice of the embedding $F'_0\hookrightarrow L$ in Lemma \ref{EHlemWDtophiNmodules}.
We will spell out this isomorphism in detail in Section~\ref{subsubsec:embedding role} below, and also give explicit descriptions of the stacks involved (at least for $F'=F$).

In order to link a $(\varphi,N,{\rm Gal}(F'/F))$-module $D\in{\rm Mod}_{d,\varphi,N,F'/F}(A)$ to equivariant vector bundles, we need to introduce a filtration on $D$. 
Let us write $G={\rm Res}_{F/\mathbf{Q}_p}{\rm GL}_d$ and $G'={\rm Res}_{F'/\mathbf{Q}_p}{\rm GL}_d$ for the moment. Moreover, we write $G'_0={\rm Res}_{F'_0/\mathbf{Q}_p}{\rm GL}_d$. 

Let $T\subseteq G$ denote the Weil restriction of the diagonal torus of ${\rm GL}_d$ and let $B\subseteq G$ denote the Weil restriction of the upper triangular matrices. 
For a choice of a dominant cocharacter $\lambdau\in X_\ast(T_{\bar L})_+$ we write $L_{\lambdau}$ for the reflex field of $\lambdau$ and $P_{\lambdau}\subseteq G_{L_{\lambdau}}$ for the parabolic subgroup (containing $B_{L_{\lambdau}}$) defined by $\lambdau$. The flag variety $G_{L_{\lambdau}}/P_{\lambdau}$ then parameterizes filtrations on $L_{\lambdau}^d\otimes_{\mathbf{Q}_p}F$ of type $\lambdau$. More precisely, using the standard identification 
$$X_\ast(T)_+=\prod_{\tau:F\hookrightarrow \bar L} \mathbf{Z}^d_+$$
we write $$\lambdau=(\lambda_{\tau,1}\geq \lambda_{\tau,2}\geq\dots\geq \lambda_{\tau,d})_{\tau:F\hookrightarrow \bar L}.$$
Then an $\bar L$-valued point of $G_{L_{\lambdau}}/P_{\lambdau}$ is given by a product of descending, exhaustive and separated $\mathbf{Z}$-filtrations on 
$$\bar L^d\otimes_{\mathbf{Q}_p}F=\prod_{\tau:F\hookrightarrow\bar L} \bar L_\tau^d$$
(with $\bar L_\tau=\bar L$) such that $(\lambda_{\tau,1}\geq \lambda_{\tau,2}\geq\dots\geq \lambda_{\tau,d})$ are precisely the weights (counted with multiplicity) of the filtration on $\bar L_\tau^d$.

Using the canonical map $G\rightarrow G'$ we can also view $\lambdau$ as a cocharacter of a torus in $G'$ (more precisely again of the Weil restriction of the diagonal torus in ${\rm GL}_d$) and we obtain a similar parabolic subgroup $P'_{\lambdau}$ and a flag variety $G'_{L_{\lambdau}}/P'_{\lambdau}$.

Note that there is a forgetful map
$${\rm WD}_{d,F}\longrightarrow B{\rm GL}_d=\ast/{\rm GL}_d$$
to the stack of rank~$d$ vector bundles. And similarly there is a canonical map 
$${\rm Mod}_{d,\varphi,N,F'/F}\longrightarrow \ast/G'_0$$
forgetting the action of $\varphi,N$ and ${\rm Gal}(F'/F)$ on a projective $A\otimes_{\mathbf{Q}_p}F'_0$-module of rank~$d$. 

The scheme $G/P_{\lambdau}$ is equipped with a (diagonal) left translation action of ${\rm GL}_d$ and we can define the fiber product
$${\rm Fil}_{\lambdau}{\rm WD}_{d,F}={\rm WD}_{d,F}\times_{\ast/{\rm GL}_d} {\rm GL}_d\backslash G_{L_{\lambdau}}/P_{\lambdau}$$
that is called the stack of filtered Weil--Deligne representations (of Hodge--Tate weight $\lambdau$). Its $S$-valued points parameterize Weil--Deligne representations on vector bundles $\mathcal{E}$ over $S$ together with a filtration on $F\otimes_{\mathbf{Q}_p}\mathcal{E}$ of type $\lambdau$.

Similarly, the scheme $G'/P'_{\lambdau}$ is equipped with a (diagonal) left translation action of $G'_0$ and we can consider the closed substack 
$${\rm Fil}_{\lambdau}{\rm Mod}_{d,\varphi,N,F'/F}\subset{\rm Mod}_{d,\varphi,N,F'/F}\times_{\ast/G'_0} G'_0\backslash G'_{L_{\lambdau}}/P'_{\lambdau}$$
parameterizing (for a given object $D\in {\rm Mod}_{d,\varphi,N,F'/F}$) those filtrations (of type $\lambdau$) on $D\otimes_{F'_0}F'$ that are stable under the action of ${\rm Gal}(F'/F)$. 
\begin{rem}
We point out that ${\rm Gal}(F'/F)$-stable filtrations can only exist if the cocharacter $\lambdau$ of (our maximal torus in) $G'$ comes from (our maximal torus in) $G$ via the map $G\rightarrow G'$.
\end{rem}

The following lemma then generalizes Lemma \ref{EHlemWDtophiNmodules}. 

\begin{lem}\label{EH:lem:filtWDrep}
Assume that $L$ is large enough such that $[F':\mathbf{Q}_p]=|{\rm Hom}(F',L)|$ and fix an embedding $F'\hookrightarrow L$. Then there is a canonical isomorphism 
$${\rm Fil}_{\lambdau}{\rm Mod}_{d,\varphi,N,F'/F} \isoto {\rm Fil}_{\lambdau}{\rm WD}_{d,F'/F}$$
lying above the isomorphism ${\rm WD}_d$ of Lemma~{\em \ref{EHlemWDtophiNmodules}}.
\end{lem}

The stack ${\rm Fil}_{\lambdau}{\rm Mod}_{d,\varphi,N,F'/F} $ plays an important role in the characterization of de Rham $(\varphi,\Gamma)$-modules and via the above isomorphism we can hence link filtered Weil--Deligne representations to $(\varphi,\Gamma)$-modules, at least if the base field is large enough. More precisely we have the following theorem: 
\begin{theorem}\label{thm: embedding of de rham stack}
For fixed $\lambdau\in X_\ast(T_{\bar L})_+$ there is a canonical closed embedding
$$({\rm Fil}_{\lambdau}{\rm Mod}_{d,\varphi,N,F'/F}) ^{\rm
  an}\longrightarrow \mathfrak{X}_{d,L_{\lambdau}}.$$
Moreover, the image of this morphism coincides with the closed substack of de Rham $(\varphi,\Gamma)$-modules of Hodge--Tate--Sen weight $\lambdau$ that become semi-stable after base change to $F'$. 
\end{theorem}

The construction of a family of $(\varphi,\Gamma)$-modules of rank~$d$ over the rigid analytic stack $({\rm Fil}_{\lambdau}{\rm Mod}_{d,\varphi,N,F'/F}) ^{\rm an}$ is the generalization of the work of Berger \cite{MR2493215} to families of the objects involved. 

To make sense of the ``moreover'' part of Theorem~\ref{thm: embedding of de rham stack} we need to define the notion of a family of (potentially) semi-stable or of de Rham $(\varphi,\Gamma)$-modules over $\mathcal{R}_{F,A}$. Given a $(\varphi,\Gamma)$-module $D_A$ over $\mathcal{R}_{F,A}$ we can consider the base change 
$$D_{\rm dif}(D_A)=D_A\otimes_{\mathcal{R}_{F,A}}{F}_\infty((t))$$
which is endowed with a derivation $$\partial: D_{\rm dif}(D_A)\longrightarrow D_{\rm dif}(D_A)$$ above $t\tfrac{d}{dt}$ that is defined by the derivative of the $\Gamma$-action. Then $D_A$ is defined to be de Rham if $D_{\rm dif}(D_A)^{\partial=0}$ is a projective $A\otimes_{\mathbf{Q}_p}{F}_\infty$-module of rank~$d$ and if the canonical map
$$D_{\rm dif}(D_A)^{\partial =0}\otimes_{{ F}_\infty}{F}_\infty((t))\longrightarrow D_{\rm dif}(D_A)$$ is an isomorphism.
To prove the ``moreover'' part of Theorem~\ref{thm: embedding of
  de rham stack} we follow closely the proof of \cite[Thms.\ 5.3.1, 5.3.2]{MR2493221} in order to show that there is a Zariski-closed substack of $\mathfrak{X}_{d,L_{\underline{\lambda}}}$ consisting of the de Rham objects of Hodge--Tate weight $\underline{\lambda}$. By the $p$-adic monodromy theorem the classical points (i.e.~the $L'$-valued points for finite extensions $L'$ of $L_{\underline{\lambda}}$) coincide with the classical points of $\bigcup_{F'} ({\rm Fil}_{\lambdau}{\rm Mod}_{d,\varphi,N,F'/F}) ^{\rm
  an}$, which implies that these stacks have to coincide. 

\begin{rem}
  Recall that the $(\varphi,\Gamma)$-modules in the analytic set up are not required to be \'etale. 
  Hence the stacks of de Rham $(\varphi,\Gamma)$-modules can be described by filtered $(\varphi,N)$-modules (with descent data) as done in Theorem~\ref{thm: embedding of
  de rham stack} without putting an additional ``weak admissibility'' condition on the filtration.
  It is however possible to pass to an open substack of $({\rm Fil}_{\lambdau}{\rm Mod}_{d,\varphi,N,F'/F}) ^{\rm
  an}$ where the $(\varphi,\Gamma)$-modules associated to the filtered $(\varphi,N,{\rm Gal}_F)$-modules are \'etale, see \cite{MR4129383} and \cite{MR3103130} for that construction in the case of a trivial ${\rm Gal}_F$-action.
  This \'etale locus then agrees with the (union over appropriate $\tau$ of the) image of the generic fibers of the stacks $\cX_d^{\semis,\lambdau,\tau}$  under the map (\ref{EHgenfibEGstacktoanalytic}).
\end{rem}
\begin{rem}
We point out that it is important to fix the cocharacter $\lambdau$ in
the theorem. If $\lambdau$ is allowed to vary arbitrarily, then the
union of the corresponding substacks $({\rm Fil}_{\lambdau}{\rm
  Mod}_{d,\varphi,N,F'/F}) ^{\rm an}$ is not closed in
$\mathfrak{X}_d$ but expected to be Zariski-dense. This should be seen
as an analogue of the density of crystalline representations in local
deformation rings, see \cite{MR3037022,MR3235551, https://doi.org/10.48550/arxiv.2206.06944}. 
\end{rem}

\subsubsection{The role of the embedding $F'_0\hookrightarrow L$}
\label{subsubsec:embedding role}
In Lemma~\ref{EHlemWDtophiNmodules} above, we have identified the stacks ${\rm Mod}_{d,\varphi,N,F'/F}$ and ${\rm WD}_{d,F'/F}$ over a field $L$ containing $F'_0$ using the choice of an embedding $\sigma_0:F'_0\hookrightarrow L$. 
We elaborate a bit further about the choice of this embedding. For simplicity (and in order not to overload the notation) we limit ourselves to the case $F'=F$. 
We refer to \cite[\S 2.6]{MR3927049} and \cite[\S 3]{MR4129383} for the explicit description of the stacks that we are using in the following. 

Assume that $L$ contains $F'_0 = F_0$ (the equality holding because we have assumed
that $F'=F$); then for an $L$-algebra $A$ we have a decomposition 
$$A\otimes_{\Q_p}F_0=A\otimes_L(L\otimes_{\mathbf{Q}_p}F_0)=A\otimes_L \prod_{\varsigma:F_0\hookrightarrow L} L_\varsigma=\prod_{\varsigma:F_0\hookrightarrow L} A_\varsigma,$$
where $L_\varsigma=L$ and $A_\varsigma=A$ for all $\varsigma$. The action of ${\rm Gal}(F_0/\mathbf{Q}_p)$ of course just permutes the factors in the product. For the remainder of this section we fix an embedding  $\varsigma_0$ and rewrite $A\otimes_{\Q_p}F_0=\prod_{i=1}^f A_i$ with $A_i=A_{\varsigma_0^i}=A$. 
We deduce that we have isomorphisms
\numequation\label{EHeqn decompoofWeilres}
\begin{aligned}
({\rm Res_{F_0/\mathbf{Q}_p}}{\rm GL}_d)_L&=\prod_{i=1}^f {\rm GL}_d \\
({\rm Res_{F_0/\mathbf{Q}_p}}{\rm Mat}_{d\times d})_L&=\prod_{i=1}^f {\rm Mat}_{d\times d}.
\end{aligned}
\end{equation}
Note that ${\rm WD}_{d,F/F}=X_d/{\rm GL}_d$, where
$$X_d=\{(\varphi,N)\in{\rm GL}_d\times{\rm Mat}_{d\times d}| 
N\varphi=q\varphi N\},$$
and where ${\rm GL}_d$ acts on ${\rm GL}_d$ and ${\rm Mat}_{d\times d}$ by conjugation. 
On the other hand ${\rm Mod}_{d,\varphi,N,F/F}=Y_d/{\rm Res}_{F_0/\mathbf{Q}_p}{\rm GL}_d$, where
\begin{align*}
Y_d=\{(\varphi,N)\in {\rm Res}_{F_0/\mathbf{Q}_p}{\rm GL}_d\times {\rm Res}_{F_0/\mathbf{Q}_p}{\rm Mat}_{d\times d}\mid N\varphi=p\varphi \sigma^{\ast}N\},
\end{align*}
where $\sigma:{\rm Res}_{F_0/\mathbf{Q}_p}{\rm Mat}_{d\times d}\rightarrow {\rm Res}_{F_0/\mathbf{Q}_p}{\rm Mat}_{d\times d}$ is the map induced by $\sigma$ on $F_0$; and ${\rm Res}_{F_0/\mathbf{Q}_p}{\rm GL}_d$ acts by conjugation on ${\rm Res}_{F_0/\mathbf{Q}_p}{\rm Mat}_{d\times d}$ and acts on ${\rm Res}_{F_0/\mathbf{Q}_p}{\rm GL}_d$ by $\sigma$-conjugation, i.e.~$g\cdot(\varphi,N)=(g^{-1}\varphi \sigma(g),g^{-1}Ng)$. 

After base change to $L$ we use (\ref{EHeqn decompoofWeilres}) to rewrite a point of $Y_d$ as a tuple $((\varphi_1,\dots,\varphi_f),(N_1,\dots, N_f))$ with $\varphi_i\in{\rm GL}_d$, $N_i\in {\rm Mat}_{d\times d}$ such that 
$$N_{i+1}\varphi_i=p\varphi_iN_i,\ \text{for}\ i=1,\dots, f$$
where we set $N_{f+1}=N_1$. 
The action of $({\rm Res}_{F_0/\mathbf{Q}_p}{\rm GL}_d)_L\cong \prod_{i=1}^f {\rm GL}_d$ is then given by 
\begin{align*}
&(g_1,\dots g_f)\cdot ((\varphi_1,\dots,\varphi_f),(N_1,\dots, N_f)\\
=&((g_1^{-1}\varphi_1g_f, g_2^{-1}\varphi_2g_1,\dots, g_f^{-1}\varphi_f g_{f-1}),(g_1^{-1}N_1g_1,\dots, g_f^{-1}N_f g_f)).
\end{align*}
We use this explicit description to write down a canonical morphism $X_d\rightarrow Y_d$ (depending on the isomorphisms (\ref{EHeqn decompoofWeilres}), hence on the choice of $\varsigma_0$) by
$$ (\varphi,N)\longmapsto ((\varphi,1,\dots, 1),(N,p\varphi N\varphi^{-1},p^2\varphi N\varphi^{-1}\dots, p^{(f-1)}\varphi N\varphi^{-1})$$
which is equivariant for the map of group schemes ${\rm GL}_d\rightarrow {\rm Res}_{F_0/\mathbf{Q}_p}{\rm GL}_d$ given by
$$g\longmapsto (g, g, \dots, g)$$
and which induces the isomorphism of stacks 
$$f_{\varsigma_0}:{\rm WD}_{d,F/F}=X_d/{\rm GL}_d\isoto {\rm Mod}_{d,\varphi,N,F/F}=Y_d/{\rm Res}_{F_0/\mathbf{Q}_p}{\rm GL}_d$$
in Lemma \ref{EHlemWDtophiNmodules}.
The inverse of this isomorphism is induced by the map $Y_d\rightarrow X_d$ given by 
$$((\varphi_1,\dots, \varphi_f),(N_1,\dots, N_f))\longmapsto (\varphi_f\cdots\varphi_2\varphi_1, N_1).$$
These formulas show that the isomorphism $f_{\varsigma_0}$ %
depends on the choice of $\varsigma_0$: for a different choice $\varsigma_0'$ the induced automorphism $f_{\varsigma_0'}\circ f_{\varsigma_0}^{-1}$ is not the identity. 
Moreover, we easily see that the isomorphism $f_{\varsigma_0}$ does not descend to a subfield $L'\subset L$ not containing $F_0$: if $\tau:L\rightarrow L$ is a field automorphism of $L$, then it commutes with $f_{\varsigma_0}$ if and only if its restriction to $F_0$ is the identity. 

On the other hand we note that the map $${\rm WD}_{d,F/F}\rightarrow{\rm GL}_d/\!/{\rm GL}_d\cong \mathbf{A}^{d-1}\times\mathbf{G}_m$$
mapping $(\varphi,N)$ to the coefficients of its characteristic polynomial induces an isomorphism on global sections, see e.g.~\cite[Remark 3.4]{hellmann2020derived}. 
We deduce that $f_{\sigma_0}$ induces an isomorphism on the (invariant) global sections of the structure sheaf 
\numequation\label{EHeqn:isoonglobalsections}
\Gamma({\rm Mod}_{d,\varphi,N,F/F},\mathcal{O}_{{\rm Mod}_{d,\varphi,N,F/F}})\longrightarrow \Gamma({\rm WD}_{d,F/F},\mathcal{O}_{{\rm WD}_{d,F/F}})
\end{equation}
that is independent of the choice of $\varsigma_0$ and that, moreover, this isomorphism  descends to any subfield $L'\subseteq L$.
We now make this descent explicit in the case $L=F_0$ and $L'=\mathbf{Q}_p$. 
Then the action $\sigma_{Y_d}$ of the generator $\sigma\in {\rm Gal}(F_0/\Q_p)$ on the $\varphi$-part in $Y_d$ is given by 
$$(\varphi_1,\dots,\varphi_f)\longmapsto (\sigma^\ast \varphi_f,\sigma^\ast\varphi_1,\dots, \sigma^\ast\varphi_{f-1}),$$
whereas the action $\sigma_{X_d}$ of $\sigma$ on the $\varphi$-part in $X_d$ is clearly given by $\varphi\mapsto \sigma^\ast\varphi$. We deduce that (spelling out only the $\varphi$-part of the formulas)
$$(f_{\varsigma_0}\circ \sigma_{Y_d})(\varphi_1,\dots,\varphi_f)=(\sigma^\ast \varphi_f)^{-1}\big(\sigma_{X_d}\circ f_{\varsigma_0}((\varphi_1,\dots,\varphi_f))\big) (\sigma^\ast \varphi_f)$$
and hence $f_{\varsigma_0}$ does not descend to $\Q_p$. On the other hand the characteristic polynomials of the $\varphi$-parts of $f_{\varsigma_0}\circ \sigma_{Y_d}$ and $\sigma_{X_d}\circ f_{\varsigma_0}$ agree, and hence the isomorphism between the global sections (\ref{EHeqn:isoonglobalsections}) descends to $\Q_p$. 

In fact this observation generalizes to ${\rm WD}_{d,F'/F}$ as follows:
\begin{prop}\label{EH:propsameglobalsections}
Let $L$ be any finite extension of $\mathbf{Q}_p$ and let $F'/F$ be a finite extension. For a field extension $M$ of $L$ containing $F'_0$ and any choice of an embedding $\varsigma_0:F'_0\hookrightarrow M$ the isomorphism
$$\Gamma(({\rm Mod}_{d,\varphi,N,F'/F})_{M},\mathcal{O}_{({\rm Mod}_{d,\varphi,N,F'/F})_{M}})\longrightarrow \Gamma(({\rm WD}_{d,F'/F})_{M},\mathcal{O}_{({\rm WD}_{d,F'/F})_{M}})$$
induced by $f_{\varsigma_0}$ is independent of the choice of $\varsigma_0$ and descends to a canonical isomorphism
$$\Gamma({\rm Mod}_{d,\varphi,N,F'/F},\mathcal{O}_{{\rm Mod}_{d,\varphi,N,F'/F}})\longrightarrow \Gamma({\rm WD}_{d,F'/F},\mathcal{O}_{{\rm WD}_{d,F'/F}}).$$
\end{prop}

\subsubsection{Invariant functions and the Bernstein center}\label{sec:EHBernsteincenters}
Similarly to the definition of $\mathcal{X}_d^{\rm ss,\underline{\lambda},\tau}$ in section \ref{subsec:  overview of EG book} for a given inertial type $\tau$ we will define stacks of de Rham objects with Hodge--Tate weight $\underline{\lambda}$ and fixed inertial type $\tau$.  
Recall that the local Langlands correspondence for ${\rm GL}_d(F)$ gives rise to a bijection $\tau\mapsto\Omega_\tau$ between the $d$-dimensional inertial types $\tau$ (defined over $L$) of $I_F$ and the Bernstein blocks $\Omega$ (defined over $L$) of the category ${\rm Rep}^{\rm sm}{\rm GL}_d(F)$ (we view our representations as representations on $\bar L$-vector spaces for an algebraic closure $\bar L$ of $L$).
More precisely, if $(r,N)$ is the Weil--Deligne representation associated to an irreducible smooth representation $\pi$ of ${\rm GL}_d(F)$, then $r|_{I_F}\cong \tau$ if and only if $\pi$ lies in the Bernstein block~$\Omega_\tau$. 
(This is a form of the {\em inertial local Langlands correspondence}, already remarked on
in Section~\ref{subsubsec:TWK patching} in the case of $\GL_2(\Q_p)$.  The relationship between 
this bijection $\tau \leftrightarrow \Omega_{\tau}$ and the formulation
of inertial local Langlands in terms of $K$-types is 
recalled in Section~\ref{subsubsec:inertial LL} below.)

Moreover, note that over an algebraic closure $\bar L$ of $L$ there is a bijection $\tau\mapsto {\rm WD}_{d,F,\tau}$ between the $d$-dimensional inertial types $\tau$ (defined over $\bar L$) and the connected components of the stack of Weil--Deligne representations, such that a Weil--Deligne representation $(r,N)$ lies in ${\rm WD}_{d,F,\tau}$ if and only if $r|_{I_F}\cong\tau$, see \cite[\S 5.2]{benzvi2020coherent}.
As we restrict ourselves to inertial types that become trivial on $I_{F'}$ we may choose $L$ to be large enough such that all inertial types $\tau$ over $\bar L$ (with $\tau|_{I_{F'}}$ being trivial) are defined over $L$ and such that all the geometric connected components ${\rm WD}_{d,F,\tau}$ are defined over $L$, and (by slight abuse of notation) we will still write ${\rm WD}_{d,F,\tau}$ for this stack on $L$-schemes.

The following proposition (which is a consequence of \cite[Theorem 5.13]{benzvi2020coherent}, see \cite{MR3867634} for an integral version) relates the Bernstein center $\mathcal{Z}_{\Omega_\tau}$ of $\Omega_\tau$ (that we view as a block of the category of smooth representations on $L$-vector spaces, so in particular $\mathcal{Z}_{\Omega_\tau}$ is an $L$-algebra) to the invariant functions on ${\rm WD}_{d,F,\tau}$.
\begin{prop}\label{EHBernsteiniso}
Let $\tau$ be an inertial type defined over $L$. Then there is a natural isomorphism
$$\mathcal{Z}_{\Omega_\tau}\isoto \Gamma({\rm WD}_{d,F,\tau},\cO_{{\rm WD}_{d,F,\tau}}).$$
\end{prop}

\begin{rem}\label{EHremBernsteiniso}
Note that there are two possibilities to normalize the isomorphism
of Proposition~\ref{EHBernsteiniso}, corresponding to two possible  normalizations of the local Langlands correspondence. In the usual (unitary) normalization the central character of a smooth representation $\pi$ matches the determinant of the corresponding Weil--Deligne representation $(\rho,N)$. 
Sometimes it is more convenient to choose a normalization such that $(\rho,N)$ corresponds to $\pi\otimes |{\rm det}|^{-\tfrac{d-1}{2}}$. This normalization is better suitable for local-global compatibility, and has the advantage that it does not depend on a choice of $q^{1/2}$ (and hence has better properties when considered over non-algebraically closed fields, compare \cite{MR2359853}).
In these notes we will from now on always use the \emph{non-unitary} normalization.
\end{rem}

We note that, given $\tau$, the isomorphism of Lemma \ref{EHlemWDtophiNmodules} allows us to define a connected component ${\rm Mod}_{d,\varphi,N,\tau}$ of ${\rm Mod}_{d,\varphi,N,F'/F}$ corresponding to ${\rm WD}_{d,F,\tau}$. We note that this connected component does not depend on the choice of the embedding $\varsigma_0$ in Lemma \ref{EHlemWDtophiNmodules}. We further use this to define the obvious connected component ${\rm Fil}_{\lambdau}{\rm Mod}_{d,\varphi,N,\tau}$ of ${\rm Fil}_{\lambdau}{\rm Mod}_{d,\varphi,N,F'/F}$.

For a fixed inertial type $\tau$ we can consider the rigid analytic generic fiber $(\cX_d^{\semis,\lambdau,\tau})_\eta^{\rm rig}$ of the formal stack $\cX_d^{\semis,\lambdau,\tau}$ defined in Section \ref{sec:moduli-stacks-phi-gamma}. 
By the very construction of the stacks $\cX_d^{\semis,\lambdau,\tau}$
this generic fiber admits a map %
\numequation\label{eqn: semistable period map}\pi_d^{\semis,\lambdau,\tau}:(\cX_d^{\semis,\lambdau,\tau})_\eta^{\rm rig}\longrightarrow ({\rm Fil}_{\lambdau}{\rm Mod}_{d,\varphi,N,\tau}) ^{\rm an}\subset \mathfrak{X}_d,\end{equation}
that is a generalization of the period map \cite[(5.37)]{MR2562795} of Pappas--Rapoport in the case $\tau=1$ and $\lambdau$ minuscule. Moreover, it fits into a commutative diagram 
$$
\begin{xy}
\xymatrix{
(\cX_d^{\semis,\lambdau,\tau})_\eta^{\rm rig} \ar[r]\ar[d]_{\pi_d^{\semis,\lambdau,\tau}} & \cX_{d,\eta}^{\rm rig} \ar[d]^{\pi_d} \\
({\rm Fil}_{\lambdau}{\rm Mod}_{d,\varphi,N,\tau}) ^{\rm an} \ar[r] & \mathfrak{X}_d
}
\end{xy}
$$
with the map $\pi_d$ from \eqref{EHgenfibEGstacktoanalytic}.
In particular the isomorphism 
\numequation\label{EH:eqnBersteiniso}
\mathcal{Z}_{\Omega_\tau}\isoto \Gamma({\rm Mod}_{d,\varphi,N,\tau},\mathcal{O}_{{\rm Mod}_{d,\varphi,N,\tau}})
\end{equation}
obtained by combining Propositions \ref{EH:propsameglobalsections} and \ref{EHBernsteiniso},
induces a map
\numequation\label{EH:eqn Bernstein map to functions on ss stack}\mathcal{Z}_{\Omega_\tau}\longrightarrow \Gamma((\cX_d^{\semis,\lambdau,\tau})_\eta^{\rm rig},\mathcal{O}_{(\cX_d^{\semis,\lambdau,\tau})_\eta^{\rm rig}})\end{equation}
that is defined over any chosen base field $L$ over which $\tau$ is defined (not necessarily containing the field $F'_0$ over which we can identify the stack of $(\varphi,N,{\rm Gal}(F'/F))$-modules with a substack of the stack of Weil--Deligne representations). 
\begin{rem}
After composing with the map to the functions on the generic fiber of a potentially semi-stable deformation ring this map coincides with the morphism of \cite[Theorem 4.1, Prop. 4.2]{Gpatch}.
In fact the result of \cite{Gpatch} is stronger: it asserts that the map takes values in the ring of bounded functions $$\Gamma(\cX_d^{\semis,\lambdau,\tau},\mathcal{O}_{\cX_d^{\semis,\lambdau,\tau}})[\tfrac{1}{p}]\subset \Gamma((\cX_d^{\semis,\lambdau,\tau})_\eta^{\rm rig},\mathcal{O}_{(\cX_d^{\semis,\lambdau,\tau})_\eta^{\rm rig}}).$$
\end{rem}

\subsection{Drinfeld-style compactifications}
\label{sec:drinf-style-comp}
We discuss variants of the stack $\mathfrak{X}_d$ parameterizing flags on a family of equivariant vector bundles of rank~$d$. These stacks are naturally associated with a choice of a parabolic subgroup $P\subseteq G={\rm GL}_d$. 
Given a parabolic subgroup $P$ with Levi quotient $M\cong {\rm GL}_{d_1}\times\dots\times{\rm GL}_{d_r}$ we define $\mathfrak{X}_P$ to be the stack of $\Gal_F$-equivariant $P$-bundles on the Fargues--Fontaine curve $X_{\bar F}$, which comes with two projections
\[\begin{xy}
\xymatrix{
& \mathfrak{X}_P\ar[dl]_{\beta_P}\ar[dr]\ar[dr]^{\alpha_P} & \\
\mathfrak{X}_d && \mathfrak{X}_M}
\end{xy}\]
induced by the maps $P\hookrightarrow {\rm GL}_d$ and
$P\twoheadrightarrow M$, where
$\mathfrak{X}_{M}\cong\mathfrak{X}_{d_1}\times\dots\times\mathfrak{X}_{d_r}$
denotes the stack of $\Gal_F$-equivariant $M$-bundles on $X_{\bar
  F}$. It will be convenient to write
$\mathfrak{X}_{G}=\mathfrak{X}_d$ in the following.

In the special case of $P=B$ the $(\varphi,\Gamma)$-modules in the image of the morphism $\mathfrak{X}_B\rightarrow \mathfrak{X}_d$ are, by definition, precisely the $(\varphi,\Gamma)$-modules that are trianguline (in the sense of Colmez, see \cite{MR2493219,MR3444228} for the definition and discussion of two dimensional trianguline representations and their role in the $p$-adic Langlands correspondence).
Beside the fact that trianguline $(\varphi,\Gamma)$-modules are very natural objects to consider (as a triangulation is nothing but a $B$-structure), their definition was inspired by the fact that the Galois representations associated to overconvergent $p$-adic modular forms of finite slope (i.e.~those $p$-adic modular forms whose eigensystems are parameterized by the Coleman--Mazur eigencurve) are trianguline at $p$, as discovered by Kisin \cite{MR1992017} (before the name ``trianguline'' was invented and using an equivalent description).  This was later extended to more general eigenvarieties (in particular to eigenvarieties for unitary groups, see e.g. \cite{MR2656025}). We will return to this point in Section \ref{sec:eigenvar} when we discuss the relation of eigenvarieties and $p$-adic automorphic forms with sheaves on $\mathfrak{X}_d$ defined using $\mathfrak{X}_B$ and its compactification.
\begin{rem}
There are of course variants, where ${\rm GL}_d$ is replaced by an arbitrary reductive group, but for simplicity of the exposition we limit ourselves to the case of ${\rm GL}_d$ here. We define the stacks in a way such that it is clear what the definitions will look like in general.
\end{rem}
Following Drinfeld, Braverman and Gaitsgory \cite{MR1933587} define two types of so-called Drinfeld compactifications of the stack of $P$-bundles on an algebraic curve, and we will define their analogues $\overline{\mathfrak{X}}_P$ and $\widetilde{\mathfrak{X}}_P$ in the context of $\Gal_F$-equivariant vector bundles or, equivalently, in the context of $(\varphi,\Gamma)$-modules, below. 
Both of these compactifications have their advantages and their disadvantages: the compactification $\overline{\mathfrak{X}}_P$ has the advantage that it is functorial with respect to inclusions $P'\subseteq P$ of parabolic subgroups, but has the disadvantage that the map $\alpha_P:\mathfrak{X}_P\rightarrow \mathfrak{X}_M$ does not extend to $\overline{\mathfrak{X}}_P$. Conversely the compactification $\widetilde{\mathfrak{X}}_P$ has the advantage that $\alpha_P$ extends to $\widetilde{\alpha}_P:\widetilde{\mathfrak{X}}_P\rightarrow \mathfrak{X}_M$, but $P\mapsto \widetilde{\mathfrak{X}}_P$ is not functorial with respect to inclusions $P'\subseteq P$ of parabolic subgroups. 

Before giving the formal definitions let briefly explain the motivation in the case $P=B$. 
A $B$-structure on a $(\varphi,\Gamma)$-module $D_A$ is a complete $(\varphi,\Gamma)$-stable flag $$0={\rm Fil}_0\subseteq {\rm Fil}_1\subseteq \dots \subseteq {\rm Fil}_d=D_A$$ such that ${\rm Fil}_i/{\rm Fil}_{i-1}$ is a projective $\mathcal{R}_{F,A}$-module of rank~$1$. Using Pl\"ucker coordinates the datum of such a filtration is the same as the datum of the $(\varphi,\Gamma)$-stable rank~$1$ subobjects 
$$\mathcal{L}_i=\bigwedge\nolimits^i {\rm Fil}_i\subseteq \bigwedge\nolimits^i D_A$$ such that the quotient $(\bigwedge^i D_A)/\mathcal{L}_i$ is a projective $\mathcal{R}_{F,A}$-module. 
The lines $\mathcal{L}_i$ satisfy some relations (called Pl\"ucker
relations) that we do not spell out explicitly, but we refer to the
more abstract description below that allows a more uniform description
of these relations (but has to pay the price that we consider the
evaluation of $D_A$ on all algebraic representations of ${\rm GL}_d$
instead of just the exterior powers $\bigwedge^i D_A$, compare also
Remark \ref{EHremexteriorpowers} below).

The following example illustrates that there are families $D_A$ over rigid spaces ${\rm Sp}\,A$ that admit a $B$-structure on a Zariski-open and dense subspace, but this $B$-structure does not extend to all of ${\rm Sp}\,A$. 
Hence we aim to somehow compactify the situation. The idea is to still parameterize lines $\mathcal{L}_i\subseteq \bigwedge^i D_A$ that satisfy the Pl\"ucker relations, but to drop the condition that $(\bigwedge^i D_A)/\mathcal{L}_i$ is projective over $\mathcal{R}_{F,A}$. This obviously fixes the issue in the following example. 

\begin{example}
Let $F=\mathbf{Q}_p$ and consider the following family of filtered $\varphi$-modules over the projective line $\mathbf{P}^1$.
Let $$M=\mathcal{O}_{\mathbf{P}^1}e_1\oplus\mathcal{O}_{\mathbf{P}^1}e_2\ \text{and}\ \varphi_M={\rm diag}(\varphi_1,\varphi_2)$$
for fixed $\varphi_1\neq\varphi_2\in\mathbf{Q}_p^\times$. 
Moreover we consider the $\mathbf{Z}$-filtration ${\rm Fil}^\bullet$ on $M$ given by 
$${\rm Fil}^i M=\begin{cases}M & i\leq -1 \\ \mathscr{L} & i=0 \\ 0 & i\geq 1\end{cases}$$
where $\mathscr{L}\subset M$ is the universal line over $\mathbf{P}^1$. 
We consider the submodule $M'=\mathcal{O}_{\mathbf{P}^1}e_1\subset M$. Then it is not possible to define a filtration ${\rm Fil}^i M'$ on $M'$ (in the sense that ${\rm Fil}^i M'\subset M'$ is locally on $\mathbf{P}^1$ a direct summand) such that for all $x\in \mathbf{P}^1$ the subobject ${\rm Fil}^iM'\otimes k(x)\subset M'\otimes k(x)$ is the intersection of ${\rm Fil}^iM\otimes k(x)$ with $M'\otimes k(x)$. 
If we insist that ${\rm Fil}^iM'\subset M'$ is (locally) a direct summand, then we can only arrange this on $\mathbf{A}^1\subset\mathbf{P}^1$ by setting 
$${\rm Fil}^i M'=\begin{cases}M' & i\leq -1 \\ 0 & i\geq 0 \end{cases}.$$
However at $x=\infty$ the specialization ${\rm Fil}^iM'\otimes k(x)$ is a proper subset of the intersection  ${\rm Fil}^iM\otimes k(x)\cap M'\otimes k(x)$.

Using the construction of $(\varphi,\Gamma)$-modules associated to filtered $\varphi$-modules (see Theorem \ref{thm: embedding of de rham stack} and the discussion below it), this translates to the following situation.
Let $$\mathcal{R}_{F,\mathbf{P}^1}^{\oplus 2}\subset D\subset
t^{-1}\mathcal{R}_{F,\mathbf{P}^1}^{\oplus 2}$$ be the $(\varphi,\Gamma)$-module over $\mathbf{P}^1$ associated to $M$ and let $\mathcal{L}$ be the $(\varphi,\Gamma)$-stable line defined by the sub-filtered $\varphi$-module $M'\subset M$. 
Then $\mathcal{L}\subset D$ is a $B$-structure on $D|_{\mathbf{A}^1}$, but after specializing to $\infty\in\mathbf{P}^1$ the cokernel of this inclusion has $t$-torsion.  

\end{example}

In order to define these compactifications in general we introduce the following notation, compare \cite{MR1933587}.
We fix a Borel subgroup $B\subseteq G={\rm GL}_d$ and write $T=B/U$, where $U$ is the unipotent radical of $B$.
As usual we write $\langle \cdot,\cdot\rangle$ for the canonical pairing between $X_\ast(T)$ and $X^\ast(T)$. 
Given a parabolic subgroup $B\subseteq P\subseteq G$ with unipotent radical $U_P$ and Levi quotient $M=P/U_P$, we write $X^\ast(T)_M^+\supset X^\ast(T)_G^+$ for the semi-groups of dominant weights for $M$ respectively for~$G$. 
The set of simple coroots $\{\alpha_i\in X_\ast(T), i\in\mathcal{I}\}$ of $G$ is indexed by the set $\mathcal{I}=\mathcal{I}_G=\{1,\dots,d-1\}$ of vertices of the Dynkin diagram of $G$. Given $M$ as above, the Dynkin diagram for $M$ gives rise to a subset $\mathcal{I}_M\subseteq \mathcal{I}$ and we set
\[X^\ast(T)_{G,P}=\{\check \lambda\in X^\ast(T)\mid \langle \alpha_i,\check \lambda\rangle=0\ \text{for}\ i\in\mathcal{I}_M\}.\]
Moreover, we write $S_M=M/[M,M]$ for the maximal torus quotient of $M$.

Given a (family of) $\Gal_F$-equivariant vector bundle(s)
$\mathcal{E}$ on $X_{\bar F}$ (or a family of
$(\varphi,\Gamma)$-modules $D$ over $\mathcal{R}_F$) and an algebraic
representation $V$ of $G$ we write $V(\mathcal{E})$ (respectively
$V(D)$) for the $\Gal_F$-equivariant ${\rm GL}(V)$-bundle
(respectively a $(\varphi,\Gamma)$-module with ${\rm
  GL}(V)$-structure) given by the push-forward of $\mathcal{E}$ along
$G\rightarrow {\rm GL}(V)$. %
We view $V(\mathcal{E})$ again as a vector bundle instead of a ${\rm GL}(V)$-bundle, and similarly we view $V(D)$ again as a $(\varphi,\Gamma)$-module.
For a character $\check \lambda\in X^\ast(T)_G^+$ we write $V^{\check \lambda}$ for the irreducible $G$-representation of highest weight $\check\lambda$ and $L^{\check \lambda}$ for the canonical representation $\check\lambda:T\rightarrow \mathbf{G}_m$ of $T$.

The following definition of $\overline{\mathfrak{X}}_P^{\rm naive}$ as well as the definition of $\widetilde{\mathfrak{X}}_P^{\rm naive}$ below are the analogues of the corresponding compactifications of ${\rm Bun}_P$ defined in \cite{MR1933587}.
\begin{defn}
The stack $\overline{\mathfrak{X}}_P^{\rm naive}$ is the category fibered in groupoids over ${\rm Rig}_L$ parameterizing triples
\[\overline{\mathfrak{X}}_P^{\rm naive}({\rm Sp}\,A)=\left\{\begin{array}{*{20}c} \mathcal{E}_G\in\mathfrak{X}_G({\rm Sp}\,A),\mathcal{E}_{S_M}\in\mathfrak{X}_{S_M}({\rm Sp}\,A)\\ (\kappa_{\check\lambda}:L^{\check\lambda}(\mathcal{E}_{S_M})\rightarrow V^{\check \lambda}(\mathcal{E}_G))_{\check\lambda\in X^\ast(T)_G^+\cap X^\ast(T)_{G,P}}\end{array}\right\}\]
such that the restriction of $\kappa_{\check\lambda}$ to fibers over ${\rm Sp}\,A$ is injective, the cokernels of the various $\kappa_{\check\lambda}$ are flat over ${\rm Sp}\,A$ and such that the following \emph{Pl\"ucker relations} are satisfied:
\begin{enumerate}
\item[-] for $\check\lambda=0$ the map $\kappa_0$ is the identity.
\item[-] for $\check\lambda,\check\mu\in X^\ast(T)_G^+\cap X^\ast(T)_{G,P}$ the composition
\begin{align*}
L^{\check\lambda}(\mathcal{E}_{S_M})\otimes L^{\check\mu}(\mathcal{E}_{S_M})=L^{\check\lambda+\check\mu}(\mathcal{E}_{S_M})\xrightarrow{\kappa_{\check\lambda+\check \mu}}&V^{\check\lambda+\check\mu}(\mathcal{E}_G)\\ \rightarrow (V^{\check\lambda}\otimes V^{\check\mu})(\mathcal{E}_G)=&V^{\check\lambda}(\mathcal{E}_G)\otimes V^{\check\mu}(\mathcal{E}_G)
\end{align*}
coincides with the map $\kappa_{\check\lambda}\otimes\kappa_{\check\mu}$.
\end{enumerate}
\end{defn} 
The stack $\mathfrak{X}_P$ of equivariant $P$-bundles on $X_{\bar F}$ can be identified with the open substack of $\overline{\mathfrak{X}}_{P}^{\rm naive}$, where the cokernels of the maps $\kappa_{\check\lambda}$ are flat over ${\rm Sp}\,A\times X_{\bar F}$ (not just over ${\rm Sp}\,A$) for all $\check\lambda$. 
It is obvious from the definition that the maps $\mathfrak{X}_P\rightarrow\mathfrak{X}_G$ and $\mathfrak{X}_P\rightarrow\mathfrak{X}_{M}\rightarrow\mathfrak{X}_{S_M}$ extend to morphisms
\begin{align*}
\overline{\beta}_P:\overline{\mathfrak{X}}_P^{\rm naive}&\longrightarrow \mathfrak{X}_G\\
\overline{\alpha}_P:\overline{\mathfrak{X}}_P^{\rm naive}&\longrightarrow \mathfrak{X}_{S_M},
\end{align*}
but $\overline{\alpha}$ does not factor though $\mathfrak{X}_M$.
Moreover $P\mapsto \overline{\mathfrak{X}}_{P}^{\rm naive}$ is functorial with respect to inclusions $P'\subseteq P$ of parabolic subgroups (compatible with the projection maps $\overline{\alpha}$ and $\overline{\beta}$).
\begin{defn}
The stack $\widetilde{\mathfrak{X}}_P^{\rm naive}$ is the category fibered in groupoids over ${\rm Rig}_L$ parameterizing triples
\[\widetilde{\mathfrak{X}}_P^{\rm naive}({\rm Sp}\,A)=\left\{\begin{array}{*{20}c} \mathcal{E}_G\in\mathfrak{X}_G({\rm Sp}\,A),\mathcal{E}_{M}\in\mathfrak{X}_{M}({\rm Sp}\,A)\\ (\kappa_{V}:V^{U_P}(\mathcal{E}_{M})\rightarrow V(\mathcal{E}_G))_V\end{array}\right\}\]
where the maps $\kappa_V$ are indexed by the finite dimensional $G$-representations,
such that the restrictions of $\kappa_{V}$ to the fibers over ${\rm Sp}\,A$ are injective, the cokernels of the various $\kappa_{V}$ are flat over ${\rm Sp}\,A$, and such that the following version of the Pl\"ucker relations are satisfied:
\begin{enumerate}
\item[-] if $V$ is the trivial representation, then the map $\kappa_V$ is the identity.
\item[-] for any  $G$-representations $V_1$ and $V_2$ the following diagram commutes:
\[
\begin{xy}
\xymatrix{
(V_1^{U_P}\otimes V_2^{U_P})(\mathcal{E}_M) \ar[d]\ar[r]^{\kappa_{V_1}\otimes\kappa_{V_2}} & (V_1\otimes V_2)(\mathcal{E}_G)\ar[d]^= \\
(V_1\otimes V_2)^{U_P}(\mathcal{E}_M)\ar[r]_{\kappa_{V_1\otimes V_2}} & (V_1\otimes V_2)(\mathcal{E}_G)
}
\end{xy}
\]
\item[-] for a morphism of $G$-representations $V_1\rightarrow V_2$ the following diagram commutes:
\[\begin{xy}
\xymatrix{
V_1^{U_P}(\mathcal{E}_M) \ar[d]\ar[r]^{\kappa_{V_1}} & V_1(\mathcal{E}_G)\ar[d] \\
V_2^{U_P}(\mathcal{E}_M) \ar[r]^{\kappa_{V_2}} & V_2(\mathcal{E}_G).
}
\end{xy}\]
\end{enumerate}
\end{defn}
Again $\mathfrak{X}_P$ can be identified with the open substack of $\widetilde{\mathfrak{X}}_P^{\rm naive}$ where the cokernels of the maps $\kappa_V$ are flat (over ${\rm Sp}\,A\times X_{\bar F}$, not just $A$-flat). Again, it is obvious from the definition that the maps $\mathfrak{X}_P\rightarrow\mathfrak{X}_G$ and $\mathfrak{X}_P\rightarrow\mathfrak{X}_{M}$ extend to morphisms
\begin{align*}
\widetilde{\beta}_P:\widetilde{\mathfrak{X}}_P^{\rm naive}&\longrightarrow \mathfrak{X}_G\\
\widetilde{\alpha}_P:\widetilde{\mathfrak{X}}_P^{\rm naive}&\longrightarrow \mathfrak{X}_{M}.
\end{align*}

Moreover, there is a forgetful map $\widetilde{\mathfrak{X}}_P^{\rm naive}\rightarrow \overline{\mathfrak{X}}_P^{\rm naive}$, inducing the identity on $\mathfrak{X}_P$, that is given by only remembering the maps $\kappa_V$ for $V=V^{\check\lambda}$ with $\check\lambda\in X^\ast(T)_G^+\cap X^\ast(T)_{G,P}$. In this case the $M$-representation $(V^{\check\lambda})^{U_P}$ factors through $S_M$.  
\begin{rem}\label{EHremexteriorpowers}
The stacks $\overline{\mathfrak{X}}_P^{\rm naive}$ and $\widetilde{\mathfrak{X}}_P^{\rm naive}$ coincide in the case $P=B$, but are distinct otherwise.
In our setting with $G={\rm GL}_d$ and in the case $P=B$ it is enough to remember the maps $\kappa_{\check\lambda}$ for the fundamental weights of~${\rm GL}_d$. That is, an ${\rm Sp}\,A$-valued point of $\overline{\mathfrak{X}}_B$ can be described by tuples $(\mathcal{E},\mathcal{L}_1,\dots,\mathcal{L}_d,\nu_1,\dots,\nu_d)$, where $\mathcal{E}$ is an equivariant vector bundle of rank~$d$, the $\mathcal{L}_i$ are equivariant line bundles and the $\nu_i$ are morphisms 
\[\nu_i:\mathcal{L}_i\longrightarrow\bigwedge\nolimits^i\mathcal{E}\]
that are fiberwise (over ${\rm Sp}\,A$) injective with $A$-flat cokernel and such that the $\nu_i$ satisfy some Pl\"ucker relations (that we do not spell out explicitly).
\end{rem}

It turns out the the equivariance condition on vector bundles adds a new phenomenon which makes the above definitions too naive (hence the name): the stack $\mathfrak{X}_P$ is not dense in $\overline{\mathfrak{X}}_P^{\rm naive}$ respectively in $\widetilde{\mathfrak{X}}_P$. 
Roughly the reason for this behavior is the existence of Hodge--Tate--Sen weights. The following example makes this phenomenon more explicit:
\begin{example}
Let $d=2$ and $D\in\mathfrak{X}_2(L)$. Then a $B$-structure on $D$ is given by a subobject $\mathcal{L}=\mathcal{R}_{F,L}(\delta_1)\subset D$ such that the quotient $D/\mathcal{R}_{F,L}(\delta_1)$ is a $(\varphi,\Gamma)$-module of rank~$1$, say $\mathcal{R}_{F,L}(\delta_2)$. It follows that ${\rm wt}_\tau(D)=\{{\rm wt}_\tau(\delta_1),{\rm wt}_\tau(\delta_2)\}$. As the labeled Hodge--Tate weights of both $D$ and the subobject $\mathcal{L}$ vary continuously on $\overline{\mathfrak{X}}_{B}^{\rm naive}$, it follows that on the closure of $\mathfrak{X}_B$ the weight of the subobject $\mathcal{L}$ must be one of the weights of~$D$. On the other hand, by definition of $\overline{\mathfrak{X}}_B$ all the objects $(D,t^m\mathcal{L})$ for $m\geq 0$ lie in $\overline{\mathfrak{X}}_B$.
\end{example}
In order to define good compactifications of $\mathfrak{X}_P$ we would like to take its closure inside the naive compactifications in a reasonable way. In the case $P=B$ this is the content of the following theorem that will be proven in \cite{HellmannHernandezSchraen}.
\begin{theorem}
\hfill\\(i) The stacks $\mathfrak{X}_B$ and $\overline{\mathfrak{X}}_B^{\rm naive}=\widetilde{\mathfrak{X}}_B^{\rm naive}$ are rigid analytic Artin stacks.\\
(ii) There is a well-defined ``scheme-theoretic" image $\overline{\mathfrak{X}}_B$ of $\mathfrak{X}_B$ in $\overline{\mathfrak{X}}_B^{\rm naive}$.
\end{theorem}
\begin{remark}\label{EHrem defofXPnotcomplete}
Of course we expect that there is a similar statement for $\mathfrak{X}_P$ that yields a definition of $\overline{\mathfrak{X}}_P$ in a similar fashion.
However, the construction of smooth surjections from a rigid analytic space to $\mathfrak{X}_B$ (respectively $\overline{\mathfrak{X}}_B^{\rm naive}, \overline{\mathfrak{X}}_B$) crucially uses the fact that $\mathfrak{X}_T$ is an Artin stack. For a more general parabolic subgroup $P$ we can prove the same result for $\mathfrak{X}_P$ and its compactification $\overline{\mathfrak{X}}_P$, if we assume that $\mathfrak{X}_M$ is an Artin stack, i.e.\ if we can prove Conjecture \ref{EHConj:Artinstack} for $M$. For this reason there is no complete definition or characterization of the stack $\overline{\mathfrak{X}}_P$. We will still use the symbol $\overline{\mathfrak{X}}_P$ in the following, for the scheme-theoretic image of $\mathfrak{X}_P$ in $\overline{\mathfrak{X}}_P^{\rm naive}$ that is conjecturally well-defined.
Similarly one expects to obtain a compactification $\widetilde{\mathfrak{X}}_P$ by taking the scheme-theoretic image of ${\mathfrak{X}}_P$ in $\widetilde{\mathfrak{X}}_P^{\rm naive}$.
The maps $\overline{\alpha}_P, \overline{\beta}_P,\widetilde{\alpha}_P,\widetilde{\beta}_P$ then conjecturally restrict to $\overline{\mathfrak{X}}_P$, respectively $\widetilde{\mathfrak{X}}_P$ and will be denoted by the same letters in the following.
\end{remark}

\begin{proof}[Sketch of proof.]
We sketch the proof of the fact that $\mathfrak{X}_B$ is an Artin stack in the $2$-dimensional case. 
The proof can be divided into two steps: in a first step we construct
the restriction $\mathfrak{X}_B^{\rm wreg}$ of $\mathfrak{X}_B$ to a
\emph{weakly regular} subset $\mathfrak{X}_T^{\rm wreg}\subset
\mathfrak{X}_T$. This weakly regular subset ensures that all ${\rm
  Ext}^2$-terms vanish which makes it possible to construct explicit
charts. In a second step we use the fact that we can twist away the
second cohomology $H_{\varphi,\Gamma}^2$ and use the Beauville--Laszlo gluing lemma to reduce the general case to the first step.

\emph{Step 1}: Let us write $\mathfrak{X}_T^{\rm wreg}\subset \mathfrak{X}_T$ for the open subset of characters $(\delta_1,\delta_2)$ such that $H^2_{\varphi,\Gamma}(\mathcal{R}(\delta_1/\delta_2))=0$. If the base field $L$ contains all Galois conjugates of $F$, then $\mathfrak{X}_T^{\rm reg}$ may be described as the set of all $(\delta_1,\delta_2)$ such that 
$$\delta_1/\delta_2\notin \{\varepsilon z^\lambda\mid \lambda\in \prod_{\tau:F\hookrightarrow L} \mathbf{Z}_+\}.$$
Let $U={\rm Sp}(A)$ be an affinoid open subset of $\mathfrak{X}_T^{\rm wreg}$. We will write $(\delta_1,\delta_2)$ for the restriction of the universal characters to $U$, and use the same notation for the pullback of $\delta_i$ along maps from rigid spaces to $U$.
As $H^2_{\varphi,\Gamma}(\mathcal{R}_{A,F}(\delta_1/\delta_2))$ vanishes pointwise on $U$, there is a quasi-isomorphism
$$[\mathcal{E}^0\xrightarrow{d}\mathcal{E}^1]\xrightarrow{f^\bullet} (C^\bullet_A,\partial^\bullet)$$
from a perfect complex concentrated in degrees $[0,1]$ to the Herr-complex $C^\bullet_A$ from Remark \ref{EHrem:propertiesphiGammacohom} (i) computing $H^\bullet_{\varphi,\Gamma}(\mathcal{R}_{A,F}(\delta_1/\delta_2))$. Moreover, this quasi-isomorphism stays a quasi-isomorphism after any pullback.
The map $f^1:\mathcal{E}^1\rightarrow \ker\partial^1\subset C_A^1$ gives rise to a universal extension 
$$0\longrightarrow \mathcal{R}_{Y,F}(\delta_1)\longrightarrow D_Y\longrightarrow \mathcal{R}_{Y,F}(\delta_2)\longrightarrow 0$$
over the geometric vector bundle $Y\rightarrow U$ associated to $\mathcal{E}^1$. This extension induces a morphism
$$g:Y\longrightarrow \mathfrak{X}_B|_{U}.$$
We claim that $g$ is a smooth surjection. Indeed, $g$ is surjective by
construction, and relatively representable as the diagonal of
$\mathfrak{X}_B$ is representable. It is left to show that $g$
satisfies the infinitesimal lifting criterion. Let $B
\twoheadrightarrow \bar B$ be a surjection of affinoid $A$-algebras and consider a commutative diagram
\[
\begin{xy}
\xymatrix{
{\rm Sp}(\bar B)\ar[r]^{\bar h}\ar[d] & Y\ar[d]^g\\
{\rm Sp}(B)\ar[r]^{h} & \mathfrak{X}_B.
}
\end{xy}
\]
We need to find a lift ${\rm Sp}(B)\rightarrow Y$ making the diagram commute.
The morphism~$\bar h$ defines an element $s_{\bar B}\in\Gamma({\rm Sp}(\bar B),\mathcal{E}^1_{\bar B})$ (where we write $\mathcal{E}^1_{\bar B}$ for the pullback of $\mathcal{E}^1$ to ${\rm Sp}(\bar B)$, and similarly for the other objects involved) and its image $f^1(s_{\bar B})$ defines an extension 
$$0\longrightarrow \mathcal{R}_{\bar B,F}(\delta_1)\longrightarrow D_{\bar B}\longrightarrow \mathcal{R}_{\bar B,F}(\delta_2)\longrightarrow 0.$$ 
The morphism $h$ gives a lift of the extension 
$D_{\bar B}$ to an extension $D_B$ of $\mathcal{R}_{B,F}(\delta_2)$ by $\mathcal{R}_{B,F}(\delta_1)$ and (after fixing appropriate bases) we may view $D_B$ as a cocycle $c^1_B\in C^1_B$ lifting $f^1(s_{\bar B})\in C^1_{\bar B}$ . 
To find a lift ${\rm Sp}(B)\rightarrow Y$ we need to find an element $s_B\in \Gamma({\rm Sp}(B),\mathcal{E}^1_{B})$ such that the extension defined by $f^1(s_B)$ is isomorphic to the extension $D_B$ via an isomorphism lifting the canonical identification of $D_{\bar B}$ with the extension defined by $f^1(s_{\bar B})$. In other words we need to find $s_B$ and $c^0_B\in C^0_B$ reducing to $0$ in $C^0_{\bar B}$ such that $c^1_B-f^1(s_B)=\partial^0(c^0_B)$.

As $\mathcal{E}^1_B\rightarrow H^1(\mathcal{E}_B^\bullet)=H^1_{\varphi,\Gamma}(\mathcal{R}_{B,F}(\delta_1/\delta_2))$ is surjective, there exists an element $s'_B\in\mathcal{E}_B^1$ lifting $s_B$ that represents the cohomology class $[c_B^1]$ defined by  the extension $D_B$. Hence there exists some element $c'^0_B\in C^0_B$ such that $c^1_B-f^1(s'_B)=\partial^0(c'^0_B)$. 
We need to modify $s'_B$ and $c'^0_B$, as a priori $c'^0_B$ does not necessarily reduce to $0$ over $\bar B$.
However, it follows from the construction that its reduction maps to zero under $\partial^0$, i.e.
$$c'^0_{\bar B}\in{\rm ker}(\partial^0:C^0_{\bar B}\rightarrow C^1_{\bar B})={\rm ker}(\mathcal{E}^0_{\bar B}\rightarrow \mathcal{E}^1_{\bar B}).$$
In particular we may regard $c'^0_{\bar B}$ as an element $e^0_{\bar B}$ of $\mathcal{E}^0_{\bar B}$. Let $e^0_B\in\mathcal{E}^0_B$ be a lift of $c'^0_{\bar B}$ and set 
\begin{align*}
s_B&=s'_B+d(e^0_B),\\
c^0_B&=c'^0_B-f^0(e^0_B).
\end{align*}
Then $c^0_{\bar B}=c'^0_{\bar B}-f^0(e^0_{\bar B})=0$ and $c^1_B-f^1(s_B)=\partial^0(c'^0_B)-f^1(d(e^0_B))=\partial^0(c^0_B)$ as desired.

\emph{Step 2}: Let $U\subset \mathfrak{X}_T$ be an open neighborhood
of a point $(\delta_1,\delta_2)$ that is not weakly regular and that
admits a smooth surjection from an affinoid space ${\rm Sp}(A)$.  Then
there exists some $N\gg 0$ such that
$H^2_{\varphi,\Gamma}(\mathcal{R}_{A,F}(z^N\delta_1/\delta_2))=0$,
compare Remark \ref{EHrem:propertiesphiGammacohom} (iii). We consider
the neighborhood $U'=\{(\eta_1,z^N\eta_2)\mid (\eta_1,\eta_2)\in U\}$
of $(\delta'_1,\delta'_2)=(\delta_1,z^N\delta_2)\in\mathfrak{X}_T^{\rm
  wreg}$ and the morphism $a:U\rightarrow U'$ given by
$(\eta_1,\eta_2)\mapsto (\eta_1,z^N\eta_2)$. This morphism fits in a
diagram
$$
\begin{xy}
\xymatrix{
& V'\ar[d]\\
\mathfrak{X}_B|_{U}\ar[d]\ar[r]^{b}&\mathfrak{X}_B|_{U'}\ar[d]\\
U\ar[r]^a& U',
}
\end{xy}
$$
where $ V'\rightarrow \mathfrak{X}_B|_{U'}$ is a smooth surjection from a rigid analytic space (that exists by Step 1) and where $b$ is the morphism defined by the pullback
\numequation\label{EHeqn:pullbackofextension}
\begin{aligned}
\begin{xy}
\xymatrix{
0\ar[r] & \mathcal{R}(\eta_1)\ar[r]\ar@{=}[d]&D\ar[r]&\mathcal{R}(\eta_2)\ar[r]&0\\
0\ar[r] & \mathcal{R}(\eta_1)\ar[r]&D'\ar@{^{(}->}[u]^{\iota}\ar[r]&z^N\mathcal{R}(\eta_2)\ar[r]\ar@{^{(}->}[u]&0.
}
\end{xy}
\end{aligned}
\end{equation}
We define $V\rightarrow \mathfrak{X}_B|_U$ as the projection to the first factor of the fiber product
$$V= V'\times_{\mathfrak{X}_B|_{U'}}\mathfrak{X}_B|_U\longrightarrow \mathfrak{X}_B|_U.$$
As $V'\rightarrow \mathfrak{X}_B|_{U'}$ is a smooth surjection we are left to show that $b$ is relatively representable.
This is the content of Lemma \ref{EHlem:modifyextensions} below.
This finishes the proof that $\mathfrak{X}_B$ is an Artin stack.

We refer to Remark \ref{EHrem: construction of compactification} below for some comments on how to adapt the strategy to prove the corresponding claims for $\overline{\mathfrak{X}}_B$ and $\overline{\mathfrak{X}}_B^{\rm naive}$.
\end{proof}

\begin{lem}\label{EHlem:modifyextensions}
The morphism $b:\mathfrak{X}_B\rightarrow\mathfrak{X}_B$ defined by mapping an extension $D$ to the pullback $D'$ defined by $(\ref{EHeqn:pullbackofextension})$ is relatively representable.
\end{lem}
\begin{proof}[Sketch of proof.]
Let $V'={\rm Sp}(A)\rightarrow \mathfrak{X}_B$ be a morphism defining the extension $D'$. 
For any $D$ as in $(\ref{EHeqn:pullbackofextension})$ the morphism $\iota$ induces an isomorphism $D'[1/t]\isoto D[1/t]$. Let us write $\hat D$ for the scalar extension 
$\hat D=D\otimes_{\mathcal{R}_F}F_\infty[[t]]$ and similarly $\hat D'$. By the Beauville--Laszlo gluing lemma, the $(\varphi,\Gamma)$-module $D$ may be reconstructed from $D'$ and the datum of the $\Gamma$-stable $(A\otimes_{\Q_p}F_\infty)[[t]]$-lattice in the free $(A\otimes_{\Q_p}F_\infty)((t))$-module $\hat D'[1/t]$. We may assume that $D'$ is free and choose a basis of $D'$ such that the first basis vector spans the subobject $\mathcal{R}_{A,F}(\eta_1)$. Using this basis we may view $\hat D$ as a $\Gamma$-stable $A$-valued point of the affine Grassmannian $${\rm Grass}_\infty=L({\rm Res}_{F_\infty/\Q_p}{\rm GL}_2)/L^+({\rm Res}_{F_\infty/\Q_p}{\rm GL}_2).$$
The lattice $\hat D'$ is then identified with the base point $e_0\in{\rm Grass}_\infty$. 
Moreover the fact that $D$ is an extension as in $(\ref{EHeqn:pullbackofextension})$ translates to the fact that the lattice $\hat D$ lies in the closed subspace $$t^\nu L^{<0}U_\infty t^{-\nu}e_0\subset {\rm Grass}_\infty,$$ where $\nu=(1,-N)$ and where $U_\infty\subset {\rm Res}_{F_\infty/\Q_p}{\rm GL}_2$ is the Weil restriction of the unipotent upper triangular matrices. \\
Conversely, the $\Gamma$-stable lattices in $t^\nu L^{<0}U_\infty t^{-\nu}e_0$ precisely describe the pullback of $V'\rightarrow \mathfrak{X}_B$ along $b$.\\
There are two technical difficulties in this proof: 
\begin{enumerate} 
\item[(a)] As $F_\infty$ is not finite over $\Q_p$, the Weil restriction ${\rm Res}_{F_\infty/\Q_p}{\rm GL}_2$ is not representable by a rigid analytic space. 
\item[(b)] The subspace $t^\nu L^{<0}U_\infty t^{-\nu}e_0$ is an ind-scheme (or ind-rigid analytic space) rather than a scheme (or a rigid analytic space). 
\end{enumerate}
By working a bit more carefully, we can overcome (a) by descending  to some finite subextension $F_n=F(\mu_{p^n})\subset F_\infty$.
To overcome (b), we need to show that the $\Gamma$-invariant lattices in fact lie inside some closed rigid analytic subspace of $t^\nu L^{<0}U_\infty t^{-\nu}e_0$.
In the case $F=\Q_p$ the main ingredients of this argument are as follows (for general $F$ a variant of this argument works as well): on $L$-valued points either $\hat D$ is uniquely determined by $\hat D'$ and the diagram $(\ref{EHeqn:pullbackofextension})$, or $D'$ is de Rham up to twist. In the latter case both $\hat D'$ and $\hat D$ lie in certain affine Schubert cells with respect to a third base point $\hat D''$ in the affine Grassmannian (in the de Rham case the lattice $\hat D''$ can be chosen such that it defines a $(\varphi,\Gamma)$-module with Hodge--Tate weights $(0,0)$ and the Schubert cells containing $\hat D'$ respectively $\hat D$ depend on their respective Hodge--Tate weights).
\end{proof}

\begin{rem}\label{EHrem: construction of compactification}
In order to prove that $\overline{\mathfrak{X}}^{\rm naive}_B$ and $\overline{\mathfrak{X}}_B$ are rigid analytic Artin stacks we follow the same idea as in the above sketch. 
In order to construct charts for $\overline{\mathfrak{X}}_B$ we have to consider the closure of $t^\nu L^{<0}U_\infty t^{-\nu}e_0$ in the affine Grassmannian. To construct charts for $\overline{\mathfrak{X}}_B^{\rm naive}$ we have to consider the union of certain $t^\nu L^{<0}U_\infty t^{-\nu}e_0$. 
\end{rem}

One important difference with the classical case of vector bundles on algebraic curves is that the map $$\overline{\beta}_B:\overline{\mathfrak{X}}_B\longrightarrow \mathfrak{X}_G$$ is neither proper nor surjective (and a similar remark applies to more general parabolic subgroups as well).  
The $(\varphi,\Gamma)$-modules in the image of this map are all trianguline, as the $(\varphi,\Gamma)$-modules over $\overline{\mathfrak{X}}_B^{\rm naive}$, and hence over $\overline{\mathfrak{X}}_B$, admit a complete $\varphi$ and $\Gamma$-stable flag after inverting $t$, by the very definition of $\overline{\mathfrak{X}}_B^{\rm naive}$, compare Remark \ref{EHremexteriorpowers} for the image of this flag under the Plücker embedding.  
This implies that $\overline{\beta}_B$ is not surjective, as not every Galois representation (or more generally not every $(\varphi,\Gamma)$-module) is trianguline (or equivalently: trianguline after inverting $t$). The easiest way to see this is to compare the dimensions of the stacks $\mathfrak{X}_d$ and $\mathfrak{X}_B$ (or of some related rigid analytic spaces of Galois representations).

The morphism $\overline{\beta}_B$ fails to be proper as well, as will be explained in Remark \ref{EH:rem density of XB} below.
Failure of properness and surjectivity are two (at the first glance surprising) differences from the classical situation of stacks of vector bundles on an algebraic curve. The following theorem however explains that $\overline{\beta}_B$ is still rather close to being proper.

\begin{theorem}\label{EHThm propernessofcompactific}
The  morphism
$$
\overline{\mathfrak{X}}_B\longrightarrow\mathfrak{X}_G\times\mathfrak{X}_{T}
$$
is representable by rigid analytic spaces and proper.
\end{theorem}
This result is proven using results on triangulations in families by Liu \cite{MR3433281} and Kedlaya--Pottharst--Xiao \cite{MR3230818}. Of course we expect similar results for the compactifications $\overline{\mathfrak{X}}_P$ and $\widetilde{\mathfrak{X}}_P$, as stated in the following conjecture.
We can also conjecturally relate the stack $\overline{\mathfrak{X}}_B$ (or rather its image in $\mathfrak{X}_d\times\mathfrak{X}_T$) to the \emph{trianguline variety} and hence to eigenvarieties and $p$-adic automorphic forms. Conjectures \ref{EHconj:localglobal1} and \ref{EHconj:localglobal2} will make the relation to spaces of $p$-adic automorphic forms more precise.
\begin{conj}
\label{conj:parabolic images}
\hfill\\(i)The  morphisms
\numalign\label{EHeqnmorphtoXtri}
\overline{\mathfrak{X}}_P&\longrightarrow\mathfrak{X}_G\times\mathfrak{X}_{S_M}\ \text{and} \\
\widetilde{\mathfrak{X}}_P&\longrightarrow\mathfrak{X}_G\times\mathfrak{X}_{M} \nonumber 
\end{align}
are representable by rigid analytic spaces and proper.
\hfill\\(ii) The scheme-theoretic images of the morphisms \eqref{EHeqnmorphtoXtri} are equidimensional of dimension  $\dim {\rm Res}_{F/\mathbf{Q}_p}P$.\\
(iii) Let $P=B$ and let $R_{\bar\rho}$ denote the framed deformation ring of a fixed residual representation $\bar\rho:\Gal_F\rightarrow {\rm GL}_d(k)$. Then the pullback of the scheme theoretic image of \eqref{EHeqnmorphtoXtri} along the smooth map
$$({\rm Spf}\,R_{\bar\rho})^{\rm rig}\times \mathcal{T}^d\longrightarrow \mathfrak{X}_G\times\mathfrak{X}_T$$
coincides with the trianguline variety of \cite{MR3623233}.
\end{conj}
Here scheme-theoretic image should of course be understood as an analogue of the usual notion of scheme-theoretic images in the context of rigid analytic Artin stacks.

\begin{rem}\label{EH:rem density of XB}
Even though Theorem \ref{EHThm propernessofcompactific} might be interpreted as $\overline{\beta}_B$ being very close to being proper, there is one aspect in which it behaves very differently from a proper map: it is expected that the image of $\overline{\mathfrak{X}}_B$ in $\mathfrak{X}_G$ is Zariski-dense (as opposed to the image of a proper map being Zariski-closed, as the direct image of the structure sheaf is coherent)
This image should be regarded as an analogue of the infinite fern of Gouvea--Mazur \cite{MR1486834} (and its generalization by Chenevier \cite{MR3037022} and Nakamura \cite{MR3235551}). The infinite fern argument implies directly that the image $\overline{\beta}_B(\overline{\mathfrak{X}}_B)$ is dense in the union of all irreducible components of $\mathfrak{X}_G$ that contain a crystalline point in their interior. We expect that all components of $\mathfrak{X}_G$ contain a crystalline point in their interior (which in fact would be a consequence of the conjectural description of the connected components of $\mathfrak{X}_G$ as a consequence of Conjecture \ref{EHconjglobalsections}).
In light of the relation with eigenvarieties the failure of properness can also be linked to the fact that the construction of eigenvarieties involves Fredholm hypersurfaces of infinite degree over the weight space.
As this failure of properness is a failure of properness over the weight space, rather than over the space of Galois representations, the implication is not immediate, but these two phenomena are related.
\end{rem}
\begin{rem}\label{EH:rem:derived}
From a purely deformation-theoretic viewpoint the expected dimension of the stack $\mathfrak{X}_P$ of equivariant $P$-bundles is $\dim{\rm Res}_{F/\mathbf{Q}_p}P$.
In fact we  expect that $\mathfrak{X}_P$ is equidimensional of this
dimension and is a local complete intersection. 
The $p$-adic situation seems to be rather different to the situation with the spaces of $L$-parameters  \cite[Remark 2.2]{hellmann2020derived}, \cite[Remark 2.3.8]{zhu2020coherent} where the stacks of $L$-parameters for parabolic subgroups can fail to be equidimensional (whereas the stacks of $L$-parameters for reductive groups are always equidimensional and complete intersections). In the situation with stacks of $L$-parameters one is hence forced to work with derived stacks due to this phenomenon. 
In the $p$-adic case we expect that the stacks $\mathfrak{X}_P$ (and similarly their compactifications $\overline{\mathfrak{X}}_P$ and $\widetilde{\mathfrak{X}}_P$) do not have any additional non-trivial derived structure.
In fact for fixed $d$ one can show that $\mathfrak{X}_B$ is equidimensional if $[F:\mathbf{Q}_p]$ is large (compared with $d$). For the time being we do not have an argument for all $F$, but it seems reasonable to expect equidimensionality independent of the degree $[F:\mathbf{Q}_p]$.
\end{rem}

Motivated by Conjecture~\ref{conj:parabolic images} we write $\mathfrak{X}_{G,\rm tri}\subset \mathfrak{X}_G\times\mathfrak{X}_T$ for the scheme-theoretic image of $\overline{\mathfrak{X}}_B$ under the proper map $\overline{\beta}_B$. 
The following theorem gives a partial description of the preimages of the morphism $\mathfrak{X}_{G,\rm tri}\rightarrow \mathfrak{X}_G$ that will turn out to be important in the construction of companion points in eigenvarieties (i.e.~in the construction of overconvergent $p$-adic automorphic forms of finite slope with prescribed associated Galois representation) in Section \ref{sec:eigenvar}. 
Before describing the preimages we need some preparation.

Let $D\in\mathfrak{X}_G(L)$ be a crystalline $(\varphi,\Gamma)$-module and assume that the eigenvalues $\varphi_1,\dots,\varphi_d$ of the Frobenius on the associated Weil--Deligne representation ${\rm WD}(D)$ lie in $L$ and satisfy $\varphi_i/\varphi_j\neq 1,q$ for all $i\neq j$. Moreover, we write $\lambdau=(\lambda_{\tau,1}\geq \dots\geq  \lambda_{\tau,d})_{\tau:F\hookrightarrow L}$ for the Hodge--Tate weights of $D$ and assume that $\lambdau$ is regular, i.e.~we assume $\lambda_{\tau,i}\neq \lambda_{\tau,j}$ for all $\tau$ and $i\neq j$. 
We fix an ordering of the Frobenius eigenvalues $\varphi_1,\dots,\varphi_d$ or equivalently a complete Frobenius stable flag $\mathcal{F}_\bullet$ on ${\rm WD}(D)$. For each embedding $\tau:F\hookrightarrow L$ we define the Weyl group element $w_{\mathcal{F},\tau}\in\mathcal{S}_n$ as the relative position of the $\tau$-part of the Hodge filtration
$${\rm Fil}_\tau^\bullet\subseteq {\rm WD}(D)={\rm WD}(D)_\tau\subseteq \prod\nolimits_{\tau':F\hookrightarrow L}{\rm WD}(D)_{\tau'}={\rm WD}(D)\otimes_{\mathbf{Q}_p}L$$
with respect to $\mathcal{F}_\bullet$ (note that by the regularity assumption on $\lambdau$ the filtration ${\rm Fil}^\bullet_\tau$ is indeed a full flag). 

We associate  to the flag $\mathcal{F}_\bullet$ (equivalently, to the
ordering of the Frobenius eigenvalues $\varphi_1,\dots,\varphi_d$) the unramified character $\delta_{\mathcal{F}_\bullet}:T=(F^\times)^n\rightarrow L^\times$ given by 
\numequation\label{EH:eqndeltasm}
(x_1,\dots, x_n)\longmapsto \prod_{i=1}^d {\rm unr}_{\varphi_i}(x_i).
\end{equation}
Moreover, associated to a tuple $(\mu_\tau)_\tau\in\mathbf{Z}^{{\rm Hom}(F,L)}$ we have the character 
\numequation\label{EH:defofztau}
z^\mu:x\longmapsto \prod_\tau \tau(x)^{\mu_\tau}.
\end{equation}
More generally for $\underline{\mu}=(\mu_1,\dots,\mu_d)$ with $\mu_i\in \mathbf{Z}^{{\rm Hom}(F,L)}$ we write
$$z^{\underline{\mu}}:(x_1,\dots, x_d)\longmapsto \prod_{i=1}^d z^{\mu_i}(x_i)=\prod_{i=1}^d\prod_\tau \tau(x_i)^{\mu_{i,\tau}}.$$
The following characterization of the points in $\mathfrak{X}_{G,\rm tri}$ above $D$ is \cite[Theorem 4.2.3]{MR4028517}.
\begin{theorem}\label{thm: local companion points}
Let $D\in\mathfrak{X}_G(L)$ be a crystalline $(\varphi,\Gamma)$-module of regular Hodge--Tate weight $\lambdau$ and such that the eigenvalues $\varphi_1,\dots,\varphi_d$ of the Frobenius on ${\rm WD}(D)$ lie in $L$ and satisfy $\varphi_i/\varphi_j\neq 1,q$ for $i\neq j$. Then
\begin{align*}
&\{\delta\in\mathfrak{X}_T\mid (D,\delta)\in \mathfrak{X}_{G,\rm tri}\}\\&=\bigcup_{\mathcal{F}_\bullet}\{\delta_{\mathcal{F}_\bullet} z^{w\lambdau}| w=(w_\tau)_\tau\in W\ \text{such that}\  w_{\mathcal{F},\tau}\preceq w_\tau w_0\ \text{for all}\ \tau\},
\end{align*}
where the union runs over all Frobenius stable flags $\mathcal{F}$ of ${\rm WD}(D)$ and where $w_0\in\mathcal{S}_d$ is the longest element. 
\end{theorem}
\begin{rem}
\hfill\\(i) Following Hansen \cite[Conjecture 6.2.3]{MR3692014}, there is a general conjectural
description of the fibers of $\mathfrak{X}_{G,\rm tri}\rightarrow \mathfrak{X}_G$ without assuming that the
$(\varphi,\Gamma)$-module $D$ is crystalline (and without the
regularity assumptions in the theorem). \\
(ii) In the case of crystabelline $(\varphi,\Gamma)$-modules this description is more or less the same as in the crystalline case. In the semi-stable case (and more generally in the semistabelline case) we only take those flags $\mathcal{F}_\bullet$ that are stable under the monodromy on the Weil--Deligne representation.\\
(iii) It is reasonable to ask for a generalization from semistabelline to all de Rham $(\varphi,\Gamma)$-modules $D$.
In this case the image of $\overline{\mathfrak{X}}_B$ in $\mathfrak{X}_G\times\mathfrak{X}_T$ should be replaced by the image of $\overline{\mathfrak{X}}_P$ in $\mathfrak{X}_G\times\mathfrak{X}_{S_M}$, where $P$ is the parabolic subgroup containing $B$ with Levi $M=\prod_{i=1}^r{\rm GL}_{d_i}$ such that the Weil--Deligne representation ${\rm WD}(D)$ associated to $D$ is a direct sum of indecomposable Weil--Deligne representations $D_1,\dots, D_r$ of respective dimension $d_i$ (compare \cite{BreuilDingeigenvar}).
\end{rem}

The following expectation is motivated by the relation of the trianguline variety with eigenvarieties and with the patched version of the eigenvariety \cite{MR3623233} (and by Conjecture \ref{EH:conjanalytic}(2) respectively Conjecture \ref{EHconj: coherentdescriptionoffiniteslop}).
\begin{conj}
Let us write $\pi_B:\overline{\mathfrak{X}}_B\rightarrow \mathfrak{X}_G\times\mathfrak{X}_T$ for the canonical {\em (}proper{\em )} projection. \\
\noindent (i) There are isomorphisms
$$R\pi_{B,\ast}\mathcal{O}_{\overline{\mathfrak{X}}_B}\cong\mathcal{O}_{\mathfrak{X}_{G,\rm tri}},\ R\pi_{B,\ast}\omega_{\overline{\mathfrak{X}}_B}\cong\omega_{\mathfrak{X}_{G,\rm tri}}.$$
\noindent (ii) The direct image
$R\pi_{B,\ast}(\mathcal{O}_{\overline{\mathfrak{X}}_B}([F:\Q_p]\rho'))$
is concentrated in degree $0$ and is a Cohen--Macaulay module {\em (}i.e.~its dual $R\pi_{B,\ast}(\omega_{\overline{\mathfrak{X}}_B}(-[F:\Q_p]\rho'))$ is also concentrated in degree $0${\em )}.
\end{conj}
Here $\rho'$ is the algebraic character $(0,-1,\dots, -d+1)$ of $T=\mathbf{G}_m^d$ that we view as a line bundle on $\ast/\mathbf{G}_m^d$ and pull it back to $\overline{\mathfrak{X}}_B$ along the canonical map
$$\overline{\mathfrak{X}}_B\rightarrow \mathfrak{X}_T=\mathfrak{X}_1^d=(\mathcal{T}/\mathbf{G}_m)^d\rightarrow \ast/\mathbf{G}_m^d,$$
see Section~\ref{sec: GL1} for the identification of $\mathfrak{X}_1$.

Moreover, it makes sense to ask whether the trianguline variety $\mathfrak{X}_{G,\rm tri}$ is normal and Cohen--Macaulay (note that Cohen--Macaulayness would be a consequence of (i) of the conjecture, if $\overline{\mathfrak{X}}_B$ is Cohen--Macaulay). We point out that the sheaf in (ii) of the conjecture should be regarded as a variant of the coherent Springer sheaf from \cite{benzvi2020coherent} (respectively \cite[4.4]{zhu2020coherent} or \cite[2.3]{hellmann2020derived}) and hence (ii) should be seen as a variant of \cite[Conjecture 4.15]{benzvi2020coherent} (respectively \cite[Conjecture 4.4.2]{zhu2020coherent} or \cite[Conjecture 2.17]{hellmann2020derived}).

In the case $F=\Q_p$ and ${\rm GL}_2$ we can prove the above conjecture and compute the singularities as well as the fiber dimensions of the Cohen--Macaulay modules in (ii) (these fiber dimensions have some interpretation in terms of the conjectural correspondence of locally analytic representations with coherent sheaves on $\mathfrak{X}_d$, Conjecture~\ref{EH:conjanalytic}). 

\begin{thm}\label{EHtheo triangvar for GL2Qp}
Let $F=\Q_p$ and $d=2$, i.e.~$G={\rm GL}_2$.\\
\noindent (i) The stack $\overline{\mathfrak{X}}_B$ is smooth.\\
\noindent (ii) The stack $\mathfrak{X}_{G,\rm tri}$ is a rigid analytic Artin stack\footnote{As we do not know whether $\mathfrak{X}_G=\mathfrak{X}_2$ is a rigid analytic Artin stack this should be interpreted as: there is a rigid analytic Artin stack $\mathfrak{X}_{G.\rm tri}$ that embeds as a closed substack into $\mathfrak{X}_G\times \mathfrak{X}_T$. The map $\overline{\mathfrak{X}}_B\rightarrow \mathfrak{X}_G\times\mathfrak{X}_T$ factors through $\mathfrak{X}_{G,\rm tri}$ and is scheme-theoretically dense.} and normal and Cohen-Macaulay.\\
\noindent (iii) The morphism $\pi_B:\overline{\mathfrak{X}}_B\rightarrow \mathfrak{X}_{G,\rm tri}$ is birational with geometrically connected fibers and the canonical map 
$$\cO_{\mathfrak{X}_{G,\rm tri}}\rightarrow R\pi_{B,\ast}\cO_{\overline{X}_B}$$
is an isomorphism, and as a consequence of duality, there also is an isomorphism $$R\pi_{B,\ast}\omega_{\overline{X}_B}\cong \omega_{\mathfrak{X}_{G,\rm tri}}.$$
\noindent (iv) Let $\mathcal{F}=R\pi_{B,\ast}(\mathcal{O}_{\overline{\mathfrak{X}}_B}(\rho'))$ and $\mathcal{G}=R\pi_{B,\ast}(\omega_{\overline{\mathfrak{X}}_B}(-\rho'))$. Then $\mathcal{F}$ and $\mathcal{G}$ are concentrated in degree zero and in particular are Cohen-Macaulay modules. \\
\noindent (v) For $x=(D,\delta_1,\delta_2)\in\mathfrak{X}_{G,\rm tri}$ the fibers $\mathcal{F}\otimes k(x), \mathcal{G}\otimes k(x)$ and $\omega_{\mathfrak{X}_{G,\rm tri}}\otimes k(x)$ are one-dimensional unless $D=\mathcal{R}(\delta)\oplus t^{-n}\mathcal{R}(\delta)$ with $n\geq 0$ and $\delta_1=\delta$. In this case
\begin{align*}
\dim \mathcal{F}\otimes k(x)&=2\\
\dim \mathcal{G}\otimes k(x)&=
\begin{cases}
2 & n=0\\
1 & n\geq 1
\end{cases}\\
\dim \omega_{\mathfrak{X}_{G,\rm tri}}\otimes k(x)&=\begin{cases}3 & n=0\\ 2 & n\geq 1.\end{cases}
\end{align*}
\end{thm}

\subsubsection{Local models}\label{subsec:localmodels}
Assume for the remainder of this section that $L$ contains all Galois conjugates of $F$ (inside $\bar F$).
We discuss the relation of stacks of $(\varphi,\Gamma)$-modules over
$\mathcal{R}_F$ with stacks of (semi-linear) $\Gamma$-representations
on $F_\infty[[t]]$. 
We don't know whether it can be expected that the latter stacks are Artin stacks (it seems that they can be rather ill-behaved), but their completions at a fixed Hodge--Tate weight turn out to be well-behaved. 
 At certain ``nice" points the morphism from a stack of $(\varphi,\Gamma)$-modules to these stacks is formally smooth and hence we obtain rather explicit local models for stacks like $\overline{\mathfrak{X}}_B$, generalizing the local models for the trianguline variety of \cite{MR4028517}.

We define the groupoid $\mathfrak{X}_d^{{\rm dR},+}$ over ${\rm Rig}_L$ as follows. For an affinoid algebra $A$ we set
\begin{align*}
\mathfrak{X}_d^{{\rm dR,+}}({\rm Sp}\,A)&=\left\{\begin{array}{*{20}c} \text{continuous, semi-linear}\ \Gamma\text{-representations}  \\ \text{on finite projective}\ A\widehat{\otimes}F_\infty [[t]]\text{-modules of rank}\ d\end{array}\right\},
\end{align*}
where (by slight abuse of notation) we write $A\widehat{\otimes}F_\infty[[t]]=(A\otimes_{\mathbf{Q}_p}F_\infty)[[t]]$.
\begin{rem}
\hfill\\(i) In a similar way we can define a groupoid $\mathfrak{X}_d^{\rm dR}$ by mapping ${\rm Sp}\,A$ to the groupoid of finite projective $A\widehat\otimes F_\infty((t))$-modules of rank~$d$ with a semi-linear $\Gamma$-action that (locally on ${\rm Sp}\,A$) admit an $A\widehat{\otimes}F_\infty[[t]]$-lattice. \\
(ii) The motivation to consider these stacks (and the explanation for their name) is Fontaine's theorem~\cite{MR2104360} that the category of semi-linear $\Gamma$-representations on finite free $F_\infty[[t]]$-modules (respectively on finite dimensional $F_\infty((t))$-vector spaces) is equivalent to the category of semi-linear $\Gal_F$-representations on finite free $B_{\rm dR}^+$-modules (respectively on finite dimensional $B_{\rm dR}$-vector spaces).  We expect that a similar equivalence holds true in families over rigid analytic spaces, but, as we do not have a proof, we restrict our attention to the semi-linear $\Gamma$-representations, as they seem to be easier to handle.
\end{rem}

\begin{expectation} \hfill\\(i)The groupoids $\mathfrak{X}_d^{\rm dR,+}$ and $\mathfrak{X}_d^{\rm dR}$ are stacks. \\
(ii) The diagonal of $\mathfrak{X}_d^{\rm dR, +}$ is representable.\\
(iii) The morphism $\mathfrak{X}_d^{\rm dR,+}\rightarrow \mathfrak{X}_d^{\rm dR}$ is representable.
\end{expectation}

Though we do not know whether one should expect that $\mathfrak{X}_d^{\rm dR,+}$ and $\mathfrak{X}_d^{\rm dR}$ are Artin stacks, both stacks admit versal rings at rigid analytic points, and more generally their completions along a fixed Hodge--Tate--Sen weight are formal rigid analytic Artin stacks (a notion that we will not formally introduce here). In order to make sense of this statement, we note that, similarly to the definition of the weight map $\omega_d$ in \eqref{EHeqnweightmap}, we have a weight map
$$\omega_d^{\rm dR}:\mathfrak{X}_d^{\rm dR,+}\longrightarrow {\rm WT}_{d,F}\cong\prod_{F\hookrightarrow L} \mathbf{A}^d/\mathcal{S}_d$$
that maps, on $A$-valued points, an $A\widehat{\otimes}F_\infty[[t]]$-module $D_A$ to the characteristic polynomial of the derivative of the $\Gamma$-action on $D_A/tD_A$ (recall that we have fixed $L$ large enough such that $[F:\mathbf{Q}_p]=|{\rm Hom}(F,L)|$).

\begin{rem}
The Hodge--Tate--Sen weights of an object $D_A\in \mathfrak{X}_d^{\rm dR}({\rm Sp}\,A)$ are only well-defined mod $\mathbf{Z}$. More precisely, under the identification of $A\otimes_{\mathbf{Q}_p}F=\prod_{F\hookrightarrow L}A$ we can define the weight  of $D_A$ as an element of 
$$\prod\nolimits_{F\hookrightarrow L}\big(\big(A/\mathbf{Z}\big)^d/\mathcal{S}_d\big).$$
\end{rem}

Fontaine's theory of almost de Rham representations gives a description of the stacks $\mathfrak{X}_d^{\rm dR}$ and $\mathfrak{X}_d^{\rm dR,+}$ after completion at a fixed integral Hodge--Tate weight, and similar results hold true for the completion at any fixed Hodge--Tate--Sen weight. 
We briefly recall Fontaine's classification \cite{MR2104360} of $B_{\rm dR}$-representations (or equivalently $F_\infty((t))$-representations). 
Let us write ${\rm Mod}_{B_{\rm dR}}^{\Gal_F}$ for the category of continuous semi-linear $\Gal_F$-representations on finite dimensional $B_{\rm dR}$ vector spaces, and similarly we define ${\rm Mod}_{B_{\rm dR}^+}^{\Gal_F}$, ${\rm Mod}_{F_\infty((t))}^{\Gamma}$ and ${\rm Mod}_{F_\infty[[t]]}^{\Gamma}$ in the obvious way. 
For a Hodge--Tate--Sen weight $a\in F$ we denote by ${\rm Mod}_{B_{\rm dR}}^{\Gal_F,a}$ the full subcategory of all objects whose Hodge--Tate--Sen weights are congruent to $a$ modulo $\mathbf{Z}$. Then 
$${\rm Mod}_{B_{\rm dR}}^{\Gal_F}=\bigoplus_{a\in F/\mathbf{Z}} {\rm Mod}_{B_{\rm dR}}^{\Gal_F,a}$$
is a block decomposition of ${\rm Mod}_{B_{\rm dR}}^{\Gal_F}$,
i.e.~the higher Ext-groups between objects in different summands vanish (this
is a consequence of  \cite[Thm.\ 3.19]{MR2104360}). 

The block ${\rm Mod}_{B_{\rm dR}}^{\Gal_F,0}$ is called the block of almost de Rham representations (or \emph{presque de Rham} representation in French). It has the following description in terms of linear algebra objects (see \cite[Prop 3.1.1]{MR4028517}):
\begin{theorem}[Fontaine]
There is an equivalence of categories
\[D_{\rm pdR}:{\rm Mod}_{B_{\rm dR}^+}^{\Gal_F,0}\longrightarrow {\rm Fil}{\rm Rep}(\mathbf{G}_a)\]
from the category of almost de Rham representations over $B_{\rm dR}^+$ to the category of finite dimensional $F$-vector spaces equipped with an algebraic representation of the additive group $\mathbf{G}_a$ and with a separated and exhaustive $\mathbf{Z}$-filtration ${\rm Fil}^\bullet$ stable under the $\mathbf{G}_a$-action.
Under this equivalence, passing from ${\rm Mod}_{B_{\rm dR}^+}^{\Gal_F,0}$ to ${\rm Mod}_{B_{\rm dR}}^{\Gal_F,0}$ corresponds to forgetting the filtration. 
\end{theorem}
We recall that (over fields of characteristic $0$) a $\mathbf{G}_a$-representation is nothing but the action of a nilpotent operator. 
Moreover, an almost de Rham representation $V$ over $B_{\rm dR}^+$ (or over $B_{\rm dR}$) is de Rham if and only if the nilpotent operator $N:D_{\rm pdR}(V)\rightarrow D_{\rm pdR}(V)$ (defined by the $\mathbf{G}_a$-action) is zero. 

There is a definition of the functor $D_{\rm pdR}$ in the spirit of period rings \cite[3.1]{MR4028517} that we do not recall here.
On the equivalent category ${\rm Mod}_{F_\infty[[t]]}^{\Gamma}\cong {\rm Mod}_{B_{\rm dR}^+}^{\Gal_F}$ one can characterize almost de Rham representations and the functor $D_{\rm pdR}$ as follows: 
Given an $F_\infty[[t]]$-module $\Lambda$ with continuous $\Gamma$-action, we can consider (similarly to the definition of Hodge--Tate--Sen weights) the derivation of the $\Gamma$-action at $1$ which gives a derivation 
$$\partial_\Lambda:\Lambda\longrightarrow\Lambda$$
above the derivation $t\tfrac{d}{dt}$ on  $F_\infty[[t]]$. 
We can then define $$D_{\rm pdR,\infty}(\Lambda)=\bigcup_{N\geq 1}{\rm ker}(\partial_\Lambda^N:\Lambda[1/t]\rightarrow\Lambda[1/t])$$ as the sub-$F_\infty$ vector space of $\Lambda[1/t]$ on which $\partial_\Lambda$ is nilpotent. 
In particular $\partial_\Lambda$ induces an $F_\infty$-linear, nilpotent endomorphism $N_\infty$ of $D_{\rm pdR,\infty}(\Lambda)$.
The representation $\Lambda$ is almost de Rham if and only if the canonical map
\[D_{\rm pdR,\infty}(\Lambda)\otimes_{F_\infty}F_\infty((t))\longrightarrow \Lambda[1/t]\]
is an isomorphism. In that case we write $\Lambda_0=D_{\rm pdR,\infty}(\Lambda)\otimes_{F_\infty}F_\infty[[t]]$ and refer to $\Lambda_0$ as the standard lattice. 
Then the lattice $\Lambda\subset\Lambda_0[1/t]$ defines a separated and exhaustive $\mathbf{Z}$-filtration ${\rm Fil}^\bullet_{\infty}$ on $D_{\rm pdR_\infty}(\Lambda)=\Lambda_{0}/t\Lambda_0$ by sub-$F_\infty$ vector spaces that are stable under the action of $N_\infty$, as $\Lambda$ is stable under $\partial_\Lambda$.
Again in a similar fashion to the definition of Hodge--Tate--Sen weights, the triple $(D_{\rm pdR,\infty}(\Lambda),N_\infty,{\rm Fil}_\infty^\bullet)$ is equipped with a $\Gamma$-action which allows us to descend\footnote{We caution the reader that this descent is not completely formal and requires some work. Recall that we pointed out in the discussion of Hodge--Tate weights that $(D_{\rm Sen}(D_A),\Theta)$ does not descend to an object over $F\otimes_{\Q_p}A$. However, this can be done after completion at a fixed integral Hodge--Tate weight. The issue here is a similar one.} $(D_{\rm pdR,\infty}(\Lambda),N_\infty,{\rm Fil}_\infty^\bullet)$ from $F_\infty$ to  the desired triple $(D_{\rm pdR}(\Lambda),N,{\rm Fil}^\bullet)$ over $F$.

The equivalence in Fontaine's theorem lies at the heart of the following description of the completions of the stacks of $B_{\rm dR}^+$-representations. Note that for an $L$-algebra $A$ we have $$F_\infty[[t]]\widehat{\otimes}_{\mathbf{Q}_p}A=\prod_{\tau:F\hookrightarrow L} F_\infty[[t]]\widehat{\otimes}_FA,$$
and hence we can  apply Fontaine's description for each embedding $\tau$ and for each labeled Hodge--Tate--Sen weight separately.
Let us write $G={\rm Res}_{F/\mathbf{Q}_p}{\rm GL}_d$ and $\mathfrak{g}={\rm Lie}\,G$ for the remainder of this subsection. Moreover, we fix an integral labeled Hodge--Tate weight 
$$\lambdau\in\prod_{\tau:F\hookrightarrow L}(\mathbf{Z}^d)/\mathcal{S}_d\subseteq {\rm WT}_{d,F}(L)$$
 We denote by $[\lambdau]$ its class modulo $\mathbf{Z}$ (which of course only remembers that the weight was integral but does not distinguish between various $\lambdau$).

We describe the formal completions $(\mathfrak{X}_d^{\rm dR,+})^{\widehat{}}_{[\lambdau]}$ and $(\mathfrak{X}_d^{\rm dR})^{\widehat{}}_{[\lambdau]} $ of the groupoids $\mathfrak{X}_d^{\rm dR.+}$ and $\mathfrak{X}_d^{\rm dR}$ along the fixed weight $\lambdau$ (respectively along the locus where the Hodge--Tate--Sen weight is congruent to $\lambdau$ modulo $\mathbf{Z}$). This result is very similar to \cite[Cor. 3.1.6, Cor. 3.1.9]{MR4028517}. 
Let 
$$\widetilde{\mathfrak{g}}_{P_{\lambdau}}=\{(A,\mathcal{F}^\bullet)\in \mathfrak{g}\times G/P_{\lambdau}\mid A\mathcal{F}^\bullet\subset\mathcal{F}^\bullet\}\subseteq \mathfrak{g}\times G/P_{\lambdau}$$
be the variant of Grothendieck's simultaneous resolution of singularities that parameterizes pairs $(A,\mathcal{F}^\bullet)$ consisting of an endomorphism $A$ and an $A$-stable flag $\mathcal{F}^\bullet$ of type $P_{\lambdau}$ (compare Section \ref{EHsub:sub:deRham} for the notation).
Recall that $\widetilde{\mathfrak{g}}_{P_{\lambdau}}$ is a subvariety of $\prod_\tau \mathfrak{gl}_d\times {\rm GL}_d/P_{\lambda_\tau}$
and again we think of the filtration on $G/P_{\lambdau}$ as an $[F:\mathbf{Q}_p]$-tuple of $\mathbf{Z}$-filtrations that jump in degrees $\lambdau$. 

\begin{prop}\label{EH:local model XG}
The functor  $D_{\rm pdR,\infty}$ induces isomorphisms
\begin{align*}
 (\mathfrak{X}_d^{\rm dR})^{\widehat{}}_{[\lambdau]} &\isoto\mathfrak{g}^{\widehat{}}_{0}/G\\
 (\mathfrak{X}_d^{\rm dR,+})^{\widehat{}}_{\lambdau}& \isoto(\widetilde{\mathfrak{g}}_{P_{\lambdau}})^{\widehat{}}_{0}/G,
\end{align*}
where on the right hand side we complete at the closed subspace of all nilpotent endomorphisms (i.e.~at the nilpotent cone, respectively its partial Springer resolution). 
\end{prop}

\begin{rem}
\noindent \hfill\\(a) The proposition implies that the
completion of $\mathfrak{X}_d^{\rm dR,+}$ at an integral Hodge--Tate weight is a formal rigid analytic Artin stack (with an appropriate definition of such an object).\\
\noindent (b) In fact, applying  Fontaine's theory for arbitrary (i.e.~not necessarily integral) weights  $\lambdau$ one can obtain a similar description of the completion of $\mathfrak{X}_d^{\rm dR,+}$ respectively $\mathfrak{X}_d^{\rm dR}$ at $\lambdau$. 
\end{rem}

In the neighborhood of certain ``nice" points, $\mathfrak{X}_d^{\rm dR,+}$ can be thought of as a local model of the stack $\mathfrak{X}_d$. A sample result for this is the following theorem, compare \cite[Theorem 3.1.1]{HMS}. We will discuss the related situation for $B$-structures (resp.~$P$-structures) and their Drinfeld compactifications below.
\begin{theorem}
Let $x\in\mathfrak{X}_d$ be a crystalline point with regular Hodge--Tate weight and Frobenius eigenvalues $\varphi_1,\dots,\varphi_d$ that satisfy $\varphi_i/\varphi_j\notin \{1,q\}$. Then the morphism $\mathfrak{X}_d\rightarrow \mathfrak{X}_d^{\rm dR,+}$ is formally smooth at $x$.
\end{theorem}

Given a parabolic subgroup $P\subseteq {\rm GL}_d$ we can consider $P$-structures (i.e.~partial flags) on objects in $\mathfrak{X}_d^{\rm dR}$ and $\mathfrak{X}_d^{\rm dR,+}$, which allows us to define stacks $\mathfrak{X}_P^{\rm dR}$  and $\mathfrak{X}_P^{\rm dR,+}$. 
Parallel to the above discussion of Drinfeld compactifications, we would like to define ``compactifications" $\overline{\mathfrak{X}}_P^{\rm dR,+}$ and $\widetilde{\mathfrak{X}}_P^{\rm dR,+}$ of the stack $\mathfrak{X}_P^{\rm dR,+}$ that fit into the diagram
\[\begin{xy}
\xymatrix{
& \widetilde{\mathfrak{X}}_P^{\rm dR,+} \ar[d]\ar[ddl] \ar[dr] \\ & \overline{\mathfrak{X}}_P^{\rm dR,+}\ar[dl] \ar[dr]& \mathfrak{X}_M^{\rm dR,+}\ar[d]\\
\mathfrak{X}_d^{\rm dR,+} && \mathfrak{X}_{S_M}^{\rm dR,+}.
}
\end{xy}
\]
 
We will define such compactifications, and prove finiteness statements about them after formal completion at an integral weight $\lambdau$.

As above we fix an integral labeled  Hodge--Tate weight $\lambdau=((\lambda_{\tau,1},\dots, \lambda_{\tau,d})_{\tau:F\hookrightarrow L}$. 
We restrict ourselves to the easiest case in which $\lambdau$ is regular. As $\lambdau$ is only well-defined up to permutation, we choose its dominant representative, i.e.~we assume $\lambda_{\tau,1}>\lambda_{\tau,2}>\dots>\lambda_{\tau,d}$ for all $\tau$.
 In particular the parabolic subgroup $P_{\lambdau}$ equals the Weil restriction of the Borel subgroup $({\rm Res}_{F/\mathbf{Q}_p}B)_L$. We write $\widetilde{\mathfrak{g}}=\widetilde{\mathfrak{g}}_{P_{\lambdau}}$ for simplicity. 

We now give a description similar to Proposition \ref{EH:local model XG} of the completion of $\mathfrak{X}_P^{\rm dR,+}$ and use this description to define ``compactifications" of the formal completion. 
Consider the space 
\numequation\label{EH:defXBP}
X_{B,P}=\widetilde{\mathfrak{g}}\times_{\mathfrak{g}}\widetilde{\mathfrak{g}}_P=\bigcup_{w\in W_P} V_w\subseteq \mathfrak{g}\times G/B\times G/P,
\end{equation}
where $W=\prod_\tau \mathcal{S}_d$ is the Weil group of ${\rm Res}_{F/\mathbf{Q}_p}{\rm GL}_d$ and $W_P\subseteq W$ is the set of minimal length representatives of elements in $W/W_M$ ($W_M\subseteq W$ being the Weyl group of the Levi ${\rm Res}_{F/\mathbf{Q}_p}M$ of ${\rm Res}_{F/\mathbf{Q}_p}P$). For $w=(w_\tau)_\tau\in W_p$, the space $V_w$ is defined to be the preimage of the Bruhat cell $G(1,w)\subseteq G/{\rm Res}_{F/\mathbf{Q}_p}B\times G/{\rm Res}_{F/\mathbf{Q}_p}P$
(recall that we write $G={\rm Res}_{F/\mathbf{Q}_p}{\rm GL}_d$).  
Over $V_w$, the intersection of the universal flag on $G/{\rm Res}_{F/\mathbf{Q}_p}B$ and over $G/{\rm Res}_{F/\mathbf{Q}_p}P$ defines a well-defined flag in ${\rm Res}_{F/\mathbf{Q}_p}M/{\rm Res}_{F/\mathbf{Q}_p}B_M$. Given $w$, we write $w\lambdau$ for the permutation of $\lambdau$ by $w$. Moreover, we can view $w\lambdau$ as a weight of $M$.
In particular for each $w\in W$, we obtain a diagram
\numequation\label{EHeqnlocalmodelalphabeta}
\begin{xy}
\xymatrix{
&V_w/G\ar[dl]\ar[dr]&\\
\widetilde{\mathfrak{g}}/G&&\widetilde{\mathfrak{m}}/{\rm Res}_{F/\mathbf{Q}_p}M,
}
\end{xy}
\end{equation}
where by abuse of notation we write $\widetilde{\mathfrak{m}}$ for the Lie algebra of ${\rm Res}_{F/\mathbf{Q}_p}M$.
\begin{prop}
The functor $D_{\rm pdR}$ induces an isomorphism of the diagram 
\[
\begin{xy}
\xymatrix{
&(\mathfrak{X}_P^{\rm dR,+})^{\widehat{}}_{ww_0\lambdau}\ar[dl]\ar[dr]\\
(\mathfrak{X}_d^{\rm dR,+})^{\widehat{}}_{\lambdau}&& (\mathfrak{X}_M^{\rm dR,+})^{\widehat{}}_{ww_0\lambdau}
}
\end{xy}
\]
with the completion of \eqref{EHeqnlocalmodelalphabeta} at the closed subspace of nilpotent endomorphisms. Here $(\mathfrak{X}_P^{\rm dR,+})^{\widehat{}}_{ww_0\lambdau}$ is the completion of $\mathfrak{X}_P^{\rm dR,+}$ along the preimage of 
 $$(\mathfrak{X}_M^{\rm dR,+})^{\widehat{}}_{ww_0\lambdau}\subseteq (\mathfrak{X}_M^{\rm dR,+})$$
 under the canonical map $\mathfrak{X}_P^{\rm dR,+}\rightarrow \mathfrak{X}_M^{\rm dR,+}$.
\end{prop}
In the diagram \eqref{EHeqnlocalmodelalphabeta} there are two ways to complete $V_w/G$ into a
stack proper over $\widetilde{\mathfrak{g}}/G$, and again only one of these compactifications will admit a morphism to $\widetilde{\mathfrak{m}}/{\rm Res}_{F/\mathbf{Q}_p}M$. 
The easiest way is to consider the Zariski-closure $\overline{X}_w$ of $V_w$ inside the space $X_{B,P}=\widetilde{\mathfrak{g}}\times_{\mathfrak{g}}\widetilde{\mathfrak{g}}_P$ defined in \eqref{EH:defXBP}. In this case the intersection of the two flags over $G/{\rm Res}_{F/\mathbf{Q}_p}B$ and $G/{\rm Res}_{F/\mathbf{Q}_p}P$ does not induce a well-defined family of filtrations on the graded pieces (with respect to the $P$-filtration), but only on their top exterior powers. We obtain a diagram
\[
\begin{xy}
\xymatrix{
&\overline{X}_w/G\ar[dl]\ar[dr]&\\
\widetilde{\mathfrak{g}}/G&&{\rm Lie}\,{\rm Res}_{F/\mathbf{Q}_p}S_M/{\rm Res}_{F/\mathbf{Q}_p}S_M.
}
\end{xy}
\]
Note that the formation of $\overline{X}_w$ is functorial with respect to inclusions $P'\subseteq P$ as can be deduced from the identification 
$$X_{B,P}/G=q^{-1}(\mathfrak{p})/{\rm Res}_{F/\mathbf{Q}_p}P,$$
where $q:\widetilde{\mathfrak{g}}\rightarrow \mathfrak{g}$ is the canonical projection and $\mathfrak{p}\subseteq \mathfrak{g}$ is the Lie algebra of ${\rm Res}_{F/\mathbf{Q}_p}P\subseteq G$.
The stack $\overline{X}_w/G$ is canonically isomorphic to $\overline{X}'_w/ {\rm Res}_{F/\mathbf{Q}_p}P$, where $\overline{X}_w'\subseteq q^{-1}(\mathfrak{p})$ is the closure of the locally closed subset $V'_w\subseteq q^{-1}(\mathfrak{p})$ such that $V_w/G=V'_w/{\rm Res}_{F/\mathbf{Q}_p}P$.

Note that for fixed $w$ the intersection of the universal flag on $G/{\rm Res}_{F/\mathbf{Q}_p}B$ with the partial standard flag (of type $P$) gives rise to a well-defined flag in the graded pieces of the standard flag. This induces a map $V'_w\rightarrow \widetilde{\mathfrak{m}}$ and  we can  define $\widetilde{X}_w$ to be the closure of $V'_w$ inside $q^{-1}(\mathfrak{p})\times \widetilde{\mathfrak{m}}$. Then we obtain a diagram
\[
\begin{xy}
\xymatrix{
&\widetilde{X}_w/{\rm Res}_{F/\mathbf{Q}_p}P\ar[dl]\ar[dr]&\\
\widetilde{\mathfrak{g}}/G&&\widetilde{\mathfrak{m}}/{\rm Res}_{F/\mathbf{Q}_p}M,
}
\end{xy}
\]
as well as a canonical projection 
$$\pi_P:\widetilde{X}_w/{\rm Res}_{F/\mathbf{Q}_p}P\longrightarrow \overline{X}_w/G$$
compatible with the projections to $\widetilde{\mathfrak{g}}$.

\begin{prop}
\hfill\\(i) If $P=B$, then the scheme $\overline{X}_w=\widetilde{X}_w$ is normal and Cohen--Macaulay.\\
(ii) For a general parabolic subgroup $P$ the scheme $\overline{X}_w$ is unibranched.
\end{prop}
These results about the local geometry of $\overline{X}_w$ are proven
in \cite[Prop.\ 2.3.3, Thm.\ 2.3.6]{MR4028517} respectively in \cite[Thm.\ 2.12]{ZhixiangWu}.
\begin{expectation}
\hfill\\(i) The schemes $\overline{X}_w$ and $\widetilde{X}_w$ are normal and Cohen--Macaulay. \\
(ii) The morphism  
$$\pi_P:\widetilde{X}_w/{\rm Res}_{F/\mathbf{Q}_p}P\longrightarrow\overline{X}_w/G$$
satisfies $R\pi_{P,\ast}\mathcal{O}=\mathcal{O}$ and $R\pi_{P,\ast}\omega=\omega$.
\end{expectation}

We can use these objects to define compactifications of the completion $\mathfrak{X}_P^{\rm dR,+}$ along fixed weights $w\underline{\lambda}$ by setting
\[(\overline{\mathfrak{X}}_P^{\rm dR,+})^{\widehat{}}_{w\lambdau}=(\overline{X}_w/G)^{\widehat{}}_0\ \text{and}\ (\widetilde{\mathfrak{X}}_P^{\rm dR,+})^{\widehat{}}_{w\lambdau}=(\widetilde{X}_w/P)^{\widehat{}}_0.\]
\begin{rem}
Of course, as the notation suggests, we would like to be able to define stacks $\overline{\mathfrak{X}}_P^{\rm dR,+}$ and $\widetilde{\mathfrak{X}}_P^{\rm dR,+}$ before completion. 
\end{rem}
By construction these compactifications fit into a diagram
\[\begin{xy}
\xymatrix{
& (\widetilde{\mathfrak{X}}_P^{\rm dR,+})^{\widehat{}}_{w\lambdau} \ar[d]\ar[ddl] \ar[dr] \\ & (\overline{\mathfrak{X}}_P^{\rm dR,+})^{\widehat{}}_{w\lambdau}\ar[dl] \ar[dr]& (\mathfrak{X}_M^{\rm dR,+}\ar[d])^{\widehat{}}_{w\lambdau}\\
(\mathfrak{X}_d^{\rm dR,+})^{\widehat{}}_{\lambdau} && (\mathfrak{X}_{S_M}^{\rm dR,+})^{\widehat{}}_{\det(w\lambdau)},
}
\end{xy}
\]
where we write $\det(w\lambda)$ for the weight of $S_M$ that is given by the top exterior power on each ${\rm GL}_m$-block of the Levi $M$. 

Finally the following result reconnects these stacks, and their explicit local description (respectively: their definition), to the Drinfeld compactifications of $\mathfrak{X}_P$. Part (ii) of the theorem is a generalization of the local model for the trianguline variety in \cite{MR4028517}.
\begin{theorem}
\hfill\\(i) The canonical map
$$(\mathfrak{X}_P)^{\widehat{}}_{w\lambdau} \longrightarrow (\mathfrak{X}_P^{\rm dR,+})^{\widehat{}}_{w\lambdau}$$
induced by $\mathfrak{X}_P\rightarrow \mathfrak{X}_P^{\rm dR,+}$ on completions extends to morphisms
\begin{align*}
(\overline{\mathfrak{X}}_P)^{\widehat{}}_{w\lambdau} &\longrightarrow (\overline{\mathfrak{X}}_P^{\rm dR,+})^{\widehat{}}_{w\lambdau}\\
(\widetilde{\mathfrak{X}}_P)^{\widehat{}}_{w\lambdau} &\longrightarrow (\widetilde{\mathfrak{X}}_P^{\rm dR,+})^{\widehat{}}_{w\lambdau}.
\end{align*}
(ii) Let $\widetilde{x}\in\widetilde{\mathfrak{X}}_P$ and write $\overline{x}$ for its image in $\overline{\mathfrak{X}}_P$. Assume that the image $x_M=(D_1,\dots, D_r)$ of $\widetilde{x}$ in $\mathfrak{X}_M$ is regular in the sense that
\begin{enumerate}
\item[-] the map $\mathfrak{X}_M\rightarrow \mathfrak{X}_M^{\rm dR,+}$ is formally smooth at $x_M$.
\item[-] for $i\neq j$ the Ext-groups ${\rm Ext}^0(D_i,D_j)=0={\rm Ext}^2(D_i,D_j)=0$ vanish.
\end{enumerate}
Then the morphism $$ (\overline{\mathfrak{X}}_P)^{\widehat{}}_{w\lambdau} \longrightarrow (\overline{\mathfrak{X}}_P^{\rm dR,+})^{\widehat{}}_{w\lambdau}$$ is formally smooth at $\overline{x}$.
\end{theorem}
\begin{remark}
In the situation of (i) of the theorem there is an obvious commutative diagram that we expect to be cartesian.
\end{remark}
\begin{conj}
In the situation of (ii) of the theorem the morphism 
 $$ (\widetilde{\mathfrak{X}}_P)^{\widehat{}}_{w\lambdau} \longrightarrow (\widetilde{\mathfrak{X}}_P^{\rm dR,+})^{\widehat{}}_{w\lambdau}$$ is formally smooth at $\widetilde{x}$.
\end{conj}
Of course in combination with the results on the local structure of the varieties $\overline{X}_{w}$ this gives an explicit description of the local geometry of the Drinfeld compactifications.

In the case $P=B$ the following theorem makes the situation even a bit more precise:
Let $\varphi_1,\dots,\varphi_d\in L^\times$ such that $\varphi_i/\varphi_j\neq 1,q$ for $i\neq j$. 
We write $\delta_{\rm sm}$ for the unramified character $(t_1,\dots, t_d)\mapsto {\rm unr}_{\varphi_1}(t_1)\cdots{\rm unr}_{\varphi_d}(t_d)$. 
For $w\in W$ let us write $$\overline{\mathfrak{X}}_{B,z^{w\lambdau}\delta_{\rm sm}}\subseteq \overline{\mathfrak{X}}_B$$ for the preimage of $z^{w\lambdau}\delta_{\rm sm}\in \mathfrak{X}_T$ and let 
$$\overline{\mathfrak{X}}_{B,z^{w\lambdau}\delta_{\rm sm}}^{\widehat{}}\subseteq \overline{\mathfrak{X}}_B$$
denote the formal completion of $\overline{\mathfrak{X}}_B$ along this closed substack. 
We write $$\mathfrak{X}_{d,(\lambdau,\delta_{\rm sm}){\rm -tri}}=\bigcup_{w\in W} \overline{\beta}_B(\overline{\mathfrak{X}}_{B,z^{w\lambdau}\delta_{\rm sm}})\subseteq \mathfrak{X}_d$$
for the union of the images of $\overline{\mathfrak{X}}_{B,z^{w\lambdau}\delta_{\rm sm}}$  in $\mathfrak{X}_d$. Note that the assumptions on $\delta_{\rm sm}$ imply that the restriction of $\overline{\beta}_B$ to (the formal completion of) $\overline{\mathfrak{X}}_{B,z^{w\lambdau}\delta_{\rm sm}}$ is a closed embedding. 

By definition this is the closed substack of $\mathfrak{X}_d$ consisting of all $(\varphi,\Gamma)$-modules $D$ such that $D$ admits a filtration ${\rm Fil}_\bullet$ with graded pieces $$\mathcal{R}_F({\rm unr}_{\varphi_i}z^{(\lambda_{w_{\tau}(i)})_{\tau}})$$ (see \ref{EH:defofztau} for the notation) for some $w\in W$. Over the closed substack $\mathfrak{X}_{d,(\lambdau,\delta_{\rm sm}){\rm -tri}}$ this filtration does not glue to a filtration of the universal $(\varphi,\Gamma)$-module $\tilde D$ (as the Weyl group element $w$ might vary). However, this turns out to be true after inverting $t$, and we write $\mathcal{M}_\bullet\subseteq \tilde D[1/t]$ for this filtration. With this notation we can consider the formal substack
$$\mathfrak{X}_{d,(\lambdau,\delta_{\rm sm}){\rm -tri}}^{\widehat{}}\subset\mathfrak{X}_d$$
parameterizing deformations $\tilde D'$ of $D$ to infinitesimal neighborhoods of $\mathfrak{X}_{d,(\lambdau,\delta_{\rm sm}){\rm -tri}}$ such that $\mathcal{M}_\bullet$ lifts to a (by our choice of $\delta_{\rm sm}$: necessarily unique) filtration $\mathcal{M}'_\bullet\subset\tilde D'[1/t]$ stable under $\varphi$ and $\Gamma$. The following theorem is basically proven in \cite{MR4028517}.
\begin{theorem}\label{EH:Theorem:localmodelB=P}
\hfill\\(i) The morphism 
$$\overline{\mathfrak{X}}_{B,z^{w\lambdau}\delta_{\rm sm}}^{\widehat{}}\longrightarrow (\overline{\mathfrak{X}}_{B}^{\rm dR,+})^{\widehat{}}_{w\lambdau}\cong(\overline{X}_{ww_0}/G)^{\widehat{}}_0$$
is formally smooth with relative dualizing sheaf $\mathcal{L}_{2\rho}^{\otimes[F:\Q_p]}$.\\ %
(ii) There is a formally smooth morphism 
$$f_{\lambdau,\delta_{\rm sm}}:\mathfrak{X}_{d,(\lambdau,\delta_{\rm sm}){\rm -tri}}^{\widehat{}}\longrightarrow (\tilde{\mathfrak{g}}\times_{\mathfrak{g}}\tilde{\mathfrak{g}})^{\widehat{}}_0$$
with relative dualizing sheaf $\mathcal{L}_{2\rho}^{\otimes [F:\Q_p]}$, whose restriction to $\overline{\beta}_B(\overline{\mathfrak{X}}_{B,z^{w\lambdau}\delta_{\rm sm}}^{\widehat{}})$ is given by the morphism in (i). 
\end{theorem}
We define the sheaf $\mathcal{L}_{2\rho}$ that occurs in the statement of the theorem:
The stack $\overline{\mathfrak{X}}_{B,z^{w\lambdau}\delta_{\rm
    sm}}^{\widehat{}}$ has a canonical map to $\ast/\mathbf{G}_m^d$
(as it sits inside $\overline{\mathfrak{X}}_B$ which maps to
$\mathfrak{X}_T$ and then further projects to
$\ast/\mathbf{G}_m^d$). Hence it makes sense to pull back the line bundle that is associated to the character given by the sum of the positive roots $2\rho$. This defines the line bundle in (i).
To make sense of this line bundle in (ii) we note that we can write $\mathcal{L}_{2\rho}$ as a pullback of a $G$-equivariant line bundle on $\overline{X}_{ww_0}$ that is in fact defined on all of $\tilde{\mathfrak{g}}\times_{\mathfrak{g}}\tilde{\mathfrak{g}}\subset \mathfrak{g}\times G/B\times G/B$, as it is given by the pullback of a line bundle (the line bundle corresponding to $2\rho$) on one of the $G/B$-factors. 

\begin{example} Let us spell out the above result in a bit more detail in the case of ${\rm GL}_2(\mathbf{Q}_p)$. \\
(i) 
Let $\delta_1={\rm unr}_{\varphi_1}z^{\lambda_1}$ and $\delta_2={\rm
  unr}_{\varphi_2}z^{\lambda_2}$ be two locally algebraic characters
with unramified smooth parts, which we assume satisfy $\varphi_1/\varphi_2\notin\{1,p^{\pm 1}\}$
and $\lambda_1\neq \lambda_2$. We consider a $(\varphi,\Gamma)$-module that is an extension 
\numequation\label{EHeqntriangulineext}
0\rightarrow \mathcal{R}_L(\delta_1)\rightarrow D\rightarrow \mathcal{R}_L(\delta_2)\rightarrow 0.
\end{equation}
There are basically two cases to consider: either $D$ is crystalline, or the extension is non-split and $\lambda_1<\lambda_2$. In the crystalline case we also need to distinguish between $\lambda_1<\lambda_2$ and $\lambda_1>\lambda_2$.
\begin{enumerate}
\item[(a)] The extension is non-split and $\lambda_1<\lambda_2$. In this case $D$ is Hodge--Tate but not de Rham and \eqref{EHeqntriangulineext} is the only way to write $D$ as an extension of two rank one $(\varphi,\Gamma)$-modules. In particular $(\delta_1,\delta_2)$ is the only tuple of characters such that $(D,\delta_1,\delta_2)\in\mathfrak{X}_B$. 
But if we set $\delta'_1={\rm unr}_{\varphi_1}z^{\lambda_2}$ and $\delta'_2={\rm unr}_{\varphi_2}z^{\lambda_1}$, then $(D,\delta'_1,\delta'_2)\in \overline{\mathfrak{X}}_B$. The points 
$$x_1=(D,\delta_1,\delta_2)\ \text{and}\ x_s=(D,\delta_1,\delta_2)$$
are the only points in $\overline{\mathfrak{X}}_B$ above
$D$. Consequently, locally at $D$ the space
$\mathfrak{X}_{d,(\lambdau,\delta_{\rm sm}){\rm -tri}}^{\widehat{}}$
has two components: one component is the image of a neighborhood
(inside $\overline{\mathfrak{X}}_B$) of $x_1$, and the other is the image of a neighborhood (inside $\overline{\mathfrak{X}}_B$) of $x_s$ under the projection $\overline{\mathfrak{X}}_B\rightarrow \mathfrak{X}_d$.
\item[(b)] The extension $D$ is crystalline. In this case there are again two cases to distinguish:
\begin{itemize}
\item[(b1)] We have $\lambda_1>\lambda_2$, and the extension is non-split. In this case the extension
  is necessarily crystalline. Moreover, the
  space $\mathfrak{X}_{d,(\lambdau,\delta_{\rm sm}){\rm
      -tri}}^{\widehat{}}$ is smooth at $D$.
  
There are exactly two points in
$\overline{\mathfrak{X}}_B$ above $D$,  and both actually  belong to $\mathfrak{X}_B$. The second point does not arise as in (a) above, but from the fact that $D$ can also be written as an extension 
$$0\rightarrow \mathcal{R}_L(\eta_1)\rightarrow D\rightarrow \mathcal{R}_L(\eta_2)\rightarrow 0,$$
with $\eta_1={\rm unr}_{\varphi_2} z^{\lambda_1}$ and $\eta_2={\rm
  unr}_{\varphi_1}z^{\lambda_2}$. That is: we have changed the
\emph{refinement}, i.e.\ the order of the $\varphi_1$ and $\varphi_2$
(whereas in (a) we have changed the ordering of the weights
$\lambda_1,\lambda_2$), i.e.\ in Theorem \ref{EH:Theorem:localmodelB=P}
these points would be addressed by a different choice of the smooth
character $\delta_{\rm sm}$.

\item[(b2)] We have $\lambda_1>\lambda_2$, and the extension is split; or $\lambda_1<\lambda_2$. Then there are three points in $\overline{\mathfrak{X}}_B$, two of which correspond to the ordering $(\varphi_1,\varphi_2)$, i.e.~belong to $\mathfrak{X}_{d,(\lambdau,\delta_{\rm sm}){\rm -tri}}^{\widehat{}}$. Their description is exactly as in (a). The only difference with (a) is that there is a third point in $\overline{\mathfrak{X}}_B$ mapping to $D\in\mathfrak{X}_{2}$ corresponding to the refinement given by the ordering $(\varphi_2,\varphi_1)$. Again this point is of course not covered by Theorem \ref{EH:Theorem:localmodelB=P}.  
\end{itemize}
\end{enumerate}
(ii) Let us describe the local models in the case $d=2$. The local model $\tilde{\mathfrak{g}}\times_{\mathfrak{g}}\tilde{\mathfrak{g}}$ of Theorem \ref{EH:Theorem:localmodelB=P} has two irreducible components both of which are non-singular. 
To see this, and in order to give explicit equations for the local models, we note that $\tilde{\mathfrak{g}}\times_{\mathfrak{g}}\tilde{\mathfrak{g}}$ is smoothly equivalent to $q^{-1}(\mathfrak{b})$, where $q:\tilde{\mathfrak{g}}\rightarrow \mathfrak{g}$ is the canonical projection. More precisely 
$$(\tilde{\mathfrak{g}}\times_{\mathfrak{g}}\tilde{\mathfrak{g}})/G=(q^{-1}(\mathfrak{b}))/B$$
as stacks. Obviously $$X_1=\mathfrak{b}\times \{[1:0]\}\subset q^{-1}(\mathfrak{b})\subset \mathfrak{b}\times\mathbf{P}^1$$ is one of the irreducible components. 
The other component $X_s$ is the closure of 
$$V_{s}=\left\{\left.\left(\begin{pmatrix}a& b\\ 0&c\end{pmatrix}, [x:1] \right)\right| (a-c)x+b=0\right\}$$
where $s$ denotes the non-trivial element of the Weyl group $\mathcal{S}_2$ of ${\rm GL}_2$. Note that $V_s$ is smooth, open in $X_s$ and disjoint from $X_1$. 
The intersection $X_1\cap X_s$ is contained in the preimage of the origin $0=[1:0]\in\mathbf{P}^1$ and in the standard affine neighborhood of this point $q^{-1}(\mathfrak{b})$ is described by
$$\left\{\left.\left(\begin{pmatrix}a& b\\ 0&c\end{pmatrix}, [1:y] \right)\right| y((a-c)+by)=0\right\}$$
whereas $X_s$ is given by the equations
$$\left\{\left.\left(\begin{pmatrix}a& b\\ 0&c\end{pmatrix}, [1:y] \right)\right| (a-c)+by=0\right\}.$$
Let us write $$Z=X_1\cap X_s=\left\{\left(\begin{pmatrix}a& b\\ 0&a\end{pmatrix}, [1:0] \right)\right\}$$ for the intersection of the two components. 
Under the formally smooth morphisms in Theorem \ref{EH:Theorem:localmodelB=P}, the points of type (a) and (b2) in the first part of the example are mapped to $Z$, whereas the points of type (b1) are mapped to
$$\left\{\left(\begin{pmatrix}a& 0\\ 0&a\end{pmatrix}, [x:1] \right)\right\}\subset V_s$$
More precisely the points of type (a) are mapped to $$\begin{pmatrix}a & b \\ 0 & a\end{pmatrix}$$ with $b\neq 0$, whereas the points of type (b2) are mapped to the diagonal matrix ${\rm diag}(a,a)$.\\
(iii) The image $\mathfrak{X}_{d,\rm tri}$ of $\overline{\mathfrak{X}}_B$ in $\mathfrak{X}_d\times\mathfrak{X}_T$ should be seen as a local Galois-theoretic counterpart of an eigenvariety. We refer to \cite{MR3623233} for a mathematically meaningful formulation of this slogan. In the case of the eigencurve we briefly indicate how the points discussed in (i) arise. 
Let $f$ be a classical modular eigenform of weight $k\geq 2$ level $\Gamma_1(N)$ for some $N$ not divisible by $p$. Then the Galois representation $\rho$ attached to $f$ is crystalline at $p$.
At level $\Gamma_1(N)\cap\Gamma_0(p)$ there are two \emph{stabilizations} $f_\alpha$ and $f_\beta$ of $f$ which are eigenforms for the same Hecke eigenvalues away from $p$ but with different $U_p$-eigenvalues, say $\alpha$ and $\beta$. These two different $U_p$-eigenvalues correspond to the two possible choices of a refinement above. 
There are at least two branches of the eigencurve that meet at $\rho$
inside the generic fiber of a deformation space of Galois
representations: the branch containing the classical modular form
$f_\alpha$ and the branch containing the classical form $f_\beta$.

In the case that $\rho$ is non-split at $p$ these are the only branches containing the Galois representation $\rho$. This corresponds to the non-split extensions in case (b1) above (note that both $f_\alpha$ and $f_\beta$ are of dominant algebraic weight). 
But if $\rho$ is split crystalline at $p$ there is third branch of the eigencurve containing a non-classical $p$-adic modular form $f'$ (of weight $1-k$). This is precisely the case (b2) discussed above. 
Let us mention explicitly that the two branches of $\mathfrak{X}_{d,(\lambdau,\delta_{\rm sm}){\rm -tri}}^{\widehat{}}$ in Theorem \ref{EH:Theorem:localmodelB=P} meeting at a given $(\varphi,\Gamma)$-module $D$ as in (b2) are not the branches containing the classical forms $f_\alpha$ and $f_\beta$, but one of these branches should be thought of as containing a classical modular form, whereas the other branch contains the non-classical form $f'$ (the points on that branch being of non-classical weight).
Finally we note that the points in (a) correspond to non-classical $p$-adic modular forms (that may or may not have classical weight).
\end{example}

\begin{rem}
In \cite{BreuilDingeigenvar} Breuil and Ding define a variant of the trianguline variety for paraboline $(\varphi,\Gamma)$-modules. 
In the stacky context discussed above their construction can be described as follows (the notations differ from the ones of Breuil--Ding):
Instead of $\mathfrak{X}_M$ we consider its closed substack of $\mathfrak{X}'_M$ of objects that are de Rham of prescribed weight and inertial type, up to a twist with an arbitrary rank $1$ object. 
Then we can define the preimage $\mathfrak{X}'_P$ of $\mathfrak{X}'_M$ in $\mathfrak{X}_{P}$ and its closure $\overline{\mathfrak{X}}'_{P}$ in the compactification $\overline{\mathfrak{X}}_{P}$.
Then, in our language, the analogue of the paraboline variety of Breuil--Ding is the image of $\overline{\mathfrak{X}}'_{P}$ inside $\mathfrak{X}_d\times\mathfrak{X}_{S_M}$. 
Breuil and Ding also construct local models for their spaces. These local models are given by the
Zariski-closure (inside the space $X_{B,P}$ of \eqref{EH:defXBP}) of the closed subspace in $V_w$ of elements whose image in $\mathfrak{m}$ is a central element of the Lie algebra. 
One might wonder whether there is variant for $\widetilde{\mathfrak{X}}_P$ instead of $\overline{\mathfrak{X}}_P$.
This variant does not appear in the work of Breuil--Ding, but there is a discussion of this in work of Huang \cite{huang2023familiesparabolinevarphigammakmodules}.
\end{rem}

\section{Categorical \texorpdfstring{$p$}{p}-adic local Langlands conjectures}\label{sec:
  categorical conjectures}In this section we state our main
conjectures, and explain some motivation for them (in particular, we
explain some motivations coming from the Taylor--Wiles method). %
\subsection{Expectations: The Banach Case}\label{sec:expectations Banach case}
We now make precise our
expectations in the Banach case. In order to do so, we introduce some
notation. Fix as usual a finite extension $F/\Qp$ and a coefficient
ring $\cO$. Set $G=\GL_d(F)$, $K=\GL_d(\cO_F)$, $Z=Z(G)$.

Our goal, as we've already explained, is to make a conjecture relating
(the derived category of) of smooth representations to coherent sheaves,
or, more generally, Ind-coherent complexes, on the formal algebraic stack~$\cX_d$.
In order to make precise statements, we require some preliminary
discussion of the relevant categories of $G$-representations.

\subsubsection{Preliminaries on categories of representations}
\label{subsubsec:rep cats}
We begin by recalling some material from Appendix~\ref{sec: reps
  of p adic analytic groups}. We let $\smG$
denote the abelian category of smooth %
representations of $G$ on $\cO$-modules,
where we define an $\cO$-linear $G$-representation 
to be smooth if any element of the representation is fixed by an open subgroup of~$G$,
and annihilated by some power of~$p$.
(See Definition~\ref{def:smooth}.) 
We let $D(\smG)$ denote the derived stable $\infty$-category of~$\smG$.

Any smooth $G$-representation admits a canonical structure of~$\cO[[G]]$-module,
where~$\cO[[G]]$ is the ring defined in Definition~\ref{defn: completed group ring of G},
and it's frequently useful to study smooth $G$-representations
from the perspective of their~$\cO[[G]]$-module structures.
To this end, Proposition~\ref{prop:smooth equivalence}
yields a fully faithful $t$-exact functor
\numequation
\label{eqn:smooth into modules}
D(\smG) \hookrightarrow D(\cO[[G]]) 
\end{equation}
(the target being the stable $\infty$-category of complexes of $\cO[[G]]$-modules)
with essential image equal to $D_{\sm}(\cO[[G]])$ (the full subcategory
of objects with smooth cohomologies).

Since there are many $\cO[[G]]$-modules
that are not smooth (e.g.\ $\cO[[G]]$ itself),
the functor~\eqref{eqn:smooth into modules} is far from being an equivalence.
On the other hand, via a consideration of pro-objects, we can construct
a functor which is close to being left-inverse to~\eqref{eqn:smooth into modules}.
Indeed, if we let $D_{\cg}^-(\cO[[G]])$ denote
the full subcategory of $D(\cO[[G]])$ consisting of complexes whose cohomologies
are countably generated as $\cO[[G]]$-modules and vanish in sufficiently high degree,
then~\eqref{eqn:functor to pro complexes}
provides a functor
\numequation
\label{eqn:modules into pro-smooth}
D_{\cg}^-(\cO[[G]]) \to \Pro D^b(\smG).
\end{equation}
(For example, $\cO[[G]]$  itself maps to the formal
pro object $\quoteslim{H,n} \cInd_H^G \cO/p^n$
of $\Pro \smG$; here $H$ runs over all compact open subgroups of~$G$
and $n$ over all positive integers, and $\cO/p^n$ is given the trivial $H$-action.)
Furthermore, the composite 
$$D_{\cg}^b(\smG) \buildrel \text{\eqref{eqn:cg smooth equivalence two}} \over \iso D_{\smcg}^b(\cO[[G]])  \hookrightarrow D_{\cg}^-(\cO[[G]])
\buildrel \text{\eqref{eqn:modules into pro-smooth}} \over \longrightarrow \Pro D^b(\smG)$$
is simply the
canonical fully faithful functor $D_{\cg}^b(\smG) \hookrightarrow \Pro D^b(\smG).$
(This provides the sense in which we regard~\eqref{eqn:modules into pro-smooth}
as close to being left-inverse to~\eqref{eqn:smooth into modules}.)

In order to properly state our categorical Langlands conjecture, 
we make a preliminary conjecture regarding the abelian category~$\smG$.

\begin{conj}
\label{conj:smooth reps are locally coherent}
The abelian category $\smG$ of smooth $\cO$-linear $G$-representations 
is {\em locally coherent}.
\end{conj}

\begin{remark}
Recall that this means that $\smG$
is compactly generated, 
and that the compact objects form an abelian subcategory of~$\smG$.
We have seen in Proposition~\ref{prop:compacts in smG}
that the compact objects in $\smG$  are precisely those representations
that are of finite presentation, in the sense of Definition~\ref{def:finitely
presented},
and that $\smG$ {\em is} generated  by these  objects.
Thus the content of Conjecture~\ref{conj:smooth reps are locally
coherent}
is that the full subcategory of $\smG$ consisting
of representations of finite presentation is
closed under the formation of kernels in~$\smG$.  (It is obviously
closed under the formation of cokernels.) 
\end{remark}

\begin{rem}\label{rem: is O[[G]] coherent?}
Since the smooth representations of finite presentation
are precisely those that are finitely  presented  as~$\cO[[G]]$-modules
(Lemma~\ref{lem:fp characterization}),
Conjecture~\ref{conj:smooth reps are locally coherent}
certainly holds if the category  $\cO[[G]]\text{-}\Mod$ of left $\cO[[G]]$-modules is locally
coherent, i.e.\ if $\cO[[G]]$ is a coherent ring, i.e.\
if the  category  of finitely  presented $\cO[[G]]$-modules  is
stable under the formation of kernels.  

  In the case~$d=1$ the ring $\cO[[G]]$ is even Noetherian, while if $d=2$ it is coherent by a result of Timmins~\cite[Thm.\ 9.7]{https://doi.org/10.48550/arxiv.2103.13765}, building on Shotton's proof~  \cite[Cor.\ 4.4]{MR4106887} of the coherence of
  the ring~$\cO[[\SL_2(F)]]$.
  In either case it thus follows
  that the category of finitely presented smooth $\cO[[G]]$-modules is
  abelian. %

  No other cases of Conjecture~\ref{conj:smooth reps are locally coherent}
  are known.   It is known that
  if~$d>2$, then $\cO[[G]]$ is not coherent, by a result of
  Timmins~\cite{https://doi.org/10.48550/arxiv.2103.13765}.
  Our main reason for believing Conjecture \ref{conj:smooth reps are locally coherent}
  is that we believe that something like Conjecture~\ref{conj: Banach functor} 
  must hold, and it is difficult to formulate such a conjecture consistently
  unless Conjecture~\ref{conj:smooth reps are locally coherent} holds.
\end{rem}

Suppose now that  Conjecture~\ref{conj:smooth reps are locally coherent} holds,
and let $\Dfp^b(\smG)$ denote the full
subcategory of $D(\smG)$ consisting of (cohomologically) bounded complexes whose
cohomologies are finitely presented $\cO[[G]]$-modules.
As noted  in the discussion  following Proposition~\ref{prop:compacts
in smG},
$\Dfp^b(\smG)$ is then precisely 
the full subcategory of $D(\smG)$ consisting of coherent objects.

We have a duality theory for
$\Dfp^b(\smG)$ defined as follows. %
We identify $\cO[[G]]^{\op}$ with $\cO[[G]]$
via $g \mapsto g^{-1}$.
In this way, we may canonically convert right~$\cO[[G]]$-modules into
left~$G$-modules.  
Then we  define
an antiequivalence %
$$\mathbf{D}: \Dfp^b(\smG) \to \Dfp^b(\smG)$$
via
$$\mathbf{D}(\text{--}) := \RHom_{\cO[[G]]}(\text{---},
\cO[[G]])[d^2[F:\Qp] + 1].$$
We anticipate that this is of amplitude $[0,d]$. %
It is
equipped with a natural isomorphism $\mathbf{D}\circ \mathbf{D} \iso \Id$,
and has the property that
\numequation\label{eqn: duality on compact induction}\mathbf{D}(\cInd_H^G V) = \cInd_H^G V^{\vee}\end{equation}
for any finite length smooth representation~$V$  of a compact open subgroup
$H$ of~$G$; here $V^{\vee}  = \Hom_{\cO}(V,E/\cO)$ denotes the Pontryagin
dual of~$V$, equipped with its contragredient $H$-action.

\subsubsection{Serre duality on $\cX_d$}
On the Galois side, we have the derived categories of (Ind-)coherent sheaves
$D_{\coh}(\cX_d)$ and $\IndCoh(\cX_d)$ (as defined in Appendix~\ref{app:
IndCoh on Ind algebraic stacks}). There is an involution $\imath$ of $\cX_d$ given in moduli-theoretic
terms by $\rho \mapsto \rho^{\vee}$. We let $\D_{\cX_d}$ denote the
antiequivalence of $D_{\coh}(\cX_d)$ %
given by composing Grothendieck--Serre duality with $\imath^*$
(i.e.\ $\imath^*\circ\sRHom_{\cO_{\cX_d}}(\text{--},\omega_{\cX_d})$ where $\omega_{\cX_d}$ is the dualizing
sheaf %
of~$\cX_d$, thought of as pro-coherent sheaf on $\cX$ (as discussed in Example~\ref{ex:dualizing
complexes}) and placed in degree $0$; recall that we anticipate that
$\cX_d$ is lci).

\subsubsection{Hodge types, inertial types, components of the Bernstein centre, and Hecke algebras}
\label{subsubsec:inertial LL}
By results of Schneider--Zink~\cite{MR1728541} (the ``inertial local Langlands
correspondence''), we may associate to any inertial type
$\tau:I_F\to\GL_d(L)$ a finite-dimensional smooth irreducible
$K$-representation $\sigma^{\crys}(\tau)$ as in~\cite[Thm.\
8.2.1]{emertongeepicture}.
 We let~$\sigma^{\crys,\circ}(\tau)$ denote
an arbitrary choice of $K$-stable
$\cO$-lattice in~$\sigma^{\crys}(\tau)$.
As is explained in~\cite{Gpatch},
the type $\tau$ determines a Bernstein component~$\Omega_{\tau}$,
and there is a canonical isomorphism
\numequation
\label{eqn:six author iso}
\mathcal{Z}_{\Omega_{\tau}} \iso  \End_G\bigl(\cInd_K^G \sigma^{\crys}(\tau)\bigr).
\end{equation}
(Here and below, for any $K$-representation~$U$,
 $\cInd_K^G U$ has the usual meaning: it is the space
of compactly supported functions $G \to U$ which are $K$-equivariant
for the left $K$-action on $G$, equipped with the $G$-action given by right
translation.)
If~$\lambdau$ is a regular Hodge type, then we have the
algebraic representation~$W_{\lambdau}$ of~$K$ (see Definition~\ref{defn: W
  lambda}).
We write
$$\Pi(\lambdau,\tau)^{\circ} := 
\cInd_K^G \bigl(W_{\underline{\lambda}} \otimes_{\cO} \sigma^{\crys,\circ}(\tau)\bigr)$$
(note that this depends on the choice 
of lattice~$\sigma^{\crys}(\tau)$),
and
$$\Pi(\lambdau,\tau) := 
\Pi(\lambdau,\tau)^{\circ}[1/p].$$
Since $W_{\underline{\lambda}}[1/p]$ is naturally a
$G$-representation, ``push-pull'' for compact induction shows that
\nummultline\label{eqn:push-pull-Pi-lambda-tau}
\Pi(\lambdau,\tau) := 
\Bigl(\cInd_K^G \bigl(W_{\underline{\lambda}} \otimes_{\cO} \sigma^{\crys,\circ}(\tau)\bigr)\Bigr)[1/p]
=  \cInd_K^G \bigl(W_{\underline{\lambda}}[1/p] \otimes_L \sigma^{\crys}(\tau)\bigr)
\\
= W_{\underline{\lambda}}[1/p] \otimes_L \cInd_K^G  \sigma^{\crys}(\tau).
\end{multline}

As $W_{\underline{\lambda}}[1/p]$ is furthermore an irreducible algebraic
representation of~$G$, while the compact induction is a smooth $G$-representation,
we then find that
\nummultline
\label{eqn:Hecke = Bernstein}
\End_G
\Bigl(
\Pi(\lambdau,\tau)
\Bigr)
= \End_G
\bigl( W_{\underline{\lambda}}[1/p] \otimes_L \cInd_K^G  \sigma^{\crys}(\tau)\bigr)
\\
=  \End_G\bigl( \cInd_K^G  \sigma^{\crys}(\tau)\bigr)
\buildrel \text{\eqref{eqn:six author  iso}} \over =  \mathcal{Z}_{\Omega_\tau}.
\end{multline}

As recalled in Section~\ref{sec:EHBernsteincenters}, the Bernstein
centre $\mathcal{Z}_{\Omega_\tau}$ can also be identified with the global
functions on a stack of Weil--Deligne representations, and as in
~\eqref{EH:eqn Bernstein map to functions on ss stack}, there is a morphism
\numequation\label{eqn: Bernstein map crystalline case}\mathcal{Z}_{\Omega_\tau}\longrightarrow \Gamma((\cX_d^{\crys,\lambdau,\tau})_\eta^{\rm rig},\mathcal{O}_{(\cX_d^{\crys,\lambdau,\tau})_\eta^{\rm rig}}).\end{equation}

\subsubsection{Statement of the conjecture}
We are now ready to state our categorical Langlands conjecture. 
\begin{conj}%
  \label{conj: Banach functor}Assume Conjecture~{\em \ref{conj:smooth reps are locally coherent}}.
Then there exists an $\cO$-linear  %
  exact fully faithful functor
$\fA:\Dfp^b(\smG) \to D_{\coh}^b(\cX_d)$
satisfying the following properties: %
\begin{enumerate}
 \item\label{item: bounded amplitude} The functor $\fA$ has bounded
cohomological amplitude. 
 \item\label{item: Linfty is a sheaf} \emph{($L_\infty$ is a sheaf.)} $L_\infty:=\fA(\cO[[G]])$ is a
   pro-coherent sheaf, concentrated in degree~$0$. Furthermore it is
   flat over~$\cO[[K]]$.  %
    \item\label{item: compatible with duality} \emph{(Compatibility with
     duality.)}There is a natural equivalence
\numequation\label{eqn: general duality expectation}\fA\circ \mathbf{D} \iso (\mathbf{D}_{\cX_d} \circ \fA)[ [F:\Qp]d(d+1)/2 + 1]\end{equation}
of contravariant functors from $\Dfp^b(\smG)$ to
$D_{\coh}(\cX_d)$. %
 \item\label{item: loc alg vectors} \emph{(Locally algebraic vectors.)} %
   For any regular Hodge type~ $\lambdau$ and %
   inertial type~$\tau$, the scheme-theoretic support of
$\fA\bigl(\widehat{\Pi(\lambdau,\tau)^{\circ}}\bigr)$
is equal to  $\cX_d^{\crys,\lambdau,\tau}$.
{\em (}The $\widehat{\phantom{X}}$ here denotes $p$-adic completion.{\em )}
Furthermore, the action of the Bernstein centre~$\mathcal{Z}_{\Omega_\tau}$ on
$\fA\bigl(\widehat{\Pi(\lambdau,\tau)^{\circ}}\bigr)[1/p]$
induced by the identification~{\em \eqref{eqn:Hecke = Bernstein}}
and the functoriality of~$\fA$
coincides with the action given by~\eqref{eqn: Bernstein map crystalline case}. %

 \end{enumerate}
\end{conj}%

\begin{rem}\label{rem:notation}
The notation $\fA$ for the functor of Conjecture~\ref{conj: Banach functor}
is drawn from the paper~\cite{zhu2020coherent}, where the same notation is used
in the $\ell \neq p$ case. %
\end{rem}

\begin{rem} 
\label{rem:Ind extension}
Assuming that it exists, the functor $\fA$ has natural Ind and Pro extensions whose considerations are important
(and indeed necessary to make sense of the conjecture as stated).
We begin with the Ind-extension of~$\fA$; this is a continuous functor
$$\Ind \Dfp^b(\smG) \to \Ind\Coh(\cX),$$
which we again denote by~$\fA$.

The category $\Ind\Dfp^b(\smG)$
is the analog, in our context, of the renormalized categories of smooth
representations considered in~\cite[\S 4.1]{zhu2020coherent}.
It admits a canonical functor
\numequation
\label{eqn:Ind smooth to smooth}
\Ind\Dfp^b(\smG) \to D(\smG),
\end{equation}
which becomes an equivalence if we restrict to  
bounded below complexes on each side.  
Indeed, assuming Conjecture~\ref{conj:smooth reps are locally coherent},
the $t$-structure on $D(\smG)$ is {\em coherent} in the
sense of Definition~\ref{def:coherent t-structure},
and $\Dfp^b(\smG)$ is then precisely the subcategory of coherent objects in~$D(\smG)$, 
so that the functor is particular instance of~\eqref{eqn:functor from IndCoh}.
In particular, the heart of $\Ind\Dfp^b(\smG)$ coincides with the heart of~$D(\smG)$,
i.e.\ with the abelian category~$\smG$.  Thus, extending $\fA$ to all of
$\Ind \Dfp(\smG)$ allows us in particular to evaluate $\fA$
on all smooth $G$-representations, not just finitely presented ones.

As an aside, we note that
in general the $t$-structure on~$D(\smG)$ need not be regular (in the sense
of Definition~\ref{def:regular t-structure}), and so~\eqref{eqn:Ind smooth to smooth}
need not be an equivalence in general.   Indeed, we expect that the $t$-structure on $D(\smG)$
{\em is} regular, so that~\eqref{eqn:Ind smooth to smooth}  {\em is} an equivalence,
precisely when $I_1$ (the pro-$p$-Iwahori subgroup of~$K$) is torsion-free.

Since $D^b(\smG)$ is a full subcategory of~$\Ind\Dfp^b(\smG),$
we may restrict (the Ind-extension of) $\fA$ to $D^b(\smG)$ so as to
obtain a functor $D^b(\smG) \to \Ind\Coh(\cX)$.  (In fact,
taking into account Remark~\ref{rem: expect amplitude d-1} below,
we see that this functor should take values in $D_{\coh}^b(\cX)$.)
We may Pro-extend this functor to obtain a functor
$$\Pro D^b(\smG)\to \Pro \Ind \Coh(\cX).$$
Composing this with the functor~\eqref{eqn:modules into pro-smooth},
we find that it makes sense to evaluate $\fA$ on objects of $D_{\cg}^-(\cO[[G]]).$
The discussion of Section~\ref{subsubsec:rep cats}
shows that the restriction of this functor to
$\Dfp^b(\smG) \iso D_{\mathrm{f.p.}\sm}^b(\cO[[G]])$ coincides with the
original functor~$\fA$
(here the target has the obvious meaning, namely it is the full subcategory
of $D_{\sm}^b(\cO[[G]])$ consisting of complexes whose cohomologies are finitely
presented; the indicated equivalence follows directly 
from Proposition~\ref{prop:smooth equivalence}
by restricting the  equivalence  stated there
to the full subcategories of  its source  and target consisting
of complexes whose cohomologies are finitely generated).
\end{rem}

\begin{rem}
\label{rem:explaining pro-nature of Linfty}
To make sense of  Conjecture~\ref{conj: Banach functor}~\eqref{item: Linfty is a sheaf},
note that, as discussed in Remark~\ref{rem:Ind extension}, 
we may evaluate $\fA$ on  $\cO[[G]]$ (thought of as a module over itself).
{\em A priori} this gives rise to an object of $\Pro \Ind (\cX_d)$, but 
part~\eqref{item: Linfty is a sheaf} asserts that in fact
we  obtain an object of $\Pro\Coh(\cX_d)$.
\end{rem}

\begin{rem}
  \label{rem: Banach FS conj}As was already mentioned in Section~\ref{subsec: FS},
it seems likely that Conjecture~\ref{conj: Banach functor}
can be viewed as a consequence of a stronger conjecture, to the effect
that there should an equivalence of (derived)
  categories as in the conjecture of Fargues--Scholze \cite[Conj.~I.10.2]{fargues--scholze},
but with $p$-adic rather than $\ell$-adic coefficients.
However it does not seem to be clear at the time of writing
  precisely which category of sheaves on~$\Bun_G$ should be considered
  in order to obtain such an equivalence.
\end{rem}%
\begin{rem}
  \label{rem: fully faithful implies spectral Hecke}The full
  faithfulness of~$\fA$ implies in particular that there should be
  identifications of  $E_1$-rings between $p$-adic derived
  Hecke algebras and certain endomorphism rings of coherent sheaves
  on~$\cX_d$ (e.g.\ $p$-adic analogues of Feng's spectral Hecke algebra~\cite{https://doi.org/10.48550/arxiv.1912.04413}); see \cite[Conj.\ 4.3.1]{zhu2020coherent} for an
  analogous conjecture in the case $\ell\ne p$.
\end{rem}

\begin{rem}
  \label{rem: semistable version if we wanted}We could equally well
  formulate Conjecture~\ref{conj: Banach functor}~\eqref{item: loc alg vectors} for potentially semistable
  representations, by modifying the representation $\sigma^{\crys}(\tau)$ (see e.g.\ \cite{MR4243650}), but for simplicity of exposition we have restricted
  to the potentially crystalline case.
\end{rem}
\begin{rem}\label{rem: expect Linfty is a kernel}
  The functor of Conjecture~\ref{conj: Banach functor} %
  is necessarily of the form considered in
  Appendix~\ref{subsec:Morita}. That is, it is of the form
$$\pi \mapsto L_{\infty}\otimes^{L}_{\cO[[G]]}\pi,$$ 
where, as in %
Conjecture~\ref{conj: Banach functor}~\eqref{item: Linfty is a
  sheaf}, $L_{\infty}:=\fA(\cO[[G]])$ is a pro-coherent sheaf on the
moduli stack~$\cX_d$, and is in addition a right $\cO[[G]]$-module
(because the endomorphism algebra of $\cO[[G]]$ is
$\cO[[G]]$).
(In the case $d=1$, we see below that $L_\infty$ is the structure
sheaf on~$\cX_1$, by~\eqref{eqn: Linfty for d=1}, and in the case of
$\GL_2(\Qp)$ we give a construction of $L_\infty$, but in general we do
not have a conjectural candidate for~$L_\infty$.)%

Indeed, %
assuming~Conjecture~\ref{conj:smooth reps are locally coherent},
we may resolve
any finitely presented smooth $\cO[[G]]$-module $\pi$ by a (possibly infinite) complex of direct sums of copies
of~$\cO[[G]]$
\[\dots\to \cO[[G]]^{\oplus m_n}\dots\to \cO[[G]]^{\oplus m_0}\to
  \pi\to 0,\] and we find that $\fA(\pi)$ is computed by
\[\dots\to \fA(\cO[[G]])^{\oplus m_n}\dots\to \fA(\cO[[G]])^{\oplus m_0}\]
and thus by
\[L_\infty\otimes_{\cO[[G]]}(\dots\to \cO[[G]]^{\oplus m_n}\dots\to
  \cO[[G]]^{\oplus m_0})=L_\infty\otimes^L_{\cO[[G]]}\pi.\] %
(Assumption~\eqref{item: bounded amplitude}, i.e.\ the boundedness of the amplitude of~$\fA$,
justifies the manipulations with these possibly unbounded below complexes.)
Since~$\Ind\Dfp^b(\smG)$ is
generated under shifts and colimits by the representations of finite
presentation, we have 
\numequation
\label{eqn:kernel for fA}
\fA(\pi)=L_\infty\otimes^L_{\cO[[G]]}\pi
\end{equation}
in general.
\end{rem}

\begin{rem}
\label{rem:explaining our intuition}
We note that the only reason for making an {\em a priori} assumption
of bounded amplitude for~$\fA$ (i.e.\ the only reason for imposing condition~\eqref{item: bounded
amplitude} in the conjecture)
is to justify the manipulations that lead to formula~\eqref{eqn:kernel for fA},
and we should emphasize that the existence
of such a formula for~$\fA$
is basic to our own point of view on Conjecture~\ref{conj: Banach functor}.
If we grant the validity of this formula, then the boundedness 
of the amplitude of~$\fA$ should follow, as we explain in Remark~\ref{rem: expect
amplitude d-1} below.
\end{rem}

\begin{rem}%
  \label{rem: expect compact inductions to go to sheaves} %
If $V$ is a topological $\cO[[K]]$-module, of finite type as an $\cO$-module
(so that it is $p$-adically complete),
then we we may regard the $p$-adically completed compactly induced representation
$(\cInd_K^G V)^\wedge$
as an object of $\Pro \smG$, and hence form
$\fA\bigl((\cInd_K^G V)^\wedge\bigr).$
(And this gives the sense
in which we form
$\fA\bigl(\widehat{\Pi(\lambdau,\tau)^{\circ}}\bigr)$
in property~\eqref{item: loc alg vectors} of Conjecture~\ref{conj: Banach functor}.)

It then follows from Conjecture~\ref{conj: Banach functor}~\eqref{item: Linfty is a
    sheaf} that %
  $\fA\bigl((\cInd_K^GV)^\wedge\bigr)$ is concentrated in degree $0$. Indeed, by
  Remark~\ref{rem: expect Linfty is a kernel}, we
  have \[\fA(\cInd_K^GV)=L_\infty\cotimes^L_{\cO[[G]]}(\cInd_K^GV)^\wedge
    =L_\infty\otimes^L_{\cO[[K]]}V \] (here $\cotimes$ indicates the
$p$-adically completed tensor product) and the claim follows from the
  hypothesized flatness of~$L_\infty$ over~$\cO[[K]]$.

  In particular, %
  the sheaves
$\fA\bigl(\widehat{\Pi(\lambdau,\tau)^{\circ}}\bigr)$
  considered in property~\eqref{item: loc alg vectors}
  are concentrated in degree $0$.
\end{rem}

\begin{rem}
  \label{rem: expect amplitude d-1}We can be much more precise
about the bounded amplitude conjectured in statement~\eqref{item: bounded amplitude}
of the conjecture.
Namely, we anticipate that the
  functor~$\fA$ has amplitude $[1-d,0]$.
Indeed, in the case~$d=2$ and $p>2$ this follows from Remark~\ref{rem: expect compact inductions to go to
    sheaves} and the resolution
  explained in Remark~\ref{rem: how we use the resolution} (at least
  if we ignore the action of the centre, which we do not expect to
  contribute to the amplitude), and we anticipate that analogous
  resolutions exist for all~$d$.
\end{rem}

\begin{rem}
  \label{rem: compact objects and fg representations}Since (by
  construction) %
the compact objects in $\Ind\Dfp^b(\smG)$ are the
objects of $\Dfp^b(\smG)$, and the compact objects in
$\Ind\Coh(\cX_d)$ are given by $D_{\coh}^b(\cX_d)$, it follows from
Lemma~\ref{lem:compact objects}~(3) that if $\pi$ is an object of
$\Ind\Dfp^b(\smG)$, then $\pi$ actually belongs to %
$\Dfp^b(\smG)$ if and
only if $\fA(\pi)$ is an object of $D_{\coh}^b(\cX_d)$.

In particular, if~$\pi$ is concentrated in degree zero, then we see
that ~$\pi$ is of finite presentation if and only if $\fA(\pi)$ is a
bounded complex with coherent cohomology sheaves. (In fact, we expect that
$\fA(\pi)$ is automatically bounded, by Remark~\ref{rem: expect
  amplitude d-1}.)

If~$\pi$ is concentrated in degree zero, and is finitely generated but
not of finite presentation, then $\fA(\pi)$ is not coherent. However,
we claim that
$H^0(\fA(\pi))$ is coherent. Indeed, since $\pi$ is finitely
generated, it admits a surjection $\cInd_K^GU\to\pi$ for some finite
length $K$-representation~$U$. Since $\fA$ is right $t$-exact (see
Remark~\ref{rem: expect amplitude d-1}),
we
have a surjection $H^0(\fA(\cInd_K^GU))\to H^0(\fA(\pi))$; and
$\fA(\cInd_K^GU)$ is coherent (because $\cInd_K^GU$ is of finite
presentation).
\end{rem}

\begin{rem}
  \label{rem: lots of ways the conjecture could be stated}We have
  tried to formulate Conjecture~\ref{conj: Banach functor} with a
  minimal set of useful properties; we explain some motivation for
  these properties in the following remarks. It is unclear to what
  extent these properties uniquely determine~$\mathfrak{A}$, but in practice they
  seem to seriously constrain it. %
\end{rem}
\begin{rem}
  \label{rem: conjectural compatibility with TW patching}%
We can think of 
$L_{\infty}$ as a ``universal'', purely local version of completed
cohomology. More precisely, under the assumption that $p\nmid 2d$, for
any point $x\in\cX_d(\Fpbar)$ 
the
paper ~\cite{Gpatch} constructs a versal morphism $f:\Spf
R_\infty\to\cX_d$ at~$x$, together
with an $R_\infty$-module $M_\infty$ with a commuting action
of~$G$. Here~$R_\infty$ is a power series ring over
a universal lifting ring for the $\Gal_F$-representation corresponding
to~$x$, with the power series variables being the ``patching
variables'' occurring in the Taylor--Wiles--Kisin method,
and~$M_\infty$ is obtained by applying the patching method to the
completed cohomology of certain unitary groups. (See Section~\ref{sec:
  TW patching} for a very similar construction for modular curves.)

 The construction of~$R_{\infty}$ and $M_\infty$ is global and depends
 on many choices, but we expect that $M_\infty=f^*L_\infty$ (see also Expected Theorem~\ref{expected
    EGZ thm}). This would confirm the
 conjecture in~\cite[\S6]{Gpatch} that~$M_\infty$ should be purely
 local; indeed the hypothetical $L_\infty$ in loc.\ cit.\ is just the
 pullback of our $L_\infty$ to a local deformation ring. It also
 explains some of our expectations for~$L_\infty$; for example,
 $M_\infty$ is flat over $\cO[[K]]$ by construction, whence our
 expectation that $L_\infty$ is flat over $\cO[[K]]$. %
\end{rem}

\begin{rem}
  \label{rem: sheaf property is an analogue of
    benzvi2020coherent}Property~\eqref{item: Linfty is a sheaf} of
  Conjecture~\ref{conj: Banach functor} is an analogue of the
  expectation  ~\cite[Conj.\ 4.15]{benzvi2020coherent} %
  that
  the coherent Springer sheaf is concentrated in degree~$0$. In our
  case it is motivated by the expectation of Remark~\ref{rem:
    conjectural compatibility with TW patching} (which implies that
  the pullback of~$L_\infty$ to any versal ring is concentrated in
  degree~$0$).
\end{rem}

\begin{rem}\label{rem: duality explicated on crystalline stack}
    By
  property~\eqref{item: loc alg vectors}, for any regular Hodge type~ $\lambdau$ and
   inertial type~$\tau$,
the sheaf
$\fA\bigl(\widehat{\Pi(\lambdau,\tau)^{\circ}}\bigr)$
 is supported on $\cX_d^{\crys,\lambdau,\tau}$. We have the closed
   immersion $i:\cX_d^{\crys,\lambdau,\tau}\into \cX_d$, which we
   expect to be pure of codimension $[F:\Qp]d(d+1)/2$ (we say
   ``expect'' because we are not aware of a dimension theory for
   formal algebraic stacks in the literature; but note that by
   \cite{MR2827797, BIP} we know the corresponding claim on versal
   rings, i.e.\ on Galois deformation rings). Granting this %
   expectation, we have $i^!\omega_{\cX_d} 
= \omega_{\cX_d^{\crys,\lambdau,\tau}}[-[F:\Qp]d(d+1)/2]$, where we
abusively write $ \omega_{\cX_d^{\crys,\lambdau,\tau}}$ for the
dualizing complex shifted by the dimension
of~$\cX_d^{\crys,\lambdau,\tau}$ (so that if
~$\cX_d^{\crys,\lambdau,\tau}$ is Cohen--Macaulay then
$\omega_{\cX_d^{\crys,\lambdau,\tau}}$ is a sheaf in degree zero). Thus
 for any $\cF\in D_{\Coh}(\cX_d^{\crys,\lambdau,\tau})$ we have %
$$i_*\sRHom_{\cO_{\cX_d^{\crys,\lambdau,\tau}}}(\cF,\omega_{\cX_d^{\crys,\lambdau,\tau}})
= \sRHom_{\cO_{\cX_d}}(i_*\cF,\omega_{\cX_d})[[F:\Qp]d(d+1)/2].$$ 
Write \[\mathbf{D}_{\cO_{\cX_d^{\crys,\lambdau,\tau}}}:D_{\Coh}(\cX_d^{\crys,\lambdau,\tau})\to
  D_{\Coh}(\cX_d^{\crys,-\lambdau,\tau^\vee})\]for the composite of
Grothendieck--Serre duality and~$\imath^*$,
where \[\imath:\cX_d^{\crys,-\lambdau,\tau^\vee}\isoto
  \cX_d^{\crys,\lambdau,\tau} \] is induced by the involution~
$\imath$ of~$\cX_d$.

Take
$\cF = \fA\bigl(\widehat{\Pi(\lambdau,\tau)^{\circ}}\bigr).$
Note that it follows from~\eqref{eqn: duality on compact induction} that
\begin{multline*}
\mathbf{D}\bigl(\widehat{\Pi(\lambdau,\tau)^{\circ}}\bigr)
:=
\mathbf{D}\Bigl(\Bigl(\cInd_{K}^{G}\bigl(W_{\lambdau}\otimes_{\cO}\sigma^{\crys,\circ}(\tau)\bigr)\Bigr)^\wedge\Bigr)
\\
=
\Bigl(\Bigl(\cInd_{K}^{G}\bigl(W^*_{-\lambdau}\otimes_{\cO}\sigma^{\crys,\circ}(\tau^\vee)\bigr)
\Bigr)^\wedge\Bigr)[1],
\end{multline*}
where we abusively write $W^*_{-\lambdau}$ for the (not dual!) Weyl module,
and choose $\sigma^{\crys,\circ}(\tau^\vee)$ to be a lattice in
$\sigma^{\crys}(\tau^\vee)$ that is dual to
our chosen lattice
$\sigma^{\crys,\circ}(\tau)$. %
Bearing in mind~\eqref{eqn: general duality expectation}, we
obtain %
\numequation\label{eqn: dual Weyl module crystalline duality}
\fA\Bigl(\Bigl(\cInd_{K}^{G}\bigl(W^*_{-\lambdau}\otimes_{\cO}\sigma^{\crys,\circ}(\tau^{\vee})\bigr)\Bigr)^\wedge\Bigr)
\iso
\mathbf{D}_{\cO_{\cX_d^{\crys,\lambdau,\tau}}}
\bigl(\widehat{\Pi(\lambdau,\tau)^{\circ}}\bigr).
\end{equation}%
\end{rem}

\begin{rem}%
  \label{rem: expect duality to imply maximal CM}Property~\eqref{item: compatible with duality} of Conjecture~\ref{conj: Banach functor}
  can be thought of explaining the fact that the patched modules (constructed by
  Taylor--Wiles--Kisin patching at finite level) are maximal
  Cohen--Macaulay over their supports. In particular we expect that we can strengthen
  property~\eqref{item: loc alg vectors} to say that for any regular Hodge type~ $\lambdau$ and
   inertial type $\tau$,
   $\fA(\widehat{\Pi(\lambdau,\tau)^{\circ}}\bigr)$
   is concentrated in degree~$0$, 
   is maximal Cohen--Macaulay over
   $\cX_d^{\crys,\lambdau,\tau}$, and
   $\fA(\widehat{\Pi(\lambdau,\tau)^{\circ}}\bigr)^{\rm rig}_{\eta}$
is locally free of rank one over %
   $(\cX_d^{\crys,\lambdau,\tau})^{\rm rig}_\eta$. We now sketch an
   explanation for this. %

Since we are assuming that
   $\fA(\widehat{\Pi(\lambdau,\tau)^{\circ}}\bigr)$
is concentrated in degree~$0$ and is supported on all of
$\cX_d^{\crys,\lambdau,\tau}$, and similarly for
$\fA\Bigl(\Bigl(\cInd_{K}^{G}\bigl(W^*_{-\lambdau}\otimes_{\cO}\sigma^{\crys,\circ}(\tau^\vee)\bigr)\Bigr)^\wedge\,\Bigr)$,\footnote{While
we have not made an explicit conjecture involving $W^*_{-\lambdau}$,
this statement follows from the version for $W_{-\lambdau}$, since
these representations are isomorphic after inverting~$p$.} it follows from~\eqref{eqn: dual Weyl module crystalline duality}
 that
   $\fA(\widehat{\Pi(\lambdau,\tau)^{\circ}}\bigr)$
is maximal Cohen--Macaulay over $\cX_d^{\crys,\lambdau,\tau}$ (see~\cite[\href{https://stacks.math.columbia.edu/tag/0B5A}{Tag
  0B5A}]{stacks-project}).

It follows that
   $\fA(\widehat{\Pi(\lambdau,\tau)^{\circ}}\bigr)^{\rm rig}_{\eta}$
is maximal Cohen--Macaulay %
   over
   $(\cX_d^{\crys,\lambdau,\tau})^{\rm rig}_\eta$. Since
   $(\cX_d^{\crys,\lambdau,\tau})^{\rm rig}_\eta$ is smooth, %
 this implies
   that
   $\fA(\widehat{\Pi(\lambdau,\tau)^{\circ}}\bigr)^{\rm rig}_{\eta}$
   is locally free, and we anticipate that 
   the full faithfulness of~$\gA$, together with \cite[Cor.\ 3.12]{Gpatch}, implies that it is furthermore
   locally free of rank one. %

Recall that
  the Cohen--Macaulay property of patched modules is usually
  established by completely different means, using the action of the
  ring~$S_\infty$ of diamond operators; see~\cite{MR1440309}. %
\end{rem}

\begin{rem}
\label{rem:interchanging involutions}
As already noted when defining our modified version $\mathbf{D}_{\cX_d}$ 
of Grothendieck--Serre duality on $\cX_d$,  
we have the ``main involution'' $\imath$ on $\cX_d$ defined as $\rho
\mapsto \rho^{\vee}$. %
This induces an auto-involution $\imath^*$ on $D^b_{\coh}(\cX_d)$.
We also have the usual main involution $g \mapsto (g^{-1})^t$ of $\GL_d(F)$,
and precomposing representations of $\GL_d(F)$ with this involution  
induces an auto-involution $\imath'$ of $\Dfp^b(\smG)$.  We expect that
$\fA$ intertwines the involutions $\imath'$ and~$\imath$.

For example, in the statement of Property~\eqref{item: compatible with
  duality} of Conjecture~\ref{conj: Banach functor},
rather than replacing Grothendieck--Serre duality $\mathbf{D}$ on~$\cX_d$
by its twist $\mathbf{D}_{\cX_d} := \imath^* \mathbf{D}$,
we could replace the duality $\mathbf D$ on $\Dfp^b(\smG)$ by
its twist ${\mathbf D}' := \imath'\circ {\mathbf D}$.  Then Property~\eqref{item: compatible
with duality} could be rephrased as saying that  $\fA$ interchanges $\mathbf D'$
and~$\mathbf D$ (up to a shift). 
Now~\eqref{eqn: duality on compact induction} implies (ignoring possible shifts) that
$\mathbf{D}'(\cInd_K^G V) = \cInd_K^G V$
(respectively that $\mathbf{D}'\bigl((\cInd_K^G V^{\circ})^{\wedge}\bigr) = 
(\cInd_K^G V^{\circ})^{\wedge}$)
if $V$ is
an irreducible mod $p$ representation
of~$K$
(respectively 
if $V^{\circ}$ is a $K$-invariant lattice
in an irreducible locally algebraic representation $V$ of~$K$).
Thus we expect that the corresponding coherent sheaves $\fA(\cInd_K^G V)$ 
(respectively $\fA\bigl((\cInd_K^G V^{\circ})^{\wedge}\bigr)$
are Grothendieck--Serre self-dual (up to a shift).  
\end{rem}

\begin{rem}%
  \label{rem: expect Banach case to have compatibility with parabolic
    induction}We also expect a compatibility with parabolic induction,
  just as in the $\ell\ne p$ and analytic cases (see Section~\ref{EH
    subsection:paraboliccompatibility} for the latter).
\end{rem}

\subsubsection{The geometric Breuil--M\'ezard conjecture}\label{subsubsec:
  geometric BM Banach}
The  geometric Breuil--M\'ezard conjecture for the stack~$\cX_d$ was formulated in
~\cite{emertongeepicture}; it was motivated by (and indeed is
equivalent to) the earlier geometric Breuil--M\'ezard conjectures
for Galois deformation rings of~\cite{MR3248725, emertongeerefinedBM},
and the reader who is unfamiliar with the Breuil--M\'ezard conjecture
is referred to these papers (and the references therein) for the
original motivations for the conjecture, which come from automorphy
lifting theorems. In this section we explain how
Conjecture~\ref{conj: Banach functor} gives a further refinement of
the Breuil--M\'ezard conjecture
(lifting the statement from an equality of cycles to an equality in a
Grothendieck group of coherent sheaves). %

By a \emph{Serre
weight} $\underline{k}$ we mean a tuple of \index{Serre weight} \index{\underline{k}}
integers~$\{k_{\sigmabar,i}\}_{\sigmabar:k\into\Fpbar,1\le i\le d}$
with the properties that \begin{itemize}
\item $p-1\ge k_{\sigmabar,i}-k_{\sigmabar,i+1}\ge 0$ for each $1\le i\le
  d-1$, and
\item $p-1\ge k_{\sigmabar,d}\ge 0$, and not every~$k_{\sigmabar,d}$
  is equal to~$p-1$.
\end{itemize}
The set of Serre weights is in bijection with the set of irreducible
$\Fpbar$-representations of $\GL_d(k_F)$, via passage to highest
weight vectors (see for example the appendix to~\cite{MR2541127});
for
each Serre weight~$\underline{k}$, we write~$F_{\underline{k}}$ for
the corresponding irreducible $k$-representation of~$\GL_d(k_F)$.

We write
$\sigmabar^{\crys}(\lambdau,\tau)$ %
for the semisimplification of the $k$-representation of~$\GL_d(k_F)$
given by
$W_{\lambdau}\otimes_\cO\sigma^{\crys,\circ}(\tau)\otimes_\cO k$.
There are unique
integers $n_{\underline{k}}^\crys(\lambdau,\tau)$ %
  such
    that \[\sigmabar^{\crys}(\lambdau,\tau)\cong\oplus_{\underline{k}}F_{\underline{k}}^{\oplus
        n_{\underline{k}}^{\crys}(\lambdau,\tau)}.\] %
Write  $Z_{\crys,\lambdau,\tau}$ for the cycle (i.e.\ an element of
the free abelian group $\Z[\cX_{d,\red}]$ on the
irreducible components of~$\cX_{d,\red}$) corresponding to the
special fibre of $\cX_d^{\crys,\lambdau,\tau}$.  The following is ~\cite[Conj.\
8.2.2]{emertongeepicture}.
\begin{conj}
  \label{conj: geometric BM}There are cycles~$Z_{\underline{k}}$ with
  the property that for each regular Hodge type~$\lambdau$ and each
  inertial type~$\tau$, we have
  $Z_{\crys,\lambdau,\tau}=\sum_{\underline{k}}n_{\underline{k}}^\crys(\lambdau,\tau)\cdot
  Z_{\underline{k}}$. %
\end{conj}

Assume Conjecture~\ref{conj: Banach functor}. For each regular Hodge type~$\lambdau$ and each
inertial type~$\tau$, after making a choice of lattice
$\sigma^{\crys,\circ}(\tau)$ in $\sigma^{\crys}(\tau)$ as above\footnote{Since we are passing to
  underlying cycles, the precise choice of lattice is irrelevant.} and (again as above) writing
$\Pi(\lambdau,\tau)^{\circ} := \cInd_K^G \bigl(W_{\lambdau}\otimes_{\cO}\sigma^{\crys,\circ}(\tau)\bigr),$ 
we set %
\[\cF_{\lambdau,\tau}=
\fA\bigl(\widehat{\Pi(\lambdau,\tau)^{\circ}}\bigr),\]
which by Remark~\ref{rem: expect compact inductions to go to sheaves} %
is
a coherent sheaf (concentrated in degree zero) on $\cX_d$,
For each Serre
weight~$\underline{k}$ we set
\[\cF_{\underline{k}}:=\mathfrak{A}(\cInd_K^GF_{\underline{k}}).\] An
easy induction using~\cite[Lem.\ 4.1.1]{emertongeerefinedBM} shows that $\cF_{\underline{k}}$ is a
coherent sheaf (again, concentrated in degree zero) on the special fibre of~$\cX_d$.

Since $\mathfrak{A}$ is a functor, it follows from the definitions
that we have an equality in $K_0(\Coh(\cX_d))$
\numequation
  \label{eqn: Banach BM lifted to K0}[\cF_{\lambdau,\tau}\otimes_{\cO}k]=\sum_{\underline{k}}n_{\underline{k}}^\crys(\lambdau,\tau)\cdot
  [\cF_{\underline{k}}].
\end{equation}
As explained in Remark~\ref{rem: expect duality to imply maximal CM},
$\cF_{\lambdau,\tau}$ is maximally Cohen--Macaulay over its support
$\cX_d^{\crys,\lambdau,\tau}$, so that the support of
$\cF_{\lambdau,\tau}\otimes_\cO k$ is equal to
$Z_{\crys,\lambdau,\tau}$.  Again, an easy induction using~\cite[Lem.\
4.1.1]{emertongeerefinedBM} shows that the support of each
$\cF_{\underline{k}}$ is a cycle~$Z_{\underline{k}}$ in
$\Z[\cX_{d,\red}]$, so that Conjecture~\ref{conj: geometric BM} is an
immediate consequence of~\eqref{eqn: Banach BM lifted to K0}.

\subsection{Expectations: the analytic case} \label{subsec: analytic expectations}

Let $G={\rm GL}_d(F)$ considered as a $\mathbf{Q}_p$-analytic group.  Moreover, we assume in this section that the fixed finite extension $L$ of $\mathbf{Q}_p$ contains the normal closure of $F$. 
We want to formulate a conjecture parallel to Conjecture~\ref{conj: Banach functor} involving locally analytic representations and sheaves on $\mathfrak{X}_d$. 
In order to do so we first have to discuss derived categories in the context of locally analytic representations.

\subsubsection{Derived categories of locally analytic representations} We write $\mathcal{D}(G)$ for the distribution algebra of $G$ with coefficients in $L$. This distribution algebra is defined as the dual space of the space of locally analytic functions $G\rightarrow L$, see \cite{MR1887640} for the precise definition (in our case the field of coefficients will always be understood and we suppress it from the notation).

By a theorem of Schneider--Teitelbaum \cite[Corollary 3.3]{MR1887640} the category of locally analytic representations of $G$ on $L$-vector spaces of compact type with continuous $G$-morphisms is anti-equivalent to the category of separately continuous $\mathcal{D}(G)$-modules on nuclear Fr\'echet spaces with continuous $\mathcal{D}(G)$-module homomorphisms. 
This anti-equivalence is given by mapping an $L$-vector space $V$ of compact type to its strong dual $V'=V'_b$.
Under this anti-equivalence of categories the admissible locally
analytic $G$-representations (by definition) correspond to
coadmissible $\mathcal{D}(G)$-modules, i.e.\ to those $\mathcal{D}(G)$-modules $M$ that are coadmissible over the distribution algebra $\mathcal{D}(K)$ of the maximal compact subgroup $K={\rm GL}_n(\mathcal{O}_F)$ of $G$.
More precisely this means that $M$ is an inverse limit of finite type modules over certain Banach completions $\mathcal{D}_r(K)$ of $\mathcal{D}(K)$, compare \cite[\S3, 4]{MR1990669}.
In the following we will write ${\rm Rep}_L^{\rm an}G$ for the category of locally analytic $G$-representations on $L$-vector spaces of compact type. 

We note moreover that the $G$-action on a locally analytic representation $\pi$ automatically extends to a separately continuous action of the distribution algebra $\mathcal{D}(G)$, i.e.~the representation $\pi$ can be regarded as a $\mathcal{D}(G)$-module~\cite[Prop.~3.2]{MR1887640}.
Similarly to the category $D(\smG)$ we consider the derived category $D({\rm an}\,G)$ of complexes of locally analytic representations. 
The definition of this category (and of derived functors from this category to other derived categories) meets some difficulties, as we are working with topological objects in a homological algebra context. Problems like these are best dealt with by using the condensed structures and solid modules introduced by Clausen--Scholze, and a theory of locally analytic representations in the framework of solid modules was developed by Rodrigues Jacinto and Rodr\'iguez Camargo in \cite{solidlocan} and \cite{solidlocan2}.
We briefly outline\footnote{The reader who is not familiar with the formalism of solid modules may take the category  $D({\rm an}\,G)$ as a black box, and just view it as the derived category of ${\rm Rep}_L^{\rm an}G$ in an appropriate sense.} the construction of the ($\infty$-)category $D({\rm an}\,G)$:
The ($\infty$-)category $D_{L_\blacksquare}(L_\blacksquare[G])$ of $L_\blacksquare[G]$-modules in the derived category of $L_\blacksquare$-modules has a full subcateogry $D_{L_\blacksquare}(\mathcal{D}(G))$ consisting of the $\mathcal{D}(G)$-modules in that category, see \cite[Cor.4.3.4]{solidlocan} and \cite[Remark 1.2.1 (2)]{solidlocan2}. 
The category $D_{L_\blacksquare}(\mathcal{D}(G))$ has an endofunctor 
$$(-)^{\rm Rla}:D_{L_\blacksquare}(\mathcal{D}(G))\longrightarrow D_{L_\blacksquare}(\mathcal{D}(G))$$
that maps a (complex of) $\mathcal{D}(G)$ module(s) to its (derived) locally analytic vectors \cite[Def. 3.2.3]{solidlocan2}. 
This functor is a (derived) generalization of the functor of Schneider--Teitelbaum that maps an admissible Banach space representation of $G$ to its locally analytic vectors. 
Rodrigues Jacinto and Rodr\'iguez Camargo then define $D({\rm an}\,G)$ as the full subcategory of $D_{L_\blacksquare}(\mathcal{D}(G))$ consisting of complexes $\pi$ such that the natural map $\pi^{\rm Rla}\rightarrow \pi$ is an isomorphism \cite[Def. 3.3.1]{solidlocan2}.
As being derived locally analytic can be tested on cohomology groups \cite[Theorem A (1)]{solidlocan2} using non-derived locally analytic vectors, the (heart of the natural $t$-structure of the) category $D({\rm an}\,G)$ contains the category ${\rm Rep}_L^{\rm an}G$ of (non-derived) locally analytic representations. 
An application of the same result shows that $D({\rm an}\,G)$ contains $D({\rm Rep}^{\rm sm}\ G)$, the derived category of the category of smooth $G$-representations (though not as a full subcategory, as the ${\rm Ext}$-groups are not the same).

In order to discuss functors from $D({\rm an}\,G)$ to coherent sheaves on the stacks~$\mathfrak{X}_d$ (and their variants) we have to impose a finiteness condition on $D({\rm an}\,G)$ in order to be able to map to coherent sheaves (as opposed to solid quasi-coherent or ind-coherent sheaves, that we do not want to discuss here). In the following we will hence write 
$$D^b_{\rm f.p.}({\rm an}\,G)\subset D({\rm an}\,G)$$
for the full subcategory of those (derived) locally analytic representations that are expected to satisfy this condition. The first guess would be that $D^b_{\rm f.p.}({\rm an}\,G)$ consists exactly of the compact objects in $D({\rm an}\,G)$, but there are some reasons to expect that this is not the case. The correct guess seems to be that $D^b_{\rm f.p.}({\rm an}\,G)$ is the full subcategory of objects that are \emph{prim} with respect to the $6$-functor formalism on locally analytic
representations. \footnote{The reader can however approximately think of the category of compact objects in  $D({\rm an}\,G)$ instead.} 
This category contains the category $D^b_{\rm f.p.}({\rm Rep}^{\rm sm}\ G)$ of bounded complexes of smooth $G$-representations on $L$-vector spaces with finitely presented cohomology (though again not as a full subcategory).

\begin{remark}
The notation $D^b_{\rm f.p.}({\rm an}\,G)$ is rather non-standard, but we want to use a notation that is reminiscent of our notation in the Banach case. 
\end{remark}

 \subsubsection{Formulation of the conjecture}
Recall that the center $Z(\mathfrak{g}_L)$ of the enveloping algebra $U(\mathfrak{g}_L)$ of the Lie algebra $\mathfrak{g}_L$ of the $\mathbf{Q}_p$-analytic group $G$ embeds into $\mathcal{D}(G)$, \cite[Prop.3.7]{MR1887640}. This center can be identified with the global sections $\Gamma({\rm WT}_{d,L},\mathcal{O}_{{\rm WT}_{d,L}})$ of the space of Hodge--Tate--Sen weights of $\mathfrak{X}_d$ using the Harish--Chandra isomorphism.\footnote{Recall that in this section the group $G$ is ${\rm GL}_d$ and hence we can easily compare the center of the enveloping algebra with the space of functions on the space of Hodge--Tate--Sen weights. The case of an arbitrary reductive group is harder and would rely on similar constructions as in \cite[4.6]{inficharinfamilies}. But anyway, we did not define spaces of $(\varphi,\Gamma)$-modules with $G$-structures for arbitrary reductive groups.}

Further recall \cite[Theorem 5.13]{benzvi2020coherent},
\cite[Conjecture 4.5.1]{zhu2020coherent}, \cite[Conjecture
3.6]{hellmann2020derived} (see also Section~\ref{subsec: categorical LL
  l not p}) that for smooth representations we have a functor 
$$\mathfrak{A}^{\rm sm}_G:D^b_{\rm f.p.}({\rm Rep}^{\rm sm} G)\longrightarrow D^b_{\rm coh}({\rm WD}_{d,F}).$$
(For a general reductive group this functor is of course conjectural but for $G={\rm GL}_d(F)$ the functor can be constructed from the case of the Iwahori--Hecke algebra using type theory, see \cite{benzvi2020coherent}.)

We formulate a conjecture about locally analytic
representations and coherent sheaves on $\mathfrak{X}_d$ that roughly
parallels Conjecture~\ref{conj: Banach functor} in the case of smooth
mod $p$ representations and sheaves on $\mathcal{X}_d$.

\begin{conj}\label{EH:conjanalytic}
There exists an exact $Z(\mathfrak{g}_L)$-linear functor
$$\mathfrak{A}_G^{\rm rig}: {D}^b_{\rm f.p.}({\rm an}\,G)\rightarrow {D}^b_{\rm coh}(\mathfrak{X}_d)$$
satisfying the following conditions:
\begin{enumerate}
\item[(1)] {\em (}Compatibility with the smooth case.{\em )} Let
  $\underline{\xi}$ be a dominant algebraic character associated to a
  regular Hodge--Tate weight $\lambdau$, and write $W_{\lambdau}[1/p]$ for
  the corresponding irreducible algebraic representation of $G$ of
  highest weight~ $\xi$ as in Definition~\ref{defn: W lambda}.
Let ${\rm Rep}_{\Omega_\tau}^{\rm sm}G$ be the Bernstein block associated to an inertial type $\tau$ such that $\tau|_{I_{F'}}$ is trivial, for some finite extension $F'$ of $F$. 
After base change to $L'\supset F'$ there is a natural isomorphism 
\numequation\label{EHanconj smoothcomparison}
\mathfrak{A}_G^{\rm rig}( W_{\lambdau}[1/p] \, \otimes \text{--} ) \cong ({\rm pr_{WD}}^\ast(\mathfrak{A}_G^{\rm sm}(-))((-\rho')_{\sigma:F\hookrightarrow L}),
\end{equation}
where
\begin{align*}
{\rm pr_{WD}}:{\rm Fil}_{\lambdau}{\rm Mod}_{d,\varphi,N,\tau}\cong {\rm Fil}_{\lambdau}{\rm WD}_{d,F,\tau}&\longrightarrow {\rm WD}_{d,F,\tau}%
\end{align*}
is the canonical projections. 

\item[(2)] {\em (}Compatibility with parabolic induction.{\em )} Let $P\subseteq G$ be a parabolic subgroup with Levi quotient $M$ and write $\hat P$ respectively $\hat M$ for the corresponding dual groups. Then locally analytic parabolic induction $${\rm Ind}_{\bar P}^G(\text{--})^{\rm an}:D({\rm an.}\, M)\rightarrow D({\rm an.} G),$$ where $\bar P\subset G$ is the opposite parabolic, satisfies the compatibility 
\numequation\label{EHanconj parabinduction}\mathfrak{A}_G^{\rm rig}({\rm Ind}_{\bar P}^G(\text{--})^{\rm an})\cong \widetilde{\beta}_{\hat P,\ast}\circ \widetilde{\alpha}_{\hat P}^\ast(\mathfrak{A}_{M}^{\rm rig}(\text{--})([F:\Q_p]\rho'_M))
\end{equation}
whenever this formula makes sense, see {\em Section~\ref{EH subsection:paraboliccompatibility}} for a more precise formulation.
\end{enumerate}
\end{conj}
\begin{remark}\label{EHrem dualityonlocanreps}
We also expect a compatibility with duality similar to (3) of Conjecture \ref{conj: Banach functor}.
This would involve the definition of a \emph{Bernstein--Zelevinsky} type duality on $D({\rm an.}\, G)$. 
It should be possible to define such a duality by interpreting  $D({\rm an.}\, G)$ as the category of quasi-coherent sheaves on an analytic stack and using 6-functor formalisms for analytic stacks. However, we will not go into this direction in these notes. 
\end{remark}

\begin{rem}\label{EHrem definitionoftwists}
Let us comment on the twists by $-\rho'$ respectively by $\rho'_M$ that occur in (\ref{EHanconj smoothcomparison}) respectively (\ref{EHanconj parabinduction}).\\
(i) Let us, as usual, write $\rho$ for the half sum of the positive roots in the dual torus $\hat T$. To fix notations let us fix $\hat T$ to be the diagonal torus in ${\rm GL}_d$ and we canonically identify $\hat T=\mathbf{G}_m^d$. Then $$\rho=(\frac{d-1}{2},\frac{d-3}{2},\dots,\frac{-d+3}{2},\frac{-d+1}{2})$$ in general is not a character of $\hat T$, but only lies in $X^\ast(\hat T)\otimes\Q$. 
Instead we consider the character $\rho'=(0,-1,\dots, -(d-1))\in X^\ast(\hat T)$. 
Note that $\rho'$ is a shift of $\rho$ by an element of $X^\ast(\hat T)\otimes\Q$ that factors through the determinant and agrees with $\rho$ when evaluated on the coroots of $\hat T$.

If $M$ is a Levi subgroup of $G$, then we write $\rho_M$ for the half sum of the roots that occur in the $\hat T$ representation on ${\rm Lie}\,\hat U$, where $\hat U\subset \hat P$ is the unipotent radical.  
Then again we can apply a shift to $\rho_M$ to obtain a character $\rho'_M$ of $\hat M=\prod_{i=1}^r{\rm GL_{r_i}}$ that factors through the maximal torus quotient $S_{\hat M}=\prod_{i=1}^r\mathbf{G}_m$ such that $\rho'_M$ is trivial on the first $\mathbf{G}_m$ factor of $S_{\hat M}$.  \\
(ii) Let us now describe the twist in (\ref{EHanconj smoothcomparison}). 
The space $ {\rm Fil}_{\lambdau}{\rm WD}_{d,F,\tau}$ has a canonical projection to 
$${\rm GL}_d\backslash ({\rm Res}_{F/\mathbf{Q_p}}{\rm GL}_d)_L/P_{\lambdau}\cong {\rm GL}_d\backslash\big(\prod_{\sigma:F\hookrightarrow L}{\rm GL}_d/\hat B\big).$$
Here $P_{\lambdau}$ denotes the parabolic defined in Section~{\em\ref{EHsub:sub:deRham}} that agrees with the Weil restriction of the Borel, as $\lambdau$ is regular.
Then we can twist a sheaf on ${\rm Fil}_{\lambdau}{\rm WD}_{d,F,\tau}$ by the pullback of a line bundle on the product of flag varieties that we choose to be the line bundle on ${\rm GL}_d/\hat B$ associated with the character $-\rho'$ of $\hat T$ in each factor.
We point out that twisting with $(-\rho')_{\sigma}$ is compatible with
duality (in the spirit of the compatibility in Remark \ref{rem: duality explicated on crystalline stack}).
Again we have an auto-duality $\mathcal{D}_{\mathfrak{X}_d}$ on $D^b_{\rm coh}(\mathfrak{X}_d)$ induced by composing Grothendieck--Serre duality with the pullback of a map $\mathfrak{X}_d\rightarrow \mathfrak{X}_d$ induced by the main involution on ${\rm GL}_d$. 
If $\mathcal{F}$ denotes a complex of sheaves on ${\rm WD}_{d,F,\tau}$, then 
$$\mathcal{D}_{\mathfrak{X}_d}({\rm pr}_{\rm WD}^\ast(\mathcal{F})(-\rho')_{\sigma})={\rm pr}^\ast_{\rm WD}(\mathcal{D}_{{\rm WD}_{d,F,\tau}}(\mathcal{F}))(-\rho')_\sigma,$$
where $\mathcal{D}_{{\rm WD}_{d,F,\tau}}$ is defined the same way as $\mathcal{D}_{\mathfrak{X}_d}$.
 \\
(iii) Finally let us describe the twist in (\ref{EHanconj parabinduction}). The determinant induces  a canonical morphism $\mathfrak{X}_d\rightarrow \mathfrak{X}_1\rightarrow \ast/\mathbf{G}_m$ (see 
(\ref{EHeqn X1analytic}) for a description of $\mathfrak{X}_1$ and the canonical projection to $\ast/\mathbf{G}_m$). In particular we obtain a projection 
$$\mathfrak{X}_{\hat M}\rightarrow \ast/S_{\hat M}$$
and hence it makes sense to twist sheaves on $\mathfrak{X}_{\hat M}$ with a character of $S_{\hat M}$ (that we view as a line bundle on $\ast/S_{\hat M}$ and pull back to $\mathfrak{X}_{\hat M}$).
\end{rem}

There should be further compatibilities, namely a compatibility with the Banach case discussed in Section \ref{EH subsection: Banachcompatibility} below and a local-global compatibility that we will talk about in Section \ref{sec:eigenvar}.
We now comment a bit on the compatibilities in the conjecture.
\begin{rem}\label{EH rem:analyticconjecture}
\noindent \hfill\\(a) Property~(1) prescribes the functor $\mathfrak{A}_G^{\rm
  rig}$ on locally algebraic representations. This means that in some
sense we expect that the extension of the functor $\mathfrak{A}_G^{\rm
  sm}$ in smooth representations to locally algebraic representations
(which takes values in coherent sheaves supported on de Rham loci) can
be interpolated to a functor that takes values in coherent sheaves on
$\mathfrak{X}_{d}$ (without restrictions on the support). In
some sense this parallels the point of view that (Hecke eigenvalues on) $p$-adic automorphic forms interpolate (Hecke eigenvalues on) classical automorphic forms. \\
\noindent (b) The natural isomorphism in Property~(1) depends on the choice of the isomorphism
\[{\rm Fil}_{\lambdau}{\rm Mod}_{d,\varphi,N,\tau}\cong {\rm Fil}_{\lambdau}{\rm WD}_{d,F,\tau}\]
in Lemma~\ref{EH:lem:filtWDrep}, i.e.~on the choice of an embedding $F'\hookrightarrow L'$. 
In fact we expect that there is a variant of the functor $\mathfrak{A}_G^{\rm sm}$ that is a functor 
$$D^b_{\rm f.p.}({\rm Rep}^{\rm sm} G)\longrightarrow D^b_{\rm coh}\big(\bigcup_{F'/F} {\rm Mod}_{d,\varphi,N,F'/F}\big)$$
which agrees with $\mathfrak{A}_G^{\rm sm}$ after appropriate scalar extensions, see Remark \ref{EHrem: formofsmoothfunctor}. This version of $\mathfrak{A}_G^{\rm sm}$
 would allow for a smoother formulation of Property~(1) (i.e.~a version without an auxiliary scalar extension to a field containing $F'$). 
\\
\noindent (c) We cannot hope that the functor $\mathfrak{A}_G^{\rm rig}$ is fully faithful without putting any additional conditions on the source $D^b_{\rm f.p.}({\rm an}\,G)$.
 In the ${\rm GL}_1$-case this is discussed in Remark \ref{EH rem: ffGL1} below. But the phenomenon is also visible by looking at locally algebraic representations. The functor $\mathfrak{A}_G^{\rm sm}$ is (conjecturally) fully faithful. Given two smooth representations $\pi_1$ and $\pi_2$ (concentrated in degree zero) the module of homomorphisms ${\rm Hom}_G(\pi_1,\pi_2)$ is the same whether we compute it in the category of smooth representations or in the category of locally analytic representations (of course this is not true for higher Ext-groups). By full faithfulness this module coincides with ${\rm Hom}(\mathfrak{A}_G^{\rm sm}(\pi_1),\mathfrak{A}_G^{\rm sm}(\pi_2))$, but in general the space of homomorphisms between the analytifications of these sheaves is strictly larger.
It seems that this issue will not occur if we either add a \emph{finite slope} condition, or ask that the representations are admissible; compare the discussion in the ${\rm GL}_1$-case.
Unfortunately, simultaneously imposing the assumptions of being finitely presentated over $\mathcal{D}(G)$ and being admissible is rather restrictive, and would rule out a lot of interesting representations.
\\ 
\noindent (d) We expect that the functor $\mathfrak{A}_G^{\rm rig}$ can be defined on a larger category like ${D}({\rm an}\,G)$ instead of ${D}^b_{\rm f.p.}({\rm an}\,G)$ if we allow more general (complexes of) sheaves than just coherent sheaves. Similarly to the case of $p$-power torsion representations we expect that in this case admissibility of a representation $\pi$ will not ensure that $\mathfrak{A}_G^{\rm rig}(\pi)$ is a coherent sheaf without assuming some additional finiteness assumption on the locally analytic representation regarded as a module over $\mathcal{D}(G)$. In the $p$-torsion case this is discussed in Section~\ref{subsubsec: irredquots}, %
and one knows \cite[Thm.~1.1]{MR2608966} that,  when $G =\GL_2(F)$,
there are unitary Banach completions of locally algebraic $G$-representations
whose reduction contains any given irreducible admissible mod $p$ $G$-representation, 
so that we expect the mod $p$ phenomena of that discussion to have analogues
in characteristic zero.
\end{rem}

\begin{rem}\label{EH:rem functor is completed tensor}
We expect that, similarly to Remark \ref{rem: expect Linfty is a
  kernel}, the functor $\mathfrak{A}_{G}^{\rm rig}$ of Conjecture
\ref{EH:conjanalytic} should be given by
$$\pi\longmapsto \mathcal{L}_\infty\widehat{\otimes}^L_{\mathcal{D}(G)}\pi$$
for a certain derived completed tensor product $-\widehat{\otimes}^L_{\mathcal{D}(G)}-$  and a family $\mathcal{L}_\infty$ of $\mathcal{D}(G)$-modules over the stack $\mathfrak{X}_d$.  This \emph{derived completed} tensor product certainly only makes sense in the world of solid modules and solid locally analytic representations. 
We expect that, using the map $\pi_d$ from (\ref{EHgenfibEGstacktoanalytic}), the pullback of $\mathcal{L}_\infty$ to the generic fiber of $\cX_d$ should be given by
$$\pi_d^\ast \mathcal{L}_\infty\cong L_\infty\otimes_{\mathcal{O}_{\cX_d}[[G]]}\mathcal{D}(G,\cO_{\mathfrak{X}_d})$$
for a relative version $\mathcal{D}(G,\cO_{\mathfrak{X}_d})$ of the distribution algebra defined in similar terms as in \cite[3.1]{MR3623233} or \cite[1.3]{MR3874944}.

In the case of ${\rm GL}_2(\Q_p)$ we expect that for a point $x\in \mathfrak{X}_2(L)$, given by a $(\varphi,\Gamma)$-module $D$, the $\mathcal{D}(G)$-module $\mathcal{L}_\infty\otimes k(x)$ is the dual of the locally analytic representation associated to $D$ by Colmez's extension of the $p$-adic Langlands correspondence to the case of (not necessarily \'etale) $(\varphi,\Gamma)$-modules over $\mathcal{R}_F$, \cite{MR3522263}. More precisely: the family $\mathcal{L}_\infty$ should interpolate Colmez's $p$-adic Langlands correspondence for $(\varphi,\Gamma)$-modules over $\mathcal{R}_L$.
Using the object $L_\infty$ from Remark \ref{rem: expect Linfty is a kernel} we can construct $\mathcal{L}_\infty$ (in the case of ${\rm GL}_2(\Q_p)$) on the open substack of $\mathfrak{X}_2$, where the $(\varphi,\Gamma)$-modules are \'etale up to twist. The complement of this open subset consists only of trianguline $(\varphi,\Gamma)$-modules.
On the other hand Gaisin and Rodrigues Jacinto \cite{MR3874944} construct an analogue of $\mathcal{L}_\infty$ over (the regular part of) the stack $\mathfrak{X}_B$ parameterizing trianguline $(\varphi,\Gamma)$-modules (together with a triangulation).
At least for the time being it is not clear how to glue $\mathcal{L}_\infty$ from these two cases.
\end{rem}

\subsubsection{Compatibility with the Banach case}\label{EH subsection: Banachcompatibility}
We expect a strong link between the conjecture in the Banach case (Conjecture \ref{conj: Banach functor}) and the conjecture in the analytic case (Conjecture \ref{EH:conjanalytic}).
On the representation theoretic side, this compatibility should
involve  passage from a lattice $\Lambda$ in a Banach space
representation $V=\Lambda[1/p]$ %
(that we view as a pro-object in
$\Dfp^b(\smG)$) to the locally analytic vectors
(cf.~\cite[\S7]{MR1990669}) in $V$. On the side of (coherent) sheaves
on stacks of $(\varphi,\Gamma)$-modules, we consider the generic fiber
$\mathcal{F}_\eta$ of a pro-object in $D_{\coh}^b(\cX_d)$. %
This  (at
least heuristically, but see the discussion following Lemma
\ref{EH:lemtorsionanalyticcomp d=1} for a discussion in the case of
${\rm GL}_1$) gives rise to a complex of sheaves on $\cX_{d,\eta}^{\rm
  rig}$, and we compare this complex of sheaves to the pullback along (\ref{EHgenfibEGstacktoanalytic}) of a complex of sheaves on $\mathfrak{X}_d$. 
More precisely, for a pro-object $\Lambda \in \Dfp^b(\smG)$ with $V=\Lambda[1/p]$, the comparison should look like
\numequation\label{EHBanachanalyticcompatibility}
\mathfrak{A}(\Lambda)_{\eta}^{\rm rig}\cong \pi_d^\ast\mathfrak{A}_G^{\rm rig}(V^{\rm an}).
\end{equation}

This formula meets with some difficulties due to the finiteness conditions that we have imposed. We expect that these difficulties can be dealt with by using solid modules on both sides (i.e.~solid locally analytic representations and solid quasi-coherent sheaves) but we will not pursue this direction here. Instead we discuss some of the difficulties that one has to deal with.

\begin{enumerate}
\item[-] We have restricted ourselves to coherent sheaves on rigid analytic spaces (and rigid analytic Artin stacks), though there should be a more flexible ambient category of sheaves on these spaces (defined in terms of solid modules), compare Remark \ref{EHrem:sheavesnotsolid}. Though it is possible to put conditions on pro-coherent sheaves $\mathcal{F}$ on $\mathcal{X}_d$ that will assure that  $\mathcal{F}$ has a well-defined generic fiber $\mathcal{F}_\eta^{\rm rig}$ that is a coherent sheaf in $\mathcal{X}_{d,\eta}^{\rm rig}$, it is less clear which finiteness conditions on the pro-object $\Lambda$ of $D(\smG)$ would assure that $\mathfrak{A}(\Lambda)^{\rm rig}_\eta$ is an object of $D_{\rm coh}(\mathcal{X}_{d,\eta}^{\rm rig})$. Given a single representation $\Lambda$ (say in cohomological degree zero) the finiteness assumption should at least involve finiteness over $\cO[[G]]$.
\item[-] In our context of derived categories the passage to the locally analytic vectors $V\mapsto V^{\rm an}$ should be derived. It turns out that this is not necessary on the category of admissible Banach spaces representations. More precisely, if $V$ is an admissible Banach space representation (i.e.~its dual is finitely generated over $L[[K]]$), then $V^{\rm an}$ is an admissible locally analytic representation and its dual may be described as follows: By \cite[\S4]{MR1990669} there is a canonical flat map $\cO[[K]]\rightarrow \mathcal{D}(K)$ from the completed group ring $\cO[[K]]$ to the distribution algebra $\mathcal{D}(K)$. Given a finitely generated $\cO[[K]]$-module $\pi$  we can consider the extension of scalars $\pi\otimes_{\cO[[K]]}\mathcal{D}(K)$. If $\pi$ is equipped with a $G$-action (compatible with the $\cO[[K]]$-module structure), then $\pi\otimes_{\cO[[K]]}\mathcal{D}(K)$ naturally becomes a coadmissible $\mathcal{D}(G)$-module. 
As $O[[K]]\rightarrow \mathcal{D}(K)$ is flat, the functor $(-)^{\rm an}$ is exact on the category of admissible Banach space representations, and one can in fact deduce that $V^{\rm an}$ is even derived locally analytic, see \cite[Prop. 4.48]{solidlocan}, \cite[Theorem 2.2.3]{pan2020locallyanalytic}.
In general, i.e.~without invoking admissibility, we have to consider a derived variant of $V\mapsto V^{\rm an}$ (that is the derived locally analytic vectors of \cite{solidlocan}) that can only be defined in the framework of solid modules.
\item[-] As discussed in Remark \ref{EH rem:analyticconjecture}, the admissibility condition is not a good finiteness condition when dealing with the functors $\mathfrak{A}$ and $\mathfrak{A}_G^{\rm rig}$, i.e.~we should not expect that admissible representations are mapped to coherent sheaves in general. On the other hand imposing both finiteness conditions, admissibility and finiteness over $\cO[[G]]$, respectively $\mathcal{D}(G)$, is rather restrictive and excludes a lot of interesting representations. Hence it seems that a general formulation of the compatibility \eqref{EHBanachanalyticcompatibility} can only be formulated if we enlarge our ambient categories relying on the theory of solid modules.
\end{enumerate}

One instance of the compatibility (\ref{EHBanachanalyticcompatibility}) of the Banach and the analytic case should be the following. Given a regular Hodge--Tate weight~ $\lambdau$ and an inertial type~ $\tau$, Conjecture \ref{conj: Banach functor} (4) predicts that $$\mathcal{F}:=\mathfrak{A}\Bigl(\bigl(\widehat{\Pi(\lambdau,\tau)^{\circ}}\Bigr)$$
is a pro-coherent sheaf concentrated in degree $0$ and supported on the $p$-adic stack $\cX_d^{{\rm crys},\lambdau,\tau}$. We expect that its pullback along any morphism ${\rm Spf}\,A\rightarrow \cX_d^{{\rm crys},\lambdau,\tau}$ for a $p$-adically complete $\cO$-algebra topologically of finite type over $\cO$ is given by a finitely generated $A$-module. In this case there is a well-defined generic fiber $\mathcal{F}_\eta^{\rm rig}$ of $\mathcal{F}$ which is a coherent sheaf on $(\cX_d^{{\rm crys},\lambdau,\tau})_\eta^{\rm rig}$.
We then expect that there is
an isomorphism
\numequation\label{EHeqnlocalgcompatibility}
\mathcal{F}_\eta^{\rm rig}\cong \pi_d^\ast\Bigl(\mathfrak{A}_G^{\rm rig}\big(\Pi(\lambdau,\tau)\big)\Bigr).
\end{equation}
As already indicated, this is presumably an instance of %
the conjectural compatibility~(\ref{EHBanachanalyticcompatibility}),
but a fuller understanding of exactly how 
would seem to require a description of 
the subspace of locally analytic vectors in~$\widehat{\Pi(\lambdau,\tau)^{\circ}}.$
The problem of giving such a description is an interesting question in its own right,
which we unfortunately do not currently know the answer to.

Finally let us point out that Colmez \cite[Rem.\ 0.2]{MR3908766}
has described the Jordan--H\"older factors of the locally analytic vectors $\Pi(D)^{\rm an}$ in the Banach space representation $\Pi$ associated to an \'etale $(\varphi,\Gamma)$-module $D$ of rank $2$ under the $p$-adic local Langlands correspondence for ${\rm GL}_2(\Q_p)$.
It would be interesting to predict where the representations in \cite[Corollaire 0.4]{MR3908766} are mapped to under the functor $\mathfrak{A}_G^{\rm rig}$. In the trianguline case this is basically covered by the expected compatibilities with the smooth case and with parabolic induction. But for the remaining representations we do not know how to (conjecturally) describe the image.

\subsubsection{Compatibility with parabolic induction}\label{EH subsection:paraboliccompatibility}
We fix a Borel subgroup $B\subseteq G$ with split maximal torus $T$. Further let $P\supseteq B$ be a parabolic subgroup with Levi $M$.
Let $(\pi,V)$ be a locally analytic representation of $M$ on an
$L$-vector space $V$ of compact type. The locally analytic parabolic
induction of~$\pi$ is defined by $$({\rm Ind}_P^G\pi)^{\rm an}=\{f: G\rightarrow V\ \text{locally analytic}\mid f(gp)=\pi(p)^{-1}f(g)\ \text{for all}\ p\in P,g\in G\}.$$
In particular we obtain functors 
\begin{align*}
{\rm Ind}_P^G(-)^{\rm an}:{\rm Rep}_L^{\rm an}M&\longrightarrow {\rm Rep}_L^{\rm an}G,\\
{\rm Ind}_P^G(-)^{\rm an}:D({\rm an}\,M)&\longrightarrow D({\rm an}\,G),
\end{align*}
where by abuse of notation we use the same notation for the functor on the (informally defined) derived categories. In order to make the discussion less involved we focus on the case $P=B$ in the following.
Let us write 
$$s_G:D({\rm an}\,T)\longrightarrow D({\rm an}\,T)$$
for the endofunctor\footnote{We index this shift and twist by $G$ instead of $T$, as the shift and twist clearly depends on $G$, not just on $T$.} that twists a $T$-representation by 
\numequation\label{EHeqn:twistedparabolicinduction}
\delta_B\cdot (1,\dots, (\varepsilon\circ{\rm rec})^{i-1},\dots,(\varepsilon\circ{\rm rec})^{d-1}).
\end{equation}
(where $\delta_B$ is the smooth modulus character of $B$ and $\varepsilon$ is the cyclotomic character and ${\rm rec}$ is the isomorphism of local class field theory). 

\begin{expectation}\label{EHexp:parabolicinudction}
There are isomorphisms
\begin{align*}
\mathfrak{A}_G^{\rm rig}\big( {\rm Ind}_{\overline{B}}^G(\pi)^{\rm an}\big) & \cong \tilde{\beta}_{\hat B,\ast}\tilde\alpha^\ast_{\hat B}\big(\mathfrak{A}_T^{\rm rig}(s_G(\pi)[[F:\mathbf{Q}_p]\rho'])\big)%
\end{align*}
functorial in $\pi$,
whenever $\pi$ is a representation such that ${\rm Ind}_{\overline{B}}^G(\pi)^{\rm an}$ lies in $D^b_{\rm f.p.}({\rm an.}\,G)$ and such that the support of $\tilde \alpha_{\hat B}^\ast(\mathfrak{A}_T^{\rm rig}(\pi))$ is proper over $\mathfrak{X}_d$.
\end{expectation}

\begin{rem}
The appearance of the twist with the character \eqref{EHeqn:twistedparabolicinduction} is due to the chosen normalizations. In fact the normalizations we choose depend on the normalization of the functor $\mathfrak{A}_G^{\rm sm}$, or (which is more or less the same), the normalization of the isomorphism in Proposition \ref{EHBernsteiniso}, compare Remark \ref{EHremBernsteiniso}.
In fact here we use a different normalization than in \cite{hellmann2020derived}, where compatibility with parabolic induction is stated without involving a twist (but using normalized parabolic induction).
The normalization in \cite{hellmann2020derived} corresponds to the
unitary normalization of the local Langlands correspondence, rather than the non-unitary normalization that we use here. The price we have to pay is the appearance of the twist in Expectation \ref{EHexp:parabolicinudction}.
\end{rem}

Let us stress that locally analytic parabolic induction should not be expected to preserve the finiteness conditions, i.e~for $\pi\in D^b_{\rm f.p.}({\rm an.}\,T)$ we do not expect that ${\rm Ind}_{\overline{P}}^G(\pi)^{\rm an}$ lies in $D^b_{\rm f.p.}({\rm an.}\, G)$. 
However this should be true if $\pi$ is of bounded slope. We expect that the phenomenon is directly related to the question whether the support of $\alpha_{\hat B}^\ast(\mathfrak{A}_T^{\rm rig}(\pi))$ is proper over $\mathfrak{X}_d$ or not.
As already remarked above we expect that there is a variant of $\mathfrak{A}_{G}^{\rm rig}$ defined on $D({\rm an.}\, G)$ that takes values in a larger category of sheaves on $\mathfrak{X}_d$ (most probably: solid quasi-coherent sheaves). In this context one should have a full 6-functor formalism for these sheaves at hand and the compatibility with parabolic induction \ref{EHexp:parabolicinudction}  should rather be formulated with $\tilde\beta_{\hat B,!}$ instead of $\tilde\beta_{\hat B,\ast}$ (of course we don't see the difference if we assume that the support of the corresponding sheaf is proper).

We mention that Orlik and Strauch \cite{MR3264764} have described the (topologically) irreducible subquotients of some locally analytic principal series representations. More precisely, they describe the irreducible subquotients of  $({\rm Ind}_B^G\delta)^{\rm an}$ for locally algebraic characters $\delta:T\rightarrow L^\times$. 
The main tool of their study is the construction of certain bi-functors $\mathcal{F}_P^G(-,-)$ that map a pair $(M,\pi)$ consisting of a Lie-algebra representation $M$ in the BGG category $\mathcal{O}$ (such that the action of ${\rm Lie}\, P$ on $M$ can be integrated to a $P$-action) and a smooth representation $\pi$ of the Levi of $P$ to a locally analytic representation of $G$. We recall part of their construction. In fact it makes sense to speculate that the construction of Orlik--Strauch extends to functors
  \begin{align*}
\mathcal{F}_P^G:{\rm Rep}_L^{\rm an}M&\longrightarrow {\rm Rep}_L^{\rm an}G,\\
\mathcal{F}_P^G:D({\rm an}\,M)&\longrightarrow D({\rm an}\,G),
\end{align*}
that behave similarly to parabolic induction, see Expectation \ref{EH:existenceofFBG} below.

Again we focus on the case $P=B$. In the following we write $U(\mathfrak{g}_L)$ for the universal enveloping algebra of $\mathfrak{g}_L=({\rm Lie}\,G)\otimes_{\mathbf{Q}_p}L$ (note that here we consider $G$ as a $\mathbf{Q}_p$-analytic group and hence $\mathfrak{g}_L$ is isomorphic to the product $\prod_{\tau:F\hookrightarrow L}{\mathfrak{gl}}_d$). 
Similarly, we write  $\mathfrak{t}_L$ and $\mathfrak{b}_L$ for the Lie algebras over $L$ of the $\Q_p$-analytic groups $T$ respectively $B$. 
 Let us write $\mathcal{O}(U(\mathfrak{g}_L))$ for the BGG category $\mathcal{O}$ of finitely generated $U(\mathfrak{g}_L)$-modules on which the Cartan $\mathfrak{t}_L$ acts semi-simply and which are locally finite as $U(\mathfrak{b}_L)$-modules. 
We write $\mathcal{O}(U(\mathfrak{g}_L))^{\rm alg}$ for the subcategory of objects $M$ on which the action of $\mathfrak{t}_L$ is algebraic. On this subcategory one can integrate the $\mathfrak{b}_L$-action to obtain a $B$-action on $M$.  
Then every object $M$ of $\mathcal{O}(U(\mathfrak{g}_L))^{\rm alg}$ can be regarded as a module over the sub-$L$-algebra $\mathcal{D}(\mathfrak{g}_L,B)\subseteq \mathcal{D}(G)$ generated by $U(\mathfrak{g}_L)$ and $\mathcal{D}(B)$. 
We obtain an exact \emph{contravariant} functor
\[F_B^G(-):\cO(U(\mathfrak{g}_L))^{\rm alg}\longrightarrow {\rm Rep}^{\rm an}G\]
given by 
$$V\longmapsto (\mathcal{D}(G)\otimes_{\mathcal{D}(\mathfrak{g},B)}V)';$$
compare \cite[3.4]{MR3264764}. 
\begin{rem}
On the dual side, in terms of $\mathcal{D}(G)$-modules, the locally analytic parabolic induction is given by
\[\big(({\rm Ind}_B^G\delta)^{\rm an}\big)'_b=\mathcal{D}(G)\otimes_{\mathcal{D}(B)} (L_\delta)'_b,\]
where $L_\delta$ denotes the $1$-dimensional locally analytic representation of $T$ on $L$ via the character $\delta$, see \cite[Lemma 2.3]{MR3264764}. 
Now let $\xi$ be an algebraic character of $T$. The Verma module $M(\xi)=U(\mathfrak{g}_L)\otimes_{U(\mathfrak{b}_L)}L_\xi$ is an object of $\mathcal{O}(U(\mathfrak{g}_L))^{\rm alg}$ and we find that 
\[\mathcal{D}(G)\otimes_{\mathcal{D}(\mathfrak{g}_L,B)}M(\xi)=\mathcal{D}(G)\otimes_{\mathcal{D}(B)}L_\xi.\]
In particular, for an algebraic character $\xi$ we find $$F_B^G(M(-\xi))={\rm Ind}_B^G(z^\xi)^{\rm an}.$$
\end{rem}

In fact Orlik and Strauch \cite[4.4]{MR3264764} show that the functor $F_B^G(-)$ can in some sense be twisted by a smooth character of $T$: for any choice $\delta_{\rm sm}$ of a smooth character of $T$ there is an exact contravariant functor
$$F_B^G(-,\delta_{\rm sm}):\mathcal{O}(U(\mathfrak{g}_L))^{\rm alg}\longrightarrow {\rm Rep}^{\rm an}G$$
such that for an algebraic character $\xi$ the representation $F_B^G(M(-\xi),\delta_{\rm sm})$ is the locally analytic induced representation $({\rm Ind}_B^G z^{\xi}\delta_{\rm sm})^{\rm la}$. 

We point out the relation with locally algebraic representations whose
smooth part is parabolically induced. Fix an integral, dominant and
regular Hodge--Tate weight $\lambdau$ and write $\xi$ for the
corresponding character of $\mathfrak{t}_L$ (via Definition~\ref{defn: lambda versus xi}), and $\mathcal{O}(U(\mathfrak{g}_L))_{(\xi)}\subseteq \cO(U(\mathfrak{g}_L))$ for the block of the category $\mathcal{O}$ containing the Verma module $M(\xi)$. Then $\mathcal{O}(U(\mathfrak{g}))_{(\xi)}\subseteq \cO(U(\mathfrak{g}))^{\rm alg}$.
As $\lambdau$ is regular the character $\xi$ is automatically dot-regular (i.e.~the Weyl group orbit $W\cdot\xi$ of $\xi$ under the dot-action $w\cdot\xi=w(\xi+\rho)-\rho$ consists of $|W|$ elements).
The simple objects of $\mathcal{O}(U(\mathfrak{g}_L))_{(\xi)}$ are the simple quotients 
$L(w\cdot\xi)$ of the Verma modules $M(w\cdot\xi)$. Moreover,
$L(\xi)=W_{\lambdau}[1/p]$ is the irreducible algebraic representation
of $G$ defined in Definition~\ref{defn: W lambda}. %
\begin{lem}\label{EH:lemOSonGreps}
There is a canonical isomorphism
$$F_B^G(L(-\xi),\delta_{\rm sm})\cong L(\xi)\otimes_L ({\rm Ind}_B^G\delta_{\rm sm})^{\rm sm}$$
with the locally algebraic representation given by the tensor product of the irreducible algebraic representation $L(\xi)$ and the smooth parabolic induction of the smooth character $\delta_{\rm sm}$.
\end{lem}

Recall that the BGG category has an internal notion of duality $M\mapsto M^\vee$ (usually referred to as BGG duality) that is roughly given by mapping $M$ to the direct sum of the duals of its (finite dimensional) $\mathfrak{t}$-eigenspaces and using the main involution on $\mathfrak{g}$ (in order to pass from left to right modules and to obtain a $U(\mathfrak{g}_L)$-module that is locally $U(\mathfrak{b}_L)$-finite instead of $U(\overline{\mathfrak{b}}_L)$-finite, where $\overline{\mathfrak{b}}$ denotes the Lie-algebra of the opposite Borel $\overline{B}$). For a simple module $M$ in $\cO(U(\mathfrak{g})_L)$ we have $M^\vee\cong M$, but as $(-)^\vee$ is contravariant the functor exchanges subobjects and quotients. 
Let $\delta= z^{\xi}\delta_{\rm sm}$ be a locally algebraic character with smooth part $\delta_{\rm sm}$ and algebraic part $\xi$. We then define the representation 
\numequation\label{EH:OSoflocalg}
\mathcal{F}_B^G(\delta)=F_B^G(M(-\xi)^\vee,\delta_{\rm sm}). 
\end{equation}

\begin{rem}\label{EHrem localgquotofOS}It follows from \cite{MR3264764} that ${\rm Ind}_B^G(\delta)^{\rm an}$ and $\mathcal{F}_B^G(\delta)$ have the same irreducible subquotients (even with the same multiplicities). Moreover, note that Lemma \ref{EH:lemOSonGreps} together with the fact that there is a canonical surjection $M(-\xi)^\vee\rightarrow L(-\xi)$ implies that for dominant $\xi$ the representation $\mathcal{F}_B^G(z^{\xi}\delta_{\rm sm})$ has the locally algebraic representation $L(\xi)\otimes_L ({\rm Ind}_B^G\delta_{\rm sm})^{\rm sm}$ as a quotient, whereas this locally algebraic representation appears as a subrepresentation in ${\rm Ind}_B^G(z^\xi \delta_{\rm sm})^{\rm an}$. 
\end{rem}
We expect that the Orlik--Strauch functors can be used to define a variant of parabolic induction (whose derived versions of course will depend on the definition of the categories $D({\rm an}\,T)$ and $D({\rm an}\,G)$).
\begin{expectation}\label{EH:existenceofFBG}
There exist natural functors 
\begin{align*}\label{EH:existenceofFBG}
\mathcal{F}_B^G:{\rm Rep}_L^{\rm an}T&\longrightarrow {\rm Rep}_L^{\rm an}G,\\
\mathcal{F}_B^G:D({\rm an}\,T)&\longrightarrow D({\rm an}\,G),
\end{align*}
such that for a locally algebraic character $\delta$ the representation $\mathcal{F}_B^G(\delta)$ is given by \eqref{EH:OSoflocalg}.
\end{expectation}

\begin{rem}\label{EHrem OSdualtoInd}
It seems reasonable to speculate that there are internal dualities $\mathbf{D}_G$ resp. $\mathbf{D}_T$ of $D({\rm an}\,G)$ resp. $D({\rm an}\,T)$ (similar to the dualities discussed in Section \ref{sec:expectations Banach case}) such that (up to shift)
$$\mathcal{F}_B^G(-)=\mathbf{D}_G\circ {\rm Ind}_B^G(-)^{\rm an} \circ \mathbf{D}_T.$$
The restriction of $\mathbf{D}_G$ to smooth representations should be given by the Zelevinsky involution, and the above expectation suggests that $\mathbf{D}_G$ should be related to Serre-duality under the hypothetical functor $\mathfrak{A}_G^{\rm rig}$, compare Remark \ref{EHrem dualityonlocanreps}. In the case of principal series representations such a duality is defined and discussed by Strauch and Wu in \cite{strauch2025bernsteinzelevinskydualitylocallyanalytic}. 
\end{rem}

\begin{remark}
It was explained to us by Juan Esteban Rodríguez Camargo that it is possible define a generalization of the functor $\mathcal{F}_B^G$ on a category of $(\mathfrak{g},B)$-modules using the theory of analytic stacks developed by Clausen--Scholze. This generalization, and a six functor formalism on analytic stacks, can be used to construct the functors expected in Expectation \ref{EH:existenceofFBG} and to prove the expected compatibilities with duality in Remark \ref{EHrem OSdualtoInd}. As we do not intend to discuss the theory from the point of view of analytic stacks (due to many reasons, including personal insufficiencies in the handling of analytic stacks) we stick to the formulation in the expectation above. 
\end{remark}

Using the conjectural functor $\mathcal{F}_B^G$ we can formulate a variant of Expectation \ref{EHexp:parabolicinudction}. This expectation is motivated by the properties of Bezrukavnikov's functor in Theorem \ref{EH:thmBezfunctor} below.
\begin{expectation}\label{EHexp:OrlikStrauchcompatibility}
There are isomorphisms 
\begin{align*}
\mathfrak{A}_G^{\rm rig}\big( \mathcal{F}_{\overline{B}}^G(\pi)\big) & \cong \tilde{\beta}_{\hat B,\ast}\tilde\alpha^!_{\hat B}\big(\mathfrak{A}_T^{\rm rig}(s_G(\pi)[-[F:\Q_p]\rho'])\big)%
\end{align*}
functorial in $\pi$, whenever 
whenever $\pi$ is a representation such that $\mathcal{F}_{\overline{B}}^G(\pi)$ lies in $D^b_{\rm f.p.}({\rm an.}\,G)$ and such that the support of $\tilde\alpha_{\hat B}^!(\mathfrak{A}_T^{\rm rig}(\pi))$ is proper over $\mathfrak{X}_d$.
\end{expectation}

\begin{rem}
Let us consider the case $F=\Q_p$ and $d=2$ and recall from Remark \ref{EH:rem functor is completed tensor} that we expect $\mathfrak{A}_G^{\rm rig}$ to be of the form 
$$\pi\mapsto \pi\widehat{\otimes}^L_{\mathcal{D}(G)}\mathcal{L}_\infty,$$
where $\mathcal{L}_\infty$ is a family of $\mathcal{D}(G)$-modules over $\mathfrak{X}_2$. For a point $x\in\mathfrak{X}_2(L)$ represented by a $(\varphi,\Gamma)$-module $D$ we expect that 
$$\mathcal{L}_\infty\otimes k(x)=(\Pi(D)^{\rm an})'.$$
Here if~ $D$ is \'etale, $\Pi(D)^{\rm an}$ denotes the locally analytic vectors in the
representation corresponding to~ $D$ via the usual $p$-adic local Langlands correspondence for ${\rm GL}_2(\Q_p)$ (compare the discussion in section \ref{subsubsec:p-adic LL}), %
and if $D$ is not \'etale we use Colmez's extension \cite{MR3522263}
of the correspondence.
Given a character $\delta:T\rightarrow L$,  we are lead to expect that
\begin{align*}
H^0\big(\mathfrak{A}_G^{\rm rig}({\rm Ind}_{\bar B}^G (\delta))\big)\otimes k(x)&\cong\big({\rm Hom}_G({\rm Ind}_{\bar B}^G (\delta),\Pi(D)^{\rm an})\big)',\\
H^0\big(\mathfrak{A}_G^{\rm rig}(\mathcal{F}_{\bar B}^G (\delta))\big)\otimes k(x)&\cong\big({\rm Hom}_G(\mathcal{F}_{\bar B}^G (\delta),\Pi(D)^{\rm an})\big)'.
\end{align*}
The dimensions of the vector spaces on the right hand side can be
determined explicitly, as Colmez has given a precise description of
the locally analytic representations $\Pi(D)^{\rm an}$ (see
\cite{MR3444228,MR3522263}), and in particular they are zero unless $D$ is trianguline.
Using Expectations \ref{EHexp:parabolicinudction} and \ref{EHexp:OrlikStrauchcompatibility} the computation of these fiber dimensions precisely matches the computation of fiber dimensions in Theorem \ref{EHtheo triangvar for GL2Qp} (v).
\end{rem}

\subsubsection{Bezrukavnikov's functor on the BGG category}
\label{subsec:bezrukavnikov-orlik-strauch}
We elaborate a bit further on the Orlik--Strauch functor and compatibility with parabolic induction. As on the Galois side parabolic induction is defined in terms of the compactification $\overline{\mathfrak{X}}_B$ of $\mathfrak{X}_B$ and as we can explicitly control the geometry of $\overline{\mathfrak{X}}_B$ at certain regular points (Theorem \ref{EH:Theorem:localmodelB=P}), we can conjecturally link the evaluation of the functor $\mathfrak{A}_G^{\rm rig}$ on certain Orlik--Strauch representations to a functor on the BGG category constructed by Bezrukavnikov \cite{MR3502096}.

Fix $\lambdau$ and $\xi$ as above. As every object in
$\mathcal{O}(U(\mathfrak{g}_L))_{(\xi)}$ is an extension of objects
that are (sub-)quotients of Verma modules $M(w\cdot\xi)$ for $w\in W$,
the expectations on the functor $\mathfrak{A}^{\rm rig}_G$ imply that given $\delta_{\rm sm}$, there should be a functor
$$\mathfrak{A}^{\rm rig}_G\circ F_{\overline{B}}^G((-)^\vee,\delta_{\rm sm}):\cO(U(\mathfrak{g}_L))_{(\xi)}\longrightarrow D^b_{\rm coh}(\mathfrak{X}_d)$$
such that the sheaves in the image of this functor are supported on the closed substack $\mathfrak{X}_{d,(\lambdau,\delta_{\rm sm}){\rm -tri}}$ defined just before Theorem~\ref{EH:Theorem:localmodelB=P}.
We use an additional duality $(-)^\vee$ in this functor to make the functor covariant.
Using the local description of the stack $\mathfrak{X}_{d,(\lambdau,\delta_{\rm sm}){\rm -tri}}^{\widehat{}}$ in Theorem \ref{EH:Theorem:localmodelB=P} we can construct this functor, relying on work of Bezrukavnikov.

Let us write 
$$H=({\rm Res}_{F/\mathbf{Q}_p}{\rm GL}_d)_L=\prod_{\sigma:F\hookrightarrow L}{\rm GL}_d$$ for the remainder of this section. 
We note that $\mathfrak{g}_L\cong {\rm Lie} H=\mathfrak{h}$
as Lie algebras over~$L$. In particular $\mathcal{O}(U(\mathfrak{g}_L))$ and $\mathcal{O}(U({\mathfrak{h}}))$ are isomorphic. 
Recall from Section~\ref{subsec:localmodels} (where we write $G$ instead of $H$) that we write 
$$\overline{X}_w\subseteq \widetilde{{\mathfrak{h}}}\times_{{\mathfrak{h}}}\widetilde{{\mathfrak{h}}}$$
for the irreducible component that dominates the closure of the Schubert cell $ H (1,w)\subseteq H/B_H\times H/B_H$,
where $B_H\subseteq H$ is the Weil restriction of the upper triangular matrices in ${\rm GL}_d$.
 
Moreover, we note that there is a canonical projection
$$\omega: \tilde{{\mathfrak{h}}}\times_{{\mathfrak{h}}}\tilde{{\mathfrak{h}}}\xrightarrow{{\rm pr}_1} \tilde{{\mathfrak{h}}}\longrightarrow \widehat{\mathfrak{t}},$$
(where $\widehat{t}$ is the dual Lie algebra of $\mathfrak{t}_L$)
and we write $\overline{X}_{w,0}$ for the intersection of $\overline{X}_w$ with the preimage of $0$. 
Finally note that the zero section gives rise to a copy 
\begin{align*}
({\rm Res}_{F/\mathbf{Q}_p}{\rm GL}_d/P_{\lambdau})_L\times({\rm Res}_{F/\mathbf{Q}_p}{\rm GL}_d/P_{\lambdau})_L&=H/B_H\times H/B_H\\ &\subseteq \overline{X}_{w_0,0}\subseteq \overline{X}_{w_0}\subseteq
\tilde{{\mathfrak{h}}}\times_{{\mathfrak{h}}}\tilde{{\mathfrak{h}}}.
\end{align*}
Here $P_{\lambdau}$ denotes the parabolic subgroup associated to $\lambdau$ as in Section \ref{EHsub:sub:deRham} which agrees with the Borel $B$ as we assume $\lambdau$ to be regular.
The following functor has been constructed by Bezrukavnikov in \cite[11.4]{MR3502096}, see also \cite[Theorem 13.0.1]{bezrukavnikovkanstrup}.
\begin{theorem}[Bezrukavnikov]\label{EH:thmBezfunctor}
There is an exact functor
$$F_{\xi}:\mathcal{O}(U(\mathfrak{g}_L))_{(\xi)}\longrightarrow {\rm Coh}^{H}((\tilde{{\mathfrak{h}}}\times_{{\mathfrak{h}}}\tilde{{\mathfrak{h}}})^{\widehat{}}_0)$$
such that \\
\noindent \hfill\\(i) For $w\in W$ the image of the of the Verma module
$M(ww_0\cdot\xi)$ is the dualizing sheaf of $\bar{X}_w$: %
$$F_\xi(M(ww_0\cdot\xi))=\omega_{\bar{X}_w}.$$
\noindent \hfill\\(ii) For $w\in W$ the image of the of the dual Verma module
$M(ww_0\cdot\xi)^\vee$ is the structure sheaf of $\bar{X}_w$: %
$$F_\xi(M(ww_0\cdot\xi)^\vee)=\cO_{\bar{X}_w}.$$
\noindent (iii) The image of the irreducible algebraic representation $L(\xi)$ is the line bundle %
on 
$$({\rm Res}_{F/\mathbf{Q}_p}{\rm GL}_d/P_{\lambdau})_L\times({\rm Res}_{F/\mathbf{Q}_p}{\rm GL}_d/P_{\lambdau})_L=H/B_H\times H/B_H\subseteq \overline{X}_{w_0,0}$$
associated to the character $((-\rho')_\sigma,(-\rho')_\sigma)$ (compare Remark \ref{EHrem definitionoftwists} for the definition of $\rho'$).
\end{theorem}
We expect that, using the local description of the stack $\mathfrak{X}_{d,(\lambdau,\delta_{\rm sm}){\rm -tri}}^{\widehat{}}$ via the map
$$f_{\lambdau,\delta_{\rm sm}}:\mathfrak{X}_{d,(\lambdau,\delta_{\rm sm}){\rm -tri}}^{\widehat{}}\longrightarrow (\tilde{\mathfrak{h}}\times_{\mathfrak{h}}\tilde{\mathfrak{h}})^{\widehat{}}_0$$
of Theorem \ref{EH:Theorem:localmodelB=P}, the pullback of Bezrukavnikov's functor gives rise to the restriction of the conjectural functor $\mathfrak{A}_G^{\rm rig}$ via the functor of Orlik--Strauch:
\begin{expectation}\label{EH:Bezcompatibility}
Let $\varphi_1,\dots,\varphi_d\in L^\times$ such that $\varphi_i/\varphi_j\neq 1,q$ for $i\neq j$ and let $\delta_{\rm sm}$ be as in {\em \ref{EH:eqndeltasm}}. 
Then the functor
$$\mathfrak{A}^{\rm rig}_G\circ \mathcal{F}_B^G(-,\delta_{\rm sm}):\cO(U(\mathfrak{g}_L))_{(\xi)}\longrightarrow D^b_{\rm coh}(\mathfrak{X}_d)$$
factors through $\mathfrak{X}_{d,(\lambdau,\delta_{\rm sm}){\rm -tri}}^{\widehat{}}$ and is given by the composition 
$$(f_{\lambdau,\delta_{\rm sm}}^\ast \circ F_\xi(-))([F:\Q_p]\rho').$$
\end{expectation}
\begin{rem}
Here the twist by $\rho'$ makes sense on $\mathfrak{X}_{d,(\lambdau,\delta_{\rm sm}){\rm -tri}}^{\widehat{}}$, compare the discussion following Theorem \ref{EH:Theorem:localmodelB=P}.

Note that using Lemma \ref{EH:lemOSonGreps} and the fact that $f_{\lambdau,\delta_{\rm sm}}$ is formally smooth with relative dualizing sheaf the line bundle associated to $-2\rho$ this expectation is compatible with the twist by $(-\rho')_\sigma$ in (\ref{EHanconj smoothcomparison}) and with the twists in Expectation \ref{EHexp:parabolicinudction} and \ref{EHexp:OrlikStrauchcompatibility}.
\end{rem}

\subsubsection{The trianguline Breuil--M\'ezard conjecture}
\label{subsec:triang-breu-mezard}

We describe a shadow of the compatibility with parabolic induction and of the compatibility with smooth categorical Langlands in Conjecture \ref{EH:conjanalytic}. This shadow can be regarded as a version of the Breuil--M\'ezard conjecture (Conjecture~\ref{conj: geometric BM}) for locally analytic principal series representations that involves Grothendieck groups rather than groups of cycles (compare the discussion in Section \ref{subsubsec: geometric BM Banach}). In certain regular cases, a version of this conjecture involving cycles was proven in \cite{MR3623233}, see Theorem \ref{BHS:locanBM} below.

As above we write $G={\rm GL}_d(F)$ and $T\subseteq B\subset G$ for the diagonal torus respectively the upper triangular Borel. We keep the assumption that $L$ is large enough.
Note that we may regard $G$ as the group of $\Q_p$-valued points of $H={\rm Res}_{F/\Q_p}{\rm GL}_d$ and similarly $T$ and $B$ are the $\Q_p$-valued points of $T_H$ and $T_B$ the Weil restrictions of the diagonal torus respectively of the upper triangular Borel. 
We fix a regular Hodge--Tate weight $\lambdau$ and as usual write $\xi$ for the corresponding shifted character (which is a character of $T_{H}$ defined over $\bar F$). We regard $\xi$ as an (algebraic) character $z^\xi:T\rightarrow L^\times $ (compare \eqref{EH:defofztau}) and fix a smooth character ${\delta}_{\rm sm}:T\rightarrow L^\times$. 

Let us write $W\cong \prod_{K\hookrightarrow L}\mathcal{S}_n$ for the absolute Weyl group of $H$ and $W_0\cong \mathcal{S}_n\subseteq W$ for the relative Weyl group of $H$. Then for every $w\in W$ and $w'\in W_0$ we obtain well-defined characters $z^{w\cdot\xi}$ and $w'{\delta}_{\rm sm}$. 

We write $K_0(\underline{\lambda}, {\delta}_{\rm sm})$ for the free abelian group generated by the irreducible constituents of the locally analytically induced representations ${\rm Ind}_B^G(z^{w\cdot\xi}\, w'{\delta}_{\rm sm})$, for $w\in W$ and $w'\in W_0$. 
\begin{remark}
The notation $K_0$ should remind the reader of the fact that this free abelian group is the Grothendieck group of the category of finite length locally analytic representations whose irreducible subquotients occur in the above locally analytic principal series representations.
\end{remark}

On the Galois side, recall the definition of the stack $\mathfrak{X}_{d,(\lambdau,\delta_{\rm sm}){\rm -tri}}$ defined just before Theorem \ref{EH:Theorem:localmodelB=P}.
We consider the union
$$\mathfrak{X}_{d,[\lambdau, \delta_{\rm sm}]{\rm -tri}}\bigcup_{w'\in W_0}\mathfrak{X}_{d,(\lambdau,w'\delta_{\rm sm}){\rm -tri}}$$
 and the Grothendieck group $K_0({\rm Coh}(\mathfrak{X}_{d,[\lambdau,\delta_{\rm sm}]{\rm -tri}}))$ of coherent sheaves on this stack. 
 Note that $\mathfrak{X}_{d,[\lambdau, \delta_{\rm sm}]{\rm -tri}}$ is precisely the set of $(\varphi,\Gamma)$-modules that admit a filtration with graded pieces $\mathcal{R}_F(z^{w\lambdau}w'\delta_{\rm sm})$ for some $w\in W$ and $w'\in W_0$.
 The smooth character ${\delta}_{\rm sm}$ defines a Bernstein
 component $\Omega=[T,{\delta}_{\rm sm}]$ of the category of smooth
 representations of $G$ and we write $\tau$ for the inertial type such
 that $\Omega=\Omega_\tau$. As recalled in Section~ \ref{sec:EHBernsteincenters}, the component $\Omega_\tau$ defines connected components 
\begin{align*}
{\rm WD}_{d,F,\tau}&\subset {\rm WD}_{d,F},\\
{\rm Fil}_{\lambdau}{\rm WD}_{d,F,\tau}&\subset {\rm Fil}_{\lambdau}{\rm WD}_{d,F}
\end{align*}
of the stacks of rank~$d$ Weil--Deligne representations, respectively of rank~$d$ filtered Weil--Deligne representations of Hodge--Tate weight $\lambdau$.
The connected component ${\rm WD}_{d,F,\tau}$ admits a canonical map 
\numequation\label{EHmaptoSpecZ}
{\rm WD}_{d,F,\tau}\longrightarrow {\rm Spec}\, \mathcal{Z}_{\Omega_{\tau}},
\end{equation}
see \eqref{EH:eqnBersteiniso}, and ${\delta}_{\rm sm}$ defines a point $[{\delta}_{\rm sm}]\in {\rm Spec}\,\mathcal{Z}_{\Omega_\tau}$.  We write $${\rm WD}_{d,F}({\delta}_{\rm sm})\subseteq {\rm WD}_{d,F,\tau}$$
for the preimage of $[\delta_{\rm sm}]$ under \eqref{EHmaptoSpecZ}. 
\begin{remark}
In terms of (smooth) categorical local Langlands this closed substack can be characterized as the union of the supports of the coherent sheaves $\mathfrak{A}_G^{\rm sm}(\pi)$ attached to the irreducible smooth representations $\pi$ that occur in $\iota_B^G({\delta}_{\rm sm})$, where $\iota_B^G$ denotes normalized parabolic induction. Note that we do not need to take all the possible $\iota_B^G(w'{\delta}_{\rm sm})$ for $w'\in W_0$ into account as $\iota_B^G({\delta})_{\rm sm}$ and $\iota_B^G(w'{\delta})_{\rm sm}$ have the same Jordan--H\"older factors. 
\end{remark}

\begin{lemma} The isomorphism 
$${\rm Fil}_{\lambdau}{\rm WD}_{d,F,\tau}\cong {\rm Fil}_{\lambdau}{\rm Mod}_{d,\varphi,N,\tau}$$ 
of Lemma \ref{EH:lem:filtWDrep} induces a closed embedding 
$${\rm Fil}_{\lambdau}{\rm WD}_{d,F,\tau}({\delta}_{\rm sm})^{\rm an}\subset \mathfrak{X}_{d,(\lambdau,\delta_{\rm sm}){\rm -tri}}.$$
\end{lemma}
\begin{proof}
This is straightforward from the fact that the right hand side is the stack of all $(\varphi,\Gamma)$-modules that are trianguline with parameters $z^{w\lambdau}\, w'\delta_{\rm sm}$. 
\end{proof}

The map in the following conjecture should be thought of as the map induced by $\mathfrak{A}_G^{\rm rig}$ on Grothendieck groups.

\begin{conj}\label{EH:locanBM}
There is a unique injective group homomorphism
\[\mathfrak{a}_{\lambdau,\delta_{\rm sm}}:K_0(\underline{\lambda},\underline{\delta}_{\rm sm})\longrightarrow K_0({\rm Coh}(\mathfrak{X}_{d,[\lambdau,\delta_{\rm sm}]{\rm -tri}}))\]
such that 
\begin{enumerate}
\item[(i)] For $w\in W$ and $w'\in W_0$ one has 
\[\mathfrak{a}_{\lambdau, \delta_{\rm sm}}\big([{\rm Ind}_{\overline{B}}^G(z^{ww_0\cdot\xi}\,w'\delta_{\rm sm})^{\rm an}]\big)=[\widetilde{\beta}_{\hat B,\ast}\circ \widetilde{\alpha}_{\hat B}^\ast(\mathfrak{A}_M^{\rm rig}(s_G(z^{ww_0\cdot\xi}\,w'\delta_{\rm sm})([F:\Q_p]\rho')))].\]
\item[(ii)] For an irreducible smooth representation $\pi$ in the block ${\rm Rep}_{\Omega_{\tau}}^{\rm sm}G$ that is a constituent of $\iota_B^G{\delta}_{\rm sm}$ we have
\[\mathfrak{a}_{\lambdau,\delta_{\rm sm}}([\pi\otimes L(\xi)])=[({\rm pr_{WD}}^\ast(\mathfrak{A}_G^{\rm sm}(\pi))(-\rho')_\sigma].\]
\end{enumerate}
\end{conj}
Of course the formulas (i) and (ii) should be compared to Properties~(1) and~(2) in Conjecture \ref{EH:conjanalytic}.
A weaker version of this conjecture is proved in \cite[Theorem
4.3.8]{MR3623233} under more restrictive hypothesis. In fact the proof
of loc.\ cit.\ globalizes from spectra of deformation rings to the
stack $\mathfrak{X}_{d,[\lambdau,\delta_{\rm sm})]{\rm -tri}}$). More precisely, 
let ${\rm CH}^0(\mathfrak{X}_{d,([\lambdau,\delta_{\rm sm}]{\rm -tri}})$ denote the free abelian group on the irreducible components of $\mathfrak{X}_{d,[\lambdau,\delta_{\rm sm}]{\rm -tri}}$. Mapping a coherent sheaf to its support (and forgetting components of non-maximal dimension) defines a morphism
\[{\rm supp}: K_0({\rm Coh}(\mathfrak{X}_{d,[\lambdau,\delta_{\rm
      sm}]{\rm -tri}}))\longrightarrow {\rm
    CH}^0(\mathfrak{X}_{d,[\lambdau,\delta_{\rm sm}]{\rm -tri}}).\]
Then the methods of proof of \cite[Theorem
4.3.8]{MR3623233} can be used to prove the following theorem.
\begin{theorem}\label{BHS:locanBM}
Let ${\delta}_{\rm sm}=({\rm unr}_{\varphi_1},\dots, {\rm unr}_{\varphi_n})$ be an unramified character with $\varphi_i/\varphi_j\notin \{1,q\}$ for $i\neq j$. Then there is a unique homomorphism 
\[\mathfrak{a}'_{\lambdau,\delta_{\rm sm}}:K_0(\underline{\lambda},{\delta}_{\rm sm})\longrightarrow {\rm CH}^0(\mathfrak{X}_{d,[\lambdau,\delta_{\rm sm}]{\rm -tri}})\]
such that the conditions in (i) and (ii) hold after replacing $\mathfrak{a}_{\lambdau,\delta_{\rm sm}}$ by $\mathfrak{a}'_{\lambdau,\delta_{\rm sm}}$ and after composing the right hand side of the equalities with the support map ${\rm supp}$.
\end{theorem}

\begin{remark}
\hfill\\(a) In fact, in Conjecture \ref{EH:locanBM}, we could ask for more compatibilities, using parabolic induction from more general Levi subgroups than $T$. However, this would make the formulation of the conjecture a bit more involved.\\
(b) While Conjecture~\ref{EH:conjanalytic} is rather hard to attack as there is no candidate for a functor satisfying the constraints, Conjecture \ref{EH:locanBM} is much more explicit: as in the Breuil--M\'ezard conjecture, the group homomorphism $\mathfrak{a}_{\lambdau,\delta_{\rm sm}}$ is determined by the conditions (i) and (ii) in the conjecture and the conjecture is rather about checking relations than about the construction of the map (though in this case, contrary to the Breuil--M\'ezard conjecture, there are only finitely many relations to check). For example, under the additional assumptions in Theorem \ref{BHS:locanBM}, Bezrukavnikov's functor \ref{EH:thmBezfunctor} and Expectation \ref{EH:Bezcompatibility} already determine $\mathfrak{a}_{\lambdau,\delta_{\rm sm}}$ (once one can verify the corresponding relations).\\
(c) The version of Conjecture \ref{EH:locanBM} involving
Grothendieck groups is much finer than its version involving
cycles. For example, in the case $d=2$, if $\pi$ is the (parabolically
induced) smooth representation that is the non-split extension of the
trivial representation $\mathbf{1}$ and the (smooth) Steinberg representation ${\rm St}$, then the cycle underlying
$\mathfrak{a}_{\lambdau,\delta_{\rm sm}}(L(\xi)\otimes\pi)$ is
trivial, whereas the class $\mathfrak{a}_{\lambdau,\delta_{\rm
    sm}}(L(\xi)\otimes\pi)$ in the Grothendieck group of coherent sheaves is non-trivial
(compare e.g.~\cite[Remark 4.43]{hellmann2020derived}). 
More precisely, the functor $\mathfrak{A}_G^{\rm sm}$ maps ${\rm St}$ to the structure sheaf $\cO_{X_1}$ (in degree $0$) of the irreducible component $X_1\subset {\rm WD}_{2,F,1}$ (the stack of $(\varphi,N)$-modules of rank~$2$) where the monodromy $N$ is generically non-trivial, whereas $\mathbf{1}$ is mapped to a line bundle $\mathcal{L}$ on $X_1$ shifted to cohomological degree $-1$. The line bundle $\mathcal{L}$ and the structure sheaf $\mathcal{O}_{X_1}$ are non-isomorphic as coherent sheaves on the stack $X_1$, but they have of course the same underlying cycle. In the formula for the cycle defined by $\mathfrak{A}_G^{\rm sm}(\pi)$ these two cycles cancel, as they live in cohomological degrees with distinct parities. On the other hand the class of $\mathfrak{A}_G^{\rm sm}(\pi)$ in the Grothendieck group of coherent sheaves is still non-trivial.
 \\
(d) The construction of $\mathfrak{a}'_{\lambdau,\delta_{\rm sm}}$
lies at the heart of the construction of companion points on
eigenvarieties in Section~\ref{sec:companionpoints}. In a sloppy way this can
be compared to the Breuil--M\'ezard conjecture as follows: proving the
Breuil--M\'ezard conjecture implies that one can (in some nice
situations) construct automorphic forms with a prescribed Galois
representation. Similarly, the locally analytic Breuil--M\'ezard
multiplicity formula of Theorem~ \ref{BHS:locanBM} implies that one can construct (overconvergent finite slope) $p$-adic automorphic forms with prescribed Galois representation.
\end{remark}

\section{Some known cases of categorical \texorpdfstring{$p$}{p}-adic
  Langlands}\label{sec: known cases}
\subsection{\texorpdfstring{$\GL_1$}{GL1}}
\label{sec: GL1}
Unsurprisingly, the $\GL_1$ case of Conjecture~\ref{conj: Banach
  functor} is a consequence of local class field theory. Furthermore
in this case we can describe the essential image of the functor, and
we do not need to pass to the derived level, because the functor comes
from an equivalence of abelian categories. The key point is that by
local class field theory we have an isomorphism
$W_F^\ab\isoto F^\times = \GL_1(F)$, while by ~\cite[Rem.\
7.2.18]{emertongeepicture}, the stack $\cX_1$ admits a description as
a moduli stack of $1$-dimensional continuous representations of~$W_F$,
and thus of representations of $W_F^\ab\isoto \GL_1(F)$.

More precisely, if we choose a uniformizer $\varpi_F\in F$ then we can
write  $F^\times=\cO_F^\times\times \varpi_F^{\Z}$, and so write  %
$\cO[[F^\times]]=\cO[[\cO_F^\times]]\widehat{\otimes}_{\cO} \cO[X,X^{-1}]^{\wedge}$, where
$\cO[X,X^{-1}]^{\wedge}$ is the $p$-adic completion of $\cO[X,X^{-1}]$, 
and the tensor product is $p$-adically completed.
Then by~\cite[Prop.\
7.2.17]{emertongeepicture} (see also~\cite[Cor.\
1.2]{https://doi.org/10.48550/arxiv.2206.02888} for an alternative
proof by Dat Pham, which is similar to the approach of
Kedlaya--Pottharst--Xiao mentioned below), %
there is an isomorphism
	\numequation\label{eqn: description of 1-d stack} \Bigl[ \Bigl ( \Spf \cO[[F^{\times}]]\Bigr) / \widehat{\mathbf G}_m \Bigr] \iso \cX_{1}\end{equation}%
      (where $\widehat{\mathbf G}_m=\Spf \cO[X,X^{-1}]^{\wedge}$ is the
      $p$-adic completion of~$\Gm$, and in the formation of the
      quotient stack, the $\widehat{\mathbf G}_m$-action is taken to
      be trivial).

It
      follows immediately from ~\eqref{eqn: description of 1-d stack}
      that there is an equivalence of categories between the (abelian)
      category $\cO[[F^\times]]\text{-}\Mod$ %
      and the category of quasicoherent sheaves
      on~$\cX_{1}$ whose $\widehat{\mathbf G}_m$-action is trivial.
Another way to formulate this equivalence is that the structure sheaf
$\cO_{\cX_{1}}$ is naturally an $\cO[[F^\times]]$-module, and that our
functor is given by \numequation\label{eqn: Linfty for d=1} \pi\mapsto
\cO_{\cX_{1}}\otimes_{\cO[[F^\times]]}\pi.\end{equation} In other
words in the case $d=1$, the sheaf $L_\infty$ of Remark~\ref{rem:
  expect Linfty is a kernel} is~$\cO_{\cX_{1}}$.

There is a parallel description of the stack $\mathfrak{X}_1$ and an interpretation of coherent sheaves thereon in terms of locally analytic representations of $F^\times$.
By a result of Kedlaya--Pottharst--Xiao \cite[Thm.\ 6.2.14]{MR3230818} 
every rank one $(\varphi,\Gamma)$-module over $\mathcal{R}_{F,A}$, for an affinoid algebra $A$, is of the form $\mathcal{R}_{F,A}(\delta)\otimes\mathcal{L}$ for a continuous character $\delta:F^\times\rightarrow A^\times$ and a line bundle $\mathcal{L}$ on ${\rm Sp}\,A$. It follows that 
\numequation\label{EHeqn X1analytic}
\mathfrak{X}_1=\mathcal{T}/\mathbf{G}_m=\big(({\rm Spf}\,\mathcal{O}[[\mathcal{O}_F^\times]])^{\rm rig}\times\mathbf{G}_m\big)/\mathbf{G}_m,
\end{equation}
where $\mathcal{T}={\rm Hom}_{\rm cont}(F^\times,\mathbf{G}_m(-))$ is the rigid analytic space of continuous characters of $F^\times$ equipped with the trivial action of the rigid analytic multiplicative group\footnote{Note that we usually write $\mathbf{G}_m$ for the rigid analytic space associated with the scheme $\mathbf{G}_m$. This should not cause confusion as we will always explicitly mention when $\mathbf{G}_m$ should be seen as a scheme.}~$\mathbf{G}_m$. As above the choice of a uniformizer allows us to write $\mathcal{T}=\mathcal{W}\times\mathbf{G}_m$, where  
$$\mathcal{W}={\rm Hom}_{\rm cont}(\mathcal{O}_F^\times,\mathbf{G}_m(-))=({\rm Spf}\,\mathcal{O}[[\mathcal{O}_F^\times]])^{\rm rig}$$
is the space of continuous characters of $\mathcal{O}_F^\times$.
In particular we see that we have a canonical map
$$\mathcal{X}_{1,\eta}^{\rm rig}=(\mathcal{W}\times\widehat{\mathbf{G}}_m^{\rm rig})/\widehat{\mathbf{G}}_m^{\rm rig}\longrightarrow (\mathcal{W}\times\mathbf{G}_m)/\mathbf{G}_m=\mathfrak{X}_1.$$
The image of this map is an open substack of $\mathfrak{X}_1$, and we note that this map is not representable. Its non-trivial fibers are zero-dimensional and given by the stack quotient $\mathbf{G}_m/\widehat{\mathbf{G}}_m^{\rm rig}$. 

There is an
isomorphism $$\mathcal{D}(\cO_F^\times)\cong \Gamma(\mathcal{W},\mathcal{O}_\mathcal{W})$$
of the locally analytic distribution algebra of $\cO_F^\times$ with the
global sections of the structure sheaf on the space $\mathcal{W}$ of characters of $\mathcal{O}_F^\times$.
We note that $\mathcal{D}(F^\times)=\mathcal{D}(\cO_F^\times)[T^{\pm 1}]$ is a Laurent polynomial ring over $\mathcal{D}(\cO_F^\times)$, i.e.~the ring of functions of the algebraic $\mathbf{G}_m$ over $\mathcal{D}(\cO_F^\times)$ whereas 
$$\Gamma(\mathfrak{X}_1,\mathcal{O}_{\mathfrak{X}_1})=\Gamma(\mathcal{T},\cO_{\mathcal{T}})\cong \Gamma(\mathcal{W},\cO_{\mathcal{W}})\widehat{\otimes}_{L}\Gamma(\mathbf{G}_m,\cO_{\mathbf{G}_m})$$
is the space of functions on the rigid analytic $\mathbf{G}_m$ over $\mathcal{D}(\cO_F^\times)$.

 As in \eqref{eqn: Linfty for d=1} we hence want to define a functor
\numequation\label{eqn: analyticLinfty for d=1}
\pi\mapsto\mathcal{O}_{\mathfrak{X}_1}\widehat{\otimes}^L_{\mathcal{D}(F^\times)}\pi,
\end{equation}
where, after restriction to affinoids in $\mathfrak{X}_1$, the symbol $\mathcal{O}_{\mathfrak{X}_1}\widehat{\otimes}^L_{\mathcal{D}(F^\times)}\pi$ should  be a derived completed tensor product (and again: at least the derived aspect only makes sense in the world of solid modules).
Unlike in the case of $p$-power torsion representations it turns out that the functor \eqref{eqn: analyticLinfty for d=1} is not defined on the level of abelian categories, and hence the categorical $p$-adic Langlands correspondence for locally analytic representations is only a derived statement already for ${\rm GL}_1$!
In order to make this more precise it is better to consider the adjoint functor of \eqref{eqn: analyticLinfty for d=1}. To make the discussion easier we consider the case $F=\Q_p$ and even restrict to the case of $\Z_p^\times$ instead of $\Q_p^\times$. The discussion of this case can again easily be reduced to the case of $(\Z_p,+)\cong (1+p\Z_p,\cdot)$.

Recall that the distribution algebra $\mathcal{D}(\Z_p)$ is isomorphic to the ring of analytic functions $\Gamma(\mathbf{U},\cO_{\mathbf{U}})$ on the open unit disc $\mathbf{U}$ over $\Q_p$.
Here we identify $({\rm Spf}\,\Z_p[[\Z_p]])^{\rm rig}$, the space of continuous characters of $\Z_p$, with the open unit disc $\mathbf{U}=({\rm Spf}\,\Z_p[[T]])^{\rm rig}$ using an identification $\Z_p[[\Z_p]]\cong\Z_p[[T]]$.

The rigid analytic space $\mathbf{U}$ is a smooth Stein space. For coherent sheaves on such spaces Chiarelotto \cite{MR1094850} has constructed cohomology with compact support $R\Gamma_c(\mathbf{U},-)$ as well as a Serre duality
\numequation\label{EHeqn SerredualityStein}
H^0(\mathbf{U},\mathcal{F})'\cong H^1_c(\mathbf{U},\mathbf{D}(\mathcal{F})),
\end{equation}
where $(-)'$ denotes the strong dual (as above) and $\mathbf{D}(\mathcal{F})=\mathbf{D}_{\cO_{\mathbf{U}}}(\mathcal{F})$ denotes the Serre dual of the coherent sheaf $\mathcal{F}$.
By functoriality $H^1_c(\mathbf{U},\mathbf{D}(\mathcal{F}))$ (or more generally $R\Gamma_c(\mathbf{U},\mathcal{F}^\bullet)$) is a module over $\mathcal{D}(\Z_p)$ (or a complex of such modules). 
In fact the induced $\Z_p$-action on the $H^1_c$ is locally analytic (and admissible), as by Serre duality, it is the strong dual of a coadmissible module over $\mathcal{D}(\Z_p)$.
The functor \eqref{eqn: Linfty for d=1} now should be defined as the right adjoint functor to the functor 
$$\mathcal{F}\longmapsto R\Gamma_c(\mathbf{U},\mathcal{F})$$
from the derived category of coherent sheaves $\mathbf{D}_{\rm coh}^b(\mathbf{U})$ to a certain derived category of $\mathcal{D}(\Z_p)$-modules in solid $\Q_p$-vector spaces (in the sense of Clausen--Scholze).
We note that the functor $R\Gamma_c(\mathbf{U},-)$ does not preserve the cohomological degree of coherent sheaves: if $\mathcal{F}$ is a coherent sheaf concentrated in degree $0$ that has finite support, then $R\Gamma_c(\mathbf{U},\mathcal{F})$ agrees with the global sections $R\Gamma(\mathbf{U},\mathcal{F})$ and is concentrated in degree $0$ (and is a finite dimensional $\Q_p$ vector space). 
On the other hand $R\Gamma_c(\mathbf{U},\cO_{\mathbf{U}})$ is concentrated in degree $1$ where it is given by the space of locally analytic functions $\mathcal{C}^{\rm an}(\Z_p,\Q_p)$ on $\Z_p$ (as $\mathcal{C}^{\rm an}(\Z_p,\Q_p)$ is identified with the strong dual of $\mathcal{D}(\Z_p)=R\Gamma(\mathbf{U},\mathcal{O}_{\mathbf{U}})$).

Conversely the adjoint functor $-\widehat{\otimes}^L_{\mathcal{D}(\Z_p)}\cO_{\mathbf{U}}$ takes finite dimensional locally analytic representations to coherent sheaves in degree $0$ while it takes $\mathcal{C}^{\rm an}(\Z_p,\Q_p)$ to $\mathcal{O}_{\mathbf{U}}[1]$, the structure sheaf shifted to degree $-1$. 
\begin{remark}
If we restrict to modules or representations concentrated in degree $0$, then the underived version $\pi\widehat{\otimes}_{\mathcal{D}(\Z_p)}\cO_{\mathbf{U}}$ can also be computed (at least for some $\pi$) without passing to the condensed and solid world and gives the sheaf on the Stein space $\mathbf{U}$ associated to a module over the ring of global sections of the structure sheaf $\mathcal{D}(\mathbf{Z}_p)=\Gamma(\mathbf{U},\mathcal{O}_{\mathbf{U}})$. 
The reader should note however that this completed tensor product does not commute with taking global sections, i.e.
$$\Gamma(\mathbf{U},\pi\widehat{\otimes}_{\mathcal{D}(\Z_p)}\cO_{\mathbf{U}})\neq \pi\widehat{\otimes}_{\mathcal{D}(\Z_p)}\Gamma(\mathbf{U},\mathcal{O}_{\mathbf{U}})= \pi .$$
For example in the case of $\pi=\mathcal{C}^{\rm an}(\Z_p,\Q_p)$ the left hand side is $0$, as the sheaf $\pi\widehat{\otimes}_{\mathcal{D}(\Z_p)}\mathcal{O}_{\mathbf{U}}$ is zero: the completed tensor product of $\mathcal{C}^{\rm an}(\Z_p,\Q_p)$ with the ring of functions on every open affinoid of $\mathbf{U}$ vanishes.
However, as pointed out above, the derived version $\pi\widehat{\otimes}^L_{\mathcal{D}(\Z_p)}\cO_{\mathbf{U}}$ is non-zero, and applying $R\Gamma_c(\mathbf{U},-)$ instead of $R\Gamma(\mathbf{U},-)$ gives back the representation $\pi$. 
\end{remark}

\begin{rem}
Using solid condensed modules and analytic rings as introduced by
Clausen--Scholze \cite{Scholzecondensed2019} (and in the context of
rigid analytic spaces studied further by Andreychev
\cite{andreychevpseudocoherent2021}) it is indeed possible to make the
above precise. In order not to get  too involved with condensed structures we will not spell this out here. 
However we note that one should compare the above situation to the compactification of the affine line 
$$\iota:{\rm Spa}(\Z[x],\Z[x])\rightarrow {\rm Spa}(\Z[x],\Z)$$
in the context of discrete adic spaces, see \cite[Theorem 8.1]{Scholzecondensed2019}. In this context the derived completion 
$$\iota^\ast=-\otimes^L_{(\Z[x],\Z)_\blacksquare}(\Z[x],\Z[x])_\blacksquare:\mathbf{D}_\blacksquare(\Z[x],\Z)\rightarrow \mathbf{D}_\blacksquare(\Z[x],\Z[x])$$
has a fully faithful left adjoint $\iota_!$ (the cohomology with compact support) such that $\iota^\ast\circ\iota_!$ is the identity, see \cite[Observation 8.11]{Scholzecondensed2019}.
Moreover $\iota_!(\mathbf{Z}[x])=(\mathbf{Z}((x^{-1}))/\Z[x])[-1]$ lives in cohomological degree $1$ (which implies that the derived completion $\iota^\ast$ maps $\mathbf{Z}((x^{-1}))/\Z[x]$ to the structure sheaf shifted to degree $-1$), see \cite[Observation 8.12]{Scholzecondensed2019}.
\end{rem}

The above discussion can be summarized in the following Theorem, which is a special case of \cite[Theorem B]{solidlocan2}. Here we denote by $D_{\rm qc}(\mathbf{U})$ the derived category of solid quasi-coherent sheaves on $\mathbf{U}$ in the sense of Clausen--Scholze. 
\begin{theorem}
The functor $R\Gamma_c(\mathbf{U},-)$ induces an equivalences of categories
\begin{align*}
D_{\rm qc}(\mathbf{U})\xrightarrow{\cong} &D({\rm an}\,\mathbf{Z}_p)\\
D^b_{\rm coh}(\mathbf{U})\xrightarrow{\cong} &D_{\rm f.p.}({\rm an}\,\mathbf{Z}_p)
\end{align*}
whose inverse functor is given by $-\widehat{\otimes}^L_{\mathcal{D}(\Z_p)}\cO_{\mathbf{U}}$.
\end{theorem}

It may seem surprising that in the Banach case the functor \eqref{eqn:
  Linfty for d=1} is (the derived functor of) a functor between
abelian categories, while in the locally analytic setup the functor
\eqref{eqn: analyticLinfty for d=1} only makes sense on the derived
level. In particular one of these functors cannot shift the
cohomological degree, while the other one can.
Still, these two functors are compatible in the sense of Conjecture~\ref{EH:conjanalytic},
as we now explain.

Consider the pro-system of ind-coherent sheaves $$\mathcal{F}^{(n)}_m=F^{(n)}_m\otimes_{\Z_p[[T]]}\mathcal{O}_{{\rm Spf}\Z_p[[T]]}$$ on ${\rm Spf}(\Z/p^{n}[[T]])={\rm colim}_m\, {\rm Spec}((\Z/p^n)[T]/T^{p^m})$, where $F^{(n)}_m=(\Z/p^n)[\Z_p/p^m\Z_p]$ denotes the representation of $\Z_p$ on the space of $\Z/p^n$-valued functions on $\Z_p$ that are constant on $p^m\Z_p$-cosets (here we of course use the isomorphism $\Z_p[[\Z_p]]\cong \Z_p[[T]]$).
On the level of representations this pro-ind system gives rise (after
taking the limits $\lim_{\leftarrow, n}\lim_{\rightarrow, m}$) to the representation of $\Z_p$ on the space of continuous maps $\mathcal{C}^{\rm cont}(\Z_p,\Z_p)$ and hence, after inverting $p$ and passing to locally analytic vectors we obtain $\mathcal{C}^{\rm an}(\Z_p,\Q_p)$. 
The compatibility of the Banach and the analytic context is hence settled by the following lemma.
\begin{lem}\label{EH:lemtorsionanalyticcomp d=1}
The generic fiber $(\mathcal{F}_m^{(n)})_\eta$ of the pro-ind-system $(\mathcal{F}^{(n)}_m)$ (concentrated in degree $0$) is the structure sheaf on the open unit disc $\mathcal{O}_{\mathbf{U}}[1]$ shifted to degree $-1$.
\end{lem}
Before sketching the proof, we have to define what we mean by the generic fiber $(\mathcal{F}_m^{(n)})_\eta$ of a pro-ind system $(\mathcal{F}_m^{(n)})$.
Recall that the open unit disc has an admissible cover
$$\mathbf{U}=\bigcup_{k\geq 1} U_k$$ by the closed discs $U_k$ of radius $p^{-1/k}$ and each of the closed discs has a nice affine formal model $\mathcal{U}_k$. 
Let $j_k:U_k\hookrightarrow \mathbf{U}$ denote the open embedding and let us write $\pi_k:\mathcal{U}_k\rightarrow {\rm Spf}\,\Z_p[[T]]$ for the canonical map. Then $j_k^\ast((\mathcal{F}_m^{(n)})_\eta)$ is by definition the generic fiber of the coherent sheaf on $\mathcal{U}_k$ given by 
$$\lim_{\leftarrow, n}\lim_{\rightarrow, m} \pi_k^\ast \mathcal{F}_m^{(n)},$$
and we define $(\mathcal{F}_m^{(n)})_\eta$ by gluing these sheaves on the cover $\mathbf{U}=\bigcup_{k\geq 1} U_k$.
\begin{proof}[Sketch of proof of Lemma~\ref{EH:lemtorsionanalyticcomp d=1}.]
We sketch the case $k=1$ in which case $\mathcal{U}_1={\rm Spf}\,\Z_p\langle X \rangle$ with $X=p^{-1}T$.
For simplicity we only prove that 
$$\lim_{\rightarrow, m}\bar\pi_1^\ast \mathcal{F}^{(1)}_m$$ is the structure sheaf on $\bar{\mathcal{U}_1}={\rm Spec}\,\F_p[X]$ shifted to degree $-1$ and leave the general case as an exercise . Here $$\bar\pi_1: {\rm Spec}\,\F_p[X]\rightarrow {\rm Spf}\,\F_p[[T]]$$ is the mod $p$ fiber of $\pi_1$, i.e. the morphism given by $T\mapsto 0$. 
One can compute that $\mathcal{F}^{(1)}_m$ is the coherent sheaf corresponding to the module $\F_p[T]/T^{p^m}$ with transition maps induced by $T\mapsto T^p$ (a quick way to see this is to note that $\lim_{\rightarrow,m}F_m^{(1)}$ is the Pontryagin-dual of $\F_p[[\Z_p]]$ and $\lim_{\rightarrow,m}\F_p[T]/T^{p^m}$ is the Pontryagin dual of $\F_p[[T]]$, but we have fixed an identification $\F_p[[\Z_p]]\cong\F_p[[T]]$).

Computing the derived pullback $\bar\pi_1^\ast$ using the resolution $$\F_p[[T]][X]\xrightarrow{\cdot T}\F_p[[T]][X]$$ of $\F_p[X]$, we find that $\bar\pi_1^\ast \mathcal{F}^{(1)}_m$ is a two-term complex with cohomology 
$\F_p[X]=T^{p^m-1}(\F_p[T]/T^{p^m})[X]$ in degree $-1$ and $\F_p[X]=(\F_p[T]/T)[X]$ in degree $0$. 
However the transition maps $T\mapsto T^p$ induce isomorphisms in degree $-1$ whereas they induce the zero maps in degree $0$. It follows that $\lim_{\rightarrow, m}\bar\pi_1^\ast\mathcal{F}_m^{(1)}$ is $\F_p[X]$ concentrated in degree $-1$ as claimed.
\end{proof}

Let us return to the functor (\ref{eqn: analyticLinfty for d=1}), i.e.~to the case of locally analytic representations of $F^\times$ rather than $\mathbf{Z}_p$. In this case we obtain a functor 
\numequation\label{EHeqn analyticfunctorGL1}
\mathfrak{A}^{\rm rig}_{F^\times}:D^b_{\rm f.p.}({\rm an}\,F^\times)\longrightarrow D^b_{\rm coh}(\mathfrak{X}_1).
\end{equation}
From the point of view of solid locally analytic representations this functor is discussed in more detail in \cite[4.4]{solidlocan2}.

\begin{rem}\label{EH rem: ffGL1}
We note however that the functor (\ref{EHeqn analyticfunctorGL1}) cannot be fully faithful on the category $D^b_{\rm f.p.}({\rm an}\,F^\times)$: Roughly the functor involves analytification from the algebraic $\mathbf{G}_m$ to the analytic $\mathbf{G}_m$. To give an explicit example, let $\pi$ be the universal unramified representation $$F^\times \rightarrow F^\times/\cO_F^\times\cong \mathbf{Z}\rightarrow \Gamma(\mathbf{G}_m^{\rm alg},\cO^\times_{\mathbf{G}_m^{\rm alg}})=L[T,T^{-1}]^\times$$
mapping $1$ to $T$, where we write $\mathbf{G}_m^{\rm alg}$ for the scheme $\mathbf{G}_m$.
Then (\ref{EHeqn analyticfunctorGL1}) maps $\pi$ to the structure sheaf on the quotient of the rigid analytic space $\mathbf{G}_m=\{1\}\times\mathbf{G}_m\subset\mathcal{W}\times\mathbf{G}_m=\mathcal{T}$, where $1\in\mathcal{W}$ denotes the trivial character of $\cO_F^\times$, by the trivial $\mathbf{G}_m$-action .
It follows that ${\rm End}_{F^\times}(\pi,\pi)=L[T^{\pm 1}]$, whereas the endomorphism ring of $\cO_{\mathfrak{X}_1}\widehat{\otimes}^L_{\mathcal{D}(F^\times)}\pi$ is the ring of functions on the rigid analytic $\mathbf{G}_m$. 
In fact the functor (\ref{EHeqn analyticfunctorGL1}) can be seen as some kind of localization of a representation $\pi$ on $\mathfrak{X}_1$. 
The functor however is fully faithful on \emph{tempered} representations $\pi$. Roughly these are the representations $\pi$ such that the canonical map
$$\pi\rightarrow R\Gamma_c(\mathfrak{X}_1,\cO_{\mathfrak{X}_1}\widehat{\otimes}^L_{\mathcal{D}(F^\times)}\pi)$$
 is an isomorphism, or equivalently  those representations $\pi$ on which the $\mathcal{D}(F^\times)$-action extends to $\Gamma(\mathcal{T},\cO_{\mathcal{T}})$. This is the case if $\pi$ has bounded slope, for example if $\pi$ is an admissible representation. We refer to \cite[Theorem 4.4.4]{solidlocan2} for more details, in particular for a representation-theoretic account of the category of tempered $F^\times$-representations. 
 It is worth pointing out that \cite[Theorem 4.4.4]{solidlocan2} provides a second way to fix the issue that (\ref{EHeqn analyticfunctorGL1}) is not fully faithful: if the target category is replaced by the category of coherent sheaves on the analytic stack $(\mathcal{W}\times\mathbf{G}_m^{\rm alg})/\mathbf{G}_m^{\rm alg}$, then one obtains a fully faithful functor  
 $$D^b_{\rm f.p.}({\rm an}\,F^\times)\longrightarrow D^b_{\rm coh}((\mathcal{W}\times\mathbf{G}_m^{\rm alg})/\mathbf{G}_m^{\rm alg})$$
 whose essential image is given by those (complexes of) coherent sheaves on which the $\mathbf{G}_m^{\rm alg}$-action is trivial.
 This suggests that the moduli stack $\mathfrak{X}_d$ in fact should be defined on a larger test category of analytic spaces that contain rigid analytic spaces as well as algebraic varieties like $\mathbf{G}_m^{\rm alg}$ (and in which $\mathbf{G}_m^{\rm alg}$ does not coincide with the rigid analytic space $\mathbf{G}_m$), as in work of Mikami \cite{mikami2024varphigammamodulesrelativelydiscretealgebras}.
 \end{rem}

\subsection{The Banach case for \texorpdfstring{$\GL_2(\Qp)$}{GL2(Qp)} --- I. The structure of~\texorpdfstring{$\cX$}{the stack}}
\label{sec:banach-case-gl_2qp}%
In this section and the following two,
we summarise the results of the work in
progress~\cite{DEGcategoricalLanglands} of Andrea Dotto with M.E.\ and T.G., which
establishes a version of Conjecture~\ref{conj: Banach functor} for
$\GL_2(\Qp)$. We caution that since this is work in progress, there
may be imprecisions in what follows, and the sketched proofs do not always
reflect those in the final version of~\cite{DEGcategoricalLanglands}.

For simplicity of exposition, we work throughout these two sections with representations having
trivial central character. (The papers ~\cite{DEGcategoricalLanglands,
  DEGlocalization} allow an arbitrary fixed central character, and the
arguments are essentially the same, with some notational overhead.) %
Accordingly we set~$G:=\PGL_2(\Qp)$, and $K=\PGL_2(\Zp)$, so that
$\smG$ is the abelian category of smooth $\GL_2(\Qp)$-representations with
trivial central character. %

On the geometric side,
we work %
with the moduli stack~$\cX$
of rank~$2$ $(\varphi,\Gamma)$-modules for~$\GL_2/\Qp$ with
determinant $\varepsilon^{-1}$. %
As explained following Remark~\ref{rem: don't make dimension of
  whole EG stack precise}, this is a formal algebraic stack
over~$\Spf\cO$, whose underlying reduced substack~$\cX_{\red}$ is equidimensional
of dimension~$1$. We make the geometry of~$\cX$ much more explicit in
the discussion below.

The following expected theorem verifies Conjecture~\ref{conj: Banach
  functor} in the case of $\GL_2(\Qp)$.  %
\begin{expectedthm}[Dotto--Emerton--Gee]
  \label{expectedthm: DEG results}Suppose that $p\ge 5$, and let
  $G=\PGL_2(\Qp)$. %
  Then there is a
  pro-coherent sheaf \emph{(}concentrated in degree~$0$\emph{)} of  $\cO[[G]]$-modules $L_\infty$ over~$\cX$,
  such that \[\fA(\pi):= L_{\infty}\otimes^{L}_{\cO[[G]]}\pi\] defines
 an $\cO$-linear exact
  fully faithful functor
$\fA:\Dfp^b(\smG) \to D_{\coh}^b(\cX)$, which then
extends to a continuous fully faithful functor
$\fA:\Ind \Dfp^b(\smG) \to \Ind\Coh(\cX)$. %
\end{expectedthm}
\begin{rem}We also expect to prove the various properties conjectured
  in Conjecture~\ref{conj: Banach
  functor}. In particular, the statement about the support of
$\fA\Bigl(\bigl(\cInd_{K}^{G}W_{\lambdau}\otimes_{\cO}\sigma^{\crys,\circ}(\tau)\bigr)^{\widehat{p}}
\Bigr)$
should follow easily from the construction and the classification of
locally algebraic vectors in the existing $p$-adic local Langlands
correspondence for~$\GL_2(\Qp)$. 

  The expected compatibility with duality is discussed in Section~\ref{subsubsec:duality for }.   %
\end{rem}

The following theorem is an ingredient in the proof of Expected
Theorem~\ref{expectedthm: DEG results}. It is essentially the main result
of~\cite{DEGlocalization} (together with the construction of the
morphism $\cX_{\red}\to X$, which will appear
in~\cite{DEGcategoricalLanglands} and is explained below). It is special
to the case of~$\GL_2(\Qp)$, in the strong sense that we do not expect it to
generalise in any obvious way, even to $\GL_2(\Q_{p^2})$.

\begin{thm}[Dotto--Emerton--Gee]%
  \label{expectedthm: Bernstein centre}Assume that $p\ge 5$ and
  $G=\PGL_2(\Qp)$. There is a reduced scheme~$X$ given 
  by an explicit chain of $\Pone$s over~$k$ with ordinary double
  points, and a morphism $\cX_{\red}\to X$ inducing a bijection on
  closed $\Fpbar$-points. For each open subset~$U$ of~$X$, there is a certain
  localization $(\smG)_U$ of $\smG$
  such that $(\smG)_X=\smG$, and the collection $\{(\smG)_U \}$ forms a
  stack \emph{(}of abelian categories\emph{)} over the Zariski site of~$X$.
\end{thm}

\begin{rem}
  \label{rem: temporary remark about Bernstein centre}We anticipate
  that Theorem~\ref{expectedthm: Bernstein centre} can be
  refined to give much more precise information about the
  compatibility of localization with the functor~$\mathfrak{A}$, and about
  the Bernstein centres of the categories $(\smG)_U$. Indeed, we
  expect that there is a  formal scheme~$\widehat{X}$ with underlying
  reduced scheme~$X$, which is equipped with a morphism
  $\cX\to\widehat{X}$, and which simultaneously realises $\widehat{X}$ as
  a moduli space for~$\cX$ (in the sense that $\cX \to \widehat{X}$ should
  be initial for maps from $\cX$ to formal algebraic spaces; indeed 
  we expect that the map $\cX \to \widehat{X}$ will satisfy the (formal
  algebraic analogues) of the conditions
  of \cite[Prop.~7.1.1]{MR3272912}), and whose ring of functions over each~$U$
  is identified with the Bernstein centre of $(\smG)_U$.  We
  hope to prove such results in a sequel to~\cite{DEGcategoricalLanglands}. %
\end{rem}

\begin{rem}
  \label{rem: p at least 5} We expect that the statements of Expected
  Theorem~\ref{expectedthm: DEG results} and 
Theorem~\ref{expectedthm:
    Bernstein centre} continue to hold without the assumption that
  $p\ge 5$, and that with some effort the proofs sketched below can be
  extended to this case. The difficulty when $p=2$ or $3$ is that there are
  many more special cases to consider (or rather, the existing special
  cases overlap in more complicated ways); on the Galois side, this is
  because the mod~$p$ cyclotomic character is quadratic (if $p=3$) and
  trivial (if $p=2$), and on the representation theory side, there are
  no very generic Serre weights in the sense of Definition~\ref{defn:
    generic Serre weight} below (and indeed no generic Serre weights
  if $p=2$). %
\end{rem}

\begin{rem}
  \label{rem: compatibility of DEG with patching}Under mild
  hypotheses, the expected compatibility with Taylor--Wiles--Kisin
  patching explained in Remark~\ref{rem:
    conjectural compatibility with TW patching} is a consequence of
  the construction of our functor below, and local-global
  compatibility for $p$-adic local Langlands for~$\GL_2(\Qp)$ of
  \cite{emerton2010local} (see also ~\cite{MR3732208} for a more
  direct connection between the patching method for modular curves and
  $p$-adic local Langlands for~$\GL_2(\Qp)$, under slightly different
  hypotheses). See also Expected Theorem~\ref{expected EGZ thm}.
 \end{rem}
\subsubsection{Motivation}\label{subsubsec: motivation for DEG
  construction}One reason that we are able to prove
Conjecture~\ref{conj: Banach functor} for $\GL_2(\Qp)$ (but not in any
other cases) is that we are able to build upon the results of
Colmez~\cite{MR2642409} and Pa\v{s}k\={u}nas~\cite{MR3150248}. Roughly
speaking, these results can be reinterpreted as proving
Conjecture~\ref{conj: Banach functor} ``pointwise''; that is, after
completing $\cX$ at a closed $\Fpbar$-point on the spectral side, and
restricting to a block of locally admissible representations on the
automorphic side.
Since writing these notes, this reinterpretation has been made precise in the work of Johansson--Newton--Wang-Erickson~\cite{johansson2024modulistacksgaloisrepresentations}.

There does not seem to be any hope of directly passing from these
results to Expected Theorem~\ref{expectedthm: DEG results}, and morally
our strategy is to carry out constructions over~$\cX$ which are
analogues of some of the constructions  of Colmez and
Pa\v{s}k\={u}nas, and where possible to check properties of these
constructions by passing to completions and using their results. In
particular, $\cX_{\red}$ generically consists of reducible
representations, which (for $\GL_2(\Qp)$, but not more generally)
correspond via the mod $p$ local Langlands correspondence for
$\GL_2(\Qp)$ to (extensions of) principal series representations, and we are able to
interpolate various explicit calculations in~\cite{MR3150248} over the
reducible locus.

In order to carry out this strategy we need to have an analogue on the
automorphic side of familiar operations on sheaves such as restriction
to open subsets and passage to completion at closed
points. In~\cite{DEGlocalization} we develop such a localization
theory for $\smG$, proving (for the most part by purely
representation-theoretic arguments, together with some use of the
$p$-adic local Langlands correspondence for $\GL_2(\Qp)$) that it
localizes over the projective scheme~$X$ of Theorem~\ref{expectedthm: Bernstein centre}. %

We then generalize Colmez's construction in \cite{MR2642409} of the
$\GL_2(\Qp)$-representation he calls $D^\natural\boxtimes\Pone$, to
give a direct construction of our candidate for~$L_\infty$. We
therefore have an explicit functor~$\mathfrak{A}$, whose behaviour
after completion at closed points can be understood using the results
of \cite{MR2642409,MR3150248}. %

Our arguments rely on having a good understanding of the geometry of the stack
$\cX$, as well as an interpretation of this geometry in Galois representation-theoretic
terms,
and in the remainder of this section we describe this geometry.  In the following
section we explain the construction of $L_\infty$ and sketch the main
points in the proof of Expected Theorem~\ref{expectedthm: DEG results}.

\subsubsection{The geometry of~$\cX_{\red}$}
\label{subsubsec:expl-open-subst} %
We begin by describing the
irreducible components of~$\cX_{\red}$ and their intersections. These
are described by the weight part of Serre's conjecture, and we begin
by recalling this. %

For any~$x\in\Fpbartimes$, we let~$\lambda_x:\Gal_{\Qp}\to\Fpbartimes$ be the
unramified character taking a geometric Frobenius element to~$x$. We
write~$\omega=\varepsilonbar$ for the mod~$p$ cyclotomic character, so
that every character $\Gal_{\Qp}\to\Fpbartimes$ is of the
form~$\lambda_x\omega^i$ for some uniquely determined~$x\in\Fpbartimes$
and $0\le i<p-1$.

We write~$\omega_2:I_{\Qp}\to\Fpbartimes$ for a choice
of fundamental character of niveau~$2$; writing~$\Q_{p^2}$ for the
quadratic unramified extension of~$\Qp$, we can extend~$\omega_2$ to a
character $\Gal_{\Q_{p^2}}\to\Fpbartimes$ in such a way
that~$\omega_2^{p+1}=\omega|_{\Gal_{\Q_{p^2}}}$. Then the 2-dimensional absolutely irreducible
representations $\Gal_{\Qp}\to\GL_2(k)$ are precisely those of the
form~$\lambda_x\otimes\Ind_{\Gal_{\Q_{p^2}}}^{\Gal_{\Qp}}\omega_2^n$ for
some~$x\in\Fpbartimes$ and $n\in\Z/(p^2-1)\Z$ with $(p+1)\nmid n$. In this case~$x$ is
uniquely determined up to multiplication by~$-1$, and~$n$ is uniquely
determined up to multiplication by~$p$. The determinant of this
representation is an unramified twist of~$\omega^n$.

  Serre weights  for  $\GL_2(\Fp)$  are the (isomorphism classes of)
irreducible $\Fp$-representations of~$\GL_2(\Fp)$, which are
explicitly given by
  $\sigma_{a,b}=\det^{a}\Sym^b\Fp^2$, where $0\le
  a<p-1$ and~$0\le b\le p-1$. It is sometimes convenient to view~$a$
  as an element of~$\Z/(p-1)\Z$, and we will do so without further
  comment. 

  Following Serre~\cite{SerreDuke} (and as reinterpreted by
  Buzzard--Diamond--Jarvis \cite{BuzzardDiamondJarvis}) one can
  associate to each
  representation~$\rhobar:\Gal_{\Qp}\to\GL_2(\Fpbar)$ a set of Serre
  weights~$\sigma$. %
  One way of defining this is to say that
  $\sigma_{a,b}$ is a Serre weight for~$\rhobar$ if and only
  if~$\rhobar$ admits a crystalline lift with Hodge--Tate
  weights~$1-a,-(a+b)$. This description can be made completely
  explicit as follows.

If~$\rhobar=\lambda_x\otimes\Ind_{\Gal_{\Q_{p^2}}}^{\Gal_{\Qp}}\omega_2^n$ is
irreducible, then~$\sigma_{a,b}$ is a Serre weight for~$\rhobar$ if
and only if $n\equiv (p+1)a+b-p\pmod{p^2-1}$ or $n\equiv
(p+1)a+pb-1\pmod{p^2-1}$. 

If~$b\ne 0$, then~$\sigma_{a,b}$ is a Serre weight for a
reducible representation~$\rhobar$ %
if and only if %
\numequation\label{eqn: form of reducible repn Serre weight
  fixed}\rhobar\cong
\begin{pmatrix}
  \lambda_x \omega^{a+b} &* \\0&\lambda_y\omega^{a-1}
\end{pmatrix}\end{equation}for some~$x,y\in \Fpbartimes$. If~$b=0$, then we further
demand that~$\rhobar$ is finite flat. If~$x\ne y$ then this is
automatic, while if $x=y$ it is equivalent to asking that the
extension be peu ramifi\'ee (which by definition is equivalent to the
corresponding Kummer class being one associated to an integral unit).

Note that if~$\det \rhobar=\omega^{-1}$ then any Serre weight
for~$\rhobar$ has trivial central character. %
We assume that all Serre weights 
have  trivial central character from now on; that is, we only consider
Serre weights~$\sigma_{a,b}$ with $2a+b\equiv 0\pmod{p-1}$. In
particular, $b$ is always even.

Then the irreducible components of~$\cX_{\red}$ are as follows: for
each Serre weight~$\sigma=\sigma_{a,b}$ with $b\ne p-1$, there is a
unique irreducible component $\cX(\sigma)$, whose $\Fpbar$-points are
precisely those $\rhobar:\Gal_{\Qp}\to\GL_2(\Fpbar)$ admitting~$\sigma$
as a Serre weight. %
Furthermore if~$\sigma=\sigma_{a,p-1}$ for $a=0,(p-1)/2$, then there are two corresponding irreducible
components~$\cX(\sigma)^{\pm}$, which differ by a twist by the
unramified quadratic character~$\lambda_{-1}$. In this case, the
representations which admit~$\sigma_{a,p-1}$ as a Serre weight are
precisely the union of the $\Fpbar$-points of the
stacks~$\cX(\sigma_{a,p-1})^+$, $\cX(\sigma_{a,p-1})^-$, and
$\cX(\sigma_{a,0})$. %

We now describe the stacks $\cX(\sigma)$ more explicitly. %

\begin{defn}\label{defn: generic Serre weight}
A Serre weight~$\sigma_{a,b}$ is \emph{generic} if
$0 \leq b \leq p-3$ (equivalently, if $b\ne p-1$).
It is {\em very generic} if furthermore $b \neq 0.$
\end{defn}

The following result is essentially a consequence of
Fontaine--Laffaille theory (see Section~\ref{subsubsec:
  components X sigma} for slightly more explanation).
\begin{prop}
  \label{prop: explicit description of FL stack}If~$\sigma=\sigma_{a,b}$ is
  generic, then we have an isomorphism $\cX(\sigma) \iso [(\A^2 \setminus \{0\}) / \Gm],$ where
if we write~$(t,x)$ for the
coordinates on~$\A^2$, then the $\Gm$-action is given
by \[u\cdot(t,x)=(t,u^2x). \]
The reducible locus is the locus~$t\ne 0$, and the split locus
is the locus~$x=0$. More precisely, 
the locus~$t\ne 0$ parameterizes the universal extension \[
  \begin{pmatrix}\lambda_{t}\omega^{a+b}&*\\
0&    \lambda_{t^{-1}}\omega^{a-1}
  \end{pmatrix},
\]
with the variable~$x$ parameterizing the %
extension class.
\end{prop}

\[\includegraphics[width=0.8\textwidth]{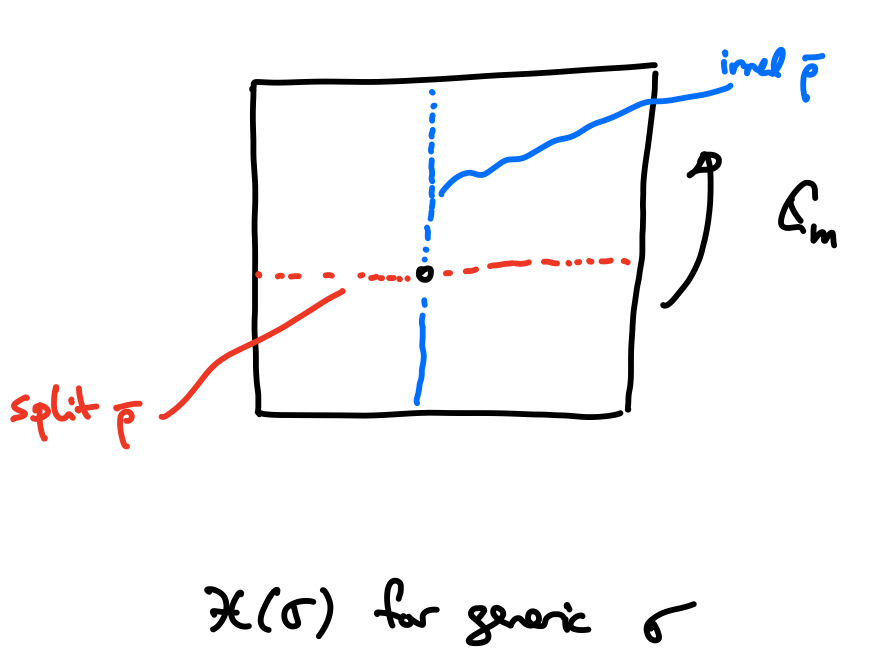}\]

\begin{remark}
We can describe and interpret 
the orbits of $\Gm$ on the $\Fbar_p$-valued points of $\A^2\setminus \{0\}$ as follows: %
\begin{itemize}
\item
If $t = 0$ then $x \neq 0$, and so the locus $t = 0$ is a single orbit,
with associated residual gerbe equal to $[\Spec \Fbar_p /\bm{\mu}_2];$
the stabilizer $\bm{\mu}_2$ appears here because the points
$u \in \Gm$ act on $x$ via multiplication by~$u^2$.
This point corresponds to the unique irreducible  
representation $\rhobar: \Gal_{\Q_p} \to \GL_2(\Fbar_p)$ having determinant~$\omega^{-1}$
and admitting $\sigma$ as a Serre weight.
(In this latter optic, the stabilizer 
$\boldsymbol{\mu}_2$ appears because it is the centralizer of $\rhobar$ in~$\SL_2$.)
\item If we fix a value of $t \neq 0$,
then there are two orbits, namely $x = 0$ and $x \neq 0.$ 
The corresponding residual gerbes are $[\Spec \Fbar_p /\Gm]$ and $[\Spec \Fbar_p / \boldsymbol{\mu}_2].$
The ``$x$-line'' through $t$ can be interpreted as the $1$-dimensional space
$\Ext^1(\lambda_{t^{-1}}\omega^{a-1} ,
\lambda_t \omega^{a+b}),$
with the $\Gm$-action corresponding to its action on this $\Ext^1$ induced
by its scaling action on these characters.
(If we scale $\lambda_t \omega^{a+b}$ by~$u$, then we scale 
$\lambda_{t^{-1}}\omega^{a-1}$ by~$u^{-1}$, so that the isomorphism
of the determinant of the extension 
--- which is the (tensor) product 
$(\lambda_t \omega^{a+b}) ( \lambda_{t^{-1}}\omega^{a-1})$ ---
with $\omega^{-1}$ stays fixed. This is why $u$ acts on elements
$x$ in the $\Ext^1$ via $u^2$.) 
The fact that there are two $\Gm$-orbits corresponds to the  
fact that there are two isomorphism classes of two-dimensional $G_{\Q_p}$-representations
arising from the elements of this $\Ext^1$: the split extension (corresponding
to~$x = 0$), and the non-split extensions (all of which give rise to isomorphic
$\Gal_{\Q_p}$-representations).
(Note that the centralizers in $\SL_2$ of these two representations 
are the diagonal torus --- a copy of $\Gm$ --- and the centre $\boldsymbol{\mu}_2$,
matching with the stabilizer groups appearing in the preceding description of
the associated residual gerbes.)
\end{itemize}
\end{remark}

 In the non-generic case, we have the following description of
  the stacks $\cX(\sigma_{a,p-1})^{\pm }$, which
can be proved by constructing universal families of extensions of
rank~1 $(\varphi,\Gamma)$-modules.
 \begin{prop}
  \label{prop: explicit description weight $p-1$ stacks}If~$\sigma=\sigma_{a,p-1}$, then  we have isomorphisms
$[\A^2/\Gm] \iso \cX(\sigma)^{\pm },$ 
where
if we write~$(x,y)$ for the
coordinates on~$\A^2$, then the $\Gm$-action is given 
by \[u\cdot(x,y)=(u^2 x,u^2 y). \]
 More precisely, $\cX(\sigma)^{\pm }$ 
is the universal extension \[
  \begin{pmatrix}\lambda_{\pm 1}\omega^{-a}&*\\
0&    \lambda_{\pm 1}\omega^{-(a+1)}
  \end{pmatrix}.
\]
The variables~$x,y$ parameterize the extension class in the
$2$-dimensional $k$-vector space $H^1(\Gal_{\Qp},\omega)$.
\end{prop}

\begin{defn}\label{defn: companion Serre weights}
  The \emph{companion} of a generic Serre weight
  ~$\sigma=\sigma_{a,b}$ is the Serre weight
  $\sigmacomp:=\sigma_{a+b+1,p-3-b}$ (which is again generic). We
  refer to the unordered pair of~$\sigma,\sigmacomp$ as a
  \emph{companion pair}.

  (We
do not define a companion for a Serre weight of the
form~$\sigma_{a,p-1}$.)%

We say that a companion pair~$\sigma,\sigmacomp$ 
is %
very generic if $\sigma$ and~$\sigmacomp$ are %
not of the form $\sigma_{a,0},\sigma_{a+1,p-3}$.
\end{defn}

\begin{rem}%
  \label{rem: companion weights and cuspidal types}The companion pairs
  $\sigma,\sigmacomp$ are precisely the sets of Jordan--H\"older
  constituents of the reductions modulo~$p$ of tame cuspidal
  types (with trivial central character).
\end{rem}

The following decomposition of~$\cX_{\red}$ will be useful
below. %
\begin{defn}
  \label{defn: X sigma sigmaco}If~$\sigma=\sigma_{a,b}$ is a generic Serre
  weight, %
  then we define a closed substack
  $\cX(\sigmasigmacomp)$ of~$\cX_{\red}$ as follows: if $b\ne 0,p-3$
  then we set \[\cX(\sigmasigmacomp)=\cX(\sigma)\cup\cX(\sigmacomp),\] if
  $b=0$ then we set \[\cX(\sigmasigmacomp)=\cX(\sigma_{a,0})\cup\cX(\sigma_{a,p-1})^+\cup\cX(\sigma_{a,p-1})^-\cup\cX(\sigma_{a+1,p-3}),\]
  and if $b=p-3$ then we set \[\cX(\sigmasigmacomp)=\cX(\sigma_{a-1,0})\cup\cX(\sigma_{a-1,p-1})^+\cup\cX(\sigma_{a-1,p-1})^-\cup\cX(\sigma_{a,p-3}).\] Note
  in particular
  that~$\cX(\sigmasigmacomp)$ depends only on the unordered pair
  $\{\sigma,\sigmacomp\}$, and that $$\cX_{\red} = \bigcup_{\sigmasigmacomp} \cX(\sigmasigmacomp),$$
where$\{\sigma,\sigmacomp\}$ runs over all companion pairs of Serre
weights.\end{defn}

\[\includegraphics[width=0.8\textwidth]{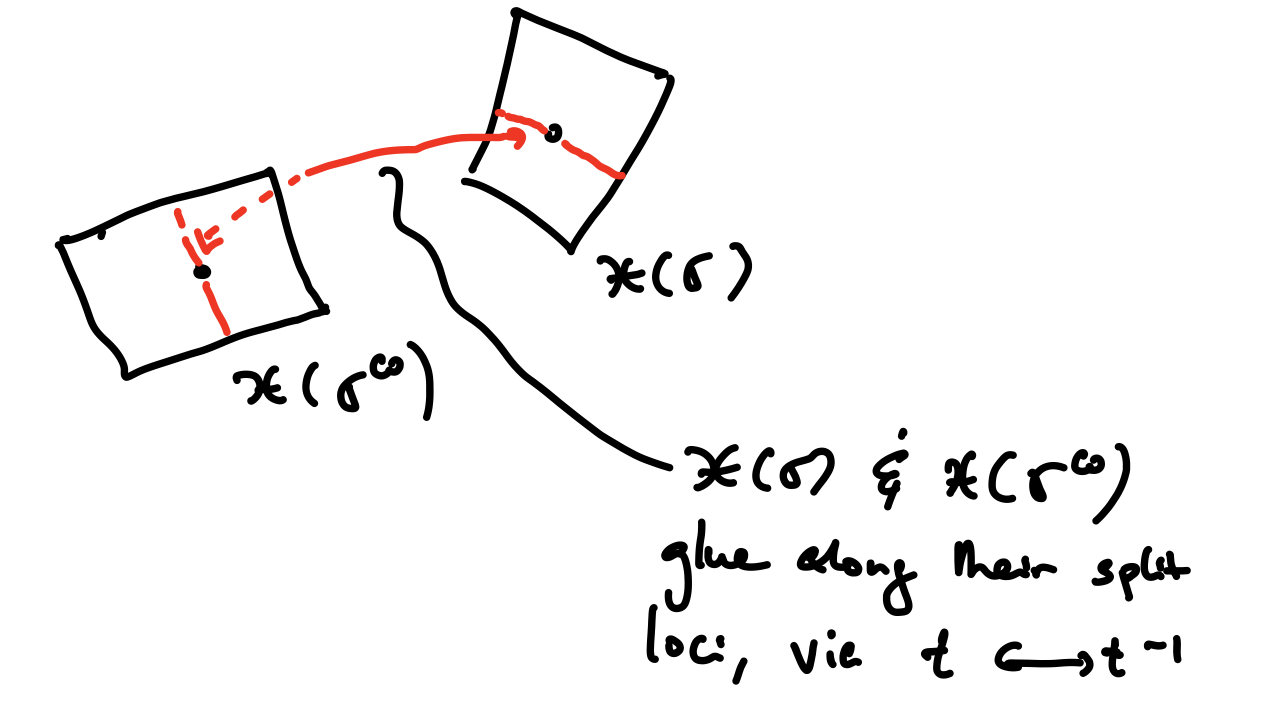}\]

We let $X$ denote a chain of $\mathbf{P}^1$'s over $k$ with ordinary
double points, of length $(p -1)/2$. %
We choose a coordinate~$T$ on each irreducible
component of~$X$ in such a way that each singular point corresponds
to~$0$ on one intersecting component and~$\infty$ on the other. We
will refer to the points~$0$ and~$\infty$ as \emph{marked
  points}.

We label the components of~$X$ by pairs $\sigma,\sigmacomp$ of
companion Serre weights, and denote 
the corresponding~$\Pone$ by~$X(\sigmasigmacomp)$. We order the
components by demanding that if~$0<b< p-1$, then the
component~$X(\sigma_{a,b}|\sigma_{a+b+1,p-3-b})$ meets the
component~$X(\sigma_{a+b,p-1-b}|\sigma_{a+1,b-2})$. It is easy to check that there
are precisely two such labellings of the components (and the choice of a labelling amounts to the choice of an end of the chain). %

\[\includegraphics[width=0.5\textwidth]{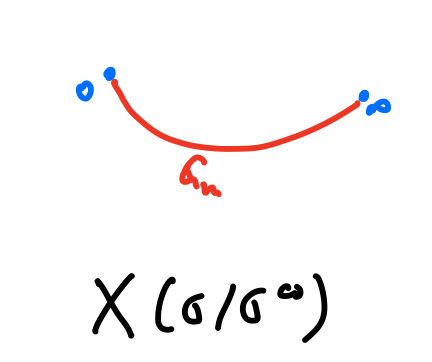}\]

The point of this construction is that using our explicit descriptions
of the stacks~$\cX(\sigma)$ we can define a morphism
$\pi:\cX_{\red}\to X$. Again we content ourselves with this
description in the case of a very generic pair $\sigma,\sigmacomp$. We define a
morphism $\pi_\sigma:\cX(\sigmasigmacomp)\to X(\sigmasigmacomp)$ as
follows. We identify $\cX(\sigma)$ with~$[(\A^2\setminus\{0\})/\Gm]$
as in Proposition~\ref{prop: explicit description of FL stack}, and
with the coordinates~$(t,x)$ there, we send $(t,x)\mapsto t$. We also
identify $\cX(\sigmacomp)$ with~$[(\A^2\setminus\{0\})/\Gm]$, and send
$(t,x)\mapsto t^{-1}$.  Note that by definition we
have $\pi_{\sigmacomp}(x,t)=\pi_\sigma(x,t)^{-1}$. %

For each companion pair~$\sigma,\sigmacomp$ of Serre weights, we make
a choice of one of
$\pi_{\sigma}:\cX(\sigmasigmacomp)\to X(\sigmasigmacomp)$ or
$\pi_{\sigmacomp}:\cX(\sigmasigmacomp)\to X(\sigmasigmacomp)$, in such
a way that the two maps agree at the intersection points of
adjacent components of~$X$. (There are two possible such choices of maps,
corresponding to the two choices of labelling of the components
of~$X$; one can convert between the two possible labellings by
twisting by $\det^{(p-1)/2}$.)%

\begin{rem}
\label{rem:description of cX(sigmasigmacomp) in the generic case}
If$\{\sigma,\sigmacomp\}$ is very generic,
then
\begin{multline*}
\cX(\sigmasigmacomp) \cong
\\
\bigl(
\Pone_t \times \bigl[ \bigl(\Spec k[x,y]/(xy)\bigr)/(\Gm)_u\bigr] 
\bigr)
\setminus ( 0 \times [\{x = 0\}/\Gm] \cup \infty \times [\{y = 0\}/\Gm]),
\end{multline*}
where $\Gm$ acts via 
$u\cdot (x,y)
= (u^2 x, u^{-2}y).$
The closed substack $\cX(\sigma)$ is cut
out by the equation~$y = 0$ (with $x$, $t$, and $u$ then corresponding 
to the variables with the same names in
the description of~$\cX(\sigma)$ given by
  Proposition~\ref{prop: explicit description of FL stack}),
while $\cX(\sigmacomp)$ is  cut out  by~$x = 0$
(with $y$, $t^{-1}$, and $u^{-1}$ then corresponding to the variables $x$, $t$, and $u$
appearing in the description of $\cX(\sigmacomp)$ given by 
  Proposition~\ref{prop: explicit description of FL stack}).
In terms of this explicit description,
the map $\cX(\sigmasigmacomp) \to X(\sigmasigmacomp)$ is simply projection to the $\Pone$
factor.
\end{rem}

\begin{rem}\label{rem: justification of existence of morphism}
It can be checked that if~$\sigma_{a,b}$ is a Serre weight
  with~$b\ne 0,p-1$, then we have an isomorphism
\begin{multline*}
  \cX(\sigma_{a,b})\cup\cX(\sigma_{a+b,p-1-b}) 
\cong
\\
 \bigl(\Spec k[t_1,t_2]/(t_1t_2)\times [\Pone_x/(\Gm)_u]\bigr) 
\setminus ( \{t_1 = 0\} \times [0/\Gm] \cup  \{t_2 = 0\} \times [\infty/\Gm]) ,
\end{multline*}
where $\Gm$ acts on $\Pone$ via 
$u\cdot x = u^2 x$.

Here we identify the locus~$t_2=0$ (resp.\ $t_1=0$) with
  $\cX(\sigma_{a,b})$ (resp.\ $\cX(\sigma_{a+b,p-1-b})$) via
  Proposition~\ref{prop: explicit description of FL stack}, with $t_1$ corresponding
 to $t$ and $x$ corresponding to~$x$ (resp.~$t_2$ corresponding to $t$ and $x$
  corresponding to~$x^{-1}$).

\[\includegraphics[width=0.8\textwidth]{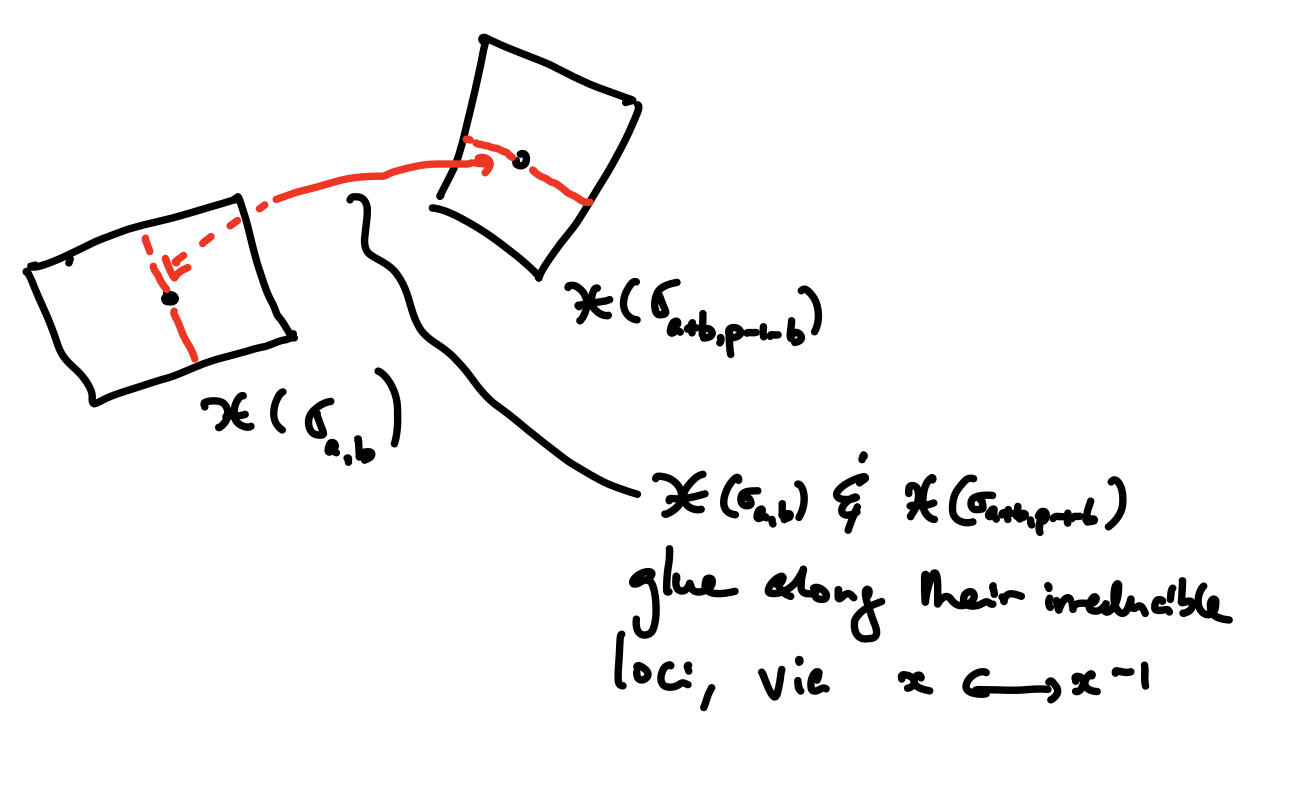}\]

  From this, it is easy to see
  that the above prescriptions for the morphisms $\pi_{\sigma}$ and $\pi_{\sigmacomp}$
  glue to give a morphism $\cX_{\red}\to X$.
\end{rem}

\begin{rem}\label{rem: what is the meaning of cX to X}%
  The motivation for the definition of~$X$ and of
  the morphism $\pi:\cX_{\red}\to X$ is as follows. Either by general
  principles, or by our explicit descriptions of the stacks
  $\cX(\sigma)$, one easily sees that the \emph{closed}
  $\Fpbar$-points of $\cX$ are precisely the semisimple
  $\rhobar:\Gal_{\Qp}\to\GL_2(\Fpbar)$ (with determinant
  $\varepsilonbar^{-1}$). It follows from the definition of
  ~$\pi$ that it induces a bijection on closed $\Fpbar$-points, and we
  think of~$X$ as a ``coarse'' moduli space\footnote{We put ``coarse'' in quotes
  because the morphism $\cX_{\red} \to X$ is not a coarse
  moduli space morphism in the technical sense of e.g.~\cite{MR3272912},
  or even an adequate moduli space morphism (for example, because it is not
  a closed map on underlying topological spaces).  It does, however, satisfy 
  the conditions of \cite[Prop.~7.1.1]{MR3272912}, and in particular
  is initial for morphisms from $\cX_{\red}$ to locally separated algebraic spaces.}
  for~$\cX_{\red}$.
  On $\Fpbar$-points, we can think
  of the morphism~$\pi$ as sending a representation to its
  semisimplification; more explicitly, the marked points of the
  component $X(\sigmasigmacomp)$ correspond to irreducible~$\rhobar$,
  and the~$\Gm$ obtained by deleting the marked points parameterizes
  sums of characters with fixed product, and with the Frobenius
  eigenvalue determined by the~$\Gm$.
\end{rem}

\begin{rem}
  \label{rem: existence of coarse moduli is unique to this case} For any choice of~$F$,
  it is easy to construct associated moduli space morphisms 
  \numequation
  \label{eqn:rank 1 moduli space morphism}
  (\cX_{F,1})_{\red} \to Y
  \end{equation}
  satisfying satisfying the conditions of~\cite[Prop.~7.1.1]{MR3272912}
  Indeed,  since $(\cX_{F,1})_{\red}$
  is a global quotient stack of the form $[\Spec A/\Gm]$, where $\Gm$ acts trivially,
  one sets $Y = \Spec A$, and the resulting morphism~\eqref{eqn:rank 1 moduli space morphism}
  is even a good moduli space morphism, in the sense of \cite{MR3237451}.

  Furthermore, since $\cX_{F,1}$ can be written in the form $[\Spf \widehat{A}/\Gmhat]$
  for some adic affine formal scheme $\Spf \widehat{A}$ which thickens $\Spec A$, endowed
  with the trivial $\Gmhat$-action,
  if we then set $\widehat{Y} := \Spf \widehat{A}^{\Gmhat}$,
  we obtain a morphism $\cX_{F,1} \to \widehat{Y}$
  which satisfies the formal algebraic analogue of the conditions of~\cite[Prop.~7.1.1]{MR3272912},
  and realizes $\widehat{Y}$ as the formal moduli space associated to~$\cX_{F,1}$ itself. 

Returning to the case of~$\GL_2(\Qp)$, as already mentioned, we anticipate that the morphism $\cX_{\red} \to X$
  constructed above also admits a formal thickening $\cX \to \widehat{X}$
  satisfying the formal algebraic analogue of the conditions of~\cite[Prop.~7.1.1]{MR3272912}.
  On the other hand, we do not expect there to exist interesting moduli space
  morphisms $\cX_{F,d} \to \widehat{Y}$, or $(\cX_{F,d})_{\red} \to  Y$,
  if $F \neq \Qp$ and $d \geq 2$, or even if $F = \Qp$ when $d > 2$.  More precisely,
  in these cases we expect that any morphism from $\cX_{F,d}$ to a locally separated formal algebraic space factors through the determinant map $\cX_{F,d}\to\cX_{F,1}$. %
\end{rem}

\subsubsection{Some open substacks}%
\label{subsubsec:the open substacks U}
Fix a very generic companion pair of
Serre weights $\sigma$, $\sigmacomp$, and write~$U(\sigmasigmacomp)$
to denote the open subset
of~$|X(\sigmasigmacomp)|$ obtained by deleting the marked points
$0,\infty$. Then $U(\sigmasigmacomp)$ is also open in $|X|$,
and we write
$\cU(\sigmasigmacomp)$, or sometimes just~$\cU$, to denote the preimage
of $U(\sigmasigmacomp)$, thought of an open substack of~$\cX$. 
 Then $\cU_{\red} = \cU(\sigmasigmacomp)_{\red}$
is the (dense) open substack
of~$\cX(\sigmasigmacomp)$ obtained by deleting the locus (consisting of two closed points)
corresponding to irreducible~$\rhobar$.
More precisely,  $\cU_{\red}$  is the union of its two closed substacks
$\cU\cap \cX(\sigma)$ and $\cU\cap \cX(\sigmacomp)$.

In terms of the description of~$\cX(\sigmasigmacomp)$ given in
Remark~\ref{rem:description of cX(sigmasigmacomp) in the generic case},
$\cU(\sigmasigmacomp)_{\red}$ is given by
$\bigl[ \big(\Gm \times \Spec k[x,y]/(xy) \bigr) / \Gm \bigr]
 = [\Spec k[t^{\pm 1},x,y]/(xy) / \Gm ]$ with $\Gm$-action given by
\[  u\cdot (x,y)   =   (u^2 x, u^{-2} y).\] %

\[\includegraphics[width=0.8\textwidth]{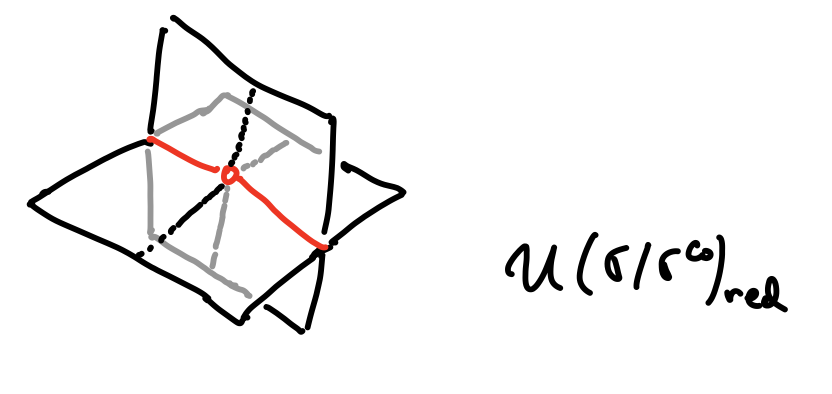}\]

By a deformation theory argument building on
this description of~$\cU_{\red}$,
we can explicitly describe~$\cU$
as %
\numequation\label{eqn: explicit description of U sigmasigmacomp}\cU(\sigmasigmacomp)=
    [\Spf \cO [s,t^{\pm 1},x,y]^{\wedge} / \Gmhat],\end{equation}
 where the hat denotes $(p,s,xy)$-adic completion, and
the $\Gmhat$-action is given by
\[  u\cdot (s,t,x,y)   =   (s,t, u^2 x, u^{-2} y).\]%
Here $\cU\cap\cX(\sigma)$ is the locus $\varpi=s=y=0$, and $\cU\cap\cX(\sigmacomp)$
is given by $\varpi=s=x=0$ (with the identification with the
description in Proposition~\ref{prop:
  explicit description of FL stack} being given by $y\mapsto x^{-1}$,
$t\mapsto t^{-1}$).  %

\[\includegraphics[width=0.8\textwidth]{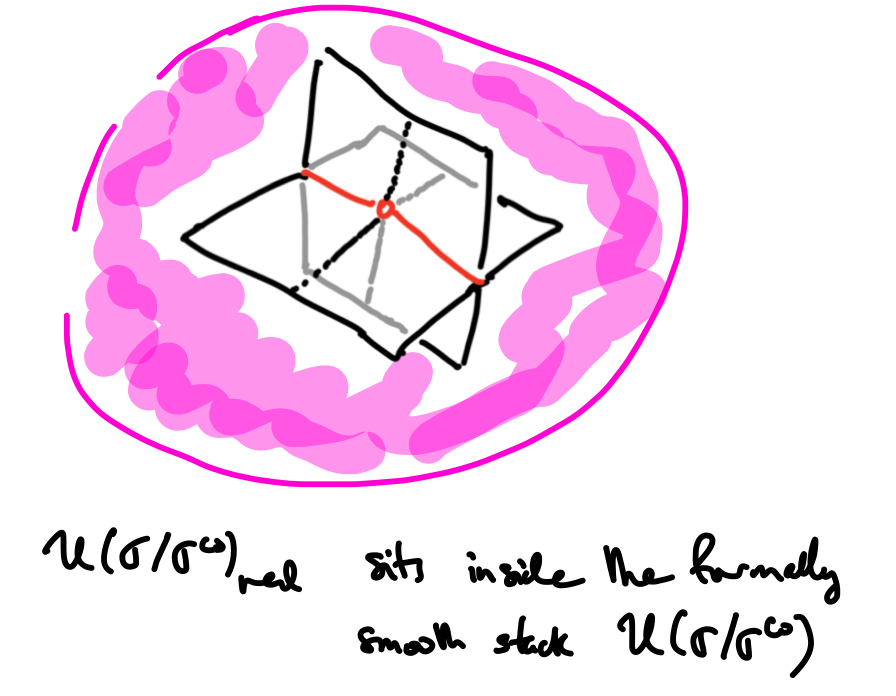}\]

\subsubsection{More open substacks}
\label{subsubsec:more opens}
As in Remark~\ref{rem: justification of existence of morphism},
consider a Serre weight~$\sigma_{a,b}$ with $b \ne 0,p-1$.
The union $\cX(\sigma_{a,b}) \cup \cX(\sigma_{a+b,p-1-b})$ then 
contains a unique closed point whose associated Galois representation
$\rhobar$ is irreducible.  If we delete the split locus from this union
(which, in terms of the explicit description of
Remark~\ref{rem: justification of existence of morphism},
means deleting the loci $x = 0$ and $x=\infty$),
we obtain an open subset $\cV_{\red}$ of $\cX_{\red}$ containing~$\rhobar,$
which admits the explicit description
$\cV_{\red} = [\Spec k[t_1,t_2]/(t_1t_2)/\boldsymbol{\mu}_2],$
where the $\boldsymbol{\mu}_2$-action is trivial.

If we let~$\cV$ denote the open substack of~$\cX$ corresponding
to  the open subset $\cV_{\red}$ of $\cX_{\red},$
then a deformation theory argument
shows that $\cV$ admits the description
$$\cV \cong [\Spf \cO[a,t_1,t_2]^\wedge/\boldsymbol{\mu}_2],$$
where the hat indicates $(p,a,t_1t_2)$-adic completion, 
and where the $\boldsymbol{\mu}_2$-action is again trivial.

\subsection{The Banach case for \texorpdfstring{$\GL_2(\Qp)$}{GL2(Qp)} --- II. Proof sketches} 
\label{subsec:GL2 Qp proof sketches}
In this section we outline a proof of Expected Theorem~\ref{expectedthm: DEG results}.
\subsubsection{The construction of $L_\infty$}%
\label{subsubsec:construction-l_infty}We begin by defining our
pro-coherent sheaf~$L_\infty$ of $\GL_2(\Qp)$-representations. In
essence this is a natural generalization of Colmez's
$\GL_2(\Qp)$-representation $D^\natural\boxtimes\Pone$. Some of the
arguments of~\cite{MR2642409} go through essentially unchanged in our
setting, but others need more serious adaptation. (Some of the main difficulties
that arise are that Colmez's $D^\sharp$ and $D^{\nr}$ are not compatible with flat
base change, so that all arguments involving them have to be reworked;
and that, when working with $(\varphi,\Gamma)$-modules with general coefficients,
quotients of lattices need not be literally finite modules,
so that arguments invoking such finiteness statements must also be reworked.)

To begin with, we note that
by definition there is a universal rank 2 $(\varphi,\Gamma)$-module $D$
over~$\cX$. We find it convenient to introduce a twist into its
definition, in the following way. If~ $A$ is a finite type
$\cO/\varpi^a$-algebra, then a morphism $\Spec A\to\cX$ corresponds to
a rank 2 $(\varphi,\Gamma)$-module  with $A$-coefficients and
determinant the inverse cyclotomic character, and the
formation of this $(\varphi,\Gamma)$-module is compatible with base
change in~$A$. We write~$D_A$ for the twist of this
$(\varphi,\Gamma)$-module by the cyclotomic character, so that~$D_A$
has determinant equal to the cyclotomic character (rather than its
inverse).

We work with $(\varphi,\Gamma)$-modules for the full cyclotomic extension
$\Qp(\zeta_{p^\infty})/\Qp$, so that $D_A$ is a rank 2 projective
module over the ring \[\A_{A}=A((T)),\] equipped with commuting($A$-linear)
semilinear actions of $\varphi$ and $\Gamma$,
where \[\varphi(1+T)=(1+T)^p,\]\[\gamma(1+T)=(1+T)^{\varepsilon(\gamma)}\]for
any $\gamma\in\Gamma$ (where $\varepsilon$ is as usual the $p$-adic
cyclotomic character). There is also a left inverse $\psi:D_A\to D_A$
to~$\varphi$; this is an $A$-linear surjection which commutes with~$\Gamma$,
and satisfies  \[\psi(\varphi(a)m)=a\psi(m), \] \[\psi(a\varphi(m))=\psi(a)m\]
  for any $a\in \A_{A}$, $m\in D_A$.

Following Colmez, we use the actions of~$\varphi$, $\psi$ and~$\Gamma$ on~$D_A$
to explicitly build a $\GL_2(\Qp)$-representation. %
Write $P^+$ for the monoid $
\big(\begin{smallmatrix}
  \Zp\setminus\{0\}&\Zp\\0&1
\end{smallmatrix}\big)
$, and~$P$ for the group~$
\big(\begin{smallmatrix}
  \Qptimes &\Qp\\0&1
\end{smallmatrix}\big)
$.
Then there is an
action of~$P^+$ on~$D_A$ given by the continuous maps\[
  \begin{pmatrix}
    p^ka&b\\0&1
  \end{pmatrix}z=(1+T)^b\varphi^k\circ\sigma_a(z)
\] for $a\in\Zptimes$, $b\in\Zp$,
$k\in\Z_{\ge 0}$, where we write~$\sigma_a\in\Gamma$ for the element with
$\varepsilon(\sigma_a)=a$.

There is a natural action of~$P^+$ on~$\Zp$ (respectively
of~$P$ on~$\Qp$) via $
\big(\begin{smallmatrix}
  a &b\\0&1
\end{smallmatrix}\big)x=ax+b $, and as explained in~\cite[\S
3.1]{ColmezMirabolique}, this action can be used to define %
a $P^+$-equivariant sheaf of
$A$-modules $U\mapsto D_A\boxtimes U$ on~$\Zp$. %
This is arranged by letting $D_A\boxtimes\Zp$ be~$D_A$ itself, %
and for any $i\in\Zp$ and~$k\ge 0$,
the sections of $D_A\boxtimes (i + p^k \Z_p)$
on $i+p^k\Zp$ are defined to be $
\big(\begin{smallmatrix}
  p^k&i\\0&1
\end{smallmatrix}\big)D_A\subseteq D_A
$. The restriction map $\Res_{i+p^k\Zp}:D_A\boxtimes\Zp\to
D_A\boxtimes(i+p^k\Zp)$ is given by \[
  \begin{pmatrix}
    1&i\\0&1
  \end{pmatrix}\circ\varphi^k\circ\psi^k\circ
  \begin{pmatrix}
    1&-i\\0&1
  \end{pmatrix}.\] We thus have in particular~$D_A\boxtimes\Zptimes=D_A^{\psi=0}$.

Given an affine linear map $f:U\to V$ between compact open subsets
of~$\Zp$ we have an induced map $f_*:M\boxtimes U\to M\boxtimes
V$. This definition extends immediately to the case of continuous maps
which are locally affine linear, and by taking limits, to functions
which can be suitably approximated by affine linear maps. In
particular this applies to the map $z\to 1/z$ on $\Zptimes$, and we
can then glue $D_A \boxtimes\Zp := D_A$ to itself over $D_A\boxtimes\Zptimes$ via
this map, giving a sheaf over~$\Pone(\Zp) (= \Pone(\Q_p) = \Q_p \cup  \{\infty\})$.
We denote the global
sections of this sheaf by $D_A \boxtimes\Pone$; by functoriality, they
have an action of $\GL_2(\Qp)$ with trivial central character.

Since $\Q_p$ is an open subset of~$\Pone(\Z_p)$, we can also take the
sections of $D_A\boxtimes\Pone$ over $\Q_p$, to obtain a~$B$-representation
that is denoted~$D_A\boxtimes\Q_p$.  Unwinding the definitions, one obtains
an explicit description of $D_A\boxtimes\Qp$ as the inverse limit
$\varprojlim_\psi D_A$.

The action of $\GL_2(\Qp)$ on $D_A\boxtimes\Pone$ is
compatible with flat base change in~$A$, but $D_A\boxtimes\Pone$ is not the
sheaf~$L_\infty$ that we are seeking. This is already apparent in the
original setting of~\cite{MR2642409}, where $D_A\boxtimes\Pone$ does not
realise the local Langlands correspondence, but rather an extension of
the correspondence by its (twisted) dual.
(See Section~\ref{subsubsec:duality for }
below for an expected extension of this phenomenon to our setting.) We
will therefore construct  $L_\infty$ (again following Colmez) as a
subrepresentation of $D_A\boxtimes\Pone$. 

In order to do so, we set~$\A_A^+=A[[T]]$, which inherits actions
of~$\varphi$, $\psi$ and~$\Gamma$ from~$\A_A$. It can be shown that
$D_A$ contains a minimal $\psi$-stable $\A_A^+$-lattice $D_A^\natural$
(it is a lattice in the sense that it is finitely generated
over~$\A_A^+$ and spans $D_A$ over $\A_A$), whose formation is
compatible with flat base change. %
There is a natural %
restriction map
\[\Res_{\Qp}:D_A\boxtimes\Pone\to D_A\boxtimes\Qp,\]
and we set \[D_A^\natural\boxtimes\Pone :=\{z\in D_A\boxtimes\Pone,\
  \Res_{\Qp}(z)\in D_A^\natural\boxtimes\Qp\},\] where
$D_A^\natural\boxtimes\Qp:=\varprojlim_\psi D_A^\natural$. It turns
out that $D_A^\natural\boxtimes\Pone$ is a
$\GL_2(\Qp)$-stable lattice in $D_A\boxtimes\Pone$, but this is
not at all obvious; indeed it is not at all obvious (even in the
original setting of~\cite{MR2642409}) that
$D_A^\natural\boxtimes\Pone$ is open in $D_A\boxtimes\Pone$. (Here ``lattice''
is in the sense of~\cite{MR2181808}: a lattice in a Tate module
over~$A$ is an open $A$-submodule~$L$ with the property that if
$U\subseteq L$ is open, then $L/U$ is a finitely generated $A$-module.)

The proof that $D_A^\natural\boxtimes\Pone$ is a $\GL_2(\Qp)$-stable lattice is complicated; it
involves a reduction (via consideration of complete local rings) to
the case that~$A$ is Artinian, and then
as in ~\cite{MR2642409} it is ultimately proved via a further reduction to the
explicit description of the $p$-adic local Langlands correspondence
for 2-dimensional crystalline representations of $\Gal_{\Qp}$ proved by
Berger--Breuil ~\cite{MR2642406}. Having established this, we have the
following useful corollary.

\begin{cor}[Dotto--E.--G.]%
  \label{cor: D natural is f.g. over K group ring} Suppose that~$p\ge 3$. Assume that~$A$ is a finite type
  $\cO/\varpi^a$-algebra for some~$a\ge 1$, and that~$D_A$ is a
  projective \'etale $(\varphi,\Gamma)$-module of rank~$2$ with
  $A$-coefficients and determinant~$\varepsilon$. Then $D_A^\natural\boxtimes\Pone$ is a
 finite $A[[K]]$-module.
\end{cor}
\begin{proof}
 This is deduced from the stability of $D_A^\natural\boxtimes\Pone$ under
 the action of~$\GL_2(\Qp)$, together with the Noetherianness of $A[[K]]$, and
 the fact  that $D_A^\natural$,
  being a lattice, is a finitely generated $\A^+_{A}$-module.
\end{proof}

It can then be shown that the formation of
$D_A^\natural\boxtimes\Pone$ is compatible with flat base change
in~$A$, and we use this to construct a pro-coherent sheaf~$L_\infty$
over~$\cX$, equipped with a continuous
  $\GL_2(\Qp)$-action. %

\subsubsection{A reinterpretation of some results of Pa\v{s}k\={u}nas}%
\label{subsubsec:vytas}We now recall some of the main results
of~\cite{MR3150248}, and reinterpret them in terms of our
functor~$\mathfrak{A}$. Let  $(\smG)^{\ladm}$ denote the full subcategory
of~$\smG$ 
consisting of {\em locally admissible}
representations. (By definition, a
representation~$\pi\in \smG$ is locally admissible if every vector~$v\in\pi$
is smooth and generates an admissible representation; in fact, for
$\GL_2(\Qp)$, locally admissible representations are automatically locally
finite, so that for every~$v$, the representation generated by~$v$ is of
finite length.)

Unlike~$\smG$, the category $(\smG)^{\ladm}$ admits a block decomposition,
as we now explain. By definition, a block of~$(\smG)^{\ladm}$
is an equivalence class of (isomorphism classes of) irreducible
objects 
under the equivalence relation generated by
\[
\pi_1 \sim \pi_2 \text{ if } \Ext^1_{\smG}(\pi_1, \pi_2) \ne 0 \text{ or } \Ext^1_{\smG}(\pi_2, \pi_1) \ne 0.
\]
Strictly speaking, we should have written ``irreducible admissible objects'',
but in fact the irreducible objects of $\smG$ are all contained in $(\smG)^{\ladm}$.
Indeed,
the absolutely irreducible smooth representations of~$\GL_2(\Q_p)$
in characteristic~$p$
admitting a central character\footnote{Subsequently, Berger showed~\cite{MR3076828}
that irreducible smooth representations of  $\GL_2(\Q_p)$ over $\Fbar_p$
necessarily admit central characters.  In our particular context, we are interested
in representations of $\PGL_2(\Q_p)$, which is to say of $\GL_2(\Q_p)$ having
trivial central character, and so this issue with central characters does not arise.}
have been classified by
Barthel--Livn\'e  and Breuil~\cite{BarthelLivneDuke, BreuilGL2I}, %
and in particular have been shown to be admissible.
The non-absolutely
irreducible case is rather easily reduced to the absolutely
irreducible case (see e.g.\ \cite[Prop.~1.2]{MR3076828} and~\cite[\S 5.13]{MR3150248}).

The blocks~$\mathfrak{B}$ are as follows:
\begin{enumerate}
\item $\fB = \{\pi\}$ for a supersingular irreducible representation~$\pi$
having trivial central character,
\item $\fB = \{\Ind_B^G(\chi \otimes \chi^{-1}), \Ind_B^G(\omega\chi^{-1} \otimes \omega^{-1} \chi)\}$, for some character~$\chi: \bQ_p^\times \to k^\times$ such that $\chi^{2} \ne 1, \omega^{2}$, %
\item $\fB = \{\chi, \chi \otimes \St_G, \Ind_B^G(\omega\chi \otimes
  \omega^{-1} \chi)\}$ for a quadratic character $\chi : \bQ_p^\times \to k^\times$, and
\item\label{item: not absolutely irreducible blocks} blocks that do not contain absolutely irreducible objects.
\end{enumerate}

A key observation for us is that there is a bijection between the blocks
of types (1)--(3) and $X(k)$. (For the sake of exposition we ignore the blocks of
type~\eqref{item: not absolutely irreducible blocks}; they are easily
handled by allowing finite extensions of our coefficient field~$k$,
and indeed correspond to points of~$X$ over extensions of~$k$.) To
describe this bijection, recall that for any Serre weight~$\sigma$,
the Hecke algebra \[\cH(\sigma) := \End_G(\cInd_{K}^G \sigma
  )\] is
equal to $k[T]$ for an explicit Hecke operator~$T$, and the
absolutely irreducible representations are all subquotients
of the representations \[(\cInd_{K}^G\sigma)/(T-\lambda)\]for ~$\lambda\in
k$. Furthermore these representations are (absolutely) irreducible
unless~$\sigma=\sigma_{a,0}$ or~$\sigma=\sigma_{a,p-1}$ for some~$a$
and~$\lambda=\pm 1$, in which case there are two Jordan--H\"older
factors, which are of the form $\chi, \chi \otimes \St_G$ for
some~$\chi$.

The bijection between the blocks and~$X(k)$ relies on regarding the
Hecke operator~$T$ as a coordinate on an~$\A^1$ in the
component~$X(\sigmasigmacomp)$ of~$X$. We do this as follows,
restricting as above to the generic case for the sake of exposition. %

\begin{defn}\label{defn: coordinate functions on X from Hecke algebra}
Let~$\sigma$ be a generic Serre weight. %
Then we defined above a morphism $\cX(\sigmasigmacomp)\to
X(\sigmasigmacomp)\subset X$ to be given by one of~$\pi_\sigma$
or~$\pi_{\sigmacomp}$, and we fixed a coordinate~$T$
on~$X(\sigmasigmacomp)$ (a copy of~$\Pone$). %
In the case that this morphism is given
 by~$\pi_\sigma$,
we define $f_\sigma: \Spec \cH(\sigma) \to X$ by
sending $T\to T$; and in the case that this morphism is given
by~$\pi_{\sigmacomp}$, we send $T\mapsto T^{-1}$.
\end{defn}

Then if $x$ is a $k$-point of~$X$, the corresponding block is given
by 
\[
\fB_x =
\bigcup_{\sigma \text{ s.t.\ } x \in f_{\sigma}(\bA^1)} \JH\left ( \cInd_{K}^G(\sigma)
\otimes_{\cH(\sigma), x} k  \right ). 
\]
More precisely, $\fB_x$ is a block of~$(\smG)_k$, and the map
$x \mapsto \fB_x$ (suitably extended to the case of non-generic Serre weights) is a bijection from~$X(k)$ to the set of blocks
of~$\smG$ containing absolutely irreducible representations.

If $\fB$ is a block of $(\smG)^{\ladm},$ then we write
$(\smG)^{\ladm}_{\fB}$ to denote the full subcategory of $(\smG)^{\ladm}$ whose
objects  are those representations  all of whose irreducible subquotients  lie in the
given block~$\fB$.  We then have the  block decomposition
$$(\smG)^{\ladm} \iso  \prod_{\fB} (\smG)^{\ladm}_{\fB}.$$

In~\cite{MR3150248} Pa\v{s}k\={u}nas establishes equivalences of
categories between the various categories~$(\smG)^{\ladm}_{\fB}$,
and categories of modules over certain pseudodeformation rings.
We use Pa\v{s}k\={u}nas's results (and their
proofs) as an input to the proof of Theorem~\ref{expectedthm: DEG
  results}, but in order to do so it is necessary to reformulate these results
(and in some cases slightly extend them) as statements about fully
faithful functors to categories of sheaves on stacks. %

Let $\thetabar:\Gal_{\Qp}\to\GL_2(k)$ be a semisimple representation
(equivalently, a pseudorepresentation), 
corresponding to a point~$x \in X(k)$.
Assume (for simplicity of exposition) that the associated block~$\fB_x$
lies in one of cases~(1) or (2) above; that is, we assume that~$x$  is
not a point~$T=\pm 1$ of~$X(\sigmasigmacomp)$ for some
$\{\sigma,\sigmacomp\}$ of the form
$\{\sigma_{a,0},\sigma_{a+1,p-3}\}$.
We write 
$$(\smG)^{\ladm}_{\thetabar} := (\smG)^{\ladm}_{\fB_x}.$$

We let  $\cX_{\thetabar}$ denote Carl Wang-Erickson's stack of Galois
representations associated to~$\thetabar$. By definition,
$\cX_{\thetabar}$ is characterised by the following property: if~$A$ is an
$\cO$-algebra in which $p$ is nilpotent, then $\cX_{\thetabar}(A)$ is the
groupoid of continuous morphisms \[\rho:\Gal_{\Qpbar}\to\GL_2(A)\]with
$\det\rho=\varepsilon^{-1}$ and are such that the pseudorepresentation
associated to $\rho$ modulo~$\varpi$ is equal to~$\thetabar$ (here~$A$
has the discrete topology, and~$\Gal_{\Qp}$ its natural profinite
topology).

Wang-Erickson showed in~\cite{MR3831282} that~$\cX_{\thetabar}$ is a
formal algebraic stack, and by~\cite[Thm.\ 6.7.2]{emertongeepicture}
and~\cite[Thm.\ 10.2.2]{emerton2020moduli}, there is a natural
monomorphism $\cX_{\thetabar}\into\cX$ which induces a closed
immersion on underlying topological spaces (with image corresponding
to those representations $\rhobar:\Gal_{\Qp}\to\GL_2(\Fpbar)$ with
semisimplification isomorphic to~$\thetabar$). This morphism is versal
at the finite type points of $\cX_{\thetabar}$, so we naturally expect
that it is some kind of completion of the formal algebraic stack~$\cX$
at the closed substack $\cZ_{\thetabar}$ of~$\cX_{\red}$ whose
$\Fpbar$-points are those~$\rhobar$ with
semisimplification~$\thetabar$.

This expectation can be made precise by using the language of coherent
completeness~\cite[Defn.\ 2.1]{MR4088350} in the following
way. %
Write
$\cI_{\thetabar}$ for the ideal sheaf of~$\cZ_{\thetabar}$ in~$\cX$,
and $\cX_{\thetabar}^{[n]}$ for the closed substack of $\cX$ cut out
by $\cI_{\thetabar}^{n+1}$. Then it can be shown that for each~$n$,
$\cX_{\thetabar}^{[n]}$ is a closed substack of~$\cX_{\thetabar}$, and
that the natural functor (given by
restriction)
\[\Coh(\cX_{\thetabar})\to\varprojlim_n\Coh(\cX_{\thetabar}^{[n]})\] is
an equivalence of categories; that is, a coherent sheaf
on~$\cX_{\thetabar}$ is the same thing as a compatible system of
coherent sheaves on the infinitesimal neighbourhoods
of~$\cZ_{\thetabar}$ in~$\cX$.

Given the above constructions, the following theorem is essentially a
matter of unwinding some of the proofs of the main results
of~\cite{MR3150248}.
Let $(\smG)_{\thetabar}^{\fgadm}$ denote the full
subcategory of $(\smG)^{\ladm}_{\thetabar}$ consisting of finitely generated representations. 
\begin{thm}[Pa\v{s}k\={u}nas+$\varepsilon$]
  \label{thm: Vytas theorem}Suppose as above that
  $\thetabar:\Gal_{\Qp}\to\GL_2(k)$ is a pseudorepresentation, 
corresponding to a block %
which
lies in one of cases~(1) or (2). Then the functor \numequation\label{eqn:
    functor restricted to Vytas}(L_{\infty})_{| \cX_{\thetabar}} \otimes_{\cO[[G]]} \text{--}\end{equation}
is an exact fully faithful functor from $(\smG)^{\fgadm}_{\thetabar}$ %
to the category of coherent $\cO_{\cX_{\thetabar}}$-modules. %
\end{thm}%
\begin{proof}%
Let $R_{\thetabar}^{\ps}$ denote the complete Noetherian local
$\cO$-algebra parameterizing deformations of $\thetabar$ of
fixed  determinant~$\varepsilon^{-1}$ over finite type
Artinian local $\cO$-algebras, and let $\theta^u$ be the universal deformation of $\thetabar$
over $R_{\thetabar}^{\ps}.$  
We write~$\widetilde{R}_{\thetabar}$ for the
Cayley--Hamilton algebra defined as follows: 
let $J$ denote the closure of the two-sided ideal of $R_{\thetabar}^{\ps}[[\Gal_{\Q_p}]]$ generated
by the elements $g^2  - \theta^u(g)g + (\varepsilon^{-1})(g),$ for
$g \in \Gal_{\Q_p}$,  %
and then set
$$\widetilde{R}_{\thetabar}
:= R_{\thetabar}^{\ps}[[\Gal_{\Q_p}]]/J.$$

  As in Section~\ref{sec: TW patching}, we let~$V$ denote Colmez's Montreal functor. %
Pa\v{s}k\={u}nas shows in~\cite{MR3150248} that the ring $\widetilde{R}_{\thetabar}$
is free of rank~$4$ as an $R_{\thetabar}^{\ps}$-module,
and the results of~\cite{MR3150248} can be interpreted as showing that
the functor $\pi \mapsto V(\pi)$ induces an equivalence
between $(\smG)_{\thetabar}^{\fgadm}$ and the category of finite length~$\widetilde{R}_{\thetabar}$-modules.

  There is a universal two-dimensional $\Gal_{\Q_p}$-representation~$\cV$
lying over~$\cX_{\thetabar}$,
which can be regarded as a rank $2$ locally free sheaf of $\cO_{\cX_{\thetabar}}$-modules
equipped with an %
action
of $\widetilde{R}_{\thetabar}$. %
Then by using the relationships between $D^\natural\boxtimes\Pone$
  and the Montreal functor~$V$ established in~\cite{MR2642409}, one
  can show that the functor~\eqref{eqn:
    functor restricted to Vytas} is naturally equivalent to the
  functor %
  $$\cV
   \otimes_{\widetilde{R}_{\thetabar}}
  V(\text{--}),$$ and it is now easy to deduce the result from the
  theorems proved in~\cite{MR3150248}.
\end{proof}

The case of a block of type~(3) (i.e.\ the case that~$x$ is a
point~$T=\pm 1$ of~$X(\sigmasigmacomp)$ for 
$\{\sigma,\sigmacomp\}$ of the form
$\{\sigma_{a,0},\sigma_{a+1,p-3}\}$) is more subtle, as in contrast to
the other blocks we see derived phenomena. This case has recently been
studied by Christian Johansson, James Newton and Carl Wang-Erickson~\cite{johansson2024modulistacksgaloisrepresentations},
and their results can be used to prove the following theorem.

\begin{thm}[Johansson--Newton--Wang-Erickson+$\varepsilon$]
  \label{thm: JNWE theorem}Suppose as above that
  $\thetabar:\Gal_{\Qp}\to\GL_2(k)$ is a pseudorepresentation,
  corresponding to a block $\fB_{\thetabar}=\fB_x$ which lies in
  case~(3). Then the functor \numequation\label{eqn: functor
    restricted to Vytas two}(L_{\infty})_{| \cX_{\thetabar}}
  \otimes_{\cO[[G]]}^{L} \text{--}\end{equation} is a fully
faithful functor with amplitude~$[-1,0]$ %
from $D(\fB_{\thetabar}^{\fg})$ to the bounded
derived category of coherent
$\cO_{\cX_{\thetabar}}$-modules. %
\end{thm}%

\subsubsection{A useful resolution}
\label{sec:useful-resolution}
Let~$N$ be the normalizer of the usual (upper triangular) Iwahori
subgroup~$\Iw$, and let~$\delta$ be the nontrivial quadratic character
of~$N/Z\Iw$. If $\pi$ is any smooth $G$-representation, then we have a 
complex of $\cO[[G]]$-modules
\numequation\label{tree}
0 \to \cInd_N^G(\delta\pi) \to \cInd_{K}^G(\pi) \to \pi \to 0
\end{equation}
where the third arrow is the tautological surjection
$\cInd_{K}^G \pi \to \pi$. The second arrow is defined as follows.  Write $\Pi=
\big(\begin{smallmatrix}
  0&1\\p&0
\end{smallmatrix}\big)
$. 
Then the morphism
$$v \mapsto  \Pi \otimes \Pi^{-1} v - 1\otimes v \in \cO[G]\otimes_{\cO[K]}
\pi = \cInd_{K}^G \pi$$
is an $N$-equivariant 
morphism $\delta \pi \to \cInd_{K}^G \pi$
which induces the second arrow $\cInd_{N}^G \delta\pi \to
\cInd_{K}^G \pi.$ It is easy to show that the
 complex~\eqref{tree} is acylic.

\begin{rem}\label{rem: how we use the resolution}
The resolution~\eqref{tree} is valid for $\PGL_2(F)$ for any $F/\Qp$,
and can be used  to reduce certain assertions
  about general~$\pi$ to the case $\pi=\cInd_{K}^G \sigma$ for a
  Serre weight~$\sigma$. %
  In particular, we can reduce the evaluation of continuous functors
  (such as our hypothetical functor~$\fA$) on general~$\pi$ to the
  evaluation on $\cInd_{K}^G \sigma$ in the following way.
  Since we are assuming
  that $p>2$, $\cInd_{N}^G\delta\pi$ is a direct summand of 
$\cInd_{\Iw}^G \pi = \cInd_{K}^G (\cInd_{\Iw}^{K} \pi)$, so we
can reduce to studying representations  of the form $\cInd_{K}^GV$.

Since compact induction is compatible with passing to filtered colimits, 
we then reduce to the case when $V$ is finitely generated.  Since
compact induction is exact, we furthermore reduce to the case when $\pi$ is of the
form $\cInd_{K}^G \sigma$ for some Serre weight~$\sigma$, as claimed.
\end{rem}

\subsubsection{The definition of the functor}
\label{sec:definition-functor}%
We are now in a position to define our functor. Since we have already
defined the pro-coherent sheaf~$L_\infty$ of
$\cO[[G]]$-representations over~$\cX$, %
we would simply like to define this via
\numequation\label{eqn: defn of the functor for GL2 Qp}\pi \mapsto
L_{\infty}\otimes^{L}_{\cO[[G]]}\pi.\end{equation} However, it
is not a priori
obvious that this definition is well-behaved, and we proceed in
stages. Firstly, if~$\pi$ is of the form~$\cInd_{K}^GV$ with~$V$
finitely generated, then we note
that \[L_{\infty}\otimes^{L}_{\cO[[G]]}\cInd_{K}^GV=L_{\infty}\otimes^{L}_{\cO[[K]]}V,\]
which can be shown (using Corollary~\ref{cor: D natural is f.g. over K
  group ring}) to be a bounded-above complex of coherent
$\cO_\cX$-modules. We claim that it is in fact a sheaf, concentrated in
degree zero. To see this, it is enough (by descending induction on the
cohomological degree, and using the coherence of the cohomology
groups) to check after restriction to each substack~$\cX_{\thetabar}$,
where the result follows from Theorem~\ref{thm: Vytas
  theorem}. %

Passing to colimits, we see that for any~$V$ (finitely generated or
otherwise), $L_{\infty}\otimes^{L}_{\cO[[G]]}\cInd_{K}^GV$ is a
quasicoherent sheaf in degree zero. By Remark~\ref{rem: how we use
  the resolution}, we deduce that for any smooth representation~$\pi$,
$L_{\infty}\otimes^{L}_{\cO[[G]]}\pi$ is a complex of
quasicoherent sheaves concentrated in degrees~$[-1,0]$, and that
if~$\pi$ is finitely presented (equivalently, finitely generated),
then this complex has coherent cohomology sheaves. We conclude that
$L_{\infty}\otimes^{L}_{\cO[[G]]}-$ defines a functor with
amplitude in $[-1,0]$ and takes bounded complexes of finitely
presented representations to bounded coherent complexes of $\cO_\cX$-modules.

\subsubsection{Localization and gluing}
\label{subsubsec:localization}%

If~$\sigma$ is a Serre weight, we heuristically regard $\cInd_{K}^G \sigma$ as ``living over'' a 
copy of $\A^1 = \Spec k[T] = \Spec \cH(\sigma)$ in~$X$ via
Definition~\ref{defn: coordinate functions on X from Hecke
  algebra}. This intuition is the basis of the localization theory
for~$\smG$ developed in~\cite{DEGlocalization}, and we now explain how
to make it precise. The results described in this section, and their
proofs, are for the most part 
purely representation-theoretic (with some occasional uses of the
classification of blocks, which is proved in~\cite{MR3150248} using
the $p$-adic local Langlands correspondence), although they are of course motivated
by the geometric picture.

If $Y$ %
is a closed subset
of~$X$,
we let $(\smG)_Y$ denote the subcategory of 
$\smG$ 
consisting of those representations all of whose irreducible
subquotients lie in a block corresponding to a closed point of $Y$. If $U$ is an open subset of $X$, then we write \[(\smG)_U := \smG/(\smG)_Y,\]
where $Y := X \setminus U$.
It is not hard to check that the subcategory~$(\smG)_Y$ is localizing in
the usual sense: namely, there is a fully faithful right adjoint \[(j_U)_*: (\smG)_U \to \smG\]
to the canonical quotient functor $(j_U)^*:\smG \to
(\smG)_U$. Where~$U$ is understood, we write~$j_*$ for~$(j_U)_*$,
and~$j^*$ for~$(j_U)^*$. %

The following result, which is~\cite[Prop.\ 3.1.7]{DEGlocalization}, confirms some basic intuitions about how the
functors $j_*,j^*$ should behave.

\begin{prop}%
\label{prop:localizing compact inductions} 
Let $Y$ be a closed subset of $X$, and write 
$U~:=~X~\setminus~Y$.  

\begin{enumerate}
\item If $\sigma$ is a Serre weight, 
and if $f_{\sigma}^{-1}(Y) = V(g)$
{\em (}a closed subset of $\Spec \cH(\sigma)${\em )}
for some $g \in \cH(\sigma)$,
then the natural map
$$\cInd_{K}^G \sigma \rightarrow
(\cInd_{K}^G \sigma)[1/g]$$
can be identified with the morphism
$$\cInd_{K}^G \sigma \to
j_*j^* \cInd_{K}^G \sigma.$$
\item The functor $j_*$ is exact and commutes with filtered colimits.
\end{enumerate}
\end{prop}

If $\{U_0,\ldots,U_n\}$ is any finite open cover of $X$,
then for any object $\pi$ of $\smG$,
we obtain a functorial \v{C}ech resolution  %
\numequation
\label{eqn:Cech resolution again}
0 \to \pi \to \prod_i (j_i)_*(j_i)^* \pi \to \cdots \to
(j_{0,\ldots,n})_*(j_{0,\ldots,n})^* \pi \to 0
\end{equation}
where as usual we write $U_{i\dots k}$
for~$U_i\cap\dots\cap U_k$, we have written $j_{i\dots k}$
for $j_{U_{i\dots k}}$, and the differentials are given by
the usual formulas. The following is ~\cite[Prop.\ 3.2.3]{DEGlocalization}. %

\begin{prop}%
\label{prop:Cech acyclicity}
For any object $\pi$ of~$\smG$,
the resolution~{\em \eqref{eqn:Cech resolution again}}
is acyclic.
\end{prop}
\begin{proof}Using the resolution~\eqref{tree}, we reduce as in
  Remark~\ref{rem: how we use the resolution} to the case
  that $\pi=\cInd_{K}^G\sigma$ for a Serre weight~$\sigma$. In this
  case Proposition~\ref{prop:localizing compact inductions} allows us
  to reduce to the standard assertion that the \v{C}ech complex
  associated to a finite open cover of an affine scheme by
  distinguished opens is acyclic.
\end{proof}

From Proposition~\ref{prop:Cech acyclicity}, we deduce that~$\smG$ is a
stack in abelian categories over~$X$, as claimed in Theorem~\ref{expectedthm:
  Bernstein centre}. %
We need a gluing statement for categories of representations that corresponds to Beauville--Laszlo gluing
of sheaves. Write $(\smG)^{\fg}$ %
for the full subcategory of $\smG$ consisting
of finitely generated objects.  Similarly, for any closed subset $Y$ of~$X$,
we write $(\smG)_Y^{\fg}$ to denote the full subcategory of $(\smG)_Y$ consisting
of finitely generated objects. We may form the category of formal pro-objects $\pro((\smG)_Y^{\fg})$.

\begin{df} If $\pi$ is an object of~$(\smG)^{\fg}$, 
and if $Y$ is a %
closed subset of~$X$,
then we write $\widehat{\pi}_Y$ for the object  $\quoteslim{I} \pi'$
 of $\pro((\smG)_Y^{\fg})$, where $\pi'$ runs over all quotients of $\pi$ which lie in $(\smG)_Y^{\fg}$.
\end{df}

The following is ~\cite[Thm.\ 3.8.1]{DEGlocalization}.
\begin{thm}%
\label{thm:BL gluing}If~$Y$ is finite, then the functor $\pi \mapsto \widehat{\pi}_Y$ is
exact, and the canonical functor
$$(\smG)^{\fg} \, \longrightarrow \, \pro((\smG)_Y^{\fg}) \times_{(\pro((\smG)_Y^{\fg})_U)} (\smG)_U^{\fg},$$
induced by completion along $Y$ and by localization over~$U$,
is an equivalence of categories.
\end{thm}

\subsubsection{A description of $\fA$ generically} 
\label{subsubsec:restriction to U}
Given an open subset~$U$ of~$|X|$, we
have the localized category $(\smG)_U$, and an open
substack ~$\cU:=\pi^{-1}(U)$ of~$\cX$. We now take~$U$ to be the
complement of the finite set~$Y$ consisting of points corresponding to blocks of types (1) or (3);
in particular we have deleted the marked points.
Then~$U$ is dense in~$X$, and as we now explain, %
it is in fact possible to explicitly describe the
restriction of~$\fA$ to $(\smG)_U$, and to show that it
is a fully faithful functor to~$D(\cU)$. The rough idea is to follow the strategy of~\cite{MR3150248} for blocks of
type~(2), but to ``algebraize in the unramified direction''.  %

The open set $U$ is a disjoint union of its intersections with the various
irreducible components of~$|X|$, and we handle each of these intersections in turn.
In the present discussion we restrict ourselves to describing the case
of $U \cap X(\sigmasigmacomp)$ when $\sigma,\sigmacomp$ is a companion pair of very
generic Serre weights.
This intersection was denoted $U(\sigmasigmacomp)$ in Section~\ref{subsubsec:the open substacks U},
but to simplify notation, we now redefine $U$ as $U := U(\sigmasigmacomp)$.
The preimage of $U$ in $\cX$ is the open substack $\cU: = \cU(\sigmasigmacomp)$.

Our aim, then, is to describe the composition of $\fA$ with restriction to~$\cU$,
which turns out to factor through~$(\smG)_U$. 
Our description of this composition will be Morita-theoretic, in terms
of a projective pro-generator for the category $(\smG)_U$,
  and the simplest way to construct this pro-generator in terms of the
sheaf $L_{\infty}$
that we have used to define~$\fA$.
Actually, we find it convenient to use a dual version~$L^\vee_\infty$ of~$L_\infty$, which can
be constructed in a similar way to~$L_\infty$ itself (see also
Section~\ref{subsubsec:duality for } below). %

We now write \[\cU\setminus \cX(\sigmacomp)= [ (
\Spf \cO [s,t^{\pm 1},x,y^{\pm 1}]^{\wedge} / \Gmhat
]=[\Spf A/\boldsymbol{\mu}_2],\]
where $A=\cO [s,t^{\pm 1},x]^{\wedge}$ (the hat denoting
$(p,s,x)$-adic completion), and the action of~$\boldsymbol{\mu}_2$ is
trivial. %
We can regard $L_\infty^\vee|_{\cU\setminus
  \cX(\sigmacomp)}$ as topological $A[[G]]$-module $P_\alpha$ with an action
of~$\boldsymbol{\mu}_2$. With some work, and a comparison to the results
of~\cite{MR3150248}, we can show that $P_\alpha$ is a ``pro-projective''
object of $\pro\bigl((\smG)_U\bigr)$,
meaning that if $\pi \to \pi'$ is an epimorphism in~$(\smG)_U$,
then $\Hom(P_{\alpha},\pi) \to \Hom(P_{\alpha},\pi')$ is surjective.
(Because projective limits of surjections of abelian groups need not
be surjective, this does not imply that  $P_{\alpha}$ is projective
in~$\pro\bigl((\smG)_U\bigr).$)

 Reversing the roles of $\sigma$ and
$\sigmacomp$, we obtain another pro-projective object~$P_\beta$. 
Together, $P_\alpha$ and $P_\beta$ generate~$(\smG)_U$.
By computing their completions at closed points using Theorem~\ref{thm: Vytas theorem},
we can show that $\mathfrak{A}(P_\alpha)$ and
$\mathfrak{A}(P_\beta)$ are certain explicit line bundles on~$\cU$. %
In fact, with an appropriate interpretation of the notation, these are the line bundles corresponding
to the graded modules~$B(1)$ and~$B(-1)$ considered in
Example~\ref{ex: Morita example for GL2 Qp}, 
 and the full 
 faithfulness of the functor
$\mathfrak{A}:(\smG)_U\to D(\cU)$ 
 is a straightforward
consequence of Proposition~\ref{prop:morita}.  %

\subsubsection{Full faithfulness of the functor}
\label{subsubsec:full faithful GL2 Qp}
As at the beginning of Section~\ref{subsubsec:restriction to U},
let~$U$ denote the
complement of the finite set~$Y$ consisting of points corresponding to blocks of types (1) or (3).
Ignoring as usual the issue of blocks of type (3)
(which require some additional arguments), the full faithfulness
of~$\mathfrak{A}$ can be deduced
from the full faithfulness of $\mathfrak{A}:(\smG)_U\to
D(\cU)$, together with  Theorem~\ref{thm:BL gluing} and
Theorem~\ref{thm: Vytas theorem}. 

\subsubsection{Compatibility with duality}
\label{subsubsec:duality for }%
For any finite type
$\cO/\varpi^a$-algebra~$A$, there is an obvious short exact sequence
of $\GL_2(\Qp)$-representations
\numequation\label{eqn: Colmez ses}0 \to D_A^{\natural} \boxtimes \Pj^1 \to D_A\boxtimes \Pj^1 \to
(D_A\boxtimes \Pj^1)/(D_A^{\natural}\boxtimes \Pj^1) \to 0.\end{equation}
This short exact sequence plays an important role in~\cite{MR2642409},
and using Colmez's results (via a reduction to the case that~$A$ is Artinian), one can prove that there is a perfect
duality pairing (ultimately coming from the residue pairing in Tate
local duality) between $D_A^{\natural} \boxtimes \Pj^1$ and $(D_A\boxtimes \Pj^1)/(D_A^{\natural}\boxtimes \Pj^1)$.

We write $\mathbf{D}: \Dfp^b(\smG) \to \Dfp^b(\smG)$ for the
contravariant functor defined via
\[\mathbf{D}(\text{--}) = \RHom(\text{--},
\cO[[G]])[4].\] Using the resolution~\eqref{tree}, it is
easy to check that~$\mathbf{D}$ has amplitude~$[0,1]$.
We also write $\mathbf{D}_{\cX}$ for  the antiequivalence of $D_{\coh}(\cX)$
given by Grothendieck--Serre duality,
i.e.\ $\RHom_{\cO_{\cX}}(\text{--},\omega_{\cX})$ where $\omega_{\cX}$ is the dualizing
sheaf of~$\cX$,  placed in degree~$0$.

By the construction of~$L_\infty$, %
we can use~\eqref{eqn: Colmez ses} to construct an extension
\[0 \to L_{\infty} \to \text{?} \to 
  \mathbf{D}_{\cX}L_{\infty} \to 0.\] This extension class is an
element of
\[ \Hom_{\cO[[G]]}\bigr( \mathbf{D}_{\cX}L_{\infty}, L_{\infty}[1]\bigr) = \Ext^1_{\cO[[G]]}\bigr(
  \mathbf{D}_{\cX}L_{\infty},
  L_{\infty}\bigr)[1] ,\] and composition with this extension class
gives us a natural morphism of functors on $\Dfp^b(\smG)$
\[\RHom_{\cO[[G]]}(\text{--}, \mathbf{D}_{\cX} L_{\infty})
\to
\RHom_{\cO[[G]]}\bigl(\text{--}, L_{\infty}\bigr)
[1]
\]
which we expect can be shown to be an isomorphism, so that we have the
expected compatibility with duality %
\[\fA\circ \mathbf{D} \iso (\mathbf{D}_{\cX} \circ \fA)[3].\]

As discussed in Remark~\ref{rem: expect duality to imply maximal CM},
this implies that the various sheaves
$\fA\bigl((\cInd_K^G W_{\lambdau}\otimes \sigma^{\crys,\circ}(\tau))^{\wedge}\bigr)$
are Cohen--Macaulay.
Since this is a property that can be checked by pulling back
over versal rings, 
it also follows more directly from the construction of the functor~$\fA$,
its relationship to Pa\v{s}k\={u}nas's results discussed in Section~\ref{subsubsec:vytas},
and the results of~\cite{paskunasBM}. %

\subsection{The Banach case for \texorpdfstring{$\GL_2(\Qp)$}{GL2(Qp)} --- III.  Examples}
\label{subsec:examples}
  We now describe some explicit values of the functor~$\fA$,
  beginning with its values on compact inductions of Serre weights and
  on (absolutely) irreducible representations, and then turning to some cases
  that are related to tamely potentially Barsotti--Tate Galois representations.
  These latter examples also illustrate the compatibility between the Banach and the analytic case posited in
  Section~\ref{EH subsection: Banachcompatibility}. 

  \subsubsection{Compact inductions of Serre weights}
  \label{subsubsec:fA on Serre weights}%
We consider a Serre weight~$\sigma_{a,b}$.
Assume
  firstly that $b\ne p-1$, i.e.\ that $\sigma_{a,b}$ is generic. Then
  $\fA(\cInd_{K}^G\sigma_{a,b})$ is a line bundle on
  $\cX(\sigma_{a,b})$. This is to be expected, for example by
  comparison to the corresponding Taylor--Wiles patched modules, which
  are free over their support as a consequence of Fontaine--Laffaille
  theory, but it is less obvious precisely which line bundle on
  $\cX(\sigma_{a,b})$ we obtain. To describe this, we note that by
  Proposition~\ref{prop: explicit description of FL stack} (and the
  triviality of the Picard group of $\A^2\setminus\{0\}$) we have
  $\Pic(\cX(\sigma_{a,b}))=\Z$, and it can be checked that
  $\fA(\cInd_{K}^G\sigma_{a,b})$ is the  line bundle corresponding to
  $-1\in\Z$. %
  (This is in fact forced by the compatibility with duality; the
  discussion of Remark~\ref{rem: duality explicated on crystalline
    stack} shows that   $\fA(\cInd_{K}^G\sigma_{a,b})$ must be the
  unique self-dual line bundle on $\cX(\sigma_{a,b})$. See also Section~\ref{subsec: Banach Qp2 stuff} for a related
  computation for~$\GL_2(\Q_{p^2})$.)

  If $b=p-1$ then  $\fA(\cInd_{K}^G\sigma_{a,p-1})$ is a line bundle
  on the union %
  \numequation\label{eqn: support of sym
    p-1}\cX(\sigma_{a,0})\cup\cX(\sigma_{a,p-1})^+\cup\cX(\sigma_{a,p-1})^-.\end{equation}
We anticipate that this line bundle can be described explicitly as
follows: By Proposition~\ref{prop: explicit description weight $p-1$
  stacks}, the stacks $\cX(\sigma_{a,p-1})^{\pm }$ also have Picard
group isomorphic to~$\Z$, and we expect that it can be checked that
the restriction of
$\fA(\cInd_{K}^G\sigma_{a,p-1})$ %
to each irreducible component in~\eqref{eqn: support of
  sym p-1} corresponds to $-1\in\Z$. %
\subsubsection{Absolutely irreducible representations}  
The results of Section~\ref{subsubsec:fA on Serre weights}
allow us to compute the values of~$\fA$ on absolutely irreducible
representations. As discussed in Section~\ref{subsubsec:vytas}, these
are subquotients of the representations \[(\cInd_{K}^G\sigma)/(T-\lambda)\]for ~$\lambda\in
k$; and these representations are in fact absolutely irreducible
unless~$\sigma=\sigma_{a,0}$ or~$\sigma=\sigma_{a,p-1}$ for some~$a$
and~$\lambda=\pm 1$. Supposing that we are not in this case, we can
deduce from the injectivity of $T-\lambda$ on the line bundle
$\fA(\cInd_{K}^G\sigma)$ that $\fA((\cInd_{K}^G\sigma)/(T-\lambda))$
is a rank 1 sheaf supported on the closed locus
~$\pi_\sigma^{-1}(\lambda)\cap\cX(\sigma)$. (In the case~$\lambda\ne
0$ we describe this rank 1 sheaf more precisely in Section~\ref{subsubsec: type 2 block} below.) 

Unsurprisingly, the trivial and Steinberg representations exhibit more
  complicated behaviour, which mirrors the corresponding picture for
  $l$-adic representations in characteristic zero (see Remark~\ref{rem: how Morita example relates to Iwahori Hecke}). %
There are short exact sequences of $G$-representations
\[0\to\cInd_K^G\sigma_{0,0}\to\cInd_K^G\sigma_{0,p-1}\to
  \St\oplus(\nr_{-1}\circ \det)\otimes\St\to
  0\] \[0\to\cInd_K^G\sigma_{0,p-1}\to\cInd_K^G\sigma_{0,0}\to
  1\oplus\nr_{-1}\circ \det \to 0\]
which in combination with the descriptions of the
$\gA(\cInd_K^G\sigma_{a,b})$ above
can be used to show that $\fA(\St)$ is a line bundle on
$\cX(\Sym^{p-1})^+$, while $\fA(1)$ is a shift to cohomological degree~$-1$ of 
another 
line bundle on $\cX(\Sym^{p-1})^+$. %
\subsubsection{Type (2) blocks}\label{subsubsec: type 2 block}
We consider the case of a block of type~(2) in more detail.  In this case  
$\fB = \{\pi_1,\pi_2\}$, with 
$\pi_1 =\Ind_B^G(\chi \otimes \chi^{-1})$ and $\pi_2= \Ind_B^G(\omega \chi^{-1} \otimes \omega^{-1}
\chi)\}$ for some suitably generic characters $\chi: \bQ_p^\times \to k^\times$.
In this context
there are non-split extensions $0 \to \pi_2 \to \kappa  \to \pi_1 \to 0$
and $0 \to \pi_1 \to \kappa' \to \pi_2 \to 0$ (each unique up to isomorphism),
and  we will describe the various  sheaves $\fA(\pi_1),$  $\fA(\pi_2)$,
$\fA(\kappa)$, and $\fA(\kappa')$.

To begin with, we note that we may write
$\pi_1 = (\cInd_K^G  \sigma)/(T-\lambda)$
and 
$\pi_2 = (\cInd_K^G  \sigmacomp)/(T-\lambda^{-1})$
for some very generic $\sigma$ and its companion $\sigmacomp$, 
and some appropriately chosen $\lambda$.
All the sheaves under consideration will then have support contained
in~$\pi_{\sigma}^{-1}(\lambda).$

The locus $\pi_{\sigma}^{-1}(\lambda)$ is a closed substack
of $\cU(\sigmasigmacomp)_{\red}$.  Following the notation of~\eqref{eqn:
explicit description of U sigmasigmacomp}, we have
$$ \cU(\sigmasigmacomp)_{\red} := [\Spec \left(k[t^{\pm 1},x,y]/(xy)\right) / \Gm],$$
with the $\Gm$-action being given by $u\cdot (x,y) = (u^2 x, u^{-2} y).$
The map $\pi_{\sigma}: \cU(\sigmasigmacomp)_{\red} \to \A^1$ is then
simply projection to the  $t$-coordinate,
and so
$$ \pi_{\sigma}^{-1}(\lambda) = [\left(\Spec k[x,y]/(xy)\right) / \Gm]$$
(the locus where $t= \lambda$).

We may give $k[x,y]/(xy)$ its grading induced by the $\Gm$-action
(so $x$  has degree $2$ and $y$ has  degree $-2$),  and may then describe
coherent sheaves on $\pi_{\sigma}^{-1}(\lambda)$ as graded $k[x,y]/(xy)$-modules.
We then have the following descriptions of the sheaves under consideration:
$$\fA(\pi_1) = \bigl(k[x,y]/(y)]\bigr)(-1);$$
$$\fA(\pi_2) = \bigl(k[x,y]/(x)]\bigr)(1);$$
$$\fA(\kappa) = \bigl(k[x,y]/(xy))(-1);$$
$$\fA(\kappa') = \bigl(k[x,y]/(xy))(1).$$
(Here the twists indicate shifts in grading, as usual.)
Note the evident short exact sequences
$$0 \to \bigl(k[x,y]/(x) \bigr)(1)
\buildrel y \cdot \over \longrightarrow\bigl( k[x,y]/(xy) \bigr)(-1)
\longrightarrow \bigl( k[x,y]/(y) \bigr)(-1) \to 0$$
and
$$0 \to \bigl(k[x,y]/(y) \bigr)(-1)
\buildrel x \cdot \over \longrightarrow\bigl( k[x,y]/(xy) \bigr)(1)
\longrightarrow \bigl( k[x,y]/(x) \bigr)(1) \to 0;$$
these can be interpreted as arising from applying $\fA$ to the short exact  sequences
that defined~$\kappa$ and~$\kappa'$.

\subsubsection{Tame principal series types}
\label{subsubsec: tame PS}
In this example and the next, we treat some $p$-adically completed
compactly induced representations that are related to tamely
potentially Barsotti--Tate representations.  We content ourselves with
merely sketching the various
computations required to verify all our claims.  

We begin with the tame principal series case,    
and to this end, let
 $\{\sigma,\sigmacomp\}$ be very generic, and let~$\tau$ be a tame principal
  series type such that~$\sigma$ is a Jordan--H\"older factor
  of ~$\sigmabar^{\crys,\circ}(\tau)$. Note that this implies that $\sigmacomp$ is \emph{not} a Jordan--H\"older factor
  of~$\sigmabar^{\crys,\circ}(\tau)$; indeed, if $\sigma  = \sigma_{a,b}$,
  then the Jordan--H\"older factors of $\sigmabar^{\crys,\circ}(\tau)$ are $\sigma$ itself,
  together with $\sigma_{a + b, p-1-b}$. (Explicitly,
  the type $\sigmabar^{\crys,\circ}(\tau)$ is inflated from a principal series
  representation of $\GL_2(\Fp)$, as in Section~\ref{subsubsec:rep theory interlude}.) %

We may then consider the stack $\cX^{\crys,(1,0),\tau}$ of potentially
Barsotti--Tate representations of type~$\tau$.  We can describe this
closed substack of $\cX$ quite explicitly, and we will proceed to do so.
The first point to note is that the mod $\varpi$ fibre $\overline{\cX^{\crys,(1,0),\tau}}$
is equal to
$\cX(\sigma_{a,b})\cup \cX(\sigma_{a+b,p-1-b}).$ 
(This is a precise
form of the geometric Breuil--M\'ezard conjecture discussed
in Section~\ref{subsubsec: geometric BM Banach} for the particular type~$\tau$.)
This union admits the explicit description given in
Remark~\ref{rem: justification of existence of morphism}, 
which embeds it as an open substack of
$\Spec k[t_1,t_2]/(t_1t_2) \times [\Pone/\Gm].$
This latter stack may be identified with the underlying reduced substack
of the $p$-adic formal algebraic stack
    $(\Spf (\cO[t_1,t_2]/(t_1t_2 - p))^{\wedge}\times [\widehat{\Pone}/\Gmhat],$
and $\cX^{\crys,(1,0),\tau}$ is isomorphic to the open substack of this
formal algebraic stack having
$\cX(\sigma_{a,b})\cup \cX(\sigma_{a+b,p-1-b})$ 
as its underlying reduced substack.

We now turn to describing the sheaf
$\fA\bigl((\cInd_{K}^G\sigma^{\crys,\circ}(\tau))^\wedge\bigr)$.
Of course, this depends on the particular choice of lattice
$\sigma^{\crys,\circ}(\tau)$.  Among all such, there are two particular lattices
$\Lambda_1$ and $\Lambda_2$, determined uniquely (up to scaling) by the condition
that the cosocle of $\Lambda_1$ (resp.\ $\Lambda_2$) is equal to $\sigma$
(resp.\ $\sigma_{a+b,p-1-b}$).
We may and do normalize these lattices so that $p \Lambda_1 \subset \Lambda_2 \subset \Lambda_1$.
Then (up to scaling) any lattice $\Lambda$ is of the form $\varpi^a\Lambda_1 + \Lambda_2,$ %
with $0 \leq a \leq e$ (where
$e$ denotes the absolute ramification degree of our coefficient
field~$L$). (This description of lattices in tame principal series types
for~$\GL_2$ is due to Breuil~\cite[\S2]{MR3274546}, and is a special
case of \cite[Prop.\ 4.1.4]{MR3323575}.)

If $\chi:\F_p^{\times}\times \F_p^{\times} \to \Z_p^{\times} \subseteq \cO^{\times}$
denotes the tame lift of the highest weight of~$\sigma$, which we inflate
to a character of the (usual upper triangular) Iwahori~$I$, then $\Lambda_1 = \Ind_I^K \chi,$ 
while if $\chi^s$ denotes the conjugate of 
$\chi$ by the non-trivial Weyl group element~$s$,
then $\Lambda_2 = \Ind_I^K \chi^s$.
Conjugation by the element $\big(\begin{smallmatrix} 0 & 1 \\ p & 0 \end{smallmatrix}\big)$
in the normalizer (in $G$) of $I$ interchanges $\chi$ and $\chi^s$,
and so there is an isomorphism
$$\cInd_K^G \Lambda_1 = \cInd_K^G \chi \cong \cInd_I^G \chi^s = \cInd_K^G \Lambda_2.$$
Following this isomorphism with the inclusion 
$\cInd_K^G \Lambda_2 \hookrightarrow \cInd_K^G \Lambda_1$
induced by the inclusion of $\Lambda_2$ in $\Lambda_1$
gives rise to an endomorphism $U_p$ of $\cInd_K^G \Lambda_1$. 
Interchanging the roles of $\Lambda_1$ and $\Lambda_2$ gives
rise to another endomorphism $U_p'$ of $\cInd_K^G \Lambda_1$,
with the property that $U_pU_p'= p;$
and indeed
$$\End_G(\cInd_K^G\Lambda_1)= \cO[U_p,U_p']/(U_pU_p'  -  p).$$
Passing to $p$-adic completions, we find that
$$\End_G\bigl((\cInd_K^G\Lambda_1)^{\wedge}\bigr)= \cO[U_p,U_p']^{\wedge}/(U_pU_p'  -  p).$$

If $\Lambda = \varpi^a\Lambda_1 + \Lambda_2,$
then the inclusion of $\Lambda$ in $\Lambda_1$ induces an identification
$\cInd_K^G \Lambda  = (\varpi^a,U_p) \cInd_K^G \Lambda_1,$
and hence a corresponding identification of $p$-adic completions 
$(\cInd_K^G \Lambda)^{\wedge}  = (\varpi^a,U_p) (\cInd_K^G \Lambda_1)^{\wedge}.$
This reduces
the computation of 
$\fA\bigr((\cInd_{K}^G\Lambda)^\wedge\bigl)$
for any lattice $\Lambda$ to the computation of~$\fA\bigr((\cInd_{K}^G\Lambda_1)^\wedge\bigl)$.
 
Now we choose our coefficient field $L$ to be unramified over~$\Q_p$ (indeed the character~$\chi$,
and hence the lattice $\Lambda_1$, are defined over~$\Z_p$), and with this hypothesis
we see that   
$\cX^{\crys,(1,0),\tau}$
is regular.   Thus, since 
(as noted in Section~\ref{subsubsec:duality for })
the sheaf
$\fA\bigr((\cInd_{K}^G\Lambda_i)^\wedge\bigl)$  
is Cohen--Macaulay, 
the Auslander--Buchsbaum formula shows that it is locally free, 
and in fact it must be invertible (e.g.\ for the reason sketched
in Remark~\ref{rem: expect duality to imply maximal CM}).
Since $\Lambda_1^{\vee}  = \Lambda_2,$
and since $(\cInd_K^G\Lambda_1)^\wedge \cong (\cInd_K^G \Lambda_2)^{\wedge},$
we further see that this line bundle is Grothendieck--Serre
self-dual on~$\cX^{\crys,(1,0),\tau}$.

Now one can check that
$\Pic([\widehat{\Pone}/\Gmhat])
= \Z \times\Z$, 
with the first copy of $\Z$ coming from $\Pic(\widehat{\Pone}) = \Z$,
and the second copy coming from twisting by characters of~$\Gmhat$.
We let $\cO(m)_{\widehat{\Pone}}(n)$ denote the invertible sheaf corresponding to
a pair $(m,n) \in \Z\times \Z$.  

Pulling back along the composite
\begin{multline*}
f:
\cX^{\crys,(1,0),\tau} \buildrel \text{open}  \over \hookrightarrow
    (\Spf (\cO[t_1,t_2]/(t_1t_2 - p))^{\wedge}\times [\widehat{\Pone}/\Gmhat]
\buildrel  \text{proj.} \over \longrightarrow
 [\widehat{\Pone}/\Gmhat]
\end{multline*}
induces an isomorphism 
$$
\Pic([\widehat{\Pone}/\Gmhat])
\iso
\Pic(\cX^{\crys,(1,0),\tau}).$$
Furthermore, the dualizing sheaf of $[\widehat{\Pone}/\Gmhat]$,
namely $\cO(-2)_{\widehat{\Pone}},$
pulls back to the dualizing sheaf of $\cX^{\crys,(1,0),\tau}.$
Thus (since it is self-dual), we find that
$$
\fA\bigr((\cInd_{K}^G\Lambda_1)^\wedge\bigl)  
\cong
f^*\cO(-1)_{\widehat{\Pone}}.$$
Since $\fA$ is fully faithful, we also obtain an isomorphism
$$\End_G\bigl((\cInd_K^G\Lambda_1)^{\wedge}\bigr)
\iso \Gamma(\cX^{\crys,(1,0),\tau}, \cO_{\cX^{\crys,(1,0),\tau}}) = \cO[t_1,t_2]^{\wedge}/(t_1 t_2 - p),$$
which is given (after choosing the $t_i$ suitably) 
by mapping 
$U_p$ to $t_1$ and $U_p'$ to $t_2$.
For any lattice $\Lambda = \varpi^a \Lambda_1 + \Lambda_2$,
we then find that
$$\fA\bigr((\cInd_{K}^G\Lambda)^\wedge\bigl)  
= (\varpi^a,t_1) f^*\cO_{\widehat{\Pone}}(-1).$$

In order to compare these values of $\fA$ 
with those we computed above for the compact inductions of Serre
weights, it helps to note first that
pullback to the special fibre induces an 
isomorphism
$$\Pic([\widehat{\Pone}/\Gmhat])\iso\Pic([\Pone/\Gm]).$$
Then,
if we let $x_0$ and $x_1$ denote the homogeneous coordinates on $\Pone$
(so that the inhomogeneous coordinate
$x$ of Remark~\ref{rem: justification of existence of morphism} is equal to $x_1/x_0$),
and let $u$ denote the coordinate on $\Gm$,
each of
$[\A^1_{x_1/x_0}/(\Gm)_u]$
and
$[\A^1_{x_0/x_1}/(\Gm)_{u^{-1}}]$
embeds as an open substack of $[\Pone/\Gm]$,
and restriction to these substacks  induces an embedding
$$\Pic( [\Pone/\Gm] )
\hookrightarrow 
\Pic([\A^1_{x_1/x_0}/(\Gm)_u])  \times
\Pic([\A^1_{x_0/x_1}/(\Gm)_{u^{-1}}])  =  \Z \times  \Z,$$
given by
$$(m,n)  \mapsto (m+n,  m - n).$$
From this, we compute the restrictions
$$(f^*\cO(-1)_{\widehat{\Pone}}\bigr)_{|\cX(\sigma_{a,b})} =
\cO_{\cX(\sigma_{a,b})}(-1)$$
and
$$(f^*\cO(-1)_{\widehat{\Pone}}\bigr)_{|\cX(\sigma_{a+b,p-1-b})} =
\cO_{\cX(\sigma_{a+b,p-1-b})}(-1).$$

Thus,
if we apply $\fA$ to the short exact sequence
$$0 \to (\cInd_K^G \Lambda_1)^{\wedge}  \buildrel U_p \over \longrightarrow
(\cInd_K^G \Lambda_1)^{\wedge}
\to \cInd_K^G \sigma_{a,b} \to 0,$$
respectively\
$$0 \to (\cInd_K^G \Lambda_1)^{\wedge}  \buildrel U_p' \over \longrightarrow
(\cInd_K^G \Lambda_1)^{\wedge}
\to \cInd_K^G \sigma_{a+b,p-1-b} \to 0,$$
we obtain the short exact sequence of sheaves
$$0 \to f^*\cO_{\widehat{\Pone}}(-1) \buildrel t_1 \over \longrightarrow f^*\cO_{\widehat{\Pone}}(-1)
\to \cO_{\cX(\sigma_{a,b})}(-1) \to 0,$$
respectively\
$$0 \to f^*\cO_{\widehat{\Pone}}(-1) \buildrel t_2 \over \longrightarrow f^*\cO_{\widehat{\Pone}}(-1)
\to \cO_{\cX(\sigma_{a+b,p-1-b})}(-1) \to 0.$$

\subsubsection{Tame cuspidal types}
\label{subsubsec: tame cuspidal}
  One can make a similar analysis to that of Section~\ref{subsubsec: tame PS}
  in the case of a tame cuspidal type.
  Namely, there is a tame cuspidal type $\tau$ uniquely determined by the condition 
  that the Jordan--H\"older factors of $\sigma^{\crys,\circ}(\tau)$  
  are $\sigma$ and~$\sigmacomp$. 
  The potentially Barsotti--Tate substack
$\cX^{\crys,(1,0),\tau}$
  is then a $p$-adic formal lift of $\cX(\sigmasigmacomp) := \cX(\sigma)\cup \cX(\sigmacomp).$  
  As noted in Remark~\ref{rem:description of cX(sigmasigmacomp) in the generic case},
$\cX(\sigmasigmacomp)$ embeds as an open substack of
$\Pone \times \bigl[\bigl(\Spec k[x,y]/(xy)\bigr)/\Gm\bigr],$
and a deformation theory argument shows that
$\cX^{\crys,(1,0),\tau}$
is isomorphic to the open formal algebraic substack of 
the $p$-adic formal algebraic stack 
$\widehat{\Pone}\times\bigl[ \bigl(\Spf (\cO[x,y]/(xy - p)\bigr)^{\wedge}/\Gmhat\bigr]$ 
lifting~$\cX(\sigmasigmacomp)$.

Just as in the principal series case,
there are distinguished lattices $\Lambda_1$  and $\Lambda_2$  in $\sigma^{\crys}(\tau)$,
uniquely determined up to scaling,
whose cosocles are
$\sigma$ and $\sigmacomp$ respectively.
Fixing a choice of $\Lambda_1$, we can then scale $\Lambda_2$ so that
$p \Lambda_1  \subset \Lambda_2  \subset \Lambda_1$, and any lattice is (up to  scaling)
of the form $\Lambda = \varpi^a\Lambda_1 + \Lambda_2$,
with $0 \leq a \leq e$.
Thus, to describe the various sheaves
$\fA\bigr((\cInd_{K}^G\Lambda)^\wedge\bigl),$  
it suffices to describe the sheaves
$\fA\bigr((\cInd_{K}^G\Lambda_1)^\wedge\bigl)$  
and
$\fA\bigr((\cInd_{K}^G\Lambda_2)^\wedge\bigl)$,
as well as the morphism 
\numequation
\label{eqn:lattice inclusion on sheaves}
\fA\bigr((\cInd_{K}^G\Lambda_2)^\wedge\bigl)  
\rightarrow
\fA\bigr((\cInd_{K}^G\Lambda_1)^\wedge\bigl)  
\end{equation}
induced by the inclusion $\Lambda_2 \subset  \Lambda_1$.

As in the principal series we choose our coefficient field $L$ to be unramified over~$\Q_p$ (tame types
may always be defined over unramified extensions), and with this hypothesis
we see that   
$\cX^{\crys,(1,0),\tau}$
is again regular, so that   %
the sheaves
$\fA\bigr((\cInd_{K}^G\Lambda_i)^\wedge\bigl)$  
are again invertible. %
Further, since $\Lambda_1^{\vee}=  \Lambda_2$,
we find that 
$\fA\bigr((\cInd_{K}^G\Lambda_1)^\wedge\bigl)$  
and
$\fA\bigr((\cInd_{K}^G\Lambda_2)^\wedge\bigl)$  
must be interchanged under Grothendieck--Serre duality.

Now pullback induces an isomorphism 
$$\Pic(
\widehat{\Pone}\times\bigl[ \bigl(\Spf (\cO[x,y]/(xy - p)\bigr)^{\wedge}/\Gmhat\bigr] )
\iso
\Pic(\cX^{\crys,(1,0),\tau}),$$
and so any invertible sheaf on
$\cX^{\crys,(1,0),\tau}$
is necessarily of the form $\cL(n),$ where $\cL$ is pulled back from $\widehat{\Pone},$
and $n$ is a an integer indicating the weight of the $\Gmhat$-action.
Furthermore, the dualizing sheaf 
on $\cX^{\crys,(1,0),\tau}$
is obtained as the pull back of the dualizing sheaf
$\cO(-2)_{\widehat{\Pone}}$  on
$\widehat{\Pone}\times\bigl[ \bigl(\Spf (\cO[x,y]/(xy - p)\bigr)^{\wedge}/\Gmhat\bigr].$ 

The short exact sequences
$$0 \to (\cInd_K^G \Lambda_2)^{\wedge}  \rightarrow
(\cInd_K^G \Lambda_1)^{\wedge}
\to \cInd_K^G \sigma \to 0$$
and
$$0 \to (\cInd_K^G \Lambda_1)^{\wedge}  \buildrel p \over \longrightarrow
(\cInd_K^G \Lambda_2)^{\wedge}
\to \cInd_K^G \sigmacomp \to 0$$
induce corresponding short exact sequences of sheaves after applying~$\fA$,
and from the known values of $\fA$
on $\cInd_K^G \sigma$ and $\cInd_K^G \sigmacomp$,
we find that
$$\fA\bigr((\cInd_{K}^G\Lambda_1)^\wedge\bigl)
\cong \cL(-1),$$
that  %
$$\fA\bigr((\cInd_{K}^G\Lambda_2)^\wedge\bigl)
\cong \cL(1),$$
and that the morphism~\eqref{eqn:lattice inclusion on sheaves}
is induced by multiplication by~$y$.  
This analysis doesn't determine the invertible sheaf~$\cL$
(since the restriction of any $\cL$ pulled back from $\Pone$
to either $\cX(\sigma)$ or $\cX(\sigmacomp)$ is necessarily
trivial). 
However, the mutual duality  of $\cL(\pm 1)$  then
shows  that $\cL$ must be the pullback of~$\cO(-1)_{\widehat{\Pone}}.$

\begin{rem}
  \label{rem: difference between PS and cuspidal}One difference
  between the cuspidal and principal series cases is that in the
  cuspidal case, we have $\Gamma(\cX^{\crys,(1,0),\tau},\cO_{\cX^{\crys,(1,0),\tau}}) = \cO$
  (just as we have $\Gamma\bigl(\cX(\sigmasigmacomp),\cO_{\cX(\sigmasigmacomp)}\bigr)  = k$).
  This reflects (via consideration of the functor~$\fA$) the fact 
  that the endomorphism ring of the compactly 
  induced representation $\bigl(\cInd_K^G \sigma^{\crys,\circ}(\tau)\bigr)^{\wedge}$
  (for any choice of lattice~$\sigma^{\crys,\circ}(\tau)$)
  is trivial in this case.
\end{rem}

\subsubsection{Illustrating the compatibility with the analytic case}%
Allowing, for a moment, $\tau$ to denote either a tame principal series
or tame cuspidal type,
then, as indicated in (\ref{eqn: semistable period map}) in the semi-stable case,
there is a morphism
$$\pi: (\cX^{\crys,(1,0),\tau})^{\rm rig}_{\eta}
\to ({\rm Fil}_{(1,0)}{\rm WD}_{\tau})^{\rm an},$$
whose target is the (rigid analytic) moduli stack 
of two-dimensional Weil--Deligne representations having inertial type~$\tau$,
of fixed inverse cyclotomic determinant, %
equipped with a filtration by a one-dimensional subspace. 

Suppose now that $\tau$ is a principal series type.
Then 
\begin{multline*}
(\cX^{\crys,(1,0),\tau})^{\rm rig}_{\eta}
\cong
\\
\bigl[ \bigl(
\{(t,x) \in (\Gm)_t\times \Pone_x\, | \,  1\geq |t| \geq |p| \} \setminus
( \{(x,t) \, | \, 1 > |t|, |x|\}
\cup \{(x,t)\, | \, |t| > |p| , |x| > 1\}) 
\bigr) / (\Gmhat)^{\rm rig}_{\eta} \bigr],
\end{multline*}
where $(\Gmhat)^{\rm rig}_{\eta} = \{ u \, | \, |u| = 1\} \subset \Gm,$
acting via $u \cdot (t,x) = (t, u^2 x),$
while
$$({\rm Fil}_{(1,0)}{\rm WD}_{\tau})^{\rm an}
\cong 
[ (\Gm)_t \times \Pone_x / (\Gm)_u],$$
where\footnote{Here we write just $\mathbf{G}_m$ and $\mathbf{P}^1$ for the multiplicative group and the projective line viewed as rigid analytic spaces.} again the action is given by 
$u \cdot (t,x) = (t, u^2 x).$
The map $\pi$ is then the obvious one.

In terms of this description,
the ``forget the filtration'' map $\rm pr_{WD}$  
is the morphism
$$
[ (\Gm)_t \times \Pone_x / (\Gm)_u]
\to [(\Gm)_t/(\Gm)_u]$$
(the right hand stack being formed with respect to the trivial action
of $(\Gm)_u$ on~$(\Gm)_t$).
The twist $((-\rho')_\sigma)$ which appears in
\eqref{EHanconj smoothcomparison}
is, in this case, the twist by the pullback of the sheaf $\cO(-1)_{\mathbf{P}^1}$ along the the ``forget the Weil--Deligne representation'' map $${\rm pr_{Fil}}:
[ (\Gm)_t \times \Pone_x / (\Gm)_u]
\to  \Pone_x/\SL_2,$$
that is induced by projection to the 
second factor together with the inclusion of $(\Gm)_u$ as the diagonal torus in~$\SL_2$. 

Now the smooth categorical Langlands correspondence associates
to $\cInd_K^G \sigma^{\crys}(\tau)$
the trivial invertible sheaf, with trivial $(\Gm)_u$-equivariant structure,
on $[(\Gm)_t/(\Gm)_u]$.
The expected compatibility (\ref{EHeqnlocalgcompatibility}) in Section \ref{EH subsection: Banachcompatibility} 
(taking into account also Conjecture~\ref{EH:conjanalytic}~(1))
then predicts that
$$\fA\bigl((\cInd_{K}^G\sigma^{\crys,\circ}(\tau))^\wedge\bigr)^{\rm rig}_{\eta}
\cong {\rm pr_{\Fil}^*}\bigl(\cO(-1)_{\Pone}\bigr).$$
The line bundle $\cO(-1)_{\Pone}$ on $\Pone_x$ pulls back to the trivial
line bundle on the open subset
$$\bigl \{(t,x) \, \mid \,  1\geq |t| \geq |p| \bigr \} \setminus
\Bigl( \bigl\{ (x,t) \, \mid \, 1 > |t| , |x|\bigr \}
\cup \bigl \{(x,t)\, \mid \, |t| > |p| , |x| > 1 \bigr\}\Bigr) 
$$
of $(\Gm)_t \times\Pone_x$
under the projection to $\Pone_x,$
but it has non-trivial $(\Gm)_u$-action: in fact, $(\Gm)_u$ acts on it with weight~$-1$.
Thus (\ref{EHeqnlocalgcompatibility}) predicts an isomorphism
$$\fA\bigl((\cInd_{K}^G\sigma^{\crys,\circ}(\tau))^\wedge\bigr)^{\rm rig}_{\eta}
\cong \cO_{(\cX^{\crys,(1,0),\tau})_{\eta}^{\rm rig}}(-1),$$
in accordance with the computation of Section~\ref{subsubsec: tame PS}.

Now suppose that $\tau$ is cuspidal.  In this case
\begin{multline*}
(\cX^{\crys,(1,0),\tau})^{\rm rig}_{\eta}
\cong
\\
\Bigl[\Bigl(
\bigl\{(t,x) \in \Pone \times \Gm  \, \mid \,  1  \geq |x| \geq |p|\bigr\} \setminus
\Bigl( \bigl \{(t,x) \, \mid \, 1 > |t| , |x| \bigr\}  \cup  \bigl\{ (t,x) \, \mid \, |t| >  1, |x| > |p|
\bigr \} ) \Bigr) /
(\Gmhat)^{\rm rig}_\eta
\Bigr] 
,
\end{multline*}
where
$(\Gmhat)^{\rm rig}_\eta$ acts via
$u\cdot (t,x) = (t,u^2 x),$
while
$$({\rm Fil}_{(1,0)}{\rm WD}_{\tau})^{\rm an}
\cong 
[\Pone_t/\boldsymbol{\mu}_2],$$
where $\boldsymbol{\mu}_2$ acts trivially.
The map $\pi$ is then the obvious one, given by forgetting~$x$.

The map $\rm pr_{WD}$ is the projection
$$
[\Pone_t/\boldsymbol{\mu}_2] \to
[*/\boldsymbol{\mu}_2],$$
while  we again have the map $${\rm  pr_{Fil}}:
[\Pone_t/\boldsymbol{\mu}_2] \to
[\Pone_t/\SL_2]
$$
that forgets the Weil--Deligne representations and which is
induced by the central embedding of $\boldsymbol{\mu}_2$ into $\SL_2$.
The expected 
compatibility (\ref{EHeqnlocalgcompatibility})
then predicts an isomorphism
$$\fA\bigl((\cInd_{K}^G\sigma^{\crys,\circ}(\tau))^\wedge\bigr)^{\rm rig}_{\eta}
\cong (\pi^*{\rm \pr_{Fil}^*}\cO(-1)_{\Pone})(-1)$$
(the outer twist by  $-1$
indicating the equivariant action of~$(\Gmhat)^{\rm rig}_{\eta}$),
which is in accordance with the computations of Section~\ref{subsubsec: tame cuspidal}.
(Note that the twists by $1$ and by $-1$ give non-isomorphic invertible sheaves
on~$\cX^{\crys,(1,0),\tau}$,
but that these line bundles become isomorphic upon being pulled back
to~$(\cX^{\crys,(1,0),\tau})^{\rm rig}_{\eta}$,
since the invertible sheaf
$\cO_{(\cX^{\crys,(1,0),\tau})_{\eta}^{\rm rig}}(2)$
is trivial, admitting as it does the coordinate $x$ as a nowhere zero global section.)

\subsection{Semiorthogonal decompositions}
\label{subsec:semi-decomp}We now explain an explicit semiorthogonal
decomposition of $\IndCoh(\cX)$ (where $\cX$ is as in Section~\ref{sec:banach-case-gl_2qp}), %
which allows us to describe the
essential image of the functor~$\gA$. We also discuss the essential
image of the analogous expected functor for representations of
$D^\times$, where~$D$ is the nonsplit quaternion algebra with centre
$\Qp$, as well as the case of $\GL_2(\Q_{p^2})$. 

We begin by noting that we have a morphism $f: \boldsymbol{\mu}_2\times\cX\to\cI_{\cX}$ (the inertia stack of~$\cX$) given
by the action of $\zeta\in\boldsymbol{\mu}_2$ on a $(\varphi,\Gamma)$-module by
multiplication by~$\zeta$ (note that since~$\zeta^2=1$, this acts
trivially on $\wedge^2D$). This morphism is furthermore {\em central}
(as one sees from its construction) with respect to the groupoid structure
on $\cI_{\cX}$ (which is a groupoid over~$\cX$).

The morphism $f$ induces a canonical action of~$\boldsymbol{\mu}_2$
on~$\Coh(\cX)$, and thus on~$\Ind\Coh(\cX)$, as follows:
We may find a formal scheme $U$ with a smooth surjective morphism $U \to \cX$
(in other words, a smooth chart for~$\cX$),
and if $R := U\times_{\cX} U,$ then we may write $\cX = [U/R].$
The morphism $f$ above then induces a morphism $U\times \boldsymbol{\mu}_2 \to R$
of groupoids over~$U$, where $\boldsymbol{\mu}_2$ acts trivially on~$U$.  A coherent sheaf
on $\cX$ may be regarded as a coherent sheaf on $U$ with an $R$-equivariant
structure (encoding the descent data back down to $\cX$); in particular
any such sheaf is equipped with a $\boldsymbol{\mu}_2$-equivariant structure,
which (by the centrality of~$f$) commutes with the $R$-equivariant structure.
Thus coherent sheaves on~$\cX$  are indeed equipped with  a canonical  $\boldsymbol{\mu}_2$-equivariant structure.

We let~$\IndCoh(\cX)^{\pm}$ be the subcategory of~$\IndCoh(\cX)$, on which
~$-1\in\boldsymbol{\mu}_2$ acts as~$\pm$. We sometimes refer to $\IndCoh(\cX)^+$ as
the \emph{even} subcategory, and $\IndCoh(\cX)^-$ as the \emph{odd}
subcategory. 

As a first step towards a semiorthogonal decomposition of
$\IndCoh(\cX)$, it is easy to see that we have an orthogonal
decomposition into the two pieces $\IndCoh(\cX)^+$  and
$\IndCoh(\cX)^-$. It follows from the construction of~$L_\infty$ that
$L_\infty$ is odd, so the essential image of~$\fA$ is contained in ~$\IndCoh(\cX)^{-}$.

\subsubsection{Motivation}We refer to Appendix~\ref{sub:
  semiorthogonal generalities} for some generalities on semiorthogonal
decompositions. Our semiorthogonal decomposition is motivated by the
possibility of a ``$p$-adic Fargues--Scholze conjecture'' as in
Remark~\ref{rem:why only fully faithful functors}, and by the
structure of the stack $\Bun_G$ of $G$-bundles on the
Fargues--Fontaine curve in the case $G=\GL_2(F)$. %
This stack admits
a stratification with strata corresponding to particular isocrystals,
and as in the $\ell$-adic case, it is reasonable to expect that
whatever categories of $p$-adic sheaves on $\Bun_G$ are considered,
there will be a corresponding semiorthogonal decomposition into pieces
which admit a description in terms of the representation theory of the
automorphism groups of isocrystals. We warn the reader that we do not,
however, expect that the semiorthogonal decomposition on
$\Ind\Coh(\cX)$ that we consider in this section will match precisely with the
(anticipated) semiorthogonal decomposition coming from~$\Bun_G$; rather we expect
that the two are related ``by an upper triangular unipotent change of
basis matrix'' which may be quite complicated.\footnote{Xinwen Zhu
has remarked to the authors that the situation is somewhat analogous
to that of a sheaf-theoretic Radon transform as considered in~\cite[\S 4]{MR2480718}.}

\begin{rem}
  \label{rem: ignoring nilpotent singular support}It is likely that the
 formulation of $p$-adic local Langlands as an equivalence of categories
 should involve a nilpotent singular support condition on
 $\Ind\Coh(\cX)$. We ignore this throughout this section; note that for
 the most part we confine ourselves to the generic part of the
 stack~$\cX$, where we expect this condition to be automatic.
\end{rem}

\subsubsection{Considerations from $\Bun_G$}
As already indicated, the structure of the semiorthogonal decomposition
of $\Ind\Coh(\cX)$ that we aim to exhibit is motivated by a parallel structure
on the category of constructible sheaves on $\Bun_G$ (the
case $G=  \GL_2(F)$ being the one of interest to us), 
arising from a natural stratification on $\Bun_G$ that we now recall.

Namely, %
following~\cite[\S I.4]{fargues--scholze} and \cite[\S
2.6]{fargues--scholze-IHES}, the stack $\Bun_{\GL_2}$ admits a
stratification into locally closed strata as follows. The connected
components of~$\Bun_{\GL_2}$ are indexed by elements
$\alpha\in\frac{1}{2}\Z$, and each connected component contains a
unique open
stratum. If~$\alpha\in\Z$, this open stratum is isomorphic to
$[*/\GL_2(F)]$, while if $\alpha\notin\Z$, it is isomorphic to
$[*/D^\times]$, where $D$ is the non-split quaternion algebra with
centre~$F$. (Here and below we informally confuse the $p$-adic points
of group schemes with the corresponding sheaves.) The other strata are indexed by the pairs of integers
$(m,n)$ with $m>n$ and $m+n=2\alpha$, and are identified with the
quotient stacks $[*/G_{m,n}]$, where \[G_{m,n}=
  \begin{pmatrix}
    F^\times &\mathcal{BC}(\cO(m-n))\\0& F^\times
  \end{pmatrix}
\]is the semidirect product of a unipotent group diamond (in fact a
Banach--Colmez space)
$\mathcal{BC}(\cO(m-n))$ by the torus $F^\times\times F^\times$. (The
stack $\Bun_{\PGL_2}$ can be described by similar data, by identifying
those isocrystals which differ by a twist by $\cO(1)$; for example,
the connected components are now
indexed by $\alpha\in\frac{1}{2}\Z/\Z$.)%

In the $\ell$-adic case there is a corresponding (infinite)
semiorthogonal decomposition of the category of sheaves
on~$\Bun_{\GL_2}$, by \cite[Thm.\ I.5.1]{fargues--scholze}, and we
expect that the same will be true in the $p$-adic case. (Indeed, %
the semiorthogonal decomposition of $\IndCoh(\cX)$ that we exhibit %
below was motivated by this expectation, and retrospectively justifies
it.) In the $\ell$-adic case one can furthermore ignore the
Banach--Colmez spaces $\mathcal{BC}(\cO(m-n))$, but it is \emph{a
  priori} unclear that this should be the case for $p$-adic
coefficients. As we explain below, we expect that they affect our
semiorthogonal decomposition if and only if $[F:\Qp]\ge m-n$. (In
particular, for $\GL_2(\Qp)$, we find that they give  rise to a single
constituent of the semiorthogonal decomposition of~$\IndCoh(\cX)^+$,
and do not intervene at all in the decomposition of~$\IndCoh(\cX)^-$ ---
and so in particular do not intervene in the description of the
essential image of~$\fA$.)

As we have already indicated, in the case $F=\Qp$, we construct below
a semiorthogonal decomposition of $\IndCoh(\cX)$, and in particular we
describe the image of~$\fA$ as the left orthogonal of a category of
sheaves pushed forward from the reducible locus in~$\cX$. We
furthermore give a semiorthogonal decomposition of the sheaves on the
residually reducible locus. As indicated above, we anticipate that
this latter decomposition should ultimately be related to the
representation theory of~$\Qptimes$. Furthermore, it should have a
description in terms of a ``spectral Eisenstein functor'' involving
the analogues in the Banach setting of the stacks of Borel
$(\varphi,\Gamma)$-modules as in Remark~\ref{EH:rem:derived}. %
\subsubsection{Restriction to the reducible locus}We return to the
setting $F=\Qp$. The residually
reducible locus~$\cU\subset\cX$ is a union of connected components
$\cU(\sigmasigmacomp)$, where as always for the sake of exposition we
assume that $\sigma$ is very generic. As we saw in~\eqref{eqn: explicit
  description of U sigmasigmacomp}, we have
$\cU(\sigmasigmacomp)=[\Spf \cO[s,t^{\pm  1}, x,y]^\wedge/\Gmhat]$,
where the wedge denotes $(p,s,xy)$-adic completion,
and $x$ and $y$ have weights $2$ and $-2$ respectively.

Write $A := \cO[s,t^{\pm 1}],$
and write $B = A[x,y]$, regarded as a graded $A$-algebra by giving $x$ weight~$2$
and $y$ weight~$-2$.  
Let $I := (p,s)\subset A$ 
and $J := (I,xy) \subset B$ denote the indicated ideals,
and let $\cB_J$ denote the full subcategory of the category of graded $B$-modules
consisting of those modules each element of which is annihilated by a power
of~$J$. 
There is then a canonical equivalence %
$$\Ind\Coh\bigl(\cU(\sigmasigmacomp)\bigr) \iso \Ind\Coh(\cB_J)$$ 
(for the definition of the right hand side, see~\eqref{eqn:functor from IndCoh}
as well as the discussion of Section~\ref{subsubsec:set of generators}),
which respects the $\pm$-decomposition on its source and
target (with respect to the $\boldsymbol{\mu}_2$-action on the source, and
with respect to parity of grading on the target).

It follows straightforwardly from Proposition~\ref{prop:semiorthogonal decomposition}
that each $\pm$-component
$\Ind\Coh(\cB_J^{\pm})$ 
admits a semiorthogonal decomposition
into semiorthogonal pieces
$\cD_0^{\pm}$ and $\cE_0^{\pm} = (\cD_0^{\pm})^{\perp},$
where 
$$\cD_0^+ = \langle B/J\rangle,$$
$$\cE_0^+ = \langle (A/I)[x](-2), (A/I)[x](-4),
\cdots, (A/I)[y](2), (A/I)[y](4),\cdots \rangle,
$$
$$\cD_0^- = \langle B/J(1), B/J(-1) \rangle,$$ 
and
$$\cE_0^- = \langle (A/I)[x](-3), (A/I)[x](-5),
\cdots, (A/I)[y](3), (A/I)[y](5),\cdots \rangle.
$$
By our Morita-theoretic construction of the restriction
to~$\cU(\sigmasigmacomp)$ of the functor~$\fA$
given in Section~\ref{subsubsec:restriction to U},
we see that the image of this restriction 
is precisely equal to $\cD_0^-$.

\subsubsection{The essential image of~$\fA$ in the case of $\PGL_2(\Qp)$}\label{subsubsec: image of GL2Qp}
We
now consider the pushforward along $j:\cU\to\cX$ of our semiorthogonal
decomposition of~$\IndCoh(\cU)$ (ignoring as ever the non-generic
components, and therefore the issue of nilpotent singular support), as
in Appendix~\ref{subsec: restricting semiorthogonal to open substacks}. By~\eqref{numeqn:
  orthogonal of pushforward}, we have a semiorthogonal decomposition
  $\cA_1,\cA_2$ of $\IndCoh(\cX)^{-}$, where $\cA_2=j_*\cE_0^{-}$, and
  $\cA_1=(j^*)^{-1}(\cD_0^{-})$.

  We claim that $\cA_1$ is the essential image of~$\fA$. %
To
  see this, note that for any $x$ in~$\cA_1$, %
  we have the morphism \[x\to
    j_*j^*x\] (from the unit of the
  adjunction between  $j_*$ and $j^*$), whose cofiber is supported on the
  complement of~$\cU$. This cofiber is therefore in the essential image
  of~$\fA$, because this essential image contains the (odd parity) skyscraper
  sheaves at the supersingular points (these are the images of the
  supersingular irreducible representations, via Theorem~\ref{thm:
    Vytas theorem}). It therefore suffices to show that the
  restriction to~$\cU$ of the essential image of~$\fA$ is
  $\cD_0^-$, which, as we already observed, follows from the explicit description of the
  restriction to~$\cU$ in Section~\ref{subsubsec:restriction to U}.

\subsubsection{Representations of~$D^\times$}We now consider the sheaves of even
parity.  We continue to assume that $p\ge 5$ for simplicity.

We begin by describing some features of a conjectural functor
$\fA_{D^\times}$ taking smooth representations of~$D^\times$ with trivial
central character to the even category $\IndCoh(\cX)^{+}$. We anticipate that
such a functor exists, and enjoys similar properties to those of our
functor~$\fA$ for~$\PGL_2(\Qp)$. In particular, we expect it to be
continuous and fully faithful, and we expect it to satisfy a
compatibility with Taylor--Wiles patching. %

Using this latter expectation, we can predict the values of
$\fA_{D^\times}$ on a set of generators. In order to do so, we briefly
recall some of the basic representation theory
of~$D^\times$. To avoid having to discuss rationality issues, we work
with $\Fpbar$-coefficients. Writing~$\Pi$ for a uniformizer of~$D$, and~$\cO_D$ for
the usual maximal order, we see that $\cO_D/\Pi$ is a quadratic
extension of~$\Fp$, which we identify with $\F_{p^2}$. Thus the
irreducible $\Fpbar$-representations of~$\cO_D^\times$ are given by the
characters of~$\F_{p^2}^\times$. Fixing an embedding
$\F_{p^2}\into\Fpbar$ defines a character
$\etabar:\cO_D^\times\to\Fpbartimes$, and every other such character is of
the form $\etabar^i$ for some~$i\in\Z/(p^2-1)\Z$. These
representations have trivial central character if and only if
$(p-1)$ divides~$i$, which we assume from now on; then we can extend them to representations of
$\Pi^{2\Z}\cO_D^\times$ with trivial central character, which we also
denote by~$\etabar^i$.

If $(p+1)\nmid i$ then
$\sigma_i:=\Ind_{\Pi^{2\Z}\cO_D^\times}^{D^\times}\etabar^i$ is an
irreducible 2-dimensional representation of~$D^\times$ (with trivial
central character); we have $\sigma_i\cong\sigma_{-i}$, and
$\sigma_i\not\cong\sigma_j$ if $i\ne \pm j$. %
If $i=0$ or $i=(p^2-1)/2$, then
$\Ind_{\Pi^{2\Z}\cO_D^\times}^{D^\times}\etabar^i$ is the direct sum of
two characters~$\sigma_i^{\pm}$ of~$D^{\times}$, satisfying
$\sigma_i^{\pm}(\Pi)=\pm 1$. These representations give a set of generators of the
category of smooth $D^\times$ representations with trivial central
character. %

We anticipate that the values of~$\fA_{D^\times}$ on these
representations take a particularly simple form. The mod~$p$
Jacquet--Langlands correspondence of T.G.\ and David Geraghty~\cite{MR3449190}
associates to $\etabar^i$ a companion pair of Serre weights
$\sigma,\sigmacomp$ (in the sense of Definition~\ref{defn:
  companion Serre weights}). Indeed, by Remark~\ref{rem: companion
  weights and cuspidal types}, the pair $\sigma,\sigmacomp$ is the set
of Jordan--H\"older constituents of a cuspidal tame type, and this
cuspidal tame type is the one corresponding to (the tame lift to
characteristic zero of the representation of~$I_{\Qp}$ corresponding
via local class field theory to) $\etabar^i\oplus\etabar^{-i}$.

Write $\cX(\sigma_i):=\cX(\sigma)\cup\cX(\sigmacomp)$. (Note that this
does not always agree with $\cX(\sigmasigmacomp)$ as defined in
Definition~\ref{defn: X sigma sigmaco}, because we are omitting
components of the form $\cX(\sigma_{a,p-1})^{\pm}.$) 

 Then we expect
 that (with $j=0$ or~$1$) \numequation\label{eqn: expectation for Dx
   2}\fA_{D^\times}(\sigma_{(p^2-1)j/2}^{\pm})=\cO_{\cX(\sigma_{(p-1)j/2,p-1})^{\pm}},\end{equation}
while if $(p+1)\nmid i$
we expect that \numequation\label{eqn: expectation for Dx 1}\fA_{D^\times}(\sigma_i)=\cO_{\cX(\sigma_i)}.\end{equation} 

Our justification for this expectation is as follows. By multiplicity
one results for mod~$p$ quaternion modular forms, i.e.\ freeness
results for the corresponding patched modules (see Section~\ref{subsec: mult one for D}), and the expected
compatibility of $\fA_{D^\times}$ with Taylor--Wiles patching, we
expect that all of these sheaves are invertible sheaves over their supports (which
are the indicated closed substacks).
(As far as we know, these multiplicity one results are not in
the literature, so we sketch a proof of them below.) Since they are also in $\IndCoh(\cX)^{+}$, it is natural to
guess that they are given by the structure sheaves of their supports.

\begin{rem}
  \label{rem: global justification for trivial twist for D}It seems
  plausible that the compatibility of the functor $\fA_{D^\times}$ with  Conjecture~\ref{conj:cohomology}
  in  Eisenstein situations
  (for a quaternion algebra over~$\Q$ which is ramified at~$p$) would
  imply that these invertible sheaves are indeed just the structure
  sheaves of their support, adding further evidence for our expectation.  
\end{rem}%

\begin{rem}
  \label{rem: extending companion weights to p-2}If we drop the
  assumption that our representations have trivial central character,
  then the above statements extend in the obvious fashion, once we
  decree that a Serre weight $\sigma_{a,p-2}$ is its own companion (so
  the ``companion pair'' is a singleton in this case). 
\end{rem}

We can now repeat the analysis of Section~\ref{subsubsec: image of
  GL2Qp}. We again have a semiorthogonal decomposition
  $\cA_1,\cA_2$ of $\IndCoh(\cX)^{+}$, where $\cA_2=j_*\cE_0^{+}$, and
  $\cA_1=(j^*)^{-1}(\cD_0^{+})$. Using expectations~\eqref{eqn: expectation for Dx
   2} and~\eqref{eqn: expectation for Dx 1}, we see that the essential
 image of $\fA_{D^\times}$ is contained in~$\cA_1$, and that $\cA_1$
 is generated by this essential image, together with the (even parity)
 skyscraper sheaves at the irreducible points. %

 However, in contrast to the case of~$\GL_2(\Qp)$,  the essential
 image of~$\fA_{D^\times}$ does not contain these
 skyscraper sheaves! Instead, we expect that these skyscraper accounts
 are accounted for by the first non-open stratum of~$\Bun_{\PGL_2}$,
 i.e.\ by the representation theory of $\mathcal{BC}(\cO(1))$. We are
 therefore suggesting that  $\mathcal{BC}(\cO(1))$ behaves differently
 from  $\mathcal{BC}(\cO(n))$ for~$n>1$; as some
 justification for this, we note that  $\mathcal{BC}(\cO(n))$ is
 representable by a perfectoid space if~$n=1$, but not if $n>1$.

\subsubsection{The case of $\GL_2(\Q_{p^2})$}\label{subsec:
  semiorthogonal Qp2}%
It is possible to make a similar analysis in the case
of~$\GL_2(\Q_{p^f})$. We can give generators for the essential image
of the hypothetical functor~$\fA$ in a similar fashion to our
predictions for~$D^\times$ above, and we can describe
the irreducible components of the corresponding
stack~$\cX_{\Q_{p^f},2}$ (for generic Serre weights, we do
this in Section~\ref{subsec: Banach Qpf stuff}). One can construct a
semiorthogonal decomposition on the reducible locus by
an analysis similar to the one given above.

The case $\GL_2(\Q_{p^2})$ is very similar to the case of~$D^\times$
considered above. We again find that the essential image of (the hypothetical) $\fA$ is contained
in the even part of the corresponding category of sheaves, and is
contained in the left orthogonal of an explicit category of sheaves
pushed forward from the reducible locus. We again find that the
skyscraper sheaves at irreducible points are contained in this left
orthogonal, but seem not to be in the essential image of~$\fA$. %

We again anticipate that these skyscraper sheaves are accounted for by
the representation theory of a Banach--Colmez space, in this case
$\mathcal{BC}(\cO(2))$; again, this Banach--Colmez space is
represented by a perfectoid space. (For a general~$F$,
$\mathcal{BC}(\cO(n))$ is represented by a perfectoid space if and
only if $0<n\le [F:\Qp]$, \cite[Prop.\ II.2.5(iv)]{fargues--scholze}.)

\begin{rem}
  \label{rem: do BC spaces explain non finitely presented ss}The
  failure of the skyscraper sheaves at irreducible points to be contained in the essential
  image of $\fA$ is related to the failure of supersingular
  representations of~$\GL_2(\Q_{p^2})$ to be finitely presented; see
  Section~\ref{subsubsec: irredquots} below. We do not know if there
  is a purely representation-theoretic interpretation of this failure
  of finite presentation in terms of $\mathcal{BC}(\cO(2))$. %
\end{rem}

\subsubsection{Multiplicity one}\label{subsec: mult one for D}We end this section by sketching a proof of the multiplicity one
result for~$D^\times$ mentioned above. For the sake of completeness we drop the
assumption that our representations have trivial central
character. Let~$D_{\Q}$ be a quaternion algebra over~$\Q$ which is
ramified precisely at $p,\infty$, and let
$\rbar:\Gal_{\Q}\to\GL_2(\Fpbar)$ be such that
$\rbar|_{\Gal_{\Q(\zeta_p)}}$ is irreducible. Let~$\m$ be the
corresponding maximal ideal of an appropriate prime-to-$p$ Hecke
algebra. (For simplicity the reader might imagine that $\rbar$ is
unramified outside~$p$.) For each character $\etabar^i$ of
$\cO_D^\times$ we have a corresponding space of quaternionic modular
forms~$M(\etabar^i)$ for the group $D_{\Q}^\times$, and a patched module $M_\infty(\etabar^i)$, which
is the reduction modulo~$p$ of a patched module $M_\infty(\eta^i)$,
where $\eta$ is the lift of~$\etabar$ to an $\cO^{\times}$-valued
character of~$\F_{p^2}^{\times}$. 

These patched modules have been studied in the papers~\cite{MR2822861,
  MR3449190, MR3411395}, and by modularity lifting theorems, one can
show that $M_\infty(\eta^i)$ is maximal Cohen--Macaulay over a ring
$R_\infty(\eta^i)$. The ring $R_\infty(\eta^i)$ is formally smooth
over a framed deformation ring $R_{\rhobar}(\eta^i)$
for~$\rhobar:=\rbar|_{\Gal_{\Qp}}$, which can be described as follows.

Firstly, if~$(p+1)|i$ then after twisting we may
suppose that~$i=0$, in which case $R_{\rhobar}(\eta^i)$ is a framed
deformation ring for non-crystalline semistable representations of
Hodge--Tate weights 0,1. Secondly, if $(p+1)\nmid i$, then
$R_{\rhobar}(\eta^i)$ is a framed deformation ring for potentially
crystalline representations of Hodge--Tate weights 0,1 and inertial
type $\eta^i\oplus\eta^{-i}$.

With this description at hand, one can easily check that the supports
of the patched modules agree with the prescriptions~\eqref{eqn:
  expectation for Dx 1}, \eqref{eqn: expectation for Dx 2}. (This is
the weight part of Serre's conjecture for quaternionic modular forms,
originally proved by T.G. and Savitt in~\cite{MR2822861}; the
refinement to distinguish between the $\sigma_i^{\pm}$ was proved by
Rozensztajn in~\cite[Lem.\ 6.1.1]{MR3411395}, by considering
``extended types'', i.e.\ by keeping track of an action of a
Frobenius.)

Now, if $R_{\rhobar}(\eta^i)$ is regular, $M_\infty(\eta^i)$ is
automatically free (this is the usual argument for deducing
multiplicity one results from Taylor--Wiles patching, due
independently to Diamond and Fujiwara). In particular, this holds if
$R_{\rhobar}(\eta^i)$ is either formally smooth or is formally smooth
over $\cO[[x,y]]/(xy-p)$ where $\cO$ is the ring of integers in an
appropriate unramified extension of~$\Qp$. The deformation rings
$R_{\rhobar}(\eta^i)$ have been computed in many cases. In the case $(p+1)|i$, they are
formally smooth (and can be described very explicitly, as all the
lifts are ordinary)  while if $(p+1)\nmid i$
then from \cite[Thm.\ 1.4]{SavittDuke} and \cite[Thm.\ 7.2.1]{MR3323575}, we know that they are regular except possibly for the cases
that $\rhobar$ is a twist of an unramified representation, or $\rhobar|_{I_p}$ is a twist of
$1\oplus\varepsilonbar$; we refer to these below as the exceptional
cases.

It can be proved
by ``pure thought'' (via a comparison to patched modules for
crystalline representations of Hodge--Tate weights $0,p-1$) that in
the case that $\rhobar$ is a direct sum of distinct unramified
characters, %
$R_{\rhobar}(\eta^i)$ is not regular. In fact by the results of~\cite{lehungmezardmorra}  it is formally smooth over
$\cO[[x,y]]/(xy-p^2)$ (at least if~$p$ is sufficiently large). In this case it seems unlikely that
one can prove the freeness of $M_\infty(\eta^i)$ without using some
geometric input. In the remaining cases the  deformation
rings   are again computed (for~$p$ sufficiently large) in~\cite{lehungmezardmorra}; see~\cite[\S 5.5]{lehungmezardmorra} for the details. %

The following argument for freeness in these exceptional cases
is taken from
unpublished work of George Boxer, Frank Calegari and T.G. (It was
independently discovered by Chengyang Bao, Andrea Dotto and Yulia
Kotelnikova.) %
By
Nakayama's lemma, it is enough to prove that the corresponding eigenspace
of quaternionic modular forms $M(\etabar^i)[\m]$ is 
one-dimensional over $\Fpbar$. Now, we interpret $M(\etabar^i)$ as
a space of functions on the supersingular points~$SS$ of a modular
curve~$X$. More precisely, after twisting, we can identify
$M(\etabar^i)_\m$ with the term $H^0(SS,\omega^{n+p-1})_\m$ in the short exact
sequence \[0\to H^0(X,\omega^n)_\m\to H^0(X,\omega^{n+p-1})_\m\to
  H^0(SS,\omega^{n+p-1})_\m\to 0\]where either $n=3$ and
$\rhobar|_{I_{\Qp}}\cong \varepsilonbar\oplus \varepsilonbar$, or $n=4$ and
$\rhobar|_{I_{\Qp}}\cong \varepsilonbar\oplus \varepsilonbar^2$. In either case it
follows from the weight part of Serre's conjecture (and the assumption
that $\rhobar$ is reducible) that
$H^0(X,\omega^n)_\m=0$, so $H^0(X,\omega^{n+p-1})_\m\isoto
  H^0(SS,\omega^{n+p-1})_\m$, and we need only show that
  $H^0(X,\omega^{n+p-1})[\m]$ is 1-dimensional. Suppose that this is
  not the case; then there is an eigenform with $a_1=0$, and by the
  $q$-expansion principle, such an eigenform is in the kernel of
  the~$\theta$ operator. But this kernel is trivial by a theorem of
  Katz~\cite{MR0463169}, so we are done.

\subsection{Preliminaries in the Banach case for  \texorpdfstring{$\GL_2(\Q_{p^{f}})$}{GL2(Qpf)}}
\label{subsec: Banach Qpf stuff}
Following the case of~$\GL_2(\Q_p)$, a natural next case to consider
is that of $G:= \GL_2(\Q_{p^f})$, for $f > 1$; here, as usual, $\Q_{p^f}$
denotes the unramified extension of~$\Q_p$ with residue
field~$\F_q := \F_{p^f}$, and we write $K:=\GL_2(\Z_{p^f})$ and $Z=Z(G)$.
For example, the mod $p$ local Langlands correspondence has been
studied quite intensively in this case,
as has the problem of mod $p$ local-global compatibility. (We recall
some of the history of these investigations as we go along; 
 the reader who is unfamiliar with it may wish to consult~\cite{MR2827792} and the
 introduction to~\cite{https://doi.org/10.48550/arxiv.2102.06188}.) 
On the other hand, the structure of the category of mod $p$ (or, more 
generally, mod $p^n$) smooth $G$-representations remains rather mysterious,
the classification of supersingular irreducible admissible representations
remains unknown, and as of yet there is no precise formulation of 
a mod $p$ or $p$-adic Langlands correspondence in the ``traditional'' mode
as exemplified
by \cite{MR2642409} and \cite{MR3150248} in the $\GL_2(\Q_p)$-case.

We believe that the categorical perspective sheds quite a bit of light 
on this case, helping both to illuminate known results and to indicate
why a traditional mode of formulation 
for the $p$-adic Langlands correspondence has proved so elusive. %
In the following discussion we elaborate on
these points, %
by describing some of our\footnote{Many of these ideas have been worked out
in joint work of M.E.\ and T.G.\ with Ana Caraiani, Michael
Harris, Bao Le Hung, Brandon Levin, 
and David Savitt, under the auspices of the NSF grant 1952556 ``Collaborative Research: Geometric Structures in the p-adic Langlands Program'', and we thank them for their permission to describe
these results here.}
preliminary results and
expectations for the categorical Langlands
program for the groups~$\GL_2(\Q_{p^f})$.  To ease notation
and fix ideas, we primarily focus on the case when~$f =2$, but we
begin in this section with some considerations that are valid for
any~$f$.

\begin{rem}
  \label{rem: where to see the highlights of GL2 Qp2}The reader who is
  familiar with the various difficulties and mysteries in the $p$-adic
  local Langlands program for $\GL_2(\Q_{p^2})$, and who is anxious to
  see what light we can shed on them, may wish to jump to
  Section~\ref{subsubsec: irredquots} and Remark~\ref{rem: we expect a
    unique supersingular pi with skyscraper H0}.
\end{rem}

\begin{rem}
  \label{rem: inherent derivedness}The similarities and differences
  between the case of $\GL_2(\Qp)$ and $\GL_2(\Q_{p^2})$ can be
  illustrated by the following example. In either case, for a generic
  Serre weight~$\sigma$, the functor~$\fA$ takes $\cInd_{KZ}^G\sigma$
  to an explicit line bundle supported on the corresponding
  component~$\cX(\sigma)$ of~$\cX_{\red}$. In the case of
  $\GL_2(\Qp)$, the correspondence is ``not very derived''; in
  particular, a simple calculation shows that the supersingular irreducible admissible representation
  $(\cInd_{KZ}^G\sigma)/T$ is taken to a (twisted) skyscraper sheaf supported at the
  corresponding irreducible Galois representation.

  For~$\GL_2(\Q_{p^2})$, the same calculation shows that
  $(\cInd_{KZ}^G\sigma)/T$ is taken to a sheaf supported on the closed
  locus on~$\cX(\sigma)$ corresponding to~$T=0$. However, this support
  is now one-dimensional, and the representation
  $(\cInd_{KZ}^G\sigma)/T$ has infinite length. In fact, we know (see
  \eqref{eqn:embedding}) that it admits as subrepresentations the full
  compact inductions $\cInd_{KZ}^G\sigma'$ for certain weights
  $\sigma'\ne \sigma$, and it follows easily that the 
  supersingular irreducible admissible representations of $\GL_2(\Q_{p^2})$ cannot possibly
  be taken to actual sheaves (as opposed to complexes of sheaves), in
  complete contrast to the case of~$\GL_2(\Qp)$  (or the case $\ell\ne p$).
\end{rem}

\subsubsection{The generic structure of $\cX_{\red}$}
As always we fix our coefficients~$\cO$ to be the ring 
of integers in some finite extension~$L$ of~$\Q_p$,
which we furthermore assume admits an embedding
of~$\Q_{p^f}$. (In other words, we may regard~$L$
as an extension of~$\Q_{p^f}$.)
We let~$k$ denote the residue field of~$\cO$.
We also fix a character $\zeta: \Q_{p^f}^{\times}
\to \cO^{\times}$,
and let $\cX$ denote the Noetherian formal algebraic stack 
classifying rank $2$ projective \'etale
$(\varphi,\Gamma)$-modules 
for~$\Gal_{\Q_{p^f}}$ over $p$-adically complete $\cO$-algebras,
having their determinant fixed to equal~$\zeta \varepsilon^{-1}$.

\subsubsection{The components~$\cX(\sigma)$}\label{subsubsec:
  components X sigma}
Just as in the case of $\GL_2(\Q_p)$,
the underlying 
reduced algebraic substack~$\cX_{\red}$ is of finite presentation
over~$\Spec k$,
and is a union of finitely many irreducible
components $\cX(\sigma)$ labelled 
by Serre weights~$\sigma$ having central character $\zeta_{| \cO^{\times}}$,
with the
exception that there is a pair of components $X(\sigma)^{\pm}$ when $\sigma$
is a twist of the Steinberg weight.
And again, as in that case,
for a {\em generic} $\sigma$
(one of the form 
$$\sigma = (\Sym^{a_0} \otimes \det{}^{b_0} ) \otimes
\cdots\otimes (\Sym^{a_i}\otimes \det{}^{b_i})^{\Fr^i} \otimes \cdots
\otimes (\Sym^{a_{f-1}} \otimes \det{}^{b_{f-1}})^{\Fr^{f-1}}$$
with all the $a_i$ satisfying $0 \leq a_i \leq p-3$) 
we have that $\cX(\sigma)$ can be described
as a moduli stack of Fontaine--Laffaille modules with coefficients
in $k$-algebras.

In general, for a reductive group $G$ over~$\Fp$,
the corresponding Fontaine--Laffaille moduli space can be described as
the quotient  $(U\backslash G)/T$, where $U$ is the unipotent radical of some
fixed Borel in~$G$, acting by left translation, and $T$ is the maximal
torus of the the same Borel, acting by Frobenius twisted conjugation. 
In our case, the group $G$ is the restriction of scalars of~$\SL_2$ from
$\F_{p^f}$ to~$\Fp$.
Since~$k$ is assumed to contain~$\F_{p^f}$,
we have $$G = G_{/k} = \SL_2^f,$$
$$U = U_{/k} = 
(\begin{pmatrix} 1 & \Ga \\ 0 & 1 \end{pmatrix})^f,
$$
and
$$T = T_{/k} = \Gm^f$$
with Frobenius cyclically permuting the factors.

Now %
$\big(\begin{smallmatrix} 1 & \Ga \\ 0 & 1 \end{smallmatrix}\big)\backslash \SL_2
= \A^2\setminus \{0\}$
(via $\big(\begin{smallmatrix} * & * \\ x & y \end{smallmatrix}\big) \mapsto (x ,y)$).
Thus we have
$$U\backslash G := \prod_{i=0}^{f-1}  (\A^2\setminus \{0\} )_{(x_i,y_i)} $$
(where the subscripts indicate our notation for the coordinates).
The torus $T$ then acts via
$$(t_0,\ldots,t_{f-1}) \cdot (x_i, y_i) = 
(t_{i-1} t_{i} x_i, t_{i-1} t_i^{-1} y_i).$$
So (for generic $\sigma$) we 
find that $\cX(\sigma)$ is equal to the quotient stack
$$\cX(\sigma)
= 
[\bigl(\prod_{i=0}^{f-1}  (\A^2\setminus \{0\} )_{(x_i,y_i)}\bigr) / \Gm^f] $$
with the torus action being given by the preceding formula.

The coordinates $(x_i,y_i)$ on~$\cX(\sigma)$ are evidently not
completely canonical, but they are not completely without meaning
either; in particular, we will shortly see that their zero loci describe the loci of intersection 
of~$\cX(\sigma)$ with other components. 
A special role
is played by the product
$$T := y_0\cdots y_{f-1}.$$
Evidently, $T$ is $\Gm^f$ invariant, and so defines a regular function
on~$\cX(\sigma)$.  
In fact, it generates the ring of functions on~$\cX(\sigma)$,
as we record in the following lemma.

\begin{lemma}\leavevmode
\begin{enumerate}
\item
$\Gamma\bigl(\cX(\sigma),\cO_{\cX(\sigma)}\bigr)
= k[T].$
\item
The distinguished open subset $D(T)$ {\em (}i.e.\ the non-vanishing
locus of~$T$\emph{)} in $\cX(\sigma)$ 
is precisely the niveau $1$ ordinary locus in~$\cX(\sigma)$.
\end{enumerate}
\end{lemma}

Concretely, we see that
$$D(T) = [\A^f_{(x_0,\ldots,x_{f-1})}/(\Gm)_t]
\times (\Gm)_T,$$
where $t$ acts via
$$t\cdot (x_0,\ldots,x_{f-1}) = (t^2 x_0, \ldots, t^2 x_{f-1});$$ 
so $(\Gm)_t$ can be identified with the diagonal copy of $\Gm$ in $\Gm^f$.
There is a character $\chi: 
\Gal_{\Q_{p^f}} \to k^{\times}$, depending on~$\sigma$, with
the property 
that the point $(T,x_0,\ldots,x_{f-1})$ of $D(T)$ corresponds to a two-dimensional
representation $\rhobar$ of $\Gal_{\Q_{p^f}}$ of the form
$$\begin{pmatrix} \ur_T \chi & * \\ 0 &
\ur_{T^{-1}} \zeta \varepsilon^{-1} \chi^{-1} \end{pmatrix},$$
where $*$ is an element of $\Ext^1_{\Gal_{\Q_{p^f}}}( 
\ur_{T^{-1}} \zeta \varepsilon^{-1} \chi^{-1},
\ur_T \chi),$
an $f$-dimensional vector space with coordinates~$x_0,\ldots,x_{f-1}$.

The zero locus $Z(T):= \cX(\sigma) \setminus D(T)$
is more complicated in its geometry; we describe it
in detail in the case~$f= 2$ in Section~\ref{subsec: Banach Qp2 stuff}. %

Suppose now that $\sigma$ is \emph{very generic} in the sense
that each $a_i$ satisfies $0 < a_i < p-3$. %
Then there are $2f$ Serre weights $\sigma'$ such 
$\Ext^1_{\GL_2(\F_{q})}(\sigma,\sigma') \neq 0,$ and each such~$\sigma'$ is generic. %
One then finds that $\cX(\sigma)$ and $\cX(\sigma')$ intersect in 
one of the codimension $1$ closed substacks of $\cX(\sigma)$
cut out by either $x_i = 0$ or $y_i = 0$ (for some choice of~$i$).

\subsubsection{Picard groups}Continue to assume that~$\sigma$ is generic.
From the description of $\cX(\sigma)$ as a quotient stack,
it is easy to see that $\Pic\bigl(\cX(\sigma)\bigr) \cong \Z^f$
(the character lattice of $\Gm^f$).
This isomorphism is not evidently canonical, but we now explain how to make it %
canonical up to permuting the coordinates of~$\Z^f$,
under the additional assumption that $\sigma$ is very generic.

To begin with, we consider the open substack $D(T)$ of $\cX(\sigma)$.
As already noted, this admits the description $D(T) = (\A^f/\Gm)\times \Gm$,
from which one deduces that
$\Pic\bigl( D(T) \bigr) = \Z$ (the character lattice of~$\Gm$).
Furthermore, this identification {\em is canonical}, since the 
$\Gm$ that appears is naturally oriented (corresponding to the
fact that $D(T)$ parameterizes extensions of $\ur_{T^{-1}}\zeta \varepsilon^{-1}
\chi^{-1}$ by $\ur_T \chi$ in a fixed direction).  

For each $i = 0,\ldots,f-1$,
write $\sigma_i$ for the weight such that $\cX(\sigma) \cap
\cX(\sigma_i) = \cX(\sigma)^{y_i = 0};$
since $\sigma$ is very generic each $\sigma_i$ is again generic.
We write $T_i$ for the
analogue of the function $T$ on $\cX(\sigma_i)$;
then 
$$(T_i)_{|\cX(\sigma_i)^{y_i  = 0}} = \dfrac{x_{i-1}}{x_i} \prod_{j \neq i-1,i} y_j
$$
Thus
$$ \cX(\sigma)^{y_i = 0} \cap D(T_i) = \cX(\sigma)^{y_i = 0,(x_{i-1}\prod_{j \neq i-1,i} y_j) \neq 0}.$$
Now
$$D(x_{i-1}x_i\prod_{j \neq i-1,i}y_j)
=
[\A^{f}_{(x_0,\ldots,x_{i-2},y_{i-1},y_i,x_{i+1}, \ldots,x_{f-1})}/(\Gm)_t],$$
where $t$ acts via
\begin{multline*}t \cdot (x_0,\ldots,x_{i-2},y_{i-1},y_i,x_{i+1}, \ldots,x_{f-1})
\\
= (t^2 x_0, \ldots, t^2 x_{i-2}, t^2y_{i-1}, t^{-2} y_i, t^2 x_{i+1},
\ldots, t^2 x_{f-1});
\end{multline*}
the torus $(\Gm)_t$ can be regarded as the copy of $\Gm$ embedded into $\Gm^f$
via $$t\mapsto (t,\ldots,t,t^{-1} \text{ in the } i-1 \text{  slot}, t,\ldots,t).$$
The closed substack $\cX(\sigma)^{y_i  = 0}\cap D(T_i)$
is then the copy of $[\A^{f-1}/\Gm]$ cut out by $y_i = 0.$
Of course the embedding
$$\cX(\sigma)^{y_i = 0} \cap D(T_i) \hookrightarrow D(T_i)$$ embeds
$[\A^{f-1}/\Gm]$ into {\em a different copy} of $[\A^f/\Gm]$.

We've already observed that $\Pic\bigl(D(T_i)\bigr) = \Pic([\A^{f}/\Gm]) = \Z$,
and the restriction map $\Pic\bigl(D(T_i)\bigr) \to \Pic([\A^{f-1}/\Gm]) = \Z$
is an isomorphism. 
Thus $\Pic\bigl(\cX(\sigma)^{y_i = 0} \cap D(T_i)\bigr) = \Z$  (canonically).
We  now consider the product of the restriction morphisms
\nummultline 
\label{eqn:Pic restriction}
\Z^f \cong \Pic\bigl(\cX(\sigma)\bigr)
\to \Pic\bigl(D(T)\bigr) \times \prod_{i = 0}^{f-1}  \Pic\bigl( \cX(\sigma)^{y_i = 0}\cap
D(T_i)\bigr) 
\\
 = \Z \times \Z^{f} = \Z^{f+1}.
\end{multline}
Recalling the descriptions of the $D(T)$  and the various $\cX(\sigma)^{y_i = 0}
\cap D(T_i)$ as quotients $[\A^f/(\Gm)_t]$ and $[\A^{f-1}/(\Gm)_t]$, and of the embeddings 
of these various copies of $(\Gm)_t$ into $(\Gm)^f$,
we find that this map is given by
$$(a_0,\ldots, a_{f-1}) \mapsto 
\bigl(\sum_j a_j, (\sum_{j  \neq {i-1}} a_j  - a_{i-1})_{i=0}^{f-1} \bigr)
=  
\bigl(\sum_j a_j, (\sum_{j=0}^{f-1} a_j  - 2 a_{i-1})_{i=0}^{f-1} \bigr)
.$$
So we see that~\eqref{eqn:Pic restriction} is an  embedding, and that we can recover
the various $a_i$ from the coordinates in $\Z^{f+1}$;
explicitly, the (left) inverse is given by 
$$(b,b_0,\ldots,b_i, \ldots, b_{f-1}) \mapsto \dfrac{1}{2}(b-b_1, \ldots, b - b_{i+1},
\ldots,  b- b_{0}).$$ 
Now the coordinates $(b,b_0,\ldots,b_{f-1})$ are
canonical, up to permuting the $b_i$ (corresponding to a relabelling of the $\sigma_i$),
and so we see
that the isomorphism $\Pic\bigl(X(\sigma)\bigr) \cong \Z^f$
is canonical, up to permuting the coordinates~$a_i$.

\subsubsection{Tame types}
If $\tau$ is a tame inertial type for~$\Q_{p^f}$,
then the associated $K$-representation $\sigma^{\crys}(\tau)$
factors through $\GL_2(k_F) = \GL_2(\F_{p^f})$. 
If $\sigma^{\crys}(\tau)^{\circ}$ is a $K$-invariant  lattice
in~$\sigma^{\crys}(\tau)$, then the multi-set of Jordan--H\"older
factors of its reduction is in fact a set of distinct Serre
weights which is independent of the particular choice of lattice,
and which we denote by
$\JH\bigl(\overline{\sigma^{\crys}(\tau)}\bigr).$
If $\tau$ is furthermore sufficiently generic, then 
$\JH\bigl(\overline{\sigma^{\crys}(\tau)}\bigr)$
consists of exactly $2^f$ weights.

We write $\cX(\tau) = \bigcup_{\sigma \in \JH\bigl(\overline{\sigma^{\crys}(\tau)}\bigr)}
\cX(\sigma);$
this is a closed substack of~$\cX_{\red}$.

\begin{thm}
If $\tau$ is a sufficiently generic tame type, then
$\Pic\bigl(\cX(\tau)\bigr) \cong \Z^{2f}.$
\end{thm}

In fact, the product of the restriction maps
$$\Pic\bigl(\cX(\tau)\bigr)
\hookrightarrow 
\prod_{\sigma \in \JH\bigl(\overline{\sigma^{\crys}(\tau)}\bigr)}
\Pic\bigl(\cX(\sigma)\bigr)$$
is injective, and one can identify the image of
$\Pic\bigl(\cX(\tau)\bigr)$ as a certain specific subgroup of the target.

\subsubsection{Values of the functor}
The results of \cite{MR3323575} show that if $\sigma$ is a  generic Serre
weight, or an appropriately chosen lattice
$\sigma^{\crys}(\tau)^{\circ}$
for a sufficiently generic tame type $\tau$ (having determinant
equal to $\zeta|_{\cO^\times_F}$),
then the patched modules $M_{\infty}(\sigma)$
and
$M_{\infty}\bigl(\overline{\sigma^{\crys}(\tau)^{\circ}}\bigr)$
(constructed from an appropriate
Shimura curve context) are free of rank one over their support. The
weight part of Serre's conjecture and the (geometric) Breuil--M\'ezard
conjecture for potentially Barsotti--Tate representations imply that
the support of $M_{\infty}(\sigma)$ is precisely the pullback to the
corresponding versal ring of $\cX(\sigma)$, and
similarly the support of
$M_{\infty}\bigl(\overline{\sigma^{\crys}(\tau)^{\circ}}\bigr)$ is the
pullback of
$\cX(\tau)$. (See~\cite[Thm.\ 8.6.2]{emertongeepicture}.)
Thus we expect that the value of the functor $\fA$ on
$\cInd_{KZ}^G \sigma$ or
$\cInd_{KZ}^G \overline{\sigma^{\crys}(\tau)^{\circ}}$
will be a line bundle on $\cX(\sigma)$ or $\cX(\tau)$ respectively.
This explains our interest in computing the Picard groups
of these stacks.

More generally, if $V$ is a representation 
of $\GL_2(k_F)$ over $F$ which is multiplicity free (i.e.\ each of whose
Jordan--H\"older factors occurs with multiplicity exactly one),
and which has central character equal to~$\overline{\zeta}$,
then the value of $\fA$ on $\cInd_{KZ}^G V$ should be
scheme-theoretically supported on 
$\bigcup_{\sigma \in \JH(V)} \cX(\sigma)$, and in particular on
$\cX_{\red}$. (This follows from the expected compatibility with the
Breuil--M\'ezard conjecture, and the expected maximal Cohen--Macaulay property
of $\fA(\cInd_{KZ}^G V)$.)

\begin{expectedthm}[Caraiani--E.--G.--Harris--Le Hung--Levin--Savitt]
\label{expthm:FRG}%
There is a fully faithful %
functor $\fA$ from the additive category of finitely presented smooth representations generated by the objects
$\cInd_{KZ}^G V$ \emph{(}with~$V$ multiplicity free, and having all
Jordan--H\"older factors being sufficiently generic\emph{)} to the abelian category 
$\Coh(\cX_{\red})$.
\end{expectedthm}

The functor $\fA$ of this theorem, which we certainly believe to be
the restriction of the functor $\fA$ of Conjecture~\ref{conj: Banach
  functor} for~$\GL_2(\Q_{p^f})$, is produced by an explicit
construction, building on the description of the various Picard groups
given above. In particular, we have the following expectation.

\begin{expectation}
\label{exp:Serre weights}
For $G  = \GL_2(\Q_{p^{f}})$,
and for a very generic
Serre weight~$\sigma$, the sheaf $\fA(\cInd_{KZ}^G\sigma)$ is the line
bundle $\cO(-1,\ldots,-1)$ on~$\cX(\sigma)$, i.e.\ the line bundle
corresponding to the element
$(-1,\dots,-1)\in\Pic(\cX(\sigma))\cong\Z^f$. 
\end{expectation}
(Note that $(-1,\dots,-1)$  is invariant under permutation,
and so does indeed give a canonically defined element of~$\Pic\bigl(\cX(\sigma)\bigr).$)

The functor $\fA$
of Expected Theorem~\ref{expthm:FRG} is defined by imposing Expectation~\ref{exp:Serre
weights} as a definition; it turns out that
the extension of $\fA$ to the additive subcategory in its statement is then essentially forced.

\begin{rem}
  \label{rem: we're forced to make this guess for the line bundles}In the case $f=1$ 
Expectation~\ref{exp:Serre weights} agrees with the results of~\cite{DEGcategoricalLanglands},
as explained in Section~\ref{subsubsec:fA on Serre weights}. In general the
expected compatibility with patching, and the weight part of Serre's
conjecture, imply that  $\fA(\cInd_{KZ}^G\sigma)$ will be a line
bundle on~$\cX(\sigma)$. One can then check that the full faithfulness
of~$\fA$, and the possible extensions between the various
$\cInd_{KZ}^G\sigma$, imply that if the corresponding element of the
Picard group --- identified
with $\Z^f$ as above --- is independent of~$\sigma$,
then it is necessarily
$(-1,\dots,-1)$. The interested reader is invited to check this in the
case~$f=2$ using the calculations of Section~\ref{subsec: Banach Qp2
  stuff} below. (Alternatively, this can be seen from the expected
compatibility with duality, as mentioned in Section~\ref{subsubsec:fA
  on Serre weights} in the case $F=\Qp$.)

We also note that the discussion of Remark~\ref{rem: guessed the twists following Breuil}
below shows that Expectation~\ref{exp:Serre weights} is compatible %
with the philosophy introduced by Breuil in~\cite{MR2783977}.  
\end{rem}

\subsection{Partial results and examples in the Banach case for \texorpdfstring{$\GL_2(\Q_{p^{2}})$}{GL2(Qp2)}}
\label{subsec: Banach Qp2 stuff}
In order to continue, we restrict to the case $f = 2$, 
in which case we can make much of the preceding discussion quite explicit.
To begin with, we simplify the notation by writing
$(u,v),(x,y)$
(rather than~$(x_0,y_0), (x_1,y_1)$)
for the coordinates on $G/U$,
and $(s,t)$ (rather than $t_0,t_1$) for the coordinates on~$T$.

Then, if $\sigma$  is a generic Serre weight,
we have
\numequation\label{eqn:X description}
\cX(\sigma) =  
[
(\A^2\setminus\{0\})_{(u,v)}
\times
(\A^2\setminus\{0\})_{(x,y)}
/ (\Gm\times \Gm)_{(s,t)} ],\end{equation}
where
$$(s,t) \cdot  
\bigl( (u,v), (x,y) \bigr) = 
\bigl((stu,s^{-1}tv), (stx, st^{-1}y) \bigr).$$
Because the diagonal copy of $\boldsymbol{\mu}_2$ in $T$ acts trivially in this formula, and because (as it
furthermore turns out) these copies of $\boldsymbol{\mu}_2$ (for the various generic Serre weights~$\sigma$)
act trivially on all the sheaves that we will consider, 
it is notationally convenient to work instead on stacks that we
will denote~$\cC_{\sigma}$,
defined as follows.

We set $a = st,$ $b = s^{-1} t,$
so that \numequation
\label{eqn:torus isogeny}
(s,t) \mapsto (a,b)
\end{equation} is a self-isogeny of $\Gm\times \Gm$ having
the diagonal copy of $\boldsymbol{\mu}_2$ as its kernel.
As already noted, the action of  $\Gm\times\Gm$ on 
$(\A^2\setminus\{0\})_{(u,v)}
\times
(\A^2\setminus\{0\})_{(x,y)}$
occurring in the definition of~$\cX(\sigma)$ then factors through
this isogeny, 
and we write $\cC_{\sigma}$ for the quotient of
$(\A^2\setminus\{0\})_{(u,v)}
\times
(\A^2\setminus\{0\})_{(x,y)}$
by the corresponding isogenous copy of $\Gm\times\Gm$.
Concretely,  we write
\numequation\label{eqn:C description}\cC_{\sigma} =  
[
(\A^2\setminus\{0\})_{(u,v)}
\times
(\A^2\setminus\{0\})_{(x,y)}
/ (\Gm\times \Gm)_{(a,b)} ],\end{equation}
where
$$(a,b) \cdot  
\bigl( (u,v), (x,y) \bigr) = 
\bigl((au,bv), (ax, b^{-1}y) \bigr).$$
There is then a natural map $\cX(\sigma) \to \cC_{\sigma}$,
and the various stacks $\cC_{\sigma}$ glue together,
in a manner that emulates the gluing together of the~$\cX(\sigma)$ in $\cX_{\red}$,
to produce a reduced algebraic stack that we will denote~$\cC_{\red}$.

There are two open covers of $\cC_{\sigma}$ that we wish to
introduce. %
To describe them, we first remark that,
as usual, if $f$ is a function on a stack, or more generally, a section
of a line bundle, then we let $D(f)$ denote the non-vanishing locus of~$f$; 
this is an open subset of the stack in question.
We then set $\cU_{\sigma} := D(vy)$, and $\cV_{\sigma} := D(u) \cup D(x),$
and note that $\cC_{\sigma} = \cU_{\sigma} \cup \cV_{\sigma}$.
We also write $\cZ_{\sigma} := \cC_{\sigma} \setminus \cU_{\sigma}.$

There is an evident isomorphism
$$\cU_{\sigma} = [(\A^2_{(u,x)}/(\Gm)_a]  \times (\Gm)_{vy}$$ 
(here the subscripts denote the coordinates), where of course
$a \cdot (u,x) = (a  u, a  x).$
The Galois-theoretic interpretation of this locus is straightforward:
the points correspond to extensions of a character~$\chi_2$ by a character~$\chi_1 $, with $\chi_1|_{I_{\Q^{p^{2}}}},\chi_2|_{I_{\Q^{p^{2}}}}$ depending only on~$\sigma$; the quantity $vy$ corresponds to a certain Frobenius eigenvalue appearing
in these characters,
and the $\A^2$ is the $2$-dimensional $\Ext^1$ between these characters.
More precisely, if (after possibly twisting) we write $\sigma = \Sym^{a_0} \otimes (\Sym^{a_1})^{\Fr}$ with
$0 \leq a_1,a_2  \leq p-3$, then (for an appropriate choice of conventions) the character~$\chi_2 $ is unramified, while $\chi_1 |_{I_{\Q_{p^{2}}}}=\omega_0 ^{a_0+1}\omega_1 ^{a_1 +1}$, where~$\omega_0$ and $\omega_1=\omega_0 ^{p} $ are fundamental characters of~$\Gal_{\Q_{p^{2}}}$ of niveau one (see for example \cite[\S2.1]{GLSII} for more details).

\[\includegraphics[width=0.8\textwidth]{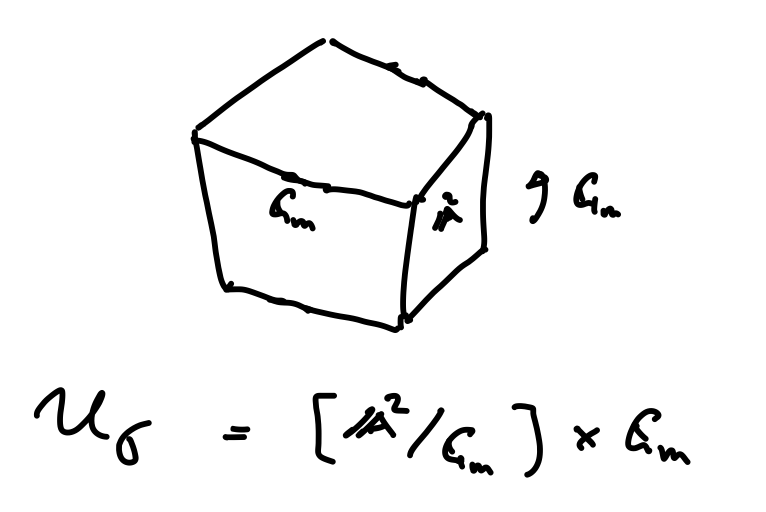}\]

There is an isomorphism %
$$\cZ_{\sigma} =  \Pone_{[u:x]} \times  
[\set{ (v,y) \, | \, vy = 0 }/(\Gm)_b]
\setminus ( \set{u = v = 0} \cup \set{x = y = 0} )
,$$
where $b$ acts via $b \cdot (v,y)  = (b v, b^{-1} y).$
If we consider the open substack of $\cZ_{\sigma}$ where $ux \neq 0,$
we obtain the stack
$$ (\Gm)_{x/u} \times [\set{ (v,y) \, | \, vy = 0 }/(\Gm)_b],$$
which is a family of extensions of characters~$\chi_3 ,\chi_4 $ of certain fixed inertial weights,
the two lines corresponding to extensions in either  possible  directions, but which
are ``split along one of the two possible embeddings''.
The variable $x/u$ again parameterizes a 
Frobenius eigenvalue appearing in these characters.  (In the notation of the previous paragraph, we have $\chi_3 |_{I_{\Q_{p^{2}}}}=\omega_0 ^{a_0 +1}$, $\chi_4 |_{I_{\Q_{p^{2}}}}=\omega_1 ^{a_1 +1}$.)
The points $\set{u  = y = 0}$ and $\set{v = x = 0}$ each correspond to an irreducible
Galois representation (there being two such representations up to isomorphism with~$\sigma$ as a Serre weight; see \cite[Defn.\ 2.4]{GLSII}).

\[\includegraphics[width=0.8\textwidth]{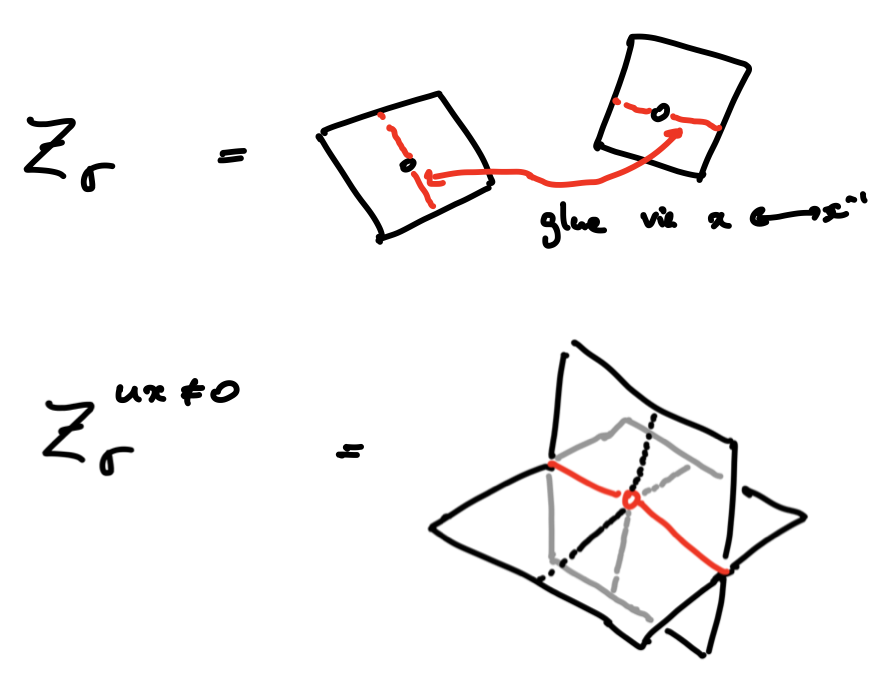}\]

The open cover $\set{\cU_{\sigma}, \cV_{\sigma}}$
of $\cC_{\sigma}$ is useful when we want to think about points of $\cX(\sigma)$
in terms of the Galois representations that they correspond to.  However,
there is another open cover we will work with, which at first
seems less conceptually meaningful,
but which leads to simple computations (especially when  we consider how the
various components $\cX(\sigma)$ intersect one another);
this  is  the open cover $\set{ D(u), D(v)}$.  

Concretely, we see that
\numequation
\label{eqn:D(u) description}
D(u) = [ \bigl(  \A^1_{v} \times  (\A^2\setminus  \set{0})_{x/u,y}\bigr) /(\Gm)_b ],
\end{equation}
where $(\Gm)_b$ acts via $b\cdot (v, x/u, y) = (b v, x/u, b^{-1} y)$,
while
\numequation
\label{eqn:D(v) description}
D(v) = [ \bigl( \A^1_u \times  (\A^2 \setminus \set{0})_{x, vy}\bigr) / (\Gm)_a],
\end{equation}
where $(\Gm)_a$ acts via $a \cdot (u, x, vy) = (a u, a x, vy).$
We can also describe their intersection:
$$D(u) \cap D(v) = D(uv) = (\A^2\setminus \set{0})_{x/u, vy}.$$

Since $(\A^2\setminus \set{0} ) \times (\A^2 \setminus \set{0})$ 
has trivial Picard group, we see from~\eqref{eqn:X description} that
$$\Pic(\cX(\sigma)) = \X^{\bullet}\bigl( (\Gm\times \Gm)_{s,t}\bigr) = \Z \times \Z,$$
and from \eqref{eqn:C description} that 
$$\Pic(\cC_{\sigma}) = \X^{\bullet}\bigl( (\Gm\times\Gm)_{a,b}\bigr) = \Z \times \Z.$$
The pull back map $\Pic(\cC_{\sigma} ) \to \Pic(\cX(\sigma))$ is just
the map of character lattices induced by the isogeny~\eqref{eqn:torus  isogeny},
and so is given by the matrix $\big(\begin{smallmatrix} -1 & 1 \\ -1 &
  -1 \end{smallmatrix}\big).$%

Similarly, we see that
$$\Pic\bigl(D(u)\bigr) = X^{\bullet}\bigl((\Gm)_b\bigr) = \Z,$$
while
$$\Pic\bigl(D(v)\bigr) = X^{\bullet}\bigl((\Gm)_a\bigr) = \Z.$$
Thus the product of the restriction  maps
$$\Pic(\cC_{\sigma}) \to \Pic\bigl(D(v)\bigr) \times \Pic\bigl(D(u)\bigr)$$
is an isomorphism, indeed is just the identity map from $\Z\times \Z$
to itself. 
This is one reason that the open cover $\set{D(u),D(v)}$ is useful: we
can describe a line bundle on  $\cC_{\sigma}$ simply by describing its
restrictions to each of $D(u)$ and $D(v)$, and line bundles on each of these
open substacks are given by prescribing a single twist (i.e.\ a single
integer).

Following Expectation~\ref{exp:Serre weights}
(and taking into account our change of coordinates~\eqref{eqn:torus isogeny}),
for a generic Serre weight $\sigma$ we define 
$$ \gA(\cInd_{KZ}^G \sigma) := \cO_{\cC_{\sigma}}(-1,0).$$
The endomorphism $T$ of $\cInd_{KZ}^G \sigma$ will then
correspond to the endomorphism of $\cO_{\cC_{\sigma}}(-1,0)$ given
by multiplication by~$vy$.

\subsubsection{Intersections and unions of generic components}
Two generic components $\cX(\sigma)$ and $\cX(\sigma')$ intersect in a codimension one locus precisely if $\sigma$ and $\sigma'$ admit a non-trivial extension as representations of $\GL_2(\F_{p^2})$.
For a sufficiently generic $\sigma$,  there are %
four such $\sigma'$. The corresponding loci
of intersection on $\cX(\sigma)$  are respectively  $\set{u=0},$ $\set{v = 0}$, $\set{x=0},$
and~$\set{y=0}$.   On each of the corresponding components~$\cX_{\sigma'}$, these same
loci are given by $\set{v' = 0}$, $\set{u' = 0}$, $\set{y' = 0}$, and~$\set{x' = 0}$.
(Here we have used primes to denote the variables on~$\cX_{\sigma'}$.)
These loci are compatible with passing from the $\cX(\sigma)$ to the~$\cC_{\sigma}$,
and indeed the whole picture of~$\cX_{\red}$ that we will describe descends 
to the $\cC_{\sigma}$ level, and this is where we work from now on.

\subsubsection{A brief interlude on representation theory}
\label{subsubsec:rep theory interlude}
If $\xi$ is a character of the torus $T(\F_q) :=\F_q^{\times}\times \F_q^{\times}$
in $GL_2(\F_q)$ (with values in $k^{\times}$, say, although the following discussion
is independent of the particular field or ring of coefficients),
then we can define the principal series 
$\PS(\xi) := \Ind_{B(\F_q)}^{\GL_2(\F_q)} \xi$
(where $B$ denotes the subgroup scheme of upper triangular matrices in $\GL_2$). 
We may inflate the $\GL_2(\F_q)$-action on~$\PS(\xi)$ to a~$K$-action,
and, assuming that $\xi$ is suitably compatible with our central character $\zeta$,
we can then extend this $K$-action to a $KZ$-action by  having $Z$ act via~$\zeta$.
We may then form the compact induction
$\cInd_{KZ}^G \PS(\xi)$.
If we let $s$ denote
the non-trivial Weyl group involution of $T(\F_q)$,
then we may similarly consider $\PS(\xi^s)$ and~$\cInd_{KZ}^G \PS(\xi^s)$.

Now if $\xi\neq \xi^s$,then $\PS(\xi)$ and $\PS(\xi^s)$
are not isomorphic as 
$\GL_2(\F_q)$-representations (and hence are not isomorphic
as $KZ$-representations).
However, there {\em is} an isomorphism 
\numequation
\label{eqn:PS intertwiner}
\cInd_{KZ}^G \PS(\xi)
\iso \cInd_{KZ}^G \PS(\xi^s).
\end{equation}
To see this, write 
$$\cInd_{KZ}^G \PS(\xi) = \cInd_{IZ}^G \xi$$
and
$$\cInd_{KZ}^G \PS(\xi^s) = \cInd_{IZ}^G \xi^s,$$
and note that 
the matrix $\big(\begin{smallmatrix}0 & 1 \\p &0\end{smallmatrix}\big) \in G$
normalizes~$I$
and 
induces the action of~$s$ on~$T(\F_q)$.

The isomorphism~\eqref{eqn:PS intertwiner} is perhaps the most fundamental
example of a morphism in the category of smooth $G$-representations 
which is not obtained by 
inducing a morphism of $KZ$-representations, and plays an important (if
sometimes only implicit) role in our analysis below.

 \subsubsection{A union of two components, and extensions of Serre weights}
\label{subsubsec:two components}
We begin by describing $\cC_{\sigma} \cup \cC_{\sigma'}$
in the case when these two components meet in codimension~$1$ (i.e.\
in the case that $\sigma,\sigma'$ admit a non-trivial extension).
Up to switching $\sigma$  and $\sigma'$ and swapping the variables $u$ and $v$
or $x$ and~$y$, we may assume that the intersection is along the locus
$u = 0$ (in $\cC_{\sigma}$), which coincides with the locus $v' = 0$ (in
$\cC_{\sigma'}$).  This locus is evidently disjoint from (indeed, equal
to the complement in $\cC_{\sigma}$ of)~$D(u)$, 
and so is contained in $D(v)$.  Similarly,
this locus is disjoint from $D(v')$ in $\cC_{\sigma'},$
and so is contained in $D(u')$.
Thus
\numequation
\label{eqn:open cover of two components}
\set{ D(u), D(v) \cup D(u'), D(v')}
\end{equation}
is an open cover of $\cC_{\sigma} \cup \cC_{\sigma'}$.  

Recall from~\eqref{eqn:D(v) description} that
$$
D(v) = [ \bigl( \A^1_u \times  (\A^2 \setminus \set{0})_{x, vy}\bigr) / (\Gm)_a],
$$
while
applying~\eqref{eqn:D(u) description} to $\cC_{\sigma'}$ shows that
$$
D(u') = [ \bigl(  \A^1_{v'} \times  (\A^2\setminus  \set{0})_{x'/u',y'}\bigr) /(\Gm)_{b'} ].
$$
The union $D(v) \cup D(u')$ is then equal to
$$
[ \set{(u,x,vy,v') \, | \,   (x,vy) \neq (0,0), u v' = 0} / (\Gm)_a ],
$$
where the locus $\set{v' = 0}$ is identified with $D(v)$ in the evident way,
while the locus $\set{u = 0}$ is identified with $D(u')$ via
$b'  = a^{-1}, x'/u' = vy, y' = x.$  

Since $\A^1\times (\A^2 \setminus \set{0})$ has trivial Picard,
we see that restriction to either $D(v)$ or $D(u')$ induces isomorphisms
$$\Pic \bigl( D(v) \cup D(u') \bigr) \iso \Pic\bigl(D(v)\bigr), \Pic\bigl( D(u') \bigr)
\iso \Z,$$
the two different identifications with $\Z$ differing by a sign (since
$b' = a^{-1}$).
A consideration of the open cover~\eqref{eqn:open cover of two components}
then shows that the product of the restriction maps induces an isomorphism
\begin{multline*}
\Pic( \cC_{\sigma} \cup \cC_{\sigma'})
\iso
\Pic\bigl( D(v') \bigr) \times \Pic\bigl( D(v) \cup D(u') \bigr) \times \Pic\bigl( D(u) \bigr)
\\
\iso 
\Pic\bigl( D(v') \bigr) \times \Pic\bigl( D(v) \bigr) \times \Pic\bigl( D(u) \bigr)
=
\Z \times \Z \times \Z.
\end{multline*}
In terms of these coordinates,
we find that the restriction
$$ \Pic( \cC_{\sigma} \cup \cC_{\sigma'}) \to
\Pic(\cC_{\sigma}) = \Pic\bigl( D(v) \times D(u) \bigr) = \Z\times \Z$$
is given by
$$(l,m,n) \mapsto (m,n),$$
while 
$$ \Pic( \cC_{\sigma} \cup \cC_{\sigma'} ) \to
\Pic(\cC_{\sigma'}) = \Pic\bigl( D(v') \times D(u') \bigr) = \Z\times \Z$$
is given by
$$(l,m,n) \mapsto (l,-m).$$

We now present some examples of morphisms of coherent sheaves
on $\cC_{\sigma}\cup \cC_{\sigma'}$ which incorporate some aspects
of the preceding discussion.  These examples will play an important
role in our  discussion of the functor~$\mathfrak{A}$.

\begin{example}
Consider the short exact sequence
$$0 \to \cI_{\cC_{\sigma}} \to \cO_{\cC_{\sigma} \cup \cC_{\sigma'}} 
\to \cO_{\cC_{\sigma}} \to 0;$$
here, of course, $\cI_{\cC_{\sigma}}$ denotes the ideal sheaf
cutting out $\cC_{\sigma}$ as a closed substack of the union
$\cC_{\sigma} \cup \cC_{\sigma'}.$
The middle term is the trivial element of 
$ \Pic( \cC_{\sigma} \cup \cC_{\sigma'} ) $,
the right-hand term is the trivial element
of $\Pic( \cC_{\sigma})$,
and the ideal sheaf is the element $(0,-1)$ in $\Pic(\cC_{\sigma'})$
(as one sees by considering its restriction to the members of the
open cover~\eqref{eqn:open cover  of two components}: its restriction
to $D(u)$ is the zero ideal, its restriction to $D(v) \cup D(u')$
is the ideal sheaf spanned by $v'$, which is an invertible sheaf
on $D(u')$ on which the $(\Gm)_{b'}$-action has been twisted by
weight~$-1$, and its restriction to $D(v')$ is the structure sheaf).
Thus we may rewrite this short exact sequence as 
\numequation
\label{eqn:s.e.s. from two components}
0 \to \cO_{\cC_{\sigma'}} (0,-1) 
\buildrel v' \cdot \over \longrightarrow
\cO_{\cC_{\sigma} \cup \cC_{\sigma'} } \to
\cO_{\cC_{\sigma}} \to 0.
\end{equation}
(Note that $v'$ is a well-defined global section of $\cO_{\cC_{\sigma} \cup
\cC_{\sigma'}}(0,-1,n)$ for any~$n$ --- and that
$\cO_{\cC_{\sigma} \cup \cC_{\sigma'}}(0,-1,n)$ restricts
to $\cO_{\cC_{\sigma'}}(0,1)$ on~$\cC_{\sigma'}$
--- so that  multiplication by $v'$
is a well-defined morphism 
$\cO_{\cC_{\sigma'}} (0,-1) \to \cO_{\cC_{\sigma} \cup \cC_{\sigma'} } .$)
\end{example}

\begin{example}
We also have the short exact sequence
$$0 \to \cI_{\cC_{\sigma'}} \to \cO_{\cC_{\sigma} \cup \cC_{\sigma'}}  
\to \cO_{\cC_{\sigma'}} \to 0,$$
which can be rewritten as the short exact sequence
\numequation
\label{eqn:other s.e.s. from two components}
0 \to \cO_{\cC_{\sigma}} (-1,0) 
\buildrel u \cdot \over \longrightarrow
\cO_{\cC_{\sigma} \cup \cC_{\sigma'} } \to
\cO_{\cC_{\sigma'}} \to 0.
\end{equation}
\end{example}

\begin{example}
The element $x= y'$ is a well-defined global section
of $\cO_{\cC_{\sigma} \cup \cC_{\sigma'}}(0,1,0)$.
Thus multiplication by $y'$, followed by restriction to~$\cC_{\sigma'}$,
induces a morphism
\numequation
\label{eqn:mult. by y'}
\cO_{\cC_{\sigma} \cup  \cC_{\sigma'}} \buildrel y' \cdot \over \longrightarrow
\cO_{\cC_{\sigma'}}(0,-1).
\end{equation}
The composite of this morphism with the first  non-trivial arrow of~\eqref{eqn:s.e.s.
from two components} (i.e.\ with the ``multiplication by $v'$'' map) gives
a morphism
\numequation
\label{eqn:mult. by v'y'}
\cO_{\cC_{\sigma'}}(0,-1) \buildrel (v'y') \cdot \over \longrightarrow
\cO_{\cC_{\sigma'}}(0,-1).
\end{equation}
\end{example}

By assumption,   $\sigma$ and~$\sigma'$  admit a nonzero extension
as~$\GL_2(\Fq)$-representations, and in fact the corresponding
$\Ext^1$ is one-dimensional. We let 
 $U$ denote the corresponding non-split extension of $\sigma$ by $\sigma'$, %
and we define $\gA(\cInd_{KZ}^G U) := \cO_{\cC_{\sigma}\cup \cC_{\sigma'}}(-1,-1,0).$
The short exact sequence
$$0 \to \cInd_{KZ}^G \sigma' \to \cInd_{KZ}^G U \to \cInd_{KZ}^G \sigma \to 0$$
corresponds to the short exact sequence
\numequation
\label{eqn:twisted s.e.s. from two components}
0 \to \cO_{\cC_{\sigma'}} (-1,0) 
\buildrel v' \cdot \over \longrightarrow
\cO_{\cC_{\sigma} \cup \cC_{\sigma'} }(-1,-1,0) \to
\cO_{\cC_{\sigma}}(-1,0) \to 0,
\end{equation}
which is a twist  of~\eqref{eqn:s.e.s. from two components}.

There is a well-known morphism\footnote{We may choose a character $\xi$ so
that $\sigma'$ is the socle of $\PS(\xi)$ and such that the  extension
$U$ embeds as a $\GL_2(\F_q)$-subrepresentation of~$\PS(\xi)$.  
Then $\sigma'$ is the cosocle of $\PS(\xi^s)$,
and so we may form the composite morphism
$\cInd_{KZ}^G U\hookrightarrow \cInd_{KZ}^G \PS(\xi) \buildrel \text{\eqref{eqn:PS intertwiner}}
\over \iso
\cInd_{KZ}^G \PS(\xi^s)  \to \cInd_{KZ}^G \sigma'.$}
$\cInd_{KZ}^G U \to \cInd_{KZ}^G \sigma'$ (which is {\em not}  the compact induction
of a morphism $U \to \sigma'$; indeed, there is no non-zero such morphism)
with the property that the composite of this morphism
with the obvious inclusion $\cInd_{KZ}^G \sigma' \hookrightarrow \cInd_{KZ}^G U$
is equal to the endomorphism $T$ of $\cInd_{KZ}^G \sigma'$.
This well-known morphism corresponds under $\mathfrak A$
to the twist
$$\cO_{\cC_{\sigma} \cup \cC_{\sigma'}} (-1,-1,0) \buildrel y'\cdot \over
\longrightarrow \cO_{\cC_{\sigma'}}(-1,0)$$
of~\eqref{eqn:mult. by y'}. 
Composing this with the first non-trivial arrow in~\eqref{eqn:twisted s.e.s.
from two components} 
then gives the twist
$$
\cO_{\cC_{\sigma'}}(-1,0) \buildrel (v'y') \cdot \over \longrightarrow
\cO_{\cC_{\sigma'}}(-1,0)
$$
of~\eqref{eqn:mult. by v'y'}.
This does indeed correspond to the  endomorphism $T$ of $\cInd_{KZ}^G \sigma'$,
and so we have verified an instance of the functoriality of~$\gA$.

Similarly, if we let $V$ be the non-split extension of $\sigma'$ by $\sigma$,
then we define $\gA(\cInd_{KZ}^G V) := 
\cO_{\cC_{\sigma} \cup \cC_{\sigma'} }(-1,0,0).$ 
Then the short exact sequence
$$0 \to \cInd_{KZ}^G \sigma \to \cInd_{KZ}^G V \to \cInd_{KZ}^G \sigma' \to 0$$
corresponds to the short exact sequence
$$
0 \to \cO_{\cC_{\sigma}} (-1,0) 
\buildrel u \cdot \over \longrightarrow
\cO_{\cC_{\sigma} \cup \cC_{\sigma'} }(-1,0,0) \to
\cO_{\cC_{\sigma'}}(-1,0) \to 0
$$
obtained by twisting~\eqref{eqn:other s.e.s. from two components}.

\subsubsection{Union of four components --- case one (principal series)}
\label{subsubsec:four components PS case}
Typically (i.e.\ for the induction of a suitably generic $k^{\times}$-valued
character),
the parabolic induction $\PS(\xi)$ of a character $\xi$ of the torus in $\GL_2(\F_q)$
can be denoted schematically as follows:
$$\xymatrix{ & \sigma_3 \ar@{-}[dl]\ar@{-}[dr]& \\ \sigma_2 \ar@{-}[dr] & &
\sigma_0 \ar@{-}[dl] \\ & \sigma_1 & }
$$
Here the various $\sigma_i$ are the Jordan--H\"older factors of~$\PS(\xi)$,
the rows of the diagram are the layers of its socle filtration (starting
from the bottom of the diagram) or equivalently of the cosocle filtration (starting
from the top of the diagram), and the lines indicate non-trivial extensions.

We write $\cC_i := \cC_{\sigma_i};$ our aim now is to describe
$\cC_0 \cup \cC_1 \cup \cC_2 \cup \cC_3,$
and then to describe its Picard group.
We begin by noting that this union admits the open cover
\nummultline
\label{eqn:PS open cover}
\{ D(u_0x_0), D(v_1y_1), D(u_2 x_2), D(v_3 y_3),
D(v_0 x_0) \cup D(u_1 y_1),
\\
D(v_1 x_1 ) \cup D(u_2 y_2),
D(v_2  x_2)  \cup D(u_3 y_3), D(v_3 x_3) \cup D(u_0 y_0),
\\
D(v_0y_0) \cup D(u_1 x_1) \cup D(v_2 y_2) \cup D(u_3 x_3)  \}.
\end{multline}
(Here and below we use the notation $u,v,x,y,a,b$ as before,
but add a subscript $i \in \set{0,1,2,3}$
to indicate the component to which the variable pertains.)
Here are explicit descriptions of these various open subsets:
$$
D(u_i x_i) = (\Gm)_{u_i/x_i} \times [\A^2_{v_i,y_i}/(\Gm)_{b_i}];$$
$$D(v_iy_i) = [\A^2_{u_i,x_i} /(\Gm)_{a_i}]\times (\Gm)_{v_i y_i};$$ 
$$D(v_i x_i) = \A^1_{u_i/x_i} \times \A^1_{v_i y_i},$$
$$D(u_{i+1} y_{i+1}) = \A^1_{x_{i+1}/u_{i+1}} \times \A^1_{v_{i+1} y_{i+1}},$$
and for $i = 0,2$, we have
$$D(v_i x_i) \cup D(u_{i+1} y_{i+1}) =
\set{ (u_i/x_i, v_i  y_i, v_{i+1} y_{i+1} ) \, | \, (u_i/x_i) (v_{i+1}  y_{i+1}) = 0 },$$
with the identification
$v_i y_i = x_{i+1}/u_{i+1}$,
while for $i = 1, 3$ we have
$$D(v_i x_i) \cup D(u_{i+1} y_{i+1}) =
\set{ (u_i/x_i, v_i  y_i, x_{i+1}/u_{i+1}) \, | \, (v_i y_i) (x_{i+1}/ u_{i+1}) = 0 },$$
with the identification
$u_i/x_i =
v_{i+1} y_{i+1} $;
and finally
\begin{multline*}
D(v_0y_0) \cup D(u_1 x_1) \cup D(v_2 y_2) \cup D(u_3 x_3) 
\\
=
[ \set{ (u_0, x_0, v_1,y_3, v_0 y_0) \, | \, u_0 v_1 = x_0 y_3 = 0, v_0 y_0 \neq 0 }
/ (\Gm)_{a_0}],
\end{multline*}
where we have the identifications
\begin{multline*}
x_0 = y_1, x_2 = y_3, u_0 = v_3, u_2 = v_1,
\\
v_0  y_0 = x_1/u_1 = (v_2 y_2)^{-1}  = u_3/x_3,
\\
a_0 = b_1^{-1} = a_2^{-1} = b_3.
\end{multline*}
Using the open cover~\eqref{eqn:PS open cover}, we compute
\numequation
\label{eqn:PS global sections}
\Gamma(\cO_{\cC_0 \cup \cC_1 \cup \cC_2 \cup \cC_3})  
= k[ v_1y_1,  v_3y_3]/( v_1y_1 v_3  y_3 ) .
\end{equation}
(To give slightly more detail: the function $v_1y_1$ is  well-defined
on~$\cC_1$, and extends by zero over the union of all four components.  Similarly,
the function $v_3y_3$ is well-defined on~$\cC_3$, and again extends by zero.)

One can show 
that there is an embedding
\numequation
\label{eqn:PS Pic embedding}
\Pic( \cC_0 \cup \cC_1 \cup \cC_2 \cup \cC_3) \hookrightarrow
\Z \times \Z \times \Z \times \Z
\times \Z,
\end{equation}
where the factors are given by
$
\Pic\bigl( D(v_1y_1) \bigr)$,
$\Pic\bigl( D(v_3y_3) \bigr),$
$\Pic\bigl(
D(v_0y_0) \cup D(u_1 x_1) \cup D(v_2 y_2) \cup D(u_3 x_3) 
\bigr),$
$ \Pic\bigl( D(u_0x_0) \bigr)$,
and
$\Pic\bigl( D(u_2x_2) \bigr),$
respectively,
with the morphisms being defined by restriction.  
The image of this embedding is equal to
\numequation
\label{eqn:PS Pic image}
\set{ (j,k,l,m,n) \, | \,  j - k + m - n = 0};
\end{equation}
in particular, 
$\Pic( \cC_0 \cup \cC_1 \cup \cC_2 \cup \cC_3) \cong \Z^4.$

The simplest way to fix our choice of coordinates precisely is to describe the
corresponding embedding
\numequation
\label{eqn:PS Picard}
\Pic( \cC_0 \cup \cC_1 \cup \cC_2 \cup \cC_3) \hookrightarrow
\Pic(\cC_0) \times
\Pic(\cC_1) \times
\Pic(\cC_2) \times
\Pic(\cC_3);
\end{equation}
it is given by
\numequation
\label{eqn:PS Picard map}
(j,k,l,m,n) \mapsto
 \bigl( (l,m) , (j,-l), (-l,n), (k,l) \bigr).
\end{equation}

\begin{example}
\label{ex:PS one}
We have the  following  diagram of ideal subsheaves of
$\cO_{\cC_0 \cup \cC_1 \cup \cC_2 \cup \cC_3}$: 
$$\xymatrix{
&\cO_{\cC_0 \cup \cC_1 \cup \cC_2 \cup \cC_3} & && \\
& \cI_{\cC_3}
\ar@{=}[r]\ar@{^{(}->}[u]^-{\cO_{\cC_3}} &
\cI_{\cC_2 \cup \cC_3} + \cI_{\cC_0 \cup \cC_3}  &&\\
\cI_{\cC_2 \cup \cC_3}  
\ar@{^{(}->}[ur]^-{\cO_{\cC_2}(-1,0)} &
 &  \cI_{\cC_0 \cup \cC_3}  
\ar@{_{(}->}[ul]_-{\cO_{\cC_0}(-1,0)} && \\
&
\cI_{\cC_0 \cup \cC_2 \cup \cC_3} \ar@{=}[r]
\ar@{_{(}->}[ul]
\ar@{^{(}->}[ur]
& \cI_{\cC_2 \cup \cC_3}
\cap \cI_{\cC_0  \cup \cC_3} 
\ar@{=}[r] &  v_1y_1 
\cO_{\cC_0 \cup \cC_1 \cup \cC_2 \cup \cC_3} & \ar[l]_-{\genfrac{}{}{0pt}{}{\scriptstyle v_1 y_1 \cdot}{\sim}}
 \cO_{\cC_1}  }
$$
where the sheaf labels on the arrows are indicate the corresponding cokernel.
Thus we may describe
$\cO_{\cC_0 \cup \cC_1 \cup \cC_2 \cup \cC_3}$
via the following extension diagram:
\numequation
\label{eqn:four components first extension diagram}
\xymatrix{ & \cO_{\cC_3} \ar@{-}[dl]\ar@{-}[dr]& \\ \cO_{\cC_2}(-1,0) \ar@{-}[dr] & &
\cO_{\cC_0}(-1,0) \ar@{-}[dl] \\ & \cO_{\cC_1} & }
\end{equation}
Similarly, we have the diagram 
$$\xymatrix{
&\cO_{\cC_0 \cup \cC_1 \cup \cC_2 \cup \cC_3} & && \\
& \cI_{\cC_1}
\ar@{=}[r]\ar@{^{(}->}[u]^-{\cO_{\cC_1}} &
\cI_{\cC_0 \cup \cC_1} + \cI_{\cC_1 \cup \cC_2}  &&\\
   \cI_{\cC_0 \cup \cC_1}  
\ar@{^{(}->}[ur]_-{\cO_{\cC_0}(-1,0)}
&&
\cI_{\cC_1 \cup \cC_2}  
\ar@{_{(}->}[ul]^-{\cO_{\cC_2}(-1,0)} \\
&
\cI_{\cC_0 \cup \cC_1 \cup \cC_2} \ar@{=}[r]
\ar@{_{(}->}[ul]
\ar@{^{(}->}[ur]
& \cI_{\cC_0 \cup \cC_1}
\cap \cI_{\cC_1  \cup \cC_2} 
\ar@{=}[r] &  v_3y_3 
\cO_{\cC_0 \cup \cC_1 \cup \cC_2 \cup \cC_3} & \ar[l]_-{\genfrac{}{}{0pt}{}{\scriptstyle v_3 y_3 \cdot}{\sim}}
 \cO_{\cC_3}  }
$$
which allows us to describe
$\cO_{\cC_0 \cup \cC_1 \cup \cC_2 \cup \cC_3}$
via the extension diagram
\numequation
\label{eqn:four components second extension diagram}
\xymatrix{ & \cO_{\cC_1} \ar@{-}[dl]\ar@{-}[dr]& \\ \cO_{\cC_0}(-1,0) \ar@{-}[dr] & &
\cO_{\cC_2}(-1,0) \ar@{-}[dl] \\ & \cO_{\cC_3} & }
\end{equation}
\end{example}

\begin{example}
\label{ex:PS two}
We also have the  following  diagrams of ideal subsheaves of
$\cO_{\cC_0 \cup \cC_1 \cup \cC_2 \cup \cC_3}$,
namely 
$$\xymatrix{
&\cO_{\cC_0 \cup \cC_1 \cup \cC_2 \cup \cC_3} & && \\
& \cI_{\cC_0}
\ar@{=}[r]\ar@{^{(}->}[u]^-{\cO_{\cC_0}} &
\cI_{\cC_0 \cup \cC_3} + \cI_{\cC_0 \cup \cC_1}  &&\\
\cI_{\cC_0 \cup \cC_3}  
\ar@{^{(}->}[ur]^-{\cO_{\cC_3}(0,1)} &
 &  \cI_{\cC_0 \cup \cC_3}  
\ar@{_{(}->}[ul]_-{\cO_{\cC_1}(0,-1)} && \\
&
\cI_{\cC_0 \cup \cC_1 \cup \cC_3} \ar@{=}[r]
\ar@{_{(}->}[ul]
\ar@{^{(}->}[ur]
& \cI_{\cC_0 \cup \cC_3}
\cap \cI_{\cC_0  \cup \cC_1} 
&
\ar[l]_-{\sim} \cO_{\cC_2}(-2,0)}
$$
and
$$\xymatrix{
&\cO_{\cC_0 \cup \cC_1 \cup \cC_2 \cup \cC_3} & && \\
& \cI_{\cC_2}
\ar@{=}[r]\ar@{^{(}->}[u]^-{\cO_{\cC_2}} &
\cI_{\cC_1 \cup \cC_2} + \cI_{\cC_2 \cup \cC_3}  &&\\
\cI_{\cC_1 \cup \cC_2}  
\ar@{^{(}->}[ur]^-{\cO_{\cC_1}(0,1)} &
 &  \cI_{\cC_2 \cup \cC_3}  
\ar@{_{(}->}[ul]_-{\cO_{\cC_3}(0,-1)} && \\
&
\cI_{\cC_1 \cup \cC_2 \cup \cC_3} \ar@{=}[r]
\ar@{_{(}->}[ul]
\ar@{^{(}->}[ur]
& \cI_{\cC_1 \cup \cC_2}
\cap \cI_{\cC_2  \cup \cC_3} 
&
\ar[l]_-{\sim}   \cO_{\cC_0}(-2,0)}
$$
which allows us to describe
$\cO_{\cC_0 \cup \cC_1 \cup \cC_2 \cup \cC_3}$
via the extension diagrams
\numequation
\label{eqn:four components third extension diagram}
\xymatrix{ & \cO_{\cC_0} \ar@{-}[dl]\ar@{-}[dr]& \\ \cO_{\cC_3}(0,1) \ar@{-}[dr] & &
\cO_{\cC_1}(0,-1) \ar@{-}[dl] \\ & \cO_{\cC_2}(-2,0) & }
\end{equation}
and
\numequation
\label{eqn:four components fourth extension diagram}
\xymatrix{ & \cO_{\cC_2} \ar@{-}[dl]\ar@{-}[dr]& \\ \cO_{\cC_1}(0,1) \ar@{-}[dr] & &
\cO_{\cC_3}(0,-1) \ar@{-}[dl] \\ & \cO_{\cC_0}(-2,0) & }
\end{equation}
\end{example}

We now explain the representation-theoretic interpretation of the
above calculations.
We write $U:= \PS(\xi)$, and recall
the description
\numequation
\label{eqn:PS first extension diagram}
\xymatrix{ & \sigma_3 \ar@{-}[dl]\ar@{-}[dr]& \\ \sigma_2 \ar@{-}[dr] & &
\sigma_0 \ar@{-}[dl] \\ & \sigma_1 & }
\end{equation}
of $U$ in terms of its Jordan--H\"older factors.
We then define
$$\gA(\cInd_{KZ}^G U) := \cO_{\cC_0 \cup \cC_1 \cup \cC_2 \cup \cC_3}(-1,-1,0,0,0),$$
where the twisting is understood to
be in terms of the embedding~\eqref{eqn:PS Pic embedding};
note that $-1 - (-1) + 0 = 0 = 0,$ so that
the constraint of~\eqref{eqn:PS Pic image} is satisfied.
If we form the corresponding twist of the
diagram~\eqref{eqn:four components first extension diagram},
and take into account the formula of~\eqref{eqn:PS Picard map},
we obtain the diagram
\numequation
\label{eqn:four components first extension diagram twisted}
\xymatrix{ & \cO_{\cC_3}(-1,0) \ar@{-}[dl]\ar@{-}[dr]& \\ \cO_{\cC_2}(-1,0) \ar@{-}[dr] & &
\cO_{\cC_0}(-1,0) \ar@{-}[dl] \\ & \cO_{\cC_1}(-1,0) & }
\end{equation}
which shows that  our definition is compatible with the functoriality
of $\gA$ and diagram~\eqref{eqn:PS first extension diagram}.

If we induce the character $\xi^s$ (the non-trivial Weyl group conjugate of~$\xi$)
instead, we obtain a representation %
$U^s := \PS(\xi^s)$, which can be described
by the extension diagram
\numequation
\label{eqn:PS second extension diagram}
\xymatrix{ & \sigma_1 \ar@{-}[dl]\ar@{-}[dr]& \\ \sigma_0 \ar@{-}[dr] & &
\sigma_2 \ar@{-}[dl] \\ & \sigma_3 & }
\end{equation}
Although $U$ and $U^s$ are {\em not} isomorphic representations of $\GL_2(\F_{p^2})$,
their compact inductions to $\GL_2(\Q_{p^2})$  {\em are}
isomorphic (via the isomorphism~\eqref{eqn:PS intertwiner}),  and so we also have that
$\gA(\cInd_{KZ}^G U^s) = \cO_{\cC_0 \cup \cC_1 \cup \cC_2 \cup \cC_3}(-1,-1,0,0,0)$.
Note that if we twist
diagram~\eqref{eqn:four components second extension diagram}
by~$(-1,-1,0,0,0)$, we obtain the diagram
\numequation
\label{eqn:four components second extension diagram twisted}
\xymatrix{ & \cO_{\cC_1}(-1,0) \ar@{-}[dl]\ar@{-}[dr]& \\ \cO_{\cC_0}(-1,0) \ar@{-}[dr] & &
\cO_{\cC_2}(-1,0) \ar@{-}[dl] \\ & \cO_{\cC_3}(-1,0) & }
\end{equation}
This verifies  the functoriality of $\gA$ when applied to
the diagram~\eqref{eqn:PS second extension diagram}.

The composite $\cInd_{KZ}^G \sigma_1 \hookrightarrow \cInd_{KZ}^G U \iso
\cInd_{KZ}^G U^s \to \cInd_{KZ}^G \sigma_1$ is known to coincide with the Hecke operator~$T$.
A consideration of the diagrams in Example~\ref{ex:PS one}
shows that this map is given by multiplication by~$v_1y_1$,
which is indeed $\gA(T)$. %
Related to this computation, we note that~\eqref{eqn:PS global sections}
is compatible with the full faithfulness of~$\gA$
and the known structure of the endomorphism
ring of $\cInd_{KZ}^G U$.

As is well-known, by reducing lattices in integral principal series,
we can also obtain representations $V$ and $W$ of $\GL_2(\F_q)$, which
can be described by the extension diagrams
\numequation
\label{eqn:PS third extension diagram}
\xymatrix{ & \sigma_0 \ar@{-}[dl]\ar@{-}[dr]& \\ \sigma_3 \ar@{-}[dr] & &
\sigma_1 \ar@{-}[dl] \\ & \sigma_2 & }
\end{equation}
and
\numequation
\label{eqn:PS fourth extension diagram}
\xymatrix{ & \sigma_2 \ar@{-}[dl]\ar@{-}[dr]& \\ \sigma_1 \ar@{-}[dr] & &
\sigma_3 \ar@{-}[dl] \\ & \sigma_0 & }
\end{equation}
respectively.

We define
$$\gA(\cInd_{KZ}^G V) :=  \cO_{\cC_0 \cup \cC_1 \cup \cC_2 \cup \cC_3}(-1,-1,-1,0,0)$$
and
$$\gA(\cInd_{KZ}^G W) :=  \cO_{\cC_0 \cup \cC_1 \cup \cC_2 \cup \cC_3}(-1,-1,1,0,0).$$
Then twisting the diagrams~\eqref{eqn:four components third extension diagram}
and~\eqref{eqn:four components fourth extension diagram}
by $(-1,-1,-1,0,0)$ and $(-1,-1,1,0,0)$ respectively, we obtain diagrams
\numequation
\label{eqn:four components third extension diagram twisted}
\xymatrix{ & \cO_{\cC_0}(-1,0) \ar@{-}[dl]\ar@{-}[dr]& \\ \cO_{\cC_3}(-1,0) \ar@{-}[dr] & &
\cO_{\cC_1}(-1,0) \ar@{-}[dl] \\ & \cO_{\cC_2}(-1,0) & }
\end{equation}
and
\numequation
\label{eqn:four components fourth extension diagram twisted}
\xymatrix{ & \cO_{\cC_2}(-1,0) \ar@{-}[dl]\ar@{-}[dr]& \\ \cO_{\cC_1}(-1,0) \ar@{-}[dr] & &
\cO_{\cC_3}(-1,0) \ar@{-}[dl] \\ & \cO_{\cC_0}(-1,0) & }
\end{equation}
confirming the functoriality of $\gA$ when applied
to~\eqref{eqn:PS third extension diagram} and~\eqref{eqn:PS fourth extension diagram}.

\subsubsection{Union of four components --- case two (cuspidal)}
\label{subsubsec:four components cuspidal case}
The computations in this case follow similar lines to those in the
principal series case.  For a tame cuspidal type, there is no preferred lattice,
though.  Rather, the reduction (generically) contains four weights $\sigma_0$,  
$\sigma_1,$ $\sigma_2,$ and $\sigma_3$, each with multiplicity one,
and one can arrange for them to appear in the extension diagram
$$
\xymatrix{ & \sigma_0 \ar@{-}[dl]\ar@{-}[dr]& \\ \sigma_3 \ar@{-}[dr] & &
\sigma_1 \ar@{-}[dl] \\ & \sigma_2 & }
$$
or more generally in any cyclic  permutation of this diagram.

For an appropriate choice of labelling of the weights, the
union $\cC_{\sigma_0} \cup \cC_{\sigma_1} \cup \cC_{\sigma_2} \cup \cC_{\sigma_3}$
admits the open cover
\nummultline
\label{eqn:cuspidal open cover}
\{ D(v_0x_0), D(u_1y_1), D(v_2 x_2), D(u_3 y_3),
D(v_0 y_0) \cup D(u_1 x_1),
\\
D(v_1 y_1 ) \cup D(u_2 x_2),
D(v_2  y_2)  \cup D(u_3 x_3), D(v_3 y_3) \cup D(u_0 x_0),
\\
D(u_0y_0) \cup D(v_1 x_1) \cup D(u_2 y_2) \cup D(v_3 x_3)  \}.
\end{multline}
We have
$$D(u_iy_i) = \A^1_{x_i/u_i} \times \A^1_{v_i y_i};$$
$$D(v_i x_i) = \A^1_{u_i/x_i} \times \A^1_{v_i y_i};$$
$$D(u_i x_i) = (\Gm)_{x_i/u_i} \times [\A^2_{v_i,y_i} / (\Gm)_{b_i}],$$
$$D(v_i y_i) = [\A^2_{u_i,x_i}/(\Gm)_{a_i}] \times (\Gm)_{v_i y_i},$$
and for $i = 0,2$, we have
$$D(v_i y_i) \cup D(u_{i+1} x_{i+1}) = 
[\set{(u_i,x_i,v_i y_i,y_{i+1} ) \, | \,
v_i y_i \neq 0, x_i  y_{i+1}= 0 }/ (\Gm)_{a_i}],$$
with the identifications $u_i = v_{i+1}$, $v_i y_i = u_{i+1}/x_{i+1}$,
and $a_i = b_{i+1}$,
while for $i  = 1,3$ we have
$$D(v_i y_i) \cup D(u_{i+1} x_{i+1}) = 
[\set{(u_i,x_i,v_i y_i,v_{i+1} ) \, | \,
v_i y_i \neq 0, u_i  v_{i+1}= 0 }/ (\Gm)_{a_i}],$$
with the identifications $x_i = y_{i+1}$, $v_i y_i = x_{i+1}/u_{i+1}$,
and $a_i = b_{i+1}^{-1}$;
and finally
$$
D(u_0y_0) \cup D(v_1 x_1) \cup D(u_2 y_2) \cup D(v_3 x_3)  
= \set{( v_0 y_0 ,  v_1 y_1, v_2 y_2, v_3 y_3 ) \, | \,
(v_0 y_0) (v_2 y_2) = (v_1y_1) (v_3y_2) = 0 },$$
with the identifications
$x_0/u_0 = v_3 y_3$,
$u_1/x_1 = v_0y_0$,
$x_2/u_2 = v_1y_1$,
and
$u_3/x_3 = v_2 y_2$.

Using the  open cover~\eqref{eqn:cuspidal open cover},
we compute
\numequation
\label{eqn:cuspidal global sections}
\Gamma(\cO_{\cC_0 \cup \cC_1 \cup \cC_2 \cup \cC_3})  
= k.
\end{equation}

We also compute that
$$\Pic( \cC_0 \cup \cC_1 \cup \cC_2 \cup \cC_3) = \Z \times \Z \times \Z \times \Z,$$
where the factors are given by
$
\Pic\bigl( D(v_2y_2) \cup D(u_3 x_3)\bigr)$,
$\Pic\bigl( D(v_1y_1) \cup D(u_2 x_2)\bigr),$
$ \Pic\bigl( D(v_0y_0) \cup D(u_1 x_1)\bigr)$,
and
$\Pic\bigl( D(v_3 y_3) \cup D(u_0x_0) \bigr),$
respectively.
We fix our choice of coordinates by describing the embedding
\numequation
\label{eqn:cuspidal Picard}
\Pic( \cC_0 \cup \cC_1 \cup \cC_2 \cup \cC_3) \hookrightarrow
\Pic(\cC_0) \times
\Pic(\cC_1) \times
\Pic(\cC_2) \times
\Pic(\cC_3);
\end{equation}
it is given by
\numequation
\label{eqn:cuspidal Picard map}
(k,l,m,n) \mapsto
 \bigl( (m,n) , (-l,m), (k,l), (-n,k) \bigr).
\end{equation}

\begin{example}
\label{ex:cuspidal}
We have the  following  diagrams of ideal subsheaves of
$\cO_{\cC_0 \cup \cC_1 \cup \cC_2 \cup \cC_3}$: 
$$\xymatrix{
&\cO_{\cC_0 \cup \cC_1 \cup \cC_2 \cup \cC_3} & && \\
& \cI_{\cC_0}
\ar@{=}[r]\ar@{^{(}->}[u]^-{\cO_{\cC_0}} &
\cI_{\cC_0 \cup \cC_3} + \cI_{\cC_0 \cup \cC_1}  &&\\
\cI_{\cC_0 \cup \cC_3}  
\ar@{^{(}->}[ur]^-{\cO_{\cC_3}(-1,0)} &
 &  \cI_{\cC_0 \cup \cC_3}  
\ar@{_{(}->}[ul]_-{\cO_{\cC_1}(0,1)} && \\
&
\cI_{\cC_0 \cup \cC_1 \cup \cC_3} \ar@{=}[r]
\ar@{_{(}->}[ul]
\ar@{^{(}->}[ur]
& \cI_{\cC_0 \cup \cC_3}
\cap \cI_{\cC_0  \cup \cC_1} 
&
\ar[l]_-{\sim} \cO_{\cC_2}(-1,-1)}
$$
\medskip
$$\xymatrix{
&\cO_{\cC_0 \cup \cC_1 \cup \cC_2 \cup \cC_3} & && \\
& \cI_{\cC_1}
\ar@{=}[r]\ar@{^{(}->}[u]^-{\cO_{\cC_1}} &
\cI_{\cC_0 \cup \cC_1} + \cI_{\cC_1 \cup \cC_2}  &&\\
   \cI_{\cC_0 \cup \cC_1}  
\ar@{^{(}->}[ur]_-{\cO_{\cC_0}(-1,0)}
&&
\cI_{\cC_1 \cup \cC_2}  
\ar@{_{(}->}[ul]^-{\cO_{\cC_2}(0,-1)} \\
&
\cI_{\cC_0 \cup \cC_1 \cup \cC_2} \ar@{=}[r]
\ar@{_{(}->}[ul]
\ar@{^{(}->}[ur]
& \cI_{\cC_0 \cup \cC_1}
\cap \cI_{\cC_1  \cup \cC_2} 
& \ar[l]_-{\sim}
 \cO_{\cC_3}(-1,1)  }
$$
\medskip
$$\xymatrix{
&\cO_{\cC_0 \cup \cC_1 \cup \cC_2 \cup \cC_3} & && \\
& \cI_{\cC_2}
\ar@{=}[r]\ar@{^{(}->}[u]^-{\cO_{\cC_2}} &
\cI_{\cC_1 \cup \cC_2} + \cI_{\cC_2 \cup \cC_3}  &&\\
\cI_{\cC_1 \cup \cC_2}  
\ar@{^{(}->}[ur]^-{\cO_{\cC_1}(-1,0)} &
 &  \cI_{\cC_2 \cup \cC_3}  
\ar@{_{(}->}[ul]_-{\cO_{\cC_3}(0,1)} && \\
&
\cI_{\cC_1 \cup \cC_2 \cup \cC_3} \ar@{=}[r]
\ar@{_{(}->}[ul]
\ar@{^{(}->}[ur]
& \cI_{\cC_1 \cup \cC_2}
\cap \cI_{\cC_2  \cup \cC_3} 
&
\ar[l]_-{\sim}   \cO_{\cC_0}(-1,-1)}
$$
\medskip
$$\xymatrix{
&\cO_{\cC_0 \cup \cC_1 \cup \cC_2 \cup \cC_3} & && \\
& \cI_{\cC_3}
\ar@{=}[r]\ar@{^{(}->}[u]^-{\cO_{\cC_3}} &
\cI_{\cC_2 \cup \cC_3} + \cI_{\cC_0 \cup \cC_3}  &&\\
\cI_{\cC_2 \cup \cC_3}  
\ar@{^{(}->}[ur]^-{\cO_{\cC_2}(-1,0)} &
 &  \cI_{\cC_0 \cup \cC_3}  
\ar@{_{(}->}[ul]_-{\cO_{\cC_0}(0,-1)} && \\
&
\cI_{\cC_0 \cup \cC_2 \cup \cC_3} \ar@{=}[r]
\ar@{_{(}->}[ul]
\ar@{^{(}->}[ur]
& \cI_{\cC_2 \cup \cC_3}
\cap \cI_{\cC_0  \cup \cC_3} 
& \ar[l]_-{\sim} \cO_{\cC_1}(-1,1)  }
$$
\end{example}

The following representations (described via their extension diagrams)
arise as the reduction of lattices in cuspidal types:
\label{subsec:cusp}
\numequation
\label{eqn:first cuspidal extension diagram}
\xymatrix{ & \sigma_0 \ar@{-}[dl]\ar@{-}[dr]& \\ \sigma_3 \ar@{-}[dr] & &
\sigma_1 \ar@{-}[dl] \\ & \sigma_2 & }
\end{equation}
\numequation
\label{eqn:second cuspidal extension diagram}
\xymatrix{ & \sigma_1 \ar@{-}[dl]\ar@{-}[dr]& \\ \sigma_0 \ar@{-}[dr] & &
\sigma_2 \ar@{-}[dl] \\ & \sigma_3 & }
\end{equation}
\numequation
\label{eqn:third cuspidal extension diagram}
\xymatrix{ & \sigma_2 \ar@{-}[dl]\ar@{-}[dr]& \\ \sigma_1 \ar@{-}[dr] & &
\sigma_3 \ar@{-}[dl] \\ & \sigma_0 & }
\end{equation}
and
\numequation
\label{eqn:fourth cuspidal extension diagram}
\xymatrix{ & \sigma_3 \ar@{-}[dl]\ar@{-}[dr]& \\ \sigma_2 \ar@{-}[dr] & &
\sigma_0 \ar@{-}[dl] \\ & \sigma_1 & }
\end{equation}
The corresponding values of $\gA$ are 
$\cO_{\cC_0 \cup \cC_1 \cup \cC_2 \cup \cC_3} ( 0,1,-1,0)$,
$\cO_{\cC_0 \cup \cC_1 \cup \cC_2 \cup \cC_3} ( -1,1,0,0)$,
$\cO_{\cC_0 \cup \cC_1 \cup \cC_2 \cup \cC_3} ( -1,0,0,1)$,
and
$\cO_{\cC_0 \cup \cC_1 \cup \cC_2 \cup \cC_3} ( 0,0,-1,1)$.
If we twist each of the diagrams  in Example~\ref{ex:cuspidal}
by the corresponding twist (and compute using the formula~\eqref{eqn:cuspidal
Picard map}), 
we obtain an extension diagram constituted from the $\gA(\sigma_i)$ in the
appropriate order; this confirms the functoriality of $\gA$ on
the various diagrams
\eqref{eqn:first cuspidal  extension diagram},
\eqref{eqn:second cuspidal  extension diagram},
\eqref{eqn:third cuspidal  extension diagram},
and~\eqref{eqn:fourth cuspidal  extension diagram}
in turn.
 
\subsubsection{A non-compactly induced example}
\label{subsubsec:non CM structure sheaf example}
The isomorphism $\cInd_{KZ}^G \overline{\PS(\xi)} \iso \cInd_{KZ}^G
\overline{\PS(\xi^{s})}$
allows us to construct an interesting representation.  
(Here
$\overline{\PS(\xi)}$ is the representation denoted $U$
in Section~\ref{subsubsec:four components PS case}, 
and 
$\overline{\PS(\xi^s)}$ is the representation denoted $U^s$.)
If we let $\sigma_{\xi}$  denote the socle of $\overline{\PS(\xi)}$,
and $\sigma_{\xi^{s}}$ denote the socle of $\overline{\PS(\xi^{s)}}$,
then there are embeddings
$$\cInd_{KZ}^G \sigma_{\xi} \hookrightarrow
\cInd_{KZ}^G \overline{\PS(\xi)}$$
and
$$\cInd_{KZ}^G \sigma_{\xi^{s}} \hookrightarrow
\cInd_{KZ}^G \overline{\PS(\xi^{s})}\iso
\cInd_{KZ}^G \overline{\PS(\xi)},$$
inducing
$$
\cInd_{KZ}^G \sigma_{\xi} \oplus
\cInd_{KZ}^G \sigma_{\xi^{s}} \hookrightarrow
\cInd_{KZ}^G \overline{\PS(\xi)}.$$
Define $\Pi$ to be the cokernel of this embedding.

We can compute $\mathfrak{A}(\Pi)$ using
the calculations of Section~\ref{subsubsec:four components PS case}.  Indeed, noting 
that $\sigma_{\xi}$  and $\sigma_{\xi^{s}}$ can be taken to
be $\sigma_1$ and $\sigma_3$ in the notation of that discussion,
we find that  $\mathfrak{A}(\Pi)$ is the cokernel of 
the embedding
$$\cO_{\cC_1}(-1,0)\oplus \cO_{\cC_3}(-1,0) \hookrightarrow
\cO_{\cC_0 \cup \cC_1 \cup \cC_2 \cup \cC_3}(-1,-1,0,0,0).$$
The image of this embedding is identified with the 
twist  by~$(-1,-1,0,0,0)$ of the ideal sheaf $\cI_{\cC_0 \cup \cC_2},$
and one finds that $\mathfrak A(\Pi)$ 
is a certain invertible sheaf on  $\cC_0 \cup \cC_2$.
(In particular, it is {\em not} Cohen--Macaulay.)

\subsubsection{An example --- weight cycling}
In the 2006  AIM meeting, the participants found an argument which has come
to be known as ``weight  cycling''.   Here it is in its original
form. (For more details, the interested reader could consult the
contemporary reports
\url{https://aimath.org/WWN/padicmodularity/emertonreport.pdf} and \url{https://www.ma.ic.ac.uk/~buzzard/maths/research/notes/notes_on_mod_p_local_langlands.pdf}.)

\begin{prop}
\label{prop:weight cycling}
 Suppose that $\pi'$ is a representation of $\GL_2(\Q_{p^2})$
with the following properties:

\begin{enumerate} 
\item The constituents of the $K$-socle of $\pi'$ are among
the constituents  of the reduction of some sufficiently generic cuspidal
type~$\tau$.
\item Each constituent of the $K$-socle of $\pi'$ appears with multiplicity one.
\item Each constituent of the  $K$-socle of $\pi'$  is supersingular.
\end{enumerate}
Then there is a  {\em (}necessarily irreducible{\em )} subrepresentation $\pi$  of  $\pi'$  with
the following properties:
\begin{enumerate}
\item[(a)] $\pi$ admits a central character.
\item[(b)] Every $G$-subrepresentation of $\pi'$ contains $\pi$.
\item[(c)] 
The $K$-socle of  $\pi$ coincides
with the $K$-socle of $\pi'$.
\item[(d)]
All four constituents of the reduction of $\tau$ actually appear
in the common $K$-socle of $\pi'$  and $\pi$.
\end{enumerate}
\end{prop}
\begin{proof}
Any $G$-subrepresentation of $\pi'$ contains at least one of the constituents
of the $K$-socle  of $\pi'$, and thus contains the $G$-subrepresentation
that this constituent  generates.  Thus if we show
that the $G$-subrepresentation $\pi$  generated by any one  of these
constituents contains {\em all} of these constituents --- that 
is, if we show that  $\pi$ satisfies~(c) --- then we will also
have shown that $\pi$ is independent of the initial choice
of generating constituent, and thus that
$\pi$ also satisfies condition~(b).
(And the fact that $\pi$  satisfies~(b) implies that it is irreducible.)

So we turn to  proving that $\pi$ does indeed satisfy~(c).
Along the way, we will also show that it satisfies (a) and~(d).

Let $\sigma_0$ denote the chosen constituent of the reduction of $\tau$  that
appears in the $K$-socle of $\pi'$ and generates~$\pi$.
Because $\sigma_0$ appears in the $K$-socle of $\pi$ with multiplicity one,
the centre $Z$ of $G$ acts on $\sigma_0$, and hence
also on the representation $\pi$  that it generates, through a character. 
This proves~(a).

Recalling that  $\sigma_0$ is one of the weights appearing in the reduction
of~$\tau$,
we then label the remaining weights in the reduction of $\tau$ as
$\sigma_1,\sigma_2,\sigma_3$.
Since $\tau$ has a central character as a representation
of~$K$, each $\sigma_i$  has a central character  as a representation
of~$K$, and these characters all coincide (they are just the reduction
of the central character of~$\tau$).  In particular
the central character of $\pi$  is compatible with  the
central character of the $\sigma_i$  (i.e.\ the restriction of
the former to $K\cap Z$ coincides with the latter), and so we can and do extend each
$\sigma_i$  to a representation of $KZ$ by having $Z$ act through
the central character of~$\pi$.

We can then choose our labels so that
for each  $i = 0,\dots,3$, we have an embedding (originally
constructed by Pa\v{s}k\={u}nas)
\numequation
\label{eqn:embedding}
(\cInd_{KZ}^G \sigma_{i+1} \oplus \cInd_{KZ}^G \sigma_{i+1}'  )
\hookrightarrow \cInd_{KZ}^G \sigma_i  / T,
\end{equation}
where $i+1$  is to be computed mod~$4$, 
and $\sigma_{i+1}'$  is another Serre weight
which is not isomorphic to any of the~$\sigma_j$ ($j = 0,\dots,3$).
(See Remark~\ref{rem:weight facts} below for a slight elaboration
on this.)

Let $\pi'(\sigma_i)$  denote the  cokernel of~\eqref{eqn:embedding}.  
Then a basic fact is that  $\pi'(\sigma_i) = \pi'(\sigma_i^{(s)}),$
where $\sigma_i^{(s)}$ denotes the Serre weight
whose highest  weight
is obtained from that  of $\sigma_i$ by  applying the Weyl group involution~$s$.
In particular, $\pi'(\sigma_i)$ is not  only a quotient
of $\cInd_{KZ}^G \sigma_i$, but is also a quotient of~$\cInd_{KZ}^G  \sigma_i^{(s)}.$
Another basic fact is that none of the  $\sigma_i^{(s)}$ (for $i = 0, \dots,3$)
is isomorphic to any of the $\sigma_j$;  i.e.\ none of the former
representations are constituents of  the reduction of~$\tau$.
(Again, see Remark~\ref{rem:weight facts} below for an elaboration on this.)
By our assumption on the $K$-socle of $\pi'$,
this means that  any morphism
$\cInd_{KZ}^G  \sigma_i^{(s)} \to \pi'$ (for any value of~$i$) must vanish,
and consequently  any morphism $\pi'(\sigma_i) \to \pi'$ must vanish.
Since $\pi \subseteq \pi',$  the same applies with $\pi$ in place of~$\pi'$.

Now consider the  morphism $\cInd_{KZ}^G \sigma_0 \to  \pi$
induced by the  inclusion $\sigma_0 \subseteq  \pi$.  The supersingularity
assumption implies (indeed, more-or-less means) that this factors through
$\cInd_{KZ}^G \sigma_0/T.$  Since it is non-zero (it contains $\sigma_0$ in
its image) it must  not factor through $\pi'(\sigma_i)$, and so 
it must not vanish identically on %
the domain of~\eqref{eqn:embedding}.
However, it {\em must} vanish on the  second such summand (by our assumption
on the $K$-socle of $\pi'$), and hence it is non-zero on the first summand.
Thus we obtain a non-zero morphism
$\cInd_{KZ}^G \sigma_{i+1} \to \pi.$  Equivalently,
we obtain an embedding $\sigma_{i+1} \hookrightarrow \pi$.

Continuing  in  this manner, we cycle through all the weights
$\sigma_i$, and find that the $\pi'$ contains all of them in its $K$-socle,
and that in fact  they are all contained in~$\pi$.  This completes
the proof of the proposition.
\end{proof}

\begin{remark}
\label{rem:weight facts}
Just to say slightly more about the facts regarding Serre weights that
we used above:  if $\chi_i$  denotes the highest weight  of $\sigma_i$,
then the weights in the reduction of $\PS(\chi_i)$  are $\sigma_i, \sigma_{i+1},
\sigma_{i+1}',$  and $\sigma_i^{(s)}$.  And the reduction of a (sufficiently
generic) cuspidal type
and of a (sufficiently generic) principal series type have either no or
two weights  in common.  So, since $\tau$  contains $\sigma_i$ and $\sigma_{i+1}$
in its reduction, it  does not contain $\sigma_{i+1}'$  or $\sigma_i^{(s)}$.
\end{remark}

\begin{rem}\label{rem: BP}The original aim of the participants at AIM
  in 2006 had been to apply Proposition~\ref{prop:weight
    cycling} with~$\pi'$ equal to the $\m$-torsion in the completed
  cohomology of a Shimura curves (where~$\m$
  denotes a maximal ideal of an appropriate Hecke algebra). In this
  case the hypotheses of the proposition were known to hold under
  appropriate genericity hypotheses by~\cite{geebdj}. Much more
  recently it has been proved that in this case we even
  have~$\pi'=\pi$ (see~\cite[Cor.\ 1.3.10]{https://doi.org/10.48550/arxiv.2102.06188}). %

  Shortly after the 2006 AIM meeting,
  Breuil and Pa\v{s}k\={u}nas~\cite{BreuilPaskunas} gave a direct
  construction of
  an irreducible admissible representation $\pi$ as in Proposition~\ref{prop:weight
    cycling}. %
  In
  fact, they exhibited infinitely many pairwise non-isomorphic such
  representations.
\end{rem}

The proof of Proposition~\ref{prop:weight cycling} shows that
there is actually a scalar  invariant attached to $\pi'$ --- or, really,
to its irreducible  subrepresentation~$\pi$.
Namely, if we begin with  the  inclusion $\sigma_0 \subseteq \pi$
and cycle four times, we return to an embedding
$\sigma_0  \hookrightarrow \pi$, which must then (by our multiplicity
one assumption) be equal to a non-zero scalar multiple of the original
embedding.    In order to make this well-defined, one has to pin down
the embeddings~\eqref{eqn:embedding}, but with a bit of care
one can do that (we explain below how to do this after applying~$\fA$),
and one obtains a scalar invariant of $\pi$.
More precisely, we choose a non-zero scalar $\lambda = \lambda_0
\in k^{\times}$,
write $\lambda_i  = 1 $ if $i  =  1,2,3,$,
and consider the two step complex
\numequation
\label{eqn:diagram}
\bigoplus_{i = 0}^3 \bigl( \cInd_{KZ}^G  \sigma_i \oplus 
\cInd_{KZ}^G \sigma_i' \bigr) 
\to
\bigoplus_{i=0}^3 \cInd_{KZ}^G  \sigma_i /T ,
\end{equation}
where the arrow is defined
on the summand 
$\cInd_{KZ}^G \sigma_{i} \oplus \cInd_{KZ}^G \sigma_{i}'$ via
\begin{multline*}
 \cInd_{KZ}^G \sigma_i \oplus \cInd_{KZ} \sigma_i' \ni (v,w)
\\  \mapsto 
( v + w, - \lambda_i v \bmod T) \in  
\cInd_{KZ}^G  \sigma_{i-1}/T
 \oplus  
\cInd_{KZ}^G  \sigma_i/T;
\end{multline*}
here we regard $v$  and $w$ as elements of $\cInd_{KZ} \sigma_{i-1}/T$
via~\eqref{eqn:embedding}.

\begin{lemma}
\label{lem:injectivity}
The arrow of~{\em \eqref{eqn:diagram}} is injective.
\end{lemma}
\begin{proof}
If the  kernel is non-zero, then it has a non-zero $K$-socle,
and so there is a non-zero morphism from $\cInd_{KZ}^G \sigma'$ 
to this  kernel, for some weight~$\sigma'$.
But if $\cInd_{KZ}^G \sigma' \to 
\bigoplus_{i = 0}^3 \bigl( \cInd_{KZ}^G  \sigma_i \oplus 
\cInd_{KZ}^G \sigma_i' \bigr) 
$
is  non-zero, then $\sigma'$ is  equal to one of the $\sigma_i$ or
$\sigma_i'$ (we are assuming that these weights are sufficiently generic),
and  the morphism is just an embedding into the corresponding summand.
The lemma then follows from the fact that the arrow
of~\eqref{eqn:diagram} is injective on each summand in its domain.
\end{proof}

Let $\Pi_{\lambda}$ denote the cokernel of~\eqref{eqn:diagram}.  
Then giving a morphism $\Pi_\lambda \to \pi$, for some other representation~$\pi$,
is equivalent to giving a copy of $\sigma_0$ (or indeed, any one of
the $\sigma_i$) in $\pi,$  such that the  weight cycling of
Proposition~\ref{prop:weight cycling} ``works'' in~$\pi$, and such
that the corresponding
scalar %
is equal to~$\lambda$.

We can compute the coherent sheaf (or rather, complex of coherent sheaves) associated to $\Pi_\lambda$.  The answer depends
on~$\lambda$.  Indeed, the geometry  allows us to normalize the
embeddings~\eqref{eqn:embedding}. 
Namely, each of the arrows
$\cInd_{KZ}^G \sigma_{i+1} \to \cInd_{KZ}^G \sigma_i/T$
and
$\cInd_{KZ}^G \sigma'_{i+1} \to \cInd_{KZ}^G \sigma_i/T$
corresponds, on the coherent sheaf side,
to an arrow
$\cO_{\cC_{\sigma_{i+1}}}(-1,0) \to \cO_{\cC_{\sigma}}(-1)/(vy),$
respectively
$\cO_{\cC_{\sigma'_{i+1}}}(-1,0) \to \cO_{\cC_{\sigma}}(-1,0)/(vy),$
and we can take these to simply be the restriction maps induced by
the inclusions $\cC_{\sigma}^{vy  = 0} \hookrightarrow \cC_{\sigma_{i+1}}, \cC_{\sigma'_{i+1}}$.  
Applying $\gA$ to the morphism~\eqref{eqn:diagram}, 
with the morphisms being interpreted in the manner just described,
one can (with a little effort) explicitly compute the complexes
$\fA(\Pi_{\lambda})$.

We summarize the results of this computation
in the following statement,
which shows in particular that the value $\lambda = 1$ is distinguished.

\begin{conc} \leavevmode
  \begin{enumerate}
  \item $H^0\bigl(\mathfrak A(\Pi_1)\bigr) = \text{ skyscraper supported
      at } x,$ where $x$ is the closed point obtained as the
    intersection of the components $\cC_{\sigma_i}$ ($i = 0, \dots,3$),
    while $H^{-1}\bigl(\mathfrak A(\Pi_1)\bigr)$ is a line bundle
    supported on
    $\bigcup_{i = 0}^{3} \cC_{\sigma_i} \cup \cC_{\sigma_i'}.$

  \item If $\lambda \neq 1$, then
    $H^0\bigl(\mathfrak A(\Pi_\lambda)\bigr) = 0,$ while
    $H^{-1}\bigl(\mathfrak A(\Pi_\lambda)\bigr)$ is a torsion-free
    sheaf supported on
    $\cC:= \bigcup_{i = 0}^{3} \cC_{\sigma_i} \cup \cC_{\sigma_i'}.$
    More precisely, it is the pushforward to $\cC$ of a line bundle on
    $\cC\setminus \set{x}$ which depends on the parameter~$\lambda$.
  \end{enumerate}
\end{conc}

\subsubsection{Irreducible quotients}\label{subsubsec: irredquots}

Consider the short exact sequence
$$0  \to \mathcal K  \to \Pi_\lambda \to \pi \to 0$$
(where $\pi$ is an irreducible admissible representation as above, and $\mathcal K$ is defined to be the kernel of the surjection).
Passing to the (now merely hypothetical) sheaves, and bearing in mind
that~$\fA$ has amplitude $[-1,0]$ (see Remark~\ref{rem: expect amplitude d-1}), we obtain a long
exact sequence %
\begin{multline*}
0 \to
H^{-1}\bigl(\mathfrak A(\mathcal K)\bigr)
\to
H^{-1}\bigl(\mathfrak A(\Pi_\lambda)\bigr)
\to
H^{-1}\bigl(\mathfrak A(\pi)\bigr)
\\
\to
H^0\bigl(\mathfrak A(\mathcal K)\bigr)
\to
H^0\bigl(\mathfrak A(\Pi_\lambda)\bigr)
\to
H^0\bigl(\mathfrak A(\pi)\bigr)
\to 0.
\end{multline*}

If $\lambda \neq 1$, then we find that 
$H^0\bigl(\mathfrak A(\pi)\bigr) = 0$
(since
$H^0\bigl(\mathfrak A(\Pi_{\lambda})\bigr) = 0$),
and so $\mathfrak A(\pi)$ is supported in degree~$-1$.
If $\lambda = 1$, then 
since $H^0\bigl(\mathfrak A(\Pi_1)\bigr)$
is
the skyscraper supported at~$x$,
we see either that
$H^0\bigl(\mathfrak A(\pi)\bigr)$  is again
the skyscraper supported at~$x$,
or else that
$H^0\bigl(\mathfrak A(\pi)\bigr) = 0$.

On the other hand, since $\pi$ is not finitely presented~\cite{MR3365778},
Lemma~\ref{lem:compact objects}~(3) implies that
$H^{-1}\bigl(\mathfrak A(\pi)\bigr)$ is {\em not} coherent,
 no matter the value of $\lambda$.

What can we say about the sheaf $H^{-1}\bigl(\mathfrak A(\pi)\bigr)$?
Here is one property it has:  since $\Hom(\cInd_{KZ}^G \sigma' , \pi) = 0$
if $\sigma' \neq \sigma_i$ for some $i = 0,\dots, 3$,
we find that
$$\Ext^1\Bigl( \cO_{\cC_{\sigma'}}, H^{-1}\bigl(\mathfrak A(\pi)\bigr)\Bigr)
= 0$$
for any such $\sigma'$. %
If $\lambda = 1$, %
a representation~
$\pi$ as above (that is, an irreducible quotient of~$\Pi_\lambda$
whose $K$-socle consists precisely of the $\sigma_i$, each with multiplicity one) can also be constructed globally,
via the completed cohomology of definite quaternion algebras and/or
Shimura curves. (That $\lambda$ necessarily equals~$1$ in this case was proved by
Dotto--Le~\cite{MR4283560}. The  representations in
completed cohomology have been studied by many authors, using
Taylor--Wiles patching as a key ingredient, and we refer to the
introduction to~\cite{https://doi.org/10.48550/arxiv.2102.06188} for a
brief overview of what is known.)

In the case when $\pi$ is constructed via completed cohomology,
the expected relationship between the functor $\mathfrak A$
and completed cohomology
suggests that the map
$H^0\bigl(\mathfrak A(\Pi_1)\bigr)
\to
H^0\bigl(\mathfrak A(\pi)\bigr)$
should be non-zero,
and so at least in this case we anticipate that
$H^0\bigl(\mathfrak A(\pi)\bigr)$ is equal to the skyscraper
at~$x$. %
We
explore the sheaf $\mathfrak A(\pi)$ further in Section~\ref{subsec: adjoint
  functor}, where we see that our conjectures imply that the
irreducible representation~$\pi$ should be uniquely determined by the
requirement that $H^0\bigl(\mathfrak A(\pi)\bigr)\ne 0$ (see
Remark~\ref{rem: we expect a unique supersingular pi with skyscraper
  H0}).

\subsection{The adjoint functor}\label{subsec: adjoint
  functor}Return to  the situation of Conjecture~\ref{conj: Banach
  functor}, and the continuous functor
$\fA:\Ind \Dfp^b(\smG) \to \Ind\Coh(\cX_d)$.  By Lemma~\ref{lem:compact
  objects}~(1) (i.e.\ the adjoint functor theorem), ~$\fA$ admits a right
adjoint~$\fB: \Ind\Coh(\cX_d)\to \Ind \Dfp^b(\smG)$. 
  By part~(2) of this same lemma (which applies, because~$\fA$
  preserves compact objects by assumption; see Remark~\ref{rem: compact objects and fg representations}), we see that $\fB$ is again continuous. %

\begin{rem}%
  \label{rem: interpretation of adjoint in FS Mann context}In the
  setting of Remark~\ref{rem:why only fully faithful functors} (i.e.\
  of a $p$-adic version of the Fargues--Scholze conjecture), %
  the
  functors~$\fA$ and~$\fB$ will presumably be induced by the functors
  $i_!$ and $i^!$, where $i$ is the inclusion of the open stratum $[*/G(F)]\into \Bun_G$. %
\end{rem}

We can use the calculations of the previous section to investigate the
functor~$\fB$ in the cases of~$\PGL_2(\Qp)$
and~$\PGL_2(\Q_{p^2})$. (Of course, this is not precisely the setting
of Conjecture~\ref{conj: Banach
  functor}, but the above remarks go over unchanged to this case.) This will allow us to
connect the conjectures of this paper to the more traditional approach to the $p$-adic
and mod~ $p$
local Langlands correspondences via the cohomology of Shimura curves,
and in particular to the hope that it is possible to associate to each representation
$\rhobar:\Gal_{\Q_{p^2}}\to\GL_2(\Fpbar)$ a specific %
admissible representation of~$\GL_2(\Q_{p^2})$, %
which is determined by its occurrence in any (that is to say, every!) global
context related to a global lift of~$\rhobar$
(see for example \cite[Q.\ 3.11]{MR2827792}). We return to this theme
in Remark~\ref{rem: we expect a unique supersingular pi with
  skyscraper H0} below. %

We can study the expected properties of~$\fB$ by means of the expected
compatibility of~$\fA$ with Taylor--Wiles
patching. Let~$F/\Qp$ be arbitrary, and let~$\cX$ be the stack of rank~2
$(\varphi,\Gamma)$-modules with determinant~$\varepsilon^{-1}$; and
let~$L_\infty$ be the (hypothetical, for $F\ne \Qp$) pro-coherent sheaf which is a
kernel for~$\fA$, in the sense that %
\[\fA(\pi):= L_{\infty}\otimes^{L}_{\cO[[G]]}\pi.\] Then the adjoint functor is given by \numequation\label{eqn: formula
  for adjoint in terms of L infty}\fB(\cF):=
\RHom_{\pro-\cO_{\cX}}(L_{\infty},\cF).\end{equation}%

Suppose now that $x:\Spec k\to\cX$ is a finite type %
point, corresponding
to a continuous representation $\rhobar:\Gal_F\to \GL_2(k)$, and let $f:\Spf
R_\infty\to\cX$ be the versal morphism arising from patching at some
globalization of~ $\rhobar$. %
By Remark~\ref{rem: conjectural compatibility with TW patching} we expect that  $f^*L_\infty$ is equal
to the patched module~$M_\infty$. For generic choices of~$\rhobar$,
the ring~$R_\infty$ is formally smooth, and it follows
from~\cite[Prop.\ 4.3.1]{MR4386819} and the results
of~\cite{https://doi.org/10.48550/arxiv.2009.03127, MR4461834, https://doi.org/10.48550/arxiv.2209.09639} %
that~$M_\infty$ is pro-flat over~$R_\infty$. We
expect that it is even pro-free. %
Accordingly, we expect that we can compute $\fB(\cF)$ via the formula
\begin{align*}
  \fB(\cF) &= \RHom_{\pro-\cO_{\cX}}(L_\infty,\cF)\\ &=
                                    R\Gamma\circ\sRHom_{\pro-\cO_{\cX}}(L_\infty,\cF)
  \\ &=R\Gamma(\cX, \sHom_{\pro-\cO_{\cX}}(L_\infty,\cF)).
\end{align*}

Write
$G_x=\Aut_{\Gal_F}(\rhobar)$. Assume further that~$\rhobar^{\semis}$ is not
scalar (this is automatically the case if as above~$\rhobar$ is assumed to be generic); so
$G_x=\boldsymbol{\mu}_2$ if~$\rhobar$ is irreducible or non-split reducible, and equals ~$\Gm$ if
$\rhobar$ is a sum of two distinct characters. Then the residual gerbe
at $x$ is an affine immersion\footnote{If $\rhobar$ is semisimple
then $i$ is a closed immersion; otherwise $i$ factors as
$[\Spec k/\boldsymbol{\mu}_2] \buildrel \text{open}\over \hookrightarrow 
[\A^1_{k}/\Gm] \buildrel \text{closed}\over\hookrightarrow \cX,$ which
exhibits it as a composite of affine immersions.}
\[i:[\Spec k/G_x]\into\cX,\]and we can
form the quasicoherent sheaf $\delta_x:=i_*\cO_{[\Spec k/G_x]}$. If, as above, $f:\Spf
R_\infty\to\cX$ is a versal morphism at~$x$ arising from Taylor--Wiles
patching, and~$\m_\infty$ denotes the maximal ideal of~$R_\infty$, then we see that %
\[\sHom_{\pro-\cO_{\cX}}(L_\infty,\delta_x)=M_\infty^\vee[\m_\infty]\]equipped
with some $G_x$-action. (To see this, note that %
$L_\infty \otimes_{\cO_{\cX}} \delta_x = M_\infty / \m_\infty M_\infty$
equipped with some $G_x$-action,
and that the continuous $k$-dual of $M_\infty/\m_\infty M_\infty$ 
equals~$M_\infty^\vee[\m_\infty]$.) Note that by the construction of~
$M_\infty$ by Taylor--Wiles patching, $M_\infty^\vee[\m_\infty]$ is
equal to the $\m$-torsion in completed cohomology, where~$\m$ is a
certain maximal ideal of an appropriate Hecke algebra. %
Since~$G_x$ is linearly reductive (being equal to either $\boldsymbol{\mu}_2$ or
$\Gm$),
we then see that \[\fB(\delta_x)=(M_\infty^\vee[\m_{\infty}])^{G_x}.\]  %

We can recover all of $M_{\infty}^\vee[\m_{\infty}]$ by considering the values
\[\fB\bigl(\delta_x(n)\bigr)=\bigl(M_\infty^\vee[\m_{\infty}](n)\bigr)^{G_x}\]  
as $n$ varies.
Or, rather than considering the various twists $\delta_x(n)$ individually,
if we let $g:\Spec k \to \cX$ denote  the composite
$$\Spec k \to [\Spec k/G_x] \buildrel i \over \hookrightarrow 
\cX,$$
then we can apply $\fB$  to  the quasicoherent  sheaf $\Delta_x := g_*\cO_{\Spec k}.$
(This is the quasicoherent sheaf over the residual gerbe of~$x$
which corresponds to the regular representation of~$G_x$.)
Indeed, we compute that
\begin{multline*}
\fB(\Delta_x)  = \RHom_{\pro-\cO_{\cX}}(L_{\infty},g_*\cO_{\Spec \cF})
\\
= \RHom_{\pro-{\cO_{\Spec k}}}(g^*L_{\infty},\cO_{\Spec  \cF})
= \Hom_{\cont}(M_{\infty}/\mathfrak \m_{\infty}, k) = M_{\infty}^{\vee}[\m_{\infty}].
\end{multline*}
Thus  $x \mapsto \fB(\Delta_x)$ should recover the 
hoped-for ``traditional'' 
mod $p$ local Langlands correspondence alluded to above.
The $G_x$-action on $M_{\infty}^{\vee}[\m_{\infty}]$ that
arises by considering 
the sheaf-Hom
$\sHom_{\pro-\cO_{\cX}}(L_\infty,\delta_x)$
as above is some additional structure in the mod $p$ Langlands correspondence
that (as far as we are aware) has not been previously considered
(although one could say that it is implicit in the ideas of Breuil and his collaborators
expressed in papers such as
\cite{MR2783977} and~\cite{MR3347316}).

\subsubsection{Some examples}
In the cases that $F=\Qp$ or $F=\Q_{p^2}$, we know something
about $M_\infty^\vee[\m]$ (in the former case by M.E.'s local-global
compatibility theorem~\cite{emerton2010local}, and in the latter case
by~\cite{https://doi.org/10.48550/arxiv.2102.06188}), and we either know or can guess something
about how~$G_x$ acts. We have the following examples.
\begin{enumerate}
\item
If $F=\Qp$ and $\rhobar=\chibar_1\oplus\chibar_2$ with
  $\chibar_1 \chibar_2^{-1} \neq \triv$ or $\varepsilon^{\pm 1}$,
then $M_\infty^\vee[\m]$ is a direct sum
  $\pi_1\oplus\pi_2$ of two irreducible principal series
  representations, and in fact for an appropriate ordering of $\pi_1,\pi_2$ we
  have %
  \[\sHom_{\pro-\cO_{\cX}}(L_\infty,\delta_x)=\pi_1(1)\oplus\pi_2(-1)\]
  (where the twists denote the weight of the $\Gm$-action; note that
  the apparent asymmetry is due to the need to choose an ordering of
  $\chibar_1,\chibar_2$ in order to identify $G_x$ with~$\Gm$).
  
  In this case $\Delta_x := \oplus_{n = -\infty}^{\infty} \delta_x(n),$
  and so %
$$\fB(\Delta_x) = %
\fB\bigl(\delta_x(-1)\bigr) \oplus 
\fB\bigl(\delta_x(1)\bigr) 
=  \pi_1 \oplus \pi_2.$$
  If~ $\sigma_1,\sigma_2$ are (a chosen ordering of) the two
  Serre weights of~ $\rhobar$, we find that $(\fA\circ\fB)\bigl(\delta_x(-1)\bigr)=\fA(\pi_1)$
is supported on
  $\cX(\sigma_1)$ while $(\fA\circ\fB)(\delta_x(1))=\fA(\pi_2)$ is
  supported on $\cX(\sigma_2)$. 
  As a side remark,
  we note that although $\delta_x$ is supported on the intersection
  $\cX(\sigma_1)\cap\cX(\sigma_2)$,
  these representations are not;
  the representations $\fA(\pi_i)$ are supported on the
  closed locus in the respective component $\cX(\sigma_i)$
corresponding to Galois representations with
  semisimplification~$\rhobar$.
  Thus $(\fA\circ \fB)(\Delta_x) = \fA(\pi_1) \oplus \fA(\pi_2)$ is supported
  on the union of these two loci, which is the locus of all Galois representations
  having $\rhobar$ as their semisimplification.
\item
If $F=\Qp$ and $\rhobar$  is a non-split  extension
of characters $\chibar_1$  and $\chibar_2$, again with
  $\chibar_1 \chibar_2^{-1} \neq \triv$ or $\varepsilon^{\pm 1}$,
then $M_\infty^\vee[\m]$ is a non-split extension of the two irreducible principal series
  $\pi_1$ and $\pi_2$ appearing in case~(1),
and 
 \[\sHom_{\pro-\cO_{\cX}}(L_\infty,\delta_x)= M_\infty^{\vee}[\m](1),\]
the twist indicating that $\boldsymbol{\mu}_2$ acts non-trivially. 

In this case
$\Delta_x := \delta_x \oplus \delta_x(1)$,
and so 
$$\fB(\Delta_x) = \fB\bigl(\delta_x(1)\bigr)=M_\infty^{\vee}[\m]$$
(the non-split extension of $\pi_1$ and $\pi_2$).
Furthermore, 
$ (\fA\circ\fB)(\Delta_x) = (\fA\circ\fB)\bigl(\delta_x(1)\bigr)$
is a non-split extension of $\fA(\pi_1)$ and $\fA(\pi_2)$,
and in particular has support equal to 
the closed locus corresponding to Galois representations with
semisimplification equal to $\chi_1\oplus \chi_2$.
\item If $F=\Qp$ and $\rhobar$ is irreducible, then
  $M_\infty^\vee[\m]=\pi$ is a supersingular irreducible 
  representation, and \[\sHom_{\pro-\cO_{\cX}}(L_\infty,\delta_x)=\pi(1)\]
  (the twist indicating that the action of~$\boldsymbol{\mu}_2$ is nontrivial).
  
  As in case~(2), we have
$\Delta_x := \delta_x \oplus \delta_x(1)$,
and so 
$$\fB(\Delta_x) = \fB\bigl(\delta_x(1)\bigr)=\pi,$$
and 
$ (\fA\circ\fB)(\Delta_x) = (\fA\circ\fB)\bigl(\delta_x(1)\bigr)
=\delta_x(1)$.
\item If $F=\Q_{p^2}$ and $\rhobar=\chibar_1\oplus\chibar_2$ with
  $\chibar_1,\chibar_2$ suitably generic, then by~\cite[Cor.\
  1.3.10]{https://doi.org/10.48550/arxiv.2102.06188}, we have
  $M_\infty^\vee[\m]=\pi_1\oplus\pi\oplus\pi_2$, where $\pi_1,\pi_2$
  are explicit irreducible principal series representations, and $\pi$
  is a (non-explicit!) supersingular irreducible representation. We
  expect %
  that \[\sHom_{\pro-\cO_{\cX}}(L_\infty,\delta_x)=\pi_1(2)\oplus\pi\oplus\pi_2(-2).\]

  Assuming that this holds,  
  we see that %
$$\fB(\Delta_x) = %
\fB\bigl(\delta_x(-2)\bigr) \oplus 
\fB(\delta_x) 
\oplus
\fB\bigl(\delta_x(2)\bigr) 
=  \pi_1 \oplus \pi \oplus \pi_2.$$ (Here we used that (as in case~(1)) we have  $\Delta_x := \oplus_{n = -\infty}^{\infty} \delta_x(n)$.)
In particular, we find that 
$\fB(\delta_x)=\pi$ is irreducible and supersingular. (A similar
analysis could be performed in the case of a nonsplit extension of
characters, using the results of~\cite{MR4461834}.)
\item If $F=\Q_{p^2}$ and $\rhobar$ is irreducible and suitably
  generic, then by~\cite[Thm.\
  1.3.11]{https://doi.org/10.48550/arxiv.2102.06188}, we have
  $M_\infty^\vee[\m]=\pi$, where ~$\pi$ is a supersingular irreducible representation. We
  expect %
  that \[\sHom_{\pro-\cO_{\cX}}(L_\infty,\delta_x)=\pi\](with
  trivial $\boldsymbol{\mu}_2$-action).  Assuming that this holds, 
we see that
$$\fB(\Delta_x) = \fB(\delta_x)=\pi$$
  is irreducible and supersingular. (Here we used that (as in case~(3)) we have
$\Delta_x := \delta_x \oplus \delta_x(1)$.)
\end{enumerate}

\begin{rem}
  \label{rem: guessed the twists following Breuil}
  The twists in
  examples~(4) and~(5) above were, at first,  arrived at heuristically, based on
  examples~(1), (2), and~(3) and the philosophy suggested by Breuil's
  paper~\cite{MR2783977}.  However, one can confirm them (within our conjectural
  framework, i.e.\ assuming that $\fA$ exists with its conjectured properties,
  and that the values of $\fA(\cInd_{KZ}^G \sigma)$ for sufficiently generic Serre weights
  $\sigma$ are as conjectured in Expectation~\ref{exp:Serre weights}) in the following manner:  
  The Serre weights $\sigma$ of $\rhobar$ are known, and are precisely the
  irreducible  constituents of the $K$-socle of $M_{\infty}^{\vee}[\m_{\infty}]$,
  and are known to appear with multiplicity one.
  It is also known how they distribute themselves among the various direct summands 
  of $M_{\infty}^{\vee}[\m_{\infty}].$ (See~\cite{GLSII, geekisin,
    MR3323575, MR3274546} for these results.)
 Thus, as $\sigma$ ranges over the Serre weights of~$\rhobar$, the multiplicity spaces
 $\Hom_G(\cInd_{KZ}^G \sigma, M_\infty^{\vee}[\m_\infty])$ are one-dimensional,
 and it is known which direct summand of $M_\infty^{\vee}[\m_\infty]$ any non-zero morphism
 factors through.

On the other hand, we can write this Hom space as
\begin{multline*}
 \Hom_G(\cInd_{KZ}^G \sigma, M_\infty^{\vee}[\m_\infty])
= \Hom_G\bigl(\cInd_{KZ}^G \sigma, \fB(\Delta_x)\bigr)
\\
= \Hom_G\Bigl(\cInd_{KZ}^G \sigma, \fB\bigl(\oplus_n \delta_x(n)\bigr)\Bigr)
= \Hom_G\Bigl(\cInd_{KZ}^G \sigma, \oplus_n \fB\bigl(\delta_x(n)\bigr)\Bigr)
\\
= \oplus_n \Hom_G\Bigl(\cInd_{KZ}^G \sigma, \fB\bigl(\delta_x(n)\bigr)\Bigr),
\end{multline*}
the sum ranging over either $n \in \Z$
or $n \in \Z/2\Z$, depending on whether $\rhobar$ is a direct sum of
characters  or is irreducible.
Now the adjunction between $\fA$ and $\fB$  shows that
\begin{multline*}
\Hom_G\Bigl(\cInd_{KZ}^G \sigma, \fB\bigl(\delta_x(n)\bigr)\Bigr)
= \Hom_G\bigl(  \gA(\cInd_{KZ}^G \sigma), \delta_x(n)\bigr)
\\
= \Hom_G\bigl( \cO_{\cX(\sigma)}(-1,-1), \delta_x(n)\bigr)
=  i^*\bigl(\cO_{\cX(\sigma)}(-1,-1)\bigr)(n)^{G_x},
\end{multline*}
where, as above, 
$i:[\Spec k/G_x] \hookrightarrow \cX$ is the %
immersion of the residual gerbe
at $x$ into~$\cX$.

Since $\cO_{\cX(\sigma)}(-1,-1)$ is an invertible sheaf, we see that
$i^*\bigl(\cO_{\cX(\sigma)}(-1,-1)\bigr)$ is a character of~$G_x$,
and  so has the form $k(m)$ for some  integer~$m$ (which we see is determined by $\sigma$ and~$x$). 
Thus 
$i^*\bigl(\cO_{\cX(\sigma)}(-1,-1)\bigr)(n)^{G_x}  = k(m+n)^{G_x}$
is non-zero for exactly one value of $n$ (namely $n  = - m$),
and this determines the twist of the corresponding constituent 
of $M_{\infty}^{\vee}[\m_\infty]$; namely, it appears with a twist by~$m$.
Computing the values of $m$ for the various possible choices of $x$ (which is to say,
of~$\rhobar$) and of $\sigma$ then yields the twists appearing in  examples~(4)  and~(5).     
\end{rem}

\begin{rem}
As already recalled in the discussion of Section~\ref{subsubsec: irredquots},
it is shown in~\cite{MR3365778} that supersingular representations of~$\GL_2(\Q_{p^2})$
are never finitely presented.  Comparing this fact with the preceding (conjectural)
computations, we see that we should not expect the adjoint functor $\fB$ 
to preserve compact objects.  (We saw that it should take certain skyscraper sheaves,
which are coherent and hence compact as objects of the stable $\infty$-category
of Ind-coherent sheaves on~$\cX$, to non-finitely presented representations
of~$\GL_2(\Q_{p^2})$, which are not compact as objects of~$\Ind \Dfp^b(\smG)$.) 

To see why this might be, note that 
the adjoint functor~$\fB$ should be computable in terms of a semiorthogonal 
decomposition of~$\Ind\Coh(\cX)$ of the kind conjectured in Section~\ref{subsec: semiorthogonal Qp2}.
Namely, the conjectural functor $\fA$ should induce an equivalence  
between $\Ind \Dfp^b(\smG)$ and the ``first piece'' of a semiorthogonal decomposition
of~$\Ind\Coh(\cX)$, while the other pieces of this decomposition should be obtained via
``geometric Eisenstein series''-type constructions, which involve pushforward from
the reducible locus in~$\cX$ (as well as more mysterious contributions from 
$\mathcal{BC}(\cO(1))$ and $\mathcal{BC}(\cO(2))$).  The composite functor~$\fA \circ \fB$ can 
then be described  as ``projection away from'' these other pieces so as to land
in the essential image of~$\fA$ (and so
$\fB$ itself can be described as this projection,
 followed by an application of $\fA^{-1}$ on the resulting
object in this essential image). 

Since the 
inclusion of the reducible locus of a typical irreducible component $\cX(\sigma)$
into $\cX(\sigma)$ is an open, but not closed, immersion (and so not proper),
the objects in these ``other pieces'' of the semiorthogonal decomposition
 are typically not coherent. Thus, evaluating
$\fA \circ \fB$ on a coherent sheaf will typically involve forming cones of morphisms from
these non-coherent sheaves to the given coherent sheaf, and so we shouldn't expect
the values of $\fA \circ \fB$ on coherent sheaves to themselves be coherent in general.
Correspondingly, we shouldn't expect the value of $\fB$ on a coherent sheaf
to be a compact object of $\Ind\Dfp(\smG)$ in general. 
\end{rem}

\begin{rem}
  \label{rem: we expect a unique supersingular pi with skyscraper H0}
Since the results of Breuil and Pa\v{s}k\={u}nas~\cite{BreuilPaskunas}
it has been clear that the formulation of the $p$-adic local Langlands
correspondence for $\GL_2(\Q_{p^2})$ would have to be very different
from the (classical) correspondence for $\GL_2(\Q_p)$, for the reason
that there appear to be many more representations of $\GL_2(\Q_{p^2})$
than there are 2-dimensional representations of $\Gal_{\Q_{p^2}}$. In
particular, as already discussed in Section~\ref{subsubsec:
  irredquots}, there are infinite families of supersingular irreducible admissible
mod $p$ representations of~$\GL_2(\Q_{p^2})$, all of which
appear to correspond to the same irreducible representation~$\rhobar$ of
$\Gal_{\Q_{p^2}}$.

Our conjectural functor~$\fA$ gives one explanation of this
phenomenon: namely, the apparent lack of Galois representations is
made up for by an abundance of coherent sheaves on~$\cX$, and in
particular, these supersingular representations all correspond to
Ind-coherent sheaves on~$\cX$, which are distinguished by their
(quasicoherent, but not coherent) $H^{-1}$ sheaves.

There is, however, an earlier approach to reconciling the results of~\cite{BreuilPaskunas} with the expectation that a $p$-adic local Langlands correspondence for $\GL_2(\Q_{p^2})$ should exist; namely one can hope for the existence of a uniquely determined, canonical element~$\pi_{\rhobar}$ within the infinite family of supersingular representations constructed in~\cite{BreuilPaskunas}. Moreover, this representation~$\pi_{\rhobar}$ should be characterised by local-global compatibility: it
should appear in the cohomology of Shimura curves for any globalisation of~$\rhobar$. In this sense, the existence of~$\pi_{\rhobar}$ amounts to showing that the representations arising in the cohomology of Shimura curves depend only on~$\rhobar$ -- and that this representation is indeed a supersingular irreducible representation of the type constructed in~\cite{BreuilPaskunas}. The question has been studied by many authors, and there is substantial evidence supporting this viewpoint; we refer again to the introduction of~\cite{https://doi.org/10.48550/arxiv.2102.06188} for a concise summary of these investigations.

We strongly believe in the existence of the canonical representations~$\pi_{\rhobar}$. Indeed, this perspective motivated the proposed construction of a $p$-adic local Langlands correspondence via Taylor–Wiles
patching in~\cite{Gpatch}, which in turn inspired the conjectures in these notes
(see Remark~\ref{rem: conjectural compatibility with TW patching}).
We discuss local-global compatibility more fully in Section~\ref{subsec:computing-cohomology-Shimura} below. The computations in this section are also consistent with this %
expectation, and indeed we see from examples~(4) and~(5) above that if~$x$ is
the closed point corresponding to a suitably generic irreducible~$\rhobar$, then
we should have $\pi_{\rhobar}=\fB(\delta_x)$.

Furthermore, within our conjectural framework, we claim that
$\pi_{\rhobar}$ can be characterised among all irreducible
representations ~$\pi$ of $\GL_2(\Q_{p^2})$ by the following property:
it is the unique irreducible~$\pi$ for which the support of $H^0(\fA(\pi))$ contains~$x$. %
Indeed, suppose that $\pi$ is irreducible and that the support of $H^0(\fA(\pi))$ contains~$x$. Since~$\pi$ and~$\pi_{\rhobar}=\fB(\delta_x)$ are irreducible,
it is enough to show that
\[\Hom_{\GL_2(\Q_{p^2})}(\pi,\fB(\delta_x))\ne 0.\]By adjunction, this
is equivalent to showing that
\[\Hom_{\pro-\cO_{\cX}}(\fA(\pi),\delta_x)\ne 0.\] To see that this holds, we
consider the canonical morphism \[\fA(\pi)\to\tau^{\ge 0}\fA(\pi),\] and note
that $\tau^{\ge 0}\fA(\pi)$ is coherent (by Remark~\ref{rem: compact
  objects and fg representations}), and thus admits a morphism to~$\delta_x$ by the
assumptions that the support of $H^0(\fA(\pi))$ contains~$x$ (here we
use that $\fA(\pi)$ is contained in the even part of $\Ind\Coh(\cX)$, see
Section~\ref{subsec:
  semiorthogonal Qp2}).

Note that since the supersingular representation~$\pi_{\rhobar}$
admits a surjection from a representation $\cInd_{KZ}^G\sigma$
(where~$\sigma$ is some Serre weight), the coherent sheaf
$\tau^{\ge 0}\fA(\pi)$ admits a surjection from a line bundle. In
particular in case~(5) we see that $H^0(\fA(\pi))$ should equal~$\delta_x$.
\end{rem}

\begin{rem}
  \label{rem: in theory we could compute F pi rho}
In cases (4) and (5) above we do not have an explicit description of
$\fA(\pi)$, %
except that in case~(5)  its $H^0$ is equal to~$\delta_x$ (as discussed in
Remark~\ref{rem: we expect a unique supersingular pi with skyscraper
  H0}). (Similarly, we expect that in case (4) the ~$H^0$ is the
structure sheaf of the closed locus of all points whose
semisimplification is~$x$.) %
However, $\fA(\pi)$ itself could in principle be computed
via the machinery of Section~\ref{sec:gener-adjunct}, because it is
equal to $(\fA\circ\fB)(\delta_x)$, and we expect to be able to write down a
collection of generators for the image of the hypothetical
functor~$\fA$ (namely, the various coherent sheaves %
corresponding to the $\fA(\cInd_{KZ}^G\sigma)$, which in the generic case we made explicit
in Section~\ref{subsec: Banach Qp2 stuff}). %
\end{rem}

\begin{remark} One anticipated feature of the
the traditional mod $p$ local 
Langlands correspondence for  $\GL_2(F)$
is that, for reducible~$\rhobar$,
and writing (in terms of the notation introduced above) $\pi(\rhobar) := \fB(\Delta_x)$, 
one should have 
\numequation
\label{eqn:compatibility with semisimplification}
\pi(\rhobar^{\semis})\cong \pi(\rhobar)^{\semis}.
\end{equation}
It is natural to ask if this can be  seen from our present perspective.

To this end, we note that if $\rhobar^{\semis}$ is 
a  direct sum  of  two suitably  generic characters $\chi_1$ and $\chi_2$ of $G_{\Q_{p^f}}$,
then the locus of $\rhobar$ which are extensions of $\chi_1$ by $\chi_2$
is equal to a copy of $[\A^f/\Gm]$.
If we let $\cF$ denote the structure sheaf of $[\A^f/\Gm]$,
pushed forward to a coherent sheaf on~$\cX$,
then $\sHom(L_{\infty},\cF)$ can be interpreted as a $\Gm$-equivariant
flat family of smooth $\GL_2(F)$-representations over~$\A^f$.  
If we forget the $\Gm$-equivariant structure, this family interpolates
the various $\pi(\rhobar)$ as $\rhobar$ runs over the Galois representations
(i.e.\ extensions of characters) parameterized by~$\A^f$.  It 
seems plausible that the existence of this flat family (with its $\Gm$-equivariant structure!)
could imply 
the anticipated isomorphism~\eqref{eqn:compatibility with semisimplification}.
\end{remark}

  \subsection{Global functions and Hecke algebras}\label{subsec:
    global functions and Hecke algebras}The full faithfulness of the
  putative functor $\gA$ has some concrete consequences for the
  geometry of the stacks~$\cX_d$. For example, we saw above that if
  $G=\GL_2(F)$ with $F/\Qp$ unramified, then for a generic Serre
  weight $\sigma$ we expect $\gA(\cInd_K^G\sigma)$ to be a line bundle
  on the component $\cX(\sigma)$; so that in particular the
  (underived) Hecke algebra $\End_G(\cInd_K^G\sigma)$ should be
  identified with the global functions on~$\cX(\sigma)$.

  This expectation follows immediately from the explicit description
  of $\cX(\sigma)$ in terms of Fontaine--Laffaille theory, and 
  similar results can be proved for tame types by explicit integral
  $p$-adic Hodge theory (again in the generic case for
  $\GL_2(F)$ with $F/\Qp$ unramified).

More generally, the following result is proved in the paper~\cite{hjelle2024functionsirreduciblecomponentsemertongee} of Eivind Otto Hjelle, Louis
  Jaburi, Rachel Knak,  Hao Lee, and Shenrong
  Wang: %
  if $\sigma$ is a suitably generic Serre weight for~$\GL_d(F)$,
  with~$F/\Qp$ unramified, and $\cX(\sigma)$ is the irreducible component
  of the stack~$\cX_{F,d}$ labelled by $\sigma$, then the global
  functions on $\cX(\sigma)$ are given by a ring of the form
  $k[T_1,\dots,T_{d-1},T_d,T_d^{-1}]$, which is isomorphic to the
  Hecke algebra $\End_G(\cInd_K^G\sigma)$.  Although the paper~\cite{hjelle2024functionsirreduciblecomponentsemertongee} %
  does not define a specific map between this Hecke algebra and the ring of global functions, %
such a map has subsequently been constructed by Heejong Lee~\cite{lee2024spectralmodpsatake}, where it is interpreted as a geometric analogue of the
mod $p$ Satake isomorphism of~\cite{MR2771132}.

  In the non-generic case, the
  paper~\cite{https://doi.org/10.48550/arxiv.2209.09439} of Anthony Guzman, Kalyani
  Kansal, Jason Kountouridis, Ben Savoie, and Xiyuan Wang %
  proves under mild hypotheses on~$\sigma$ that,
  for $F/\Qp$ unramified and $\sigma$ a not-necessarily-generic Serre
  weight of $\GL_2(F)$, the global functions on the normalization of
  $\cX(\sigma)$ can be identified with the corresponding Hecke algebra.
These results are further refined by Kansal and Savoie in their paper~\cite{kansal2025nongenericcomponentsemertongeestack}.

\section{Categorical local Langlands in the case \texorpdfstring{ $\ell \neq p$}{l not p}}\label{sec: categoric LL for l not p}
Although our focus so far has been on the categorical
formulation of the $p$-adic local Langlands
conjecture for the group~$\GL_d$ over a $p$-adic local field,
we also wish to address the global context; indeed, there is
a case to be made that the main point of developing the local
aspects of number theory is to provide tools for investigating global
problems. We will explain the global theory in Section~\ref{sec:
  global stacks and cohomology of Shimura varieties}; this naturally
involves the local $\ell\ne p$ case, which we study in this section.

The main global problem we will focus on is the conjectural
description of the cohomology
of Shimura varieties (and, even more speculatively, of general congruence quotients) 
in terms of the sheaves arising from the categorical local conjecture.  On the one hand,
this is related to the completed cohomological approach to the Fontaine--Mazur conjecture
(as in \cite{emerton2010local,pan2019fontainemazur}).  On the other,
it is related to Taylor--Wiles patching (and in particular, provides a conjectural
explanation for the expected relationship between patched modules and the coherent
sheaves of the categorical conjecture). 

In order to study global problems, especially those involving Shimura varieties,
we have to be able to consider reductive groups that are not necessarily split.
Although currently our understanding of $p$-adic local Langlands in
the non-$\GL_d$ case
is even more fragmentary than our understanding for~$\GL_d$, %
the local situation at primes away from $p$ is better understood; indeed, 
there are precise conjectures in this case,
on which significant progress has been made, in particular in  the papers \cite{fargues--scholze,
    zhu2020coherent, hellmann2020derived, benzvi2020coherent, zhu2025tamecategoricallocallanglands}. %
In this section we briefly recall some of these local conjectures and
results, referring to these papers and to the notes~\cite{fargues--scholze-IHES} for further details.

We begin by briefly recalling the formalism of $C$-groups; these will be the targets
of the Galois representations arising in the local context when non-split groups
are involved.
We then describe the moduli stacks of local Langlands parameters 
in the $\ell \neq p$ case, and very briefly recall the categorical
local Langlands conjectures that we need for our global discussion. %

\subsection{\texorpdfstring{$C$}{C}-groups}\label{subsec: C groups}%
As always, we have fixed a finite
extension $L$ of~$\Q_p$, whose ring of integers $\cO$ 
serves as our coefficients.

If $F$ is a 
finite extension of either $\Q_{\ell}$ (for some prime~$\ell$)
or of $\Q$,
and $G$ is
a connected reductive linear algebraic group over~$F$,  
then we can form the $C$-group ${}^cG$ of $G$  %
as in~\cite{MR3444225}; see also ~\cite[\S 3.0]{zhu2020coherent}.
This is a variant of the $L$-group ${}^LG$ which is better adapted to 
issues relating to fields of definition.  
By construction,
it is a semi-direct product %
\[{}^cG := \widehat{G} \rtimes \bigl( \Gm \times \Gal(\tF/F)
  \bigr),\]
where $\tF$ is a certain
finite Galois extension of~$F$,
chosen ``large enough'' for the purposes
at hand, and $\widehat{G}$ is the usual Langlands dual group. 
(In particular, $\tF$ should contain a splitting field for the quasi-split inner form of~$G$; %
the definitions involving ${}^cG$
that appear below are easily seen
to be independent  of the particular choice of~$\tF$.)
Here $\Gm$ acts
via   the
cocharacter $\rho_{\ad}: \Gm \to \widehat{G}_{\ad}$  defined
by one-half the sum  of the positive  coroots  of $\widehat{G}$.
We regard ${}^cG$ as a group scheme
over the ring of integers~$\cO$.

\begin{rem}
  \label{rem: equivalent definitions of the C group}The $C$-group can
  also be defined as the $L$-group of a certain $\Gm$-extension $\widetilde{G}$
  of~$G$, see~\cite[Defn.\ 5.3.2]{MR3444225}. 
\end{rem}

\subsection{Langlands parameters in the local case ---\texorpdfstring{ $\ell \neq p$}{l not p}}
If $F$ is an $\ell$-adic field  with $\ell \neq p,$
then one can define a stack $\cX_{G,F}$ which parameterizes
representations of a Weil--Deligne-type group of $F$
(suitably interpreted)
into~${}^cG.$ 
There are
in fact various ways of defining $\cX_{G,F}$, depending on the base
over  which one wishes to have it live.  We will consider $\cX_{G,F}$ as
a stack over $\Spf \cO$, and then we can make the following  definition:
we
first choose elements $\sigma, \tau \in \Gal_F:=
\Gal(\overline{F}/F)$  %
that respectively lift Frobenius
and a topological generator of tame inertia, and let
$\Gamma_F$ denote the subgroup of $\Gal_F$
generated by
 $P_F$ (the wild inertia subgroup of~$\Gal_F$),
 $\sigma$, and $\tau$.  
We have a short exact sequence of groups
$$1  \to P_F \to \Gamma_F \to \Gamma_F^{\tame}:= \langle \sigma, \tau \, \ \,
\sigma \tau \sigma^{-1}  = \tau^q\rangle \to 1$$
(where $q$ denotes the order of the residue field of $F$).
Note that $P_F$ has a neighbourhood basis of the identity 
consisting of  open subgroups $Q$ which are normal in $\Gamma_F$, so
that we may write $\Gamma_F = \varprojlim_Q \Gamma_F/Q$
as a projective limit of discrete quotients,
with each $Q$ being
a subgroup of $P_F$ (and hence pro-$\ell$).

For any $\cO$-algebra $A$ in which $p$ is  nilpotent,
and for any of the precedingly mentioned subgroups $Q$, %
we define $\cX_{G,F}^Q(A)$ to be the groupoid of
homomorphisms $\Gamma_F/ Q \to {}^cG(A)$ which are
\emph{$C$-homomorphisms}, i.e.\ which satisfy the $C$-group version
of the usual ``$L$-homomorphism'' condition:  the composite %
\numequation\label{eqn: C group hom condition}\Gamma_F \to \Gamma_F/Q \to {}^cG(A) \to  A^{\times}\times \Gal(\tF/F)\end{equation}
(the  second arrow being the tautological projection)
should equal the product of the inverse cyclotomic
character $\Gamma_F \to \Gal_F \to \Z_p^{\times}$ and the canonical surjection
$\Gamma_F \to \Gal_F \to \Gal(\tF/F).$
We then set $\cX_{G,F} := \varinjlim_Q \cX^Q_{G,F};$ each transition map 
is an open and closed embedding (because each $Q$ is pro-$\ell$ 
and $\ell \neq p$), and in this way $\cX_{G,F}$ 
is given the structure of an Ind-algebraic stack.

Actually, while this discussion defines $\cX_{G,F}$ as a classical
stack, we really want to define it as a {\em derived} stack,
which we can do by following the ideas explained in Appendix~\ref{app:derived moduli}. %
A key point is that, even when  $\cX_{G,F}$  is  defined 
as a derived  stack, it turns out to be classical, and 
even a local  complete intersection (\cite{zhu2020coherent}, Prop.~3.1.6).  
Since local Galois cohomology
is supported in cohomological degrees $[0,2]$, 
this amounts to the
classical stack underlying $\cX_{G,F}$ being of the  expected dimension,
which can be proved by an explicit computation, which we now briefly recall.

The first point is  that, for an $A$-valued point
$\rho: \Gamma_F \to {}^cG(A)$ of $\cX_{G,F}$ (here $A$ is a (classical!) $\cO$-algebra
in which $p$ is nilpotent), the tangent complex to  $\cX_{G,F}$  at $\rho$ 
is computed by the shifted continuous cochain complex
$C^{\bullet}(\Gamma_F, \mathfrak{g}\otimes_{\cO} A)[1]$, where 
$\Gamma_F$ acts on $\mathfrak{g}\otimes_{\cO} A$ via the composite of $\rho$
with the adjoint action of ${}^cG$ on its Lie algebra.
This complex has cohomology supported in degrees $[-1,1]$,
and so we find that the cotangent complex of $\cX_{G,F}$ at any point
has cohomology vanishing in degrees $\leq  -2$.
This implies that $\cX_{G,F}$ is  {\em quasi-smooth}, i.e.\ a (possibly
derived) local complete intersection.
Furthermore,
$\cX_{G,F}$ is of its expected dimension (the Euler characteristic
of $C^{\bullet}(\Gamma_F,\mathfrak{g}\otimes_{\cO} A)[1]$  equals zero,
essentially by the local Euler characteristic formula,
and $\cX_{G,F}$ is flat over $\Spf \cO$
of relative dimension zero --- see ~\cite[Prop.~3.2.7]{zhu2020coherent}).
This implies that $\cX_{G,F}$ is  actually a classical local complete intersection,
as claimed
(see the proof of~\cite[Prop.~2.2.13]{zhu2020coherent},
on which the proof of ~\cite[Prop.~3.1.6]{zhu2020coherent} ultimately depends).

\subsection{Categorical local Langlands}\label{subsec: categorical LL
  l not p}
Let $F$ be an $\ell$-adic local field, for some prime~$\ell \neq p$, 
and let $G$  be  a connected reductive group over~$F$.
We let
$D\bigl(G(F)\bigr)$
denote the stable $\infty$-category 
of complexes of smooth representations of  $G(F)$  on torsion
$\cO$-modules.  %
We let
$D\bigl(G(F)\bigr)'$
denote a ``renormalization'' of
$D\bigl(G(F)\bigr)$ analogous to the one  
described in \cite[\S 4.1]{zhu2020coherent}
(in which
$\Rep\bigl(G(F),\cO\bigr)$
is replaced by
$\Rep^{\mathrm{ren}}\bigl(G(F),\cO\bigr)$,
the key difference between that context and ours being that no torsion
hypothesis is applied there); 
in $D\bigl(G(F)\bigr)'$, all compactly induced representations are
compact objects (by definition).

We then have the following categorical
form of the local Langlands conjecture in the $\ell \neq p$ context.

\begin{conj}
\label{conj:categorical local l not p}
There is a fully faithful continuous functor of stable $\infty$-categories
$$\mathfrak A_G: 
D\bigl(G(F)\bigr)'
\to
\Ind\Coh(\cX_{G,F}).$$
\end{conj}

\begin{rem}
  \label{rem: expect that l not p functor is given by a kernel}As in
  Remark~\ref{rem: expect Linfty is a kernel}, we expect that the
  functor $\mathfrak A_G$ of Conjecture~\ref{conj:categorical local l
    not p} should be given by the (derived) tensor product with a
  kernel. In the case ~$G=\GL_n$ this expectation was first made
  explicit by E.H.\ in~\cite{hellmann2020derived}; here for the
  Iwahori block of $L$-representations of~$\GL_n(F)$, the functor is
  expected to be given by the derived tensor product over the Iwahori
  Hecke algebra with the Iwahori invariants in the Emerton--Helm--Moss
  local Langlands correspondence in families.
\end{rem}

Conjecture~\ref{conj:categorical local l not p}, as stated, is rather  imprecise. %
A much more precise version is stated by Zhu as~\cite[Conj.\ 4.5.1]{zhu2020coherent}.
(In that statement there is also no restriction to torsion
representations, and the target category is the  Ind-coherent sheaves
on an algebraic $\cO$-stack, of which our stack $\cX_{G,F}$ is the $p$-adic
completion.)
The  analogous conjecture in the case of $L$-coefficients
has been made independently by
E.H.~\cite{hellmann2020derived}.  When $G$ is split, and continuing to work
with $L$-coefficients,
an analogous conjecture was made, yet again independently,
by Ben-Zvi--Chen--Helm--Nadler~\cite{benzvi2020coherent}, who furthermore
construct the restriction of the conjectural
functor to the Iwahori block. 
This restriction (and an extension of it to the general case of
unipotent parameters) was also constructed independently
by Hemo--Zhu~\cite{zhu2025tamecategoricallocallanglands}.  
We remark that when $G$ is not of the form~$\GL_d$
(for some~$d$),
the functor~$\mathfrak A_G$ depends on the choice of auxiliary data
(e.g.\  a pinning and associated Whittaker data, in the case 
that $G$ is quasi-split).
There should be line bundles on the stack  $\cX_{G,F}$
that allow one to twist between the functors arising from
different choices of auxiliary data (see~\cite[Rem.~4.5.7]{zhu2020coherent}). 

In the case  $G = \GL_d$, the construction
of \cite{benzvi2020coherent} extends to the entire category of smooth
$\GL_d(F)$-representations (essentially because all the Hecke algebras
appearing in the theory of the Bernstein centre for  $\GL_d$ are
Iwahori Hecke algebras). As already explained in Remark~\ref{rem: expect that l not p functor is given by a kernel}, E.H. gives another conjectural
description of this functor (still in the  case of $L$-coefficients)
in terms of the ``local Langlands in families''
whose existence was conjectured in~\cite{MR3250061}, and proved by Helm and
Moss~\cite{MR3867634}.
It seems likely that E.H.'s conjectural description should apply also
in the case of $\cO$-coefficients; see e.g.~\cite[Rem.\ 4.4.4]{zhu2020coherent},
where Zhu also gives yet another conjectural description of the functor in
the case of~$\GL_d$, inspired by a conjecture made by
Braverman--Finkelberg~\cite{braverman2019quasicoherent}
in the context of geometric Langlands.  (As far as we know, even for $L$-coefficients,
the construction of~\cite{benzvi2020coherent} has not  been shown to match
with either E.H.'s conjectural description via local Langlands in families,
or Zhu's conjectural description \`a la Braverman--Finkelberg.)

In fact, the fully faithful functor
of~Conjecture~\ref{conj:categorical local l not p}
should actually arise as the restriction of an equivalence of categories
coming from an arithmetic analogue of the geometric Langlands correspondence
in its categorical form.
Such a conjecture has been made by Fargues and Scholze~\cite{fargues--scholze}.
Another such conjecture has also been made by
Zhu~\cite[\S 4.6, esp.\ Conj.~4.6.4 and Rem.~4.6.7]{zhu2020coherent}.

\begin{rem}\label{EHrem: formofsmoothfunctor}
Let $G={\rm GL}_d(F)$. In characteristic zero (over a fixed $p$-adic base field $L$) there should be a version of Conjecture \ref{conj:categorical local l not p} involving smooth $G$-representations on $L$-vector spaces and (ind-)coherent sheaves on the stack 
$$\mathcal{Y}_{G,F}=\bigcup_{F'}{\rm Mod}_{d,\varphi,N,F'/F}$$ of
Section~\ref{EHsub:sub:deRham} instead of the stack of $L$-parameters
$\cX_{G,F}$. Note that the stack of $L$-parameters for ${\rm GL}_d$ agrees with the moduli stack of $d$-dimensional Weil--Deligne representations ${\rm WD}_{d,F}$.
The point is that we expect ``local Langlands in families'' to make sense over $\mathcal{Y}_{G,F}$. Indeed, following \cite{MR3513570}, %
there is a family of representation $\mathcal{V}$ over $\cX_{G,F}$ constructed using a choice of a Whittaker datum $(N,\psi)$ and an isomorphism from the ring of invariant functions on $\cX_{G,F}$ to the Bernstein center $\mathcal{Z}$ of the category of smooth $G$-representations on $L$-vector spaces. 
Then $\mathcal{V}$ is the unique quotient of $$({c\text{-}\rm Ind}_N^G\psi)\otimes_{\mathcal{Z}}\mathcal{O}_{\cX_{G,F}}$$ with a certain prescribed quotient at the generic points of $\cX_{G,F}$. 
Using Proposition \ref{EH:propsameglobalsections} and Proposition \ref{EHBernsteiniso} the same construction should carry over to $\cY_{G,F}$. There are two difficulties in this strategy: we need to fix a Whittaker datum $(N,\psi)$, which is only possible if $L$ contains all the $p$-power roots of unity (and hence in particular $L$ can't be finite over $\Q_p$). Additionally, the stacks in question are not geometrically irreducible (and their connected components are not geometrically connected). Hence one needs to prove some rationality results in order to define the family $\mathcal{V}$ over a given base field $L$. \end{rem}

\section{Moduli stacks of global Langlands parameters and the
  cohomology of Shimura varieties}\label{sec: global stacks and cohomology of Shimura
varieties}%
In this section we briefly describe the global aspects of the theory.
We sketch the construction of moduli stacks of global Galois representations,
and of the associated  global-to-local restriction maps, and then explain
how these maps can be used, along with the conjectural 
categorical local Langlands correspondences, to give a conjectural
formula for the cohomology of Shimura varieties. We next explain some
results in the case of~$\GL_2/\Q$ obtained in~\cite{emertongeezhu},
which relates these conjectures to the
Taylor--Wiles method and to M.E.'s results on local-global compatibility
in the $p$-adic Langlands program~\cite{emerton2010local}.

\subsection{Global moduli stacks}
We begin with some basic facts regarding moduli stacks of global Langlands parameters.   
Although one can usually reduce to the case of considering reductive groups
defined over~$\Q$ (via an application of restriction of scalars),
it is no more difficult to consider the case of general number fields,  
and we do this.
Thus we fix a number  field $F$  and a connected  reductive group $G$  defined
over~$F$.  We also 
fix a finite set $S$  of places of $F$ containing all places above $p$ and~$\infty$.

We can define a derived stack $\cX_{G,F,S}$ which lies over $\Spf \cO$, and 
whose groupoid $\cX_{G,F,S}(A)$  of $A$-valued points,
for a classical $\cO$-algebra $A$ in which $p$ is nilpotent,
coincides with the groupoid of continuous representations
$\rho: \Gal_{F,S} \to {}^cG(A)$ (the target  being  given its discrete topology)
which are $C$-homomorphisms
(i.e.\ those which satisfy ~\eqref{eqn: C group hom condition}). 
(This defines the points of $\cX_{G,F,S}$ on classical rings,
i.e.\ it  characterizes the underlying classical stack $\cX^{cl}_{G,F,S}$;
the definition
of the points of $\cX_{G,F,S}$ on derived rings is more involved --- see 
\cite[\S 2.4]{zhu2020coherent} for the general formalism,
the  discussion at the end of \cite[\S 3.4]{zhu2020coherent} for 
a summary of the specifics in this context,
and forthcoming paper \cite{emertonzhu} for the precise details.)

In the case when $G = \GL_d$  for  some $d \geq 1,$
the underlying classical stack $\cX^{cl}_{\GL_d,F,S}$ has already been introduced by
Wang-Erickson \cite{MR3831282}.  He shows that it is a disjoint  union
of formal algebraic stacks, each of finite type over the formal Spec
of complete Noetherian local ring (in fact, an appropriate {\em pseudodeformation
ring}).  In general that  $\cX_{G,F,S}$  has an analogous 
structure,
and in particular
is a formal derived algebraic stack.  More precisely, if $\thetabar$
  is a residual pseudorepresentation of $\Gal_{F,S}$, then we let
  $X_{\thetabar}$ denote the associated formal deformation space of
  $\thetabar$.  Passing from a Galois representation to its associated
  pseudorepresentation gives rise to a morphism
  $\cX_{G,F,S} \to \coprod_{\thetabar} X_{\thetabar}$, and
  we let $\cX_{\thetabar}$ denote the preimage of~$X_{\thetabar}$.  We
  thus obtain a decomposition
\numequation\label{eqn: global stack pseudorep decomposition}\cX_{G,F,S} := \coprod_{\thetabar}\cX_{\thetabar},\end{equation}
  which is precisely the decomposition of $\cX_{G,F,S}$
  into its connected components.

Just as in the local case, if  $\rho: \Gal_{F,S} \to {}^cG(A)$ is an  $A$-valued point
of $\cX_{G,F,S}$, for some classical $\cO$-algebra~$A$ in which $p$  is nilpotent,
then the tangent complex to $\cX_{G,F,S}$ at~$\rho$ is computed via group cohomology,
in this  case the complex of continuous 
cochains~$C^{\bullet}(\Gal_{F,S}, \mathfrak g)[1]$.  
From now on we assume that  $p > 2$.  Then, 
just like in the  local case,
the cohomology of $\Gal_{F,S}$ vanishes in degrees~$>~2$,
and so we  find that $\cX_{G,F,S}$ is formally quasi-smooth over $\Spf \cO$.
However, unlike in the local case, $\cX_{G,F,S}$ need not be of its
expected dimension (the expected dimension at a particular point
being computable as the Euler characteristic of~$C^{\bullet}(\Gal_{F,S},\mathfrak g)[1]$,
which is given by the global Euler characteristic formula,
and so depends on the local behaviour at the real places of $F$
of the Galois representation corresponding to the point under consideration); more precisely, it can have components that are  of greater than
the expected  dimension.  Because it is formally quasi-smooth,
though, it {\em is} a derived local complete intersection,
i.e.\ can  be written as the quotient by an appropriately defined conjugation
action of ${}^cG\,\,\widehat{}$  (the  $p$-adic completion of ${}^cG$) 
 on a derived formal affine scheme of the  form
$$\Spf \cO[[X_1,\ldots,X_d]]/ (f_1,\ldots, f_r),$$
where $1 + d - r - \dim {}^cG$ is  equal to the expected dimension of~$\cX_{G,F,S}.$
(Here the first summand ``$1$'' corresponds to the dimension of~$\cO$.)
The relations $f_1,\ldots,f_r$ need not form a regular sequence, though,
and so the quotient  ring must be understood in a derived sense.  
The derived structure of $\cX_{G,F,S}$ is thus supported on those components
of $\cX^{cl}_{G,F,S}$  which are of greater-than-expected dimension.

\subsection{The global-to-local map}
\label{sec:localtoglobal}
The restriction of Galois representations is expected to give rise to  a morphism
\numequation\label{eqn: global to local restriction map of Galois
  stacks}f:  \cX_{G,F,S} \to \prod_{v \in S}
\cX_{G,F_v}.\end{equation} We only define this morphism when  $G
= \GL_d$ locally at places~$v|p$. %
Assuming that
this is the case, the  definition of this morphism is non-trivial, %
since the local stack $\cX_{G,F_v}$ %
is defined in terms of $(\varphi,\Gamma)$-modules rather than directly
in terms of local Galois representations. %
However, there is a
tautological restriction morphism to a derived version of Wang-Erickson's $p$-adic
local stacks~\cite{MR3831282}, and  the results of
~\cite{BIP} show that these stacks are actually classical (\emph{cf.}\
Remark~\ref{rem: expect that X is lci}); and these classical stacks are
substacks of our local stacks $\cX_{G,F_v}$ by Remark~\ref{rem: comparison to CWE stacks}.

In principle, we should also consider local stacks $\cX_{G,F_v}$ for
the real archimedean primes $v$ of~$F$.  However, since we are assuming
that $p  > 2$,  these will simply be  classical stacks~\cite[Prop.~2.3.2]{zhu2020coherent}, and classify conjugacy classes  of involutions %
lying above the
complex conjugation $c_v$ under the projection of ${}^cG$ onto its Galois factor.
In fact, we only want to consider a particular conjugacy class,
namely the ``odd'' class in ${}^cG$ --- see \cite[Prop.\ 6.1]{MR3028790}.
We let $\cX^{\odd}_{G,F,S}$ denote the union of connected components of 
$\cX_{G,F,S}$  classifying Galois representations  which are odd at every
real prime of~$F$.
\begin{rem}\label{rem: expect that the restriction morphism is representable}
  As will be explained in \cite{emertonzhu}, standard
  Fontaine--Mazur-type conjectures imply that the morphism $f$ is
  representable by algebraic (derived) stacks.
\end{rem}

\begin{rem}\label{rem: when is the global to local map quasismooth}
  The tangent complex of $f$ (i.e.\ the relative tangent complex of
  $\cX^{\odd}_{G,F,S}$ over $\prod_{v\in S} \cX_{G,F_v}$) at a point
  is computed (up to a shift) by the cone of the restriction morphism
  from global-to-local Galois cohomology; namely we have a
  distinguished triangle
$$\cdots \to \mathbf{t}_f \to C^{\bullet}(\Gal_{F,S},\mathfrak g)[1] \to 
\prod_{v \in S} C^{\bullet}(\Gamma_{F_v}, \mathfrak g)[1] \to \cdots
,$$ so that the relative tangent
complex %
$\mathbf{t}_f$ is  a Selmer-type complex of the sort considered
in \cite[App.\ B]{MR3762000}, which can possibly admit a non-zero $H^3$
term. This~$H^3$ term is in turn (by global duality) dual to another Galois
$H^0$, namely $H^0\bigl(\Gal_{F,S}, \widehat{\mathfrak g}(1)\bigr)$.  (This
local analysis of $f$ is a more geometric version of the traditional
method of presenting a global deformation ring or lifting ring over a
tensor product of local lifting rings. A similar analysis, in the more traditional language
of lifting rings, is made in~ \cite[\S4]{MR2459302}.  Such presentations are
important in modularity lifting theorems (and in their applications to
potential modularity); an
example in our earlier presentation is the morphism~\ref{eqn: restriction morphism with q added}.) %

The condition of $f$ being relatively quasi-smooth at $\rhobar$ is an
interesting one.  For example, if $G = \GL_d$, so that
$\widehat{\mathfrak g} = \End(\rhobar)$, then
$H^0\bigl(\Gal_{F,S}, \widehat{\mathfrak g}(1)\bigr)
= \End_{\Gal_{F,S}}\bigl(\rhobar, \rhobar(1)\bigr),$ and so we are
considering whether or not $\rhobar$ admits a non-trivial
$\Gal_{F,S}$-equivariant map to its Tate twist~$\rhobar(1)$.  If
$F = \Q$ and $G = \GL_2$, then (recalling $p > 2$) this is impossible
unless either $\rhobar^{\mathrm{ss}}$ is a twist of the direct sum
$1 \oplus \chibar$ (where $\chibar$ is the mod $p$ cyclotomic
character) or $p = 3$ and $\rhobar$ is induced from
$\Gal_{\Q(\sqrt{-3})}$. These are precisely the cases excluded by the
local-to-global principle for liftings of~$\rhobar$ of Diamond--Taylor
\cite[Thm.\ 1]{MR1272977}, and is known that such a local-to-global
principle fails in this case. In the case that $\rhobar$ is induced from
$\Gal_{\Q(\sqrt{-3})}$, this goes back to work of Carayol and Serre,
and an explicit example is given in~\cite[\S6.4]{Allen2019}. In the
reducible case the question of the possible levels has been
investigated by Ribet and Yoo, see in particular~\cite{MR3917209}. (See also \cite{MR3542487} for
analogous results in the case $p=2$.)

More generally, the condition that
$H^0\bigl(\Gal_{F,S}, \widehat{\mathfrak g}(1)\bigr) = 0$ is closely
related to the ``big image'' conditions in the Taylor--Wiles method,
and to the {\em genericity} condition of \cite{MR3702677}.
\end{rem}

\subsection{Computing the cohomology of Shimura varieties
via categorical Langlands}\label{subsec:computing-cohomology-Shimura}
We state a conjectural formula for the cohomology of Shimura varieties,
which is related to the various conjectures of \cite{zhu2020coherent},~\S 4.7.
(This conjecture will be discussed in more detail in~\cite{emertonzhu}.)
While the general conjecture seems far out of reach (since,
for example, it implies all kinds of automorphy theorems!), we discuss
the case of the modular curve in Section~\ref{subsec:
    GL2 Q}. %

\subsubsection{Cohomology of Shimura varieties}\label{subsubsec: cohomology of Shimura}
Suppose now that $F = \Q$,
and that the reductive group $G$ admits a Shimura datum.  We let
$F$ (which is now freed up as a piece of notation!) denote the reflex 
field of this Shimura datum, and let
$\mu: (\Gm)_{/F} \to G_{/F}$ denote the associated Hodge cocharacter
(well-defined up to conjugacy).
In the definition of ${}^cG$,  we choose the extension
$\tF$ of $\Q$ large  enough that it  contains~$F$,
and then write 
\[{}^cG_F :=  \widehat{G} \rtimes \bigl( \Gm \times \Gal(\tF/F) \bigr)\subseteq
\widehat{G} \rtimes \bigl( \Gm \times \Gal(\tF/\Q) \bigr) =: {}^cG.\]
The Hodge cocharacter~$\mu$, thought of as a weight of $\widehat{G}$,
gives rise to an irreducible (highest weight) representation $V_{\mu}$
of $\widehat{G}$, which naturally extends to a representation of~${}^cG_F.$
For any set $S$ of  primes (containing~$p$ and~$\infty$ as usual),
the moduli stack
$\cX_{G, \Q,S}$  parameterizes a universal family of (derived)
$C$-homomorphisms %
$\Gal_{\Q}  \to {}^cG$,
which we can restrict to~$\Gal_F$ so as to
obtain a family of representations $\Gal_F \to {}^cG_F.$
Composing this family with the representation $V_{\mu}$ of ${}^cG_F,$
we obtain a vector  bundle $\cV_{\mu}$  over  $\cX_{G,\Q,S}$, which is endowed
with an action of $\Gal_F$.

Choose a level subgroup~$K_f$, i.e.\ a compact open subgroup of~$G(\A)$,
which we suppose may be factored as
$K_f = \prod_v K_v,$ with  $K_v$  compact open in  $G(\Q_v)$ for each
finite place~$v$. 
We choose (as we may) $S$ large enough that $G$ is unramified outside~$S$,
and so that there is an extension $\cG$ of $G$ to a reductive group
scheme over $\Z[1/S]$ such that $K_v = \cG(\Z_v)$ for $v \not\in S$.
The locally symmetric quotient
$Y(K_f) := 
G(\Q)\backslash \bigl( X \times G(\A_f)\bigr)/K_f$
is a Shimura variety.

For each  finite prime $v \in S$, we let $W_v$ be a smooth representation of $K_v$ on
a finite torsion $\cO$-module. 
We let $W := \bigotimes_{v \in S} W_v$
denote the corresponding representation of $K_f$ (which acts through its projection
onto $\prod_{v \in S} K_v$; the tensor product is taken over~$\cO$, of course),
and we let $\cW$ denote the associated local system
$\cW := G(\Q)\backslash \bigl( X \times G(\A_f) \times W\bigr)/K_f$
over~$Y(K_f)$. 
For each~$v$, the compact induction $\cInd_{K_v}^{G(\Q_v)} W_v$
is a smooth $G(\Q_v)$-representation. 
Applying the functors of Conjectures~\ref{conj: Banach functor} and~\ref{conj:categorical local l not p},
we obtain (conjecturally!)
a coherent sheaf on~$\cX_{G,\Q_v}$, which we denote simply by~$\gA_v$.
We may form the exterior tensor product  $\gA := \boxtimes_{v \in S_f} \gA_v$
on $\prod_{v\in S_f} \cX_{G,\Q_v}.$

As above, let 
$f: \cX^{\odd}_{G,\Q,S} \to  \prod_{v \in S} \cX_{G, \Q_v}$ 
denote the global-to-local restriction morphism. 
We then have the following conjectural formula for the cohomology of the
Shimura variety. 
\begin{conj}
\label{conj:cohomology}Let~$d$ be the complex dimension of~$Y(K_f)$.
Then there is an isomorphism
\numequation\label{eqn: Shimura variety conjecture}
R\Gamma_c\bigl( Y(K_f), \cW)[d] \iso
R\Gamma( \cX^{\odd}_{G,\Q,S}, f^!\gA \otimes \cV_{\mu} )
\end{equation}
compatible with the $\Gal_F$- and {\em (}derived{\em )}
Hecke actions 
on source and target.
\end{conj}
By %
taking the inverse limit of the mod $p^n$ cases,
the conjecture naturally extends to the case of cohomology with $p$-adic
coefficients. We can also take a limit over the tame levels at places
in~$S$, to study completed cohomology, as in Expected
Theorem~\ref{expected EGZ thm} below. %
\begin{rem}
\label{rem:Whittaker independence}
As briefly mentioned in
our discussion of Conjecture~\ref{conj:categorical local l not p}
 in Subsection~\ref{subsec: categorical LL l not p},
the local functors 
appearing in the statement
of Conjecture~\ref{conj:cohomology}
depend on auxiliary choices of Whittaker-type data.  
In the statement of the conjecture,
we should make a global choice of this data,
and then compute the local functors~$\fA_v$ in terms of the induced
local data.  The resulting sheaf $f^!\fA$ should then
be well-determined independent of the global choice.
(A different choice of global data should cause each value of $\fA_v$
to be twisted by a certain line bundle, but the tensor product of the pullbacks
to $\cX^{\odd}_{G,\Q,S}$ of these various line bundles should  
then be trivial; see the discussion before~\cite[Conj.\ 4.7.5]{zhu2020coherent}
for the analogous discussion in the function field context.)
\end{rem}

\begin{rem}
  \label{rem: categorical implies FML}Conjecture~\ref{conj:cohomology}
  implies in particular that various representations of $\Gal_F$ occur
  in the cohomology of Shimura varieties, so in particular it implies
  cases of the Fontaine--Mazur conjecture.
  To make this implication precise, one needs to know something about the
supports of the sheaves $\mathfrak{A}_v$.   For $v|p,$
Conjecture~\ref{conj: Banach functor}~\eqref{item: loc alg vectors} is the relevant statement:
it shows that the $\mathfrak{A}_v$ will detect all points of $\cX^{\odd}_{G,\Q,S}$ whose
restrictions to $\cX_v$ satisfy an appropriate $p$-adic Hodge theoretic condition. 
Analogous statements in the case $v \nmid  p$,
describing the support of $\mathfrak A_v$ in terms of ramification conditions,
should also hold.  We don't spell these out here, but refer to 
the discussions of the $v\nmid p$  case in
\cite{benzvi2020coherent}, \cite{fargues--scholze}, \cite{hellmann2020derived},
and~\cite{zhu2020coherent}; the  essential point in this case is  
that the categorical theory should be compatible, in a suitable sense,
with the classical local  Langlands correspondence.
\end{rem}

\begin{rem}\label{rem: where does the Hecke action on Shimura
    varieties come from}
The Hecke algebras that act naturally on the source and
  target of~ \eqref{eqn: Shimura variety conjecture} are the derived Hecke algebras
  $\RHom_G( \cInd_{K_v}^{G(\Q_v)} W_v, \cInd_{K_v}^{G(\Q_v)}W_v )$ for
  finite places~$v$ (where we interpret $W_v$ simply as the trivial
  representation of $K_v$ in the case $v \not\in~S$; to see that the
  derived Hecke algebras do indeed act at places  $v \not\in~S$, one compares
  the situations at~$S$ and $S\cup\{v\}$ just after~\cite[Conj.\ 4.7.5]{zhu2020coherent}). %
  We refer
  the reader to \cite{zhu2020coherent} for a discussion of these
  derived Hecke algebras, and several conjectures related to them, in
  the case when $v \neq p.$
\end{rem}

\begin{rem}\label{rem: general locally symmetric spaces}
While the factor~$\cV_{\mu}$ has a global
  aspect, since it carries an action of the global absolute Galois
  group $\Gal_F$, if we ignore this then it can be regarded as the
  contribution from~$\infty$ to the formula
  of~Conjecture~\ref{conj:cohomology}.  We suspect that an analogous
  isomorphism should hold for {\em any} connected reductive group~$G$ (i.e.\ for
  the cohomology of all the locally symmetric spaces~$Y(K_f)$, whether
  or not they admit the structure of Shimura variety), with an
  appropriate choice of factor coming from~$\infty$ (which is now
  presumably just a vector bundle on~$\cX^{\odd}_{G,\Q,S}$, pulled
  back from an appropriately chosen vector bundle on the odd moduli
  stack at $\infty$, and not admitting the additional structure of a
  Galois representation).  %
 In particular, we expect to recover the derived phenomena observed and conjectured
by Galatius and Venkatesh  in~\cite{MR3762000, MR4061961}. \end{rem} %

  \subsection{An example in the global case: \texorpdfstring{$(\PGL_2)_{/\Q}$}{PGL2 over Q}}\label{subsec:
    GL2 Q}

In~\cite{emertongeezhu}, M.E., T.G.\ and Xinwen Zhu use the
Taylor--Wiles--Kisin patching method to study
Conjecture~\ref{conj:cohomology} in the case of the modular curve. We
briefly explain some of this work in progress in a special case, which
is adapted to our discussion of the patching method in
Section~\ref{sec: TW patching}; see also~\cite[\S 6]{johansson2024modulistacksgaloisrepresentations}.

\subsubsection{Relationship to patching}
  To fix ideas, set
$G = \PGL_2$, so we are in the modular curve
case. %
Choose~$\rbar:\Gal_{\Q,\{p,\infty\}}\to\GL_2(k)$ as in Section~\ref{sec: very
  brief introduction to patching}, so that in particular $\rbar|_{\Gal_{\Q(\zeta_p)}}$ is
absolutely irreducible, and~$\rbar$ has
determinant~$\varepsilonbar^{-1}$. Write $\rhobar:=\rbar|_{\Gal_{\Qp}}$. We consider ~\eqref{eqn: global to local restriction map of Galois
  stacks} in the special case that~$S=\{p,\infty\}$, so that we have
the
morphism \[f:\cX_{\GL_2,\Q,\{p,\infty\}}^{\varepsilon^{-1},\odd}\to\cX_{\GL_2,\Qp}^{\varepsilon^{-1}},\]where
we are only considering representations with determinant
$\varepsilon^{-1}$, so that the target is the stack denoted~$\cX$ in
Section~\ref{sec:banach-case-gl_2qp}. Write~$\cV$ for the rank 2 vector
bundle on $\cX_{\GL_2,\Q,\{p,\infty\}}^{\varepsilon^{-1},\odd}$ given
by the universal $\Gal_{\Q}$-representation.

As in Section~\ref{sec: very
  brief introduction to patching}, we write $R(\rhobar)$ for the
universal deformation ring of~$\rhobar$ with fixed determinant~$\varepsilon^{-1}$. By definition there is  a versal morphism
$g: \Spf R(\rbar_p) \to \cX_p$ that maps
the closed point in its domain to the local restriction $f(\rbar)$. 
The patching construction gives  an $R(\rbar|_{\Gal_{\Qp}})$-module $M_{\infty}$ with
an action of~$\PGL_2(\Qp)$.

We expect to prove the following result
in~\cite{emertongeezhu}, describing the localization at the maximal
ideal $\m$ of the prime-to-$p$ Hecke algebra corresponding
to~$\rbar$ of the completed homology group $\widetilde{H}_{1}$. As
in~\eqref{eqn: global stack pseudorep decomposition} we write
$\cX_{\rbar}$ for the irreducible component of
$\cX^{\varepsilon^{-1},\odd}_{\GL_2,\Q,\{p,\infty\}}$ corresponding to~$\rbar$. %

\begin{expectedthm}[M.E.--T.G.--Zhu]\label{expected EGZ thm}Let
  $L_\infty$ be the sheaf of $\PGL_2(\Qp)$-representations of Expected
  Theorem~\ref{expectedthm: DEG results}. Then there is a
  $\PGL_2(\Qp)$-equivariant isomorphism\numequation\label{eqn: modular
    curve cohomology from L infty}
\widetilde{H}_{1,\m}\iso
R\Gamma( \cX_{\rbar}, f^!L_\infty \otimes \cV)[-1]
\end{equation}
compatible with the actions of $\Gal_\Q$ and the \emph{(}derived\emph{)} Hecke
algebras at primes $\ell\ne p$ on each side.

  In addition, there is a $\PGL_2(\Qp)$-equivariant isomorphism
  $M_\infty\isoto g^*L_\infty$.  
\end{expectedthm}
\begin{rem}\label{rem: hope to improve on EGZ}
  We anticipate that the hypotheses in Expected Theorem~\ref{expected
    EGZ thm} can be considerably relaxed. In particular, we can allow
  $\rbar$ to be at ramified at places not dividing~$p$, and we hope
  to relax the other assumptions on~$\rbar$ as much as possible in
  the final version of~\cite{emertongeezhu}. 
\end{rem}

We expect to prove Expected Theorem~\ref{expected EGZ thm} as
follows. The basic idea is that $M_{\infty}$ is ``patched'' out of the
cohomology of modular curves, and we want to show that the patching
process ``undoes'' the $f^!$ appearing in the conjectural formula for
cohomology. %
We show (using existing
modularity lifting theorems) that $\cX_{\rbar}$ is equal to
$[\Spf R_{\{p,\infty\}}(\rbar)/\boldsymbol{\mu}_2]$. We then show that 
$R_{\{p,\infty\}}(\rbar)$ is the quotient of~$R(\rhobar)$ by a certain regular sequence which is also a regular sequence 
on~$M_\infty$, in order to compute the pullback~$f^!$ via a Koszul complex, and ultimately
reduce~\eqref{eqn: modular
    curve cohomology from L infty} to the local-global compatibility
  results originally proved by M.E.\ in~\cite{emerton2010local}, or
  more precisely to their interpretation and reproof in~\cite{MR3732208}.

  \subsection{The structure of \texorpdfstring{$\cX^{\odd}_{\GL_2,\Q,S}$}{the odd stack over Q}}
Another aim of~\cite{emertongeezhu} is to use Taylor--Wiles patching
and Iwasawa theory to describe the structure of the
  stacks~$\cX^{\odd}_{\GL_2,\Q,S}$. We briefly describe some of our
  expected results. For simplicity we again restrict to the stacks
  ~$\cX^{\varepsilon^{-1},\odd}_{\GL_2,\Q,S}$ with determinant $\varepsilon^{-1}$.

Similar to (\ref{eqn: global stack pseudorep decomposition}) we have the decomposition 
  \[\cX^{\odd}_{\GL_2,\Q,S} := \coprod_{\thetabar}\cX_{\thetabar}\]
  into connected components, where $\thetabar$ runs over the  odd residual pseudorepresentation of
$\Gal_{\Q,S}$ with determinant $\varepsilonbar^{-1}$.  The simplest case is when $\thetabar$ arises from an irreducible odd
  representation $\rbar: \Gal_{\Q,S} \to \GL_2(\Fbar_p)$.  We then
  write $X_{\thetabar} = X_{\rbar} = \Spf R_{\rbar}$, where
  $R_{\rbar}$ is the deformation ring of $\rbar$, and we find that
  $\cX_{\thetabar} = \cX_{\rbar} = [\Spf R_{\rbar}/\Gmhat],$ where
  $\Gmhat$ denotes the $\mathfrak m_{R_{\rbar}}$-adic completion of
  $\Gm$, acting trivially on $\Spf R_{\rbar}$.  In this case, then,
  there is no essential difference between the theory of the moduli
  stack~$\cX_{\rbar}$ and usual Galois deformation theory.  Serre's
  conjecture tells us that $\rbar$ is modular,
  and then known ``big $R = $ big $\T $'' theorems imply (at least
  under the usual Taylor--Wiles hypotheses on~$\rbar$) that
  $R_{\rbar}$ may be identified with an appropriately defined Hecke
  algebra acting on a space of $p$-adic modular forms. %

  The case when $\thetabar$ is reducible is more subtle, and our goal
  in this case is again to describe $\cX_{\thetabar}$ %
  as much as possible in terms of modular forms.  To give
  a sense of what this means, set $\thetabar = 1 + \varepsilonbar^{-1},$ %
  and let~$\cY_{\thetabar}$ be the (derived) substack obtained by  imposing the local conditions of
  being unipotently tame at the primes
  $\ell \in S \setminus \{p, \infty\}$, as well as the local condition
  of being finite flat at~$p$.  %
  Let
  $N = \prod_{\ell \in S \setminus \{p,\infty\}}\ell$, and choose
  $\cO$ large enough that each cuspidal eigenform $f$ of weight two
  and level $N$ has coefficients in $\cO$.  Then each such cuspform
  $f$ whose associated mod $p$ pseudorepresentation $\thetabar_f$ is
  equal to $1+\varepsilonbar^{-1}$ contributes an irreducible component to
  $\cY_{\thetabar}$ of the form
  $[\Spf \cO\langle X,Y\rangle /(XY - c) / \Gmhat],$ where $c$ is a
  certain modulus of reducibility, $\cO\langle X,Y\rangle $ denotes
  the $p$-adic completion of $\cO[X,Y]$, and $\Gmhat$ denotes the
  $p$-adic completion of $\Gm$, acting via
  $t\cdot (X,Y) = (t^2 X, t^{-2} Y).$ This component parameterizes
  lattices in the $p$-adic Galois representation attached to~$f$.

  In addition to such components, there is an ``Eisenstein''
  component, whose underlying classical stack is of the form
  $[\Spf \cO\langle X_1,\ldots, X_n \rangle /\Gmhat],$ where $n$
  denotes the order of $S \setminus \{p,\infty\}$ (i.e.\ the number of
  factors of~$N$), and $\Gmhat$ acts via $t\cdot X_i = t^2 X_i$.  This
  component parameterizes extensions of $\varepsilon^{-1}$~by~$1$; the space of
  such extensions is spanned by the Kummer classes associated to the
  $n$ factors of~$N$.  If $n > 1$, then this component is of
  larger-than-expected dimension, and so (following the discussion
  above) will carry non-trivial derived structure.

  There can be yet another irreducible component of $\cX$,
  parameterizing extensions of $1$ by~$\varepsilon^{-1}$ (i.e.\ extensions ``the
  wrong way'').  Assuming $p \geq 5$ for simplicity, these are
  governed by those $\ell \in S$ for which $\ell \equiv -1 \bmod p$.
  If there are $m$ such~$\ell$ , then the family of such extensions is
  $(m-1)$-dimensional (there is an $m$-dimensional space of extensions
  in characteristic~$p$, which is then quotiented by a
  $\Gmhat$-action) and is supported (set-theoretically) purely in
  characteristic~$p$.  If $m > 1$ it is a component of
  $\cY_{\thetabar}$.  If $m > 2$ it is also of greater-than-expected
  dimension, and so admits derived structure.

  We expect  that the preceding discussion enumerates all
  possible components of~$\cY_{\thetabar}$. %
  We hope that this (and more general statements) can be proved by
  suitable applications of (known cases of the) the Fontaine--Mazur
  conjecture (especially in the residually reducible case, where it
  has been proved by Skinner--Wiles~\cite{SkinnerWiles} and
  Pan~\cite{pan2019fontainemazur}), together with the main conjecture
  of Iwasawa theory \cite{MR742853}.

\subsection{Eigenvarieties and sheaves of  overconvergent \texorpdfstring{$p$}{p}-adic automorphic forms}
\label{sec:eigenvar}
We now  state an analogue of Conjecture \ref{conj:cohomology}, giving a conjectural description of the cohomology of some (overconvergent) $p$-adic coefficient systems. This gives a conjectural description of the sheaves of overconvergent $p$-adic automorphic forms that appear on eigenvarieties.
We maintain the notation from the preceding section.

There are basically two approaches to $p$-adic automorphic forms and
eigenvarieties: one can use the cohomology of locally symmetric spaces
with coefficients in $p$-adic coefficient systems (see
e.g.~\cite{MR2075765} or \cite{MR3692014}), or one can use the
completed cohomology groups introduced by M.E.\ and further developed by M.E.\ and Frank Calegari
(see \cite{MR2207783,MR3728617,MR2905536}) and apply the locally analytic Jacquet functor introduced by M.E.~in \cite{MR2292633} (see \cite{MR2207783} for a construction of eigenvarieties using this approach).
We briefly discuss the relation between the two approaches and their relation with the functor $\mathfrak{A}_G^{\rm rig}$ of Conjecture \ref{EH:conjanalytic}.

\begin{rem}\label{EHrem:eigenvarforBernsteinblocks}
In general an eigenvariety can be associated with a tame level~ $K^p$,
Bernstein blocks $\Omega_v=[M_v,\sigma_v]$ of ${\rm Rep}_L^{\rm
  sm}G(F_v)$ at places $v|p$, and  algebraic representations $W_v$
(defined up to twist with an algebraic character of $M_v$) of the
supercuspidal support $M_v$ of $\Omega_v$. We will focus on the
classical case where all the Bernstein blocks are principal (and hence
$W_v$ is not an extra datum) and  refer to \cite{BreuilDingeigenvar}
for the treatment of arbitrary Bernstein components. We point out that
in this more general context it is also possible to state similar conjectures to Conjectures \ref{EHconj:localglobal1} and \ref{EHconj:localglobal2} below.
\end{rem}

\subsubsection{The case of Shimura sets}\label{subsec: eigenvarcompactatinfty} Let us first discuss the special case where the reductive group $G$ is compact at infinite places, and hence the locally symmetric spaces $Y(K_f)$ are just point sets. We refer to these sets as Shimura sets even though strictly speaking they are not Shimura varieties. 
In order to fix the setup we assume from now on that $F$ is totally real and the group $G$ is
\begin{enumerate}
\item[-] either a definite unitary group associated to a hermitian space over a CM extension $E$ of $F$,
\item[-] or the group associated to a definite quaternion algebra over $F$,
\end{enumerate}
and that $G$ is split at places dividing $p$.

 Let us fix for each place $v|p$ a split maximal torus and a Borel $T_v\subset B_v\subset G_v=G(F_v)={\rm GL}_{d}(F_v)$. We also fix the choice of an Iwahori subgroup $I_v$ adapted to $B_v$ and write $T_v^0=T_v\cap I_v$ for the maximal compact subgroup of the torus $T_v$, and $B_{v,0}=B_v\cap I_v$. Moreover, we fix a \emph{tame level} $K^p$ which is a compact open subgroup of $G(\mathbf{A}_F^{\infty,p})$.
Let $\kappa_v:T_v^0\rightarrow L^\times$ be a locally analytic character. Then we can define the following space of \emph{module valued overconvergent $p$-adic automorphic forms}
\numequation\label{EHeqn:modulevaluedpadicforms}
M^\dagger(G,K^p,(\kappa_v)_{v|p})=\left\{\begin{array}{*{20}c} f:G(F)\backslash G(\mathbf{A}_F^\infty)/K^p\rightarrow \widehat\bigotimes_{v|p}({\rm Ind}_{B_{v,0}}^{I_v}\kappa_v)^{\rm an}\\\text{s.~th.}\ f(gu)=u^{-1}f(g)\ \text{for all}\ g\in G(\mathbf{A}_F^\infty), u\in \prod_{v|p}I_v\end{array}\right\}
\end{equation}
compare e.g.~\cite[3.3]{MR2769113}. This space comes equipped with an
action of the so-called Atkin--Lehner ring $\mathcal{A}^+(G_v)$ whose
definition (and action) we now briefly sketch.

For each place $v|p$ the choice of the Borel subgroup $B_v$ and its subgroup $N_{v,0}=N_v\cap B_{v,0}$ cuts out a submonoid $T_v^+\subset T_v$ consisting of all $t\in T_v$ such that $t N_{v,0} t^{-1}\subset N_{v,0}$. Clearly $T_v^+$ contains the maximal compact subgroup $T_v^0$. We fix the choice of a free abelian subgroup $\Sigma_v\subset T_v$ such that $\Sigma_v\cong T_v/T_v^0$ (in the case of the diagonal torus and after choosing a uniformizer $\varpi_v$ we might for example choose the diagonal matrices all of whose entries are powers of $\varpi_v$) and set $\Sigma_v^+=\Sigma_v\cap T_v^+ $, which is a free abelian monoid. Then we define 
\begin{align*}
\mathcal{A}^+(G_v)&=L[\Sigma_v^+]\\
\mathcal{A}^+(G)&=\bigotimes_{v|p}\mathcal{A}^+(G_v).
\end{align*}
It turns out \cite[3.4]{MR2769113} that $\mathcal{A}^+(G_v)$ is canonically isomorphic to the subalgebra of the Hecke algebra $\mathcal{H}(G_v)$ consisting of functions that are $I_v$-biinvariant and that are supported on the submonoid $\Sigma_v^+I_v$ of $G_v$ generated by $\Sigma_v^+$ and $I_v$. 
The action of $I_v$ on ${\rm Ind}_{B_{v,0}}^{I_v}\kappa_v$ canonically extends to $\Sigma^+_vI_v$ \cite[Theorem 2.4.7]{MR2769113} and hence this monoid acts on $M^\dagger(G,K^p,(\kappa_v)_{v|p})$ via $(\gamma.f)=\gamma f(g\gamma)$.
Viewing elements of $\mathcal{A}^+(G)$ as locally constant functions on $G(\prod_{v|p}F_v)$ we hence obtain an action of $\mathcal{A}^+(G)$ on $M^\dagger(G,K^p,(\kappa_v)_{v|p})$.

We now explain a different characterization of the above space of overconvergent $p$-adic automorphic forms together with its $\mathcal{A}^+(G)$-action that involves completed cohomology and is in fact less involved to write down (in particular it is easier to write down the action of the Atkin--Lehner ring in this setup). We however preferred to give the above definition first as it is directly related to the cohomology with coefficients in a $p$-adic coefficient system defined by the representations ${\rm Ind}_{B_{v,0}}^{I_v}\kappa_v$ (which makes it easier to formulate Conjecture \ref{EHconj:localglobal1} below).

Consider the space
$$\hat S(K^p,L)=\{f:G(F)\backslash G(\mathbf{A}_F^{\infty})/K^p\rightarrow L\ \text{continuous}\}$$
which is an admissible Banach space representation of $G(F_p)$ and coincides with the completed cohomology (base changed to $L$) of $Y(K^p)$ as our group is compact at infinity. 
The locally analytic vectors in this Banach space representations are given by  
\numequation\label{EHeqn:locanfcts}\hat S(K^p,L)^{\rm an}=\{f:G(F)\backslash G(\mathbf{A}_F^{\infty})/K^p\rightarrow L\ \text{locally analytic}\}.
\end{equation}
As above we write $N_v \subset B_v$ for the unipotent radical of $B_v$ and $N_{v,0}=N_v\cap B_{v,0}$. 
If we set $N_0=\prod_{v|p}N_{v,0}$ we can then define an action of $T^+$ on the $N_0$-invariants $\hat S(K^p,L)^{N_0}$ by 
\numequation\label{EHeqn defofT+action} 
t.f=\sum_{n\in N_0/t N_0t^{-1}} ntf. 
\end{equation}
The following proposition compares the modules \eqref{EHeqn:modulevaluedpadicforms} and \eqref{EHeqn:locanfcts} equipped with their $\mathcal{A}^+(G)$-action (where $\mathcal{A}^+(G_v)$ acts on the $N_0$-invariants of \eqref{EHeqn:locanfcts} via the inclusion $L[\Sigma_v^+]\subset L[T_v^+]$).
\begin{prop}\label{EHprop:compareocandcompletedcohom1}
Let $N_0=\prod_v N_{v,0}$. Then there is an isomorphism 
$$M^\dagger(G,K^p,(\kappa_v)_{v|p})\isoto (\hat S(K^p,L)^{\rm an})^{N_0}[(\kappa_v)_{v|p}]$$
equivariant for the action of the Atkin--Lehner ring $\mathcal{A}^+(G)=\bigotimes_{v|p}\mathcal{A}^+(G_v)$ and for the Hecke-action away from $p$. Here the $(-)[(\kappa_v)_{v|p}]$ on the right hand side means the subspace of the $N_0$-invariants on which $T^0=\prod_v T_v^0$ acts via $\prod_v \kappa_v$.
\end{prop}
\begin{proof}
This is for example proven in \cite[Prop. 3.10.1, Prop. 3.10.2]{MR2769113}. \end{proof}

Instead of a fixed locally analytic character one can apply the same construction to a family of characters of $T^0$ or even to the universal character 
$$\kappa^{\rm univ}:T^0\rightarrow \Gamma(\hat T^0,\cO_{\hat T^0})^\times,$$
where $\hat T^0$ denotes the space of continuous characters of $T^0$. 
This means that in the definition of $M^\dagger(G,K^p,\kappa^{\rm univ})$ we replace the representation ${\rm Ind}_{B_{v,0}}^{I_v}\kappa_v$ by the $I_v$-representation ${\rm Ind}_{N_{v,0}}^{I_{v}}\mathbf{1}$ (which contains all the representations ${\rm Ind}_{B_{v,0}}^{I_v}\kappa_v$ for the various characters $\kappa_v$). 
\begin{prop}\label{EHprop: comparecompletedcohom}
There is an isomorphism
\[M^\dagger(G,K^p,\kappa^{\rm univ})\isoto (\hat S(K^p,L)^{\rm an})^{N_0},\]
equivariant for the action of $L[T^+]$ and for the action of the Hecke-operators away from $p$.
\end{prop}
\begin{proof}
This is a special case of \cite[Theorem 1.2]{Fu2021derivedeigenvar};
see also the discussion of~\cite[\S 3.2]{MR2207783}.
\end{proof}

We now conjecturally relate these spaces of overconvergent $p$-adic automorphic forms to the functor $\mathfrak{A}_G^{\rm rig}$.
On the Galois side, let us write $\mathfrak{X}_{G,F,S}=\mathcal{X}_{G,F,S}^{\rm rig}$ for
the rigid analytic generic fiber of the stack 
$\mathcal{X}_{G,F,S}$. Note that $\mathcal{X}_{G,F,S}$ in general is a derived stack (not a classical stack), and we hence have to adapt the definition of its generic fiber (which will be a derived rigid analytic Artin stack), generalizing (\ref{EHdef:genfibofstack}) to the derived setting. As we do not want to discuss derived rigid analytic geometry here, we do not spell out the precise definition of $\mathfrak{X}_{G,F,S}$.

We note that similarly to \eqref{eqn: global to local restriction map of Galois stacks}, given $v|p$ there is a morphism
$$f_{{\rm rig},v}:  \mathfrak{X}_{G,F,S}=(\cX_{G,F,S})_\eta^{\rm rig} \to  \mathfrak{X}_{G,F_v},$$
that is given by composing the corresponding map to
$(\cX_{G,F_v})^{\rm rig}_{\eta}$ with the morphism
\eqref{EHgenfibEGstacktoanalytic}.
 We write 
$$f_{\rm rig}: \mathfrak{X}_{G,F,S} \to  \prod_v \mathfrak{X}_{G,F_v}$$
for the product of the maps $f_{{\rm rig},v}$. By abuse of notation we also write $f_{\rm rig}$ for the restriction of $f_{\rm rig }$ to $\mathfrak{X}_{G,F,S}^{\rm odd}=(\mathcal{X}_{G,F,S}^{\rm odd})_\eta^{\rm rig}\subset \mathfrak{X}_{G,F,S}$.

Let us write $\sigma_{\kappa_v}={\rm Ind}_{B_{v,0}}^{I_v}\kappa_v$ and $\sigma_v={\rm Ind}_{N_{v,0}}^{I_v}\mathbf{1}$.
We can use these representation and the functors $\mathfrak{A}_{G(F_v)}^{\rm rig}$ to define the coherent sheaves
\numequation\label{EHeqn:productoflocalsheaves}
\begin{aligned}
\mathfrak{A}_{\rm
  rig}(({\kappa}_v)_{v|p})&=\big(\bigotimes_{v\nmid p} \mathfrak{A}_{G(F_v)}(\cInd_{K_v}^{G(F_v)}\mathbf{1})\big)^{\rm rig} \otimes \bigotimes \mathfrak{A}_{G(F_v)}^{\rm rig}(\cInd_{I_v}^{G(F_v)}\sigma_{\kappa_v}),\\
 \mathfrak{A}_{\rm
  rig}^{\rm univ}&=\big(\bigotimes_{v\nmid p} \mathfrak{A}_{G(F_v)}(\cInd_{K_v}^{G(F_v)}\mathbf{1})\big)^{\rm rig} \otimes \bigotimes \mathfrak{A}_{G(F_v)}^{\rm rig}(\cInd_{I_v}^{G(F_v)}\sigma_{v}). 
  \end{aligned}
\end{equation}

\begin{conj}\label{EHconj:localglobal1}
There are isomorphisms 
$$
\begin{aligned}
M^\dagger(G,K^p,(\kappa_v)_{v|p})&\cong R\Gamma_c(\mathfrak{X}_{G,F,S}^{\rm odd},f_{\rm rig}^*\mathfrak{A}_{\rm rig}(({\kappa}_v)_{v|p}),\\
M^\dagger(G,K^p,\kappa^{\rm univ})&\cong R\Gamma_c(\mathfrak{X}_{G,F,S}^{\rm odd},f_{\rm rig}^*\mathfrak{A}_{\rm rig}^{\rm univ})
\end{aligned}
$$
compatible with the Hecke action and the action of the Atkin--Lehner ring.
\end{conj}
\begin{rem}
The Atkin--Lehner ring acts on the right hand side via the canonical action of $L[T^+]$ on $\cInd_{I_v}^{G_v}\sigma_{\kappa_v}$, see for example \cite[Lemma 2.2, Lemma 2.3]{MR2910788}.
\end{rem}
We point out that this conjecture is somehow dual in nature to Conjecture \ref{conj:cohomology}: we use cohomology on the side of the symmetric space (which is just a point set in this case) and cohomology with compact support on the side of coherent sheaves (note that the stack $\mathfrak{X}_{G,F,S}^{\rm odd}$ is usually the quotient of a Stein space by the action of a reductive group and hence we can define cohomology with compact support as in \cite{MR1094850}, similar to its use in Section~\ref{sec: GL1}). This duality also explains why in this conjecture we use the $\ast$-pullback along $f_{\rm rig}$ instead of the $!$-pullback of Conjecture \ref{conj:cohomology}.

\begin{rem}
In the definition of $M^\dagger(G,K^p,(\kappa_v)_{v|p})$ we could also
have inserted non-trivial $K_v$-representations $\sigma_v$ for some
places $v\nmid p$ as coefficients. 
We then expect the same isomorphism as in Conjecture \ref{EHconj:localglobal1}, but  then of course the trivial $K_v$-representation $\mathbf{1}$ in (\ref{EHeqn:productoflocalsheaves}) has to be replaced by $\sigma_v$.\\
\end{rem}

\subsubsection{Finite slope spaces and Emerton's locally analytic Jacquet functor} We recall the definition of finite slope spaces and of the locally analytic Jacquet module \cite{MR2292633}. 

The spaces $M^\dagger(G,K^p, (\kappa_v)_{v|p})$ and $(\hat S(K^p,L)^{\rm an})^{N_0}$ are locally convex topological $L$-vector spaces (of compact type) that come equipped with an action of the monoid $L[T^+]$. 
The construction of a finite slope space $V_{\rm fs}$ attached to a locally convex $L$-vector space $V$ equipped with a topological $T^+$-action (or a $T_v^+$-action for a place $v|p$) is a canonical procedure to extend the $L[T^+]$-action to a locally analytic representation of $T\supset T^+$. Its formal definition is given by 
$$V_{\rm fs}={\rm Hom}_{L[T^+],b}(\Gamma(\hat T,\mathcal{O}_{\hat T}), V)$$
compare \cite[Definition 3.2.1]{MR2292633}. Here $\hat T$ is the rigid analytic space of characters of $T$ (note that $\Gamma(\hat T,\cO_{\hat T})$ contains the distribution algebra $\cD(T)$ but is strictly larger) and the subscript $b$ indicates that we take the space of continuous homomorphisms and equip it with the strong topology. Then the $\Gamma(\hat T,\mathcal{O}_T)$-module structure on $V_{\rm fs}$ restricts to an action of $\mathcal{D}(T)$ and hence to an action of $T$ which in fact is locally analytic \cite[Proposition 3.2.4]{MR2292633}. If this representation is (essentially) admissible, the strong dual of the finite slope space is a coadmissible module over the Frechet--Stein algebra $\Gamma(\hat T,\cO_{\hat T})$. As $\hat T$ is a Stein space this implies that there exists a (necessarily unique) coherent sheaf $\mathcal{M}^\ast_V$ on $\hat T$ such that 
\numequation\label{EHeqn:dual of Vfs}
(V_{\rm fs})'\cong\Gamma(\hat T,\mathcal{M}^\ast_V)
\end{equation}
as $\Gamma(\hat T,\mathcal{O}_{\hat T})$-modules. 
\begin{rem}\label{EHrem Jacquetandsheaves}
\hfill\\(i) If $V=\Pi^{N_{v,0}}$ for a locally analytic $G_v$-representation
$\Pi$, then $V$ is equipped with a $T^+$-action defined by (\ref{EHeqn
  defofT+action}). Then essential admissibility of $V_{\rm fs}$ (as a
locally analytic $T_v$-representation) can for example be assured by
assuming that $\Pi|_H=\mathcal{C}^{\rm an}(H,L)^m$ for some open
compact subgroup $H\subset G_v$ and some $m>0$, see \cite[Proposition
3.2.24]{MR2292633} and \cite[Proposition 5.3]{MR3623233} (which relies on \cite{MR2292633}). In this
situation \cite[Proposition 5.3]{MR3623233} moreover shows that locally on $\hat T_v$ the sheaf $\mathcal{M}^\ast_V$ is finite projective over $\hat T_v^0$. \\
(ii) Assume that the $T$-representation on $V_{\rm fs}$ is essentially admissible.
Then using Serre duality on the Stein space $\hat T$ (as done for the open unit disc $\mathbf{U}$ in \eqref{EHeqn SerredualityStein}) we can rewrite (\ref{EHeqn:dual of Vfs}) as
$$V_{\rm fs}=R\Gamma_c(\hat T, \mathcal{M}_V),$$
where $\mathcal{M}_V=\mathbf{D}(\mathcal{M}^\ast_V)$ is the Serre dual (on $\hat T$) to $\mathcal{M}^\ast_V$. A priori this is only a complex of coherent sheaves. However, in good situations (e.g.~in the situation of (i)) the sheaf $\mathcal{M}^\ast_V$ turns out to be a Cohen--Macaulay module and hence $\mathcal{M}_V$ is concentrated in one degree (though this degree might be strictly negative, depending on the dimension of the support of $\mathcal{M}^\ast_V$). \end{rem}

\begin{defn}
The locally analytic Jaquet module $J_{B_v}(\Pi_v)$ of a locally analytic $G_v$-representation $\Pi_v$ is the $T_v$-representation on $(\Pi_v^{N_{v,0}})_{\rm fs}$.
\end{defn}
With this definition at hand Proposition \ref{EHprop: comparecompletedcohom} directly implies the following:

\begin{prop} There is an isomorphism 
\[M^\dagger(G,K^p,\kappa^{\rm univ})_{\rm fs}\xrightarrow\cong J_B(\hat S(K^p,L)^{\rm an})\]
of locally analytic $T$-representations.
\end{prop}

If we restrict attention to the finite slope subspace of $M^\dagger(G,K^p,\kappa^{\rm univ})$, then we can formulate a conjecture  that gives a precise description of the coherent sheaf whose compactly supported cohomology conjecturally computes $M^\dagger(G,K^p,\kappa^{\rm univ})_{\rm fs}$ and that does not involve the conjectural functors $\mathfrak{A}_G^{\rm rig}$ (though of course the conjecture is inspired by these conjectural functors and their expected properties with respect to parabolic induction).  

As above, for $v\nmid p$, let us write $\mathfrak{A}_{{\rm rig},v}=(\mathfrak{A}_{G(F_v)}(\cInd_{K_v}^{G_v}\mathbf{1}))^{\rm rig}$. For $v|p$ consider the morphism
$$f_{{\rm rig},v}\times {\iota_v}:\mathfrak{X}_{G,F,S}^{\rm odd}\times \hat T\longrightarrow \mathfrak{X}_{G,F_v}\times\hat T_v,$$
where $\iota_v$ is given by projection to the $v$-th factor followed by the automorphism 
$$(\delta_1,\dots,\delta_n)\mapsto \delta_{B_v}\cdot (\delta_1,\delta_2 (\varepsilon\circ {\rm rec}_{F_v}),\dots, \delta_n (\varepsilon\circ{\rm rec}_{F_v})^{n-1})$$
of $\hat T_v$, where $\delta_{B,v}$ denotes the modulus character.
Moreover consider the diagram
$$
\begin{xy}
\xymatrix{
& \overline{\mathfrak{X}}_{B_v} \ar[ld]_{\pi_v}\ar[rd]^{\alpha_v} & \\
\mathfrak{X}_{G,F_v}\times\hat T_v && \mathfrak{X}_{T_v}.
}
\end{xy}
$$
and define $$\mathfrak{A}_{{\rm rig},v}=\pi_{v,\ast}\alpha_v^\ast(\mathcal{O}_{\mathfrak{X}_{T_v}}([F_v:\Q_p](0,-1,\dots,-n+1))[{\rm rk}_{\Q_p}\, T_v]).$$ 
Here, as in Conjecture \ref{EH:conjanalytic} (2), the twist $(0,-1,\dots,-n+1)$ means that we take a sheaf on $\hat T_v$ and consider it as a sheaf on $\mathfrak{X}_{T_v}=\hat T_v/\mathbf{G}_m^n$ by equipping it with the $\mathbf{G}_m^n$-equivariant structure given by the algebraic character $(0,-1,\dots, -n+1)$.
\begin{rem}
We note that it would be more canonical to talk about ${\pi_{v,!}}$ instead of $\pi_{v,\ast}$. But as $\pi_v$ is proper by Theorem \ref{EHThm propernessofcompactific} these two functors agree and we do not need to worry about the definition of $\pi_{v,!}$ (that most likely would involve condensed structures).
\end{rem}
Writing $$g_{\rm rig}=\prod_{v\not |p} f_{{\rm rig},v}\times \prod_{v|p} g_{{\rm rig},v}:\mathfrak{X}_{G,F,S}^{\rm odd}\times \hat T\rightarrow \prod_{v}\mathfrak{X}_{G,F_v}\times\hat T$$
we can formulate the following conjecture:
\begin{conj}\label{EHconj: coherentdescriptionoffiniteslop}
There is an isomorphism
$$M^\dagger(G,K^p,\kappa^{\rm univ})_{\rm fs}\cong R\Gamma_c(\mathfrak{X}_{G,F,S}\times\hat T,g_{{\rm rig}}^\ast \big(\bigotimes_v \mathfrak{A}_{{\rm rig},v}\big)).$$
equivariant for the action of the Hecke-operators away from $p$ and for the action of the Atkin--Lehner ring at the places dividing $p$.
\end{conj}
\begin{rem}
Note that the degree shift in the definition of $\mathfrak{A}_{{\rm rig}, v}$ is necessary: it is expected that the coherent sheaf in Conjecture \ref{EHconj: coherentdescriptionoffiniteslop} has support of dimension $\sum_{v|p}{\rm rk}_{\Q_p}\, T_v$. Hence with this degree shift we should expect that the compactly supported cohomology sits in degree $0$. This should also explain why we actually need some degree shift in the analytic case of Section~\ref{sec: GL1}.  
\end{rem}

\subsubsection{Description of Hecke eigenspaces of $p$-adic automorphic forms}
 We assume that $\Pi_v$ is a $G_v$-representation satisfying the
 assumptions in Remark \ref{EHrem Jacquetandsheaves}, that is we assume that $\Pi_v|_{H_v}$ is isomorphic to $\mathcal{C}^{\rm an}(H_v,L)^m$ for some $m> 0$ and some compact open subgroup $H_v\subset G_v$.
Given a locally analytic character $\delta$ of $T$ we can compute the eigenspace $${\rm Hom}_T(\delta, J_{B_v}(\Pi_v))={\rm Hom}_{T^+}(\delta,\Pi_v^{N_{v,0}})$$ using the coherent sheaf $\mathcal{M}^\ast_{\Pi_v}:=\mathcal{M}^\ast_{\Pi_v^{N_{v,0}}}$ defined as in (\ref{EHeqn:dual of Vfs}). Indeed we can consider $\delta$ as a point in $\hat T$ and compute that
\numequation\label{EHeqn compute eigenspaces}
{\rm Hom}_{T}(\delta,J_{B_v}(\Pi_v))=(\mathcal{M}^\ast_{\Pi_v}\otimes_{\cO_{\hat T}} k(\delta))'
\end{equation}
is the dual space to the fiber of $\mathcal{M}^\ast$ at the point $\delta$.
We remark that there is also a derived version of this formula given by
$$R{\rm Hom}_{T}(\delta,J_{B_v}(\Pi_v))=(R\Gamma(\hat T,\mathcal{M}^\ast_{\Pi_v}\otimes^L_{\cO_{\hat T}} k(\delta)))'.$$

\begin{rem}\label{EHrem eigenspacesandcharacters}
\hfill\\(i) At least if $\delta$ is a locally algebraic character the eigenspace ${\rm Hom}_{T}(\delta,J_{B_v}(\Pi_v))$ can be computed using Breuil's adjunction formula \cite[Remarque 4.4]{MR3319547}
 $${\rm Hom}_{T}(\delta,J_{B_v}(\Pi_v))={\rm Hom}_{G_v}(\mathcal{F}_{B_v}^{G_v}(\delta),\Pi_v),$$
 where $\mathcal{F}_{B_v}^{G_v}(-)$ is the functor defined by Orlik and Strauch, see Section~\ref{EH subsection:paraboliccompatibility} above.
Roughly this adjunction says that the finite slope part $$M^\dagger(G,K^p,\kappa^{\rm univ})_{\rm fs}=J_B(\hat S(K^p,L)^{\rm an})$$ consists of those forms $f\in \hat S(K^p,L)^{\rm an}$ that generate a locally analytic subrepresentation all of whose subquotients occur in parabolically induced representations. 
 \\
 (ii) If $\delta_v$ is a locally analytic character of $T_v$ it induces by restriction a character $\kappa_v$ of $T_v^+$ and a character $\chi:\mathcal{A}^+(G_v)=L[\Sigma_v^+]\rightarrow L$. Conversely any such pair $(\kappa_v,\chi)$ canonically defines a character $\delta_v$ of $T_v$. With these notations and using (\ref{EHeqn compute eigenspaces}) we can rewrite 
 $$
 \begin{aligned}
 {\rm Hom}_{T}(\delta,J_{B}(\hat S(K^p,L)^{\rm an}))&={\rm Hom}_{T_v^+}(\delta,(\hat S(K^p,L)^{\rm an})^{N_{0}})\\ &=M^\dagger(G,K^p,(\kappa_v)_v)^{\mathcal{A}^+(G)=\chi}.
 \end{aligned}$$
\end{rem}
\subsubsection{Eigenvarieties}
We indicate (slightly informally) the relation with the construction of eigenvarieties as e.g.~in \cite{MR2207783}.
The Hecke action in Conjecture \ref{EHconj:localglobal1} is in particular the action of a big Hecke algebra $\mathbf{T}$ that is a tensor product of spherical Hecke algebras at good places away from $p$.
Let us write $\hat{\mathbf{T}}(K^p)$ for the inverse limit of the images of $\mathbf{T}$ in the endomorphism rings of 
$$\hat S(K^p,\cO/\varpi^m)=\{f:G(F)\backslash G(\mathbf{A}_F^{\infty})/K^p\rightarrow \cO/\varpi^m\ \text{locally constant}\}.$$
Then $\hat{\mathbf{T}}(K^p)$ is a noetherian semi-local ring and we may form its attached formal scheme $\mathcal{Y}(K^p)={\rm Spf}\,\hat{\mathbf{T}}(K^p)$ as well as its rigid analytic generic fiber $\mathfrak{Y}(K^p)=\mathcal{Y}(K^p)_\eta^{\rm rig}$.
Following \cite{MR2207783} (or our Remark \ref{EHrem
  Jacquetandsheaves} above and the discussion preceding it) the dual
of the Jacquet module $J_B(\hat S(K^p,L)^{\rm an})'$ becomes in a
natural way a coherent sheaf $\mathcal{M}^\ast(G,K^p)$ on
$\mathfrak{Y}(K^p)\times \hat T$. The support of this coherent sheaf
is called the eigenvariety $\mathcal{E}(G,K^p)$ associated with $K^p$
(and the choice of the principal Bernstein blocks that we omit from the notation, compare Remark \ref{EHrem:eigenvarforBernsteinblocks} ).
\begin{rem}
Note that of course one could as well use $M^\dagger(G,K^p,\kappa^{\rm univ})_{\rm fs}$ to construct the eigenvariety. This is the point of view of \cite{MR2075765} or \cite{MR2769113}.
\end{rem}
Conjecturally, as indicated for example in Section~\ref{sec: proofs of reciprocity}, the completed big Hecke algebra $\hat{\mathbf{T}}(K^p)$ should coincide with a product of pseudodeformation rings, and these pseudodeformation rings should either be defined as or proven to be isomorphic to the invariant global sections of the classical stack $\cX_{G,F,S}^{\rm cl}$ underlying $\cX_{G,F,S}$; that is, conjecturally one has an isomorphism
$$\hat{\mathbf{T}}(K^p)\cong \Gamma(\cX_{G,F,S}^{\rm cl},\cO_{\cX_{G,F,S}^{\rm cl}}).$$

Let us write $\cX_{G,F,S}^{\rm ps}$ for the disjoint union of the
formal spectra of pseudodeformation rings corresponding to $(G,F,S)$. Then, using a map $\mathcal{Y}(K^p)\rightarrow \cX_{G,F,S}^{\rm ps}$ (constructed using the construction of Galois representations associated to automorphic forms) we may regard $\mathcal{M}^\ast(G,K^p)$ as a coherent sheaf on the Stein space $\mathfrak{X}_{G,F,S}^{\rm ps}\times \hat T$, with $\mathfrak{X}_{G,F,S}^{\rm ps}=(\cX_{G,F,S}^{\rm ps})_\eta^{\rm rig}$. Conjecture \ref{EHconj: coherentdescriptionoffiniteslop} now roughly says that the coherent sheaf $\mathcal{M}^\ast(G,K^p)$ in some sense localizes over $\mathfrak{X}_{G,F,S}\times\hat T$, namely that there is a coherent sheaf on this stack (the Serre dual of the explicitly described coherent sheaf in the conjecture) whose (derived) global sections are concentrated in degree zero and agree with the global sections of the sheaf $\mathcal{M}^\ast(G,K^p)$ on the Stein space $\mathfrak{X}_{G,F,S}^{\rm ps}\times \hat T$.

\begin{rem}\label{EHrem:fibersandadjunction}
\hfill\\(i) We remark that we can in fact strengthen the computation of eigenspaces (\ref{EHeqn compute eigenspaces}) in order to not just compute the eigenspace for the Hecke action at $p$ (i.e.~the action of the Atkin-Lehner ring $\mathcal{A}^+(G)$) but also for the Hecke action away from $p$, i.e. for the action of the Hecke operators in $\hat{\mathbf{T}}(K^p)$. 
Let $\mathfrak{m}_\rho\subset \mathbf{T}(K^p)$ be the ideal attached to a Hecke eigensystem with associated Galois representation $\rho:{\rm Gal}_{F,S}\rightarrow {\rm GL}_n(L)$. Then $\mathfrak{m}_\rho$ defines a point $x_\rho\in \mathcal{Y}(K^p)_\eta^{\rm rig}$ and we find
\numequation\label{EHeqn eigenspacesoneigenvars}
\begin{aligned}
{\rm Hom}_G(\widehat{\otimes}_v \mathcal{F}_B^G(\delta_v),\hat S(K^p,L)^{\rm an}[\mathfrak{m}_\rho])&={\rm Hom}_T(\prod_v \delta_v,J_B(\hat S(K^p,L)^{\rm an})[\mathfrak{m}_\rho])\\
&=(\mathcal{M}^\ast(G,K^p)\otimes k(x))',
\end{aligned}
\end{equation}
where $x=(x_\rho,(\delta_v)_v)\in \mathcal{Y}(K^p)_\eta^{\rm rig}\times\hat T$.\\
(ii) In good situations the spaces $\mathcal{Y}(K^p)_\eta^{\rm rig}$ and $\mathfrak{X}_{G,F,S}$ agree (up to the trivial action of a copy of $\mathbf{G}_m$) after passing to connected components where for example the residual Galois representation is a fixed irreducible and odd representation. In these cases (up to using Serre duality) Conjecture \ref{EHconj: coherentdescriptionoffiniteslop} gives a very precise description of the coherent sheaf $\mathcal{M}^\ast(G,K^p)$ that only involves stacks of Galois representations (or $(\varphi,\Gamma)$-modules). Hence we also obtain a (conjectural) purely Galois-theoretic description of the eigenspaces $(\ref{EHeqn eigenspacesoneigenvars})$.

\end{rem}

\subsubsection{More general Shimura varieties}
We return to the more general situation, i.e.~we drop the assumption that $G$ is compact at infinity. In particular the locally symmetric spaces (Shimura varieties) $Y(K_f)$ can have cohomology in several degrees, and it is better to work with complexes rather than individual cohomology groups. 

Again we fix split maximal tori $T_v$ and Borel subgroups $B_v$ at the places dividing~ $p$. Let $\kappa=\prod_{v|p}\kappa_v:T\rightarrow L$ denote a locally analytic character. We fix $K^p$ and $I_v$ as above and let $K_f=K^p\times\prod_{v|p}I_v$. Then the representation $$W(\kappa)=\widehat{\otimes}_{v|p}{\rm Ind}_{B_{v,0}}^{I_v}\kappa_v$$ defines a ``locally analytic" coefficient system again denoted by $W(\kappa)$ over $Y(K_f)$, see \cite[2.2, 3]{MR3692014} (the coefficient system $W(\kappa)$ is denoted $\mathscr{A}_\kappa$ in loc.\ cit.).
The cohomology complex $R\Gamma(Y(K_f),W(\kappa))$ is then  equipped with an action of the Atkin--Lehner ring $\mathcal{A}^+(G)$. Of course there is also a variant $W(\kappa^{\rm univ})$ of this construction involving the universal character $\kappa^{\rm univ}$.
\begin{rem}
The complex $R\Gamma(Y(K_f),W(\kappa))$ should be seen as an object in the derived category of solid $L$-vector spaces, in particular it is only well-defined up to quasi-isomorphism. Its cohomology groups are $L$-vector spaces of compact type.  The action of $\mathcal{A}^+(G)$ should also be regarded as an action in the derived category, i.e.~the endomorphism defined by a given element is only well-defined up to homotopy, etc.
\end{rem} 
Similarly we can consider the completed cohomology, compare \cite[5]{Fu2021derivedeigenvar} for the discussion of this in the derived sense. Again we consider complexes instead of individual cohomology groups, that is we consider the complex
$$R\tilde\Gamma(Y(K^p),L)=\big(\lim\limits_{\leftarrow, m}\lim\limits_{\rightarrow, K_p}R\Gamma(Y(K^pK_p),\cO/\varpi^m)\big)\otimes_{\cO}L.$$
Again we regard this as an object in the derived category of solid $L$-vector spaces. In this derived category it is equipped with an action of $G$ such that the cohomology groups are admissible $G$-representations. Even better, the dual of $R\tilde\Gamma(Y(K^p),L)$, i.e.~the complex computing completed homology, is a perfect complex, when we only consider it as a complex of $\cO[[K]]$-modules.
Similar to (\ref{EHeqn:locanfcts}) we may pass to the locally analytic vectors $R\tilde\Gamma(Y(K^p),L)^{\rm an}$. We then have a comparison similar to Proposition \ref{EHprop:compareocandcompletedcohom1}.
\begin{prop}
There is a quasi-isomorphism
$$R\Gamma(Y(K_f),W(\kappa^{\rm univ}))\cong (R\tilde\Gamma(Y(K^p),L)^{\rm an})^{N_0}$$
equivariant for the action of the Atkin--Lehner ring.
\end{prop}
\begin{proof}
This is basically \cite[Theorem 1.2]{Fu2021derivedeigenvar}.
\end{proof}
\begin{rem}
There is of course a comparison of the individual complexes $R\Gamma(Y(K_f),W(\kappa))$ with a fixed $\kappa$ instead of the universal one. For that comparison one has to take a derived eigenspace on the right hand side.
We omit the technical details here.
\end{rem}
Fu \cite[6]{Fu2021derivedeigenvar} continuous the parallels with Section~\ref{subsec: eigenvarcompactatinfty} and gives a derived construction of eigenvarieties building upon the complex $(R\tilde\Gamma(Y(K^p),L)^{\rm an})^{N_0}$, respectively its dual.
We will not discuss this here, but we point out that (at least on a heuristic level) the dual of the finite slope part $(R\tilde\Gamma(Y(K^p),L)^{\rm an})^{N_0}_{\rm fs}$ conjecturally should localize to an  object of the derived category of coherent sheaves on $\mathfrak{X}_{G,F,S}\times\hat T$.
Instead of making the construction of eigenvarieties more precise we state the following analogue of
Conjecture \ref{EHconj:localglobal1} (and hence of Conjecture \ref{conj:cohomology}) in this context. 
Similarly to Section~\ref{subsec: eigenvarcompactatinfty} we have a stack $\mathfrak{X}_{G,F,S}$ and a global-to-local map $f_{\rm rig}$. 
\begin{conj}\label{EHconj:localglobal2}
Let $d$ denote the dimension of $Y(K_f)$. Then there is an isomorphism 
$$R\Gamma(Y(K_f),W(\kappa))[d]\cong R\Gamma_c(\mathfrak{X}^{\rm odd}_{G,F,S},f_{\rm rig}^\ast\mathfrak{A}_{\rm rig}(({\kappa}_v)_{v|p})\otimes V_\mu)$$
compatible with the (derived) Hecke action and the action of the Atkin--Lehner ring.
\end{conj}
\begin{rem}
In fact, as in Section~\ref{subsubsec: cohomology of Shimura}, one should actually replace $G$ by its Weil restriction to $\mathbf{Q}$ and obtain an additional action of ${\rm Gal}_F$ ($F$ now being the reflex field of the Shimura datum) on the cohomology complexes. As in Conjecture \ref{conj:cohomology} the isomorphism in Conjecture \ref{EHconj:localglobal2} should then also be equivariant for the ${\rm Gal}_F$-action. 
\end{rem}

\subsubsection{Companion points and Breuil's socle conjecture}\label{sec:companionpoints}
We return to the setup of Section~\ref{subsec: eigenvarcompactatinfty} and assume in addition that $G$ is a definite unitary group (of course similar results are expected in a more general context).
In \cite[Conj.\ 5.3, Conj.\ 6.1]{MR3319547} Breuil has made a precise
conjecture about the socles (i.e.~the largest semisimple
subrepresentations) of the (admissible) locally analytic
representations $\hat S(G,K^p)^{\rm an}[\mathfrak{m}_\rho]$ for Galois
representations $\rho$ (with coefficients in $L$) that are crystalline
at places dividing $p$ and are such that the eigenvalues $\varphi_{1,v},\dots,\varphi_{n,v}\in L$ of the crystalline Frobenius on the Weil--Deligne representation ${\rm WD}(\rho|_{{\rm Gal}_{F_v}})$ satisfy $\varphi_i/\varphi_j\notin\{1,q_v\}$ for $i\neq j$. 
Breuil's conjecture describes which (tensor products indexed by the places dividing $p$ of) irreducible subquotients of locally analytic principal series representations $({\rm Ind}_{B_v}^{G(F_v)}\delta_v)^{\rm an}$ embed into the space $\hat S(G,K^p)^{\rm an}[\mathfrak{m}_\rho]$.
A weaker version of this conjecture is the precise description of the characters $(\delta_v)_v\in\prod_v \hat T_v$ such that $$(\rho,(\delta_v)_v)\in\mathcal{E}(G,K^p)\subset \mathfrak{X}_{G,F,S}^{\rm ps}\times \prod_v \hat T_v.$$ 
We give a precise version of the weaker conjecture and refer to \cite{MR3319547} for the precise formulation of Breuil's socle conjecture. Using Remark \ref{EHrem:fibersandadjunction} one can easily show that Breuil's conjecture implies the weaker version presented here. 

We continue to assume that $\rho_v={\rho|_{{\rm Gal}_{F_v}}}$ is
crystalline with pairwise distinct Frobenius eigenvalues on ${\rm
  WD}(\rho_v)$ that satisfy the above regularity assumption. Assume
that $\rho_v$ has regular Hodge--Tate weights $\lambdau_v$ and write $\xi_v$ for the corresponding highest weight.
To each Frobenius stable flag $\mathcal{F}_v\subset {\rm WD}(\rho_v)$ (defined over $\Qpbar$) we may associate 
\begin{enumerate}
\item[-] an unramified character 
$$\delta_{\mathcal{F}_v}={\rm unr}_{\varphi_{1,v}}\otimes\dots\otimes{\rm unr}_{\varphi_{n,v}}:T_v\rightarrow L^\times$$
where $\varphi_{1,v},\dots,\varphi_{n,v}$ is the ordering of the Frobenius eigenvalues given by $\mathcal{F}_v$.
\item[-] for each embedding $\tau:F_v\hookrightarrow \bar{\mathbf{Q}}_p$ an element $w_\tau\in \mathcal{S}_n$ giving the relative position of the ($\tau$-part of the) Hodge-Filtration on $D_{\rm dR}(\rho_v)\otimes_{F,\tau}\bar{\mathbf{Q}}_p$ with respect to $\mathcal{F}_v$.  We write $w(\mathcal{F}_v)$ for the tuple $(w_\tau)_{\tau}\in\prod_{\tau:F_v\hookrightarrow \Qpbar}\mathcal{S}_n$.
\end{enumerate}
\begin{conj}\label{conjcompanionpoints}
Let $\rho$ be a representation as above and let $(\delta_v)_{v|p}\in \hat T(L)$. 
Then $(\rho,(\delta_v)_v)\in\mathcal{E}(G,K^p)$ if and only if each $\delta_v$ is of the form 
$$\delta_{v}=z^{w_v w_0\cdot \xi_v}{\delta}_{\rm sm}$$
for some choice of a collection of Frobenius stable flags $(\mathcal{F}_v)_{v|p}$ and Weyl group elements  $w(\mathcal{F}_v)\preceq w_v=(w_\tau)_{\tau:F_v\hookrightarrow L}$. 
Here $w_0$ denotes the longest Weyl group element and $\preceq$ denotes the usual Bruhat order.
\end{conj}
We point out that this conjecture would be a direct consequence of Conjecture \ref{EHconj: coherentdescriptionoffiniteslop}: we only need to compute the support of the coherent sheaf occurring on the right hand side of Conjecture \ref{EHconj: coherentdescriptionoffiniteslop}. However, this support is basically computed using Theorem \ref{thm: local companion points}. 
\begin{rem}
The points $$x_{w_v}=(\rho,(\delta_{\mathcal{F}_v}z^{w_vw_0\cdot\xi_v})_v)$$ for $w_v$ as in the conjecture are often referred to as \emph{companion points} of the point $x_{w_0}=(\rho,(\delta_{\mathcal{F}_v}z^{\cdot\xi_v})_v)$.
Note that the latter point conjecturally corresponds to (the Hecke eigensystem of) a classical automorphic form $f$. The $p$-adic automorphic forms in the Hecke eigensystem defined by $x_{w_v}$ are then often referred to as \emph{companion forms} of $f$.
\end{rem}
Under additional assumptions, so-called Taylor--Wiles assumptions related to the patching construction in \ref{sec: TW patching}, this conjecture is proven in \cite[Theorem 5.3.3, Remark 5.3.2]{MR4028517}.
\begin{theorem}\label{EHTWassumptions}
Assume that 
\begin{enumerate}
\item[-] $p>2$
\item[-] the CM extension $E$ of $F$ is unramified and does not contain a $p$-th root of unity. 
\item[-] the group $G$ is quasi-split at finite places and $K_v$ is hyperspecial if $v$ is a place of $F$ inert in $E$.
\item[-] the residual Galois representation $\bar \rho$ of $\rho$ is adequate.
\end{enumerate}
Then Conjecture \ref{conjcompanionpoints} holds true.
\end{theorem}

\begin{rem}
\hfill\\(i) In fact  Conjecture \ref{EH:conjanalytic} (iii) and Conjecture \ref{EHconj:localglobal1} imply the full conjecture in the locally analytic socle and \cite[Theorem 5.3.3]{MR4028517} proves the full socle conjecture under the assumptions in the above theorem.
The main input into this is the (proof of) the locally analytic Breuil--M\'ezard conjecture, Theorem \ref{BHS:locanBM}.
\\
(ii) There are generalizations of the above Theorem in the case of non-regular Hodge--Tate weights by Zhixiang Wu \cite{ZhixiangWu} and to the case of eigenvarieties associated to more general Bernstein components by Breuil--Ding \cite{BreuilDingeigenvar}.
\end{rem}

\subsubsection{Classical and non-classical $p$-adic automorphic forms}
As in Section~\ref{sec:companionpoints} we continue to work in the fixed
global setup of a definite unitary group $G$.
Let $x=(\rho,(\delta_v)_v)\in \mathfrak{X}_{G,F,S}^{\rm ps}\times\prod_v \hat T_v$ be a point. 
We assume that the characters $\delta_v=z^{\xi_v}\delta_{v, \rm sm}$ are locally algebraic with unramified smooth part $\delta_{v, \rm sm}$.
Recall from Remark \ref{EHrem eigenspacesandcharacters} and Remark \ref{EHrem:fibersandadjunction} that we may view $\delta$ as a tuple $((\kappa_v)_v,\chi)$ consisting of a character $\kappa=(\kappa_v)_v:T^0=\prod_{v|p}T_v^0\rightarrow L^\times$ and a character
$$\chi=\bigotimes_v \chi_{\delta_{v,\rm sm}}:\mathcal{A}^+(G)\rightarrow L$$ of the Atkin--Lehner ring. Moreover, we have the identification of eigenspaces

\numequation\label{EHeqn:comparisonofeigenspaces}
M^\dagger(G,K^p,(\xi_v)_{v|p})[\mathfrak{m}_\rho]^{\mathcal{A}^+(G)=\chi}\cong{\rm Hom}_G(\widehat{\otimes}_v \mathcal{F}_B^G(\delta_v),\hat S(K^p,L)^{\rm an}[\mathfrak{m}_\rho]).
\end{equation}
Assume that $\xi_v$ is the algebraic character associated with a regular Hodge--Tate weight $\lambdau_v$. Then the locally analytic representation ${\rm Ind}_{B_{v,0}}^{I_v}z^{\xi_v}$ contains the algebraic representation $V_{\lambdau_v}$ of highest weight $\xi_v$. In particular we find that the space $M^\dagger(G,K^p,(z^{\xi_v})_{v|p})$ of overconvergent $p$-adic automorphic forms of weight $\xi=(z^{\xi_v})_{v|p}$ contains the space
 $$S(K^pI,(\kappa_v)_v)= \left\{\begin{array}{*{20}c} f:G(F)\backslash G(\mathbf{A}_F^\infty)/K^p\rightarrow \bigotimes_{v|p}V_{\lambdau_v}\\\text{s. th.}\ f(gu)=u^{-1}f(g)\ \text{for all}\ g\in G(\mathbf{A}_F^\infty), u\in \prod_{v|p}I_v\end{array}\right\}$$  
 of classical automorphic forms of weight $\xi$.
 \begin{lem}
 Under the identification (\ref{EHeqn:comparisonofeigenspaces}) an overconvergent $p$-adic automorphic form $$f\in M^\dagger(G,K^p,(z^{\xi_v})_{v|p})[\mathfrak{m}_\rho]^{\mathcal{A}^+(G)=\chi}$$ is classical (i.e.~lies in $S(K^pI,(z^{\xi_v})_v)$) if and only if the corresponding morphism 
 $$\widehat{\otimes}_v \mathcal{F}_B^G(\delta_v)\longrightarrow \hat S(K^p,L)^{\rm an}[\mathfrak{m}_\rho]$$
 factors through the locally algebraic quotient\footnote{See Remark \ref{EHrem localgquotofOS} for the description of that quotient map.} 
 $$\widehat{\otimes}_v \mathcal{F}_B^G(\delta_v)\longrightarrow \bigotimes_{v|p}V_{\lambdau_v}\otimes ({\rm Ind}_{B_v}^{G(F_v)}\delta_{\rm sm})^{\rm sm}.$$
 \end{lem}
 \begin{proof}
 The proof is similar to the proof of \cite[Proposition 3.4]{MR3660309}.
 \end{proof}

The following conjecture is basically implied by standard Fontaine--Mazur type conjectures. 
\begin{conj}\label{classicalityconj}
Let $\rho\in\mathfrak{X}_{G,F,S}$ be a crystalline representation with regular Hodge--Tate weights $\lambdau$ and let $\chi$ be a character of the Atkin--Lehner ring such that the eigenspace $$M^\dagger(G,K^p,(z^{\xi_v})_{v|p})[\mathfrak{m}_\rho]^{\mathcal{A}^+(G)=\chi}$$ is non-zero, i.e.~$\rho$ is associated to a $p$-adic automorphic form $f$.
Then the eigenspace 
$$S(K^pI,(z^{\xi_v})_v)[\mathfrak{m}_\rho]^{\mathcal{A}^+(G)=\chi}$$
is non-zero, i.e.~$\rho$ is associated to a classical automorphic form. 
\end{conj}
In view of our conjectures on functors from locally analytic representations to coherent sheaves on stacks of Galois representations, one can approach this conjecture by studying the morphism of coherent sheaves obtained by applying the conjectural functor $\mathfrak{A}_G^{\rm rig}$ to the map
 $$\widehat{\otimes}_v \mathcal{F}_B^G(\delta_v)\longrightarrow \bigotimes_{v|p}(V_{\lambdau_v}\otimes {\rm Ind}_{B_v}^{G(F_v)}\delta_{\rm sm})^{\rm sm}.$$
of locally analytic representations. In fact this strategy led to a proof of Conjecture \ref{classicalityconj} under the additional Taylor--Wiles assumptions as in Theorem \ref{EHTWassumptions}, \cite[Theorem 5.1.3]{MR4028517}.

It is an interesting question whether the $p$-adic automorphic form $f$ in Conjecture \ref{classicalityconj} is automatically classical. 
Note that this is \emph{not} covered by the phenomenon of companion points or companion forms, as companion forms are $p$-adic forms of non-classical weight (that hence can't be classical!), but here we have fixed the weight to be classical. 
Still our conjectures imply that one cannot expect $f$ to be
classical. Conjecture \ref{EHconj: coherentdescriptionoffiniteslop} and Remark \ref{EHrem:fibersandadjunction} give a precise conjectural description of the eigenspaces
$$S(K^pI,(z^{\xi_v})_v)[\mathfrak{m}_\rho]^{\mathcal{A}^+(G)=\chi}\ \text{and}\ M^\dagger(G,K^p,(z^{\xi_v})_{v|p})[\mathfrak{m}_\rho]^{\mathcal{A}^+(G)=\chi}$$
in terms of a coherent sheaf that can be computed using the Drinfeld compactification $\overline{\mathfrak{X}}_B$. It turns out that this compactification is not Gorenstein in general and these singularities have the consequence that the fiber dimension of the sheaf $\mathcal{M}^\ast(G,K^p)$ can jump. This basically implies that 
  $$M^\dagger(G,K^p,(z^{\xi_v})_{v|p})[\mathfrak{m}_\rho]^{\mathcal{A}^+(G)=\chi}.$$
should sometimes have larger dimension than its classical subspace. 
This expectation can be made into a concrete example in the case of a definite unitary group in three variables, see \cite{hellmann2024patchingmultiplicitiespadiceigenforms}.

\appendix
\addtocontents{toc}{\protect\setcounter{tocdepth}{1}}

\section{Background on category theory}%
\label{sec:general-nonsense}
\renewcommand{\theequation}{\Alph{section}.\arabic{subsection}.\arabic{subsubsection}}
We recall various results in ($\infty$-)category theory that will be applied
in the main body of the notes.

\subsection{Stable \texorpdfstring{$\infty$}{infinity}-categories}
The functors we consider in the categorical local Langlands program will
necessarily be derived, in so far as they won't necessarily map objects from
the abelian category of smooth representations of the $p$-adic group
$G$ under consideration
into an abelian category
of coherent or quasicoherent sheaves on the corresponding moduli stack of Langlands
parameters; rather, 
certain representations 
will be mapped to complexes of sheaves. 
(This is already forced by the desired 
compatibility with Kisin--Taylor--Wiles 
patching, and so is an inevitable part of the story.)
For this reason, it is necessary to 
work with appropriate categories of complexes of sheaves.

If we were writing these notes in the twentieth century, 
we would then proceed to explain that this obliges us to
work with various derived or triangulated categories, and we would go on
to explain which ones.  But we will take advantage of contemporary
advances in homological algebra by working with stable $\infty$-categories,
which avoid the well-known technical drawbacks of triangulated categories,
and are especially useful to work with when one wants to apply
various gluing or limiting processes
(both of which we {\em will} be applying, e.g.\ since we are working with sheaves
on formal algebraic stacks, which have both a topological structure --- 
necessitating gluing arguments --- and an Ind-stack structure, 
necessitating limiting arguments).

We very briefly recall some of these basics, mainly for the purpose of orienting those readers (and those authors!)  who are more familiar with the traditional theory of derived and triangulated categories than with the theory of stable $\infty$-categories.
In order to deal with set-theoretic issues, we follow the approach of~\cite{MR2522659}.
In particular, we fix a Grothendieck universe, %
and sets are called \emph{small} if they belong to this fixed universe.
Furthermore, all limits and colimits are assumed to be small, and we will not usually comment on this.

For the basics of stable $\infty$-categories, and of $t$-structures on them,
we refer to %
~\cite[Ch.~1]{LurieHA}. If $\cC$ is a stable $\infty$-category, then its underlying homotopy
category has a canonical triangulated category structure.
In practice one can often just imagine that one is working 
in this underlying triangulated category, but with the improvement that
cones (and thus homotopy limits and colimits) are canonical.

One important example of a stable $\infty$-category is the derived
$\infty$-category $D(\cA)$ associated to a Grothendieck abelian category$~\cA$,
whose corresponding homotopy category is the usual (unbounded) derived
category of~$\cA$. (Recall that~$\cA$ is a Grothendieck abelian
category if it is cocomplete,
 the formation of filtered colimits in $\cA$ is exact,
and~$\cA$ admits a set of generators.)
Throughout this paper, we use {\em cohomological} indexing, whereas \cite{LurieHA} (along with most of our references) uses homological indexing.
Accordingly, we work with cochain complexes, while the results we cite pertain to chain complexes; we will not usually comment on this point.
Bearing this in mind, the derived $\infty$-category $D(\cA)$  is defined in \cite[Def.\ 1.3.5.8]{LurieHA}
as the {\em differential graded nerve} of the full subcategory of the category
of cochain complexes in~$\cA$ consisting of objects which are {\em fibrant}
for a certain model structure on the category of cochain complexes.
(See~\cite[Prop.~1.3.5.3]{LurieHA} for the construction of this model structure,
and~\cite[Prop.~1.3.5.6]{LurieHA} for a partial description of the fibrant
objects.)  We can also describe $D(\cA)$ 
as the $\infty$-category localization 
of  the category of cochain complexes with respect to
quasi-isomorphisms.  (See~\cite[Prop.~1.3.5.15]{LurieHA} and its
proof.)

If $X$ and $Y$ are objects of an $\infty$-category~$\cC$,
then we have a mapping anima\footnote{In \cite{lurie2009stable} and~\cite{LurieHA},
the terminology {\em spaces} is used for the $\infty$-category arising
from simplicial sets under the usual Quillen model structure,
and for its objects,
but we follow~\cite{cesnavicius2021purity} in
using the terminology {\em anima} instead.}
$\Maps_{\cC}(X,Y)$. 
When $\cC$ is furthermore stable, this mapping anima 
will in fact be a connective spectrum.  %
In our contexts our categories will  be furthermore $\cO$-linear
(for some ring of coefficients~$\cO$),
and so $\Maps_{\cC}(X,Y)$ will be a connective $\cO$-module spectrum,
or equivalently an animated $\cO$-module,
which (by the Dold--Kan correspondence) can be thought of as an object of $D^{\leq 0}(\cO\text{-}\Mod),$
the stable $\infty$-category overlying the derived category of non-positively
graded cochain complexes of $\cO$-modules. 

Again using stability, one can furthermore define a (typically non-connective) spectrum
$\RHom_{\cC}(X,Y)$, %
from which the mapping space can be recovered via
the formula
$$\Maps_{\cC}(X,Y) = \tau^{\leq 0} \RHom_{\cC}(X,Y).$$
Again, in the contexts we consider, our categories will be $\cO$-linear,
and so the various $\RHom_{\cC}(X,Y)$ will 
be $\cO$-module spectra,
which we can regard 
as objects of $D(\cO\text{-}\Mod)$ (unbounded in both directions).

Placing a $t$-structure on $\cC$ amounts (by definition) 
to placing a $t$-structure on its underlying triangulated category.
The full subcategories\footnote{Meaning sub-$\infty$-categories, of course!}
$\cC^{\leq n}$, $\cC^{\geq n}$ (for $n \in \Z$) are then defined
as pullbacks of the corresponding subcategories of the homotopy category.
The key facts are that: %
\begin{enumerate}
\item  $\Maps_{\cC}(X,Y) = 0$ (i.e.\ is contractible) if $X\in \cC^{\leq n}$
and $Y \in \cC^{\geq n+1}$;
\item
The inclusion $\cC^{\leq n} \hookrightarrow \cC$ 
admits a right adjoint $\tau^{\leq n}$.
\item 
The inclusion $\cC^{\geq n} \hookrightarrow \cC$ 
admits a left adjoint $\tau^{\geq n}$.
\end{enumerate}

\begin{df} 
\label{def:heart}
If $\cC$ is a stable $\infty$-category endowed with a $t$-structure,
then we define the {\em heart} of $\cC$  to be the intersection
$\cC^{\heartsuit} := \cC^{\geq 0} \cap \cC^{\leq 0}.$
\end{df}

The heart $\cC^{\heartsuit}$ is an abelian category \cite[Rem.~1.2.1.12]{LurieHA}.
(More precisely, it is equivalent to the nerve of  the  heart  
of the associated $t$-structure on  the homotopy category of~$\cC$,
and this latter heart is literally an abelian category.
Recall that forming  nerves  is the technical device for converting
usual  categories into $\infty$-categories, in the framework of~\cite{MR2522659}.)

If~$\cA$  is a Grothendieck abelian category, 
then the stable $\infty$-category $D(\cA)$ is endowed with a canonical $t$-structure,
arising from the usual $t$-structure on the derived category
of~$\cA$~\cite[Prop.~1.3.5.21]{LurieHA}.  There is a fully faithful functor
$\cA \hookrightarrow D(\cA)$ which induces an equivalence between~$\cA$
and the heart of the $t$-structure on~$D(\cA)$.
(By definition the heart of $D(\cA)$ coincides with the heart of its
full subcategory~$D^+(\cA)$ discussed below, and then this claim
follows from the dual version of~\cite[Prop.~1.3.2.19]{LurieHA}.)
We sometimes refer to $\cA$ as the subcategory of~$D(\cA)$ consisting of  
{\em static} objects  (sometimes known as ``discrete'' or ``classical'' objects). %

We may form the full subcategory $D^+(\cA)  = \bigcup_{n \in \Z} D^{\geq n}(\cA).$
This category may also be constructed as the differential graded nerve
of the full subcategory of bounded below cochain complexes of injective objects of~$\cA$.
(This is because the bounded below fibrant complexes are precisely
the bounded below complexes of injectives; see~\cite[Prop.~1.3.5.6]{LurieHA}.
This construction of~$D^+(\cA)$
is then dual to the construction of $D^-(\cA)$ for abelian categories with
enough projectives given in~\cite[\S 1.3.2]{LurieHA}. Note also that
this construction of~$D^+(\cA)$ only requires that~$\cA$ has enough
injectives, and doesn't require $\cA$ to be Grothendieck.)

\subsection{Associative algebra objects}\label{subsec: E1 objects}The
notion of $E_1$-rings (sometimes known as $E_1$-ring spectra or
$A_\infty$-ring spectra) is defined in~\cite[Defn.\
7.1.0.1]{LurieHA}. By \cite[Rem.\ 7.1.0.3]{LurieHA}, static
$E_1$-rings are associative rings in the usual sense. There is a
theory of left and right modules over an $E_1$-ring~$A$, and we write
$\LMod_A$, $\RMod_A$ respectively for these stable $\infty$-categories. In
the case that $A$ is static, these identify with the derived $\infty$-categories of the usual abelian categories of modules, \cite[Rem.\
7.1.1.16]{LurieHA}.

If~$x$ is an object of a stable $\infty$-category~$\cC$, then by \cite[Rem.\
7.1.2.2]{LurieHA} the spectrum
$\RHom_{\cC}(x,x)$ acquires the structure of an
 $E_1$-ring which we denote $\End_{\cC}(x)$, so that in particular
$\pi_n\End_{\cC}(x)=\Ext^{-n}_{\cC}(x,x)$ for all $n\in\Z$.

\subsection{Pro-categories (and Ind-categories)}\label{subsec: pro
  categories}
We very briefly recall the definitions of $\Pro$- and
$\Ind$-$\infty$-categories; see~\cite[\S 5.3.5]{MR2522659} for more
details. (More precisely, this reference treats the case of
$\Ind$-$\infty$-categories, but for any $\infty$-category~$\cC$, there
is an equivalence $\Ind(\cC)\cong\Pro(\cC^\op)^\op$.)

If $\cC$ is a (small) $\infty$-category,
one can define its associated  pro-category $\Pro(\cC)$
to be the category whose objects are the diagrams
$F: I \to \cC$ indexed by a cofiltered small $\infty$-category~$I$,
with the morphisms between two diagrams $F: I \to \cC, G: J \to \cC$ being defined
by the formula
\[
\Mor_{\Pro(\cC)}(F, G) = \varprojlim_J \varinjlim_I \Mor_\cC(F(i), G(j)).
\] (In fact, as explained in the introduction to~\cite[\S
5.3]{MR2522659}, this is not the most convenient way to set up the
theory of $\Ind$- and $\Pro$- $\infty$-categories, but we ignore this point.)

In the case that $\cC$ is an ordinary category, we will sometimes denote the diagram $F: I \to \cC$ via
$\quoteslim{I}  F(i)$, 
and refer to it as a pro-object of~$\cC$.
The point of this notation is to distinguish the pro-object
from the limit $\varprojlim_I F(i)$ in $\cC$ itself,
if this limit happens to exist.
Often, when employing this notation, we write $X_i$ rather than~$F(i)$,
and so denote the  object  $F$ of $\Pro(\cC)$ as $\quoteslim{I} X_i.$

If $\cC$ is a stable $\infty$-category, so are $\Ind \cC$ and~$\Pro \cC$.
If $\cC$ is equipped with a $t$-structure, there are canonically induced
$t$-structures on $\Ind \cC$ and $\Pro \cC$ (see e.g.~\cite[Prop.~2.13]{MR3935042}).

\subsection{Continuous functors}
\begin{defn}\label{defn: compact object}
Let~$\cC$ be an $\infty$-category admitting filtered colimits. We say that an object $x\in\cC$ is
\emph{compact} if $\Maps_{\cC}(x,\text{--})$ preserves filtered colimits
(see~\cite[Defn.\ 5.3.4.5]{MR2522659}). 
\end{defn}

\begin{rem}%
 \label{rem: compact objects in stable infinity}
  If $\cC$ is a stable $\infty$-category admitting filtered colimits,
  then it is cocomplete. Then $x \in \cC$ is compact
  if and only if $\RHom_{\cC}(x,\text{--})$ preserves filtered
  colimits, %
if and only if
  $\RHom_{\cC}(x,\text{--})$ preserves all (small) colimits. (This is because
  $\RHom_{\cC}(x,\text{--})$ 
  automatically preserves all limits and is therefore exact.)%
\end{rem}

\begin{defn}
  We say that an $\infty$-category~$\cC$ is \emph{compactly generated}
  if it admits filtered colimits, if its subcategory of compact objects $\cC^c$ is small, and if the natural functor
  $\Ind(\cC^c)\to\cC$
  (sending a filtered diagram in~$\cC^c$ to its colimit in~$\cC$) is an
  equivalence.
\end{defn}

\begin{defn}
  We say that a functor between stable $\infty$-categories is \emph{continuous} if it preserves (small)
  colimits. %
\end{defn}
\begin{rem}
  \label{rem:continuous-implies-exact}In many texts, {\em continuous} signifies the preservation of filtered colimits,
and so the functors between stable $\infty$-categories
that we call continuous would more usually be called exact and continuous.
Since we have no occasion to consider non-exact functors, we have opted 
to incorporate the exactness condition into our definition of
continuous functors.
\end{rem}

The notes~\cite{indcoh1} provide a useful overview of categorical notions
related to compact objects.  
The following lemma is partly based on some of the results
discussed in these notes.

\begin{lemma}
\label{lem:compact objects}
Let $F: \cC \to \cD$ be a functor
between %
compactly generated stable $\infty$-categories.

\begin{enumerate}
\item 
$F$ is continuous if and only if it admits a right adjoint~$G$.
\item 
If $F$ is continuous, then $F$ preserves compact objects
if and  only if its right adjoint $G$ is continuous.
\item If $F$ is continuous, fully faithful, and preserves
compact objects, then $F$ also reflects compact objects,
i.e.\ $F(X)$ compact implies $X$ itself is compact. %
\end{enumerate}
\end{lemma}
\begin{proof}
The first part is the  adjoint functor theorem \cite[Cor.\ 5.5.2.9]{MR2522659}.  %
For part (2), we suppose firstly that~$F$ preserves compact objects. In order to
show that~$G$ is continuous, it is enough to show that it preserves
finite colimits and filtered colimits. For the former, recall that a
functor between stable $\infty$-categories preserves finite colimits
if and only if it preserves finite limits; and since~$G$ is a right
adjoint, it in fact preserves all limits.

Thus for any filtered colimit $\colim_j Y_j \iso Y$ in~$\cD$,
we have to show that the induced morphism $\colim_j G(Y_j) \to G(Y)$
is an isomorphism.  By Yoneda's lemma, its enough to do this
after applying $\Mor(X,\text{--})$ for any object $X$ of~$\cC$,
and since  $\cC$ is compactly generated, we can further assume
that $X$ is compact.  We then find that
\begin{multline*}
\Mor\bigl(X,\colim_j G(Y_j)\bigr) 
= \colim_j \Mor\bigl(X, G(Y_j)\bigr) \iso \colim_j  \Mor\bigl( F(X), Y_j\bigr)
\\
\iso  \Mor\bigl( F(X), \colim_j Y_j\bigr) \iso \Mor\bigl(F(X), Y\bigr)
\iso \Mor\bigl(X,G(Y)\bigr),
\end{multline*}
as required.

Conversely, if $G$ is continuous, it in particular preserves
filtered colimits. Accordingly, if~$X$ is compact, and $\colim_j Y_j$
is a filtered colimit in~$\cD$,
then we may write \begin{multline*}
  \Mor\bigl( F(X), \colim_j Y_j\bigr) 
  \iso \Mor\bigl(X,G(\colim_j Y_j)\bigr) \iso 
\Mor\bigl(X,\colim_j G(Y_j)\bigr) \\
= \colim_j \Mor\bigl(X, G(Y_j)\bigr)  \iso  \colim_j \Mor\bigl(F(X), Y_j\bigr),
\end{multline*} so that $F(X)$ is compact, as required.

Finally we turn to (3), which follows by noting that if~$F(X)$ is compact, then we have
\begin{multline*}
  \Mor\bigl( X, \colim_j Y_j\bigr) 
  \iso \Mor\bigl(F(X),F(\colim_j Y_j)\bigr) \iso 
\Mor\bigl(F(X),\colim_j F(Y_j)\bigr) \\
= \colim_j \Mor\bigl(F(X), F(Y_j)\bigr)  \iso  \colim_j \Mor\bigl(X, Y_j\bigr),
\end{multline*}so that~$X$ is compact, as required.
\end{proof}

\subsection{Complete \texorpdfstring{$t$}{t}-structures}\label{sec:Complete t-structures} %
References for this material are~\cite[\S 6]{MR3730514}
and~\cite[\S\S 1.2, 1.3]{LurieHA}.
(Important caution:  these references use homological indexing, 
rather than the cohomological indexing that we use!)

\begin{df}
\label{def:complete categories and completion}
If  $\cC$ 
is a stable  $\infty$-category equipped with
a $t$-structure, then it admits a {\em left completion}
$$\widehat{\cC} := \lim \cC^{\geq n}.$$
The limit is indexed by the totally ordered set~$\Z$, %
and the transition functors 
are given by the truncations $\tau^{\geq n}: \cC^{\geq n-1} \to \cC^{\geq n}.$

There is a canonical functor $\cC\to \widehat{\cC}$,
defined by mapping an object $X$ to  the sequence $(\tau^{\geq n} X)$, 
which has a natural interpretation as an object of~$\widehat{\cC}$.
We say that $\cC$ (with its given $t$-structure) is {\em left complete} if this functor is an
equivalence. 

There are dual notions of right completion and right completeness.  We don't
introduce notation for the right completion, but to be explicit: $\cC$  is
{\em right complete}
if the canonical functor $\cC \to \lim  \cC^{\leq n}$  defined
by $X  \mapsto (\tau^{\leq n} X)$ is an equivalence.
\end{df}

Related to these notions of completeness are corresponding notions of separatedness.

\begin{df}
\label{def:separated categories}
If  $\cC$ 
is a stable $\infty$-category equipped with
a $t$-structure, then
we say it is {\em left separated} if
$$\cC^{-\infty} :=  \bigcap_{n \in  \Z}  \cC^{\leq n}$$
consists only of zero objects.
Dually, 
we say it is {\em right separated} if
$$\cC^{\infty} :=  \bigcap_{n \in  \Z}  \cC^{\geq n}$$
consists only of zero objects.
\end{df}

\begin{remark}
\label{lem:left complete implies left separated}
If $\cC$ is a stable $\infty$-category
equipped with a $t$-structure that is left complete,
then $\cC$ is left separated. (If $X$ is an object
of $\cC^{-\infty}$, then $\tau^{\geq n}X = 0$ for all~$n$,
so that under the equivalence $\cC \iso \widehat{\cC}$,
the object $X$ is identified with a sequence of zero objects, which is a
zero object of~$\widehat{\cC}$; thus $X$ is a zero object of~$\cC$.)
\end{remark}

If $\cA$ is a Grothendieck abelian category, the stable $\infty$-category $D(\cA)$ 
is left and right separated, and also
right complete~\cite[Prop.~1.3.5.21]{LurieHA}.
It need not be  left complete in general~\cite{MR2875857}.

If $\cA$ admits countable
products, and if these are furthermore exact, then $D(\cA)$ admits
countable products (take  products level  wise on complexes),  and 
it follows from~\cite[Prop.~1.2.1.19]{LurieHA} that $D(\cA)$ is
left complete (since Condition~(2) there is satisfied).
For example,  if  $\cA$ is the abelian  category of  modules over
a ring,  then $D(\cA)$ is left complete as well as  right complete.
We will be interested in contexts in which products in $\cA$, although
they  exist, may not be exact, 
and so we note the following mild generalization of~\cite[Prop.~1.2.1.19]{LurieHA}.
(The possibility of such a generalization is noted 
in~\cite[Rmk.~6.1.5]{MR3730514}, where the hypothesis of right $t$-exactness
of products is required to hold only ``up to a finite shift''.
Essentially the same result is also proved
as~\cite[Prop.~8.14]{antieau2021uniquenessinfinitycategoricalenhancementstriangulated}.)

\begin{prop}
\label{prop:left completeness criterion}
Let~$\cC$  be  a  stable $\infty$-category equipped with a $t$-structure.
Suppose that $\cC$ admits countable products, and that the formation of countable products
is of bounded amplitude \emph{(}i.e.\ there is some $a \geq 0$
such that if $\{X_n\}$ is any sequence  of objects of~$\cC^{\leq 0},$
then $\prod_n X_n$ lies in $\cC^{\leq a}$\emph{)}.  Then $\cC$ is  
left complete  if and only if $\cC$ is left separated.
\end{prop}
\begin{proof}
One follows the proof of~\cite[Prop.~1.2.1.19]{LurieHA},
noting now that if $\cdots \to X_n \to X_{n-1} \to \cdots$
is a tower of objects in $\cC^{\leq 0},$
then $\lim X_n$ lies in $\cC^{\leq a+1}.$

Now the proof works just as written, 
except that rather than considering the factorization
$\lim f \to f(n-1)  \to f(n),$
we consider the factorization
$\lim f \to f(n-a-1) \to f(n).$
\end{proof}

\begin{remark}
\label{rem:neeman example}
We elaborate slightly on the construction of \cite{MR2875857} alluded to above:
it gives examples of Grothendieck abelian categories~$\cA$ for which both the hypothesis 
of Proposition~\ref{prop:left completeness criterion}
(that countable products in $D(\cA)$ have bounded cohomological amplitude),
and the conclusion (that $D(\cA)$ is left complete) fail. 
Example~\ref{ex:non-analytic} below gives an explicit example of this construction, 
which arises naturally in the context of smooth representation theory
of (non-analytic) pro-$p$-groups.
\end{remark}

One of the advantages of the framework of stable $\infty$-categories,
compared to that of  triangulated categories, is the following
simple and very general result (which is the dual version to~\cite[Thm.~1.3.3.2]{LurieHA})
on the  existence of right derived functors. Recall that we say that a
functor $F:\cC\to\cD$ between stable $\infty$-categories endowed with
$t$-structures is \emph{left $t$-exact} if $F(\cC^{\ge 0})\subseteq
\cD^{\ge 0}$, that it is \emph{right $t$-exact} if $F(\cC^{\le 0})\subseteq
\cD^{\le 0}$, and \emph{$t$-exact} if it is both left and right $t$-exact.

\begin{thm}%
\label{thm:right derived functors}
If $\cA$ is an abelian category with enough injectives,
and if $\cC$ is a stable  $\infty$-category equipped with 
a right complete $t$-structure,
then $F \mapsto \tau^{\leq 0} F_{|\cA}$ {\em (}the restriction
being taken by identifying $\cA$ with the heart of~$D^+(\cA)${\em )}
is an equivalence between the $\infty$-category of left $t$-exact
functors $D^+(\cA)\to\cC$ which carry injective objects of~$\cA$ into $\cC^{\heartsuit},$
and the ordinary category of left exact functors from $\cA$ to~$\cC^{\heartsuit}$.
\end{thm}

As a particular consequence we have the following result~\cite[Prop.~1.3.3.7]{LurieHA}
(which notoriously does not hold in the same level of generality in
the triangulated category context, because of the non-functoriality
of the formation of cones in that context).

\begin{prop}
\label{prop:D(heart) to C}
Suppose that $\cC$ is a stable $\infty$-category equipped with a right complete $t$-structure,
and that its heart $\cC^{\heartsuit}$ has enough injectives.
Then the inclusion $\cC^{\heartsuit} \hookrightarrow \cC$ extends
essentially uniquely %
\footnote{Meaning that the sub-$\infty$-category
of the $\infty$-category of functors $\Fun\bigl(D^+(\cC^{\heartsuit}),\cC\bigr)$
whose  objects consist of such extensions, whose $1$-morphisms
consist of natural transformations that restrict to the identity on the
restriction of $F$ to $\cC^{\heartsuit}$, and whose higher morphisms 
consist of arbitrary higher morphisms between these $1$-morphisms, is contractible.
This is the $\infty$-categorical analogue of being unique up to a unique natural isomorphism
that restricts to  the identity  on the restriction of $F$ to~$\cC^{\heartsuit}$.}
to a $t$-exact functor
$F: D^+(\cC^{\heartsuit}) \rightarrow \cC$.

Furthermore, the following properties are equivalent:

\begin{enumerate}
\item
 This functor is fully faithful.
\item For any objects $X$ and $Y$  of $\cC^{\heartsuit}$ the  natural map
$\Ext^i_{\cC^{\heartsuit}}(X,Y) \to  \Ext^i_{\cC}(X,Y)$ is an isomorphism
for $i > 0$.
\item For any objects $X$ and $Y$ of $\cC^{\heartsuit}$ with $Y$ injective,
$\Ext^i_{\cC}(X,Y) = 0$ for $i >  0$.
\item For any objects $X$ and $Y$ of $\cC^{\heartsuit}$ with $Y$ injective,
there is an epimorphism $Z \to X$ such that
$\Ext^i_{\cC}(Z,Y) = 0$ for $i > 0$. 
\end{enumerate}
Finally, if these equivalent properties hold, then 
the essential image of~$F$ is equal to~$\cC^+:=\cup_{n}\cC^{\ge n}$.
\end{prop}
\begin{proof}
Except for the insertion of item~(2)  (which we add just for expositional clarity)
this is simply the dual version of~\cite[Prop.~1.3.3.7]{LurieHA}. 
To see the equivalence of~(2) with the other conditions, note that
(1) certainly implies~(2) (even for~$i = 0$; however this case is in any
event automatic,  since $\cC^{\heartsuit}$ is full in~$\cC$ by construction), 
while~(2) clearly implies~(3).  
\end{proof}

Below we will use this proposition both as stated, and in the dual form for projective
objects and categories bounded above  
(which is the case literally treated in~\cite[Prop.~1.3.3.7]{LurieHA}).

\subsection{Compact and coherent objects and regular  \texorpdfstring{$t$}{t}-structures}
\label{subsec:compact and coherent}
A reference for this material is~\cite[\S 6]{MR3730514}.
(Reminder: this reference uses homological indexing,
rather than the cohomological indexing that we use.)

The following is~\cite[Def.~6.2.2]{MR3730514}.

\begin{df}
\label{def:coherent object}
Let $\cC$ be a stable $\infty$-category 
that admits filtered colimits, and is equipped with a $t$-structure
that is  compatible with filtered colimits (which is to say that $\cC^{\geq  0}$,
or equivalently every $\cC^{\geq n}$, is closed under the formation of
filtered colimits  in~$\cC$).  We say that an object $x \in \cC$
is {\em coherent}
if $x$ is bounded below, i.e.\  $x$ lies  in $\cC^{\geq n}$ for some~$n$,
and if $x$ is furthermore compact in~$\cC^{\geq m}$ for all $m  \leq  n$.
\end{df}

\begin{defn}
  \label{defn: Coh of C}In the context of Definition~\ref{def:coherent object}, we write $\Coh(\cC)$ for the full sub-$\infty$-category of
  coherent objects of~$\cC$. This is a stable sub-$\infty$-category by~\cite[Rem.\ 6.2.3.]{MR3730514}.
\end{defn}

\begin{remark}
In the context of Definition~\ref{def:coherent object},
if $m \geq n$, then maps from $x$ to objects of $\cC^{\geq  m}$
factor through $\tau^{\geq m}x,$ and so if $x$ is compact
in $\cC^{\geq  n}$ then we see that $\tau^{\geq m}x$
is compact as an object of $\cC^{\geq m}$ for $m \geq  n$.
\end{remark}

\begin{remark}
If $\cC$ is as in Definition~\ref{def:coherent object}
and is furthermore right complete,
and if $x \in \cC$ is coherent, 
then $x$ is bounded above as well  
as below. (Apply the compactness assumption to the identity map
$x \to x =  \colim  \tau^{\leq n} x.$)
\end{remark}

The following result combines parts of~\cite[Lems.~6.2.4,
6.2.5]{MR3730514}. %
\begin{lemma}
\label{lem:coherent objects}
Let $\cC$ be a stable $\infty$-category endowed with
a right complete $t$-structure that is also compatible with
filtered colimits.  Then the following are equivalent:
\begin{enumerate}
\item The inclusion $\cC^{\geq 0} \hookrightarrow \cC^{\geq -1}$
preserves compact  objects.
\item The full subcategory
of compact objects of $\cC^{\heartsuit}$ {\em (}the heart of the $t$-structure on~$\cC${\em )}
is an abelian subcategory of~$\cC^{\heartsuit}$. 
\end{enumerate}
If these conditions hold,
then the coherent objects in $\cC$ %
are precisely those objects $x$ which are bounded both above and below
and such that each $H^n(x)$ is a compact object of~$\cC^{\heartsuit}$.
\end{lemma}

\begin{defn}
  \label{defn: locally coherent abelian category}A Grothendieck abelian category is  {\em locally coherent} if it is compactly generated,
and if its full subcategory of compact objects
is abelian. (See for example~\cite[Thm.\ 1.6]{MR1434441} for some equivalent
characterisations of this notion.)  %
\end{defn}

The following definition is~\cite[Def.~6.2.7]{MR3730514}.

\begin{df}
\label{def:coherent t-structure}
We say that a $t$-structure on a stable  $\infty$-category admitting
filtered colimits is {\em  coherent}
if it is right complete and is compatible with filtered colimits, and
if furthermore $\cC^{\heartsuit}$ is locally coherent.
\end{df}

Since the compact objects of a locally coherent abelian category
form an abelian subcategory (by definition),
a coherent $t$-structure in particular satisfies the equivalent conditions
of Lemma~\ref{lem:coherent  objects}.

We now define what it means for a coherent $t$-structure to be regular
(following ~\cite[\S 6.3]{MR3730514}).
We first note that if $\cC$ is equipped with a  coherent $t$-structure,
then since by assumption $\cC$  admits filtered colimits, there
is a canonical functor %
\numequation
\label{eqn:functor from IndCoh}
\Ind\Coh(\cC) \to \cC.
\end{equation}
This functor is continuous (essentially by construction).
We then make the following definition. 
\begin{df}
\label{def:regular t-structure}
We say that a coherent $t$-structure on~$\cC$ is {\em regular}
if $\cC$ is compactly generated by its coherent objects %
(in the sense that the coherent objects of~$\cC$ are precisely the compact objects of~$\cC$,
and generate~$\cC$),
or if (equivalently)
the  functor~\eqref{eqn:functor from IndCoh}  is an
equivalence. %
\end{df}
\begin{rem}
  \label{rem:not-usual-defn-regular-t-structure}The definition of a regular $t$-structure in~\cite[Defn.~6.3.1]{MR3730514} does not presume that the $t$-structure is coherent, and does not obviously agree with Definition~\ref{def:regular t-structure}; however in the case that the $t$ structure is coherent, the two definitions agree by~\cite[Prop.~6.3.2]{MR3730514}.
\end{rem}

If $\cC$ is equipped with a coherent $t$-structure,
then the pullback of this  $t$-structure along~\eqref{eqn:functor
from IndCoh} endows $\Ind\Coh(\cC)$ 
with a  $t$-structure, which is again coherent and  is in fact regular~\cite[Prop.~6.3.2]{MR3730514}.
The functor~\eqref{eqn:functor from IndCoh} is $t$-exact by construction.
For any~$n$, it %
restricts to  an {\em equivalence} $$\Ind\Coh(\cC)^{\geq n}  \iso \cC^{\geq n}.$$

If $\cC$ is equipped with a coherent $t$-structure,
then the $t$-structure  on its left completion~$\widehat{\cC}$
is again  coherent~\cite[Prop.~6.3.4]{MR3730514}.   
Again,  the canonical functor 
$$\cC  \to  \widehat{\cC}$$
induces  equivalences
$$\cC^{\geq  n}  \isoto  \widehat{\cC}^{\geq  n}.$$
Thus the formation of $\Ind\Coh(\cC)$  and of $\widehat{\cC}$ are
two different ways of ``filling out'' $\cC^+=\cup_n\cC^{\ge n}$ in the leftward direction.

\begin{remark}
\label{rem:Ind Coh can be non-separated}
If the $t$-structure on $\cC$ is
left separated (e.g.\ because it is left complete;
see Remark~\ref{lem:left complete implies left separated}),
then the kernel of the functor~\eqref{eqn:functor from IndCoh}
consists precisely of $\Ind\Coh(\cC)^{-\infty}.$
In examples when~\eqref{eqn:functor  from IndCoh}
is not an equivalence, i.e.\ when the $t$-structure on $\cC$
is not regular, it is often the case that $\Ind\Coh(\cC)^{-\infty}$
is non-zero, i.e.\ that the $t$-structure on $\Ind\Coh(\cC)$ is not left separated.
\end{remark}

For later reference, we note the following general fact
about $t$-structures and Ind-categories~\cite[Prop.~2.13]{MR3935042}.

\begin{prop}
\label{prop:Ind t-structures}
If $\cC$ is a stable $\infty$-category endowed with a $t$-structure,
then $\Ind \cC$ inherits a $t$-structure, characterized by the requirements
that $\cC  \hookrightarrow  \Ind \cC$ be $t$-exact, 
and that the inclusion $\cC^{\geq 0}  \to \cC$ induce an  equivalence
$\Ind(\cC^{\geq 0})  \iso  \Ind(\cC)^{\geq  0}.$ 
If the $t$-structure on $\cC$ is furthermore  right bounded,
then the  $t$-structure  on $\Ind \cC$ is  right  complete.
\end{prop}

\subsection{Derived functors and adjoints}

If~$F:\cA\to\cB$ is an additive functor between Grothendieck
abelian categories, then~ $F$ induces a functor on the corresponding categories of cochain complexes $\CoCh(\cA)\to\CoCh(\cB)$, and
thus a functor $\overline{F}:\CoCh(\cA)\to D(\cB)$.
Write~$Q:\CoCh(\cA)\to D(\cA)$ for the localization
functor. Following~\cite[\S 7.5.23]{MR3931682}
(which is a ``lifting'' to the $\infty$-categorical setting of a definition
on the level of homotopy categories
that goes back
to Quillen; 
see~\cite[Def.~1.4.1]{MR223432},
or~\cite[Def.~2.3.1]{MR3931682} for a more recent treatment),
a \emph{right derived functor~$RF$}
of~$F$ is a functor $RF:D(\cA)\to D(\cB)$ which is equipped with a natural
transformation $\overline{F}\to RF\circ Q$, which represents the
functor~$\Maps(\overline{F},Q^{*}(-))$.
More precisely, for any other functor
$G:D(\cA)\to D(\cB)$ together with a natural transformation
$\overline{F}\to G\circ Q$, there is a unique (up to a contractible space of
choices) natural transformation $RF\to G$ giving rise to the given
$\overline{F}\to G\circ Q$. 
\begin{prop}%
\label{prop:adjoints}
Suppose that $F:\cA \to \cB$
is an exact functor between Grothendieck abelian categories 
which is compatible with colimits,  
and let $G:\cB \to \cA$ denote the right adjoint to $F$ {\em (}which
exists by the adjoint functor theorem{\em )}.
Then 
\begin{enumerate}
\item
$F$ induces a %
$t$-exact continuous functor
$D(\cA) \to D(\cB)$, which we again denote by~$F$.
\item The right derived functor $RG: D(\cB) \to D(\cA)$ exists,
and is right adjoint to $F$.
\end{enumerate}
\end{prop}
\begin{proof}
That $F$ induces a 
$t$-exact  functor
$D(\cA) \to D(\cB)$  follows immediately from the description of $D(\cA)$ 
as the $\infty$-category localization 
of  the category $\CoCh(\cA)$ of cochain complexes with respect to
quasi-isomorphisms. Once we have shown that~$F$ is a left adjoint, then~$F$ is automatically continuous. %

It thus suffices to prove part (2), which we claim follows from %
  \cite[Thm.\ 7.5.30]{MR3931682}. Indeed, $\CoCh(\cA)$ (and likewise $\CoCh(\cB)$) admits
  the model structure defined in~\cite[Prop.\ 1.3.5.3]{LurieHA}, and ~$D(\cA)$
  is the underlying $\infty$-category of this model category by~\cite[Prop.\
  1.3.5.15]{LurieHA}. Now, since~$F:\cA\to\cB$ is exact it preserves all quasi-isomorphisms, and its
extension $F:D(\cA)\to D(\cB)$ is obviously a left derived functor~$LF$ of ~$F$
(see also~\cite[Lem.\ 7.5.24]{MR3931682}). By ~\cite[Thm.\ 7.5.30]{MR3931682},
we need to show that ~$G:\CoCh(\cB)\to \CoCh(\cA)$ takes quasi-isomorphisms between
fibrant cochain complexes to quasi-isomorphisms. Now, by~\cite[Prop.\
1.3.5.14]{LurieHA}, a quasi-isomorphism between fibrant cochain complexes is
necessarily a homotopy equivalence, and is thus taken by~$G$ to a homotopy
 equivalence, and in particular to a quasi-isomorphism, as required. %
 left adjoint, it is continuous.
\end{proof}

\begin{example}
A typical context in which one might apply Proposition~\ref{prop:adjoints} 
is when $\cB$ is a Grothendieck category, and $\cA$ is a localizing subcategory,
i.e.\ a Serre subcategory which is furthermore closed under the formation of
arbitrary direct sums (and hence arbitrary colimits) in~$\cB$. 
Then $\cA$ is also a Grothendieck abelian category,
and we may take $F:\cA \to \cB$ to be the inclusion.  

Although $F$ is fully faithful by construction, its extension to 
a functor $D(\cA) \to D(\cB)$ need not be fully
faithful in general (as far as we know). %
However, we have the following positive result in this direction.
(The triangulated category analogue of
the dual version of statement~(1) of the proposition --- i.e.\
the version for $D^-$ rather than~$D^+$,
and involving 
projectives rather than injectives --- is a consequence
of~\cite[\href{https://stacks.math.columbia.edu/tag/0FCL}{Tag
   0FCL}]{stacks-project}.)

\begin{prop}
\label{prop:full faithfulness}
As in the preceding discussion, let $F:\cA \hookrightarrow \cB$
be the inclusion of a localizing  subcategory into a Grothendieck abelian category,
and consider its  $t$-exact extensions
$D^+(\cA) \to D^+(\cB)$ and  $D(\cA) \to D(\cB)$.

\begin{enumerate}
\item
If $F$ satisfies the equivalent properties in the statement of
Proposition~{\em \ref{prop:D(heart) to C}},
then the extension $D^+(\cA) \to D^+(\cB)$
is again fully faithful,
and its essential image consists of the full subcategory $D^+_{\cA}(\cB)$
consisting of  objects whose cohomologies lie in~$\cA$.

\item 
Suppose that the extension $D^+(\cA) \to D^+(\cB)$ is fully faithful.
Suppose, furthermore, that the formation of products in~$\cB$ is exact, 
and that the derived right adjoint $RG$ of Proposition~{\em \ref{prop:adjoints}}  
has finite cohomological dimension. Then the extension
$D(\cA) \to D(\cB)$
is again fully faithful,
with essential image equal to $D_{\cA}(\cB)$
{\em (}the subscript~$\cA$ again denoting the full subcategory
consisting of objects whose cohomologies lie in~$\cA${\em )}.
\end{enumerate}
\end{prop}
\begin{proof}
The $t$-structure on $D(\cB)$ induces a 
$t$-structure on $D_{\cA}^+(\cB)$, whose heart is precisely
$\cA$.  
If we assume that $F$  satisfies the equivalent properties in the statement of
Proposition~\ref{prop:D(heart) to C}, then that proposition
gives rise to a $t$-exact equivalence
$$D^+\bigl(\cA)  \iso  D^+_{\cA}(\cB)$$
which extends the identity functor on hearts.  The essential uniqueness
statement of that proposition shows that this equivalence must essentially
coincide with the canonical extension of~$F$.
This proves~(1).

Suppose now that the extension of $F$ to $D^+(\cA)$ is fully faithful;
as above, we continue to denote this extension, as well as the extension to~$D(\cA)$,
 by~$F$.
For any object $X$ of $D(\cA)$, the unit of adjunction gives 
rise to a natural morphism $X \to RG\bigl(F(X)\bigr).$
Thus if $n$ is any integer, applying this to the exact triangle
$\tau^{<n}X \to X \to \tau^{\geq n} X$ (and recalling that~$F$ is
$t$-exact),
we obtain a morphism of exact triangles
$$\xymatrix{
\tau^{<n}  X \ar[d]\ar[r] & X \ar[r]\ar[d] &  \tau^{\geq n} X \ar[d] \\ 
RG\bigl(\tau^{<n} F(X)\bigr) \ar[r] & RG\bigl(F(X)\bigr) \ar[r]
&  RG\bigl(\tau^{\geq n}F(X)\bigr)  \\ 
}
$$
Since the functor~$F$ is assumed to be fully faithful on $D^+(\cA)$,
we see that the right hand vertical arrow is an isomorphism.
Of course $\tau^{< n}X$ has vanishing cohomology in degrees $\geq n,$
and if $a\geq 0$ denote the cohomological amplitude of $RG$,
then $RG\bigl(\tau^{< n}F(X)\bigr)$ has vanishing cohomology
in degrees~$\geq n+a$.
Thus if $i \geq n+a$, we obtain a commutative square
$$\xymatrix{ H^i(X) \ar^-{\sim}[r]\ar[d]  & H^i(\tau^{\geq n}X) \ar^-{\sim}[d] \\
H^i\Bigl( RG\bigl(F(X)\bigr) \Bigr) \ar^-{\sim}[r]
&  H^i\Bigl(RG\bigl(\tau^{\geq n}F(X)\bigr) \Bigr) \\
}
$$
in which (as indicated) all but the left hand vertical arrow are known to
be isomorphisms.  Thus this arrow is an isomorphism as well.  

Since $n$ was arbitrary, we find that the unit of adjunction
$X \to RG\bigl(F(X)\bigr)$ induces an isomorphism
on cohomology in all degrees, and hence is an isomorphism
in~$D(\cA)$. 
Thus $RG$ is a left quasi-inverse to~$F$,
and so the latter is fully faithful, as claimed.

Thus $D(\cA)$ is equivalent to a full subcategory of $D(\cB)$, which
is evidently contained in~$D_{\cA}(\cB)$. At this point
there may well be a direct argument to show that $D_{\cA}(\cB)$ is actually equal
to the  essential image of~$F$, but we finish the argument in a slightly roundabout way.

First, we note that since the formation of products in $\cB$ is exact,
the category $D(\cB)$ admits products, and their formation is $t$-exact.
Proposition~\ref{prop:left completeness criterion} shows that $D(\cB)$  is furthermore
left complete.
One then sees that $D_{\cA}(\cB)$ is also left complete.
Now, if $\{X_n\}_{n=0}^{\infty}$ is a sequence of objects of $D(\cA)^{\leq 0}$, then
$\prod_{n=0}^{\infty} F(X_n)$ is an object of $D(\cB)^{\leq 0}$,
and so 
$RG\bigl(\prod_{n=0}^{\infty} F(X_n)\bigr)$
is an object of $D(\cA)^{\leq a}$ (where, as above, $a$ denotes
the cohomological amplitude of~$RG$).  Furthermore, one checks using the
adjointness of $F$ and $RG$ and the full faithfulness of~$F$
that 
$RG\bigl(\prod_{n=0}^{\infty} F(X_n)\bigr)$
is the product of the $X_n$ in $D(\cA)$.  Thus countable products 
exist in $D(\cA)$, and their formation has bounded cohomological amplitude.
Applying Proposition~\ref{prop:left completeness criterion} again,
we find that that~$D(\cA)$ is left complete as well.
 
Thus we may consider the square
$$\xymatrix{
\widehat{D(\cA)} \ar^-{\widehat{F}}[r]  &   \widehat{D_{\cA}(\cB)}\\
D(\cA) \ar^-\sim[u]\ar^-F[r] &  D_{\cA}(\cB) \ar^-\sim[u] }
$$
in which the hats denote left completions,  the vertical arrows,
which are the canonical functors from the indicated categories
to their left  completions, are (as  we have just noted) equivalences,
the bottom  horizontal arrow is the functor~$F$,
and the top horizontal  arrow  is the  functor~$\widehat{F}$ induced
by~$F$ on left completions.
We claim that this diagram  commutes (up to  natural isomorphism).
Since $\widehat{F}$ is an equivalence (as $F$ restricts to an
equivalence  $D(\cA)^{\geq  n} \to  D_{\cA}(\cB)^{\geq  n}$ for each $n$),
we conclude that the bottom horizontal arrow is an equivalence,
as required.

To see the claimed commutativity, note that if $X$ is an object
of~$D(\cA)$, mapping to the element $\lim \tau^{\geq  n}X$  in~$\widehat{D(\cA)}$,
then we have a natural isomorphism
$$\widehat{F}(\lim \tau^{\geq n}) := \lim F(\tau^{\geq  n}X)  \iso  \lim \tau^{\geq n} F(X),$$
using the $t$-exactness of~$F$,
and $\lim \tau^{\geq n} F(X)$ is precisely the image  of the object $F(X)$ of $D_{\cA}(\cB)$
in~$\widehat{D_{\cA}(\cB)}.$
\end{proof}

\end{example}

\begin{remark}
\label{rem:left completeness}
In the context of 
Proposition~\ref{prop:full faithfulness}~(2),
the proof shows that both $D(\cA)$ and $D(\cB)$ are left complete,
and that one may compute the extension of $F$ from $D^+(\cA)$ to $D(\cA)$
via the formula 
$F(X) := \lim F(\tau_{\geq n} X).$
\end{remark}

\subsection{Semiorthogonal
  decompositions}\label{sub: semiorthogonal generalities}
In this section we outline the basic concepts related to semiorthogonal decompositions
of stable $\infty$-categories.  The stable $\infty$-categories under consideration
will always be cocomplete.  Furthermore, when we speak of a cocomplete subcategory,
we will always mean this in the strongest sense, i.e.\ that the subcategory in question
is closed under the formation of colimits in the ambient category.

We begin by proving some generalities about
generators in cocomplete stable $\infty$-categories.  %

\begin{defn}
  \label{defn: orthogonals}Let $\cC$ be a stable $\infty$-category,
  and $X$ a set of objects of~$\cC$. %
  Then we define the left and right
  orthogonals \[{}^{\perp}X:=\{y\in\cC|\RHom(y,x)=0\textrm{ for all
    }x\in X\},\] \[X^{\perp}:=\{y\in\cC|\RHom(x,y)=0\textrm{ for all
    }x\in X\}.\]
\end{defn}%

\begin{lem}
  \label{lem: left orthogonal is cocomplete}If~$\cC$ is a cocomplete
  stable $\infty$-category, and~$X$ is a set of objects of~$\cC$, then
  ${}^\perp X$ is a cocomplete stable subcategory of~$\cC$.
\end{lem}
\begin{proof}
The cocompleteness of ${}^\perp X$  is immediate
from \[\RHom(\colim_iy_i,x)=\lim_i\RHom(y_i,x),\]so we only need to
check stability. By~\cite[Lem.\ 1.1.3.3]{LurieHA}, it is enough to
check that~$\cC$ is stable under the formation of cofibers and
translations. The claim for cofibers is immediate from the
definitions, and that for translations from the
equality \[\RHom(y[1],x)=\Omega\RHom(y,x)\] which shows that $y\in {}^\perp X$
if and only if $y[1]\in{}^\perp X$.
\end{proof}
\begin{lem}
  \label{lem: right orthogonal of compact is cocomplete}If~$\cC$ is a cocomplete
  stable $\infty$-category, and~$X$ is a set of compact objects of~$\cC$, then
  $X^\perp$ is a cocomplete stable subcategory of~$\cC$.
\end{lem}
\begin{proof}
  Each $x\in X$ is compact, so by Remark~\ref{rem: compact objects in
    stable infinity} we
  have \[\RHom(x,\colim_iy_i)=\colim_i\RHom(x,y_i),\] showing that
  $X^\perp$ is cocomplete. It is immediate from the definitions that
  $X^\perp$ is stable under the formation of fibers, and the result
  follows as in the proof of Lemma~\ref{lem: left orthogonal is
    cocomplete} (using a variant of \cite[Lem.\ 1.1.3.3]{LurieHA} for
  fibers rather than cofibers, which follows from \emph{loc.\ cit.} by
  passage to the opposite category). %
\end{proof}

We will use Lemma~\ref{lem: left orthogonal is cocomplete} to define
the subcategory $\langle X\rangle$ of~$\cC$ generated by~$X$. Before
doing so, we find it useful to consider the following more general situation. Let $F: \cC \to \cD$ be a fully faithful and continuous functor
between cocomplete stable $\infty$-categories.
By Lemma~\ref{lem:compact objects}, %
the functor $F$  admits a right adjoint $G: \cD \to \cC$.
Let $\cA\subseteq \cD$ denote the kernel of~$G$, and
let $\cB \subseteq \cD$ denote the essential image of~$F$. These 
are both stable %
sub-$\infty$-categories of~$\cD$.

\begin{lemma}
\label{lem:orthogonality}In the preceding situation we have the
following statements.
\begin{enumerate}
\item $F$ induces an equivalence $\cC \iso \cB$.
\item $G$  induces an equivalence $\cD/\cA \iso \cC$.
\item The composite $\cB \subseteq \cD \to \cD/\cA$
is an equivalence.
\item %
$\cB$ is closed under the formation of colimits in~$\cD$.
\item $\cA = \cB^{\perp}$.
\item $\cB = {}^{\perp}\cA$.
\end{enumerate}
\end{lemma}
\begin{proof}
(1) follows from the very definition of~$\cB$ as the essential image
of~$F$, together with the assumption that $F$ is fully faithful.
To prove~(2), note first that by the definition of~$\cA$ as the kernel of~$G$,
the functor $G$ induces a functor $\cD/\cA \to \cC.$  
Since the unit of adjunction $\id_{\cC} \iso GF$ is an 
equivalence (by  full faithfulness of~$F$), and the counit is an
equivalence modulo~$\cA$ (i.e.\ gives an equivalence $GFG\iso G$) we deduce that this induced
functor is indeed an equivalence.
(3) is a reformulation, using~(1) and~(2), of the  fact that
$\id_{\cC} \iso GF$.
(4) follows from the fact that $F$ is %
continuous.

To prove~(5), note that $Y \in \cA$ iff $G(Y) = 0$ iff $\RHom\bigl(X,G(Y)\bigr) = 0$
for all $X \in \cC$ iff $\RHom\bigl(F(X),Y\bigr) = 0$ for all $X \in \cC$
(by the adjunction between $F$  and $G$)  iff $Y \in \cB^{\perp}$.
To prove~(6), note first that~(5) gives the inclusion
$\cB \subseteq {}^{\perp}\cA.$
Now let $Y \in \cD$, %
and let $Z$ denote the cofiber of 
the counit of adjunction $FG(Y) \to Y$. %
Since $GF \iso \id_{\cC}$, the morphism $GFG(Y)\to G(Y)$ is an
equivalence, so that $Z \in \cA$. 
In particular, if $Y \in {}^{\perp}\cA,$ then as
$FG(Y)\in\cB\subseteq{}^\perp\cA$, we must also have
~$Z\in{}^\perp\cA$; and since $Z\in\cA$, we see that $\id_Z=0$,
so~$Z=0$, as required. %
\end{proof}
\begin{defn}\label{defn: subcategory generated by some objects}
  If~$\cC$ is a cocomplete
  stable $\infty$-category, and~$X$ is a set of objects of~$\cC$, then
  we define the \emph{subcategory generated by~$X$} to be $\langle
  X\rangle :={}^\perp(X^\perp)$. If $\langle X\rangle =\cC$ then we
  say that $X$ is a \emph{set of generators} of~$\cC$.
\end{defn}
This may look like a rather strange definition, but we note that by
Lemma~\ref{lem: left orthogonal is cocomplete}, $\langle
X\rangle$ is a cocomplete stable subcategory
of~$\cC$. The following lemma shows that it is the smallest cocomplete
stable subcategory containing~$X$, justifying the definition.

  \begin{lem}
    \label{lem: generated subcategory is smallest}Suppose that $\cC$ is a cocomplete
  stable $\infty$-category, that~$X$ is a set of objects of~$\cC$, and
  that~$\cE$ is a cocomplete stable subcategory of~$\cC$ which
  contains every object in~$X$. Then $\cE$ contains~$\langle X\rangle$.
\end{lem}
\begin{proof}
  Applying Lemma~\ref{lem:orthogonality} with~$F$ given by the
  inclusion~$\cE\subseteq\cC$, we see that
  $\cE={}^\perp(\cE^\perp)$. Since $\langle
  X\rangle={}^\perp(X^\perp)$ (by definition), we need to show that
  ${}^\perp(X^\perp)$ is contained in $ {}^\perp(\cE^\perp)$. By definition,
it suffices to prove that  $X^\perp$ contains $\cE^\perp$, which is
immediate from our assumption that~$\cE$ contains~$X$.
\end{proof}
The following corollary is immediate.
\begin{cor}\label{cor: equivalence of definitions of generators}Let
  $\cD$ be a cocomplete stable $\infty$-category, and let~$X$ be a set
  of objects of~$\cD$. Then the following are equivalent.
  \begin{enumerate}
  \item $X$ is a set of generators of~$\cD$.
  \item The only cocomplete stable subcategory of~$\cD$ containing every
    element of~$X$ is~$\cD$ itself.
  \item If $y$ is an object of~$\cD$ satisfying $\RHom(x,y)=0$ for
    every $x\in X$, then $y=0$. 
  \end{enumerate}
\end{cor}

We now recall what we mean by a
semiorthogonal decomposition (which may differ slightly from other
definitions in the literature, in that we allow the decomposition to
be infinite). Let~$\cC$ be a cocomplete stable $\infty$-
category, and let $\cA_1,\cA_2,\dots$ be cocomplete %
stable
subcategories (either finite or infinite in number).
\begin{defn}\label{defn: semiorthogonal decomposition}
  We say that $\cA_1,\cA_2,\dots$ is a \emph{semiorthogonal
    decomposition} of~$\cC$ if
  \begin{enumerate}
  \item The  subcategories $\cA_1,\cA_2,\dots$ generate~$\cC$, and
  \item if $x\in\cA_i$, $y\in\cA_j$ with $i<j$, then $\RHom(x,y)=0$.
  \end{enumerate}
\end{defn}
\begin{rem}
  \label{rem: equivalent conditions for semiorthogonal in terms of
    definition as orthogonals}By definition, condition~(2) in Definition~\ref{defn: semiorthogonal decomposition} is
  equivalent to asking that if $i<j$ then $\cA_i\subseteq
  {}^\perp\cA_j$, and is also equivalent to asking that $\cA_j\subseteq\cA_i^{\perp}$.
\end{rem}

We can define semiorthogonal decompositions in terms of generators, as
in the following definition and lemma.
\begin{defn}
  \label{defn: weakly exceptional collection}Let~$\cC$ be a cocomplete
  stable $\infty$-category. We say that a (possibly
  infinite) sequence $x_1,x_2,\dots$ of objects of~$\cC$ is
  \emph{weakly exceptional} if:
  \begin{enumerate}
  \item The~$x_i$ are all compact objects of~$\cC$.
  \item The~$x_i$ generate~$\cC$.
  \item $\RHom(x_i,x_j)=0$ for $i<j$.
  \end{enumerate}
\end{defn}%

\begin{lem}\label{lem: weakly exceptional gives
    semiorthogonal}Let~$\cC$ be a cocomplete
  stable $\infty$-category, let $x_1,x_2,\dots$ be a weakly
  exceptional sequence of objects of~$\cC$, and let~$\cA_i:=\langle x_i\rangle$ be the
  subcategory generated by~$x_i$. Then $\cA_1,\cA_2,\dots$ is a
  semiorthogonal decomposition of~$\cC$.
  \end{lem}
\begin{proof}By assumption the sequence $x_1,x_2,\dots$
  generates~$\cC$. %
  Suppose that $i<j$. Since~$x_i$ is compact, we see from
  Lemma~\ref{lem: right orthogonal of compact is cocomplete} that
  $x_i^{\perp}$ is  cocomplete and stable, and since it
  contains $x_j$ by assumption, we see that it must contain $\langle
  x_j\rangle =\cA_j$.

  Now let $y\in\cA_j$ be arbitrary. We have just shown that
  $\RHom(x_i,y)=0$, so that $x_i\in{}^\perp y$, which is  cocomplete
  and stable by Lemma~\ref{lem: left orthogonal is cocomplete}. Thus
  ${}^\perp y$ contains $\cA_i=\langle x_i\rangle$, and since~$y$ is
  arbitrary, we have shown that $\RHom(x,y)=0$ for any $x\in
  \cA_i,y\in\cA_j$, as required.\end{proof}

\subsubsection{Generators and adjunctions}
\label{sec:gener-adjunct}
We now return to the setting of Lemma~\ref{lem:orthogonality}. Let $F': \cC \iso  \cB$ denote the equivalence of part~(1),
and let $G': \cB \iso \cC$ denote its quasi-inverse.
Let $\iota: \cB \subseteq \cD$ denote the inclusion,
and $\pi: \cD \to \cD/\cA$ denote the projection.
Then $F \iso \iota  F'$, while $G \iso G'  (\pi\iota)^{-1}\pi$
(noting that $\pi\iota$ {\em is} an equivalence, by~(3)).
Now, by~(5) and~(6), either of the subcategories~$\cA$ and $\cB$
determines the other, and the functors $\iota$ and $\pi$ are
defined canonically in terms of~$\cA$ and~$\cB$.
Thus the original data of the functor~$F$ is determined
by giving either of the subcategories $\cA$ or~$\cB$,
together with the equivalence~$F'$.

Suppose that $\{X_i\}_{i \in I}$ is a collection of generators
of~$\cC$. %
Then we see that  $$\cB^{\perp} =  \cap_{i \in I} F(X_i)^{\perp}
= \{ Y \in \cD \, | \,  \RHom\bigl(F(X_i),Y\bigr) = 0 \text{ for all  } i \in I\}.$$
By Lemma~\ref{lem:orthogonality}~(5), we then see that
$\cA =  \cap_{i \in I} F(X_i)^{\perp}.$

Thus, suppose that we know the values of the functor~$F$ on the
generating collection~$\{X_i\}_{i\in I}$, say $Y_i := F(X_i)$.
Then we can characterize $\cB$ as the stable full  sub-$\infty$-category 
generated by the~ $Y_i$. %
Given~$\cB$, 
we can then determine~$\cA$,
and although from this data we can't determine~$F$ itself,
we can determine the composite~$FG$.  Namely, we have~$FG \iso
\iota (\pi \iota)^{-1} \pi.$

Another way to express this is as follows: any object $Y$ of $\cD$ {\em semiorthogonal decomposition}, which is a cofiber sequence (equivalently, a fiber sequence, or exact triangle)
of the form
$$X \to Y \to Z$$
with $X \in \cB$ and $Z \in \cA$.
This cofiber sequence  is unique up to equivalence, and is given by taking $X = FG(Y)$, and taking $Z$ to be the cofiber of the counit map in the adjunction $FG(Y) \to Y$. Accordingly, we recover $FG(Y)$ as the term $X$ in the distinguished triangle.

We can apply this formalism in applications where the functor~$F$ is
conjectural. More precisely, we can suppose given a collection of objects
$\{Y_i\}_{i \in I}$ of~$\cD$ (which we imagine arise as the
image of a generating set of objects $X_i \in \cC$ under a functor~$F$;
but we won't actually be given~$F$).
Then we can define $\cB$ to be the  subcategory  generated by the~$Y_i,$
just as in  the preceding paragraph.  
By construction the inclusion $\cB \subset \cD$ is fully faithful and
preserves colimits,
and so by Lemma~\ref{lem:compact objects}~(1) admits a right adjoint.
Lemma~\ref{lem:orthogonality} then applies: we may form $\cA := \cB^{\perp}$,
every object $Y$ of $\cD$ admits a functorial semiorthogonal
decomposition given by a cofiber sequence
$$  X \to Y \to Z $$ with $X  \in \cB$  and $Z  \in \cA$,
and the functor $Y \mapsto X$ is the right  adjoint to the inclusion of~$\cB$.

\subsubsection{Restriction to open substacks}\label{subsec: restricting
  semiorthogonal to open substacks}Suppose further that we have a cocomplete
stable $\infty$-category $\cD'$ with a functor $j^*:\cD\to\cD'$,
together with a fully faithful right adjoint
$j_*:\cD'\to\cD$. %
(In applications,
$\cD$ will be a category of sheaves on a stack $\cX$, $\cD'$ will
be the category of sheaves on an open substack $j:\cU\to\cX$, and
$j^*,j_*$ will have their usual meanings.)  Let $\cB'$, $\cA'$ be
sub-stable $\infty$-categories of~$\cD'$ which give a semiorthogonal
decomposition as above: so $\cA',\cB'$ together generate $\cD'$, and
we have $\cB' = {}^{\perp}\cA' $. 

Set $\cA:=j_*\cA'$. Then (since the full faithfulness of~$j_*$ implies
that $j^*j_*$ is naturally equivalent to the identity) we have
$\cA'=j^*\cA$, and by adjointness we have \numequation\label{numeqn:
  orthogonal of pushforward}{}^{\perp}\cA=(j^*)^{-1}({}^{\perp}\cA')=(j^*)^{-1}(\cB').\end{equation}

\section{Coherent, Ind-Coherent, and Pro-Coherent sheaves on
  Ind-algebraic stacks}\label{app: Ind coherent sheaves on stacks}

\subsection{Modules over \texorpdfstring{$I$}{I}-adically complete rings}
We will establish some results relating various stable $\infty$-categories of
modules.  These provide the commutative algebra background  for the theory 
of (Ind-)coherent sheaves on formal algebraic  stacks.  
In fact our results apply in certain non-commutative contexts as well,
and in this case we will use them below to deduce certain facts
about stable $\infty$-categories of smooth representations of $p$-adic analytic groups.

The setup is as follows:  we let $A$ be a (not necessarily commutative) Noetherian
ring,  and $I$ a two-sided ideal in $A$ with respect to which $A$ is $I$-adically complete. 
We  furthermore  assume that $I$ satisfies the {\em Artin--Rees property}: namely,
that any $A$-submodule $M$ of $A^{\oplus n}$ (for any $n\ge 0$) is closed with respect to the $I$-adic
topology on~$A^{\oplus n}$, and the topology on $M$ induced
by the $I$-adic topology  on $A^{\oplus n}$ coincides with the $I$-adic topology on~$M$.
It  follows from this hypothesis that  any finitely generated $A$-module
is complete and separated with respect  to its $I$-adic topology,
and that morphisms between finitely generated $A$-modules have $I$-adically closed closed image
and are strict with respect to the $I$-adic topologies on their source and target. 

If $A$ is commutative, then the Artin--Rees property holds automatically.  
If $A$ is non-commutative, then it need not hold (as far as we know), %
but it does hold in one 
particular case that we are interested in. 

\begin{example}
\label{ex:compact analytic group ring}
If $H$ is a compact $p$-adic analytic  group, if $P$ is a normal  open pro-$p$  subgroup
of $H$ (the identity has a neighbourhood basis consisting of such~$P$),  
and if $I$ is the kernel of the projection $\cO[[H]]  \to k[H/P],$
then the $I$-adic  topology on $\cO[[H]]$, and more generally on any
finitely generated $\cO[[H]]$-module, coincides with the canonical topology
(as discussed e.g.\ in Section~\ref{subsec:cg top mod} below).  Thus $I$ satisfies
the Artin--Rees property.
\end{example}

\begin{df}
We write $D(A):= D(A\text{-}\Mod)$,
and similarly for any $n$ we write $D(A/I^n) := D(A/I^n\text{-}\Mod)$.
\end{df}

Since $A$ is Noetherian, the $t$-structure on $D(A)$ is coherent.
Indeed the compact objects in its heart, which is the abelian category $A\text{-}\Mod$,
are precisely the finitely presented $A$-modules, and they form
a generating  abelian subcategory of $A\text{-}\Mod$.  The full subcategory
of $D(A)$ consisting of coherent objects is thus precisely 
the category  $\Dfp^b(A)$ consisting of (cohomologically) bounded complexes
whose cohomologies are finitely presented (or, equivalently,
finitely generated, since $A$ is Noetherian).
Similar remarks apply with $A/I^n$ in place of~$A$.

There are obvious exact inclusions
\numequation
\label{eqn:torsion inclusions}
A/I^m\text{-}\Mod \hookrightarrow A/I^n\text{-}\Mod  \hookrightarrow
A\text{-}\Mod,
\end{equation}
for  $m \leq  n$,
which induce $t$-exact functors
$$\Dfp^b(A/I^m)  \to \Dfp^b(A/I^n) \to \Dfp^b(A).$$
Since these functors  are compatible in  an evident way
as $m$ and $n$  vary,  we obtain a functor
\numequation
\label{eqn:Coh functor}
\colim_n  \Dfp^b(A/I^n) \to  \Dfp^b(A)
\end{equation}
(the colimit being taken in the  $\infty$-category of stable $\infty$-categories).
We will show 
below that~\eqref{eqn:Coh functor} is fully faithful,  and describe its essential image.
To begin with, we consider an analogous construction on the level of abelian categories.

\begin{df}
We denote by $(A\text{-}\Mod)_I$
the full abelian subcategory (indeed Serre subcategory) of $A\text{-}\Mod$
consisting of
modules each element of which is annihilated by some power of~$I$.
Equivalently, $(A\text{-}\Mod)_I$ is the closure under filtered colimits
of the  union $\bigcup_n A/I^n\text{-}\Mod,$ the union taking place in
$A\text{-}\Mod$. %
\end{df}

\begin{remark}
\label{rem:I-torsion subcat properties}
Since $(A\text{-}\Mod)_I$ 
is a Serre subcategory of the Grothendieck abelian category ~$A\text{-}\Mod$
which is furthermore closed under the formation
of arbitrary direct sums, it is itself a Grothendieck abelian category. 
In particular $D\bigl((A\text{-}\Mod)_I\bigr)$ is defined, and is equipped with
its canonical $t$-structure, having $(A\text{-}\Mod)_I$ as its  heart.

The compact objects in the category $(A\text{-}\Mod)_I$
are precisely the finitely presented modules (i.e.\ those that are compact
simply as $A$-modules), and these form a generating
 abelian subcategory of $(A\text{-}\Mod)_I$.  Thus the $t$-structure
on~$D\bigl((A\text{-}\Mod)_I\bigr)$ is coherent.
Its full subcategory of coherent objects is equal to $\Dfp^b\bigl((A\text{-}\Mod)_I),$
the subcategory consisting of complexes whose cohomologies are finitely presented.
\end{remark}

There is one more collection of categories we need to introduce.

\begin{df}
We let let
$D_I(A)$ denote the full  subcategory of $D(A)$ 
consisting of objects whose cohomologies lie in $(A\text{-}\Mod)_I$, 
and set $\DfpI^b(A) = \Dfp^b(A) \cap D_I(A).$
\end{df}

The heart of the $t$-structure on $D_I(A)$ is precisely
$(A\text{-}\Mod)_I$,  which is locally coherent,
and hence $\DfpI^b(A)$ is precisely the full  subcategory
of coherent objects of $D_I(A)$.

The inclusions $A/I^n\text{-}\Mod \hookrightarrow A\text{-}\Mod$ of~\eqref{eqn:torsion inclusions}
evidently factor through~$(A\text{-}\Mod)_I$, 
and so~\eqref{eqn:Coh functor} factors as
\numequation
\label{eqn:Coh functors}  
\colim_n  \Dfp^b(A/I^n) \to  \Dfp^b\bigl((A\text{-}\Mod)_I) \to  \DfpI^b(A) \hookrightarrow
\Dfp^b(A).
\end{equation}
We will show that the first two functors in this factorization are
equivalences.

\begin{prop}
\label{prop:first Coh equivalence}
The functor 
\numequation
\label{eqn:one direction}
\colim_n  \Dfp^b(A/I^n) \to  \Dfp^b\bigl((A\text{-}\Mod)_I)
\end{equation}
{\em (}the first arrow of~{\em \eqref{eqn:Coh functors})}
is an equivalence. 
\end{prop}
\begin{proof}
We will prove this by constructing a quasi-inverse to~\eqref{eqn:one direction}. 
But we begin by noting that~\eqref{eqn:one direction} is essentially surjective,
since we will need to use
this in our construction of the quasi-inverse.
This essential surjectivity  amounts to the fact that if $X^{\bullet}$ is
a bounded complex of objects in~$(A\text{-}\Mod)_I$, each of whose cohomologies is finitely
generated, then we can find a quasi-isomorphism $Y^{\bullet} \to X^{\bullet}$
of bounded complexes with each $Y^{\bullet}$ being a finitely presented object
of~$(A\text{-}\Mod)_I$ (so that $Y^{\bullet}$ is in fact a complex of
$A/I^n$-modules for some $n$). (This is standard, using that since $A$
is Noetherian, all submodules of finitely generated $A$-modules are finitely
presented, but for the convenience of the reader we sketch a
proof. Assume that~$X^m$ is finitely presented for~$m>n$, and
let~$(X^n)'$ be a finitely generated $A$-submodule of $X^n$ which
surjects onto both  $H^n(X^\bullet)$ and $\im(d^n:X^n\to
X^{n+1})$. Let~$Y^\bullet$ be the subcomplex of~$X^\bullet$ with
$Y^r=X^r$ for $r\ne n-1,n$, $Y^n=(X^n)'$,
$Y^{n-1}=(d^{n-1})^{-1}(Y^n)$. Then the inclusion $Y^\bullet\to
X^\bullet$ is a quasi-isomorphism, and $Y^m$ is finitely presented for
all $m\ge n$, and we may proceed by downward induction on~$n$.)

We now turn to constructing the required quasi-inverse.
Proposition~\ref{prop:Ind t-structures}
endows
$\Ind\colim_n\Dfp^b(A/I^n)$ with a right complete $t$-structure,
whose heart is
precisely~$(A\text{-}\Mod)_I$, 
so that Proposition~\ref{prop:D(heart) to C} gives rise to a $t$-exact functor
\numequation
\label{eqn:preliminary functor}
D^+\bigl((A\text{-}\Mod)_I\bigr) \to \Ind\colim_n\Dfp^b(A/I^n).      
\end{equation}
We claim that~\eqref{eqn:preliminary functor} 
restricts to a functor
\numequation
\label{eqn:other direction}
\Dfp^b\bigl((A\text{-}\Mod)_I\bigr)\to \colim_n \Dfp^b(A/I^n),
\end{equation}
which is furthermore quasi-inverse to~\eqref{eqn:one direction}.

The composite of~\eqref{eqn:preliminary functor}
with the functor $D^+(A/I^m) \to D^+\bigl((A\text{-}\Mod)_I)$
(for some fixed~$m$)
is a $t$-exact functor %
\numequation
\label{eqn:first composite}
D^+(A/I^m) \to \Ind\colim_n\Dfp^b(A/I^n)
\end{equation}
whose restriction to $A/I^m\text{-}\Mod$
coincides with the canonical functor
$$A/I^m\text{-}\Mod = \Ind\Dfp^b(A/I^m)^{\heartsuit} \hookrightarrow \Ind\Dfp^b(A/I^m) 
\to \Ind\colim_n \Dfp^b(A/I^n).$$
The same is true of the canonical functor
\numequation
\label{eqn:second composite}
D^+(A/I^m) \iso \bigl(\Ind\Dfp^b(A/I^m)\bigr)^+ \to \Ind\colim_n\Dfp^b(A/I^n).  
\end{equation}
Thus Proposition~\ref{prop:D(heart) to C}
shows that the two functors~\eqref{eqn:first composite} and~\eqref{eqn:second composite}
 essentially coincide.

Consequently, the composite of~\eqref{eqn:preliminary functor} with the canonical functor
$\Dfp^b(A/I^m)\to D^+\bigl((A\text{-}\Mod)_I\bigr),$
which we can factor as
$$\Dfp^b(A/I^m) \to D^+(A/I^m) \buildrel \text{\eqref{eqn:first  composite}} \over
\longrightarrow \Ind\colim_n \Dfp^b(A/I^n),$$
essentially coincides with the composite
$$\Dfp^b(A/I^m) \to D^+(A/I^m) \buildrel \text{\eqref{eqn:second composite}} \over
\longrightarrow \Ind\colim_n \Dfp^b(A/I^n),$$
which in turn evidently essentially coincides with the composite
$$\Dfp^b(A/I^m) \to \colim_n \Dfp^b(A/I^m) \hookrightarrow \Ind \colim_n\Dfp^b(A/I^n).$$
Varying~$m$, and relabelling it as (the dummy variable)~$n$,
we find that the composite of~\eqref{eqn:preliminary functor}
with the functor 
$$\colim_n \Dfp^b(A/I^n) \buildrel \text{\eqref{eqn:one direction}}\over\longrightarrow
 \Dfp^b\bigl((A\text{-}\Mod)_I) \hookrightarrow D^+\bigl((A\text{-}\Mod)_I\bigr)$$
coincides (up to natural equivalence) with the inclusion
$$\colim_n \Dfp^b(A/I^n) \hookrightarrow \Ind \colim_n \Dfp^b(A/I^n).$$

Since~\eqref{eqn:one direction} is essentially surjective,
we also find that the restriction of~\eqref{eqn:preliminary functor}
to $\Dfp^b\bigl((A\text{-}\Mod)_I\bigr)$ does have image lying in~$\colim_n \Dfp^b(A/I^n)$,
and  thus  does induce a functor~\eqref{eqn:other direction}.
Furthermore, what we have shown so far is that the composite of \eqref{eqn:other direction}
with \eqref{eqn:one direction} essentially coincides with the identity 
on~$\colim_n \Dfp^b(A/I^n)$.

It remains to consider the composite of \eqref{eqn:one  direction}
with \eqref{eqn:other direction}.
This composite is a $t$-exact functor 
from  $\Dfp^b\bigl((A\text{-}\Mod)_I)$
to itself, which restricts to the identity on~$\Dfp^b\bigl((A\text{-}\Mod)_I)^{\heartsuit}.$
Theorem~\ref{thm:right derived functors} shows that it essentially coincides
with the identity functor.
(Technically, we should be applying Theorem~\ref{thm:right derived functors}
to $D^+\bigl((A\text{-}\Mod)_I)$ rather than its full subcategory~$\Dfp^b\bigl((A\text{-}\Mod)_I\bigr)$.
But we can do this completely formally, by forming the canonical extension
of our functor to a $t$-exact endofunctor of $\Ind\bigl(\Dfp^b\bigl((A\text{-}\Mod)_I\bigr)$,
restricting this to
$\Bigl(\Ind\bigl(\Dfp^b\bigl((A\text{-}\Mod)_I\bigl)\Bigr)^+$,
and then using the canonical equivalence
$D^+\bigl((A\text{-}\Mod)_I\bigr) \iso
\Ind\bigl(\Dfp^b\bigl((A\text{-}\Mod)_I)^+$.)
\end{proof}

Before showing that the second functor in~\eqref{eqn:Coh functors}
is an equivalence, 
we recall what is presumably a standard result comparing injective
objects in $A\text{-}\Mod$ and $(A\text{-}\Mod)_I$.

The inclusion $(A\text{-}\Mod)_I \hookrightarrow A$ has a right adjoint,
namely $M \mapsto M[I^{\infty}].$  Since the right adjoint to  an
exact functor preserves injectives, we see that if $M$ is an
injective $A$-module, then $M[I^{\infty}]$ is  an injective object
of~$(A\text{-}\Mod)_I$. But in fact more is true, as the following lemma shows.

\begin{lemma}
\label{lem:torsion in injectives}
 If $M$ is an injective $A$-module, then $M[I^{\infty}]$ 
is again an injective $A$-module.
\end{lemma}
\begin{proof}
By Baer's criterion, we must show that if $J$ 
is an ideal of $A$, then any morphism $f:J\to M[I^{\infty}]$ extends to a
morphism $A \to M[I^{\infty}].$  
Since $J$ is finitely generated ($A$ being Noetherian),  the morphism $f$
has  image lying in $M[I^m]$  for some $m$, and  so factors through
the quotient $J/I^m  J$  of~$J$.   Now our Artin--Rees hypothesis implies that   
$I^n \cap J  \subseteq I^m J$ for  
some  $n\geq m.$  Thus $f$  factors through a  morphism $g: J/I^n \cap  J \to M[I^{\infty}].$
If we rewrite the domain  of $g$  as $I^n+J/I^n$, and regard this as a submodule
of $R/I^n$, then the hypothesis that $M$ is injective ensures
that $g$ extends to a morphism $R/I^n \to M,$ which obviously has image
lying in $M[I^n] \subseteq M[I^{\infty}].$  Composing
this extension of $g$ with the quotient map $R\to R/I^n$ 
gives the required extension of~$f$. 
\end{proof}

\begin{cor}
\label{cor:preservation of injectives}
The inclusion $(A\text{-}\Mod)_I\hookrightarrow  A\text{-}\Mod$ preserves injectives.
\end{cor}
\begin{proof}
Let $M$ be an injective object of $(A\text{-}\Mod)_I$,  and
let  $M \hookrightarrow N$ be an inclusion of $M$
into an injective $A$-module.
This inclusion evidently factors through an inclusion $M \hookrightarrow N[I^{\infty}].$
This is a morphism of injective objects of  $(A\text{-}\Mod)_I$, and hence is a
split inclusion.  On the other hand, Lemma~\ref{lem:torsion in injectives}
shows that $N[I^{\infty}]$ is again an injective $A$-module, and hence  so 
is its direct summand~$M$.
\end{proof}

\begin{prop}
\label{prop:second Coh equivalence}
The canonical functor $D^+\bigl((A\text{-}\Mod)_I) \to \D_I^+(A)$
is an equivalence, which restricts to an equivalence
$\Dfp^b\bigl((A\text{-}\Mod)_I\bigr) \iso \DfpI^b(A)$
{\em (}and this is the second arrow of~{\em \eqref{eqn:Coh functors})}.
\end{prop}
\begin{proof}%
The first claimed equivalence follows by combining Corollary~\ref{cor:preservation
of injectives} with Proposition~\ref{prop:full faithfulness}~(1).
Restricting this equivalence to the compact objects in the
truncations~$\tau^{\ge n}$ for each~$n$, we see that the equivalence then restricts to %
an equivalence
$\Dfp^b\bigl((A\text{-}\Mod)_I\bigl) \iso  \DfpI^b(A).$
\end{proof}

\subsection{Ind-coherent  sheaves  on Ind-algebraic stacks}\label{app:
IndCoh on Ind algebraic stacks}
In this section we introduce (or, really, recall from the literature)
some basic notions related to derived categories of sheaves of $\cO$-modules
on algebraic and formal algebraic stacks.
Our discussion is slightly unsatisfactory, for the reason that we wish
to work with stable $\infty$-categories, rather than merely the underlying
triangulated categories, but our basic reference for the theory of stacks
is~\cite{stacks-project}, which does not use $\infty$-categorical
notions.

For example, in the $\infty$-categorical world, equivariance is a
perfectly good notion for objects in the  (stable $\infty$-categorical version of the)
derived category of quasicoherent sheaves  on a scheme, and so by writing
an algebraic stack $\cX$ as a quotient $[U/R]$ for some scheme $U$ and groupoid~$R$,
one can define the stable $\infty$-category of quasicoherent complexes on~$\cX$
as the $R$-equivariant objects of the corresponding category for~$U$. 
(Cf.\ Remark~\ref{rem:qc via equivariance} below.) 
In this way one can even  define  the {\em derived category} of quasicoherent 
complexes on $\cX$ {\em prior} to having defined the abelian category
of quasicoherent  sheaves on~$\cX$; the latter  can then be defined as
the heart of a $t$-structure  on the former.
In this approach, the more general concept of $\cO_{\cX}$-module, or the derived
category of  complexes of such, does not even arise.   
(And the question of, e.g., the relationship between $D\bigl(\QCoh(\cX)\bigr)$
and $D_{\qc}(\cO_{\cX})$ doesn't come up.)
However, in  the framework 
of~\cite{stacks-project}, such a definition is not feasible.

We don't really try to reconcile these conflicting approaches here.
We simply state
various possible definitions, with pointers to the  literature, and do our best to indicate the
relationships between them.
Ultimately, the  $\infty$-categorical perspective will be the one best suited to our
goals.

\subsubsection{The case of algebraic stacks}
If $\cX$ is
an algebraic stack,
then we will let $D(\cX)$ denote the stable $\infty$-category of quasicoherent
complexes on $\cX$.  
Let us briefly indicate what this means.
(Note that the citations we are making typically work in the context of derived/triangulated
categories, rather than stable $\infty$-categories, but their results are readily adapted
to this latter setting.)

There are actually (at least!) two possible approaches to defining
 this category (for any algebraic stack~$\cX$):  we can first consider 
the abelian category $\QCoh(\cX)$ of quasicoherent sheaves
on~$\cX$  (and even this category admits multiple definitions:  
via small lisse-\'etale sites, small {\em fppf} sites
(this is used e.g.\ in \cite{MR3720855}),
the {\em flat-fppf} site
of~\cite[\href{https://stacks.math.columbia.edu/tag/0786}{Tag 0786}]{stacks-project}
(this is a variant of the small {\em fppf} site, 
in which we consider arbitrary flat morphisms to~$\cX$, but the coverings are taken 
to be~{\em fppf}),
or large {\em fppf} sites, which however lead to equivalent
notions (see e.g.\ ~\cite[\href{https://stacks.math.columbia.edu/tag/07B1}{Tag 07B1}]{stacks-project})),
which is a Grothendieck abelian category, and then form $D\bigl(\QCoh(\cX)\bigr).$

Alternatively, 
we can consider the larger Grothendieck abelian category $\cO_{\cX}\text{-}\Mod$ of
{\em all} sheaves of $\cO_{\cX}$-modules (again on some appropriate site of~$\cX$),
form its associated stable $\infty$-category $D\bigl(\cO_{\cX})$,
and then (attempt to) consider the full subcategory
$D_{\qc}\bigl(\cO_{\cX})$ consisting of  objects whose cohomologies 
are quasicoherent, i.e.\ lie in~$\QCoh(\cX)$.
If we work with small sites, this definition is valid and reasonable,
because the inclusion $\QCoh(\cX) \to \cO_{\cX}\text{-}\Mod$ is exact in
that context.  If we work with big sites, then the inclusion $\QCoh(\cX) \to \cO_{\cX}\text{-}\Mod$
is generally not exact
(see point
(6) of~\cite[\href{https://stacks.math.columbia.edu/tag/06VE}{Tag 06VE}]{stacks-project}),
and so this construction is {\em not} valid; but a more elaborate
version of it can be carried out 
instead~\cite[\href{https://stacks.math.columbia.edu/tag/07B5}{Tag 07B5}]{stacks-project}).
Again, the resulting category is independent, up to canonical equivalence,
of the method of construction. 

There is a canonical $t$-exact functor 
$D\bigl(\QCoh(\cX)\bigr)\to D_{\qc}(\cO_{\cX})$.
We don't know if this is an equivalence, even under reasonable finiteness hypotheses
on~$\cX$.  However, if $\cX$  is quasi-compact and {\em geometric} (meaning
that is has affine diagonal), then this functor does restrict to an equivalence
\numequation
\label{eqn:bounded below qc equiv}
D^+\bigl(\QCoh(\cX)\bigr) \iso D^+_{\qc}(\cO_{\cX});
\end{equation}
this is~\cite[Prop.\ 1.6]{MR3720855}.

The heart of the $t$-structure on either
$D\bigl(\QCoh(\cX)\bigr)$ or $D_{\qc}(\cO_{\cX})$
is the category $\QCoh(\cO_{\cX})$.  
If $\cX$ is a Noetherian algebraic stack, 
then this category is locally coherent, and its full subcategory
of compact objects is precisely the subcategory
$\Coh(\cX)$ of coherent sheaves on~$\cX$.
In particular,
the $t$-structures on both 
$D\bigl(\QCoh(\cX)\bigr)$ and $D_{\qc}(\cO_{\cX})$
are then coherent.
Suppose that $\cX$ is furthermore geometric; then
the equivalence~\eqref{eqn:bounded below qc equiv} 
induces an equivalence between the subcategories of %
coherent objects in 
$D\bigl(\QCoh(\cX)\bigr)$ and $D_{\qc}(\cO_{\cX})$.
We identify these stable $\infty$-categories of coherent objects via this equivalence,
and denote either one by $D^b_{\coh}(\cO_{\cX}).$  This is reasonable,
since these objects (being coherent objects for a coherent $t$-structure) are bounded,
and the heart of $D^b_{\coh}(\cO_{\cX})$ is
precisely~$\Coh(\cX)$ 
(again by the general properties of coherent  $t$-structures).

\begin{df}%
Suppose that $\cX$ is a Noetherian and geometric algebraic stack
(so that, by the preceding discussion, $D^b_{\coh}(\cO_{\cX})$ is unambiguously
defined). %
We then define $\Ind\Coh(\cX) := \Ind D^b_{\coh}(\cO_{\cX})$,
and sometimes refer to it as the stable  $\infty$-category of Ind-coherent
complexes on~$\cX$.
\end{df}

\begin{remark}%
\label{rem:qc via equivariance}
If one works in the $\infty$ (rather than triangulated) categorical
context from the beginning, then it is possible to give 
yet another definition of (a version of) $D\bigl(\QCoh(\cX)\bigr)$, 
via descent from the affine case;
see e.g.~\cite[\S 1.1.3]{GLQCoh}.
(Working $\infty$-categorically is crucial for this definition to work,
since the coherent higher homotopies that are implicit in that context
are required for descent to be a valid process on the derived level.)
This stable $\infty$-category is usually denoted $\QCoh(\cX)$,
and in fact is defined in much greater generality ---
namely for derived prestacks.
This perspective is reconciled with the one described above by~\cite[Prop.\ 1.3]{Hall_2017},
which proves that the resulting category is equivalent to~$D_{\qc}(\cO_{\cX})$. 

We note that the definition
given in Section~\ref{subsec:rigid stacks}
of the stable $\infty$-category  
$\mathbf{D}_{\rm coh}(\mathfrak{X})$
on a rigid analytic Artin stack~$\mathfrak{X}$ is
analogous to the $\infty$-categorical definition of~$D_{\qc}(\cX)$ that we have just described.
We also note that
one can similarly give a definition of %
$\Ind\Coh(\cX)$ in the same style, and so also 
in much greater generality,
see e.g.~\cite{gaitsgory2012indcoherent},\cite[\S 2.3.1]{https://doi.org/10.48550/arxiv.1108.1738}. and~\cite{cautis2024indgeometric}.
\end{remark}

\subsubsection{The case of Ind-algebraic stacks}
Suppose that $\cX$ is an Ind-algebraic stack which admits a description
\numequation
\label{eqn:ind description}
\cX \iso \colim_n \cX_n,
\end{equation}
with the $\cX_n$ being geometric Noetherian algebraic stacks,
and the transition morphisms $\cX_n \to \cX_{n+1}$ being
closed immersions.
These closed immersions then induce $t$-exact functors
$$D^b_{\coh}(\cX_n) \to  D^b_{\coh}(\cX_{n+1}),$$
and thus (by passing to Ind categories)
$$\Ind\Coh(\cX_n) \hookrightarrow  \Ind\Coh(\cX_{n+1}).$$

\begin{df}%
\label{def:IndCoh}
We define $\Ind\Coh(\cX)  := \colim_n \Ind\Coh(\cX_n),$
the colimit being taken in the $\infty$-category whose objects
are cocomplete stable $\infty$-categories
and whose morphisms are continuous functors.
\end{df}

\begin{remark}
\label{rem:IndCoh t-structure}
The pushforward functors $\Ind\Coh(\cX_n) \to \Ind\Coh(\cX_{n+1})$
used to define the colimit in  Definition~\ref{def:IndCoh}
are $t$-exact, and
so $\Ind\Coh(\cX)$  is  naturally equipped with a $t$-structure,
which is furthermore coherent,
and with respect  to which $\Ind\Coh(\cX)$  is left regular.

We define $\Coh(\cX)$ to be  the abelian category  of compact objects
in $\Ind\Coh(\cX)^{\heartsuit}$.  Equivalently,
$\Coh(\cX) :=  \colim_n \Coh(\cX_n)$.
\end{remark}

\begin{remark}
Any two descriptions of $\cX$ of the form~\eqref{eqn:ind description}
are mutually cofinal in an evident sense,
and thus $\Ind\Coh(\cX)$ and $\Coh(\cX)$ are well-defined up to canonical equivalence
independent of the choice of such a description. 
\end{remark}

\begin{remark}
The definition of $\Ind\Coh(\cX)$
given in \cite{gaitsgory2012indcoherent} 
and~\cite[\S 2.3.1]{https://doi.org/10.48550/arxiv.1108.1738},
which we recalled in Remark~\ref{rem:qc via equivariance} above,
can also be applied directly to the Ind-algebraic stack~$\cX$.
We've given the more concrete definition above simply because we find
it easier to think about.
\end{remark}

\subsection{Pro-coherent sheaves on formal algebraic stacks}
Suppose that $\cX$ is a geometric Noetherian formal
algebraic stack, described as an Ind-algebraic stack
via $\colim_n \cX_n \iso \cX$,
where  the $\cX_n$ are geometric  Noetherian algebraic stacks,
and the transition maps $\cX_{n} \to  \cX_{n+1}$
are thickenings, i.e.\ closed immersions inducing  isomorphisms on underlying
reduced  substacks.
We may then define $\Ind\Coh(\cX)$ and $\Coh(\cX)$
following Definition~\ref{def:IndCoh} and Remark~\ref{rem:IndCoh t-structure}.

There are natural ``sheaves'' on $\cX$
that one might want to consider, which do not
arise as objects of $\Coh(\cX)$  according to our  definition.
For example,  if $A$  is an $I$-adically complete  Noetherian ring,
then the usual definition of coherent sheaves on $\Spf A$
yields a category  equivalent (via passage to global sections)
to the category of finitely generated
$A$-modules, whereas according 
to our definition, $\Coh(\Spf A)$  is equivalent (via passage
to global sections) to the subcategory of finitely generated
$A$-modules annihilated by $I^n$ for some~$n$.  In particular,
the {\em structure sheaf} on $\Spf A$, corresponding to the ring~$A$ itself,
is {\em not} a coherent sheaf on~$\Spf A$.
Instead, we identify it with the formal pro-object 
$\quoteslim{n} A/I^n$, and regard it as a {\em pro-coherent sheaf} on~$\Spf A$.

\begin{df}
\label{df:pro-coherent}
The abelian category of {\em pro-coherent sheaves} on $\cX$
is defined to be the formal pro-category of $\Pro \Coh(\cX)$.

We may also consider the stable $\infty$-category $\Pro D^b_{\coh}(\cX)$,
the {\em stable $\infty$-category of pro-coherent complexes} on~$\cX$.
The usual $t$-structure on $D^b_{\coh}(\cX)$ extends to a $t$-structure
on~$\Pro D^b_{\coh}(\cX)$, 
whose heart is (canonically equivalent to)~$\Pro \Coh(\cX)$. 
\end{df}

\begin{example}
\label{ex:dualizing complexes}
If $A$ is a Cohen--Macaulay $I$-adically complete Noetherian ring,
with dualizing module $D$, then $D$ also gives rise to a pro-coherent sheaf
on $\Spf A$.  More generally, if a Noetherian formal algebraic stack
admits a dualizing complex, this can naturally be regarded as an object of $\Pro
D^b_{\coh}(cX).$
We adopt this point of view on dualizing sheaves/complexes in the body of our
notes, since we will want to apply Grothendieck--Serre duality in the context of pro-coherent
sheaves.

We note, though, that if we restrict our duality theory to objects
of $\Ind\Coh(\cX)$, then we can construct an alternative version of the dualizing complex,
which is itself an object of $\Ind\Coh(\cX)$.
Namely, if we write $\cX \iso \colim_n \cX_n,$
and let
$i_{n,n+1}:\cX_n \hookrightarrow \cX_{n+1}$ denote the transition maps
(assumed to be closed immersions) and
$i_n : \cX_n\hookrightarrow \cX$ denote the canonical closed immersions,
then if $\omega_n$ denotes the dualizing complex on~$\cX_n$
(assuming it exists),
the natural isomorphism $\omega_n  \iso i_{n,n+1}^!\omega_{n+1}$  induces
a morphism $(i_{n,n+1})_*\omega_n \to \omega_{n+1}$, and hence a morphism
$(i_n)_*\omega_n \to (i_{n+1})_*\omega_{n+1}$.   If we set
$\omega := \colim_n (i_n)_*\omega_n,$ then $\omega$ is an object
of $\Ind\Coh(\cX)$ which serves as a dualizing sheaf (for objects of $\Ind\Coh(\cX)$).

As a concrete example, consider $A = k[[T]]$ (with its $T$-adic topology). 
Then $k[[T]]$ is a dualizing module for $A$.  On the other hand,
if we write $A = \varprojlim_n k[T]/(T^n)$, then the construction of
the preceding paragraph gives rise to the $A$-module
$\colim_n k[T]/(T^n)$, the transition maps being given by multiplication
by~$T$, which is isomorphic to $k((T))/k[[T]]$.  This module corresponds
to an object of $\Ind\Coh(\Spf A)$, and a consideration of the short exact sequence
$$0 \to k[[T]]\to k((T)) \to k((T))/k[[T]] \to 0$$
shows that  $\RHom_A\bigl(\text{--}, k[[T]]\bigr) $  and $\RHom_A\bigl(\text{--},k((T))/k[[T]]\bigr)$
coincide (up to a shift) on objects of $\Ind\Coh(\Spf A)$
(since $\RHom_A\bigl(\text{--},  k((T))\bigr)$
vanishes on $T$-power torsion $A$-modules).
\end{example}

\section{Localization and Morita theory} %

``Localization'' refers to the idea of
interpreting modules $M$ over some ring $A$ as sheaves (of some sort; for us they
will always be quasicoherent sheaves, or some variant thereof,
such as complexes of quasicoherent sheaves or ind-coherent sheaves)
on an appropriate sort of space
(for us, a scheme or variant  thereof, such 
as an algebraic stack or a formal algebraic stack).
The basic intuition comes from spectral theory: if $A$ is a $k$-algebra,
then an $A$-module $M$ is a $k$-vector space endowed with 
an action of $A$ as a ring of operators.
We can then try to decompose $M$ in terms of its various simultaneous
eigenspaces with respect to the action of~$A$ (if $A$ is  commutative)
or more generally into a direct sum or direct integral of simple $A$-modules
(if $A$ is non-commutative).  We can hope that the collection of simple $A$-modules
has the structure of a space of some sort, and that this decomposition is
effected by spreading $M$ out over the space of simple $A$-modules, with
the fibre of the associated sheaf at a given point corresponding to the multiplicity
space that the measures the contribution of that particular simple module to~$M$.

\subsection{Localizing \texorpdfstring{$A$}{A}-modules over Spec \texorpdfstring{$A$}{A} (the paradigmatic example)}
The most basic example of localization in algebraic geometry is the
theory of quasicoherent sheaves over affine varieties: if~$X = \Spec A$ is an
affine variety, then the functor of global sections
\numequation\label{eqn: QCoh global sections
  localization}\Gamma(X,-):\QCoh(X)\to
\Gamma(X,\cO_X)\text{-}\Mod\end{equation}is an equivalence.
The quasi-inverse is,  of course, given by
the localization functor
$M  \mapsto \cO_X \otimes_A M.$

Another way to describe the global sections functor, when applied to sheaves
of~$\cO_X$-modules, is as $\Hom_{\cO_X}(\cO_X, \text{--}).$
Thus this fundamental example of localization illustrates a general paradigm:
the sheaf $\cO_X$ is naturally an $(\cO_X,A)$-bimodule (in a somewhat trivial manner,
just via the identification of~$A$  with the ring of global sections of~$\cO_X$),
and so we have the pair of adjoint functors
$\cO_X\otimes_A \text{--}$ and $\Hom_{\cO_X}(\cO_X,\text{--})$, which together
induce an equivalence between the category of $A$-modules and  the
category of quasicoherent $A$-modules.  This is an instance of a more general
structure provided by {\em Morita theory}, which we discuss more formally below.

To elaborate even further on this example, we remark that we
can regard any commutative ring $A$
as an ``abstract'' ring  of commuting operators, and an $A$-module  $M$
as a concrete realization of $A$ as a ring of operators (namely, via its operation
on~$M$).  Note then that $X = \Spec A$
can be thought of as parameterizing all the possible systems of simultaneous eigenvalues
that the elements of~$A$ can assume (since the points of $\Spec A$
correspond, via taking the kernel,
 to equivalence classes of morphisms $A \to F$ with $F$  a field,
two morphisms $A\to E$ and $A\to  F$ being equivalent if we can find a common
field extension $K$ of 
$E$ and $F$ such that the induced morphisms $A\to K$ coincide);
hence the designation of  $X$ as the ``spectrum'' of~$A$.
Roughly, the support of the quasicoherent sheaf $\cM$  associated  to $M$
consists of those  systems of simultaneous eigenvalues  that  appear in~$M$.
(More  precisely, the  fibre of $\cM$
at  a point of $X$ is precisely the co-eigenspace of $M$ for the corresponding
system of eigenvalues.)
Thus localization provides an algebraic interpretation of spectral theory.

\begin{example}
\label{ex:affine line}
If we set $A = k[x]$ (for some field~$k$) and apply 
the preceding discussion,
we recover the usual derivation of the spectral theory of a single
operator
via an analysis of the structure of~$k[x]$-modules.
\end{example}

We are interested in other, less obvious, examples of localization.
But in each case, just as in the paradigmatic
example of localizing $A$-modules,
we hope that the structure of the sheaf associated to a module
(its fibres, its support, and so on) will illuminate, and separate out,
various pieces of internal structure that are more obscure when
considered purely in terms of modules.

\subsection{Morita theory}\label{subsec:Morita}
Morita theory, at least as we intend it here,
is a slightly catch-all term for the theory that describes functors,
especially equivalences, between categories of modules
over various rings.   We recall the classical (i.e.\ abelian, i.e.\ module-theoretic)
case of the theory first, before turning to the derived context
(where the ideas of the theory easily extend to provide a flexible
and convenient context for studying functors
that are not necessarily equivalences).

\subsubsection{The abelian case}
One way to phrase the classical (i.e.\ module-theoretic)
version of Morita theory is that if~$A$ and~$B$ are (not
necessarily commutative) associative rings, with respective categories of left
modules $A\text{-}\Mod$ and $B\text{-}\Mod$, then equivalences of categories
$A\text{-}\Mod\to B\text{-}\Mod$ are all of the form \[M\mapsto P\otimes_{A}M,\]
where $P$ is a finitely generated projective $B$-module which
generates $B\text{-}\Mod$ equipped with
an isomorphism $A\isoto\End_B(P)$. A quasi-inverse is given
by \[N\mapsto \Hom_B(P,N).\]
(Of course, Morita theory applies to this adjoint functor as well;
we find that $Q := P^{\vee}:= \Hom_B(P,B)$ --- which is an $(A,B)$-bimodule
--- is necessarily projective as an $A$-module,  
and the quasi-inverse can be rewritten as $N \mapsto N\otimes_B Q.$)

More generally, any functor
$F:  A\text{-}\Mod \to B\text{-}\Mod$
which commutes with colimits is of the form
\numequation
\label{eqn:functor formula}
F: M \mapsto P\otimes_A M
\end{equation}
for some $(B,A)$-bimodule $P$.  (The bimodule $P$ is given by evaluating
the functor on~$A$.  The formula~\eqref{eqn:functor formula}  
is then obtained by evaluating $F$ on some presentation of~$M$.)
The right adjoint of $F$ (which exists since $F$ commutes with colimits) is evidently 
given by $\Hom_B(P,\text{--}).$

If $F$ is furthermore fully faithful, then we see that $A \to \End_B(P)$
must be an isomorphism.  If it is in fact an equivalence,
then $P$ must inherit all the purely categorical properties that $A$ has
in $A\text{-}\Mod$ (compactness, projectivity, being a generator), which gives
the statement of classical Morita theory. 

\begin{example}
\label{ex:matrices}
A standard example of classical Morita theory is given by the 
case $B  =  M_n(A)$  and $P  = A^{\oplus n}$ 
(thought of as length $n$ column vectors)
with its obvious $\bigl(M_n(A),A\bigr)$-bimodule structure.
Since $B \iso P^{\oplus n}$ as a module  over itself, we see
that  $P$ is finitely  generated and projective, and
so $M \mapsto P\otimes_A M =  M^{\oplus  n}$ gives an equivalence
between  the categories $A\text{-}\Mod$ and $M_n(A)\text{-}\Mod$.

Going the other way, if we set $Q = A^{\oplus n}$ (now thought of as
length $n$ row vectors), with  its obvious $\bigl(A,M_n(A)\bigr)$-bimodule structure,
then $N \mapsto Q\otimes_{M_n(A)}N$  gives the quasi-inverse equivalence
between $M_n(A)\text{-}\Mod$ and~$A\text{-}\Mod$.
\end{example}

\begin{example}
\label{ex:localizing matrix modules}
When $A$ is commutative,
the preceding example can be turned into an example of localization,
by writing $X = \Spec A$, and considering the functor
$$N \mapsto \cO_X^{\oplus n}\otimes_{M_n(A)} N,$$
which induces an equivalence of categories
between $M_n(A)\text{-}\Mod$ and  the category of quasicoherent sheaves on~$X$.
\end{example}

Of course we are interested in studying modules 
over non-commutative rings that are more complicated than matrix rings 
over commutative rings.
In such contexts, it is generally not reasonable to expect 
to obtain an equivalence of categories with the category of modules
over a commutative ring, or with the category of quasicoherent
sheaves on a scheme. 
As we will see, it is more plausible to hope for a fully faithful
embedding than an equivalence, and it also makes sense to expand 
the possibilities for the ``space'' $X$ which will carry the
sheaves, e.g.\ by allowing it to be an algebraic stack.

It is also more natural to pass to a derived setting.
Indeed, if one considers the problem of determining when
a functor of the form~\eqref{eqn:functor formula} is
fully  faithful (but  not necessarily an equivalence),
one immediately runs into  questions about the homological
nature of~$F$ which are more naturally considered in the
derived world.  We will thus briefly discuss derived Morita
theory next, with an emphasis on the extra flexibility
it provides (in comparison to the classical abelian setting)
for investigating functors which are not necessarily  
equivalences.

\subsubsection{The case of stable  $\infty$-categories}
\label{subsubsec:derived morita}
If $A,B$ are $E_1$-rings
(see Appendix~\ref{subsec: E1 objects} for the
terminology of $E_1$-rings),
then by~\cite[Prop.\ 7.1.2.4]{LurieHA}, the continuous functors
$\LMod_A\to\LMod_B$ are given by tensoring with $(A,B)$-bimodules. In
particular if 
$A$ and $B$ are static (i.e.\ are associative rings in the usual
sense), then equivalences of stable $\infty$-categories $D(A\text{-}\Mod)\to D(B\text{-}\Mod)$ are all of the form
\numequation\label{eqn:derived Morita}M\mapsto
  T\otimes^{L}_{A}M,\end{equation} where $T$ is a perfect complex of $B$-modules
which generates $D(B\text{-}\Mod)$, and is
equipped with an isomorphism  %
of  $E_1$-rings $A^{\op}\isoto\REnd_B(T)$,
where we write $\REnd_B(T)$ for the $E_1$-ring $\Hom_{D(B\text{-}\Mod)}(T,T)$. %

If we drop the assumption that the perfect complex $T$ generates $D(B\text{-}\Mod)$, but continue
to assume that $A^{\op}\isoto\REnd_B(T)$,
then~\eqref{eqn:derived Morita} still defines a continuous fully
faithful functor $F:D(A\text{-}\Mod)\to D(B\text{-}\Mod)$, with a right adjoint
$G:D(B\text{-}\Mod)\to D(A\text{-}\Mod)$ given by %
\numequation\label{eqn:right
  adjoint to derived Morita}N\mapsto \RHom_B(T,N).\end{equation} More generally, if $T$ is a
perfect complex of $B$-modules %
and $S$ is a perfect complex of
$A$-modules which generates $D(A\text{-}\Mod)$ and is equipped with an isomorphism
$\REnd_A(S)\isoto\REnd_B(T):=\cH^{\op}$, then combining the equivalence
of $D(A\text{-}\Mod)$ with $\cH\text{-}\Mod$ %
given by~$S$ and the continuous fully
faithful functor from $\cH\text{-}\Mod$ to $D(B\text{-}\Mod)$ given by~$T$, we have a continuous fully
faithful functor %
$D(A\text{-}\Mod)\to D(B\text{-}\Mod)$ given by
\[M\mapsto
T\otimes_{\cH}^L\RHom_A(S,M),\] with a right adjoint
$D(B\text{-}\Mod)\to D(A\text{-}\Mod)$ given by %
\[N\mapsto
S\otimes_{\cH}^L\RHom_B(T,N).\]

There are, of course, many generalisations of this framework. In
particular, in the body of the text we consider examples where
$D(B\text{-}\Mod)$ is replaced by an appropriate derived category of sheaves
on a stack, as in the examples below (see also Remark~\ref{rem: expect Linfty is a kernel}).
\begin{rem}
  \label{rem: notation for derived Morita has Hecke}The notation $\cH$ for the derived endomorphism ring is intended to be suggestive of a (derived) Hecke algebra. See for
  example~\cite[(1.1)]{hellmann2020derived} for a functor from a
  category of representations to a category of sheaves defined in this
  fashion.
\end{rem}
\begin{rem}\label{rem: E1 versus dga}
  In the literature, derived Hecke algebras are often thought of as
  dgas, rather than $E_1$-rings. The correspondence between these
  notions is as follows: by \cite[Prop.\ 7.1.4.6]{LurieHA}, if $A$ is
  a (static) commutative ring, then the $\infty$-category of
  $E_1$-algebras over $A$ can be identified with the $\infty$-category
  of differential graded $A$-algebras. For example, if we take $A=\Fp$,
  Schneider's derived Hecke algebra \cite{MR3437786} can equally well
  be regarded as an $\Fp$-dga or an $\Fp$-$E_1$-algebra.
\end{rem}

\subsection{Non-commutative phenomena}
\label{subsec:non-comm}
There are some phenomena that can occur in the theory
of modules over non-commutative rings that don't occur in
the commutative case, 
which have to be taken into account when trying to localize
in the non-commutative setting.
We briefly recall some of these,  and explain why allowing 
localization over algebraic stacks (rather than just schemes)
helps to ameliorate the situation. 

There is one basic phenomenon that modules over a non-commutative 
ring $A$ can exhibit which does not occur in the commutative context,
namely, if $S_1$ and $S_2$ are non-isomorphic simple $A$-modules, then
although $\Hom_A(S_1,S_2)$ evidently vanishes, the higher Exts,
i.e.\  $\Ext^i_A(S_1,S_2)$  for $i> 0$, need not vanish.  This means that if the category $A\text{-}\Mod$ is mapped functorially into a category of sheaves in a way that  not only preserves $\operatorname{Hom}$ spaces but also induces isomorphisms on higher $\operatorname{Ext}$ groups, then the images of $S_1$ and $S_2$ cannot have disjoint support.  %

Remembering that the simple objects in the category of quasicoherent sheaves
on a Noetherian\footnote{We make this assumption just to ensure  that the
discussion is  not overburdened with technical complications.} 
scheme are precisely the skyscraper sheaves at closed points,
so that non-isomorphic simple objects have disjoint support,
we see for example that if $A$ is a non-commutative ring admitting 
simple  modules $S_1$ and  $S_2$ with this property (i.e.\ such
that $\Hom_A(S_1,S_2) =  0$ but $\RHom_A(S_1,S_2) \neq 0$), %
then the category  $A\text{-}\Mod$ certainly  cannot  be {\em equivalent}
to  the category of  quasicoherent sheaves on  a Noetherian scheme.     

\begin{remark}
\label{rem:D-module note}
The preceding discussion  helps to explain why in Beilinson--Bernstein localization theory,
which does give an equivalence (between certain categories of Lie
algebra representations and the category of algebraic $\cD$-modules on a flag
variety),   
it is quasicoherent $\cD$-modules 
that appear, rather  than  merely quasicoherent $\cO$-modules. 
The point is that the Riemann--Hilbert correspondence shows
that $\cD$-modules behave like perverse sheaves,
in which (to first approximation) the simple objects are described 
by their support (a stratum in some stratification of the ambient scheme~$X$),
and there are interesting extension classes between sheaves supported on different strata. 
(The most basic example is the sequence $$0 \to j_!\C \to  \C  \to  i_*\C \to 0,$$
with $j:\A^1\setminus \{0\}  \to  \A^1$,  resp.\ $i:\{0\} \to  \A^1$ the natural open, 
resp.\ closed,  immersion.)
\end{remark}

If we only want fully faithful functors, the various (im)possibilities 
are less clear,  since then a simple $A$-module  needn't  map  to a simple 
quasicoherent sheaf (it might just be that the proper quasicoherent
subsheaves  of  its image are not themselves in  the image of the functor).
If we  pass to the  derived context, the situation becomes yet murkier,
because then the notions  of {\em simple module} and {\em submodule}  aren't
defined.
Nevertheless, the difference in behaviour  between the non-commutative
and commutative  situations creates a tension which makes it
difficult to construct interesting localization functors 
for  modules  over non-commutative 
rings. One way to lower this tension is by allowing the target category to  be
a category of quasicoherent sheaves on an algebraic stack~$\cX$, as we now explain.

If $x \in |\cX|$ is a closed point of (the  underlying topological space of) an
algebraic stack~$\cX$, then it typically does not correspond to a closed immersion
$\Spec k \hookrightarrow \cX$ for  some  field~$k$, but rather (under mild
finiteness hypotheses, such as having quasi-compact
diagonal; see~\cite[\href{https://stacks.math.columbia.edu/tag/06RD}{Tag 06RD}]{stacks-project},
as well as the  discussion of~\cite[App.~E]{emertongeepicture})
to a closed immersion
$\cZ  \hookrightarrow  \cX$ with $\cZ$ being a one-point {\em gerbe}.
If, for example,  our stack~$\cX$ is of finite type over an
algebraically  closed  field~$k$, then $\cZ = [\Spec k/G]$ for some
algebraic group over~$G$.

Suppose then, to fix ideas, that we have a closed embedding
$[\Spec k/G] \hookrightarrow \cX$  for some field $k$  and some
algebraic group $G$  over~$k$,  giving  rise to the closed  point~$x \in |\cX|$.
Representations  of~$G$ then correspond  to quasicoherent  sheaves
on~$\cX$  that  are  supported at~$x$; 
irreducible representations in particular correspond to
simple objects in the category of quasicoherent sheaves
on~$\cX$. 
In examples, we can often produce non-split extensions
of these sheaves by considering the higher infinitesimal neighbourhoods of
$x$ in~$\cX$.

\begin{example}
Consider  $\cX  = [\A^1/\Gm]$,  with $\A^1  = \Spec k[x]$,  $\Gm = \Spec k[t^{\pm 1}]$,
and $\Gm$  acting on~$\A^1$  via $t\cdot x = tx$.
Let  $\cI  \subset \cO_{\cX}$ denote the  ideal  sheaf  cutting out the closed  
point $\{ x  =0\}.$  Then we have a  non-split short exact sequence
of coherent sheaves
$$0 \to k(-1) \cong 
\cI/\cI^2  \to \cO_{\cX}/\cI^2 \to  \cO_{\cX}/\cI \cong k(0) \to 0$$ 
(where as usual we identify $\Gm$-representations with  graded $k$-modules).
\end{example}

In short, in the category of coherent sheaves on algebraic stacks,  
it is possible to find non-isomorphic simple
objects that admit non-split extensions.  This makes 
categories of coherent sheaves on algebraic stacks amenable targets 
for localization functors defined
on categories of modules 
over  non-commutative rings.

\subsection{Examples}\label{subsec: Morita examples}We now give some
examples of Morita theory and localization, building up from some very
simple ones to some which occur in the representation theory of
$\GL_2(\Qp)$. We work throughout over some commutative ring~$k$.

\begin{example}\label{ex: Morita rank 2}
  As we %
saw in Example~\ref{ex:affine line}
above, the category $k[x]\text{-}\Mod$ is equivalent to
  $\QCoh(\A^1)$; this Morita equivalence is mediated by the structure
  sheaf~$\cO_{\A^1}$ (which is of course a module over its ring of
  global sections, namely $k[x]$). Slightly less trivially, as a special
 case of Example~\ref{ex:matrices},
these categories are equivalent to the categories of left modules for the
  matrix algebra \[
    \begin{pmatrix}
      k[x] & k[x]\\k[x] & k[x]
    \end{pmatrix},
  \]via the bimodule $\cO_{\A^1}\oplus\cO_{\A^1}$.
\end{example}

\begin{example}
\label{ex: Morita model for Iwahori Hecke}
  Consider instead the algebra
  \[ A:= \begin{pmatrix} k[x] & k[x]\\xk[x] & k[x]
    \end{pmatrix}.
  \]
Its category of left modules is not Morita equivalent to the
  category of modules over any commutative ring. However, there is a
  fully faithful functor from $D(A\text{-}\Mod)$ to $D(\QCoh(\cX))$,
  where \[\cX:=[(\Spec k[x,y]/(xy))/\Gm],\] with the $\Gm$-action
  given by\footnote{The appearance of $t^2$ rather than $t$ in the formula
for the action is of no intrinsic importance for this example; its only
significance is that it yields an example that connects directly with
an instance of categorical Langlands, as we note below.}
\[t\cdot(x,y)=(x,t^2y),\] and the functor being given by
  (derived!) tensoring with \[P:=\cO_{\cX}\oplus (\cO_{\cX}/y).\] %
  In order for this
  to make sense, we need~$P$ to be a right $A$-module. We claim that
  in fact \[\REnd_{\cO_{\cX}}(P)=\End_{\cO_{\cX}}(P)=A^{\op},\] so that the
  functor $P\otimes_A-$ is %
fully faithful.

  In general, for an object $X \oplus Y$ in an abelian category,
we have $$\REnd(X\oplus Y) =
\begin{pmatrix} \REnd(X) & \RHom(Y,X) \\ \RHom(Y,X) &
\REnd(Y)\end{pmatrix},$$
acting on the left on $X\oplus Y$, thought  of as column vectors. 
Then $$\REnd(X\oplus Y)^{\op} =
\begin{pmatrix} \REnd(X) & \RHom(X,Y) \\ \RHom(Y,X) &
\REnd(Y)\end{pmatrix},$$
acting on the right on $X\oplus Y$ (now thought of as row vectors).

In our example, $\cX$ is an affine scheme modulo a linearly reductive group,
in
fact the torus~$\Gm$,
and so the global sections functor on~$\cX$ is exact; equivalently,
$\cO_{\cX}$ is projective, 
and $\RHom_{\cO_{\cX}}(\cO_{\cX},\text{--}) = \Hom_{\cO_{\cX}}(\cO_{\cX}, \text{--})$
is simply the functor of global sections.  In particular
$$\RHom_{\cO_{\cX}}(\cO_{\cX},\cO_{\cX}) = 
\RHom_{\cO_{\cX}}(\cO_{\cX},\cO_{\cX}/y) = 
k[x].$$
  To compute $\RHom_{\cO_{\cX}}(\cO_{\cX}/y,\text{--}),$
we use the %
projective resolution
  \[\dots\stackrel{y}{\to}\cO_{\cX}(-2)\stackrel{x}{\to}\cO_{\cX}(-2)\stackrel{y}{\to}\cO_{\cX}\to\cO_{\cX}/y\to
    0.\] (Here as usual our notation for invertible sheaves on the
  quotient stack $\cX$ corresponds to our notation for $\Z$-graded
  modules, with the convention that $\Gr^iM(n)=\Gr^{i+n}M$.)

  We have $\Hom_{\cO_{\cX}}(\cO_{\cX}(-n),\cO_{\cX})=k[x]$ if $n=0$
  and $=k\cdot y^n$ if $n>0$, and similarly
  $\Hom_{\cO_{\cX}}(\cO_{\cX}(-n),\cO_{\cX}/y)=k[x]$ if $n=0$ and $=0$
  if $n>0$. It follows that 
  $\RHom_{\cO_{\cX}}(\cO_\cX/y,\cO_\cX)$ is computed by the complex
  \[0\to k[x]\to k\stackrel{0}{\to} k\isoto k\stackrel{0}{\to}
    k\to\dots \](where the morphism $k[x]\to k$ is given by evaluation
  at~$0$) and that $\RHom_{\cO_{\cX}}(\cO_\cX/y,\cO_\cX/y)$ is  computed by the
  complex \[0\to k[x]\to 0\to 0\to\dots,\] %
  so that
  $$\RHom_{\cO_{\cX}}(\cO_{\cX}/y,\cO_{\cX})
  = \Hom_{\cO_{\cX}}(\cO_{\cX}/y,\cO_{\cX}) = x k[x]$$
  and that
  $$\RHom_{\cO_{\cX}}(\cO_{\cX}/y,\cO_{\cX}/y)
  = \Hom_{\cO_{\cX}}(\cO_{\cX}/y,\cO_{\cX}/y) =  k[x].$$
   Thus indeed $\REnd_{\cO_{\cX}}(P)^{\op} = A$.

Continuing with this example a little more, writing
$M = k[x]\oplus xk[x]$ (thought
of as column vectors, with the evident left action of~$A$) 
and
$N = k[x]\oplus k[x]$ (again thought
of as column vectors, with the evident left action of~$A$) 
then $A = M \oplus N$ as left $A$-modules,
so that  $M$ and $N$ are  both projective  as $A$-modules,
and  $P\otimes_A M  = \cO_{\cX}$,
while  $P\otimes_A N = \cO_{\cX}/y$.

The evident inclusion 
\numequation
\label{eqn:M into N}
M\hookrightarrow  N
\end{equation}
has cokernel $S_1 := 0 \oplus k$  (whose $A$-module structure comes by thinking
of elements as column vectors and having $x$ act via $0$),
while the the morphism
\numequation
\label{eqn:N into M}
N \buildrel x \over\longrightarrow M
\end{equation}
has cokernel $S_2:= k \oplus 0$ (which again has the $A$-module structure given by regarding
its elements as column vectors and having $x$ act via $0$).
The two $A$-modules $S_1,S_2$ are simple and not isomorphic (and, up to isomorphism,
are precisely the two simple $A$-modules on which the the element $x$ in
the centre $k[x]$ of $A$
acts via~$0$).   The quotient $M/xM$ is a non-split extension of $S_2$ by $S_1$,
while $N/xN$ is a non-split 
extension of $S_1$ by~$S_2$.

If we pass to sheaves on~$\cX$, the inclusion~\eqref{eqn:M into N}
corresponds to the surjection $\cO_{\cX}\to \cO_{\cX}/y,$ and thus $S_1$
corresponds to the cone of this morphism, which is $(\cO_{\cX}/x)(-2).$
The morphism~\eqref{eqn:N into M} corresponds to the morphism
$\cO_{\cX}/y \buildrel x \over\longrightarrow \cO_{\cX}$, and thus $S_2$
corresponds to the cone of this morphism, which is $\cO_{\cX}/x.$
Thus we see an example of the possibility discussed in
Section~\ref{subsec:non-comm}, i.e.\ of two non-isomorphic simple modules,
admitting a non-split extension between them,
corresponding to coherent sheaves which are non-isomorphic but have
the same support.
\end{example}

\begin{rem}\label{rem: how Morita example relates to Iwahori Hecke}%
Example~\ref{ex: Morita model for Iwahori Hecke}
is closely
  connected to the $\ell$-adic ($\ell\ne p$) representation theory of
  $\PGL_2(\Qp)$, see for example~\cite[Rem.\
  4.43]{hellmann2020derived}.   %
  Indeed the Iwahori
  Hecke algebra~$\cH_I$ for $\PGL_2(\Qp)$ is a matrix algebra 
 \[
 \cH_I= \begin{pmatrix}
    \Q_{\ell}[T] & \Q_{\ell}[T] \\ (T^2-(p+1)^2)\Q_{\ell}[T] & \Q_{\ell}[T]
  \end{pmatrix}
\] and there is a fully faithful functor (determined by the coherent
Springer sheaf) from the derived category of
left $\cH$-modules to the category of ind-coherent sheaves on the
moduli stack $\cX_{2}^{\unip}$ of 2-dimensional unipotently ramified
Weil--Deligne representations with inverse cyclotomic
determinant.

This stack has three irreducible components; the
unramified component, and two ``Steinberg'' components which are
related by a quadratic twist. The two Steinberg components do not
intersect, but each intersects the unramified component, at
respectively the representation $1\oplus \varepsilon^{-1}$ and its
quadratic twist. In a neighbourhood of either one of these intersection
points we can identify $\cX_2^{\unip}$ with (a corresponding neighbourhood of the singular point
on) the stack~$\cX$ above. 
The locus $y=0$ then corresponds to the Steinberg component, and~$P$
corresponds to the coherent Springer sheaf (see e.g.\ \cite[Ex.\ 4.4.4]{zhu2020coherent}).

The module $M$ corresponds to the Steinberg representation (or its quadratic twist),
whose associated
sheaf is $\cO_{\cX}/y$, the structure sheaf of the Steinberg component; while
$N$ corresponds to the trivial 
representation (or its quadratic twist),
whose associated sheaf is $(\cO_{\cX}/y)(-2)[1]$,
a shift of a twist of the structure sheaf of the Steinberg component.
\end{rem}%

\begin{rem}
  \label{rem: Iwahori Morita example shows ind coherent}Note that the
object~$P$ of Example~\ref{ex: Morita model for Iwahori Hecke}
  is not a perfect complex, but it is coherent; thus already in this
  example we see that it is natural to consider ind-coherent sheaves
  on stacks of Langlands parameters, rather than quasicoherent sheaves.
\end{rem}

The following variant on Example~\ref{ex: Morita model for Iwahori Hecke}
plays a similar role in the $p$-adic
Langlands correspondence for $\GL_2(\Qp)$, as described in
Section~\ref{subsec:GL2 Qp proof sketches}.

\begin{example}\label{ex: Morita example for GL2 Qp}
  Let $A$ be a (commutative) ring, and consider the graded ring
  \[B:= A[x,y],\] where the grading is defined by giving $x$ and $y$
  weights $2$ and $-2$ respectively.  Let $\cB$ denote the category of
  graded $B$-modules (equivalently, quasicoherent sheaves on the
  quotient stack $[\Spec B/\Gm]$). %
  If $M$ is a graded module, and $n \in \Z$, then as usual we let
  $M(n)$ denote the graded module whose underlying module coincides
  with that of~$M$, and for which $\Gr^i M(n) = \Gr^{i+n} M.$

  We consider the module $P := B(1) \oplus B(-1),$ and define
  $R := \End_{\cB}(P)^{\op};$ then $R$ acts %
{\em on the right}
  on~$P$.
  Concretely, we find that
\begin{multline*}
  R = \begin{pmatrix} \End\bigl(B(1)\bigr) &  \Hom\bigl(B(1),B(-1)\bigr) \\
\Hom\bigl(B(-1),B(1)\bigr)  & \End\bigl(B(-1)\bigr) \end{pmatrix}
=  \begin{pmatrix}  \Gr^0 B  & \Gr^0 B(-2) \\
\Gr^0 B(2) & \Gr^0 B \end{pmatrix} 
\\
=  \begin{pmatrix}  \Gr^0 B  & \Gr^{-2} B \\
\Gr^2 B & \Gr^0 B \end{pmatrix} 
    =\begin{pmatrix} A[xy] & y A[xy] \\ x A[xy] & A[xy] \end{pmatrix}
\end{multline*}
  acting on $P$ via the usual right multiplication of matrices on row
  vectors.

  Certainly $P$ is a projective object of $\cB$.  Slightly less
  obviously, we have the following lemma.

  \begin{lemma} $P$ is projective as a right $R$-module.
    \label{lem:P is projective}
  \end{lemma}
  \begin{proof}
    We first note that as a right module over itself, we have
    $R = R_1 \oplus R_2,$ where $R_1 = A[xy] \oplus y A[xy]$ (the
    first row of~$R$ regarded as above as a matrix order) and
    $R_2 = x A[xy] \oplus A[xy]$ (the second row of~$R$).  Each of
    $R_1$ and $R_2$ is projective (being a direct summand of~$R$).

    We next note that
    $P = \bigoplus_{m = -\infty}^{\infty} \Gr^{1+2m} P;$ since $R$
    acts as graded endomorphisms of~$P$, this is an isomorphism of
    $R$-modules.  If $1+2m \geq 1$ (i.e.\ $m \geq 0$), then
    $\Gr^{1+2m} P = x^{m} A[xy] \oplus x^{m+1} A[xy] \iso A[xy] \oplus
    x A[xy] = R_2,$ while if $1 + 2m \leq -1$ (i.e.\ $m \leq -1$), then
    $\Gr^{1+2m} P = y^{-m-1} A[xy] \oplus y^{-m} A[xy] \iso A[xy]
    \oplus yA[xy] = R_1.$ Thus $P$ is a direct sum of projective
    $R$-modules, and so is itself a projective~$R$-module.
  \end{proof}
  We have the following immediate consequence.

  \begin{prop}
    \label{prop:morita}
    The functor $N \mapsto P \otimes_R N$ induces a fully faithful and
    exact embedding from the abelian category of~$R$-modules
    into~$\cB$, whose essential image is equal to the full subcategory
    $\cC$ of $\cB$ generated by~$P$.
  \end{prop}

In fact, since $P$ is projective over~$B$, so that  $\RHom_B(P,P) = R$,
the preceding proposition has a derived analogue (following
the discussion of~\ref{subsubsec:derived morita}).

  \begin{prop}
    \label{prop:derived morita}
    The functor $N \mapsto P \otimes_R N$ induces a fully faithful and
    $t$-exact embedding of stable $\infty$-categories with $t$-structures
    $D(R) \to D(\cB)$,
    whose essential image is equal to the full subcategory
    of $D(\cB)$ generated by~$P$.
  \end{prop}
\end{example}

\subsubsection{Characterizing an image}
As already remarked, Example~\ref{ex: Morita example for GL2 Qp}
has an application to the categorical $p$-adic 
Langlands correspondence in the $\GL_2(\Q_p)$-case (see Section~\ref{subsec:semi-decomp}).
In fact, in this application we will study a variant of
this functor, in which we work with stable $\infty$-categories
and also impose a support condition.

We now describe our setup precisely.
To begin with, we maintain the notation of Example~\ref{ex: Morita example for GL2 Qp},
but furthermore assume that the ring $A$ is Noetherian.  
We also introduce additional notation.
Namely, 
we let $\cR$ denote the abelian category of $R$-modules,
and let $\cR_0$ denote the  full  subcategory of $\cR$ consisting 
of those modules, each element of which is annihilated by some power of~$xy$.
Similarly, we let $\cB_0$ denote the subcategory of the category $\cB$ of graded $B$-modules
consisting of those modules satisfying the same condition, i.e.\ each element
is annihilated by some power of~$xy$.

Proposition~\ref{prop:second Coh equivalence}
(and its evident variant for graded rings) shows that
the embedings $\cR_0 \hookrightarrow R\text{-}\Mod$ and $\cB_0 \hookrightarrow \cB$ 
induce equivalences
\numequation
\label{eqn:R equiv}
D^+(\cR_0) \iso D^+_{(xy)}(R)
\end{equation}
and
\numequation
\label{eqn:B equiv}
D^+(\cB_0) \iso D^+_{(xy)}(\cB).
\end{equation}
(The subscript~$(xy)$ here denotes the full subcategory
consisting of those objects each of whose cohomologies
is annihilated elementwise by some power of~$xy$.)

\begin{lemma}
\label{lem:equivalence}
The equivalences~{\em \eqref{eqn:R equiv}}
and~{\em \eqref{eqn:B equiv}}
extend to equivalences
$D(\cR_0) \iso D_{(xy)}(R)$
and
$D(\cB_0) \iso D_{(xy)}(\cB)$.
\end{lemma}
\begin{proof}
This will follow from Proposition~\ref{prop:full faithfulness}~(2),
once we show that the adjoint functors to the inclusions $\cR_0 \hookrightarrow R\text{-}\Mod$
and $\cB_0 \hookrightarrow \cB$ have bounded cohomological dimension.
However, in either case, the derived functor of this adjoint is computed 
via $X \mapsto \colim_n ( X \buildrel (xy)^n \over \longrightarrow X),$
and so has a cohomological amplitude of~$1$.
\end{proof}

The $t$-exact fully faithful functor $N\mapsto P\otimes_R N$
of Proposition~\ref{prop:derived morita}
evidently restricts to a $t$-exact functor
fully faithful functor $D_{(xy)}(R) \to D_{(xy)}(\cB),$
which by Lemma~\ref{lem:equivalence}
we may also regard as a $t$-exact fully faithful functor
$D(\cR_0) \to D(\cB_0)$, 
We will characterize the image of this functor
in terms of a certain semi-orthogonal decomposition,
which we now explain.

The adjoint functor $\Hom_{\cB}(P,\text{--}):  D_{(xy)}(\cB) \to D_{(xy)}(R)$
admits the more explicit  description (remembering that  $P := B(1) \oplus B(-1)$)
$$M \mapsto \Gr^{-1} M \oplus \Gr^1 M.$$
(The right hand side being equipped with its evident  $R$-structure:  each summand
is naturally an $A[xy]$-module; multiplication by $x$ takes the first summand to
the second; and multiplication by $y$ takes the second summand to the first.)

\begin{prop}
\label{prop:image}
Let $\cK$ denote 
the full subcategory of $D_{(xy)}(\cB)$  consisting of those objects 
$M$ such that $\Gr^{-1} M = \Gr^1  M = 0$.
Then the image of $P \otimes_R \text{--}$ is equal
to ${}^{\perp}\cK$.
\end{prop}
\begin{proof}
This follows from Lemma~\ref{lem:orthogonality}~(6).
\end{proof}

\subsubsection{A set of generators for $\cK$}
\label{subsubsec:set of generators}
For applications to categorical $p$-adic Langlands, 
we wish to have an alternative description of the kernel~$\cK$,
in terms of a set of generators. 
In order to do this, we change the context slightly,
replacing the unbounded stable infinity categories
$D(\cR_0)$ and $D(\cB_0)$ by their regularizations
$\Ind\Coh(\cR_0)$ and $\Ind\Coh(\cB_0)$, following the procedure 
described in Section~\ref{subsec:compact and coherent}.
The functor $P\otimes_R \text{--}$ and its adjoint
$\Hom_{\cB}(P,\text{--})$ are both exact and continuous (the
latter because $P$ is finitely generated and projective as
a $B$-module), and preserve compact objects.  Thus they
induce an adjoint pair of functors between $\Ind\Coh(\cR_0)$
and $\Ind\Coh(\cB_0)$,
and the results
proved above in the context of $D(\cR_0)$ and $D(\cB_0)$ carry
over to the Ind-coherent context.

The category
$\cB_0$ splits into a product
$\cB_0^+ \times \cB_0^-$, where $\cB_0^+$ (resp. $\cB_0^-$) consists of objects
whose graded pieces are supported purely in even (resp.\ odd)
degrees. 
Correspondingly, $\Ind\Coh(\cB_0)$ decomposes as a product
$\Ind\Coh(\cB_0^+) \times \Ind\Coh(\cB_0^-).$

The functor $P\otimes_R \text{--}$ evidently factors through 
$\Ind\Coh(\cB_0^-)$, and we let $\cD_0^-$ denote its essential image in this category. 
We also consider the analogous
functor $B\otimes_{A[xy]} \text{--}: \Ind\Coh(\cA_0)
\to \Ind\Coh(\cB_0^+)$, where $\cA_0$  denotes the subcategory
of $A[xy]\text{-}\Mod$ consisting of modules annihilated elementwise by some power of~$xy$.
We let $\cD_0^+$ denote the essential image in $\Ind\Coh(\cB_0^+)$ of this latter functor.

We regard $A[x]$ as a graded $B$-module via the isomorphism $B/(y) \iso A[x],$
and similarly for $A[y].$
Similarly we regard $A$ as a graded  $B$-module via the isomorphism
$B_i/(x,y) \iso A.$  All three of $A[x]$, $A[y]$, and~$A$
are then objects of~$\cB_{0}$.

Using the graded resolution
\numequation
\label{eqn:A[x] res}
0 \to B(2) \buildrel y \over \longrightarrow  B \to A[x]\to 0,
\end{equation}
we see (for any choice of twist $n$) that
$\RHom_{\cB}\bigl(A[x], A[x](n)\bigr)$ is computed by
the complex (in cohomological degrees %
$0$ and $1$)
$$\Gr_0(A[x](n)) \buildrel 0 \over \longrightarrow \Gr_0(A[x](n-2)),$$
while 
$\RHom_{\cB}\bigl(A[x], A[y](n)\bigr)$ is computed by
the complex (in cohomological degrees $0$ and $1$)
$$0 \to \Gr^0 A(n-2).$$

Similarly, we find that
$\RHom_{\cB}\bigl(A[y], A[y](n)\bigr)$ is computed by
the complex (in cohomological degrees $0$ and $1$)
$$\Gr_0(A[y](n)) \buildrel 0 \over \longrightarrow \Gr_0(A[y](n+2)),$$
while 
$\RHom_{\cB}\bigl(A[y], A[x](n)\bigr)$ is computed by
the complex (in cohomological degrees $0$ and $1$)
$$0 \to \Gr^0 A(n+2).$$

Consequently, we find that the collections of objects
$\{ A[x](-2), A[x](-4), \dots\}$ and $\{A[y](2), A[y](4), \dots\}$
each form a weakly exceptional collection in~$\Ind\Coh(\cB_{0}^+)$ in the sense
of Definition~\ref{defn: weakly exceptional collection}, 
and are furthermore mutually orthogonal.
By Lemma~\ref{lem: weakly exceptional gives
    semiorthogonal},
 they induce a semiorthogonal decomposition of
  the full subcategory $\cE_0^+$ of $\Ind\Coh(\cB_{0}^+)$ that they generate.

Similarly, the collections
$\{ A[x](-3), A[x](-5), \dots\}$ and $\{A[y](3), A[y](5), \dots\}$
each form a weakly exceptional collection in~$\Ind\Coh(\cB_{0}^-)$, 
and are mutually orthogonal.
 They similarly induce a semiorthogonal decomposition of
  the full subcategory $\cE_0^-$ of $\Ind\Coh(\cB_{0}^-)$ that they generate.

\begin{prop}
\label{prop:semiorthogonal decomposition}
\leavevmode
\begin{enumerate}
\item
$\cD_0^+$ is the cocomplete full subcategory of $\Ind\Coh(\cB_0^+)$
generated by the object $B/(xy)$.
\item
$\cE_0^+ = (\cD_0^+)^{\perp}.$
\item
$\cD_0^-$ is the cocomplete full subcategory of $\Ind\Coh(\cB_0^+)$
generated by the objects $B/(xy)(1)$ and~$B/(xy)(-1)$.
\item
$\cE_0^- = (\cD_0^-)^{\perp}.$
\end{enumerate}
\end{prop}
\begin{proof}
The cocomplete stable $\infty$-category 
$\Ind\Coh(\cA_0)$ is generated by $A := A[xy]/(xy)$.  Thus
$\cD_0^+$, which is defined to be the image of this category under
the functor $B\otimes_{A[xy]} \text{--}$, is generated by the image of~$A$,
namely~$B/(xy)$.  This proves~(1).

Lemma~\ref{lem:orthogonality},
together with the Ind-coherent analogue of Proposition~\ref{prop:image},
shows that 
$(\cD_0^+)^{\perp}$ consists of those objects $M$ for which $\Gr^0 M = 0$. 
On the other hand,
Lemma~\ref{lem: generated subcategory is smallest} shows that $\cE_0^+$
can be described as the smallest cocomplete stable subcategory of
$D(\cB_0^+)$ which contains the set of objects~$X$. 
Since $\Gr M_0 = 0$ for each object in~$X$,
we thus find that $\cE_0^+ \subseteq (\cD_0^+)^{\perp}.$

Now, a consideration of (the twist by $-2$ of) the presentation~\eqref{eqn:A[x] res},
along with the fact that $\Gr M_0 = 0$ for objects of $(\cD_0^+)$,
shows that $\RHom_{\cB}\bigl(A[x](-2), M\bigr) = \Gr^2 M$
for objects~$M$ of~$(\cD_0^+)^{\perp}$.
In particular, $\RHom_{\cB}\bigl(A[x](-2),A[x](-2)\bigr) = A.$
Thus $N \mapsto A[x](-2) \otimes_A N$ is a fully faithful functor
$\Ind\Coh(A) \to (\cD_0^+)^{\perp},$ with adjoint given by $M \mapsto \Gr^2 M$.
In particular, the cone of $A[x](-2)\otimes_A \Gr^2 M \to M$
has vanishing~$\Gr^2$.  Continuing inductively with a consideration
of the various $A[x](-n)$, and then also of the $A[y](n)$ (for positive even
values of $n$),
we find that if $M$  is an object of $(\cE_0^+)^{\perp} \cap (\cD_0^+)^{\perp}$,
then $\Gr^i M = 0$ for all
values of~$i$, and thus that~$M = 0$.  Consequently, $\cE_0^+ = (\cD_0^+)^{\perp}$,
proving~(2).

The proofs of~(3) and~(4) are entirely analogous to those of~(1) and~(2),
and so we omit the details.
\end{proof}

\section{Topological abelian groups, rings and modules}
\subsection{Functional analysis}
Let $\{G_n\}$ be a sequence of linearly topologized\footnote{Meaning that $0$ admits
a neighbourhood basis of subgroups.}  Hausdorff topological abelian groups,
equipped with closed embeddings $G_n \hookrightarrow G_{n+1}$.
We can then define the inductive limit $G := \varinjlim_n G_n,$
endowed with the inductive limit linear topology (i.e.\ the topology
for which a subset is open if and only if its inverse image in
each~$G_n$ is open);
this is again a linearly topologized topological abelian group.
We recall some basic properties that $G$ satisfies; these are analogues 
of standard properties of inductive limits of sequences of locally convex 
topological spaces in various contexts ($LF$-spaces, spaces of compact type,
\ldots), and our arguments are an adaptation of
the arguments of~\cite[Ch.~VII]{MR0350361}.

\begin{remark}
\label{rem:inductive limits}\leavevmode
\begin{enumerate}
\item
To simplify the notation, we regard each $G_n$ as a subset of $G$ 
via the canonical embedding.  We will see below that this embedding
of groups is furthermore a topological embedding, so this should cause
no confusion.

\item
By definition, a subgroup $U$ of $G$ is open if each pullback $G_n \cap U$
is an open subgroup of $G_n$ for each~$n$.
A typical way to construct such open subgroups is to be given
an open subgroup $V_n$ of $G_n$ for each~$n$,
and to define $U_n := \sum_{i=1}^n V_i$.
Then $U_n$ is again an
open subgroup of $G_n$, and evidently $U_n \subseteq U_{n+1}$,
so that $U := \bigcup_n U_n = \sum_{i=1}^{\infty} V_i$
is a subgroup of~$G$.  Since $U_n \subseteq G_n \cap U,$
we see that $U$ is an open subgroup of~$G$.
If the $V_n$ satisfy the additional condition 
$G_n \cap V_{n+1} \subseteq V_n$,
then one sees that $U_n = G_n \cap U_{n+1}$,
and thus $U_n = G_n \cap U.$

If we are given open subgroups $V_n$, then we may inductively construct
open subgroups $V_n' \subseteq V_n$ such that furthermore $G_n \cap V'_{n+1}
\subseteq V_n'$.  Namely, assuming $V'_i$
is constructed, for $i = 1,\ldots,n$,
and and since $G_n \hookrightarrow G_{n+1}$ is a topological
embedding, we may find an open subgroup $W_{n+1}$ 
such that $G_n \cap W_{n+1} \subseteq V'_n.$  We then set $V_{n+1}' := 
V_{n+1} \cap W_{n+1}$.
\end{enumerate}
\end{remark}

\begin{lemma}
\label{lem:discreteness}
If $X \subseteq G$ has the property that
$X \cap G_n$ is finite for each~$n$,
then  $X$  is closed and discrete as a subset of~$G$.
\end{lemma}
\begin{proof}
We first show that  $X$ is discrete.  To this end, choose $x \in X$;
we must find a neighbourhood of $x$ in $G$ which is disjoint from $X\setminus \{x\}$.
To ease notation, replace $X$ by $X - x,$ and $x$ by~$0$. 
Since $G_n \cap (X\setminus \{0\})$ is finite for each~$n$,  
we may inductively choose open subgroups $V_n$ of $G_n$
so that $U_n:= V_1 + \cdots + V_n$
is disjoint from 
$G_n \cap (X\setminus \{0\})$.
Then (following Remark~\ref{rem:inductive limits}~(2)),
we see that $U := \bigcup U_n$ is an open subgroup of~$G$ which is disjoint
from~$X\setminus \{0\}$.

Now suppose that $x$ is an element of $G \setminus X$.
Then $X' := X \cup \{x\}$ again satisfies the hypothesis of the
lemma, and thus, by what we have already proved, is discrete.
Thus we may find a neighbourhood of $x$ which is disjoint from~$X$,
and so $x$ does not lie in the closure of~$X$.  This proves that~$X$
is closed, as claimed.
\end{proof}

\begin{prop}
\label{prop:inductive limits}\leavevmode
\begin{enumerate}
\item Each of the canonical maps $G_n \hookrightarrow G$
is a closed embedding.
\item $G$ is Hausdorff, and if each~$G_n$  is furthermore complete,
then $G$ is complete.
\item Any compact subset of $G$ is contained in some~$G_n$.
\item If $H$ is a subgroup of~$G$,
and if we set  $H_n := G_n \cap H,$
then each of the natural maps $\varinjlim_n H_n \to H$
and
$\varinjlim_n G_n/H_n \to G/H$ 
is a topological isomorphism.
Furthermore,
$H$ is closed in $G$ if and only if $H_n$ %
is closed in $G_n$ for each~$n$.  
\end{enumerate}
\end{prop}
\begin{proof}
Fix~$n_0,$ and let
$U_{n_0}$ be 
any open subgroup of $G_{n_0}$.
If we set $V_n := U_{n_0} \cap G_n$ for $n \leq n_0$, and choose $V_n$
such that $V_n \cap G_{n-1} \subseteq V_{n-1}$ 
for $n > n_0,$ then we may follow Remark~\ref{rem:inductive limits}~(2)
to find an open subgroup $U \subset G$ such that
$U \cap G_{n_0} = U_{n_0}$.
Thus the topology of $G_{n_0}$  is induced by that 
of~$G$, proving~(1).

To see that $G$ is Hausdorff, we have to show that if $g \in G \setminus \{0\}$,
then we can find an open subgroup of $G$ not containing~$g$.
Choose $n_0$ such that $g \in G_{n_0} \setminus \{0\}$,
and let $U_{n_0}$ be 
any open subgroup of $G_{n_0}$
such that $g \not\in U_{n_0}$.
By the result of the preceding paragraph,
we may 
find an open subgroup $U$ of $G$ such that $G_{n_0} \cap U = U_{n_0}$.
Thus we see that $g \not\in U,$ as required.
This proves the first  claim of~(2).

Suppose now that each $G_n$ is complete,
and let $\widehat{G}$ denote the completion of~$G$.
The tautological morphism $G \to \widehat{G}$ is a topological embedding
(since~$G$ is Hausdorff, as we've just seen).  Since each $G_n$ is
complete, the composite embedding $G_n \hookrightarrow G \hookrightarrow
\widehat{G}$  realizes $G_n$ as a closed subgroup of~$\widehat{G}$.
If $G$ is not complete, then we may find $x \in 
\widehat{G} \setminus G$, and then,
since each $G_n$ is closed in $\widehat{G}$
(as we've just noted) 
and does not contain~$x$,
we may find an open subgroup 
$W_n$ of~$\widehat{G}$ such that $x + W_n$ is disjoint from~$G_n$;
we may furthermore choose that $W_n$ so that  $W_{n+1} \subseteq W_n$.

Now let $U := \sum_{n = 1}^{\infty} W_n \cap G_n$; by Remark~\ref{rem:inductive
limits}, this is an open subgroup of~$G$. 
By the construction of completions, its completion $\widehat{U}$ (i.e.\
its closure in~$\widehat{G}$) is an open subgroup of~$\widehat{G}$.
Since $G$ is dense in~$\widehat{G}$, we see that $(x + \widehat{U}) \cap G$ is non-empty,
and thus that $(x + \widehat{U}) \cap G_n$ is non-empty for some value of~$n$. 
On the other hand, we see that
\begin{multline*}
\widehat{U} \subseteq U + W_n = \sum_{i = 1}^{\infty} (W_i \cap G_i ) + W_n
\\
= \sum_{i \leq n} (W_i \cap G_i) + \sum_{i > n} (W_i \cap G_i) + W_n
\subseteq \sum_{i \leq n} G_n + \sum_{i> n} W_i  + W_n \\
= G_n + W_{n+1} + W_{n} = G_n + W_n,
\end{multline*}
Thus $(x + G_n + W_n) \cap G_n$ is non-empty, 
and so $(x + W_n) \cap G_n$ is non-empty, contradicting our choice of~$W_n$.
This shows that necessarily $G = \widehat{G}$, completing the proof of~(2).

Let  $X \subseteq G$ be a compact subset; since $G$ is Hausdorff (by what
we've just proved), $X$ is closed.  Suppose that $X \subsetneq G_n$ for
every~$n$; then we may find a sequence $x_n$ of elements of $X$ such
that $x_n \not\in G_n$.  Then $\bigcup_{n=1}^\infty\{x_n\}$ is an infinite
subset of~$X$ which  is closed and discrete in~$G$ (and thus in~$X$), by
Lemma~\ref{lem:discreteness}.  This contradicts the compactness
of~$X$, and so~(3) is proved.

Now let $H$ be a subgroup of~$G$.  By definition
of the induced topology,
the open subgroups of $H$ are of the form  $H \cap U$, with
$U$  an open subgroup of~$G$.  Write $U_n := G_n \cap U,$
so that $U_n$ is an open subgroup of $G$ and $U  = \bigcup_{n = 1}^{\infty} U_n.$
Then $H \cap U_n = H_n \cap U_n = H_n \cap U$ is an open subgroup of $H_n$,
and Remark~\ref{rem:inductive limits} shows that
$$\bigcup_{n=1}^{\infty} H_n \cap U_n  = H \cap \bigcup_{n=1} U_n = H \cap U$$
is also an open subgroup of~$\varinjlim_n H_n$, proving the first claim of~(4).
The topological isomorphism $\varinjlim_n G_n/H_n \iso G/H$
is immediately verified via a comparison of the universal mapping
properties of the source and target.

Certainly, if $H$ is closed in $G$ then $H_n$ is closed in
$G_n$ for each~$n$.  Conversely, suppose that 
$H_n$ is closed in~$G_n$ for each~$n$.
Then each quotient $G_n/H_n$ is again Hausdorff, and so part~(2)
shows that $\varinjlim_n G_n/H_n$ is Hausdorff.  Thus $H$, which is the
kernel of the canonical morphism $G = \varinjlim_n G_n \to \varinjlim_n G_n/H_n$,
is closed in~$G$, as claimed,
completing the proof of~(4).
\end{proof}

\subsection{Countably generated topological modules}
\label{subsec:cg top mod}
We first note the following simple lemma.

\begin{lemma}
\label{lem:cg Noeth}
If $A$ is a left Noetherian {\em (}not necessarily commutative{\em )} ring,
then
the category of countably generated $A$-modules
is an abelian subcategory of the category of all $A$-modules.
\end{lemma}
\begin{proof}
The only point that is not completely clear is closure under
passage to subobjects.  To see this, note that if $N$ is
countably generated, then we may write $N = \varinjlim_n N_n,$
with each $N_n$ finitely generated.  Then if $M$ is an $A$-submodule
of~$N$, we have $M = \varinjlim_n M_n,$ where $M_n := M \cap N_n.$
Since $A$ is left Noetherian by assumption, we see that each $M_n$ is
again finitely generated, and so $M$ is indeed countably generated.
\end{proof}

Suppose now that
$A$ is a compact Hausdorff left Noetherian 
 topological (possibly non-commutative) ring.
Any finitely generated $A$-module then
admits a unique compact Hausdorff topology --- its {\em canonical topology}
--- with
respect to which it becomes a topological  $A$-module,
and all morphisms between finitely generated $A$-modules
are continuous and strict (i.e.\ induce topological quotient maps from
their domains onto their images) with respect to the 
canonical topologies on their source and target.
(See~\cite[Prop.~2.1.3]{MR2667882} for a proof in a particular
case, which immediately generalizes.  The key point is that, given
a finite presentation
$A^{\oplus m} \to A^{\oplus n} \to M \to 0,$
the image of the first arrow is necessarily closed,
being  a compact  subset of a Hausdorff space, and so $M$ obtains a natural
compact Hausdorff quotient topology, which is easily checked to be independent
of the choice of finite presentation. Similar arguments with finite presentations
establish the claim about morphisms.)

In the cases we are interested in, the topology on $A$ will be defined
by a neighbourhood basis of $0$ consisting of open ideals,
and then one sees that the canonical topology on any finitely generated $A$-module
is similarly defined by a neighbourhood basis of $0$ consisting of open $A$-submodules.
In this case we say that $A$ and $M$ admit {\em linear} topologies,
and we assume from now on that the topology on~$A$ is indeed linear.

We let $A^{\oplus \NN}$ denote a countable direct sum of copies of~$A$.
Writing $A^{\oplus \NN} := \varinjlim_n A^{\oplus n}$ (with respect to the
evident transition maps),
we endow $A^{\oplus \NN}$ with its inductive limit linear
topology. %

\begin{prop}
\label{prop:canonical topology}\leavevmode
\begin{enumerate}
\item
If $M$ is a countably generated 
$A$-module, then $M$ admits a unique topology --- which we call
its {\em canonical topology} --- with respect to which some
{\em (}equivalently any{\em )} surjection
$A^{\oplus \NN} \to M$
is a topological quotient map.  Endowed with its canonical topology,
$M$ becomes a linear topological $A$-module. 
\item
Morphisms between countably generated $A$-modules are necessarily
continuous and strict with respect to 
the canonical topology on source and target.
\end{enumerate}
\end{prop}
\begin{proof}
If $M$ is a countably 
generated $A$-module, then we may write $M = \varinjlim_n M_n,$
where $\{M_n\}$ 
is an increasing sequence of finitely generated $A$-submodules of~$M$.
If $\{M'_n\}$ is any other increasing sequence of finitely generated
$A$-submodules of~$M$, then  evidently $\{M_n\}$ and $\{M'_n\}$
are mutually cofinal, in the sense that each $M_n$ is contained in some
$M'_{n'}$, and conversely.  
Thus, if we endow each $M_n$ with its canonical topology, the resulting
inductive limit linear topology on~$M$, coming from writing it
as the colimit of the~$M_n,$
is independent of the choice of the sequence~$\{M_n\}$.
We declare this topology to be the canonical topology of~$M$.

Since each $M_n$ admits a neighbourhood basis of $0$ consisting
of open $A$-submodules, it follows from Remark~\ref{rem:inductive limits}~(2)
that the same is true of~$M$.   One then easily deduces that $M$
is a topological $A$-module.

If $f:M \to N$ is a morphism of finitely generated $A$-modules,
then one may choose the $\{M_n\}$ and $\{N_n\}$ such that
$f$ restricts to morphisms $M_n \to N_n$.  These morphisms are necessarily
continuous with respect to the canonical topology on source and target,
and hence so is the morphism
$$f: M = \varinjlim_n M_n \to \varinjlim_n N_n = N.$$

If $f$ is injective, then we may in fact first choose the $N_n$,
and then set $M_n := M \cap N_n.$  Then $M_n$ is closed in $N_n$,
and so Proposition~\ref{prop:inductive limits}~(4)  shows that
$M \to N$ is a closed embedding.

If $f$ is surjective, then we may in fact first choose the $M_n$,
and then set $N_n := f(M_n)$.  
Again by Proposition~\ref{prop:inductive limits}~(4),
we see that $M \to N$ is a topological quotient map.

Altogether, we have now proved~(2), and~(1) is
a special case of the strictness
claim of~(2).
\end{proof}

\section{Representations of \texorpdfstring{$p$}{p}-adic analytic groups}\label{sec: reps
  of p adic analytic groups}

\subsection{Smooth
\texorpdfstring{$G$}{G}-representations and modules over the ring \texorpdfstring{$\cO[[G]]$}{OG}}%
\label{subsec:smooth  G reps}
We suppose that $G$ is a $p$-adic analytic group,
and fix a compact open subgroup $H$ of~$G$.
Then the usual completed group ring
$\cO[[H]]:=\varprojlim_{J \triangleleft H \textrm{
    open}}\cO[H/J]$ %
is a compact (two-sided) Noetherian linear topological ring by a
result of Lazard~\cite{MR209286}. %

\begin{defn}
  \label{defn: completed group ring of G}
We write\[\cO[[G]]:=\cO[G]\otimes_{\cO[H]}\cO[[H]]= \cInd_H^G \cO[[H]];
\] this is {\em a priori} a (left) representation of~$G$ with a compatible
right module structure over~$\cO[[H]]$, but following~\cite[\S 1]{MR3682662} (see also
\cite[\S 3]{MR4106887}), it actually has a natural ring
structure, uniquely determined by the requirement that each of $\cO[G]$
and $\cO[[H]]$ is a subring. The $\cO$-algebra $\cO[[G]]$ is
furthermore independent of the choice of~$H$ up to canonical
isomorphism.
\end{defn}

Any $\cO[[G]]$-module is in particular an~$\cO[G]$-module, and so can
be thought of as an $\cO$-linear $G$-representation.  Not every $\cO$-linear $G$-representation
admits a structure of $\cO[[G]]$-module compatible with its $\cO[G]$-action
(which is just to say that not every $\cO[G]$-module extends to a $\cO[[G]]$-module),
but an important fact is that smooth $G$-representations
(whose definition we now recall) do.

\begin{df}
\label{def:smooth}\leavevmode
\begin{enumerate}
\item
If $M$ is an $\cO$-linear $G$-representation,
we say that an element  $m \in M$ is {\em smooth}
if $m$  is $\cO$-torsion, and if $M$ is fixed by some
open subgroup of~$G$.
We let $M_{\sm}$ denote the subset of smooth elements of~$M$;
it is a $G$-invariant $\cO$-submodule of~$M$.
\item
A representation of $G$ on a $\cO$-module $M$ is called {\em smooth} 
if $M_{\sm} = M$.
\end{enumerate}
\end{df}

If $M$ is a smooth $G$-representation, then the $\cO[G]$-action on~$M$
extends canonically to an~$\cO[[G]]$-action.  Thus the abelian category
of smooth $\cO$-linear $G$-representations may be regarded as a full subcategory
of the category of $\cO[[G]]$-modules.
  (See~\cite[Lem.\ 3.5]{MR4106887}.)

\begin{remark}
\label{rem:smooth one}
By definition, an element of an $\cO$-linear  $G$-representation $M$ is smooth
if and only if it is smooth when $M$ is regarded as  an $H$-representation
for some (or, equivalently, any) compact open subgroup $H$ of~$G$.  
By definition of the topology on~$\cO[[H]]$,
we see that an element $m \in M$ is smooth if and only
if the cyclic $\cO[[H]]$-module generated by $m$ is discrete  
when regarded as a quotient of~$\cO[[H]]$.
\end{remark}

\begin{remark}
\label{rem:smooth two}
Just from the definitions, we see that if $M$ is an $\cO$-linear $G$-representation,
then~$M_{\sm}$  is the maximal smooth sub-$G$-representation  
of~$M$.
Furthermore, if $M$ itself is an $\cO[[G]]$-module, 
then $M_{\sm}$ (defined in terms  of the underlying $\cO[G]$-module structure)
becomes a sub-$\cO[[G]]$-module of~$M$, when endowed with its canonical~$\cO[[G]]$-module
structure.  (This follows from the fact that if $H$ is any compact $p$-adic
analytic group, e.g.\ a compact open subgroup of~$G$, 
then the augmentation ideal of $\cO[H]$ generates the augmentation ideal
of~$\cO[[H]]$.) 
\end{remark}

\begin{remark}
Most of the preceding discussion (the definitions of~$\cO[[G]]$ and of smooth $G$-representations,
and the realization of smooth representations as particular kinds of~$\cO[[G]]$-modules)
applies more generally to any topological group $G$
which contains a profinite open subgroup. 
However, many of the results 
that we prove below (for $p$-adic analytic~$G$) do not hold 
in this more general context.
(See Example~\ref{ex:non-analytic} below for an illustration of this.)
\end{remark} 

We let~$(\cO[[G]]\text{-}\Mod)^{\cg}$ denote
the full subcategory of~$\cO[[G]]\text{-}\Mod$ consisting of countably generated modules.

\begin{lemma}
\label{lem:cg mods}
$(\cO[[G]]\text{-}\Mod)^{\cg}$
is a Serre
subcategory of $\cO[[G]]\text{-}\Mod$,
which has enough  projectives.  Its projective objects are precisely
the countably generated projective~$\cO[[G]]$-modules.
\end{lemma}
\begin{proof}
The only non-trivial point to check regarding
the Serre subcategory claim is that submodules of countably generated~$\cO[[G]]$-modules
are again countably generated.  To see this,
note first that,
for any compact open subgroup~$H$ of~$G$, the ring $\cO[[G]]$ %
is countably
generated as a $\cO[[H]]$-module.  Thus any countably generated
$\cO[[G]]$-module is also countably generated as an $\cO[[H]]$-module.
Lemma~\ref{lem:cg Noeth}
shows that any $\cO[[G]]$-submodule is then countably generated over $\cO[[H]]$,
and thus also over~$\cO[[G]]$.

Evidently any countably generated $\cO[[G]]$-module is
a quotient of a countably generated free $\cO[[G]]$-module.  Since the
latter are projective even in~$\cO[[G]]\text{-}\Mod$, they are projective
in~$(\cO[[G]]\text{-}\Mod)^{\cg}.$
This proves the claim about enough projectives.
Any projective object of $(\cO[[G]]\text{-}\Mod)^{\cg}$ is then
a retract of a countably generated free module, and so is a countably
generated projective module.
\end{proof}

As already noted  in the proof of Lemma~\ref{lem:cg mods},
any countably generated
$\cO[[G]]$-module is also countably generated as an $\cO[[H]]$-module,
and so admits a canonical topology (independent of the choice of~$H$), as in Proposition~\ref{prop:canonical
topology}.

We note the following elementary lemma. %

\begin{lemma}
\label{lem:discrete}
A countably generated topological $\cO[[H]]$-module $M$ is 
discrete in its canonical topology if and only if $M$ 
is a smooth $H$-representation.
\end{lemma}
\begin{proof}
If $M$ is discrete, and if $m \in M$, then the sub-$\cO[[H]]$-module
of $M$ generated by~$m$ is discrete.  Remark~\ref{rem:smooth one} then
shows that $m$ is smooth.  Thus $M_{\sm} = M$.

Conversely, suppose that $M$ is smooth,
and choose a surjection $\cO[[H]]^{\oplus \mathbf{N}} \to M.$   For any $n \geq 0$,
the image of $\cO[[H]]^{\oplus n}$ factors through a discrete,
and so finite, quotient, and so has finite image $M_n \subseteq M$. We have
$M = \varinjlim_n M_n$, and by construction,
the canonical topology on $M$ is the inductive limit linear topology
on~$M$.   
Lemma~\ref{lem:discreteness}
implies that $M$ is discrete, as claimed.
\end{proof}

We now define a functor from the category
$(\cO[[G]]\text{-}\Mod)^{\cg}$
of countably generated
$\cO[[G]]$-modules to the category $\Pro \smG$ of pro-smooth $G$-representations
via %
$$M \mapsto \Pro(M) := \quoteslim{}{} M/U,$$
where $U$ runs over the codirected set of $\cO[[G]]$-submodules of~$M$
that are open with respect to the canonical topology on~$M$.
It follows from Lemma~\ref{lem:discrete}
that $M/U$ is in fact a smooth~$G$-representation.
The functoriality of the formation of $\Pro(M)$ follows from the fact that
if $f: M \to N$ is a morphism of countably generated $\cO[[G]]$-modules,
then $f$ is continuous with respect to the canonical topologies on its
source and target. Indeed, we then see that
 $f^{-1}(V)$  is an open $\cO[[G]]$-submodule
of $M$ whenever $V$ is an open $\cO[[G]]$-submodule of~$N$,
so that we have an induced morphism %
$$\Pro(M) := \quoteslim{} M/U \to \quoteslim{} M/f^{-1}(V) \to \quoteslim{} N/V
=: \Pro(N).$$

\begin{prop}
\label{prop:Pro}\leavevmode
\begin{enumerate}
\item
The functor $M \mapsto \Pro(M)$ is right exact.
\item
If $H$ is a compact open subgroup of~$G$,
and if $N$ is a countably generated
$\cO[[H]]$-module,
then $\cO[[G]]\otimes_{\cO[[H]]}N$
is acyclic for the left-derived functors
$L^i\Pro.$
\item
All countably generated smooth representations are
acyclic for the left-derived functors
$L^i\Pro.$
\end{enumerate}
\end{prop}
\begin{proof}
Let $0 \to M' \to M \to M'' \to 0$ be a short exact sequence
of countably generated $\cO[[G]]$-modules.  If $U$ is an open $\cO[[G]]$-submodule
of~$M$, we let $U' = M' \cap U,$ and let $U''$ denote the image of $U$ in~$M''$.
Since $M \to M''$ is an open $G$-equivariant mapping,
we see that, as $U$ runs over all $G$-invariant open 
submodules of $M,$ the images $U''$ run over all $G$-invariant
open submodules of~$M''$.
Thus we obtain a short exact sequence 
$$0 \to \quoteslim{} M'/U' \to 
\Pro(M) \to \Pro(M'')\to 0.$$
There is an evident epimorphism $\Pro(M') \to \quoteslim{} M'/U'$,
so that
\numequation
\label{eqn:right exact pro sequence}
\Pro(M') \to \Pro(M) \to \Pro(M'') \to 0
\end{equation}
is exact, as claimed.

We now prove~(2). Namely, if $N$ is a
countably generated $\cO[[H]]$-module, 
then we may find a resolution $P^{\bullet}\to N$ 
consisting of countably generated free $\cO[[H]]$-modules,
and
$\cO[[G]]\otimes_{\cO[[H]]} P^{\bullet}$ 
is then a resolution of $M:=\cO[[G]]\otimes_{\cO[[H]]}N$
by countable generated free $\cO[[G]]$-modules.  (Recall
that $\cO[[G]]$ is flat over~$\cO[[H]]$.)  Using
this resolution to compute the left derived functors $L^i\Pro,$
the claim of~(2)  follows by applying Lemma~\ref{lem: conditions for
  Pro exactness}~(2) below.

Finally, we turn to proving~(3).
Suppose to begin with that $M''$ is a countably generated smooth $G$-representation,
and let $M$ be a countably generated free~$\cO[[G]]$-module
admitting a surjection $M\to M''$.   
Lemma~\ref{lem:discrete} shows that the kernel~$M'$ of this surjection
is open in~$M$, and so Lemma~\ref{lem: conditions for
  Pro exactness}~(1) below shows that $L^1\Pro(M'')= 0;$
thus $L^1\Pro$ vanishes on countably generated smooth representations.

Now let $H$ be a compact open subgroup of~$G$. 
We may write find a surjection $M:= \cO[[G]]\otimes_{\cO[[H]]} N  \to M''$
for some smooth representation $N$ of $H$. 
The kernel $N'$ of this surjection is again countably generated and smooth,
and since~(2) shows that $M$ is $L^{\bullet}\Pro$-acyclic,
we find that $L^{i+1}\Pro(M'') \iso L^i\Pro(M')$ for $i\geq 1$.  An
obvious  dimension shifting argument,
taking into account the result of the preceding paragraph, gives~(3).
\end{proof}

It is not clear to us that~\eqref{eqn:right exact pro sequence}
is left-exact in general;
such left-exactness would be equivalent to the intersections $U'$
being cofinal in the collection of open $G$-invariant submodules
of~$M'$.  (Without the $G$-invariance condition, this would always
hold, since $M' \to M$ is necessarily a topological embedding.)
The following lemma gives some situations in which this {\em does} hold.

\begin{lem}
  \label{lem: conditions for Pro exactness}\eqref{eqn:right exact pro sequence}
is left-exact if either
\begin{enumerate}
\item $M'$ is open in $M$, or
\item the embedding $M'\into M$ is of the form \[M' :=
    \cO[[G]]\otimes_{\cO[[H]]} N' \hookrightarrow
    \cO[[G]]\otimes_{\cO[[H]]} N =: M,\] where $N' \hookrightarrow N$ is an
embedding of $\cO[[H]]$-modules (for some compact
open subgroup $H$ of~$G$).
\end{enumerate}
\end{lem}
\begin{proof}
  (1) If $M'$ is open in $M$, then any open $G$-invariant submodule
$U'$ of $M'$ is itself an open $G$-invariant submodule of~$M$,
and so we can set $U = U'$. Thus in this case~\eqref{eqn:right exact pro sequence}
is exact.  

(2) %
If $V'$ is an open $G$-invariant submodule of $M'$,  then $W':=N' \cap  V'$ is an open
$\cO[[H]]$-submodule of $N'$,
and so $U' := \cO[[G]]\otimes_{\cO[[H]]}  W'$  is  an  open $G$-invariant submodule of $M'$
contained  in $V'$;  thus submodules of this latter form 
are cofinal among all open $G$-invariant submodules of~$M'$.
Now since $N' \hookrightarrow N$ is necessarily a topological embedding,
we may find an open $\cO[[H]]$-submodule $W$ of $N$ such that
$W' = N' \cap  W.$ Then $U:= \cO[[G]]\otimes_{\cO[[H]]} W$  is
an open  $G$-invariant submodule of $M$  for which $M' \cap U =  U'$.
Thus in this case~\eqref{eqn:right exact pro sequence}
is again exact. 

\end{proof}

\subsection{Derived categories of  smooth representations}
Let $G$ be a  $p$-adic analytic group,
and let $\cO[[G]]$ be the ring of Definition~\ref{defn: completed group ring of G}.
We let $D(\cO[[G]])$  denote the stable $\infty$-category associated
to the abelian category of $\cO[[G]]$-modules.  It admits its standard $t$-structure,
with respect to which  it is both left and  right complete.  (This is true
for the stable $\infty$-category associated to the abelian category
of modules over any ring, as we noted in Section~\ref{sec:Complete t-structures}.)%

\begin{df} We let $D_{\sm}(\cO[[G]])$  denote the full subcategory
of $D(\cO[[G]])$ consisting  of objects whose cohomologies in every
degree are smooth $G$-representations.
\end{df}

The $t$-structure on $D(\cO[[G]])$ induces a $t$-structure on  $D_{\sm}(\cO[[G]])$,
whose heart is the abelian category $\smG$ of smooth $G$-representations.
The full subcategory $\smG$ of $\cO[[G]]\text{-}\Mod$
is closed under passage to subobjects, quotients, and colimits in~$\cO[[G]]$. 
Thus $\smG$ is 
again a  Grothendieck abelian category, and in particular admits enough injectives.

\begin{prop}
\label{prop:smooth equivalence}
There is a canonical equivalence $D(\smG) \iso  D_{\sm}(\cO[[G]])$.
\end{prop}
\begin{proof}
 As noted
above, the inclusion of smooth  $G$-representations
as a full subcategory of $\cO[[G]]$-modules realizes the former as a localizing
subcategory of the latter. We denote this functor by~$F$, and begin by showing
that~$F$ satisfies the equivalent conditions of Proposition~\ref{prop:D(heart)
  to  C}. To this end, it suffices to show that if $X$ is a smooth $G$-representation,
and  $Y$ is an injective smooth $G$-representation, then we may
find a smooth $G$-representation $Z$ admitting a surjection
$Z \to  X$ such  $\Ext^i_{\cO[[G]]}(Z,Y) = 0$ for $i>0$.

For this, choose a compact open subgroup $H$  of~$G$,
and set  $Z =\cInd_{H}^GX=\cO[[G]]\otimes_{\cO[[H]]}X,$
which admits a natural surjection to~$X$.
Since the forgetful functor
$\cO[[G]]\text{-}\Mod \to  \cO[[H]]\text{-}\Mod$ admits an exact left adjoint,
namely $\cO[[G]]\otimes_{\cO[[H]]} \text{--},$
it preserves injectives.
This same functor, restricted to $\smH$, where it can be interpreted as $\cInd_H^G$,
provides a  left adjoint to the forgetful functor
$\smG \to \smH;$
thus this latter functor also  preserves injectives.
Consequently $\Ext^i_{\cO[[G]]}(\cInd_H^G X,Y)  = \Ext^i_{\cO[[H]]}(X,Y)$,
with $Y$  injective in~$\smH$,
and so we reduce to verifying that
if $X$ and  $Y$ are objects of $\smH$ with  $Y$ injective,
then $\Ext^i_{\cO[[H]]}(X,Y) = 0$ for $i >  0$.  In short,
we have replaced $G$ in the original problem by its compact open
subgroup~$H$.
In this case, the required statement follows from Corollary~\ref{cor:preservation
of injectives} (applied in the context of Example~\ref{ex:compact analytic group ring}).

We next note that the inclusion $\smG \hookrightarrow \cO[[G]]\text{-}\Mod$
admits a right adjoint (as it must, by the adjoint functor theorem,
since colimits of smooth representations are smooth),
which admits a quite explicit description: it is
$$M \mapsto M_{\sm} := \colim_{H,n} M[\varpi^n]^H,$$
where in the colimit, $n$ ranges over positive  integers and $H$ ranges over 
open subgroups of~$G$.
By Proposition~\ref{prop:adjoints},  
we may then form the corresponding right derived functor
$R\Gamma_{\sm}: D(\cO[[G]]) \to  D(\smG)$,
which is right adjoint to~$F$.
Furthermore, 
$R\Gamma_{\sm}$ has finite cohomological amplitude.\footnote{If $G$
has dimension~$d$, and if $H$ is sufficiently small (so that  it is a uniform
pro-$p$-group), then  the derived functor of $M \mapsto M[\pi^n]^H$  as cohomological
amplitude~$d+1$, and so the same is true of $R\Gamma_{\sm}$ (which is the filtered
colimit of these various derived functors).} The proposition now follows from Proposition~\ref{prop:full faithfulness}.
\end{proof}

\begin{example}
\label{ex:non-analytic}
As the following (counter)example will show, Proposition~\ref{prop:smooth equivalence}
genuinely depends upon the hypothesis that $G$ is $p$-adic analytic, or at least on
various properties of the completed group ring $\cO[[G]]$ as a pro-Artinian ring which
need not hold for more general locally pro-$p$-groups.

For our example, we take $G = \prod_{i =  0}^{\infty} \Z/p\Z,$ with its natural profinite
topology.  To ease the comparison with the results we cite, we work with coefficients in
$k$ rather than~$\cO$, but the reader may easily convince themselves that
this alteration of the coefficients plays no essential role, and that
the discussion readily extends to the case of~$\cO$-coefficients.
The group ring of the product $G_n := \prod_{i = 0}^n \Z/p\Z$ admits
the description $k[G_n] = k[x_0,\ldots,x_n]/(x_0^p,\ldots,x_n^p),$
and so the completed group ring $k[[G]]$  admits the description
$$k[[G]] = k[[x_0,\ldots, x_n, \ldots]]/(x_0^p,\ldots,x_n^p,\ldots)
:=  \varprojlim_n k[[x_0,\ldots,x_n]]/(x_0^p,\ldots,x_n^p).$$
If we let $\m$ denote the maximal ideal (i.e.\ the augmentation ideal) of~$k[[G]]$,
then although $\m$ is open in $k[[G]]$, the natural profinite topology on 
$k[[G]]$ is weaker than its $\m$-adic topology.   
E.g.\  $\m/\m^2 = \varprojlim_n k\langle x_0,\ldots,x_n\rangle $ is infinite dimensional,
and the profinite topology on $k[[G]]$ endows $\m/\m^2$ with its natural profinite  topology,
rather than the discrete topology.

Let $R$ denote the polynomial subring $k[x_0,\ldots,x_n,\ldots]$ 
of~$k[[G]]$.
The category of smooth $G$-representations may be identified with the category of
$R$-modules in which each element is killed by one of the ideals $(x_i)_{i\geq n}$ for some $n\geq 0.$
This is an example of one of the categories denoted~$\cA$ that are introduced
in~\cite[Construction~1.1]{MR2875857}, and  \cite[Thm.~1.1]{MR2875857}
then shows that $D(\smG)$ is not left complete, while \cite[Rem.~1.2]{MR2875857}
shows that the formation of countable products in $D(\smG)$ is  
of unbounded cohomological amplitude
(illustrating the necessity of {\em some} bounded amplitude 
criterion in Proposition~\ref{prop:left completeness criterion}).
Thus the argument in the proof of Proposition~\ref{prop:smooth equivalence}
which extends the result from $D^+(\smG)$ to  $D(\smG)$  breaks down.

In fact even the canonical functor $F:D^+(\smG) \to D^+(k[[G]])$ 
is not fully faithful.  Indeed, the one-dimensional trivial representation of $G$  corresponds to
the $k[[G]]$-module $k := k[[G]]/\m$, and so $\Ext^1_{k[[G]]}(k,k)  =  \Hom_k(\m/\m^2,k)$
(the $k$-linear dual to $k$).
On the other hand $\Ext^1_{\smG}(k,k)$  corresponds to the proper subspace 
of $\Hom_k(\m/\m^2,k)$ consisting of functionals that are continuous when $\m/\m^2$
is endowed with its natural profinite topology.
\end{example}

\subsection{Compact and coherent objects}
We begin with some preliminaries based on the results
of~\cite[\S\S 2,3]{MR4106887}.

\begin{lem}
\label{lem:fg smooth}
If $M$ is a smooth $\cO$-linear $G$-representation,
then the following are equivalent:
\begin{enumerate}
\item $M$ is finitely generated as an $\cO[G]$-module.
\item $M$ is finitely generated as an $\cO[[G]]$-module.
\item There is a compact open subgroup $H$ of $G$,
a smooth representation of $H$ on a finitely generated torsion
$\cO$-module $N$, and a surjection $\cInd_H^G N \to  M$.
\item  For every compact open subgroup $H$ of $G$, there is
a smooth representation of $H$ on a finitely generated torsion
$\cO$-module $N$, and a surjection $\cInd_H^G N \to  M$.
\end{enumerate}
\end{lem}
\begin{proof}
Clearly (4) implies (3) implies (1) implies~(2).  Suppose  then that (2)
holds, and let $S \subset M$ be a finite generating set  for $M$
as an~$\cO[[G]]$-module.  If $H$ is any compact open subgroup of~$G$, and  if $N$ is the
$\cO[[H]]$-submodule of~$M$ generated by~$S$ (which coincides with the~$\cO[H]$-submodule
generated by $S$, since the smoothness of $M$ implies  that an open subgroup
of $H$ acts trivially on~$S$), then $\cInd_H^G = \cO[G]\otimes_{\cO[H]} N$
surjects onto~$M$,  proving~(4). 
\end{proof}

\begin{df}
\label{def:fg smooth reps}
We say that a smooth $G$-representation is finitely generated
if it satisfies the equivalent conditions of Lemma~\ref{lem:fg smooth}.
\end{df}

\begin{df}
  \label{def:finitely presented} A smooth representation of~$G$ on an
  $\cO$-module~$M$ is \emph{of finite presentation} if there is an
  exact sequence 
$$  M_1 \to M_2  \to M  \to 0$$
of smooth $G$-representations  with $M_1$ and $M_2$ finitely generated.
\end{df}

The following result is~\cite[Prop.~3.8]{MR4106887}.

\begin{lem}
\label{lem:fp characterization}
A smooth $G$-representation is  of finite presentation {\em (}in the
sense  of Definition~{\em \ref{def:finitely presented})} if and only
if it is finitely presented as an  $\cO[[G]]$-module. 
\end{lem}

We now have the following result.

\begin{prop}
\label{prop:compacts in smG}
The category $\smG$ is compactly generated.
The compact objects of $\smG$ are precisely the smooth
representations of finite presentation.
\end{prop}
\begin{proof}
We sketch the proof, which is standard.
To begin with, note that if $M$ is a smooth $\cO$-linear $G$-representation,
and $\{M_i\}$ is the (directed!) set of 
finitely generated subrepresentations of~$M$,
then evidently $M = \bigcup_i M_i = \colim M_i.$
For each value of~$i$,  choose a surjection $\cInd_{H_i}^G N_i \to M_i,$
with $H_i$ compact open in~$G$ and $N_i$ a smooth representation
of $H_i$ on a finitely generated torsion $\cO$-module.
If  $M'_i$ denotes the kernel of  this  surjection,  then,
as we have just seen, $M'_i  =  \colim (M'_i)_j,$ where  $(M'_i)_j$
runs over the finitely generated subrepresentations of~$M'_i$.
Then $M = \colim_{i,j} (\cInd_{H_i}^G N_i )/ (M'_i)_j,$
which exhibits  $M$  as  a colimit of representations
of finite presentation.  Thus, if  we prove  that
representations  of finite presentation  are compact,
we will have proved that $\smG$ is compactly generated.

Now if $M$ is compact, then applying its compactness property
to its representation as a colimit of representations
of finite presentation, we find that $M$ is a direct summand 
of a representation of finite presentation, and so is
easily seen itself to be a representation of finite presentation
(either directly, or via an application of Lemma~\ref{lem:fp characterization}).
Conversely, if $M$ is of finite presentation,
then Lemma~\ref{lem:fp characterization} shows that
$M$ is compact (even in the category of~$\cO[[G]]$-modules).
\end{proof}

We don't know in general whether the smooth $G$-representations of
finite presentation form an abelian subcategory
of~$\smG$, although it seems plausible that they do;
see Conjecture~\ref{conj:smooth reps are locally coherent},
as well as Remark~\ref{rem: is O[[G]] coherent?} for a
discussion of known results in this direction.
If they  do,
then  we find that $\smG$ is a locally  coherent abelian category,
and thus that the $t$-structure on $D(\smG)$ is coherent.
In this case we find (applying Lemma~\ref{lem:coherent objects}) that
the subcategory of coherent objects in $D(\smG)$ is precisely
the full subcategory $\Dfp^b(\smG)$ 
of (cohomologically) bounded complexes whose cohomologies are finitely presented.

\subsection{Complexes of \texorpdfstring{$\cO[[G]]$}{OG}-modules as pro-smooth complexes}
The  formal Pro-category $\Pro \smG$ admits enough projectives. 
Thus we may form  its associated stable $\infty$-category  $D^-(\Pro \smG)$,
with its natural $t$-structure.

We may also from the formal Pro-$\infty$-category $\Pro D^b(\smG)$.  The $t$-structure
on $D^b(\smG)$ induces a left complete $t$-structure on $\Pro D^b(\smG)$
(this is the dual statement to Proposition~\ref{prop:Ind t-structures}),
whose heart is equivalent to $\Pro \smG$.  
Thus (by the dual version of Proposition~\ref{prop:D(heart) to C},
i.e.~\cite[Prop.~1.3.3.7]{LurieHA})
there is a canonical functor
\numequation
\label{eqn:pro functor}
D^-(\Pro \smG) \to \Pro D^b(\smG).
\end{equation}

Recall that $(\cO[[G]]\text{-}\Mod)^{\cg}$ denotes the abelian category of countably
generated $\cO[[G]]$-modules. 
By Lemma~\ref{lem:cg mods}, this is an abelian category
with enough projectives (namely, the countably generated free~$\cO[[G]]$-modules),
and so we may consider its stable 
$\infty$-category of bounded above complexes~$D^-\bigl((\cO[[G]]\text{-}\Mod)^{\cg}\bigr).$
Since in fact (by Lemma~\ref{lem:cg mods}) 
$(\cO[[G]]\text{-}\Mod)^{\cg}$ is a Serre subcategory of~$\cO[[G]]\text{-}\Mod$,
it also makes sense to consider~$D_{\cg}^-(\cO[[G]])$,
the full sub-stable $\infty$-category of $D^-(\cO[[G]])$ consisting
of complexes with countably generated cohomologies. 
The canonical $t$-structure on $D^-(\cO[[G]])$ induces a $t$-structure
on~$D^-_{\cg}(\cO[[G]]),$
whose heart is precisely~$(\cO[[G]]\text{-}\Mod)^{\cg}$.

\begin{lemma}
\label{lem:cg equivalence}
There is a canonical equivalence  
\numequation
\label{eqn:cg equivalence}
D^-\bigl((\cO[[G]]\text{-}\Mod)^{\cg}\bigr)
\iso 
D_{\cg}^-(\cO[[G]]).
\end{equation}
\end{lemma}
\begin{proof}
This follows from the dual version of Proposition~\ref{prop:D(heart) to C},
i.e.~\cite[Prop.~1.3.3.7]{LurieHA}, 
since the projective objects in $(\cO[[G]]\text{-}\Mod)^{\cg}$
are precisely the countably generated projective $\cO[[G]]$-modules,
which can be used to compute Exts in~$\cO[[G]]\text{-}\Mod$.
\end{proof}

Now left deriving the functor of 
Proposition~\ref{prop:Pro} 
gives a  functor
$$L\Pro: D^-\bigl((\cO[[G]]\text{-}\Mod)^{\cg}\bigr)
\to D^-(\Pro \smG).$$
Precomposing with an inverse to the equivalence of~\eqref{eqn:cg equivalence},
and postcomposing with~\eqref{eqn:pro functor},
we obtain a functor
\numequation
\label{eqn:functor to pro complexes}
D_{\cg}^-(\cO[[G]])
\to \Pro D^b(\smG).
\end{equation}

Let $D_{\cg}^b(\smG)$ denote the full subcategory
of $D^b(\smG)$ consisting of objects whose cohomologies are countably generated,
and
let $D_{\smcg}^b(\cO[[G]])$ denote the full subcategory of $D^b(\cO[[G]])$ 
consisting of complexes whose cohomologies are countably generated and smooth.
The equivalence of Proposition~\ref{prop:smooth equivalence}
evidently induces an equivalence
\numequation
\label{eqn:cg smooth equivalence two}
D_{\cg}^b(\smG) \iso D_{\smcg}^b(\cO[[G]]).
\end{equation}
It follows from Proposition~\ref{prop:Pro}~(3)
that the composite of the equivalence~\eqref{eqn:cg smooth equivalence two}
with (the restriction to $D_{\smcg}^b(\cO[[G]])$ of)
the functor~\eqref{eqn:functor to pro complexes}
is equivalent to the composite of 
the canonical functors
$$D_{\cg}^b(\smG) \hookrightarrow D^b(\smG) \hookrightarrow\Pro D^b(\smG).$$ 

In summary, we have the functor $D_{\cg}^-(\cO[[G]]) \to \Pro D^b(\smG)$
of~\eqref{eqn:functor to pro complexes},
which when restricted to $D_{\cg}(\smG)^b$ (thought of as
a subcategory of $D_{\cg}^-(\cO[[G]])$ via the equivalence~\eqref{eqn:cg smooth equivalence two})
is equivalent to the evident inclusion.

\section{Derived moduli stacks of group
representations} %
\label{app:derived moduli}In this appendix we give a brief
introduction to and motivation of derived moduli stacks of Galois
representations. We refer to~\cite[\S2]{zhu2020coherent} for the
details.

Let \[\Gamma := \langle g_1,\ldots,g_m \,  | \, r_1 = \cdots =r_n = 1\rangle\]
be a finitely presented group, 
with generators $g_1,\ldots,g_m$
and relations $r_1,\ldots, r_n$. (The motivational example to bear in mind is the
Galois group of a finite extension of local or global fields; in the
body of the notes we will consider instead various topologically
finitely presented profinite groups, given for example by the absolute
Galois groups of $p$-adic local fields, but for the sake of motivation
the case of a literally finite group suffices.)

Then giving a representation $\Gamma  \to
\GL_d(A)$ (for some ring $A$) amounts to giving elements $x_1,\ldots,
x_m \in \GL_d(A)$ satisfying the relations~$r_j$. 
Let $V_{\Gamma}$  be the closed subscheme of $(\GL_d)^m$ classifying
$m$-tuples satisfying the relations $r_j$.  
There is an action of $\GL_d$ on $V_{\Gamma}$ via simultaneous conjugation, and
the quotient stack
$\cX_{\Gamma} := [ V_{\Gamma}/\GL_d]$ %
is then the moduli stack of 
$d$-dimensional representations of~$\Gamma$
(over arbitrary commutative rings,
as we have written it here, although of course we could work over
any given base ring $\cO$, and then consider representations over $\cO$-algebras).
It is an algebraic stack (over $\Z$, or whichever other
base ring $\cO$ we choose).

If we fix a particular representation $\rho: \Gamma \to \GL_d(l)$
(defined over some  field~$l$),
then this representation gives rise to a $l$-valued point of $\cX_{\Gamma}$,
and the lifting ring $R_{\rho}$ of $\rho$ is a versal ring
to $\cX_{\Gamma}$ at the point~$\rho$.
The tangent space to $\cX_{\Gamma}$ at~$\rho$ can be described
in terms of group cohomology as $H^1(\Gamma, \Ad \rho),$
while $H^0(\Gamma, \Ad \rho)$
computes the infinitesimal  automorphisms of~$\rho$.  
Passing to higher degree cohomology, one knows that
elements of $H^2(\Gamma, \Ad \rho)$ 
give rise to obstructions
to  higher  order  deformations of $\rho$.
However, 
this obstruction theory %
is an  extra piece of  data --- it is not determined intrinsically
by the stack $\cX_{\Gamma}$ itself.
For example, in general one can't 
``compute'' $H^2(\Gamma, \Ad \rho)$
in terms of~$\cX_{\Gamma}$ (unlike $H^0$ and $H^1$).

This comes up concretely when one uses group cohomology to describe
the lifting ring~$R_\rho$: if
we  consider liftings over $l$-algebras,
and write  $h^i := \dim_l H^i(\Gamma,  \Ad \rho)$,
then we can write $R_{\rho}$ as a quotient
$l[[x_1,\ldots, x_{h^1 -h_0+ d^2}]]/(f_1,\ldots,f_{h^2})$.
When $\cX_{\Gamma}$ has the ``expected dimension''
$ -h^0+ h^1- h^2$ at~$\rho$, the $f_i$ must form a regular 
sequence  in 
$l[[x_1,\ldots, x_{h^1 -h_0+ d^2}]],$
and $R_{\rho}$ is a local complete intersection.  But if
$\cX_{\Gamma}$
has larger-than-expected dimension, then the ring-theoretic structure
of $R_{\rho}$ is less clear, and less clearly related to the quantity~$h^2$.

One reason, then, for considering a {\em derived} analog of $\cX_{\Gamma}$
is that it tightens the relationship between group cohomology and
the infinitesimal structure of $\cX_{\Gamma}$ at its points~$\rho$. %
The precise definition of the derived representation stack
is somewhat involved (see~\cite[\S2]{zhu2020coherent}), 
but it proceeds roughly as follows:
as is  usual  in a derived setting,
rather than just taking  into account the relations~ $r_i$   appearing
in~ $\Gamma$, we consider  the ``higher relations'', i.e.\  the
relations between the relations and so on.  In the non-abelian 
setting of not-necessarily-commutative groups, this  can be 
expressed formally by finding a simplicial resolution $\Gamma_{\bullet}$
of $\Gamma$ by  free groups (or, in fact, for technical reasons,
by free monoids). %
The schemes $V_{\Gamma_{n}}$  
classifying homomorphisms $\Gamma_{n} \to \GL_d(A)$ (for coefficient
rings $A$) are then of the form  $\Spec R_{\Gamma_{n}}$
for certain rings $R_{\Gamma_{n}}$ that are localizations
of polynomial rings
(since the $\Gamma_{n}$
are free monoids), and which are themselves organized into a simplicial commutative
ring~$R_{\Gamma_{\bullet}}$.  
We may consider the corresponding derived affine scheme
$V_{\Gamma_{\bullet}}$, and form the derived algebraic stack
$\cX_{\Gamma_{\bullet}} := [ V_{\Gamma_{\bullet}}/\GL_d],$
whose underlying classical substack coincides with $\cX_{\Gamma}$.
In a suitable $\infty$-categorical framework (see
e.g.\ the discussion of~\cite[\S2]{zhu2020coherent})  the derived
stack $\cX_{\Gamma_{\bullet}}$ is well-defined independently
of the  choice of resolution $\Gamma_{\bullet}$.

There is  a standard  useful analogy, in which one
compares the passage from classical to derived rings (or, in the 
terminology of \cite{cesnavicius2021purity}, which we employ, ``animated rings'') to
the passage from reduced rings to arbitrary rings that occurs in passing from varieties to schemes,
and
the passage from classical to derived schemes or stacks 
to the passage from varieties to not-necessarily-reduced
schemes. Even for a classical  variety over a field~$l$, it is interesting
to consider its  points over non-reduced rings: e.g.\ its tangent spaces
are computed via
$l[\varepsilon]/(\varepsilon^2)$-valued points.
Similarly,
for a derived stack $\cX$ (even one which is actually classical),
there are non-classical analogues of
$l[\varepsilon]/(\varepsilon^2)$ whose points in a classical or derived stack $\cX$
determine
the obstruction theory of $\cX$, and using such points,
we can extend  the classical computation of tangent spaces
to compute the cotangent complex $\mathbf L_{\cX,x}$ of~$\cX$
at points  $x$ of~$\cX$.
In particular, when $\cX = \cX_{\Gamma_{\bullet}}$, 
one finds that \[\mathbf L_{\cX_{\Gamma_{\bullet}},\rho} =
C^{\bullet}(\Gamma,\Ad \rho)^*[-1]\] (a shift of the dual 
to  the cochain complex that computes the group 
cohomology $H^{\bullet}(\Gamma, \Ad \rho)$) --- see e.g.\ \cite[\S2.2]{zhu2020coherent}.
This extends the geometric interpretations
of $H^0$ and $H^1$ 
in the classical context. 

This is especially useful in cases when
$H^i(\Gamma,\Ad \rho)   =  0$ for  $i  > 2$, since
then the cohomology of $\mathbf L_{\cX_{\Gamma_{\bullet}},\rho}$  
appears only  in degrees $[-1,1]$. This implies that  $\cX_{\Gamma_{\bullet}}$
is {\em quasi-smooth} (i.e.\ derived l.c.i.) in a neighbourhood of~$\rhobar$.
Concretely, this  means that there is a formally smooth morphism
\[\Spec l[[x_1,\ldots, x_{h^1 -h_0+ d^2}]]/(f_1,\ldots,f_{h^2})
\to  \cX_{\Gamma_{\bullet}},\]
taking the closed point  of the domain to~$\rho$,
{provided that} the domain is understood in 
a {derived} sense:  if the $f_i$ do not form a regular sequence,
then the ring
$l[[x_1,\ldots, x_{h^1 -h_0+ d^2}]]/(f_1,\ldots,f_{h^2})$
will be a non-classical animated ring (its higher homotopy  groups
will be computed by the Koszul complex of the sequence $(f_i)$).
We then see that $\cX_{\Gamma_{\bullet}}$ is  in fact classical in a neighbourhood
of $\rho$ if and only if the underlying classical
stack $\cX_{\Gamma}$  is  of the expected dimension
$-h_0 +h_1 - h_2 = \chi(\mathbf L_{\cX_{\Gamma_{\bullet},\rho}})$ 
in such a neighbourhood.

In Sections~\ref{sec: categoric LL for l not p} and~\ref{sec: global stacks and cohomology of Shimura
  varieties} we apply these ideas when $\Gamma$
is (at least closely related to) a local or global absolute Galois group.
(Their extension to the context of  profinite
groups is treated in \cite[\S 2.4]{zhu2020coherent}.)
To this end, we recall that local Galois cohomology
always vanishes in degrees $> 2$, and that  the same is true
for global Galois cohomology, provided we work in odd residue characteristic
(as we  always will).

\bigskip

\emergencystretch=3em

\printbibliography

\bigskip

\emergencystretch=3em

\let\theequation\originaltheequation
\end{document}
